\documentclass{memo-l}

\usepackage{amsmath,amsfonts,amsthm,amssymb}
\usepackage{graphicx,psfrag,overpic}
\usepackage{mathrsfs}
\usepackage{bm}
\usepackage[all,cmtip,line]{xy}
\usepackage{array}

\usepackage[usenames,dvipsnames]{xcolor}

\usepackage{color}
\usepackage{geometry}

\newtheorem{theorem}{Theorem}[chapter]
\newtheorem{lemma}[theorem]{Lemma}

\theoremstyle{definition}
\newtheorem{definition}[theorem]{Definition}

\theoremstyle{remark}
\newtheorem{remark}[theorem]{Remark}

\newtheorem{problemL}[theorem]{Problem}
\newtheorem{proposition}[theorem]{Proposition}
\newtheorem{corollary}[theorem]{Corollary}

\newtheorem*{propositionnolabel}{Proposition}
\theoremstyle{definition}
\newtheorem*{proofdist}{Proof of Proposition {\rm \ref{proposition-distance}}}

\numberwithin{section}{chapter}
\numberwithin{equation}{chapter}

\newcommand \gam{\gamma}
\newcommand \R{\mathbb{R}}
\newcommand \Om{\Omega}
\newcommand \der{\partial}
\newcommand \vphi{\varphi}

\newcommand \mcl{\mathcal}
\newcommand \Gam{\Gamma}
\newcommand \alp{\alpha}
\newcommand \tx{\text}
\newcommand \til{\tilde}

\newcommand \ol{\overline}

\newcommand \eps{\varepsilon}

\newcommand \Shock{\Gamma_{{\rm{shock}}}}
\newcommand \Wedge{\Gamma_{{\rm{wedge}}}}

\newcommand \sonic{\Gamma_{\rm{sonic}}}
\newcommand \shock{\Gamma_{\rm{shock}}}

\newcommand \wedgetip{P_{\beta}}

\newcommand \lefttop{P_1}
\newcommand \righttop{P_2}

\newcommand \leftbottom{P_4}
\newcommand \rightbottom{P_3}

\newcommand \leftshockop{S_{0}}
\newcommand \leftsonicop{\Gamma_{{\rm{sonic}}}^{0}}
\newcommand \leftvecop{\bm{e}_{S_{0}}}
\newcommand \leftvphiop{\varphi_{0}}
\newcommand \leftuop{u_0}
\newcommand \leftvop{v_0}

\newcommand \oOmop{\Om_0}
\newcommand \leftcop{c_0}
\newcommand \leftPop{Q_0}

\newcommand \rightshockop{S_{1}}
\newcommand \rightsonicop{\Gamma_{{\rm{sonic}}}^{1}}
\newcommand \rightvecop{\bm{e}_{S_{1}}}
\newcommand \rightvphiop{\varphi_{1}}

\newcommand \rightrhoop{\rho_1}

\newcommand \nOmop{\Om_1}
\newcommand \netaop{\xi_2^{(1)}}
\newcommand \rightcop{c_1}
\newcommand \rightPop{Q_1}

\newcommand \nDop{\mcl{D}^1}
\newcommand \oDop{\mcl{D}^0}

\newcommand \leftshock{S_{\mathcal{O}}}
\newcommand \leftsonic{\Gamma_{{\rm{sonic}}}^{\mathcal{O}}}
\newcommand \leftvec{\bm{e}_{S_{\mathcal{O}}}}
\newcommand \leftvphi{\varphi_{\mathcal{O}}}
\newcommand \leftu{u_{\mathcal{O}}}
\newcommand \leftv{v_{\mathcal{O}}}

\newcommand \rightshock{S_{\mcl{N}}}
\newcommand \rightsonic{\Gamma_{{\rm{sonic}}}^{\mcl{N}}}
\newcommand \rightvec{\bm{e}_{S_{\mathcal{N}}}}

\newcommand \rightvphi{\varphi_{\mcl{N}}}
\newcommand \ivphi{\varphi_{\infty}}
\newcommand \vphib{\varphi_{\beta}}

\newcommand \iv{v_{\infty}}
\newcommand \iu{u_{\infty}}
\newcommand \ic{c_{\infty}}

\newcommand \irho{\rho_{\infty}}

\newcommand \leftphi{\phi_{\mathcal{O}}}

\newcommand \iphi{\phi_{\infty}}

\newcommand \leftc{c_{\mathcal{O}}}
\newcommand \leftrho{\rho_{\mathcal{O}}}
\newcommand \rightc{c_{\mcl{N}}}
\newcommand \rightrho{\rho_{\mcl{N}}}
\newcommand \ik{k_{\infty}}
\newcommand \leftk{-\iv \xi_2^{(\beta)}}

\newcommand \nrho{\rho_{\mcl{N}}}
\newcommand \orho{\rho_{\mcl{O}}}
\newcommand \iM{M_{\infty}}
\newcommand \iq{q_{\infty}}
\newcommand \oM{M_{\mcl{O}}}
\newcommand \oxi{{\xi_1^{\mcl{O}}}}
\newcommand \nxi{{\xi_1^{\mcl{N}}}}

\newcommand \neta{{\xi_2^{\mcl{N}}}}
\newcommand \oeta{{\xi_2^{\mcl{O}}}}

\newcommand \oL{L_{\mcl{O}}}
\newcommand \iL{L_{\infty}}

\newcommand \otheta{\theta_{\mcl{O}}}
\newcommand \itheta{\theta_{\infty}}

\newcommand \nD{\mcl{D}^{\mcl{N}}}
\newcommand \oD{\mcl{D}^{\mcl{O}}}
\newcommand \nOm{\Omega^{\mcl{N}}}
\newcommand \oOm{\Omega^{\mcl{O}}}
\newcommand \iOm{\Omega^{\infty}}

\newcommand \Oi{O_{\infty}}
\newcommand \Oo{O_{\mcl{O}}}
\newcommand \Onormal{O_{\mcl{N}}}
\newcommand \Qo{Q_{\mcl{O}}}
\newcommand \Qn{Q_{\mcl{N}}}

\newcommand \nLambda{\Lambda^{\mcl{N}}}
\newcommand \oLambdab{\Lambda^{\mcl{O}}_{\beta}}
\newcommand \go{g_{\mcl{O}}}
\newcommand \gn{g_{\mcl{N}}}

\newcommand \betac{\beta_{\rm s}}
\newcommand \Lbeta{\Lambda_{\beta}}

\newcommand \cone{\mathrm{Cone}^0(\leftvec, \rightvec)}

\newcommand \leftchi{\chi_{\mcl{O}}}
\newcommand \rightchi{\chi_{\mcl{N}}}

\newcommand \gshock{\mathfrak{g}_{\rm sh}}

\newcommand \tgshock{\tilde{\mathfrak{g}}_{\rm sh}}

\newcommand \gso{\mathfrak{g}_{\mcl{O}}}

\newcommand \gsn{\mathfrak{g}_{\mcl{N}}}

\newcommand \hgshock{\hat{\mathfrak{g}}_{\rm sh}}

\newcommand \un{u^{\rm (norm)}}

\newcommand \leftsonicdelta{\Gam^{\mcl{O},\delta_0}_{\rm sonic}}
\newcommand \rightsonicdelta{\Gam^{\mcl{N},\delta_0}_{\rm sonic}}

\newcommand \tleftchi{\tilde{\chi}_{\mcl{O}}}
\newcommand \trightchi{\tilde{\chi}_{\mcl{N}}}

\newcommand \iter{\mcl{Q}^{\rm iter}}
\newcommand \Q{\mcl{Q}^{\rm iter}}

\newcommand \Itr{\mcl{I}}
\newcommand \Itrf{\mcl{I}_1}

\newcommand \dershock{\der_{\rm sh}\iter}
\newcommand \derwedge{\der_{\rm w}\iter}

\newcommand \D{\mcl{D}}

\newcommand \fshock{f_{\rm sh}}
\newcommand \fshpolar{f_{\Oi,{\rm sh}}}
\newcommand \fshockn{\hat f_{\mcl{N},{\rm sh}}}
\newcommand \fshocko{\hat f_{\mcl{O},{\rm sh}}}

\newcommand \qi{q_{\infty}}
\newcommand \leftq{q_{\mcl{O}}}

\newcommand \dsonic{d_{\rm so}}

\newcommand \sbeta{s_{\beta}}
\newcommand \Lb{L_{\beta}}
\newcommand \Qbeta{Q^{\beta}}
\newcommand \cbeta{c_{\beta}}

\newcommand \mfrakO{O}

\newcommand \Kext{\mcl{K}^{\rm ext}}

\newcommand \ee{{\bm e}}

\newcommand \q{{\bf q}}
\newcommand \xib{\xi_1^{P_\beta}}

\newcommand \zetao{\zeta_{\mcl{O}}}
\newcommand \zetan{\zeta_{\mcl{N}}}

\newcommand \lsonic{l_{\rm so}}
\newcommand \betad{\beta_{\rm d}}

\newcommand \tdomain{\mathfrak{D}_{\theta_{\rm w}}} 
\newcommand \bdomain{\Lambda_{\beta}} 

\newcommand \M{\mcl{M}}
\newcommand \hfshock{\hat{f}_{\rm sh}}

\newcommand \leftshockseg{S_{\mathcal{O},{\rm seg}}}
\newcommand \rightshockseg{S_{\mathcal{N},{\rm seg}}}

\newcommand \leftshocksegop{S_{0,{\rm seg}}}
\newcommand \rightshocksegop{S_{1,{\rm seg}}}

\newcommand \dersonic{\der_{\rm so}\iter}
\newcommand \leftch{\hat{c}_{\mcl{O}}}

\newcommand \bmxi{{\bm{\xi}}}
\newcommand \xin{\xi_1}
\newcommand \etan{\xi_2}

\newcommand \om{\omega}
\newcommand \bPhi{\bar{\Phi}}
\newcommand \betadet{\beta_{\rm d}^{(\iv)}}
\newcommand \betasonic{\beta_{\rm s}^{(\iv)}}
\newcommand \cinfty{c_{\infty}}
\newcommand \Npolar{\mcl{N}_{p}}
\newcommand \rx{{\rm{\bf x}}}
\newcommand \ry{\rm{\bf y}}

\newcommand \coneclosure{{\mathrm{Cone}(\leftvec, \rightvec)}}
\newcommand{\xxi}{{\boldsymbol \xi}}
\newcommand{\betaa}{{\boldsymbol \beta}}
\newcommand{\pb}{{\bf p}}

\newcommand{\rightch}{\hat{c}_{\mcl{N}}}
\newcommand{\leftchnew}{{c}^*_{\mcl{O}}}
\newcommand{\Pmax}{P_{\rm{max}}}

\newcommand{\dd}{{\rm d}}

\makeindex

\begin{document}

\frontmatter
\title[Prandtl-Meyer Reflection Configurations]
{Prandtl-Meyer Reflection Configurations,\\
Transonic Shocks,\\
and Free Boundary Problems}

\author{Myoungjean Bae}
\address{Department of Mathematical Sciences,
Korea Advanced Institute of Science and Technology (KAIST),
291 Daehak-ro, Yuseong-gu, Daejeon, 43141, Korea}
\curraddr{}
\email{mjbae@kaist.ac.kr}

\author{Gui-Qiang G. Chen}
\address{Mathematical Institute, University of Oxford,
Andrew Wiles Building, Radcliffe Observatory Quarter, Woodstock Road,
Oxford, OX2 6GG, UK}
\curraddr{}
\email{chengq@maths.ox.ac.uk}

\author{Mikhail Feldman}
\address{Department of Mathematics\\
         University of Wisconsin\\
         Madison, WI 53706-1388, USA}
\curraddr{}
\email{feldman@math.wisc.edu}
\date{}

\subjclass[2010]{Primary: 35M10, 35M12, 35R35, 35B65, 35L65, 35L70, 35J70, 76H05, 35L67, 35B45, 35B35, 35B40, 35B36, 35B38;
Secondary: 35L15, 35L20, 35J67, 76N10, 76L05, 76J20, 76N20, 76G25}

\keywords{Prandtl-Meyer reflection, Prandtl conjecture,
supersonic flow, unsteady flow, steady flow, solid wedge,
nonuniqueness,
weak shock solution, strong shock solution,
stability, self-similar, global solution,
transonic flow, transonic shock,
sonic boundary, free boundary, existence, regularity,
long-time asymptotics, detachment angle,
admissible solutions,
elliptic-hyperbolic mixed type,
degenerate elliptic equation, nonlinear PDEs,
monotonicity, {\it a priori} estimates, uniform estimates}

\begin{abstract}
We are concerned with the Prandtl-Meyer reflection configurations
of unsteady global solutions for supersonic flow
impinging upon a symmetric solid wedge.
Prandtl (1936) first employed the shock polar analysis to
show that there are two possible steady configurations:
the steady weak shock solution and the steady strong shock solution,
when a steady supersonic flow impinges upon
the solid wedge -- the half-angle of which is less than a critical angle ({\it i.e.}, the detachment angle),
and then conjectured that
the steady weak shock solution is physically admissible
since it is the one observed experimentally.
The fundamental issue of whether one or both of the steady weak and strong
shocks are physically
admissible has been vigorously debated over the past eight decades
and has not yet been settled in a definitive manner.
On the other hand, the Prandtl-Meyer reflection configurations are
core configurations in the structure of global entropy solutions
of the two-dimensional Riemann problem,
while the Riemann solutions themselves are local building blocks
and determine local structures, global attractors, and large-time asymptotic
states of general entropy solutions of
multidimensional hyperbolic systems
of conservation laws.
In this sense, we have to understand the reflection configurations
in order to understand fully the global entropy solutions of two-dimensional
hyperbolic systems of conservation laws, including the admissibility issue
for the entropy solutions.
In this monograph, we address this longstanding open issue
and present our analysis to establish
the stability theorem for the steady weak shock solutions as
the long-time asymptotics of
the Prandtl-Meyer reflection configurations
for unsteady potential flow
for all the physical parameters up to the detachment angle.
To achieve these, we first
reformulate the problem as a free boundary problem involving transonic shocks
and then obtain appropriate monotonicity properties
and uniform {\it a priori} estimates for admissible solutions,
which allow us to employ
the Leray-Schauder degree argument to complete the theory
for all the physical parameters up to the detachment angle.
\end{abstract}
\maketitle

\tableofcontents
\mainmatter

\chapter{Introduction}

\numberwithin{equation}{chapter}
\numberwithin{figure}{chapter}

We are concerned with
unsteady global solutions for supersonic flow impinging upon a solid ramp,
which can equivalently be regarded as portraying the symmetric gas flow impinging
upon a solid wedge (by symmetry).
When a steady supersonic flow impinges upon the solid wedge  -- the half-angle $\theta_{\rm w}$
of which is less than a critical
angle ({\it i.e.}, the detachment angle $\theta_{\rm d}$),
Prandtl first employed the shock polar analysis to show that there are two possible
steady configurations:
the steady weak shock reflection with supersonic or subsonic downstream flow
(determined by the wedge angle that is less or larger than the sonic angle $\theta_{\rm s}<\theta_{\rm d}$)
and the steady strong shock reflection
with subsonic downstream flow, both of which satisfy the entropy
conditions, provided that no additional conditions are assigned downstream;
see Courant-Friedrichs \cite{CF}, von Neumann \cite{Neumann}, and Prandtl \cite{Prandtl}.

A fundamental issue is whether one or both of the steady weak and strong shocks are physically
admissible. This has been debated vigorously over the past eight decades and has not yet been settled
in a definitive manner ({\it cf.} \cite{CF,Dafermos,Liu,Neumann,Serre}).
On the basis of experimental and numerical evidence, there are strong indications to show,
as Prandtl conjectured (see \cite{Busemann,Meyer,Prandtl}),
that it is
the steady weak shock
solution that is physically admissible as the long-time asymptotics of
the Prandtl-Meyer reflection configurations.

Furthermore, the Prandtl-Meyer reflection configurations are
solutions of the lateral Riemann problem (Problem 2.6 below),
and are core configurations in the structure of global entropy solutions
of the two-dimensional Riemann problem for hyperbolic conservation laws.
On the other hand, the Riemann solutions are building blocks
and determine local structures, global attractors, and large-time asymptotic
states of general entropy solutions
of multidimensional hyperbolic systems of conservation laws
(see \cite{CCY2,CCY3,CH,CF2,GlimmMajda,KTa,LaxLiu,LZY,SCG,Zhe}
and the references cited therein).
Consequently, we have to
understand the reflection configurations in order to
fully understand global entropy solutions of the two-dimensional
hyperbolic systems of conservation laws, including the admissibility issue for
the entropy solutions.

A natural mathematical approach is to single out
steady shock reflections by
the stability analysis -- the stable ones are physically admissible.
It has been shown in the steady regime that the steady (supersonic or transonic) weak reflection
is always structurally stable in Chen-Chen-Feldman \cite{CCF1} and  Chen-Zhang-Zhu \cite{CZZ} with respect
to the steady perturbation of both the wedge slope
and the incoming steady upstream flow
(even $L^1$--stable for the supersonic weak reflection with respect to the $BV$--perturbation of both the wedge slope
and the incoming steady upstream flow as shown in Chen-Li \cite{CLi}),
while the strong reflection is also structurally stable
under conditional perturbations
({\it cf.} Chen-Chen-Feldman \cite{CCF1,CCF2} and Chen-Fang \cite{ChenFang}).
The first rigorous unsteady analysis of the steady supersonic weak shock solution
as the long-time behavior of an unsteady potential flow
was due to Elling-Liu \cite{EL2}, who dealt with a class of physical
parameters determined by an assumption for angle $\theta_{\rm w}$
less than the sonic angle $\theta_{\rm s}\in (0, \theta_{\rm d})$
(see Chapter 3).

The purpose of this monograph is
to establish the stability theorem for the steady (supersonic or transonic)
weak shock solutions as
the long-time asymptotics of
the global Prandtl-Meyer reflection configurations for unsteady potential flow
for all the admissible physical parameters,
even beyond the sonic angle $\theta_{\rm s}$, up to the detachment
angle $\theta_{\rm d}>\theta_{\rm s}$.
As a corollary, the assumption in Elling-Liu's theorem \cite{EL2}
for the case that $\theta_{\rm w}\in (0, \theta_{\rm s})$
is no longer required.
The global Prandtl-Meyer reflection configurations involve two types of transonic
flow boundaries:
discontinuous and continuous hyperbolic-elliptic phase transition
boundaries for the fluid fields (transonic shocks and sonic arcs).
To establish this theorem, we first reformulate the problem as a free boundary problem
involving transonic shocks
and then carefully establish the
required appropriate monotonicity properties and uniform
{\it a priori} estimates for admissible solutions so that
the approach developed in Chen-Feldman \cite{CF2} can be employed.
This involves several core difficulties in the theory of the underlying nonlinear PDEs:
optimal estimates of solutions of nonlinear degenerate PDEs
and corner singularities (at the corners between the transonic shock as a free boundary and the sonic arcs,
and between the transonic shock and the wedge when the wedge angle $\theta_{\rm w}$ across the sonic angle $\theta_{\rm s}$),
in addition to the involved nonlinear PDEs of mixed elliptic-hyperbolic type and free boundary problems.
Some parts of the results have been announced in Bae-Chen-Feldman \cite{BCF2}.

\begin{figure}[htp]
\centering
\begin{psfrags}
\psfrag{i}[cc][][0.7][0]{$\phantom{aa}(0,-\iv)$}
\psfrag{b}[cc][][0.7][0]{$\beta$}
\psfrag{tw}[cc][][0.7][0]{$\theta_{\rm w}$}
\psfrag{o}[cc][][0.7][0]{$O$}
\psfrag{iv}[cc][][0.7][0]{$(0,-\iv)$}
\psfrag{x}[cc][][0.7][0]{$\xi_1$}
\psfrag{y}[cc][][0.7][0]{$\xi_2$}
\psfrag{s}[cc][][0.7][0]{$\shock:\tx{free boundary}$}
\psfrag{ls}[cc][][0.7][0]{$\leftshock$} \psfrag{rs}[cc][][0.7][0]{$\rightshock$}
\psfrag{Tip}[cc][][0.7][0]{$\xi^*(\beta)$}
\psfrag{snr}[cc][][0.7][0]{$\phantom{aa}\rightsonic$}
\psfrag{snl}[cc][][0.7][0]{$\phantom{aaa}\leftsonic$}
\includegraphics[scale=.5]{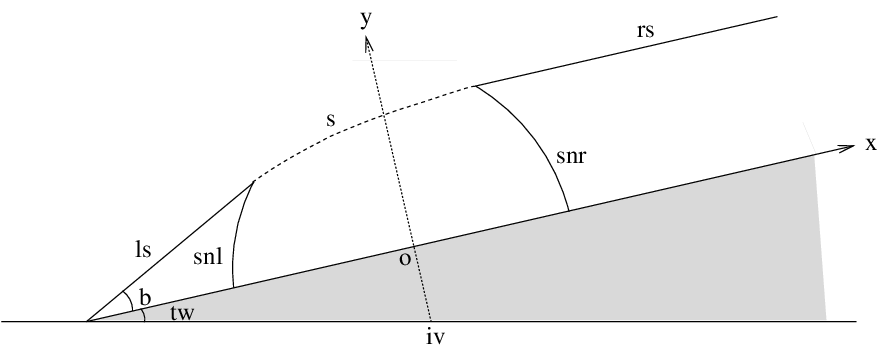}
\includegraphics[scale=.5]{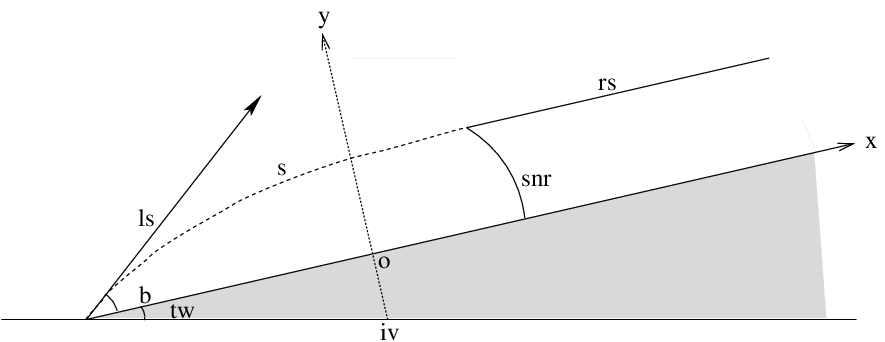}
\caption{Admissible solutions in the $(\iv, \beta)$--parameters in the rotated coordinates $(\xi_1,\xi_2)$
by angle $\theta_{\rm w}$ counterclockwise (Left: $0<\beta<\betasonic$; Right: $\betasonic\le\beta<\betadet$).}
\label{fig1_intro}
\end{psfrags}
\end{figure}

More precisely, in Chapter \ref{section:shock-polar}, we first formulate the physical problem of supersonic flow
impinging upon the solid wedge as an initial-boundary value problem.
By using the invariance under a self-similar scaling
and the physical structure of the problem (see Fig. \ref{fig1_intro}),
the initial-boundary value problem is reformulated
as a boundary value problem in an unbounded domain (Problem \ref{problem-2})
and further as a free boundary problem (Problem \ref{fbp})
for a pseudo-steady potential flow in a bounded domain in the self-similar coordinates
$\xxi=(\xi_1,\xi_2)=\frac{\rx}{t}$ for $t>0$.
Next, we introduce the notion of admissible solutions that we seek in this monograph for all the admissible
physical parameters $(\iu, u_0)\in \mathfrak{P}_{\rm weak}$, where $\iu$ represents the speed of the incoming supersonic flow
and $u_0$ represents the horizontal speed of downstream flow behind a steady weak shock which is uniquely determined
by $\iu$ and angle $\theta_{\rm w}$.
For simplicity, the density of incoming supersonic flow is normalized to be $1$ without loss of generality.
In \S 2.3, the existence of admissible solutions for all $(\iu, u_0)\in \mathfrak{P}_{\rm weak}$ is stated as one of the main theorems.

In order to prove the existence of admissible solutions for all $(\iu, u_0)\in \mathfrak{P}_{\rm weak}$
by employing the Leray-Schauder degree argument,
the first essential step is to introduce a new parameter set $\mathfrak{R}_{\rm weak}$ in \S \ref{section-change-par}.
Given $(\iu, u_0)\in \mathfrak{P}_{\rm weak}$, the half-angle $\theta_{\rm w}$
of the symmetric solid wedge is uniquely determined.
Define $\iv:=\iu\sin \theta_{\rm w}$.
As we will discuss later, $u_0>0$ represents the horizontal speed of the downstream flow behind the weak oblique shock $\leftshock$.
Then we define $\beta\in(0, \frac{\pi}{2})$ as the angle between the wedge boundary and $\leftshock$.
Parameters $(\iv, \beta)$ were first introduced in \cite{EL2}.
In Lemma \ref{lemma-parameters}, we show that there exists
a homeomorphism $\mcl{T}: \mathfrak{P}_{\rm weak}\rightarrow \mcl{T}(\mathfrak{P}_{\rm weak})=:\mathfrak{R}_{\rm weak}$.
More importantly, we show that $\mathfrak{R}_{\rm weak}$ is in the form of
\begin{equation*}
 \mathfrak{R}_{\rm weak}=\underset{\iv>0}{\bigcup}\{\iv\}\times (0,\beta_{\rm d}^{(\iv)}).
\end{equation*}
This structure of $\mathfrak{R}_{\rm weak}$ enables us to prove the existence of admissible solutions
for all $\beta\in(0, \betadet)$ for any fixed $\iv>0$ via the Leray-Schauder degree theorem.
In particular, for each $\iv>0$, there exists an admissible solution for $\beta=0$
and, in \S \ref{subsec-pf-mthm-p1}, we prove that the Leray-Schauder fixed point index of this solution is $1$.
We also show that, for each $\iv>0$, there exists a unique $\betasonic\in(0, \betadet)$,
called the {\emph{sonic angle}}, so that the structure of admissible solutions becomes different as $\beta$ increases
across $\beta=\betasonic$ (see Fig. \ref{fig1_intro}).
Finally, we restate both the definition and existence of admissible solutions
for $(\iv, \beta)\in \mathfrak{R}_{\rm weak}$ in \S \ref{subsec-mainthm}.

\smallskip
In Chapter \ref{section-unif-est-1}, we establish all the {\it a priori} estimates that are essential for solving
the free boundary problem introduced in Chapter \ref{section:shock-polar}.
Furthermore, the {\it a priori} estimates are achieved uniformly on parameters $(\iv, \beta)$.
In particular, this chapter contains the following estimates:
\begin{itemize}
\item[(i)] Strict directional monotonicity properties of $\ivphi-\vphi$;
\item[(ii)] Strict directional monotonicity properties of $\vphi-\rightvphi$ and $\vphi-\leftvphi$;
\item[(iii)] Uniform positive lower bound of the distance between $\shock$ and $\Wedge$ away from the wedge vertex;
\item[(iv)] Uniform positive lower bound of ${\rm dist}(\shock, \der B_1(0, -\iv))$;
\item[(v)] Uniform estimates of the ellipticity of equation $N(\vphi)=0$ in $\Omega$, given in \eqref{fbp-general-description} below;
\item[(vi)] Uniform weighted $C^{2,\alp}$ estimates of admissible solutions in $\Om$.
\end{itemize}
In the above, $\ivphi$, $\leftvphi$, and $\rightvphi$ represent the pseudo-velocity potential functions
for the state of incoming supersonic flow, the state behind the oblique shock $\leftshock$,
and the state behind the normal shock $\rightshock$, respectively.
Moreover, $\der B_1(0,-\iv)$ is the {\emph{sonic circle}} of the incoming supersonic flow:
\begin{equation*}
\der B_1(0, -\iv)
:=\{{\bmxi}\in \R^2\,:\,|D\ivphi(\bmxi)|=1 \}.
\end{equation*}

For fixed $\iv>0$ and $0<\beta<\betasonic$, let $\Om$ be the bounded region
enclosed by $\leftsonic$, $\shock$, $\rightsonic$,
and $\xi_2=0$ in Fig. \ref{fig1_intro}.
In order to find an admissible solution in the sense of Definition \ref{def-regular-sol},
we need to solve the following free boundary problem
for $(\vphi, \shock)$:
\begin{equation}
\label{fbp-general-description}
  \begin{split}
  N(\vphi):={\rm div}\big(\rho(|D\vphi|^2,\vphi)D\vphi\big)+2\rho(|D\vphi|^2,\vphi)=0
  \qquad&\mbox{in $\Om$},\\
  \vphi=\ivphi,\quad
  \rho(|D\vphi|^2,\vphi)D\vphi\cdot {\bm\nu}=D\ivphi\cdot{\bm\nu}\qquad&\mbox{on $\shock$},\\
  \vphi=\leftvphi\qquad&\mbox{on $\leftsonic$},\\
  \vphi=\rightvphi\qquad&\mbox{on $\rightsonic$},\\
  \der_{\xi_2}\vphi=0\qquad&\mbox{on $\der\Om\cap\{\xi_2=0\}$},
  \end{split}
\end{equation}
where $\rho=\rho(|{\bf q}|^2, z)$ is smooth with respect to $({\bf q},z)\in \R^2\times \R$
for $|{\bf q}|\le R_0$ and $|z|\le R_1$ for some positive constants $R_0$ and $R_1$.
Moreover, ${\bm\nu}$ is the inward unit normal vector to $\shock$.
In particular, we seek a solution so that equation $N(\vphi)=0$ is strictly elliptic in $\Om$,
but its ellipticity
degenerates on $\leftsonic\cup\rightsonic$.
As $\beta\in(0,\betasonic)$ tends to $\betasonic$,
$\leftsonic$ shrinks to the wedge vertex $P_{\beta}$,
and the ellipticity of $N(\vphi)=0$ degenerates at $P_{\beta}$  for $\beta=\betasonic$.
For $\beta>\betasonic$, $N(\vphi)=0$ is strictly elliptic at $P_{\beta}$.
For $\beta\ge \betasonic$, the boundary condition $\vphi=\leftvphi$ on $\leftsonic$
given in \eqref{fbp-general-description} becomes a one-point Dirichlet boundary condition.
Therefore, it is crucial to achieve estimate (v) and then employ the result to establish
the uniform {\it a priori} estimates of admissible solutions in $\Om$ by estimate (vi).

Once estimates (i)--(ii) are established,
we adjust the argument in \cite{CF2} to achieve estimates (iii)--(vi),
although there are several technical differences, due to the structural differences
of the solutions constructed in this monograph compared to those in \cite{CF2}.
We also point out that estimate (iv) is the key to achieving estimates (v)--(vi).
Using the argument in \cite{CF2}, for any fixed $\iv>0$, we are able to establish a uniform estimate of
positive lower bound of ${\rm dist}(\shock, \der B_1(0, -\iv))$ for all the admissible solutions corresponding
to $\beta\in(0,\beta_*]$ whenever $\beta_*\in (0,\betadet)$.
Owing to this property, we prove the existence of admissible solutions for all the admissible physical
parameters $(\iv, \beta)\in \mathfrak{R}_{\rm weak}$, even beyond the sonic angle $\beta^{(v_\infty)}_{\rm s}$.

\begin{figure}[htp]
\centering
\begin{psfrags}
\psfrag{i}[cc][][0.7][0]{$-\iv\;\;\;$}
\psfrag{b}[cc][][0.7][0]{$\beta$}
\psfrag{ns}[cc][][0.7][0]{$\rightshock$}
\psfrag{o}[cc][][0.7][0]{$O$}
\psfrag{ls}[cc][][0.7][0]{$\leftshock$}
\psfrag{s}[cc][][0.7][0]{$\Shock$}
\psfrag{lt}[cc][][0.7][0]{$\lefttop$}
\psfrag{rt}[cc][][0.7][0]{$\righttop$}
\psfrag{lso}[cc][][0.7][0]{$\;\;\leftsonic$}
\psfrag{rso}[cc][][0.7][0]{$\rightsonic$}
\psfrag{x}[cc][][0.7][0]{$\xi$}
\psfrag{oom}[cc][][0.7][0]{$\vphi=\leftvphi$}
\psfrag{nom}[cc][][0.7][0]{$\vphi=\rightvphi$}
\psfrag{y}[cc][][0.7][0]{$\eta$}
\psfrag{lb}[cc][][0.7][0]{$\leftbottom$}
\psfrag{rb}[cc][][0.7][0]{$\rightbottom$}
\psfrag{om}[cc][][0.7][0]{$\Om: M<1$}
\psfrag{pb}[cc][][0.7][0]{$P_{\beta}$}
\psfrag{vn}[cc][][0.7][0]{${\bf e}_{\mcl{N}}$}
\psfrag{vo}[cc][][0.7][0]{${\bf e}_{\mcl{O}}$}
\psfrag{n}[cc][][0.7][0]{$\Om_{\mcl{N}}$}
\psfrag{o}[cc][][0.7][0]{$\Om_{\mcl{O}}$}
\psfrag{e}[cc][][0.7][0]{$\xi_2$}
\psfrag{M}[cc][][0.7][0]{$\Om$}
\psfrag{c}[cc][][0.7][0]{${\rm cone}(\bm e_{\mcl{N}},\bm e_{\mcl{O}})$}
\includegraphics[scale=.5]{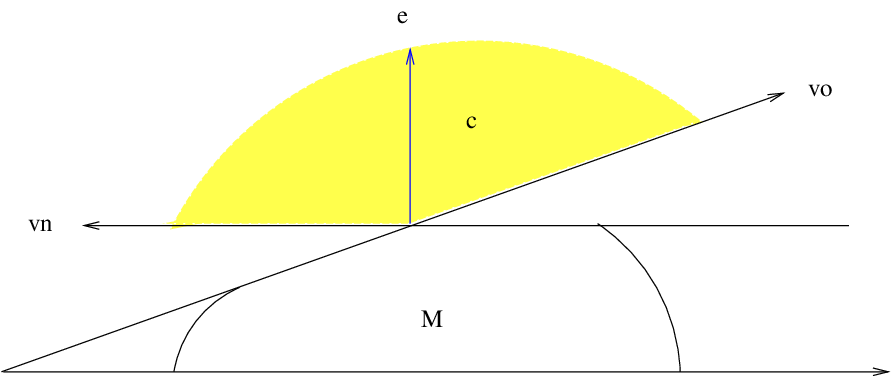}
\caption{The cone of monotonicity}
\label{fig2_intro}
\end{psfrags}
\end{figure}

Even though the overall argument follows \cite{CF2},
there are several significant differences from \cite{CF2}.
One of them is the choice of directions for the monotonicity properties of $\ivphi-\vphi$, $\vphi-\leftvphi$, and $\vphi-\rightvphi$.
For fixed $(\iv, \beta)\in \mathfrak{R}_{\rm weak}$, define ${\bf e}_{\mcl{N}}:=(0,-1)$ and ${\bf e}_{\mcl{O}}:=(\cos\beta, \sin \beta)$.
Then ${\bf e}_{\mcl{N}}$ is the unit tangent vector to the normal shock $\rightshock$,
and ${\bf e}_{\mcl{O}}$ is the unit tangent vector to the oblique shock $\leftshock$.
Moreover, we define the cone of monotonicity as shown in Fig. \ref{fig2_intro}
by
\begin{equation*}
\cone:=\{\alp_1\leftvec+\alp_2\rightvec\,:\,\alp_1, \alp_2> 0\}.
\end{equation*}
In \S \ref{subsection-directional-mono}, we show that any admissible solution $\vphi$ satisfies
\begin{equation}
\label{monotonicity-intro}
\der_{\ee}(\ivphi-\vphi)<0\qquad \tx{in $\ol{\Om}\, $ for all $\ee\in \cone$},
\end{equation}
from which many essential estimates of admissible solutions can be further obtained.
For example, \eqref{monotonicity-intro}, combined with the Rankine-Hugoniot conditions on $\shock$,
implies that $\shock$ is represented as a graph of a function $\xi_2=\fshock(\xi_1)$ with $\fshock'(\xi_1)>0$.
This property is a key ingredient in the proof of the separation of $\shock$ from {\emph{the sonic circle}} $\der B_1(0, -\iv)$
of the incoming supersonic flow. Notice that this separation property is crucial for establishing
the uniform estimate of the ellipticity of equation $N(\vphi)=0$ in $\Om$.
In addition, further monotonicity properties of $\vphi-\leftvphi$ and $\vphi-\rightvphi$ in $\cone$ are achieved,
which play important roles in the {\it a priori} estimates of admissible solutions near $\leftsonic\cup\rightsonic$.

In Chapter \ref{section-itr-set}, we define the iteration set $\mcl{K}$ consisting of approximate admissible solutions.
Note that the pseudo-subsonic region $\Om$ of each admissible solution is different.
Furthermore, as $\beta$ increases across $\betasonic$, the shape of $\Om$ changes from a rectangular domain
to a triangular domain. This is because the sonic arc $\leftsonic$ corresponding to the oblique shock $\leftshock$
shrinks to the wedge vertex  $P_{\beta}$ as $\beta\in(0,\betasonic)$ tends to $\betasonic$,
and $\leftsonic=\{P_{\beta}\}$ for $\beta\ge \betasonic$.
For this reason, it is necessary to introduce a diffeomorphism $\mathfrak{F}$ so that $\mathfrak{F}^{-1}(\Om)$ is
the fixed domain $\iter:=(-1,1)\times(0,1)$.
Moreover, $\mathfrak{F}$ should be defined so that $\mathfrak{F}$ depends continuously on $\beta\in[0, \betadet)$
and admissible solutions in an appropriately chosen norm.
In \S\ref{subsec-Q}, we define a map $\mathfrak{F}$ for each admissible solution such that
\begin{equation*}
\begin{split}
  &\mathfrak{F}(\iter)=\Om,\qquad \,\, \mathfrak{F}(\shock)=\{(s,1)\,:\,-1<s<1\},\\
  &\mathfrak{F}(\leftsonic)=\{(-1,t)\,:\,0<t<1\},\quad\,\,  \mathfrak{F}(\rightsonic)=\{(1,t)\,:\,0<t<1\}.
  \end{split}
\end{equation*}
Since the sonic arc $\rightsonic$ corresponding to the normal shock $\rightshock$ is fixed so as to be the same
for all $\beta\in[0, \betadet)$ (see Fig. \ref{fig1_intro}),
the definition of $\mathfrak{F}$ in this monograph can be given more explicitly than the one given in \cite{CF2}; see Definition \ref{definition-Gset-shocks-new}.
In \S \ref{subsec-mapping-invert}, the definition of $\mathfrak{F}$ is extended to a class of approximate admissible solutions.
Then we set up the iteration set $\mcl{K}$ and analyze its properties
in \S \ref{subsection-def-iterset}--\S \ref{subsec-prop-iterset}.
The iteration set $\mcl{K}$ is given in the form
$$
\mathcal{K}:=\bigcup_{\beta\in[0,\beta_*]}\{\beta\}\times \mcl{K}(\beta)
\qquad\,\,\mbox{for fixed $\iv>0$ and $\beta_*\in(0,\betadet)$},
$$
where each $\mcl{K}(\beta)$ is a subset of $C^{1,\alp}(\ol{\iter})$ for some $\alp\in(0,1)$.

In Chapter \ref{sec-proof-theorem5-case1}, for fixed $\iv>0$, we define an iteration
map
$$
\mcl{I}(\cdot, \beta)\,:\,\mcl{K}(\beta)\rightarrow C^{2,\alp}_{(*, \alpha_1)}(\iter)
\qquad\,\, \mbox{for $\iter:=(-1,1)\times (0,1)\subset \R^2$},
$$
where $C^{2,\alp}_{(*,\alpha_1)}(\iter)$ is a weighted $C^{2,\alp}$ space.
The iteration map $\mcl{I}$ is defined so that, if $\mcl{I}(u_*, \beta)=u_*$
for $u_*\in \mcl{K(\beta)}$, then $(\vphi, \shock)$, given by
\begin{equation*}
  \vphi=u_*\circ \mathfrak{F}^{-1}_{(u_*, \beta)}+\vphi_{\beta}^*\quad \tx{in $\Om=\mathfrak{F}_{(u_*, \beta)}(\iter)$},
  \qquad\,\shock=\mathfrak{F}_{(u_*, \beta)}(\{(s,1)\,:\,-1<s<1\}),
\end{equation*}
solves the free boundary problem \eqref{fbp-general-description}.
In the above, $\vphi_{\beta}^*$ is a smooth interpolation of $\leftvphi$ and $\rightvphi$.
The precise definition of $\vphi_{\beta}^*$ is given by \eqref{12-32}.
Finally, the existence of a fixed point of $\mcl{I}(\cdot, \beta)$ in $\mcl{K}(\beta)$ for all $\beta\in(0,\beta_*]$
is proved by employing the Leray-Schauder degree argument in \S \ref{subsec-pf-mthm-p1}.
In this way, we establish the existence of admissible solutions
for all $(\iv, \beta)\in \mathfrak{R}_{\rm weak}$ (Theorem \ref{theorem-0}),
hence the existence of admissible solutions for all $(\iu, u_0)\in \mathfrak{P}_{\rm weak}$ (Theorem \ref{theorem-0}).

Theorem \ref{theorem-2-op}, or equivalently, Theorem \ref{theorem-3}, which pertains to the optimal regularity of
admissible solutions,
is established in Chapter \ref{sec-optimal}.

To make the monograph self-contained, we also include
Appendices A--C, which contain some results required for establishing the main theorems
and a proof of the non-existence of self-similar strong shock solutions.

\smallskip
A closely related problem to the one we have solved here is
the shock reflection-diffraction problem
which was addressed in Chen-Feldman \cite{CF2}.
Even though the two problems are two different lateral Riemann problems and have different issues and features,
the approach developed in Chen-Feldman \cite{CF2} for the shock reflection-diffraction problem
has been adopted for solving our Prandtl-Meyer reflection problem in this monograph.
As discussed earlier, one of the main contributions of this monograph is to identify
appropriate monotonicity properties
and establish suitable uniform {\it a priori} estimates for admissible solutions, based on the new and careful choice
of the directions for the monotonicity properties; as a result,
the Chen-Feldman approach in \cite{CF2} can be employed.

In this monograph, we have solved the Prandtl-Meyer reflection problem up to the detachment angle
in the framework of the potential flow equation,
which has been widely used for discontinuous flows in applications in aerodynamics, especially
when the amount of vorticity is relatively small in the region of interest.
When the flow regions of interest have large amounts of vorticity,
the full compressible Euler equations are usually required.
Nevertheless, for the solutions containing a shock of small strength,
the potential flow equation and the full Euler equations match
each other well, right up to the third-order of the shock strength.
Furthermore, for the problem analyzed in this monograph,
the Euler equations for potential flow is actually {\it exact} in two important regions
of the solutions near the two sonic arcs in the subsonic domain $\Omega$.
Even in the other part of  domain $\Omega$, under the Helmholtz-Hodge decomposition for the velocity field,
the full Euler equations in the self-similar coordinates can be decomposed as the potential flow type equation,
coupled with the incompressible Euler type equations plus a transport equation for the entropy function.
These can be shown by directly following the arguments in \S 18.7 in Chen-Feldman \cite{CF2}.
In this sense,
the analysis and related methods/techniques developed in this monograph could also play an essential role
in finding a solution of the problem in the framework of the full Euler equations.
In particular, our results for the potential flow equation have provided useful insights on what will happen
for the case of the full Euler equations.

Finally, we remark in passing that, for the uniqueness/stability problems, it is necessary to consider solutions
in a restricted class.
Recent results \cite{ChiodaroliDeLellisKreml,ChiodaroliKreml,FeireislKlingenbergKremlMarkfelder,MarkfelderKlingenberg}
show the non-uniqueness of solutions with flat shocks in the class of entropy solutions
of the Cauchy problem (initial value problem) for the multidimensional compressible Euler equations
(isentropic and full).
The Prandtl-Meyer reflection problem under consideration in this monograph is different -- the problem for
solutions with non-flat shocks for potential flow on the domain with
boundaries, so these non-uniqueness results do not apply directly.
However, these results indicate that it is natural to study the uniqueness and stability problems
in a more restricted class of solutions.
Since the completion of this monograph, some progress on the uniqueness in the class of self-similar solutions
of regular shock reflection-diffraction configurations
with convex transonic shocks (which are called admissible solutions) has been made, as announced recently in \cite{CFX}.
A similar uniqueness result can also be obtained by combining the approach in \cite{CFX,CFX2} with the estimate techniques
developed in this monograph.
Technically, restricting the uniqueness to the class
of admissible solutions allows us to reduce
the problem
to a corresponding uniqueness problem for solutions of a
free boundary problem for a nonlinear elliptic equation, which is degenerate
for the supersonic case.




\chapter{Mathematical Problems and Main Theorems}
\label{section:shock-polar}

\numberwithin{equation}{section}

In this chapter, we first formulate the physical problem of a supersonic flow impinging upon the solid wedge
into an initial-boundary value problem.
Then, based on the invariance of both the problem and
the governing equations under the self-similar scaling,
we reformulate the initial-boundary value problem as a boundary value problem in an unbounded domain (Problem \ref{problem-2}),
and further as a free boundary problem  in a bounded domain (Problem \ref{fbp})
for the existence of Prandtl-Meyer reflection configurations involving two types of transonic flow boundaries:
discontinuous and continuous hyperbolic-elliptic phase transition boundaries
for the fluid fields (transonic shocks and sonic arcs).
The main theorems of this monograph are presented in \S \ref{subsection-mainthm-ori-par} and \S \ref{subsec-mainthm}.

\section{Mathematical Problems}
\label{subsec-math-prob}
The compressible potential flow is governed by the conservation law of mass and the Bernoulli law:
\begin{align}
\label{1-a}
&\der_t\rho+ \nabla_{\bf x}\cdot (\rho \nabla_{\bf x}\Phi)=0,\\
\label{1-b}
&\der_t\Phi+\frac 12|\nabla_{\bf x}\Phi|^2+h(\rho)=B,
\end{align}
where $\rho$ is the density, $\Phi$ is the velocity potential,
$B$ is the Bernoulli constant determined by the incoming flow and/or boundary conditions,
and $h(\rho)$ is given by
\begin{equation*}
h(\rho)=\int_1^{\rho}\frac{p'(\varrho)}{\varrho}\,
\dd\varrho=\int_1^{\rho}\frac{c^2(\varrho)}{\varrho}\,\dd\varrho
\end{equation*}
for the sound speed $c(\rho)$ and pressure $p$.
For an ideal polytropic gas, the sound speed $c$ and pressure $p$ are given by
\begin{equation}
\label{cont-rel-pt}
 c^2(\rho)=\kappa \gam \rho^{\gam-1}, \qquad p(\rho)=\kappa \rho^{\gam}
\end{equation}
for constants $\gam>1$ and $\kappa>0$.
If $(\rho, \Phi)(t, {\bf x})$ solves \eqref{1-a}--\eqref{1-b} with \eqref{cont-rel-pt},
then $(\til{\rho}, \til{\Phi})(t,{\bf x})=(\rho, \Phi)(\alp ^2t, \alp{\bf x})$
with $\alp:=\frac{1}{\sqrt{\kappa\gam}}$ solves
\begin{align*}
&\der_t\til{\rho}+ \nabla_{\bf x}\cdot (\til{\rho} \nabla_{\bf x}\til{\Phi})=0,\\
&\der_t\til{\Phi}+\frac 12|\nabla_{\bf x}\til{\Phi}|^2+\frac{\til{\rho}^{\gam-1}-1}{\gam-1}=\alp^2 B.
\end{align*}
Therefore, we choose $\kappa=\frac{1}{\gam}$ without loss of generality so that
\begin{equation}
\label{1-c}
h(\rho)=\int_1^{\rho} h'(\varrho) \dd\varrho=\frac{\rho^{\gam-1}-1}{\gam-1},\qquad c^2(\rho)=\rho^{\gam-1}.
\end{equation}
The case of the isothermal flow can be included as the isothermal limit $\gam\to 1+$ in \eqref{1-c}.
Therefore, we define $(h, c^2)(\rho)$ by
\begin{equation}
\label{def-h-and-c}
(h, c^2)(\rho)=\begin{cases}
(\frac{\rho^{\gam-1}-1}{\gam-1}, \rho^{\gam-1})&\qquad\mbox{for $\gam>1$},\\
(\ln \rho,1) &\qquad \mbox{for $\gam=1$}.
\end{cases}
\end{equation}
By \eqref{1-b},
$\rho$ can be expressed as
\begin{equation}
\label{1-b1}
\rho(\der_t\Phi,\nabla_{\bf x}\Phi)=h^{-1}(B-\der_t\Phi-\frac 12|\nabla_{\bf x}\Phi|^2).
\end{equation}
Then system \eqref{1-a}--\eqref{1-b} can be rewritten as
\begin{equation}
\label{1-b2}
\der_t\rho(\der_t\Phi, \nabla_{\bf x}\Phi)
+\nabla_{\bf x}\cdot\big(\rho(\der_t\Phi, \nabla_{\bf x}\Phi)\nabla_{\bf x}\Phi\big)=0,
\end{equation}
with $\rho(\der_t\Phi, \nabla_{\bf x}\Phi)$ determined by \eqref{1-b1}.

\bigskip
A steady state solution $\bPhi({\bf x})$ to \eqref{1-a}--\eqref{1-b} yields the steady potential flow equations
\begin{equation}
\label{eqn-steady-pot-new2015}
\begin{split}
&\nabla_{\bf x}\cdot (\bar{\rho} \nabla_{\bf x}\bPhi)=0,\\
&\frac 12|\nabla_{\bf x}\bPhi|^2+h(\bar{\rho})=B.
\end{split}
\end{equation}

\begin{figure}[htp]
\centering
\begin{psfrags}
\psfrag{ic}[cc][][0.8][0]{$\irho>0$, $\iu>\irho^{(\gam-1)/2}$}
\includegraphics[scale=1.0]{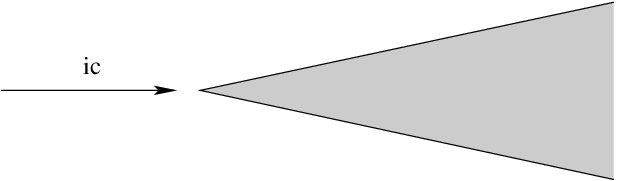}
\caption{Supersonic flow impinging upon a solid wedge}\label{Figure-1a}
\end{psfrags}
\end{figure}

A symmetric wedge $W$ of half-angle $\theta_{\rm w}\in (0, \frac{\pi}{2})$
in $\R^2$ (Fig. \ref{Figure-1a}) is defined by
\begin{equation}
\label{definition-W}
W:=\{{\bf x}=(x_1,x_2)\in \R^2\,:\,|x_2|<x_1\tan\theta_{\rm w}, x_1>0\}.
\end{equation}
On the wedge boundary  $\partial W$, $\bPhi$ must satisfy the slip boundary condition
$\der_{{\bf n}_{\rm w}} \bPhi=0\,\,\tx{on}\,\, \der W$,
where ${\bf n}_{\rm w}$ indicates the outward unit normal vector to $\der W$.
Denote $D:=\R^2\setminus W$, and consider the boundary value problem
for \eqref{eqn-steady-pot-new2015} in $D$ with
\begin{equation}
\label{bvp-steady-new2015}
\der_{{\bf n}_{\rm w}} \bPhi=0\,\qquad \tx{on}\,\, \der D=\der W.
\end{equation}

If a supersonic flow
with a constant density $\irho>0$ and a velocity ${\bf u}_{\infty}=(\iu, 0)$,
$\iu>\irho^{(\gam-1)/2}$, moves towards wedge $W$, and if $\theta_{\rm w}$ is less
than a critical angle called {the \emph{detachment angle}},
then the well-known {\emph{shock polar analysis}} shows that there are two different steady weak solutions
to the boundary value problem
\eqref{eqn-steady-pot-new2015}--\eqref{bvp-steady-new2015}:
{\emph{the steady weak shock solution}} and {\emph{the steady strong shock solution}}.
For more precise arguments, we first define a class of weak solutions of
the boundary value
problem \eqref{eqn-steady-pot-new2015}--\eqref{bvp-steady-new2015}.

\begin{definition}
\label{def-weaksol-steady-new2015}
Let $\Gamma_{\rm sh}$ be a $C^1$--curve that lies in $D$
and divides $D$ into two open subsets $D^-$ and $D^+$.
We say that $\bPhi\in W^{1,\infty}(D)$ is a steady entropy solution
with a shock $\Gamma_{\rm sh}$
of the boundary value
problem \eqref{eqn-steady-pot-new2015}--\eqref{bvp-steady-new2015}
if $\bPhi$ satisfies the following properties{\rm :}

\smallskip
\begin{itemize}
\item[(i)] $B-\frac 12|\nabla_{\bf x}\bPhi|^2 > h(0+)$ {\it a.e.} in $D${\rm ;}

\smallskip
\item[(ii)] For each $\zeta\in C^{\infty}_0(\R^2)$,
\begin{equation*}
\int_{D} \bar{\rho}(|\nabla_{\bf x}\bPhi|^2)\nabla_{\bf x}\bPhi\cdot \nabla_{\bf x}\zeta\, \dd{\bf x}=0{\rm ;}
\end{equation*}

\item[(iii)]
$\bPhi\in C^1(\ol{D^{\pm}})\cap C^2(D^{\pm})${\rm ;}

\smallskip
\item[(iv)] Entropy condition{\rm :} for $\bPhi^\pm:=\bPhi|_{D^\pm \cup\Gam_{\rm sh}}$,
\begin{equation*}
\der_{{\bf n}_{\rm sh}} {\bPhi}^{-}> \der_{{\bf n}_{\rm sh}} {\bPhi}^+ >0 \qquad\,\tx{on $\Gam_{\rm sh}$},
\end{equation*}
or equivalently, $\bar{\rho}(\nabla_{\bf x}\bPhi^-)<\bar{\rho}(\nabla_{\bf x}\bPhi^+)$ along the flow direction,
where ${\bf n}_{\rm sh}$ represents the unit normal vector to $\Gam_{\rm sh}$ pointing from $D^-$ towards $D^+$.
\end{itemize}
\end{definition}

\begin{remark}\label{remark-conormal-BC}
By performing integration by parts,
condition {\rm (ii)} of Definition {\rm \ref{def-weaksol-steady-new2015}} implies that
any entropy solution with a shock $\shock$
of problem \eqref{eqn-steady-pot-new2015}--\eqref{bvp-steady-new2015} in the sense of
Definition {\rm \ref{def-weaksol-steady-new2015}} satisfies the conormal boundary condition:
\begin{equation*}
\bar{\rho}(|\nabla_{\bf x}\bPhi|^2)\nabla_{\bf x}\bPhi\cdot{\bf n}_{\rm w}=0 \qquad\tx{on $\der W$}.
\end{equation*}
Furthermore, combining conditions {\rm (i)} and {\rm (iii)} of Definition {\rm \ref{def-weaksol-steady-new2015}}
with the conormal boundary condition stated immediately above yields that the entropy solution $\bar{\Phi}$ indeed
satisfies the boundary condition \eqref{bvp-steady-new2015} if $\bar{\rho}(|\nabla_{\bf x}\bar{\Phi}|^2)>0$ holds on $\der W$.
\end{remark}

In particular, Definition \ref{def-weaksol-steady-new2015}, via integration by parts, leads to
the following Rankine-Hugoniot jump conditions for the steady potential flow equations \eqref{eqn-steady-pot-new2015}:
\begin{equation}
\label{steady-RH-new}
[\bPhi]_{\Gam_{\rm sh}}=[\bar{\rho}(|\nabla_{\bf x} \bPhi|^2)\nabla\bPhi\cdot{\bf n}_{\rm sh}]_{\Gam_{\rm sh}}=0,
\end{equation}
where $[F({\bf x})]_{\Gam_{\rm sh}}:=F^+({\bf x})-F^-({\bf x})$ for ${\bf x}\in \Gam_{\rm sh}$.

\begin{definition}[The steady Prandtl-Meyer reflection solution]
\label{def-PM-steady-new2015}
{\emph{The steady Prandtl-Meyer reflection solution}} for potential flow is an entropy solution $\bPhi$ with a shock $\Gam_{\rm sh}$
of the boundary value problem \eqref{eqn-steady-pot-new2015}--\eqref{bvp-steady-new2015} in the sense of Definition
{\rm \ref{def-weaksol-steady-new2015}} with the following additional features{\rm :}

\smallskip
\begin{itemize}
\item[(i)] $\Gam_{\rm sh}=\{{\bf x}=(x_1,x_2)\in \R^2\setminus W\,:\, |x_2|=x_1\tan\theta_{\rm sh}, x_1\ge 0\}$
for some $\theta_{\rm sh}\in(\theta_{\rm w}, \frac{\pi}{2})${\rm ;}

\smallskip
\item[(ii)] For some constants $u_0,v_0>0$,
$$
\bPhi({\bf x})=\begin{cases}
\iu x_1\,\, &\tx{in}\,\,D^-=\{{\bf x}\in D\,:\, x_1<|x_2|\cot \theta_{\rm sh}\},\\
u_0x_1+v_0 x_2\,\quad &\tx{in}\,\,D^+:=D\setminus \ol{D^-};
\end{cases}
$$

\item[(iii)] $\tan \theta_{\rm sh}=\frac{\iu-u_0}{v_0}${\rm ;}

\smallskip
\item[(iv)] Entropy condition{\rm :} for the unit normal vector ${\bf n}_{\rm sh}$ to $\Gam_{\rm sh}$ pointing from $D^-$ towards $D^+$,
$$
\nabla \bPhi^-\cdot {\bf n}_{\rm sh}>\nabla_{\bf x} \bPhi^+\cdot {\bf n}_{\rm sh}>0 \qquad\,\,\, \mbox{on $\Gam_{\rm sh}$},
$$
or equivalently, $\bar{\rho}(|\nabla_{\bf x}\bPhi^-|^2)<\bar{\rho}(|\nabla_{\bf x}\bPhi^+|^2)$.
\end{itemize}
\end{definition}

\begin{lemma}\label{lemma-spolar-steady-new}
Given any $\gam\ge 1$ and $(\irho, \iu)$ with $\iu>c_\infty=\irho^{(\gam-1)/2}>0$,
there exist unique $\underline{u}^{(\irho, \iu)}\in (0, \iu)$ and
$\theta_{\rm d}^{(\irho, \iu)}\in(0,\frac{\pi}{2})$  such that the following properties hold{\rm :}

\smallskip
\begin{itemize}
\item[(a)] For each $\theta_{\rm w}\in(0, \theta_{\rm d}^{(\irho, \iu)})$,
there are exactly two constants $u_{\rm st}$ and $u_{\rm wk}$ with $\underline{u}^{(\irho, \iu)}<u_{\rm st}<u_{\rm wk}<\iu$
yielding two steady Prandtl-Meyer reflection configurations
in the sense that, if $(u_0, v_0)=u_{\rm st}(1,\tan\theta_{\rm w})$ or $u_{\rm wk}(1,\tan\theta_{\rm w})$
in Definition {\rm \ref{def-PM-steady-new2015}},
then the corresponding function $\bPhi$ is an entropy solution of the boundary value problem
\eqref{eqn-steady-pot-new2015}--\eqref{bvp-steady-new2015}
with shock $\Gam_{\rm sh}$ given by Definition {\rm \ref{def-PM-steady-new2015}}{\rm (i)}
with $\theta_{\rm sh}$ being determined by Definition {\rm \ref{def-PM-steady-new2015}}{\rm (iii)}{\rm ;}

\smallskip
\item[(b)] $u_{\rm st}$ and $u_{\rm wk}$ depend continuously on
$(\irho, \iu, \gam)$ and $\theta_{\rm w}\in(0, \theta_{\rm d}^{(\irho, \iu)})$,
and $u_{\rm st}=u_{\rm wk}$ at $\theta_{\rm w}=\theta_{\rm d}^{(\irho, \iu)}${\rm ;}

\smallskip
\item[(c)] For each $\theta_{\rm w}\in(0, \theta_{\rm d}^{(\irho, \iu)})$,
let $u_{\rm wk}^{(\theta_{\rm w})}$ denote the value of $u_{\rm wk}$ corresponding to $\theta_{\rm w}$.
Then there exists a unique $\theta_{\rm s}^{(\irho, \iu)}\in(0, \theta_{\rm d}^{(\irho, \iu)})$ such that
\begin{equation*}
|u_{\rm wk}^{(\theta_{\rm s}^{(\irho, \iu)})}| |(1,\tan \theta_{\rm s}^{(\irho,\iu)})|
=\big(\bar{\rho}( |u_{\rm wk}^{(\theta_{\rm s}^{(\irho, \iu)})}|^2 |(1,\tan \theta_{\rm s}^{(\irho,\iu)})|^2)\big)^{(\gam-1)/2}.
\end{equation*}
In other words, the flow behind the weak shock corresponding to $\theta_{\rm s}^{(\irho, \iu)}$ is sonic.
\end{itemize}
\end{lemma}

\begin{figure}[htp]
\centering
\begin{psfrags}
\psfrag{tc}[cc][][0.8][0]{$\phantom{aaaaaaaaaaa}\frac{v}{u}=\tan\theta_{\rm s}^{(\irho, \iu)}$}
\psfrag{tw}[cc][][0.8][0]{$\phantom{aaaaaaaaaaa}\frac{v}{u}=\tan\theta_{\rm w}$}
\psfrag{td}[cc][][0.8][0]{$\phantom{aaaaaaaaaaa}\frac{v}{u}=\tan\theta_{\rm d}^{(\irho, \iu)}$}
\psfrag{u}[cc][][0.8][0]{$u$}
\psfrag{v}[cc][][0.8][0]{$v$}
\psfrag{u0}[cc][][0.8][0]{$\iu$}
\psfrag{zeta}[cc][][0.8][0]{}
\psfrag{ud}[cc][][0.8][0]{$u_{\rm d}$}
\includegraphics[scale=.8]{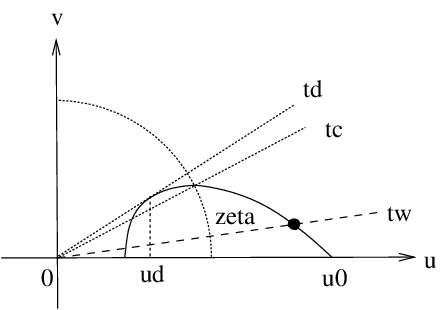}
\caption{Shock polars in the $(u,v)$--plane}\label{fig:polar}
\end{psfrags}
\end{figure}

\begin{proof}
(a) and (b) can be checked directly from Lemmas \ref{lemma-appendix} and \ref{lemma-app2}.

Define $q(\theta_{\rm w}):=|u_{\rm wk}^{(\theta_{\rm w})}| |(1,\tan\theta_{\rm w})|$.
We first observe that $|q(\theta_{\rm w})|^2=\big(\bar{\rho}(|q(\theta_{\rm w})|^2)\big)^{\gamma-1}$
if and only if $|q(\theta_{\rm w})|^2=\frac{2}{\gam+1}\big(1+(\gam-1)B\big)=:\mcl{K}_0$.
To prove (c),  it suffices to show that there exists
a unique $\theta_*\in(0,\theta_{\rm d}^{(\irho, \iu)})$ satisfying $|q(\theta_*)|^2=\mcl{K}_0$.

Condition $\iu^2>\irho^{\gam-1}$ implies that $|q(0)|^2> \mcl{K}_0$.
This can also be checked from the Bernoulli law ({\it i.e.}, $\frac 12|\nabla_{\bf x}\bPhi|^2+h(\bar{\rho})=B$)
and the conservation law of mass ({\it i.e.}, $\bar{\rho}(u_{\rm st}^{(0)}) u_{\rm st}^{(0)}=\rho_{\infty}u_{\infty}$
so that $|u_{\rm st}^{(0)}|^2<\mcl{K}_0$).
Then
there exists a unique point  $P_*=u_*(1,\tan\theta_*)$
on the shock polar $\Upsilon^{(\irho, \iu)}$
satisfying $|P_*|^2=\mcl{K}_0$ (see Lemma \ref{lemma-app2}).
It remains to verify that $u_*=u_{\rm wk}^{(\theta_*)}$;
that is, $P_*$ is the weak shock point corresponding to
$\theta_*\in(0,\theta_{\rm d}^{(\irho, \iu)})$.

In Lemmas \ref{lemma-appendix} and \ref{lemma-app2}, it is shown that the shock polar
curve $\Upsilon^{(\rho_{\infty}, u_{\infty})}$, as shown in Fig. \ref{fig:polar},
is given as the zero-level curve of $g({\bf u})$
in the first quadrant of the $(u,v)$--plane
and that $\Upsilon^{(\irho, \iu)}$ is convex.
Furthermore, $g_{\bf u}({\bf u})$ is a normal vector to $\Upsilon^{(\irho, \iu)}$
at ${\bf u}\in \Upsilon^{(\irho, \iu)}$
towards the $u$--axis.
From this observation, we see that
\begin{equation*}
\begin{split}
&g_{{\bf u}}(P_*)\cdot P_*>0\qquad\tx{if and only if}\qquad u_*=u_{\rm st}^{(\theta_*)},\\
&g_{{\bf u}}(P_*)\cdot P_*=0\qquad\tx{if and only if}\qquad \theta_*=\theta_{\rm d}^{(\irho, \iu)},\\
&g_{{\bf u}}(P_*)\cdot P_*<0\qquad\tx{if and only if}\qquad u_*=u_{\rm wk}^{(\theta_*)}.
\end{split}
\end{equation*}
Now we compute $g_{{\bf u}}(P_*)\cdot P_*$. A direct computation by using \eqref{g} gives that
\begin{equation*}
g_{\bf u}({\bf u})=\frac{1}{\bar{\rho}^{\gam-2}}\left(\bar c^2 \frac{{\bf u}_{\infty}
-{\bf u}}{|{\bf u}_{\infty}-{\bf u}|}-\big({\bf u}\cdot \frac{{\bf u}_{\infty}-{\bf u}}{|{\bf u}_{\infty}-{\bf u}|}\big){\bf u}\right)
-\frac{\bar{\rho}{\bf u}-\irho {\bf u}_{\infty}}{|{\bf u}_{\infty}-{\bf u}|},
\end{equation*}
where $\bar{\rho}=\bar{\rho}(|{\bf u}|^2)$, $\bar c^2=\bar{\rho}^{\gam-1}$,
and ${\bf u}_{\infty}=(\iu, 0)$.
Combining \eqref{steady-RH-new} with $|P_*|^2=\mcl{K}_0$ yields
\begin{equation*}
g_{\bf u}(P_*)\cdot P_*= -\big(\bar{\rho}(|P_*|^2)-\irho\big)(P_*\cdot{\bm \tau}_{\rm s})^2,
\end{equation*}
where ${\bm\tau}_{\rm s}$ represents a unit tangent vector to
shock $\leftshockop$ corresponding to state $P_*$.
Since $P_*\cdot {\bm\tau}_{\rm s}\neq 0$, we obtain from the entropy condition $\bar{\rho}(P_*)-\irho>0$
that $g_{\bf u}(P_*)\cdot P_*<0$.
From this, we conclude that $u_*=u_{\rm wk}^{(\theta_*)}$.
Choosing $\theta_{\rm s}^{(\irho, \iu)}=\theta_*$, we complete the proof.
\end{proof}

\begin{definition}
\label{def-weak-strong-steady-new2015}
Fix parameters $(\irho, \iu, \gam, \theta_{\rm w})$.
In Lemma {\rm \ref{lemma-spolar-steady-new}}, $\bPhi$ with $(u_0,v_0)=u_{\rm st}(1,\tan\theta_{\rm w})$
is called a steady Prandtl-Meyer strong reflection solution, and $\bPhi$ with $(u_0,v_0)=u_{\rm wk}(1,\tan\theta_{\rm w})$
is called a steady Prandtl-Meyer weak reflection solution in the sense that
\begin{equation*}
|(\iu,0)-u_{\rm st}(1,\tan\theta_{\rm w})|>|(\iu,0)-u_{\rm wk}(1,\tan\theta_{\rm w})|\qquad
\tx{for $0<\theta_{\rm w}<\theta^{(\irho,\iu)}_{\rm d}$};
\end{equation*}
that is, the shock strength of a steady Prandtl-Meyer weak reflection solution
is weaker than the steady strong one.
\end{definition}

The goal of this work is to prove the existence of global unsteady
Prandtl-Meyer reflection configurations for unsteady potential flow, determined by
Eq. \eqref{1-b2},
which converge to the steady Prandtl-Meyer weak reflection solution
as
$t$ tends to infinity for all possible physical
parameters $\gam\ge 1$, $\iu>c_\infty$,
and $\theta_{\rm w}\in(0, \theta^{(\irho,\iu)}_{\rm d})$.
Therefore, we consider the following initial-boundary value problem
for \eqref{1-b2}:

\begin{problemL}[Initial-boundary value problem]
\label{problem-1}
Given $\gam\ge 1$, fix $(\irho, \iu)$ with $\iu>c_\infty$.
For a fixed $\theta_{\rm w}\in(0,\theta_{\rm d}^{(\irho, \iu)})$,
let $W$ be given by \eqref{definition-W}.
Find a global weak solution $\Phi\in W^{1,\infty}_{\rm loc}(\R_+\times (\R^2\setminus W))$
of Eq. \eqref{1-b2} with $\rho$ determined by \eqref{1-b1}
and
\begin{equation}
\label{definition-bernoulli-constant}
  B=\frac{\iu^2}{2}+h(\irho)
\end{equation}
so that $\Phi$ satisfies both the initial condition at $t=0${\rm :}
\begin{equation}\label{1-d}
(\rho,\Phi)|_{t=0}=(\irho, \iu x_1) \qquad \text{for $(x_1,x_2)\in \R^2\setminus W$},
\end{equation}
and the slip boundary condition along the wedge boundary $\der W${\rm :}
\begin{equation}\label{1-e}
\nabla_{\bf x}\Phi\cdot {\bf n}_{\rm w} |_{\der W}=0\qquad\tx{for $t>0$},
\end{equation}
where ${\bf n}_{\rm w}$ is the exterior unit normal vector to $\der W$.
\end{problemL}

\begin{remark}\label{remark-for-problem1}
In particular, we seek a solution $\Phi\in W^{1,\infty}_{\rm loc}(\R_+\times (\R^2\setminus W))$
that converges to
the steady Prandtl-Meyer weak reflection solution $\bPhi$ when
$t$ tends to infinity
in the following sense{\rm :}
if $\bPhi$ is the steady Prandtl-Meyer weak reflection solution corresponding to
the fixed parameters $(\irho, \iu, \gam, \theta_{\rm w})$ in the sense of
Definition {\rm \ref{def-weak-strong-steady-new2015}} with $\bar{\rho}=h^{-1}(B-\frac 12|\nabla \bPhi|^2)$,
then, for any $R>0$, $\Phi$ satisfies
\begin{equation}
\label{time-asymp-lmt}
\lim_{t\to \infty} \left(\|\nabla_{\bf x}\Phi(t,\cdot)-\nabla_{\bf x}\bPhi\|_{L^1(B_R({\bf 0})\setminus W)}
+\|\rho(t,\cdot)-\bar{\rho}\|_{L^1(B_R({\bf 0})\setminus W)}\right)=0
\end{equation}
for $\rho(t,{\bf x})$ given by \eqref{1-b1}.
\end{remark}

\medskip
The definition of a weak solution of Problem \ref{problem-1} is given as follows:
\begin{definition}
\label{def-weaksol}
A function $\Phi\in W^{1,\infty}_{\rm loc}(\R_+\times (\R^2\setminus W))$ is called a weak solution
of Problem {\rm \ref{problem-1}}  if $\Phi$ satisfies the following properties{\rm :}

\smallskip
\begin{itemize}
\item[(i)] $B-\der_t\Phi-\frac 12|\nabla_{\bf x}\Phi|^2 >  h(0+)\, $ {\it a.e.} in $\R_+\times (\R^2\setminus W)${\rm ;}

\smallskip
\item[(ii)] $(\rho(\der_t\Phi, \nabla_{\bf x}\Phi),\rho(\der_t\Phi, \nabla_{\bf x}\Phi)|\nabla_{\bf x}\Phi|)
  \in \left(L^1_{\rm loc}(\R_+\times (\R^2\setminus W))\right)^2${\rm ;}

\smallskip
\item[(iii)] For every $\zeta\in C^{\infty}_{\rm c}(\R_+\times \R^2)$,
\begin{equation*}
\int_0^{\infty}\int_{\R^2\setminus W}\Bigl(\rho(\der_t\Phi, \nabla_{\bf x}\Phi)\der_t\zeta
+\rho(\der_t\Phi, \nabla_{\bf x}\Phi)\nabla_{\bf x}\Phi\cdot \nabla_{\bf x}\zeta\Bigr)\;\dd{\bf x} \dd t
+\int_{\R^2\setminus W} \irho \zeta(0, {\bf x}) \; \dd{\bf x}=0.
\end{equation*}
\end{itemize}
\end{definition}

\medskip
Since the initial data \eqref{1-d} does not satisfy the boundary condition \eqref{1-e},
a boundary layer is generated along the wedge boundary starting at $t=0$,
which is proved to form the Prandtl-Meyer reflection configuration in this monograph.

Notice that the initial-boundary value problem, Problem \ref{problem-1}, is invariant under the scaling
\begin{equation*}
(t, {\bf x})\rightarrow (\alp t, \alp{\bf x}),\quad (\rho, \Phi)\rightarrow (\rho, \frac{\Phi}{\alp})
\qquad\,\, \text{for}\;\;\alp\neq 0,
\end{equation*}
in the sense that, if $(\rho, \Phi)(t,{\bf x})$ is a solution,
then so is $(\til{\rho}, \til{\Phi})(t,{\bf x})=(\rho, \frac{\Phi}{\alp})(\alp t, \alp{\bf x})$.
Based on this observation, we look for {\emph{self-similar solutions}} of Problem \ref{problem-1} in the form
\begin{equation}
\label{definition-selfsimsol}
\rho(t, {\bf x})=\rho(\xxi),\quad \Phi(t, {\bf x})=t\phi(\xxi)\qquad\,\,\text{with}\;\; \xxi=(\xi_1,\xi_2)=\frac{{\bf x}}{t}\,\,\,\tx{for}\,\,t>0.
\end{equation}
For such $\phi$, introduce {\emph{the pseudo-potential function}} $\vphi$ given by
\begin{equation*}
\vphi=\phi-\frac 12|\xxi|^2.
\end{equation*}
If $\Phi$ solves \eqref{1-b2} with \eqref{1-b1}, then $\vphi$ satisfies the following {\emph{Euler equations for
the self-similar solutions}}:
\begin{align}
\label{1-r1}
&{\rm div}(\rho D\vphi)+2\rho=0,\\
\label{1-p1}
&\frac 12|D\vphi|^2+\vphi+ h(\rho)=B,
\end{align}
where the {\it divergence} ${\rm div}$ and {\it gradient} $D$ are with respect to the self-similar variables $\xxi\in \R^2$.
Solve \eqref{1-p1} first for $\rho$ and then substitute the result into \eqref{1-r1} to obtain
\begin{equation}
\label{2-1}
N(\vphi):=
{\rm div}\big(\rho(|D\vphi|^2,\vphi)D\vphi\big)+2\rho(|D\vphi|^2,\vphi)=0
\end{equation}
for
\begin{equation}
\label{1-o}
\rho(|D\vphi|^2,\vphi)
=\begin{cases}
\bigl(1+(\gam-1)(B-\frac 12|D\vphi|^2-\vphi)\bigr)^{\frac{1}{\gam-1}}&\tx{if}\;\;\gam>1,\\
\exp(B-\frac 12|\nabla\vphi|^2-\vphi)&\tx{if}\;\;\gam=1.
\end{cases}
\end{equation}
Note that the Bernoulli constant $B$ is given by \eqref{definition-bernoulli-constant}.

\medskip
The local sound speed $c=c(|D\vphi|^2, \vphi)>0$ for the pseudo-steady potential flow equation \eqref{2-1} is given by
\begin{equation}
\label{1-a1}
c^2(|D\vphi|^2,\vphi)=
1+(\gam-1)\big(B-\frac 12|D\vphi|^2-\vphi\big).
\end{equation}
Eq. \eqref{2-1} is a second-order nonlinear equation of mixed elliptic-hyperbolic type.
It is elliptic if and only if
\begin{equation}
\label{1-f}
|D\vphi|<c(|D\vphi|^2,\vphi) \,\, \Longleftrightarrow \,\,  |D\vphi|<\sqrt{\frac{2}{\gam+1}\bigl(1+(\gam-1)(B-\vphi)\bigr)}\quad \,\, \tx{(pseudo-subsonic)},
\end{equation}
and \eqref{2-1} is hyperbolic if and only if
\begin{equation*}
|D\vphi|>c(|D\vphi|^2,\vphi) \,\,\, \Longleftrightarrow \,\,\, |D\vphi|>\sqrt{\frac{2}{\gam+1}\bigl(1+(\gam-1)(B-\vphi)\bigr)}\quad\,\, \tx{(pseudo-supersonic)}.
\end{equation*}
In order to find a function $\vphi(\bmxi)$ such that $\Phi(t, {\bf x})$ with $\rho(t, {\bf x})$ given by \eqref{definition-selfsimsol}
is a solution of Problem \ref{problem-1} satisfying \eqref{time-asymp-lmt}, we make the following additional observations:

\smallskip
\begin{itemize}
\item[(i)] {\emph{Symmetric domain}}:
Since the solid wedge $W$ is symmetric with respect to the axis $x_2=0$, it suffices to consider Problem \ref{problem-1}
in the upper half-plane $\{(x_1,x_2)\in \R^2\,:\,x_2>0\}$. In the self-similar plane, define
\begin{equation}
\label{Lambda}
\tdomain:=\{\xxi\in\R^2\,:\,\xi_2>0\}\setminus\{\xxi\,:\,\xi_2\le \xi_1\tan\theta_{\rm w}, \xi_1\ge 0\}.
\end{equation}
Then Problem \ref{problem-1} is reformulated as a boundary value problem in $\tdomain$.

\smallskip
\item[(ii)] {\emph{Initial condition}}:
For each ${\bf x}\in \R^2\setminus (W\cup\{\mathbf{0}\})$,
$|{\bm\xi}|=|\frac{{\bf x}}{t}|\rightarrow \infty$ as $t\to 0+$.
This means that the initial condition \eqref{1-d} in Problem \ref{problem-1}
becomes an asymptotic boundary condition in the self-similar variables.

\smallskip
\item[(iii)] {\emph{Time-asymptotic limit}}: For each ${\bf x}\in \R^2\setminus W$,
$|{\bm\xi}|=|\frac{{\bf x}}{t}|\rightarrow 0$ as $t\to \infty$.
To find a global weak solution of Problem \ref{problem-1} satisfying \eqref{time-asymp-lmt},
we seek a self-similar weak solution $\vphi({\bm\xi})$ satisfying
\begin{equation*}
\lim_{R\to 0+} \frac{1}{|B_R({\bf 0})\cap \tdomain|}\int_{B_R({\bf 0})\cap \tdomain}|\nabla_{\bm\xi} \vphi-\nabla_{{\bf x}}\bPhi|\, \dd\xxi=0,
\end{equation*}
where $\bPhi$ is the steady Prandtl-Meyer weak reflection solution
of problem \eqref{eqn-steady-pot-new2015}--\eqref{bvp-steady-new2015},
and $|B_R({\bf 0})\cap \tdomain|$ is the area of $B_R({\bf 0})\cap \tdomain$.

\smallskip
\item[(iv)] \emph{Constant density state}: If $\rho>0$ is a constant in \eqref{1-r1}--\eqref{1-p1},
then the corresponding pseudo-potential $\vphi$ is given in the form
\begin{equation}
\label{cd-psp}
\vphi(\xxi)=-\frac 12|\xxi|^2+ (u,v)\cdot\xxi+k
\end{equation}
for some constant state $(u,v)$ and a constant $k$.
In Problem \ref{problem-1}, the initial state has a constant
density $\irho>0$ and a constant velocity $(\iu,0)$.
Then the corresponding pseudo-potential $\ivphi$ in the self-similar variables is given by
\begin{equation}
\label{1-m}
\ivphi=-\frac 12|\xxi|^2+ (\iu,0)\cdot\xxi+\ik
\end{equation}
for a constant $\ik$. It follows from \eqref{definition-bernoulli-constant} that
$\ik=0.$
\end{itemize}

Hereafter, we assume without loss of generality that $\rho_\infty=1$, so that $c_\infty=1$.
This can be achieved by
the scaling
\begin{equation*}
\xxi\mapsto c_\infty\xxi,
\qquad  (\rho, \vphi, u_\infty)\rightarrow(\frac{\rho}{\irho}, \frac{\vphi}{c^2_\infty},
\frac{\iu}{c_\infty})
\end{equation*}
for any $\gam\ge 1$.

\medskip
Given $\gam\ge 1$, $\irho=1$, and $\iu>1$,
we now reformulate Problem \ref{problem-1} in the self-similar variables.
Hereafter, we denote $(\theta_{\rm d}^{(\irho, \iu)}, \theta_{\rm s}^{(\irho, \iu)})$
by $(\theta_{\rm d}^{(\iu)}, \theta_{\rm s}^{(\iu)})$,
since $\irho$ is fixed as $1$.

Taking into account the additional observations stated above,
we reformulate Problem \ref{problem-1} as a boundary value problem in the self-similar variables.

\begin{figure}[htp]
\centering
\begin{psfrags}
\psfrag{tw}[cc][][0.7][0]{$\theta_{\rm w}$}
\psfrag{sb}[cc][][0.7][0]{$\nabla\vphi\cdot{\bf n}_{\rm w}=0$}
\psfrag{conv}[cc][][0.7][0]{$\phantom{aa}R_{\theta}=\{(\xi_1,\xi_2)\,:\,\xi_1=\xi_2\cot\theta,\xi_2>0\}$}
\psfrag{t}[cc][][0.7][0]{$\theta$}
\includegraphics[scale=0.8]{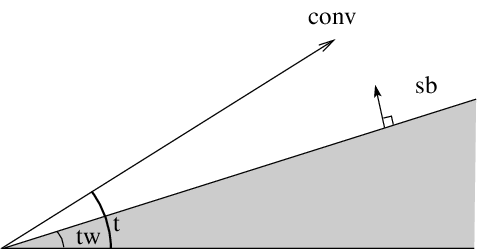}
\caption{Asymptotic boundary condition at infinity}\label{Figure-2a}
\end{psfrags}
\end{figure}

\begin{problemL}[Boundary value problem in the self-similar variables $\xxi$]\label{problem-2}
Given $\gam\ge 1$, $\iu>1$, and $\theta_{\rm w}\in(0, \theta_{\rm d}^{(\iu)})$,
find a weak solution $\vphi\in W^{1,\infty}(\tdomain)$ of Eq. \eqref{2-1} in  $\tdomain$
satisfying the following conditions{\rm :}

\smallskip
\begin{itemize}
\item[(i)] {\emph{Slip boundary condition on $\Wedge$}}{\rm :}
\begin{equation}
\label{1-k}
D\vphi\cdot{\bf n}_{\rm w}=0\qquad\;\text{on}\;\;\Wedge=\{\xxi\,:\,\xi_2=\xi_1\tan\theta_{\rm w}, \xi_1>0\},
\end{equation}
where ${\bf n}_{\rm w}$ represents the exterior unit normal vector to the wedge boundary $\Wedge${\rm ;}

\smallskip
\item[(ii)] {\emph{Time-asymptotic limit condition in the self-similar variables}}{\rm :}
\begin{equation}
\label{time-symp-lmi}
\lim_{R\to 0+} \frac{1}{|B_R({\bf 0})\cap \tdomain|}\int_{B_R({\bf 0})\cap \tdomain}|\nabla_{\bm\xi} \vphi-\nabla_{{\bf x}}\bPhi|\, \dd\xxi=0,
\end{equation}
where $\bPhi$ is the steady Prandtl-Meyer weak reflection solution corresponding to $\theta_{\rm w}${\rm ;}

\smallskip
\item[(iii)] {\emph{Asymptotic boundary condition at infinity}}{\rm :}
For each $\theta\in(\theta_{\rm w}, \pi]$,
\begin{equation}\label{1-k-c}
\lim_{r\to \infty} \|\varphi  - \varphi_\infty \|_{C(R_{\theta}\setminus B_r({\bf 0}))} = 0
\end{equation}
for each ray $R_\theta:=\{ \xi_1=\xi_2 \cot
\theta, \xi_2 > 0 \}$; see Fig. {\rm \ref{Figure-2a}}.
\end{itemize}
\end{problemL}

\begin{definition}
\label{definition-weak-sol}
A function $\vphi\in W^{1,1}_{\rm loc}(\tdomain)$ is called a weak solution of Problem {\rm \ref{problem-2}}
if $\vphi$ satisfies conditions {\rm (i)}--{\rm (iii)} of Problem {\rm \ref{problem-2}}
and the following additional properties{\rm :}

\smallskip
\begin{itemize}
\item[(i)] $\rho(|D\vphi|^2,\vphi)> 0$ $\,\,$ {\it a.e.} in $\tdomain${\rm ;}

\smallskip
\item[(ii)] $(\rho(|D\vphi|^2,\vphi),\rho(|D\vphi|^2,\vphi)|D\vphi|)\in L^1_{\rm loc}(\tdomain)${\rm ;}

\smallskip
\item[(iii)] For every $\zeta\in C^{\infty}_{\rm c}(\R^2)$,
\begin{equation}
\label{1-q}
\int_{\tdomain} \big(\rho(|D\vphi|^2,\vphi)D\vphi\cdot D\zeta-2\rho(|D\vphi|^2,\vphi)\zeta\big)\;\dd\xxi=0.
\end{equation}
\end{itemize}
\end{definition}

For $\rho >0$, note that \eqref{1-k} is equivalent to the conormal boundary condition:
\begin{equation}
\label{1-l}
\rho D\vphi\cdot {\bf n}_{\rm w}=0\;\qquad\text{on}\;\;\Wedge.
\end{equation}

Condition (ii) of Problem \ref{problem-2} indicates that a solution of Problem \ref{problem-2}
converges
to a steady potential flow with a shock near the wedge vertex.
To find such a solution, we define an entropy solution of Problem \ref{problem-2} with a shock.
The definition is given in a way similar to Definition \ref{def-weaksol-steady-new2015}.

\begin{definition}
\label{def-shocksol-self-similar-new2015}
Let $\Gam_{\rm sh}$ be a $C^1$--curve that lies in $\tdomain$ and divides
$\tdomain$ into two subdomains{\rm :} $\tdomain^-$ and $\tdomain^+$.
A weak solution $\vphi$ of Problem {\rm \ref{problem-2}} is an {\emph{entropy solution}} with a shock $\Gam_{\rm sh}$
if $\vphi$ satisfies the following properties{\rm :}

\smallskip
\begin{itemize}
\item[(i)] $\vphi\in W^{1,\infty}_{\rm loc}(\tdomain)${\rm ;}

\smallskip
\item[(ii)]
$\vphi\in C^1_{\rm loc}(\ol{\tdomain^{\pm}})\cap C^2(\tdomain^{\pm})${\rm ;}

\smallskip
\item[(iii)] For $\vphi^+:=\vphi|_{\tdomain^+\cup \Gam_{\rm sh}}$
and $\vphi^-:=\vphi|_{\tdomain^-\cup\Gam_{\rm sh}}$,
\begin{equation*}
\der_{{\bf n}_{\rm sh}} {\vphi}^{-}>\der_{{\bf n}_{\rm sh}} {\vphi}^+>0
\qquad \tx{on $\Gam_{\rm sh}$},
\end{equation*}
where ${\bf n}_{\rm sh}$ represents a unit normal vector to $\Gam_{\rm sh}$
pointing from $\tdomain^-$ towards $\tdomain^+${\rm ;}

\smallskip
\item[(iv)] $\vphi$
satisfies
the Rankine-Hugoniot jump conditions on $\Gam_{\rm sh}${\rm :}
\begin{equation}
\label{1-h}
[\vphi]_{\Gam_{\rm sh}}=[\rho(|D\vphi|^2, \vphi)D\vphi\cdot{\bf n}_{\rm sh}]_{\Gam_{\rm sh}}=0,
\end{equation}
which is similar to the steady case
\eqref{eqn-steady-pot-new2015}.
\end{itemize}

If ${\bf n}_{\rm sh}=\frac{D\vphi^--D\vphi^+}{|D\vphi^--D\vphi^+|}$ is oriented
so that $\der_{{\bf n}_{\rm sh}}\vphi^\pm>0$, and if $\der_{{\bf n}_{\rm sh}}\vphi^->\der_{{\bf n}_{\rm sh}}\vphi^+$
holds on $\Gam_{\rm sh}$,
the shock solution is said to satisfy the entropy condition. By \eqref{1-h},
the entropy condition is equivalent to $\rho(|D\vphi^-|^2, \vphi^-)< \rho(|D\vphi^+|^2, \vphi^+)$ on $\Gam_{\rm sh}$.
\end{definition}

\section{Structure of Solutions of Problem 2.9}
\label{subsec-const-density}

Given $\gam\ge 1$, $\irho=1$, and $\iu>1$, fix $\theta_{\rm w}\in(0,\theta_{\rm d}^{(\iu)})$.

\subsection{Near the origin}
\label{subsubsec-2-1}
We seek a solution $\vphi$ of Problem \ref{problem-2} so that the solution at the origin coincides with
the steady Prandtl-Meyer weak reflection solution corresponding to
parameters $(1, \iu, \gam, \theta_{\rm w})$ in the sense of Definition \ref{def-weak-strong-steady-new2015}.
For $\ivphi$ given by \eqref{1-m}, define
\begin{eqnarray}
\label{newlabel-N0}
&\leftvphiop=-\frac 12|\xxi|^2+(\leftuop,\leftvop)\cdot\xxi, \qquad
\leftshockop=\{\bmxi\in \tdomain\,:\,\leftvphiop(\bmxi)=\ivphi(\bmxi)\}.
\label{1-n}
\end{eqnarray}
Choose the constant vector $(u_0, v_0)$ as
\begin{equation}
\label{N}
(u_0, v_0)=u_{\rm wk}^{(\theta_{\rm w})}(1,\tan\theta_{\rm w}),
\end{equation}
and define
\begin{equation*}
\bar{\vphi}(\bmxi):=\max\{\ivphi(\bmxi), \leftvphiop(\bmxi)\}.
\end{equation*}
Then $\vphi:=\bar{\vphi}$ satisfies
\eqref{1-k}--\eqref{time-symp-lmi}
and \eqref{1-h} with $\shock=\leftshockop$.

For the nonlinear differential operator $N$ defined by \eqref{2-1},
equation $N(\leftvphiop)=0$ introduces {\emph{the pseudo-sonic circle}} $\der B_{c_0}(\leftuop, \leftvop)$
with $c_0^2=\rho_0^{\gam-1}$ for $\rho_0=\rho(|D\leftvphiop|^2, \leftvphiop)$ in the following-sense:

\smallskip
\begin{itemize}
\item $N(\leftvphiop)=0$ is elliptic in $B_{c_0}(\leftuop, \leftvop)$,

\smallskip
\item $N(\leftvphiop)=0$ is hyperbolic in $\R^2\setminus \ol{B_{c_0}(\leftuop, \leftvop)}$.
\end{itemize}

\begin{remark}\label{remark-origin-inclusion}
Let $\theta_{\rm s}^{(\iu)}$ be from Lemma {\rm \ref{lemma-spolar-steady-new}(c)}.
Then the wedge vertex ${O}=(0,0)$ satisfies the following{\rm :}

\smallskip
\begin{itemize}
\item ${O}\in \R^2\setminus \ol{B_{c_0}(u_0, v_0)}\,\,\tx{for}\,\,0<\theta_{\rm w}<\theta_{\rm s}^{(\iu)}${\rm ,}

\smallskip
\item ${O}\in \der B_{c_0}(u_0,v_0)\,\,\tx{at}\,\,\theta_{\rm w}=\theta_{\rm s}^{(\iu)}${\rm ,}

\smallskip
\item ${O}\in B_{c_0}(u_0, v_0)\,\,\tx{for}\,\,\theta_{\rm s}^{(\iu)}<\theta_{\rm w}<\theta_{\rm d}^{(\iu)}$.
\end{itemize}
\end{remark}

\subsection{Away from the origin }
To determine a solution $\vphi$ of Problem \ref{problem-2}, we look for a solution $\vphi$ with a piecewise constant
density $\rho(|D\vphi|^2, \vphi)$,
defined by \eqref{1-o} in $\tdomain \setminus B_R({ O})$ for some sufficiently large $R>0$,
so that such a solution $\vphi$ satisfies the asymptotic boundary condition (iii) of Problem \ref{problem-2}.
For this purpose, we introduce a straight shock solution in $\tdomain \setminus B_R({O})$.
In fact, the only straight shock solution that satisfies \eqref{1-k-c} is a normal shock solution.
This can be seen more clearly  in \S \ref{section-change-par}.
We now compute the normal shock solution and discuss its useful properties.

To compute the normal shock, denoted by $S_1$, and the corresponding pseudo-potential
$\rightvphiop$ below $\rightshockop$, it is convenient to rotate the self-similar plane
by angle $\theta_{\rm w}$ clockwise.
\begin{figure}[htp]
\centering
\begin{psfrags}
\psfrag{y}[cc][][0.8][0]{$\xi_2$}
\psfrag{x}[cc][][0.8][0]{$\xi_1$}
\psfrag{tw}[cc][][0.7][0]{$\theta_{\rm w}$}
\psfrag{pn}[cc][][0.7][0]{$(\iu\cos\theta_{\rm w},0)$}
\psfrag{pi}[cc][][0.7][0]{$\phantom{aaaaaaa}\iu(\cos\theta_{\rm w},-\sin\theta_{\rm w})$}
\psfrag{po}[cc][][0.7][0]{$(0,-\iu\sin\theta_{\rm w})$}
\psfrag{en}[cc][][0.7][0]{$\xi_2^{(1)}$}
\psfrag{sn}[cc][][0.7][0]{$\rightshockop$}
\psfrag{w}[cc][][0.7][0]{$\Wedge$}
\includegraphics[scale=.8]{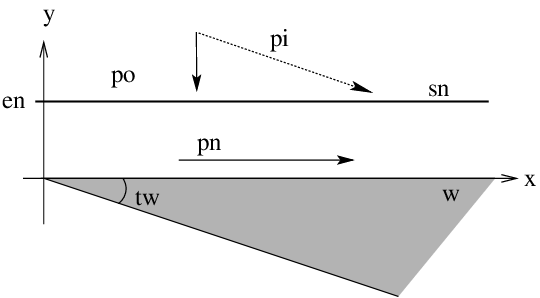}
\caption{The normal shock}
\label{fig:normalshock}
\end{psfrags}
\end{figure}
In the rotated self-similar plane, $\ivphi$ in \eqref{1-m} is written as
\begin{equation*}
\ivphi=-\frac 12|\xxi|^2+\iu (\cos\theta_{\rm w},-\sin\theta_{\rm w})\cdot\xxi.
\end{equation*}
Then $\rightvphiop$ is in the form
\begin{equation*}
\rightvphiop=-\frac 12|\xxi|^2+\iu (\cos\theta_{\rm w},-\sin\theta_{\rm w})\cdot (\xi_1, \xi_2^{(1)}),
\end{equation*}
where $\netaop$ is the distance of $\rightshockop$ from $\Wedge$.
Denote
\begin{equation}
\label{1-22}
\iv:=\iu \sin\theta_{\rm w}.
\end{equation}
It follows from \eqref{1-o} and \eqref{1-h} that density $\rightrhoop$ and distance $\netaop$ satisfy
\begin{align}
\label{1-23}
&\netaop=\frac{\iv}{\rightrhoop-1},\\
\label{1-1a}
&h(\rightrhoop)-h(1)=\frac 12\iv^2+\netaop \iv,
\end{align}
where $h(\rho)$ is defined by \eqref{def-h-and-c}.

Consider
$$
F(\rho):=\big(h(\rho)-h(1)\big)(\rho-1)-\frac 12 (\rho-1)\iv^2-\iv^2.
$$
A direct computation shows that $F(1)=-\iv^2<0$, $\displaystyle{\lim_{\rho\to \infty} F(\rho)=\infty}$, $F'(1)=-\frac 12\iv^2<0$,
and $F''(\rho)>0$ whenever $\rho\ge 1$. This implies that there exists a unique $\rightrhoop\in (1,\infty)$
such that $F(\rightrhoop)=0$. Then \eqref{1-23} yields that $\netaop>0$.
Rotating the self-similar plane back by angle $\theta_{\rm w}$ counterclockwise, we find that $\rightvphiop$ is given by
\begin{equation}
\label{1-rightvphi}
\rightvphiop=-\frac 12|\xxi|^2+\iu\cos\theta_{\rm w}(\cos\theta_{\rm w},\sin\theta_{\rm w})\cdot\xxi
-\iu\netaop\sin\theta_{\rm w},
\end{equation}
and the normal shock $\rightshockop$ by
$$
 \rightshockop=\{\xxi\,:\,\ivphi(\xxi)=\rightvphiop(\xxi)\}
 =\{\xxi\,:\,\xi_2=\xi_1\tan\theta_{\rm w}+\netaop\sec\theta_{\rm w}\}.
$$

\begin{lemma}\label{remark-1-1}
For any given $\iu>1$ and the wedge angle $\theta_{\rm w}\in(0,\theta_{\rm d}^{(\iu)})$,
$$
{\rm dist}(\rightshockop,\Wedge)<\rightcop:=\rightrhoop^{(\gam-1)/2}.
$$

\begin{proof}
By the mean value theorem,
there exists a constant $\rho_*\in(1,\rho_1)$ satisfying
\begin{equation*}
h(\rightrhoop)-h(1)=\mu(\rightrhoop-1) \qquad\,\, \tx{for $\mu=\rho_*^{\gam-2}$}.
\end{equation*}
Then $F(\rho_1)=0$ implies that
\begin{align*}
&\mu(\rightrhoop-1)^2-\frac 12\iv^2(\rightrhoop-1)- \iv^2=0
\quad\Longrightarrow \quad \rightrhoop-1=\frac{\frac 12 \iv^2+\sqrt{\iv^2(\frac 14 \iv^2+4\mu)}}{2\mu}.
\end{align*}
Since $\iv>0$, \eqref{1-23} yields that
\begin{equation*}
\netaop=\frac{4\mu}{\sqrt{16\mu+\iv^2}+\iv}\le  \sqrt{\mu}.
\end{equation*}
By the definition of $\mu$ above, it can directly be checked that
\begin{equation*}
\sqrt{\mu}<
\begin{cases}
\sqrt{\rightrhoop^{\gam-2}}< \sqrt{\rightrhoop^{\gam-1}}=\rightcop\;\;\;\;&\text{if}\;\gam\ge 2,\\[1mm]
1< \sqrt{\rightrhoop^{\gam-1}}=\rightcop\;\;&\text{if}\;1<\gam<2,\\
1=\rightcop\;\;&\text{if}\;\;\gam=1,
\end{cases}
\end{equation*}
which implies that $\netaop<\rightcop$.
\end{proof}
\end{lemma}
\medskip

\subsection{Global configurations of the solutions of Problem {\rm \ref{problem-2}}}
\label{subsubsec-global-struction-new}

Following Remark \ref{remark-origin-inclusion},
our desired solution of Problem \ref{problem-2} has two different configurations
depending on the two different intervals of the wedge angle: $\theta_{\rm w}\in(0, \theta_{\rm s}^{(\iu)})$ and
$\theta_{\rm w}\in [\theta_{\rm s}^{(\iu)}, \theta_{\rm d}^{(\iu)})$.

\smallskip
{\textbf{Case I}}. Fix $\theta_{\rm w}\in(0, \theta_{\rm s}^{(\iu)})$.
\begin{figure}[htp]
\centering
\begin{psfrags}
\psfrag{ls}[cc][][0.8][0]{$\leftshockop$}
\psfrag{sn}[cc][][0.8][0]{$\rightshockop$}
\psfrag{ol}[cc][][0.8][0]{$\oOmop$}
\psfrag{on}[cc][][0.8][0]{$\nOmop$}
\psfrag{lsn}[cc][][0.8][0]{$\phantom{aa}\leftsonicop$}
\psfrag{rsn}[cc][][0.8][0]{$\rightsonicop\phantom{aaaaaaa}$}
\psfrag{tw}[cc][][0.8][0]{$\theta_{\rm w}$}
\psfrag{o}[cc][][0.8][0]{$O$}
\psfrag{lb}[cc][][0.8][0]{$\leftbottom$}
\psfrag{rb}[cc][][0.8][0]{$\rightbottom$}
\psfrag{lt}[cc][][0.8][0]{$\lefttop$}
\psfrag{rt}[cc][][0.8][0]{$\righttop$}
\psfrag{s}[cc][][0.8][0]{$\shock$}
\psfrag{om}[cc][][0.8][0]{$\Om$}
\psfrag{q0}[cc][][0.8][0]{$Q_0$}
\psfrag{q1}[cc][][0.8][0]{$Q_1$}
\psfrag{w}[cc][][0.8][0]{$\om_0$}
\psfrag{z}[cc][][0.8][0]{$\om_1$}
\includegraphics[scale=.8]{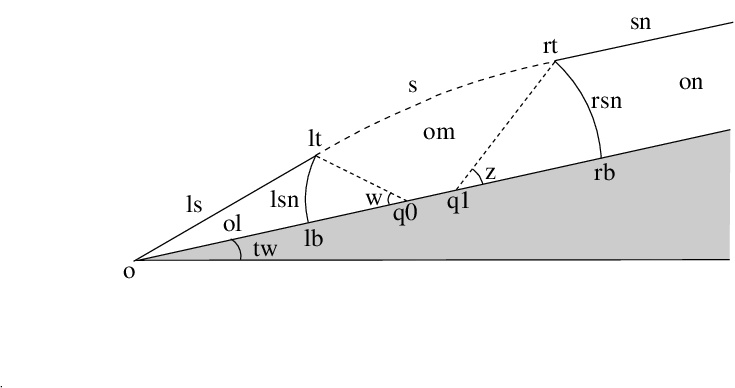}
\vspace{-18mm}
\caption{Admissible solutions for $0<\theta_{\rm w}<\theta_{\rm s}^{(\iu)}$}\label{fig:globala}
\end{psfrags}
\end{figure}
Let $\leftvphiop$ and $\rightvphiop$ be defined by \eqref{newlabel-N0} and \eqref{1-rightvphi}, respectively.
Define $\leftPop:=D\leftvphiop(O)$ and $\rightPop:=D\rightvphiop(O)$.
Consider two sonic circles $\der B_{\leftcop}(\leftPop)$ and $\der B_{\rightcop}(\rightPop)$.

{\emph{The left sonic arc}}: The sonic circle $\der B_{\leftcop}(\leftPop)$ and the straight oblique
shock $\leftshockop:=\{\xxi\,:\,\leftvphiop(\xxi)=\ivphi(\xxi)\}$
intersect at two points in $\tdomain$, which will be verified in detail
in \S \ref{section-change-par}.
Let $\lefttop$ be the intersection with the smaller $\xi_2$--coordinate.
Also, $\der B_{\leftcop}(\leftPop)$ intersects with $\Wedge$ at two points;
let $\leftbottom$ be the intersection point with the smaller $\xi_2$--coordinate.
Denote $\om_0:=\angle \leftbottom \leftPop \lefttop\in(0,\pi)$.
We define
\begin{equation*}
\leftsonicop:=\{P\in \der B_{\leftcop}(\leftPop)\,:\,0\le \angle \leftbottom \leftPop P\le \om_0\},
\end{equation*}
which is a closed subset of $\der B_{\leftcop}(\leftPop)$.
We call $\leftsonicop$ {\emph{the sonic arc corresponding to $\leftvphiop$. }}

{\emph{The right sonic arc}}:
By Lemma \ref{remark-1-1}, the sonic circle $\der B_{\rightcop}(\rightPop)$ and the normal
shock $\rightshockop=\{\xxi\,:\,\rightvphiop(\xxi)=\ivphi(\xxi)\}$
intersect at two distinct points; let $\righttop$ be the intersection point
with the larger $\xi_2$--coordinate.
Also, $\der B_{c_1}(Q_1)$ intersects with $\Wedge$ at two distinct points;
let $\rightbottom$ be the intersection point with the larger $\xi_2$--coordinate.
Denote $\om_1:=\angle \rightbottom \rightPop \righttop\in(0,\pi)$. We define
\begin{equation*}
\rightsonicop:=\{P'\in \der B_{\rightcop}(\rightPop)\,:\,0\le \angle \rightbottom \rightPop P'\le \om_1\},
\end{equation*}
which is a closed subset of $\der B_{\rightcop}(\rightPop)$,
similar to  $\leftsonicop$.
We call $\rightsonicop$ {\emph{the sonic arc corresponding to $\rightvphiop$. }}

For each $j=1,\cdots, 4$, let $\xxi^{P_j}=(\xi_1^{P_j}, \xi_2^{P_j})$ denote the $\xxi$--coordinate of point $P_j$.
Let $\leftshocksegop$ be the line segment $\overline{O\lefttop}$, and let $\oOmop\subset \tdomain$ be the open set
enclosed by $\leftshocksegop$, $\leftsonicop$, and the line segment $\ol{O\leftbottom}$.
Next, let $\rightshocksegop$ be the portion of $\rightshockop$ with the left endpoint $\righttop$,
and let $\nOmop \subset \tdomain$ be the unbounded open set enclosed by $\rightshocksegop$, $\rightsonicop$,
and
$\Wedge\cap\{\xi_2\ge \xi_2^{\rightbottom}\}$.

Our goal is to find a curved shock $\shock$ that connects $\lefttop$ with $\righttop$
and a solution $\vphi$ of Problem \ref{problem-2} to satisfy both \eqref{1-f}
in the open region $\Om$ (enclosed by $\shock$, $\rightsonicop$, $\overline{\leftbottom \rightbottom}$,
and $\leftsonicop$) and
\begin{equation*}
\vphi=\begin{cases}
\leftvphiop&\tx{in $\oOmop$},\\
\rightvphiop&\tx{in $\nOmop$},\\
\ivphi&\tx{in $\tdomain\setminus \ol{\oOmop\cup\Om\cup\nOmop}$}.
\end{cases}
\end{equation*}
Problem \ref{problem-2} is now a free boundary problem
given in a bounded region $\Om$ with a free boundary $\shock$ to be determined simultaneously with $\vphi$.

\smallskip
{\textbf{Case 2}}. Fix $\theta_{\rm w}\in[\theta_{\rm s}^{(\iu)}, \theta_{\rm d}^{(\iu)})$.
The right sonic arc $\rightsonicop$ is given in the same way as {\bf Case 1}.
By Remark \ref{remark-origin-inclusion}, since the triangular region $\oOmop$ in Fig. \ref{fig:globala} shrinks
to the origin as $\theta_{\rm w}\in(0, \theta_{\rm s}^{(\iu)})$ increases up to $\theta_{\rm s}^{(\iu)}$,
we look for a curved shock $\shock$ that connects origin $O$ with $\righttop$ for $\theta_{\rm w}\ge \theta_{\rm s}^{(\iu)}$
and
\begin{figure}[htp]
\centering
\begin{psfrags}
\psfrag{C}[cc][][0.8][0]{$c_*$}
\psfrag{TC}[cc][][0.8][0]{$\theta_{\rm c}$}
\psfrag{tw}[cc][][0.8][0]{$\theta_{\rm w}$}
\psfrag{U}[cc][][0.8][0]{$u$}
\psfrag{V}[cc][][0.8][0]{$v$}
\psfrag{U0}[cc][][0.8][0]{$u_0$}
\psfrag{oi}[cc][][0.8][0]{$\phantom{aaaaaaa}\Om_{\infty}: \rho_\infty=1, (\iu, 0)$}
\psfrag{TH}[cc][][0.8][0]{$\theta_{\rm w}$}
\psfrag{om}[cc][][0.8][0]{$\Omega$}
\psfrag{ol}[cc][][0.8][0]{$\oOmop$}
\psfrag{on}[cc][][0.8][0]{$\nOmop$}
\psfrag{ls}[cc][][0.8][0]{$\leftshockop$}
\psfrag{rsn}[cc][][0.8][0]{$\phantom{aa}\rightsonicop$}
\psfrag{rt}[cc][][0.8][0]{$\righttop$}
\psfrag{rb}[cc][][0.8][0]{$\rightbottom$}
\psfrag{s}[cc][][0.8][0]{$\shock$}
\psfrag{sn}[cc][][0.8][0]{$\rightshockop$}
\psfrag{O}[cc][][0.8][0]{$\phantom{a}O$}
\includegraphics[scale=0.9]{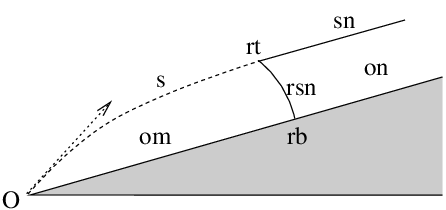}
\caption{Admissible solutions for $\theta_{\rm w}\ge\theta_{\rm s}^{(\iu)}$}\label{fig:global-subsonic-new}
\end{psfrags}
\end{figure}
a solution $\vphi$ to satisfy both \eqref{1-f} in the triangular domain $\Om$ (enclosed
by $\shock$, $\rightsonicop$, and the line segment $\ol{O\rightbottom}$) and
\begin{equation*}
\vphi=\begin{cases}
\rightvphiop&\tx{in $\nOmop$},\\
\ivphi&\tx{in $\tdomain\setminus \ol{\Om\cup\nOmop}$},
\end{cases}
\end{equation*}
with
\begin{equation}
\label{cond-op-origin-new}
\lim_{\overset {|P|\rightarrow 0}{P\in\Om}}\vphi(P)=\leftvphiop(O), \qquad\,\,
\lim_{\overset {|P|\rightarrow 0}{P\in\Om}}D\vphi(P)=D\leftvphiop(O).
\end{equation}
The condition that $\vphi=\leftvphiop$ in $\oOmop$ for $\theta_{\rm w}<\theta_{\rm s}^{(\iu)}$
is replaced by \eqref{cond-op-origin-new} so that our desired solution still satisfies \eqref{time-symp-lmi}.

\section{Main Theorems}
\label{subsection-mainthm-ori-par}
Fix $\gam\ge 1$ and  $\iu>1$. For each $\theta_{\rm w}\in(0, \theta_{\rm d}^{(\iu)})$, let $u_0$ be given by \eqref{N}.
By Lemmas \ref{lemma-appendix} and \ref{lemma-app2},
$u_0$ decreases with respect to $\theta_{\rm w}$. Define
\begin{equation*}
u_{\rm d}^{(\iu)}:=\lim_{\theta_{\rm w}\rightarrow \theta_{\rm d}^{(\iu)}-}u_0,
\qquad\,\, u_{\rm s}^{(u_{\infty})}:=\lim_{\theta_{\rm w}\rightarrow \theta_{\rm s}^{(u_{\infty})}} u_0.
\end{equation*}
For each $\iu>1$,  define an open interval $ I^{(\iu)}=(u_N^{(\iu)}, \iu)$, where $u_N^{(\iu)}$ is from Lemma \ref{lemma-app2}.
Given $\gam\ge 1$, we introduce a set of parameters
\begin{equation*}
\mathfrak{P}=\underset{\iu>1}{\cup}\{\iu\}\times I^{(\iu)}.
\end{equation*}
Then $\mathfrak{P}$ consists of three disjoint
sets $\mathfrak{P}_{\rm weak}$, $\mathfrak{P}_{\rm detach}$, and $\mathfrak{P}_{\rm strong}$:
\begin{equation}
\label{2-4-a7}
\begin{split}
&\mathfrak{P}_{\rm weak}=\underset{\iu>1}{\cup}\{\iu\}\times(u^{(\iu)}_{\rm d},\iu),\\
&\mathfrak{P}_{\rm detach}=\{(\iu, u_{\rm d}^{(\iu)})\,:\,\iu>1\},\\
&\mathfrak{P}_{\rm strong}=\underset{\iu>1}{\cup}\{\iu\}\times(u_N^{(\iu)},u^{(\iu)}_{\rm d}).
\end{split}
\end{equation}
\medskip
Our goal is to prove the existence of a global weak solution
of Problem \ref{problem-2}, satisfying the entropy condition,
for each $(\iu, u_0)\in\mathfrak{P}_{\rm weak}$
so that,
if $\theta_{\rm w}<\theta_{\rm s}^{(\iu)}$,
the solution has the configuration
of Fig. \ref{fig:globala} and,
if $\theta_{\rm w}\ge \theta_{\rm s}^{(\iu)}$,  the solution
has the configuration
of Fig. \ref{fig:global-subsonic-new}.
We first give a definition of {\emph{admissible solutions}} of Problem \ref{problem-2}
for $(\iu, u_0)\in \mathfrak{P}_{\rm weak}$.

\begin{definition}[Admissible solutions]\label{def-regular-sol-op}
Given $\gam\ge 1$, $\iu> 1$, and $(\iu,\leftuop)\in\mathfrak{P}_{\rm weak}$,
define $\theta_{\rm w}$ as
\begin{equation}\label{def-v-2015}
\tan\theta_{\rm w}=\frac{f_{\rm polar}(\leftuop)}{\leftuop},
\end{equation}
where
$f_{\rm polar}$ is determined in Lemma {\rm \ref{lemma-app2}}.
Let $\tdomain$ be the domain defined by \eqref{Lambda}, and
let $\leftvphiop$ and $\rightvphiop$ be defined by \eqref{newlabel-N0} and \eqref{1-rightvphi}, respectively.
A weak solution $\vphi\in C^{0,1}(\tdomain)$ of Problem {\rm \ref{problem-2}}
is called \emph{an admissible solution of Problem} \ref{problem-2}
if $\vphi$ satisfies the following properties{\rm :}

\medskip
\noindent
\textbf{Case {\rm I}}. $\leftuop > u_{\rm s}^{(\iu)}$, or equivalently,
$\theta_{\rm w}\in (0,\theta_{\rm s}^{(\iu)})${\rm :}

\smallskip
\begin{enumerate}
\item[\rm (i)] There exists a shock curve $\shock$ with
endpoints $\lefttop=(\xi_1^{\lefttop},\xi_2^{\lefttop})$ and $\righttop=(\xi_1^{\righttop},\xi_2^{\righttop})$
such that the following properties hold{\rm :}

\smallskip
\begin{itemize}
\item[(i-1)] Curve $\shock$ satisfies
\begin{equation}
\label{definition-1-op}
\shock\subset \tdomain\setminus\ol{B_{1}(\iu,0)},
\end{equation}
where $\partial B_{1}(\iu,0)$ is the sonic circle of the state in $\iOm:=\tdomain\setminus \ol{\oOmop\cup\nOmop\cup\Om}${\rm ;}

\smallskip
\item[(i-2)] Curve $\shock$ is $C^{2}$ in its relative interior. That is, for any $P\in \shock\setminus\{\lefttop,\righttop\}$,
 there exist $r>0$, $f\in C^{2}$, and an orthogonal coordinate system $(S,T)$ in $\R^2$
  such that $\shock\cap B_r(P)=\{S=f(T)\}\cap B_r(P)${\rm ;}

\smallskip
\item[(i-3)] Curve $\ol{\leftshocksegop\cup\shock\cup\rightshocksegop}$ is $C^1$, including at points $\lefttop$ and $\righttop${\rm ;}

\smallskip
\item[(i-4)] $\shock, \rightsonicop, \leftsonicop$, and
$\Wedge:=\{\xi_2=\xi_1\tan\theta_{\rm w},\xi_2\ge 0\}\cap\{\xxi\,:\,\xi_1^{\leftbottom}\le\xi_1\le \xi_1^{\rightbottom}\}$
 do not have common points except for $\lefttop, \righttop, \rightbottom$, and $\leftbottom$.
    Thus, ${\shock}\cup {\rightsonicop}\cup {\leftsonicop}\cup{\Wedge}$ is a closed curve without self-intersection.
    Denote by $\Om$ the bounded domain enclosed by this closed curve.
\end{itemize}

\smallskip
\item[\rm (ii)] $\vphi$ satisfies the following properties{\rm :}

\smallskip
\begin{itemize}
\item[(ii-1)] $\vphi\in C^{0,1}_{\rm{loc}}(\tdomain)\cap C^1_{\rm{loc}}\bigl(\ol{\tdomain}\setminus \ol{\leftshocksegop\cup\shock\cup \rightshocksegop}\bigr)${\rm ;}

\smallskip
\item[(ii-2)] $\vphi\in C^3(\Omega)\cap C^{2}\bigl(\ol{\Om}\setminus({\ol{\leftsonicop}}\cup{\ol{\rightsonicop}})\bigr)\cap C^{1}(\ol{\Om})${\rm ;}

\smallskip
\item[(ii-3)]
\begin{equation}\label{1-24-op}
\vphi=
\begin{cases}
\ivphi\;\;&\text{in}\;\;\ol{\tdomain}\setminus \ol{\oOmop\cup\Om\cup\nOmop},\\
\leftvphiop\;\;&\text{in}\;\;{\oOmop},\\
\rightvphiop\;\;&\text{in}\;\;{\nOmop},
\end{cases}
\end{equation}
where $\ol{\oOmop}$ shrinks to $\{O\}=\{\lefttop\}=\{\leftbottom\}$
when $\theta_{\rm w}=\theta_{\rm s}^{(\iu)}${\rm ;}

\smallskip
\item[(ii-4)] $\vphi$ satisfies

\smallskip
\begin{itemize}
\item[-] Eq. \eqref{2-1} in $\Om$ with $\rho(|D\vphi|^2, \vphi)$ defined by \eqref{1-o},

\smallskip
\item[-] the slip boundary condition{\rm :} $\der_{\xi_2}\vphi=0$ on $\Wedge\cap\partial\Omega$,

\smallskip
\item[-] the Rankine-Hugoniot conditions{\rm :} $[\vphi]_{\shock}=[\rho(|D\vphi|^2, \vphi)D\vphi\cdot{\bf n}_{\rm sh}]_{\shock}=0$
for the unit normal vector ${\bf n}_{\rm sh}$ to $\shock$ towards the interior of $\Om$.
\end{itemize}
\end{itemize}

\smallskip
\item[\rm (iii)]  Eq. \eqref{2-1} is strictly elliptic in $\ol{\Om}\setminus(\ol{\leftsonicop}\cup\ol{\rightsonicop})${\rm ;} that is,
    \begin{equation*}
    |D\vphi|<c(|D\vphi|^2,\vphi)\qquad \tx{in $\ol{\Om}\setminus(\ol{\leftsonicop}\cup\ol{\rightsonicop})$}.
    \end{equation*}

\smallskip
\item[\rm (iv)]
$\max\{\leftvphiop, \rightvphiop\}\le \vphi \le \ivphi$ in $\Om$.

\smallskip
\item[\rm (v)] Let ${\bm \tau}_{\rm w}=(\cos\theta_{\rm w},\sin\theta_{\rm w})$,
which is tangential to the wedge boundary $\Gamma_{\rm wedge}$.
Let $\leftvecop$ be the unit vector parallel to $\leftshockop$ and oriented
so that $\leftvecop\cdot{\bm \tau}_{\rm w}>0$,
and let $\rightvecop$ be the unit vector parallel to $\rightshockop$ and oriented
so that $\rightvecop\cdot{\bm \tau}_{\rm w}<0${\rm :}
    \begin{equation*}
    \leftvecop=\frac{O\lefttop}{|O\lefttop|}
    =\frac{(\leftvop, \iu-\leftuop)}{\sqrt{(\leftuop-\iu)^2+\leftvop^2}},\qquad
    \rightvecop=-(\cos\theta_{\rm w},\sin\theta_{\rm w}).
    \end{equation*}
    Then
    \begin{equation*}
    \der_{\leftvecop}(\ivphi-\vphi)\le 0, \quad
    \der_{\rightvecop}(\ivphi-\vphi)\le 0 \qquad\,\, \tx{in $\Om$}.
    \end{equation*}
\end{enumerate}

\medskip
\noindent
\textbf{Case {\rm II}.} $\leftuop\le u_{\rm s}^{(\iu)}$, or equivalently,
$\theta_{\rm w}\in [\theta_{\rm s}^{(\iu)},\theta_{\rm d}^{(\iu)})${\rm :}

\smallskip
\begin{enumerate}
\item[\rm (i)] There exists a shock curve $\shock$ with endpoints $O=(0,0)$ and $\righttop=(\xi_1^{\righttop},\xi_2^{\righttop})$
such that the following properties hold{\rm :}

\smallskip
\begin{itemize}
\item[(i-1)] Curve $\shock$ satisfies
\begin{equation}
\label{definition-1a-op}
\shock\subset(\tdomain\setminus\ol{B_{1}(\iu,0)}),
\end{equation}
where $\partial B_{1}(\iu,0)$ is the sonic circle of the state in $\iOm:=\tdomain\setminus \ol{\Om\cup\Om_1}${\rm ;}

\smallskip
\item[(i-2)] Curve $\shock$ is $C^{2}$ in its relative interior. That is, for any $P\in \shock\setminus\{O,\righttop\}$,
 there exist $r>0$, $f\in C^{2}$, and an orthogonal coordinate system $(S,T)$ in $\R^2$
  such that $\shock\cap B_r(P)=\{S=f(T)\}\cap B_r(P)${\rm ;}

\smallskip
\item[(i-3)] Curve $\ol{\shock\cup\rightshocksegop}$ is $C^1$, including at point $\righttop${\rm ;}

\smallskip
\item[(i-4)] $\shock, \rightsonicop$, and $\Wedge:=\{\xi_2=\xi_1\tan\theta_{\rm w}, \xi_2\ge 0\}\cap\{\xxi\,:\,0\le\xi_1\le \xi_1^{\rightbottom}\}$
do not have common points except for $O, \righttop$, and $\rightbottom$.
    Thus, ${\shock}\cup {\rightsonicop}\cup {\Wedge}$ is a closed curve without self-intersection.
    Denote by $\Om$ the bounded domain enclosed by this closed curve.
\end{itemize}

\smallskip
\item[\rm (ii)] $\vphi$ satisfies the following properties{\rm :}

\smallskip
\begin{itemize}
\item[(ii-1)] $\vphi\in C^{0,1}_{\rm{loc}}(\tdomain)\cap C^1_{\rm{loc}}\bigl(\ol{\tdomain}\setminus \ol{\shock\cup \rightshocksegop}\bigr)${\rm ;}

\smallskip
\item[(ii-2)] $\vphi\in C^3(\Omega)\cap C^{2}\bigl(\ol{\Om}\setminus(\{O\}\cup{\ol{\rightsonicop}})\bigr)\cap C^{1}(\ol{\Om})${\rm ;}

\smallskip
\item[(ii-3)] $D\vphi(O)=D\leftvphiop(O)$ and
\begin{equation}\label{1-24a}
\vphi=
\begin{cases}
\ivphi\;\;&\text{in $\ol{\tdomain}\setminus \ol{\Om\cup\nOmop}$},\\
\leftvphiop\;\;&\text{at $O$},\\
\rightvphiop\;\;&\text{in $\nOmop$};
\end{cases}
\end{equation}

\item[(ii-4)] $\vphi$ satisfies

\smallskip
\begin{itemize}
\item[-] Eq. \eqref{2-1} in $\Om$ with $\rho(|D\vphi|^2, \vphi)$ defined by \eqref{1-o},

\smallskip
\item[-] the slip boundary condition{\rm :} $\der_{\xi_2}\vphi=0$ on $\Wedge\cap \partial\Omega$,

\smallskip
\item[-] the Rankine-Hugoniot conditions{\rm :} $[\vphi]_{\shock}=[\rho(|D\vphi|^2, \vphi)D\vphi\cdot{\bf n}_{\rm sh}]_{\shock}=0$
for the unit normal vector ${\bf n}_{\rm sh}$ to $\shock$ towards the interior of $\Om$.
\end{itemize}
\end{itemize}

\smallskip
\item[\rm (iii)] Eq. \eqref{2-1} is strictly elliptic in $\ol{\Om}\setminus (\{O\}\cup \ol{\rightsonicop})$; that is,
    \begin{equation*}
    |D\vphi|<c(|D\vphi|^2,\vphi)\qquad \tx{in $\ol{\Om}\setminus (\{O\}\cup \ol{\rightsonicop})$}.
    \end{equation*}

\smallskip
\item[\rm (iv)] $\max\{\leftvphiop, \rightvphiop\}\le \vphi \le \ivphi$ in $\Om$.

\smallskip
\item[\rm (v)] Let ${\bm \tau}_{\rm w}=(\cos\theta_{\rm w},\sin\theta_{\rm w})$, which is tangential to the wedge boundary $\Wedge$.
Let $\leftvecop$ be the unit vector parallel to $\leftshockop$ and oriented so that $\leftvecop\cdot{\bm \tau}_{\rm w}>0$,
and let $\rightvecop$ be the unit vector parallel to $\rightshockop$ and oriented so that $\rightvecop\cdot{\bm \tau}_{\rm w}<0$. Then
    \begin{equation*}
    \der_{\rightvecop}(\ivphi-\vphi)\le 0, \,\,\,\,
    \der_{\leftvecop}(\ivphi-\vphi)\le 0\qquad\,\,\tx{in $\Om$}.
    \end{equation*}
\end{enumerate}
\end{definition}

Our two main theorems are as follows:

\begin{theorem}
\label{theorem-1-op}
Fix $\gam\ge 1$ and $\iu>1$.
For any $(\iu,\leftuop)\in\mathfrak{P}_{\rm weak}$,
there exists an admissible solution of Problem {\rm \ref{problem-2}}
in the sense of Definition {\rm \ref{def-regular-sol-op}}.
\end{theorem}

\begin{theorem}
\label{theorem-2-op}
Fix $\gam\ge 1$ and $\iu>1$. Given $(\iu,\leftuop)\in\mathfrak{P}_{\rm weak}$,
let $\vphi$ be an admissible solution with the curved shock $\shock$
of Problem {\rm \ref{problem-2}} in the sense
of Definition {\rm \ref{def-regular-sol-op}}. Then the following properties hold{\rm :}

\medskip
\noindent
\textbf{Case {\rm I}.} $\leftuop>u_{\rm s}^{(\iu)}$, or equivalently,
$\theta_{\rm w}\in (0,\theta_{\rm s}^{(\iu)})${\rm :}

\smallskip
\begin{itemize}
\item[(a)]
The curved shock $\shock$ is $C^{\infty}$ in its relative interior, and
   $\vphi\in C^{\infty}(\ol{\Om}\setminus ({\ol{\leftsonicop}}\cup{\ol{\rightsonicop}}))\cap C^{1,1}(\ol{\Om})$.

\smallskip
\item[(b)]  For a constant $\sigma>0$ and a set $\mcl{D}$ given by
\begin{equation*}
\mcl{D}=\big\{ \bm\xi\,:\,\max\{\leftvphiop(\bmxi), \rightvphiop(\bmxi)\}<\ivphi(\bmxi)\big\}\cap \mathfrak{D}_{\theta_{\rm w}},
\end{equation*}
define
\begin{equation}
\label{2-a5-op}
\begin{split}
&\oDop_{\sigma}=\mcl{D}\cap \{\bmxi\,:\,{\rm dist}\{\bmxi,\leftsonicop\}<\sigma\}\cap B_{\leftcop}(\leftPop),\\
&\nDop_{\sigma}=\mcl{D}\cap \{\bmxi\,:\,{\rm dist}\{\bmxi,\rightsonicop\}<\sigma\}\cap B_{\rightcop}(\rightPop),\\
\end{split}
\end{equation}
where $c_j=\rho_j^{(\gam-1)/{2}}$ and
$Q_j=D\varphi_j(O)$, $j=0,1$.
Then, for any $\alp\in(0,1)$ and any $\bmxi_0\in (\ol{\leftsonicop}\cup\ol{\rightsonicop})\setminus\{\lefttop, \righttop\}$,
there exist $\varepsilon_0$ depending on $(\gamma, u_{\infty})$,  and $K<\infty$ depending on $(\iu, \gam, \theta_{\rm w}, \eps_0, \alp)$,
$\|\vphi\|_{C^{1,1}(\Om\cap(\oD_{\eps_0}\cup\nD_{\eps_0}))}$,
and $d={\rm dist}\{\bmxi_0,\shock\}$ such that
    \begin{equation}
    \label{Op1-op}
    \|\vphi\|_{2,\alp,\ol{\Om\cap B_{d/2}(\bmxi_0)\cap(\nDop_{\eps_0/2}\cup\oDop_{\eps_0/2})}}\le K.
    \end{equation}

\item[(c)] For any $\bmxi_0\in \overline{\leftsonicop\cup\rightsonicop}\setminus\{\lefttop, \righttop\}$,
\begin{equation}
\label{Op2-op}
\lim_{\bmxi \to \bmxi_0 \atop \bmxi\in\Om}
\bigl(D_{rr}\vphi-D_{rr}\max\{\rightvphiop,\leftvphiop\}\bigr)(\bmxi)=\frac{1}{\gam+1},
\end{equation}
where $r=|\bmxi-\rightPop|$ near $\rightsonicop$ and $r=|\bmxi -\leftPop|$ near $\leftsonicop$.

\smallskip
\item[(d)] Limits $\displaystyle\lim_{\bmxi\to \lefttop\atop \bmxi\in\Om}D^2\vphi$
and $\displaystyle\lim_{\bmxi \to \righttop\atop \bmxi \in\Om}D^2\vphi$ do not exist.

\smallskip
\item[(e)] $\overline{S_{0, {\rm{seg}}}\cup\shock\cup S_{1,{\rm{seg}}}}$
is a $C^{2,\alp}$--curve for any $\alp\in(0,1)$,
including at points $\lefttop$ and $\righttop$.
\end{itemize}

\medskip
\noindent
\textbf{Case {\rm II}.} $\leftuop\le u_{\rm s}^{(\iu)}$, or equivalently,
$\theta_{\rm w}\in [\theta_{\rm s}^{(\iu)},\theta_{\rm d}^{(\iu)})${\rm :}

\smallskip
\begin{itemize}
\item[(a)]
The curved shock $\shock$ is $C^{\infty}$ in its relative interior,
and
$\vphi\in C^{\infty}(\ol{\Om}\setminus(\{O\}\cup \ol{\rightsonicop}))
\cap C^{1,1}(\ol{\Om}\setminus\{O\})\cap C^{1,\bar{\alp}}(\ol{\Om})$
for some $\bar{\alp}\in(0,1)$
that depends on $\iu$ and $\theta_{\rm w}$ and is non-increasing with respect to $\theta_{\rm w}$.

\smallskip
\item[(b)]  For a constant $\sigma>0$, define $\nDop_{\sigma}$ by \eqref{2-a5-op}.
Then, for any $\alp\in(0,1)$ and any $\bmxi_0 \in \ol{\rightsonicop}\setminus\{\righttop\}$,
there exist $\varepsilon_0$ depending on $(\gamma, u_{\infty})$,  and $K<\infty$ depending on $(\iu, \gam, \theta_{\rm w}, \eps_0, \alp)$,
$\|\vphi\|_{C^{1,1}(\Om\cap\nDop_{\eps_0})}$,
and $d={\rm dist}\{\bmxi_0,\shock\}$ such that
\begin{equation}\label{Op1-a-op}
    \|\vphi\|_{2,\alp,\ol{\Om\cap B_{d/2} (\bmxi_0) \cap\nDop_{\eps_0/2}}}\le K.
\end{equation}

\item[(c)] For any $\bmxi_0 \in \ol{\rightsonicop}\setminus\{\righttop\}$,
\begin{equation}
\label{Op2-a-op}
\lim_{\bmxi \to \bmxi_0 \atop \bmxi \in\Om}
\bigl(D_{rr}\vphi-D_{rr}\rightvphiop\bigr)(\bmxi)=\frac{1}{\gam+1},
\end{equation}
where $r=|\bmxi-\rightPop|$.

\smallskip
\item[(d)] Limit $\displaystyle\lim_{ \bmxi \to \righttop\atop \bmxi \in\Om}D^2\vphi$ does not exist.

\smallskip
\item[(e)] $\overline{\shock\cup\rightshocksegop}$
is a $C^{1,\bar{\alp}}$--curve for the same $\bar{\alp}$ as in statement {\rm (a)}.
Furthermore, curve $\overline{\shock\cup S_{1, {\rm{seg}}}}\setminus \{O\}$ is $C^{2,\alp}$
for any $\alp\in(0, 1)$,
including at point $\righttop$.
\end{itemize}
\end{theorem}

\medskip
\section{Change of the Parameters and Basic Properties}\label{section-change-par}

\subsection{Straight oblique shocks in the self-similar plane}
\label{subsubsec-leftshock}

Given a constant $\iv>0$, denote
\begin{equation}
\label{2-4-b6}
\ivphi:=-\frac 12|\xxi|^2-\iv \xi_2.
\end{equation}

\begin{lemma}\label{lemma-2-0}
For any given $\beta\in[0,\frac{\pi}{2})$, there exists a unique pseudo-potential function
\begin{equation*}
\leftvphi=-\frac 12|\xxi|^2+ (\leftu,\leftv)\cdot\xxi +k_{\mcl{O}}
\end{equation*}
satisfying the following properties{\rm :}

\smallskip
\begin{itemize}
\item[($\mcl{O}_1$)]
$\leftshock:=\{\bmxi\in \R^2\,:\,\ivphi(\bmxi)=\leftvphi(\bmxi)\}$
forms a line of angle $\beta$ with the $\xi_1$--axis, as shown in Fig. {\rm \ref{fig:stobliqueshock}}{\rm ;}

\begin{figure}[htp]
\centering
\begin{psfrags}
\psfrag{b}[cc][][0.8][0]{$\beta$}
\psfrag{lu}[cc][][0.8][0]{$\leftu$}
\psfrag{iv}[cc][][0.8][0]{$\phantom{aa}-\iv$}
\psfrag{x1}[cc][][0.8][0]{$\xi_1$}
\psfrag{x2}[cc][][0.8][0]{$\xi_2$}
\psfrag{o}[cc][][0.8][0]{$O$}
\psfrag{i}[cc][][0.8][0]{$\phantom{aa}\xi_2^{(\beta)}$}
\psfrag{ls}[cc][][0.8][0]{$\leftshock$}
\includegraphics[scale=0.7]{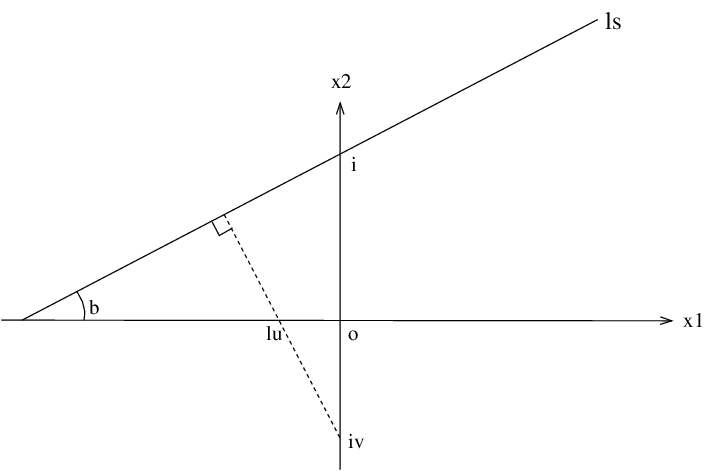}
\caption{$\leftshock$ is a line of angle $\beta$ with the $\xi_1$--axis}
\label{fig:stobliqueshock}
\end{psfrags}
\end{figure}

\smallskip
\item[($\mcl{O}_2$)] $\leftvphi$ satisfies the Rankine-Hugoniot conditions \eqref{1-h} on $\leftshock${\rm :}
    \begin{equation*}
    \leftvphi=\ivphi,\quad \rho(|D\leftvphi|^2,\leftvphi)D\leftvphi\cdot{\bm\nu}_{\rm sh}=D\ivphi\cdot{\bm \nu}_{\rm sh}
    \qquad\,\tx{on $\leftshock$}
    \end{equation*}
    for
    \begin{equation}
    \label{new-density}
    \rho(|D\vphi|^2,\vphi)=\begin{cases}
    \left(1+(\gam-1)(\mcl{B}_{\infty}-\frac 12|D\vphi|^2-\vphi)\right)^{\frac{1}{\gam-1}}&\tx{for}\;\;\gam>1,\\
    \exp\left(\mcl{B}_{\infty}-\frac 12|D\vphi|^2-\vphi\right)&\tx{for}\;\;\gam=1,
    \end{cases}
    \end{equation}
    with
    \begin{equation*}
    \mcl{B}_{\infty}=\frac 12|D\ivphi|^2+\ivphi=\frac{\iv^2}{2},
    \end{equation*}
    where ${\bm\nu}_{\rm sh}:=\frac{D(\varphi_\infty-\varphi_0)}{|D(\varphi_\infty-\varphi_0)|}${\rm ;}

    \smallskip
    \item[($\mcl{O}_3$)] Entropy condition{\rm :}
    \begin{equation*}
    \rho(|D\leftvphi|^2, \leftvphi)>1,\qquad\,\,
    0<D\leftvphi\cdot {\bm\nu}_{\rm sh}<D\ivphi\cdot {\bm\nu}_{\rm sh};
    \end{equation*}

    \item[($\mcl{O}_4$)] $\leftvphi$ satisfies the slip boundary condition on the $\xi_1$--axis{\rm :}
    \begin{equation*}
    \der_{\xi_2}\leftvphi=0\qquad\tx{on $\{\xi_2=0\}$}.
    \end{equation*}
\end{itemize}
\end{lemma}

\begin{proof}
By choosing $(\leftu,\leftv)$ as
\begin{equation}
\label{12-25}
(\leftu,\leftv)=(-\iv\tan\beta, 0),
\end{equation}
$\leftvphi$ satisfies conditions ($\mcl{O}_1$) and ($\mcl{O}_4$).
If  line $\leftshock$ has the  $\xi_2$--intercept at $(0,\xi_2^{(\beta)})$, then $\leftvphi$ can be written as
\begin{equation}
\label{2-4-a0}
\leftvphi=-\frac 12|\xxi|^2-\xi_1\iv \tan \beta-\iv\xi_2^{(\beta)}.
\end{equation}

It remains to find  the  $\xi_2$--intercept $\xi_2^{(\beta)}$ of $\leftshock$ so that $\leftvphi$
satisfies conditions ($\mcl{O}_2$)--($\mcl{O}_3$).
Define
\begin{equation*}
\leftrho:=\rho(|D\leftvphi|^2,\leftvphi).
\end{equation*}
Then $\leftrho$ satisfies
\begin{equation}
\label{2-4-a3}
h(\leftrho)+\frac 12|D\leftvphi|^2+\leftvphi=h(1)+\frac 12|D\ivphi|^2+\ivphi,
\end{equation}
where
$h(\rho)$ is defined by \eqref{1-c}.

In order to determine $\xi_2^{(\beta)}$, we follow the idea from \cite{EL2}.
Define \emph{the pseudo-Mach numbers} $\oM$ and $\iM$ by
\begin{equation}
\label{1-25}
\oM:=\frac{\der_{{\bm \nu}_{\rm sh}}\leftvphi}{\leftc}\qquad\,\tx{for $\leftc=\leftrho^{\frac{\gam-1}{2}}$
and $\iM:=\der_{{\bm \nu}_{\rm sh}}\ivphi$}.
\end{equation}
Since $\der_{{\bm\tau}_{\rm{sh}}}^k(\ivphi-\leftvphi)=0$ on $\leftshock$ for $k=0,1$,
and for a unit tangent vector ${\bm\tau}_{\rm sh}$ of $\leftshock$, it follows from \eqref{2-4-a3} that
\begin{equation}
\label{Bernoulli-RH-new}
h(\leftrho)+\frac 12 (\der_{{\bm \nu}_{\rm sh}}\leftvphi)^2
=\underset{(=0)}{\underbrace{h(1)}}+\frac 12(\der_{{\bm \nu}_{\rm sh}}\ivphi)^2\qquad\tx{on $\leftshock$.}
\end{equation}
By \eqref{1-25},  $\displaystyle{\rho(|D\leftvphi|^2,\leftvphi)D\leftvphi\cdot{\bm\nu}_{\rm sh}=D\ivphi\cdot{\bm \nu}_{\rm sh}}$
can be rewritten as
\begin{equation}
\label{leftrho-expression}
\leftrho^{\frac{\gam+1}{2}}=\frac{\iM}{\oM}.
\end{equation}
We substitute this expression into \eqref{Bernoulli-RH-new} to obtain
\begin{equation}
\label{1-1}
\begin{split}
\big(1+\frac{\gam-1}{2}\oM^2\big)\oM^{-\frac{2(\gam-1)}{\gam+1}}
=\big(1+\frac{\gam-1}{2}\iM^2\big)\iM^{-\frac{2(\gam-1)}{\gam+1}}.
\end{split}
\end{equation}

Notice that $f(M):=(1+\frac{\gam-1}{2}M^2)M^{-\frac{2(\gam-1)}{\gam+1}}$ satisfies
\begin{equation*}
\lim_{M\to 0+}f(M)=\infty,\quad \lim_{M\to \infty}f(M)=\infty,\quad
f'(M)=\frac{2(\gam-1)}{\gam+1}M^{-\frac{2(\gam-1)}{\gam+1}-1}(M^2-1).
\end{equation*}
Therefore, if $\iM=1$, then $\oM=1$ is the only solution of \eqref{1-1}.
If $\iM\in(0, \infty)\setminus \{1\}$, then \eqref{1-1} has a unique nontrivial
solution $\oM$ in $(0, \infty)\setminus \{1\}$ with $M_{\mcl{O}}\neq M_{\infty}$.
Furthermore, a direct computation from \eqref{1-1} shows that
\begin{equation}
\label{D}
\frac{\dd \oM}{\dd \iM}<0\qquad\quad\tx{for all $\iM\in(0,\infty)\setminus\{1\}$}.
\end{equation}

It follows from \eqref{leftrho-expression} that conditions ($\mcl{O}_2$)--($\mcl{O}_3$) are satisfied
if there exists $\xi_2^{(\beta)}$ so that $\iM>1$ holds.

Denote $\iq:=\iM$ and $\leftq:=\leftc \oM$.
Note that
$\leftq={\rm{dist}}(\leftshock, (\leftu,0))$ and $\iq={\rm{dist}}(\leftshock, (0,-\iv))$
for $\leftu$ given by \eqref{12-25}.
Then
\begin{equation}
\label{2-4-a4}
\iq-\leftq=\iv\sec\beta.
\end{equation}
We substitute the representations of $\iq=\iM$
and $\leftq=\oM\leftc=\oM\big(\frac{\iM}{\oM}\big)^{\frac{\gam-1}{\gam+1}}$ into \eqref{2-4-a4}
to obtain
\begin{equation}\label{2-4-a5}
\iM^{\frac{\gam-1}{\gam+1}}\big(\iM^{\frac{2}{\gam+1}}-\oM^{\frac{2}{\gam+1}}\big)
=\iv \sec\beta,
\end{equation}
where $\oM\le 1$ solves \eqref{1-1} for $\iM\ge 1$.
As a function of $\iM\ge \oM$, the left-hand side of \eqref{2-4-a5} as a function of $\iM$ has the derivative that is greater than  $\frac 2{\gamma+1}$
for $\iM\ge \oM$,
and its value at $\iM=\oM$ is $0$.
%
%
Therefore, for given constants $\iv>0$ and $\beta\in [0, \frac{\pi}{2})$, there exists a unique
\begin{equation}
\label{iM-supersonic}
\iM>1
\end{equation}
satisfying equation \eqref{2-4-a5}.
Once $\iM>1$ is decided, it follows from \eqref{12-25} and \eqref{2-4-a4} that
\begin{equation}
\label{2-4-a6}
\xi_2^{(\beta)}=\iM\sec\beta-\iv.
\end{equation}
It can be seen from $0<D\leftvphi\cdot {\bm\nu}_{\rm sh}<D\ivphi\cdot {\bm\nu}_{\rm sh}$ that
the $\xi_2$--intercept $\xi_2^{(\beta)}$ given by \eqref{2-4-a6} satisfies
\begin{equation*}
\xi_2^{(\beta)}>0.
\end{equation*}

Case $\gam=1$ can be proved similarly.
\end{proof}

\subsection{New parameters $(\iv, \beta)$}

We define $\bmxi'=(\xi_1',\xi_2')$ by
\begin{equation}
\label{2-a}
\begin{pmatrix}
\xi_1' \\
\xi_2' \end{pmatrix}
:=\begin{pmatrix}
\cos\theta_{\rm w}&\sin\theta_{\rm w}\\
-\sin\theta_{\rm w}&\cos\theta_{\rm w}
\end{pmatrix}\begin{pmatrix}
\xi_1\\ \xi_2
\end{pmatrix}
-\begin{pmatrix}
\iu\cos\theta_{\rm w}\\
0
\end{pmatrix}.
\end{equation}
In the new coordinates $(\xi_1', \xi_2')$, center $\rightPop$ of the sonic circle $\der B_{\rightcop}(\rightPop)$
becomes the origin, and $\Wedge$ lies on the horizontal axis $\xi_2'=0$.

Hereafter, for simplicity of notation, we denote $\bmxi=(\xi_1,\xi_2)$ as
the new coordinates $(\xi_1', \xi_2')$ given by \eqref{2-a}.
In the new coordinate system,
$\ivphi$, $\leftvphiop$, and $\rightvphiop$, defined by \eqref{1-m}, \eqref{newlabel-N0},
and \eqref{1-rightvphi}, are expressed respectively as
\begin{equation}
\label{1-21}
\begin{split}
&\ivphi^{\rm op}(\xxi)=-\frac 12|\xxi|^2-\xi_2 \iu\sin\theta_{\rm w}+\frac 12 \iu^2\cos^2\theta_{\rm w},\\
&\leftvphi^{\rm op}(\xxi)=
-\frac 12|\xxi|^2+(\xi_1+\iu\cos \theta_{\rm w})(\leftuop\sec\theta_{\rm w}-\iu\cos \theta_{\rm w})
+\frac 12 \iu^2\cos^2\theta_{\rm w},\\
&\rightvphi^{\rm op}(\xxi)=
-\frac 12|\xxi|^2-\iu\netaop\sin\theta_{\rm w}+\frac 12 \iu^2\cos^2\theta_{\rm w}.
\end{split}
\end{equation}
We define $(\ivphi, \leftvphi, \rightvphi)$ in the new coordinates by
\begin{equation}
\label{1-21-addition}
\begin{split}
&\ivphi(\xxi)=\ivphi^{\rm op}(\xxi)-\frac 12 \iu^2\cos^2\theta_{\rm w},\\
&\leftvphi(\xxi)=\leftvphi^{\rm op}(\xxi)-\frac 12 \iu^2\cos^2\theta_{\rm w},\\
&\rightvphi(\xxi)=\rightvphi^{\rm op}(\xxi)-\frac 12 \iu^2\cos^2\theta_{\rm w}.
\end{split}
\end{equation}
In the new coordinate system, $\leftshockop, \rightshockop, \leftsonicop$, and $\rightsonicop$
are denoted as $\leftshock, \rightshock, \leftsonic$, and $\rightsonic$, respectively.

\begin{definition}[New parameters $(\iv, \beta)$]
\label{definition-new-parameters}
For each $(\iu, u_0)\in \mathfrak{P}$, we introduce new parameters $(\iv, \beta)\in (0, \infty)\times (0, \frac{\pi}{2})$
as follows{\rm :}

\smallskip
\begin{itemize}
\item[(i)] For $\theta_{\rm w}\in (0, \theta_{\rm d}^{(\iu)})$ given by \eqref{def-v-2015}, define $\iv$ by
\begin{equation*}
\iv=\iu\sin \theta_{\rm w};
\end{equation*}

\item[(ii)] Let $\leftshockop$ be the straight oblique shock corresponding to
point $u_0(1,\tan\theta_{\rm w})$ on the shock polar {\rm (}Fig. {\rm \ref{fig:polar}}{\rm )}
with the incoming state $(\iu, 0)$.
For such $\leftshockop$, let $\theta_{\leftshockop}$ be the angle of $S_0$ from the horizontal ground
{\rm (}{\it i.e.}, $\xi_2=0$ in the coordinates $\xxi$ before \eqref{2-a}{\rm )}.
Define $\beta\in(0,\frac{\pi}{2})$ by
\begin{equation}
\begin{split}
\label{def-beta-2015}
&\beta:=\theta_{\leftshockop}-\theta_{\rm w}.
\end{split}
\end{equation}
\end{itemize}
Note that the definition of $\iv$ stated in {\rm (i)} coincides with \eqref{1-22}.
\end{definition}

The weak shock configuration in the new self-similar plane
is shown in Figs. \ref{fig:globaln1}--\ref{fig:globaln2} for $(\iv,\beta)\in (0, \infty)\times (0, \frac{\pi}{2})$.
\begin{figure}[htp]
\centering
\begin{psfrags}
\psfrag{i}[cc][][0.8][0]{$\phantom{aa}(0,-\iv)$}
\psfrag{b}[cc][][0.8][0]{$\beta$}
\psfrag{tw}[cc][][0.8][0]{$\theta_{\rm w}$}
\psfrag{o}[cc][][0.8][0]{$O$}
\psfrag{iv}[cc][][0.8][0]{$(0,-\iv)$}
\psfrag{x}[cc][][0.8][0]{$\xi_1$}
\psfrag{y}[cc][][0.8][0]{$\xi_2$}
\psfrag{s}[cc][][0.8][0]{$\shock$}
\psfrag{ls}[cc][][0.8][0]{$\leftshock$} \psfrag{rs}[cc][][0.8][0]{$\rightshock$}
\psfrag{Tip}[cc][][0.8][0]{$\xi^*(\beta)$}
\psfrag{snr}[cc][][0.8][0]{$\phantom{aa}\rightsonic$}
\psfrag{snl}[cc][][0.8][0]{$\phantom{aaa}\leftsonic$}
\includegraphics[scale=.6]{newplane2016a.eps}
\caption{Weak shock solutions in the new self-similar plane when $\theta_{\rm w}<\theta_{\rm s}^{(\iu)}$}
\label{fig:globaln1}
\includegraphics[scale=.6]{newplane2016b.eps}
\caption{Weak shock solutions in the new self-similar plane when $\theta_{\rm s}^{(\iu)}\le \theta_{\rm w}< \theta_{\rm d}^{(\iu)}$}
\label{fig:globaln2}
\end{psfrags}
\end{figure}

We define a parameter set $\mathfrak{R}$ by
\begin{equation}
\label{def-R-set}
\mathfrak{R}:=\{(\iv,\beta)\,:\,\iv>0,\,\,0<\beta<\frac{\pi}{2}\},
\end{equation}
and define a map $\mcl{T}: \mathfrak{P}\rightarrow \mathfrak{R}$ by
\begin{equation}
\label{mapping-T}
\mcl{T}(\iu, u_0)=(\iv,\beta)\qquad \,\tx{for $(\iv, \beta)$ given by Definition \ref{definition-new-parameters}}.
\end{equation}

\begin{lemma}
\label{lemma-parameters}
For any given $\gam\ge 1$, map $\mcl{T}: \mathfrak{P}\rightarrow \mathfrak{R}$ is a homeomorphism.
\end{lemma}

\begin{proof}
Fix $(\iu,u_0)\in \mathfrak{P}$. By Definition \ref{definition-new-parameters}(i),
the corresponding half-wedge angle $\theta_{\rm w}$ is given by
\begin{equation}
\label{expression-tw-op}
\theta_{\rm w}=\arctan (\frac{f_{\rm{polar}}(u_0)}{u_0}),
\end{equation}
where $f_{\rm{polar}}$ is the function introduced in Lemma \ref{lemma-app2}.

By Definition \ref{definition-new-parameters}(ii),
a unit tangent vector ${\bm\tau}_{\leftshockop}$ of the straight oblique shock $\leftshockop$ corresponding
to $(\iu, u_0)$ is ${\bm\tau}_{\leftshockop}=(\cos \theta_{\leftshockop}, \sin \theta_{\leftshockop})$
in the coordinate system introduced right before transformation \eqref{2-a}.
Substituting this expression of ${\bm\tau}_{\leftshockop}$ into one of
the Rankine-Hugoniot conditions:
$$
(\iu,0)\cdot{\bm\tau}_{\leftshockop}=(u_0, f_{\rm{polar}}(u_0))\cdot{\bm\tau}_{\leftshockop},
$$
we have
\begin{equation}
\label{expression-theta-lsop}
\tan\theta_{\leftshockop}=\frac{\iu-\leftuop}{f_{\rm{polar}}(\leftuop)}.
\end{equation}
From \eqref{def-beta-2015} and \eqref{expression-tw-op}--\eqref{expression-theta-lsop}, we obtain
\begin{equation*}
\tan\beta
=\frac{\tan \theta_{\leftshockop}-\tan\theta_{\rm w}}{1+\tan\theta_{\leftshockop}\tan\theta_{\rm w}}=
\frac{\leftuop(\iu-\leftuop)-\bigl(f_{\rm{polar}}(u_0)\bigr)^2}
{\leftuop f_{\rm{polar}}(u_0)}>0.
\end{equation*}
By Definition \ref{definition-new-parameters}(i) and \eqref{expression-tw-op}, we can express $\iv$ as
\begin{equation*}
\iv=\iu \sin(\arctan (\frac{f_{\rm{polar}}(u_0)}{u_0})).
\end{equation*}
Therefore,  map $\mcl{T}:\mathfrak{P}\rightarrow \mathfrak{R}$ is continuous.

In order to show that $\mcl{T}:\mathfrak{P}\rightarrow \mathfrak{R}$ is invertible
and its inverse is continuous,
for fixed $(\iv,\beta)\in \mathfrak{R}$,
we find  a solution $(\iu,\leftuop)\in \mathfrak{P}$ of the following equations:
\begin{align}
\label{2-4-b2}
&\iu \sin\theta_{\rm w}=\iv,\\
\label{2-4-b3}
&\iu\cos\theta_{\rm w}=\xi_2^{(\beta)}\cot\beta,\\
\label{2-4-b4}
&\leftuop\sec\theta_{\rm w}=\xi_2^{(\beta)}\cot\beta-\iv\tan\beta,
\end{align}
so that the definitions of $\leftvphi$ in  \eqref{2-4-a0} and \eqref{1-21-addition} coincide.
Combining \eqref{2-4-b2} with \eqref{2-4-b3}, we have
\begin{equation}
\label{2-4-b5}
\iu=\sqrt{\iv^2+(\xi_2^{(\beta)})^2\cot^2\beta}=:T_1(\iv,\beta).
\end{equation}
Using \eqref{2-4-b6}, we can rewrite \eqref{2-4-b5} as
\begin{equation*}
\iu=|D\ivphi(-\xi_2^{(\beta)}\cot\beta,0)|.
\end{equation*}
Then we obtain from \eqref{iM-supersonic} that
$\iu \ge \iM>1$.

Once $\iu$ is given by \eqref{2-4-b5}, we combine it with \eqref{2-4-b3}--\eqref{2-4-b4} to obtain $u_0$ as
\begin{equation}
\label{2-4-c1}
\leftuop=
\frac{\big(\xi_2^{(\beta)}\cot\beta-\iv\tan\beta\big)\xi_2^{(\beta)}\cot\beta}
{T_1(\iv, \beta)}=:T_2(\iv,\beta).
\end{equation}
Note that $(-\xi_2^{(\beta)}\cot \beta, 0)$
is the $\xi_1$--intercept of line $\leftshock$  from Lemma \ref{lemma-2-0}.
Therefore, it can be seen from Fig. \ref{fig:stobliqueshock}
that
$\xi_2^{(\beta)}\cot\beta+u_{\mathcal O}=\xi_2^{(\beta)}\cot\beta-v_{\infty}\tan\beta>0$.
This implies that $\leftuop>0$.
Since $\tan \theta_{\rm w}=\frac{\iv}{\xi_2^{(\beta)}\cot\beta}>0$
is obtained from \eqref{2-4-b2}--\eqref{2-4-b3},
we conclude that $(\iu, \leftuop)$ given by \eqref{2-4-b5}--\eqref{2-4-c1}
is contained in $\mathfrak{P}$.

Finally, the continuity of $\mcl{T}^{-1}$ follows directly from the definitions of $(T_1,T_2)$.
\end{proof}

For any given $(\iv, \beta)\in\mathfrak{R}$,
the $\xi_2$--intercept $\xi_2^{(\beta)}>0$ of the oblique shock $\leftshock$ of angle $\beta$
from the $\xi_1$--axis is uniquely defined.
Moreover, $\xi_2^{(\beta)}$ varies continuously on $\beta\in(0,\frac{\pi}{2})$,
and $\underset{\beta\to 0+}{\lim} \xi_2^{(\beta)}$ exists
and is positive.
Denote $\neta:=\xi_2^{(\beta)}|_{\beta=0}$.
Let $\rightvphi$ denote $\leftvphi$ corresponding to $\beta=0$.
Then $\rightvphi$ is given by
\begin{equation}
\label{def-rightvphi}
\rightvphi(\xxi)=-\frac 12|\xxi|^2-\iv \neta.
\end{equation}

\begin{remark}[The normal shock: Case $\beta=0$]
\label{remark-normalshock}
For fixed $\gam\ge 1$ and $\iv>0$, the straight shock of angle $\beta=0$ from the horizontal ground
{\rm (}{\it i.e.}, $\xi_2=0$ in the new coordinates $\xxi$ given by \eqref{2-a}{\rm )}
can be considered by taking limit
$\beta\to 0+$ in the argument above.
The state of $\beta=0$ is that of a normal shock, which corresponds to
the state of $\frac{\leftuop}{\iu}=1$ with $\theta_{\rm w}=0$.
Even though the case of $\beta=0$ is not physical because $\iu=\infty$,
we still put this case under our consideration
as it is useful in applying the Leray-Schauder degree argument to
prove the existence of admissible solutions of Problem {\rm \ref{problem-2}}
for all $(\iu, u_0)\in \mathfrak{P}_{\rm weak}$.
\end{remark}

\begin{remark}\label{remark-pt}
According to Lemma {\rm \ref{lemma2-appendix}},
for each $\iv>0$, there exists $\beta_{\rm d}^{(\iv)}\in(0,\frac{\pi}{2})$ such that,
if the parameter sets $\mathfrak{R}_{\rm weak}, \mathfrak{R}_{\rm detach}$, and $\mathfrak{R}_{\rm strong}$
are defined by
\begin{equation}\label{2-4-a8}
\begin{split}
&\mathfrak{R}_{\rm weak}=\underset{\iv>0}{\cup}\{\iv\}\times (0,\beta_{\rm d}^{(\iv)}),\\
&\mathfrak{R}_{\rm detach}=\underset{\iv>0}{\cup}\{\iv\}\times \{\beta_{\rm d}^{(\iv)}\},\\
&\mathfrak{R}_{\rm strong}=\underset{\iv>0}{\cup}\{\iv\}\times(\beta_{\rm d}^{(\iv)},\frac{\pi}{2}),
\end{split}
\end{equation}
then
\begin{equation}
\label{parameter-map}
\mcl{T}^{-1}(\mathfrak{R}_{\rm weak})=\mathfrak{P}_{\rm weak},\quad
\mcl{T}^{-1}(\mathfrak{R}_{\rm detach})=\mathfrak{P}_{\rm detach},\quad
\mcl{T}^{-1}(\mathfrak{R}_{\rm strong})=\mathfrak{P}_{\rm strong},
\end{equation}
for $\mathfrak{P}_{\rm weak}, \mathfrak{P}_{\rm detach}$,
and $\mathfrak{P}_{\rm strong}$ defined by \eqref{2-4-a7}.
In Lemma {\rm \ref{lemma:interval-existence}}, we will also show that, for any $\iv>0$,
there exists a unique $\beta_{\rm s}^{(\iv)}\in(0,\beta_{\rm d}^{(\iv)})$
such that $T_2(\iv,\beta)>u_{\rm s}^{(\iu)}$ if and only if $\beta<\betac^{(\iv)}$
for $\iu=T_1(\iv,\beta)$,
where $u_{\rm s}^{(u_{\infty})}$ denotes the value of $u_{\rm wk}^{(\theta_{\rm w})}$ for $\theta_{\rm w}=\theta_{\rm s}^{(\iu)}$.
\end{remark}

\begin{figure}[htp]
\centering
\begin{psfrags}
\psfrag{lt}[cc][][0.8][0]{$\lefttop$}
\psfrag{lq}[cc][][0.8][0]{$\leftq$}
\psfrag{lc}[cc][][0.8][0]{$\leftc$}
\psfrag{lu}[cc][][0.8][0]{$(\leftu,0)$}
\psfrag{ls}[cc][][0.8][0]{$\leftshock$}
\psfrag{m}[cc][][0.8][0]{$\phantom{aa}{\bm\xi}^m$}
\includegraphics[scale=1.0]{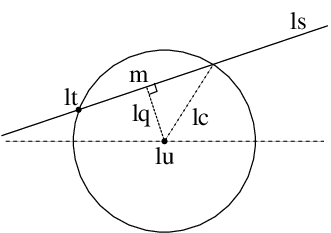}
\caption{Two intersection points of $\leftshock$ with the sonic circle $\der B_{\leftc}(\leftu,0)$}
\label{fig:global}
\end{psfrags}
\end{figure}

\smallskip
For fixed $(\iv,\beta)\in\mathfrak{R}_{\rm weak}$,  let $\oM$ be defined by \eqref{1-25}.
In the proof of Lemma \ref{lemma-2-0}, it is shown that $0<\oM<1$.
This implies that the corresponding straight oblique shock $\leftshock$ intersects
with {\emph{the sonic circle}}
$\partial B_{\leftc}(\leftu,0)=\{\bmxi\,:\,|D\leftvphi(\bmxi)|=\leftc\}$ at two distinct points.
For each $\beta\in[0,\frac{\pi}{2})$, let $\xxi^{\mathcal{O}}:=(\oxi,\oeta)$
be the intersection point $\lefttop$
with the smaller $\xi_1$--coordinate (see Fig. \ref{fig:global}).
Moreover, let $(\xi_1^{(\beta)},0)$ be the $\xi_1$--intercept of $\leftshock$.
If $\oeta>0$, then $|D\leftvphi|>\leftc$ at $(\xi_1^{(\beta)},0)$, which means that an admissible
solution in the sense of Definition \ref{def-regular-sol-op} for $(\iu,u_0)=\mcl{T}^{-1}(\iv,\beta)$
has the structure shown in Fig. \ref{fig:globaln1}.
On the other hand, if $\oeta\le 0$, then an admissible solution
for $(\iu,u_0)=\mcl{T}^{-1}(\iv,\beta)$ has the structure shown in Fig. \ref{fig:globaln2}.

\begin{lemma}
\label{lemma:interval-existence} Fix $\gam\ge 1$ and $\iv>0$.
The $\xi_2$--coordinate $\oeta$ of point $\lefttop$  satisfies
\begin{equation*}
\frac{\dd\oeta}{\dd\beta}<0\;\;\qquad\text{for all}\;\;\beta\in (0,\frac{\pi}{2})\,\,\tx{and}\,\,
\lim_{\beta\to \frac{\pi}{2}-}\oeta=-\infty.
\end{equation*}
Therefore, there exists $\betac^{(\iv)}\in(0,\frac{\pi}{2})$ such that
\begin{equation}
\label{2-4-c2}
\begin{cases}
\oeta>0\,\, \Longleftrightarrow \,\,\frac{|D\leftvphi(\xi_1^{(\beta)},0)|}{\leftc}>1&\quad\tx{for $\beta\in[0,\betac^{(\iv)})$},\\
\oeta = 0 \,\,\Longleftrightarrow\,\, \frac{|D\leftvphi(\xi_1^{(\beta)},0)|}{\leftc} = 1&\quad\tx{for $\beta=\betac^{(\iv)}$},\\
\oeta< 0 \,\,\Longleftrightarrow\,\, \frac{|D\leftvphi(\xi_1^{(\beta)},0)|}{\leftc}< 1&\quad\tx{for $\beta\in (\betac^{(\iv)},\frac{\pi}{2})$}.
\end{cases}
 \end{equation}
In addition, $\betac^{(\iv)}$ satisfies the inequality{\rm :}
\begin{equation}
\label{2-4-c3}
\betac^{(\iv)}<\betadet.
\end{equation}

\begin{proof}
For $M_{\infty}$ and $M_{\mcl{O}}$ given by \eqref{1-25}, define
\begin{equation}
\label{definition-qs}
(q_{\infty}, q_{\mcl{O}})=(M_{\infty}, M_{\mcl{O}}\leftc).
\end{equation}
For each $\beta\in(0, \frac{\pi}{2})$, let $\xxi^m=(\xi_1^m, \xi_2^m)$ be the midpoint of two intersections
of $\leftshock$ with $\der B_{\leftc}(\leftu,0)$. By \eqref{1-25}, we have
\begin{equation}
\label{1-52}
\oeta=\xi_2^m-\leftc\sqrt{1-\oM^2}\sin\beta.
\end{equation}
Since  ${(\xi_1^m-\leftu,\xi_2^m)}$ is perpendicular to $\leftshock$,
\begin{equation*}
\der_{{\bm\tau}_{\rm sh}}\leftvphi(\xxi^m)=0
=\der_{{\bm\tau}_{\rm sh}}\ivphi(\xxi^m)
=(-\xi_1^m,-\xi_2^m-\iv)\cdot{\bm\tau}_{\rm sh}
\end{equation*}
for a unit tangent vector ${\bm\tau}_{\rm sh}=(\cos\beta,\sin\beta)$ to $\leftshock$. Then we have
$$
\xxi^m
=(0,-\iv)-\iq{\bm \nu}_{\rm sh}=(0,-\iv)-\iq(\sin\beta,-\cos\beta)
$$
for the unit normal vector ${\bm\nu}_{\rm sh}$ to $\leftshock$ pointing towards the $\xi_1$--axis. This yields that
\begin{equation}
\label{1-36}
\xi_2^m=-\iv+\iq\cos\beta.
\end{equation}
We differentiate \eqref{2-4-a4} and \eqref{1-36} with respect to $\beta$ to obtain
\begin{equation}
\label{diff-qi-beta}
\frac{\dd\xi_2^m}{\dd\beta}=-\iq\sin\beta+\frac{\dd\iq}{\dd\beta}\cos\beta,\,\,\quad
\,\,\quad
\frac{\dd\iq}{\dd \beta}=\frac{\iq-\leftq}{1-\frac{\dd\leftq}{\dd \iq}}\tan\beta,
\end{equation}
and combine the results to obtain
\begin{align}\label{1-45}
\frac{\dd\xi_2^m}{\dd\beta}
=-\frac{1-\frac{\iq}{\leftq}\frac{\dd \leftq}{\dd \iq}}{1-\frac{\dd \leftq}{\dd \iq}}\xi_2^m\tan\beta.
\end{align}
If $\frac{\dd \leftq}{\dd \iq}\le 0$, then
\begin{equation*}
\frac{1-\frac{\iq}{\leftq}\frac{\dd \leftq}{\dd\iq}}{1-\frac{\dd \leftq}{\dd \iq}}>1\ge \frac{2}{\gam+1}.
\end{equation*}
A direct computation by using \eqref{1-25}--\eqref{1-1} shows that
\begin{equation}
\label{diff-lq-qi}
\begin{split}
\frac{\dd \leftq}{\dd\iq}
&=\Big(\frac{\oM}{\iM}\Big)^{\frac{2}{\gam+1}}
\Big(\frac{\gam-1}{\gam+1}+\frac{2}{\gam+1}\frac{\iM}{{\oM}}\frac{\dd {\oM}}{\dd\iM}\Big)\\
&={{\Big(\frac{\iq}{\leftrho \leftq}\Big)^{\frac{1}{\gam+1}}}}\frac{\leftq}{\iq}\Big(\frac{\gam-1}{\gam+1}
+\frac{2{\iM}}{(\gam+1){\oM}}{{\frac{\dd{\oM}}{\dd{\iM}}}}\Big)\\
&\le \frac{\gam-1}{\gam+1}\frac{\leftq}{\iq}.
\end{split}
\end{equation}
If $\frac{\dd \leftq}{\dd\iq}>0$, it follows from $0<1-\frac{\dd \leftq}{\dd \iq}<1$ that
\begin{equation*}
\frac{1-\frac{\iq}{\leftq}\frac{\dd \leftq}{\dd \iq}}{1-\frac{\dd \leftq}{\dd \iq}}
>1-\frac{\iq}{\leftq}\frac{\dd \leftq}{\dd \iq}\ge \frac{2}{\gam+1}.
\end{equation*}
We apply inequality $\displaystyle{\frac{1-\frac{\iq}{\leftq}\frac{\dd \leftq}{\dd\iq}}{1-\frac{\dd\leftq}{\dd \iq}}>\frac{2}{\gam+1}}$ to derive from \eqref{1-45} that
\begin{equation}
\label{xi2m-new1}
\frac{\dd \xi_2^m}{\dd\beta}\le -\frac{2}{\gam+1}\xi_2^m\tan\beta
\qquad\tx{for all $\beta\in(0,\frac{\pi}{2})$}.
\end{equation}

Next, we differentiate
$\leftc^2=1+\frac{\gam-1}{2}(\iq^2-q^2)$ with respect to $\beta$ and use \eqref{2-4-a4} to obtain
\begin{equation}
\label{density-mont-ox}
\begin{split}
\frac{\dd\leftc^2}{\dd \beta}
&=(\gam-1)\iq \Big(1-\frac{\leftq}{\iq}\frac{\dd \leftq}{\dd \iq}\Big)\frac{\dd \iq}{\dd\beta}\\
&\ge \frac{2(\gam-1)}{\gam+1} \iv\sec\beta\tan\beta
\qquad\,\,\text{for all $\beta\in(0,\frac{\pi}{2})$}.
\end{split}
\end{equation}
From this, we have
\begin{equation}
\label{xi2m-new2}
\lim_{\beta\to \frac{\pi}{2}-}\xi_2^m=0,\quad\, \lim_{\beta\to \frac{\pi}{2}-}\leftc=\infty,\quad \,\lim_{\beta\to \frac{\pi}{2}-}\oeta=-\infty.
\end{equation}
Notice that
\begin{equation}
\label{iq-monotonicity}
\frac{\dd\iq}{\dd\beta}>0,
\end{equation}
which can be obtained from differentiating \eqref{2-4-a5} with respect to $\beta$,
where $0<\oM<1<\iM$ is used.
From \eqref{D}, we obtain
\begin{equation}
\label{oM-monotonicity}
\frac{\dd \oM}{\dd\beta}=\frac{\dd \oM}{\dd \iM}\frac{\dd \iM}{\dd\beta}=\frac{\dd \oM}{\dd \iM}\frac{\dd \iq}{\dd\beta}<0.
\end{equation}

Therefore, we conclude from \eqref{1-52} and the monotonicity properties of $(\xi_2^m, \leftc^2, M_{\mcl{O}})$ with respect
to $\beta$ that $\frac{\der \xi_2^{\mcl{O}}}{\der\beta}<0$ for all $\beta\in(0, \frac{\pi}{2})$.
\end{proof}
\end{lemma}

\section{Main Theorems in the $(\iv,\beta)$--Parameters}
\label{subsec-mainthm}

With Lemma \ref{lemma-parameters} and Remark \ref{remark-pt},
we can restate Theorems \ref{theorem-1-op}--\ref{theorem-2-op}
by using parameters $(\iv, \beta)\in \mathfrak{R}_{\rm{weak}}$.

For fixed $\gam\ge 1$ and $(\iv,\beta)\in\mathfrak{R}$, we recall the definitions of $(\ivphi, \leftvphi, \rightvphi)$
given by \eqref{2-4-b6}, \eqref{2-4-a0}, and \eqref{def-rightvphi} as follows:
\begin{equation}\label{def-uniform-ptnl-new}
  \ivphi=-\frac 12|{\bm\xi}|^2-\iv \etan,\quad
  \leftvphi=-\frac 12|{\bm\xi}|^2+\leftu \xi_1-\iv \xi_2^{(\beta)},\quad
  \rightvphi=-\frac 12|\xxi|^2-\iv\neta,
\end{equation}
for $\xi_2^{(\beta)}$ given by \eqref{2-4-a6}.

Let
\begin{equation*}
\leftrho=\rho(|D\leftvphi|^2,\leftvphi),\qquad \rightrho=\rho(|D\rightvphi|^2,\rightvphi)
\end{equation*}
for $\rho(|D\vphi|^2,\vphi)$ defined by \eqref{new-density}.
Note that $\neta$ satisfies that $\neta<\rightc$ for $\rightc=\rightrho^{\frac{\gam-1}{2}}$.
Define
\begin{equation*}
O_{\mathcal O}:=(\leftu,0),\qquad O_{\mathcal N}=(0,0).
\end{equation*}
Since $\neta<\rightc$, $\der B_{\rightc}(O_{\mathcal N})$ intersects with $\rightshock=\{\xi_2=\neta\}$
at two distinct points.
For each $\beta\in[0, \frac{\pi}{2})$, $\xi_2=f_{\mcl{O}}(\xi_1)$, obtained by solving the equation
$\ivphi(\xi_1,\xi_2)-\leftvphi(\xi_1,\xi_2)=0$ for $\xi_2$, is given by
\begin{equation}
\label{fo}
f_{\mcl{O}}(\xi_1):=\xi_1\tan\beta+\xi_2^{(\beta)}.
\end{equation}
Note that $\leftshock=\{\xi_2=f_{\mcl{O}}(\xi_1)\}$ intersects with $\der B_{\leftc}(O_{\mathcal O})$ at two distinct points.
The $\xi_1$--intercept of $\leftshock$ is
\begin{equation}
\label{def-Pbeta}
P_{\beta}=(-\xi_2^{(\beta)}\cot\beta,0)=:(\xi_1^{(\beta)},0).
\end{equation}
The line passing through $P_{\beta}$ and $\Oi=(0, -\iv)$ is given by
\begin{equation}
\label{definition-f-w}
L_{\rm w}:=\{\xxi\,:\,  \, \xi_2= f_{\rm w}(\xi_1):=\tan \theta_{\infty}(\xi_1-\xi_1^{(\beta)})\}
\end{equation}
for
\begin{equation*}
\tan\theta_{\infty}=\frac{\iv}{\xi_1^{(\beta)}}
\qquad\tx{with $\theta_\infty\in (\frac{\pi}{2},\pi)$}.
\end{equation*}
Then $L_{\rm w}$ represents the horizontal ground in the self-similar plane
before the linear transformation \eqref{2-a} of the self-similar variables $(\xi_1, \xi_2)$.
Moreover, $\tan\theta_{\infty}$  and $L_{\rm w}$ depend continuously on $(\iv, \beta)$.

\begin{figure}[htp]
\centering
\begin{psfrags}
\psfrag{i}[cc][][0.8][0]{$-\iv\;\;\;$}
\psfrag{b}[cc][][0.8][0]{$\beta$}
\psfrag{ns}[cc][][0.8][0]{$\rightshock$}
\psfrag{o}[cc][][0.8][0]{$O$}
\psfrag{ls}[cc][][0.8][0]{$\leftshock$}
\psfrag{s}[cc][][0.8][0]{$\Shock$}
\psfrag{lt}[cc][][0.8][0]{$\lefttop$}
\psfrag{rt}[cc][][0.8][0]{$\righttop$}
\psfrag{lso}[cc][][0.8][0]{$\;\;\leftsonic$}
\psfrag{rso}[cc][][0.8][0]{$\rightsonic$}
\psfrag{x}[cc][][0.8][0]{$\xi_1$}
\psfrag{oom}[cc][][0.8][0]{$\oOm$}
\psfrag{nom}[cc][][0.8][0]{$\nOm$}
\psfrag{y}[cc][][0.8][0]{$\xi_2$}
\psfrag{lb}[cc][][0.8][0]{$\leftbottom$}
\psfrag{rb}[cc][][0.8][0]{$\rightbottom$}
\psfrag{om}[cc][][0.8][0]{$\Om$}
\psfrag{pb}[cc][][0.8][0]{$P_{\beta}$}
\psfrag{L}[cc][][0.8][0]{$L_{\rm w}$}
\includegraphics[scale=.8]{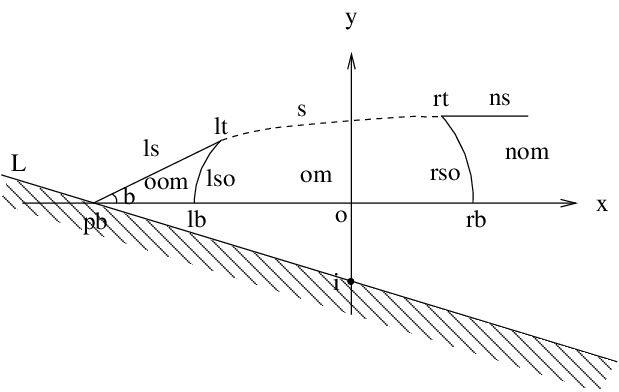}
\caption{Admissible solutions for $\beta<\betac^{(\iv)}$}
\label{fig:regularsol1}
\includegraphics[scale=.8]{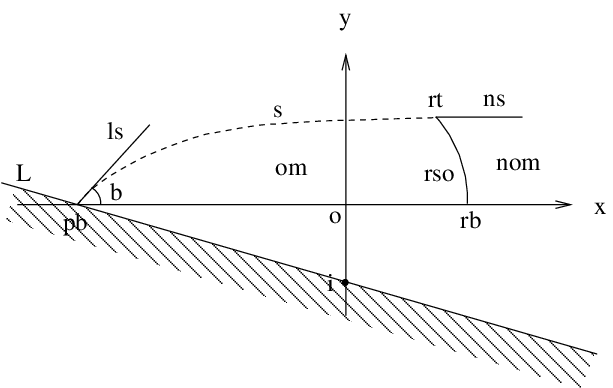}
\caption{Admissible solutions for $\betac^{(\iv)}\le \beta<\betadet$}
\label{fig:regularsol2}
\end{psfrags}
\end{figure}

\begin{definition}
\label{definition-domains-np}
For each $\iv>0$ and $\beta\in[0,\frac{\pi}{2})$, define
\begin{equation}
\label{2-m}
\begin{split}
&\Oi:=(0,-\iv),\quad
\Oo:=(\leftu,0)=(-\iv\tan\beta,0),\quad
O_{\mcl{N}}:=(0,0),\\
&\bdomain:= \R^2_+\setminus\{\bmxi\in \R^2\,:\,\xi_2 \le  f_{\rm w}(\xi_1)\},\\
&\rightsonic:=\der B_{\rightc}(O_{\mathcal{N}})\cap\{\xi_1>0, 0\le\xi_2\le\neta\},\\
&\leftvec=(\cos\beta, \sin\beta).
\end{split}
\end{equation}
For $\ivphi$, $\leftvphi$, and $\rightvphi$ given by
\eqref{def-uniform-ptnl-new},
define
\begin{equation*}
\rightshock=\{\xxi\,:\, \ivphi(\xxi)=\rightvphi(\xxi)\},\quad
\leftshock=\{\xxi\,:\,\ivphi(\xxi)=\leftvphi(\xxi)\}.
\end{equation*}
Let $\nOm$ be the unbounded open region enclosed by ${\rightshock}$, ${\rightsonic}$, and
line $\{(\xi_1,0)\,:\,\xi_1\ge \xi_1^{\rightbottom}\}$
so that $\nOm$ is a fixed domain for all $\beta\in[0, \betadet)$ for fixed $\iv>0$.
Denote the two points $\righttop$ and $\rightbottom$ by{\rm :}

\smallskip
\begin{itemize}
\item $\righttop$ -- the intersection point of line $\xi_2=\neta$ and ${\rightsonic}$,

\smallskip
\item $\rightbottom$ -- the intersection point of the $\xi_1$--axis and ${\rightsonic}$.
\end{itemize}
\smallskip

For each $\iv>0$ and $\beta\in[0,\betac^{(\iv)})$, define
\begin{equation*}
\leftsonic:=\der B_{\leftc}(\Oo)\cap\{\xi_1<0, 0\le\xi_2\le f_{\mcl{O}}(\xi_1)\}.
\end{equation*}
Set the two points $\lefttop$ and $\leftbottom$ as

\smallskip
 \begin{itemize}
\item
$\{\lefttop\}={\leftsonic}\cap \{\xi_2=f_{\mcl{O}}(\xi_1)\}$,

\smallskip
\item $\{\leftbottom\}={\leftsonic}\cap \{\xi_2=0\}$.
\end{itemize}

\smallskip
Let $\oOm$ be the bounded open region enclosed by ${\leftshock}$, ${\leftsonic}$,
and the line segment $\ol{P_{\beta}\leftbottom}$.
\smallskip

By Lemma {\rm \ref{lemma:interval-existence}}, we have
\begin{equation*}
\lim_{\beta\to \betac^{(\iv)}-}|\lefttop-P_{\beta}|=
\lim_{\beta\to \betac^{(\iv)}-}|\lefttop-\leftbottom|=0.
\end{equation*}
This implies that, as $\beta$ tends to $\beta_{\rm s}^{(\iv)}$ from the left,
${\leftsonic}$ and $\ol{\oOm}$ shrink to a single
point $P_{\beta}=\lefttop=\leftbottom$.
Therefore, the definitions of $\leftsonic$, $\lefttop$, and $\leftbottom$ for
$\beta\in[\betac^{(\iv)}, \frac{\pi}{2})$ are given by
\begin{equation}
\label{10-a3}
{\leftsonic}=\{\lefttop\}=\{\leftbottom\}:=\{P_{\beta}\}.
\end{equation}
\end{definition}

\begin{definition}[Admissible solutions with parameters $(\iv,\beta)\in \mathfrak{R}_{\rm weak}$]
\label{def-regular-sol}

Fix $\gam\ge 1$ and $(\iv,\beta)\in \mathfrak{R}_{\rm weak}$,
and let $(\ivphi, \leftvphi, \rightvphi)$ be defined by \eqref{def-uniform-ptnl-new}.
For $\leftshock$ and $\rightshock$ given in Definition {\rm \ref{definition-domains-np}},
define
\begin{equation*}
\leftshockseg:=\leftshock\cap\{-\xi_2^{(\beta)}\cot\beta \le \xi_1\le \xi_1^{\lefttop}\},\qquad
\rightshockseg:=\rightshock\cap \{\xi_1 \ge \xi_1^{\righttop}\}.
\end{equation*}

A function $\vphi\in C^{0,1}_{\rm{loc}}(\bdomain)$ is called \emph{an admissible solution}
corresponding to $(\iv, \beta)$
if $\vphi$ satisfies the following properties{\rm :}
\smallskip

\noindent
\textbf{Case {\rm I}.} $\beta\in (0,\betac^{(\iv)})${\rm :}

\smallskip
\begin{enumerate}
\item[\rm (i)] There exists a shock curve $\shock$ with endpoints $\lefttop=(\oxi,\oeta)$ and $\righttop=(\nxi,\neta)$
such that

\smallskip
\begin{itemize}
\item[(i-1)] Curve $\shock$ satisfies
\begin{equation}
\label{definition-1}
\shock\subset\,\bdomain\setminus\ol{B_{1}(\Oi)},
\end{equation}
where $\der B_{1}(0,-\iv)$ is the sonic circle of the state in $\iOm:=\bdomain\setminus \ol{\oOm\cup\Om\cup\nOm}${\rm ;}
\item[(i-2)] Curve $\shock$ is $C^{2}$ in its relative interior. That is, for any $P\in \shock\setminus\{\lefttop,\righttop\}$,
 there exist a constant $r>0$, a function $f\in C^{2}$, and an orthogonal coordinate system $(S,T)$ in $\R^2$
  such that $\shock\cap B_r(P)=\{S=f(T)\}\cap B_r(P)${\rm ;}

\smallskip
\item[(i-3)] Curve $\ol{\leftshockseg\cup\shock\cup\rightshockseg}$ is $C^1$, including at points $\lefttop$ and $\righttop${\rm ;}

\smallskip
\item[(i-4)] $\shock, \rightsonic, \leftsonic$, and $\Wedge:=\{\xi_2=0, \leftu-\leftc\le\xi_1\le\rightc\}$
do not have common points except for $\lefttop, \righttop, \rightbottom$, and $\leftbottom$.
Thus, ${\shock}\cup {\rightsonic}\cup {\leftsonic}\cup{\Wedge}$ is a closed curve without self-intersection.
Denote by $\Om$ the bounded domain enclosed by this closed curve{\rm .}
\end{itemize}

\smallskip
\item[\rm (ii)] $\vphi$ satisfies the following properties{\rm :}

\smallskip
\begin{itemize}
\item[(ii-1)] $\vphi\in C^{0,1}_{\rm{loc}}(\bdomain)\cap C^1_{\rm{loc}}\bigl(\bdomain\setminus\ol{\leftshockseg\cup\shock\cup \rightshockseg}\bigr)${\rm ;}

\smallskip
\item[(ii-2)] $\vphi\in C^3(\Omega)\cap C^2(\overline{\Omega}\setminus (\overline{\leftsonic}\cup \overline{\rightsonic}))\cap C^1(\overline{\Omega})${\rm ;}

\smallskip
\item[(ii-3)]  For $\oOm$ defined in Definition {\rm \ref{definition-domains-np}},
\begin{equation}\label{1-24}
\vphi=
\begin{cases}
\ivphi\;\;&\text{in $\bdomain\setminus \ol{\oOm\cup\Om\cup\nOm}$},\\
\leftvphi\;\;&\text{in $\oOm$},\\
\rightvphi\;\;&\text{in $\nOm$},
\end{cases}
\end{equation}
where $\ol{\Om^{\mcl{O}}}$ shrinks to $\{P_{\beta}\}=\{\lefttop\}=
\{\leftbottom\}$ when $\beta=\betac^{(\iv)}${\rm ;}

\smallskip
\item[(ii-4)] $\vphi$ satisfies

\smallskip
\begin{itemize}
\item[-] Eq. \eqref{2-1} in $\Om$ with $\rho(|D\vphi|^2, \vphi)$ defined by \eqref{new-density},

\smallskip
\item[-] the slip boundary condition $\vphi_{\xi_2}=0$ on $\Wedge$,

\smallskip
\item[-] the Rankine-Hugoniot conditions{\rm :} $[\vphi]_{\shock}=[\rho(|D\vphi|^2, \vphi)D\vphi\cdot{\bf n}_{\rm sh}]_{\shock}=0$
for the unit normal vector ${\bf n}_{\rm sh}$ to $\shock$ towards the interior of $\Om$.
\end{itemize}
\end{itemize}

\smallskip
\item[\rm (iii)] Eq. \eqref{2-1} is strictly elliptic in $\ol{\Om}\setminus(\overline{\leftsonic}\cup\overline{\rightsonic})$.

\smallskip
\item[\rm (iv)] $\max\{\leftvphi, \rightvphi\}\le \vphi \le \ivphi$ in $\Om$.

\smallskip
\item[\rm (v)] Let $\leftvec$ be the unit vector parallel to $\leftshock$ and oriented
so that $\leftvec\cdot{\bm e}_{1}>0$,
and let $\rightvec$ be the unit vector parallel to $\rightshock$ and oriented
so that $\rightvec\cdot {\bm e}_{1}<0$,
where ${\bm e}_{1}$ is the unit vector in the $\xi_1$--direction, {\it i.e.}, ${\bm e}_{1}=(1,0)$.
That is,
\begin{equation}
\label{def-monot-vecs}
\leftvec =(\cos\beta,\sin\beta), \qquad
\rightvec=(-1,0).
\end{equation}
Then
\begin{equation}\label{cond-ad}
    \der_{\leftvec}(\ivphi-\vphi)\le 0, \quad
    \der_{\rightvec}(\ivphi-\vphi)\le 0 \qquad\,\, \tx{in $\Om$}.
\end{equation}
\end{enumerate}

\medskip
\textbf{Case {\rm II}.} $\beta\in [\betac^{(\iv)},\betadet)${\rm :}

\smallskip
\begin{enumerate}
\item[\rm (i)] There exists a shock curve $\shock$ with endpoints $P_{\beta}=(-\xi_2^{(\beta)}\cot\beta, 0)$ and $\righttop=(\nxi,\neta)$
such that

\smallskip
\begin{itemize}
\item[(i-1)] Curve $\shock$ satisfies
\begin{equation}
\label{definition-1a}
\shock\subset(\bdomain\setminus\ol{B_{1}(\Oi)}),
\end{equation}
where $\der B_{1}(0,-\iv)$ is the sonic circle of the state in $\iOm:=\bdomain\setminus \ol{\Om\cup\nOm}${\rm ;}
\item[(i-2)] Curve $\shock$ is $C^{2}$ in its relative interior{\rm :} for any $P\in \shock\setminus\{P_{\beta},\righttop\}$,
 there exist $r>0$, $f\in C^{2}$, and an orthogonal coordinate system $(S,T)$ in $\R^2$
  so that $\shock\cap B_r(P)=\{S=f(T)\}\cap B_r(P)${\rm ;}

\smallskip
\item[(i-3)] Curve $\ol{\shock\cup\rightshock}$ is $C^1$, including at point $\righttop${\rm ;}

\smallskip
\item[(i-4)] $\shock, \rightsonic$, and $\Wedge:=\{\xi_2=0, -\xi_2^{(\beta)}\cot\beta\le\xi_1 \le\rightc\}$
do not have common points except for $\wedgetip, \righttop$, and $\rightbottom$.
Thus, ${\shock}\cup {\rightsonic}\cup {\Wedge}$ is a closed curve without self-intersection.
Denote by $\Om$ the bounded domain enclosed by this closed curve.
\end{itemize}

\smallskip
\item[\rm (ii)] $\vphi$ satisfies the following properties{\rm :}

\smallskip
\begin{itemize}
\item[(ii-1)] $\vphi\in C^{0,1}_{\rm{loc}}(\bdomain)\cap C^1_{\rm{loc}}(\bdomain\setminus \ol{\shock\cup \rightshockseg})${\rm ;}

\smallskip
\item[(ii-2)] $\varphi\in C^3(\Omega)\cap C^2(\overline{\Omega}\setminus (\{P_{\beta}\}\cup\overline{\rightsonic}))\cap C^1(\overline{\Omega});$

\smallskip
\item[(ii-3)] $D\vphi(\wedgetip)=D\leftvphi(\wedgetip)$ and
\begin{equation}\label{1-24ab}
\vphi=
\begin{cases}
\ivphi\;\;&\text{in $\bdomain\setminus \ol{\Om\cup\nOm}$},\\
\leftvphi\;\;&\text{at $\wedgetip$},\\
\rightvphi\;\;&\text{in $\nOm$}{\rm ;}
\end{cases}
\end{equation}

\item[(ii-4)] $\vphi$ satisfies

\smallskip
\begin{itemize}
\item[-] Eq. \eqref{2-1} in $\Om$ with $\rho(|D\vphi|^2, \vphi)$ defined by \eqref{new-density},

\smallskip
\item[-] the slip boundary condition $\vphi_{\xi_2}=0$ on $\Wedge$,

\smallskip
\item[-] the Rankine-Hugoniot conditions{\rm :} $[\vphi]_{\shock}=[\rho(|D\vphi|^2, \vphi)D\vphi\cdot{\bf n}_{\rm sh}]_{\shock}=0$
for the unit normal vector ${\bf n}_{\rm sh}$ to $\shock$ towards the interior of $\Om$.
\end{itemize}
\end{itemize}

\smallskip
\item[\rm (iii)] Eq. \eqref{2-1} is strictly elliptic in $\ol{\Om}\setminus(\{\wedgetip\}\cup \ol{\rightsonic})$.

\smallskip
\item[\rm (iv)] $\max\{\leftvphi, \rightvphi\}\le \vphi \le \ivphi$ in $\Om$.

\smallskip
\item[\rm (v)] $\vphi$ satisfies \eqref{cond-ad}.
\end{enumerate}
\end{definition}

\begin{remark} The inequalities in \eqref{cond-ad} for two directions $\leftvec$ and $\rightvec$
imply that
\begin{equation}
\label{cond-ad-1}
\der_{{\bf e}}(\ivphi-\vphi)\le 0\qquad\tx{in $\Om\,$ for all ${\bf e}\in \coneclosure$},
\end{equation}
where
\begin{equation}
\label{definition-coneclosure-1}
\coneclosure:=\{a_1\leftvec+a_2\rightvec\,:\,a_1\ge 0,\,\, a_2\ge 0\}.
\end{equation}
\end{remark}

\begin{lemma}
[Entropy condition of admissible solutions]
\label{lemma-entropycond-admsbsol}
Let $\vphi$ be an admissible solution corresponding to $(\iv, \beta)\in \mathfrak{R}_{\rm weak}$
in the sense of Definition {\rm \ref{def-regular-sol}},
and let $\shock$ be the curved shock satisfying condition {\rm (i)} of Definition {\rm \ref{def-regular-sol}}.
Let ${\bm\nu}$ be the unit normal vector to $\shock$ towards the interior of $\Om$.
Then the following properties hold{\rm :}

\smallskip
\begin{itemize}
\item[(a)]
$\der_{\bm\nu}\ivphi>\der_{\bm\nu}\vphi>0\quad\tx{on $\ol{\shock}$}${\rm ;}

\smallskip
\item[(b)] Let
\begin{equation*}
M_{\infty,\bm\nu}:=\frac{\der_{\bm\nu}\ivphi}{c(|D\ivphi|^2, \ivphi)}=\der_{\bm\nu}\ivphi,\quad
M_{\bm\nu}:=\frac{\der_{\bm\nu}\vphi}{c(|D\vphi|^2,\vphi)}
\end{equation*}
for
\begin{equation}
\label{definition-localsoundsp}
c(|{\bf q}|^2, z)=\rho^{\frac{\gam-1}{2}}(|{\bf q}|^2, z),
\end{equation}
where $\rho(|{\bf q}|^2, z)$ is defined by \eqref{new-density}. Then
\begin{equation*}
0<M_{\bm\nu}<1<M_{\infty, \bm\nu}\qquad\tx{on $\ol{\shock}$.}
\end{equation*}
\end{itemize}
\end{lemma}

\begin{proof}
Denote $w:=\ivphi-\vphi$.
From \eqref{2-1}, \eqref{new-density}, and \eqref{2-4-a3}, it can directly be checked that
\begin{equation*}
(c^2-\vphi_{\xin}^2)w_{\xin\xin}-2\vphi_{\xin}\vphi_{\etan}w_{\xin\etan}
+(c^2-\vphi_{\etan}^2)w_{\etan\etan}=0\qquad\tx{in $\Om$}
\end{equation*}
for $c^2=\rho^{\gam-1}(|D\vphi|^2, \vphi)$,
where $\rho(|D\vphi|^2,\vphi)$ is given by \eqref{new-density}.
By condition (iii) of Definition \ref{def-regular-sol}, the minimum principle applies to $w$
so that $w$ cannot attain its minimum in $\Om$,
unless it is a constant in $\ol{\Om}$.
By conditions (ii) and (iv) of Definition \ref{def-regular-sol},
we see that $w\ge 0$ in $\ol{\Om}$, and $w=0$ on $\ol{\shock}$.
Furthermore, $w$ is not a constant in $\ol{\Om}$,
because $\der_{\etan}w=-\iv$ on $\Wedge$ by \eqref{2-4-b6}
and the slip boundary condition $\der_{\etan}\vphi=0$ on $\Wedge$, stated in (ii-4) of Definition \ref{def-regular-sol}.
Then it follows from Hopf's lemma that $\der_{\bm\nu}w>0$ on ${\shock}$. This implies that
\begin{equation}
\label{ent-cond1}
\der_{\bm\nu}\ivphi>\der_{\bm\nu}\vphi\qquad\tx{on $\shock$}.
\end{equation}
If $\der_{\bm\nu}\vphi(P)=0$ for some $P\in \shock$, then it follows from
the condition: $\rho(|D\vphi|^2,\vphi)\der_{\bm\nu}\vphi(P)=\der_{\bm\nu}\ivphi(P)$ stated in (ii-4) of Definition \ref{def-regular-sol}
that $\der_{\bm\nu} \ivphi(P)=0$, which is impossible, due to \eqref{ent-cond1}. Therefore, we have
\begin{equation}
\label{ent-cond2}
|\der_{\bm\nu}\vphi|>0\qquad \tx{on $\shock$.}
\end{equation}

By conditions (ii-2)--(ii-3) of Definition \ref{def-regular-sol}, we have
\begin{equation*}
D\vphi(\righttop)=D\rightvphi(\righttop).
\end{equation*}
Then it follows from the definitions of $(\ivphi, \rightvphi)$  given in \eqref{def-uniform-ptnl-new}
and conditions (ii-4) and (iv) of Definition \ref{def-regular-sol} that
\begin{align}
&{\bm\nu}(\righttop)=\frac{D\ivphi-D\rightvphi}{|D\ivphi-D\rightvphi|}=(0,-1),\\
\label{endpoint-rt}
& \der_{\bm\nu}(\ivphi-\vphi)(\righttop)=|D\ivphi-D\rightvphi|=\iv>0,\quad
\der_{\bm\nu}\vphi(\righttop)=\der_{\bm \nu}\rightvphi(\righttop)=\xi_2^{\righttop}>0.
\end{align}
Similarly, at $\lefttop$, we have
\begin{equation*}
D\vphi(\lefttop)=D\leftvphi(\lefttop),
\end{equation*}
so that \eqref{12-25}, \eqref{1-36}, \eqref{xi2m-new1}, and \eqref{xi2m-new2} yield that
\begin{align}
&{\bm\nu}(\lefttop)=\frac{D\ivphi-D\leftvphi}{|D\ivphi-D\leftvphi|}=(\sin\beta, -\cos\beta),\notag\\
\label{endpoint-lt1}
& \der_{\bm\nu}(\ivphi-\vphi)(\lefttop)=|D\ivphi-D\leftvphi|=\iv\sec\beta>0,\\
\label{endpoint-lt2}
&\der_{\bm\nu}\vphi(\lefttop)=\der_{\bm\nu}\leftvphi(\lefttop)=\der_{\nu}\ivphi(\lefttop)-\iv \sec\beta=\xi_2^m>0.
\end{align}
Then statement (a) follows directly from \eqref{ent-cond1}--\eqref{endpoint-lt2}
and the continuity of $\der_{\bm\nu}\vphi$ along $\shock$ up to its endpoints $\lefttop$ and $\righttop$.

Note that the calculations given in \eqref{leftrho-expression}--\eqref{1-1} are still valid
when $(\leftrho, \oM, \iM)$ are replaced by $(\rho, M_{\bm\nu}, M_{\infty, {\bm\nu}})$ on $\ol{\shock}$.
Then we see that, on $\ol{\shock}$,
\begin{align}
\label{M-eqn-1}
&\rho^{\frac{\gam+1}{2}}=\frac{M_{\infty, {\bm\nu}}}{M_{\bm\nu}},\\
\label{M-eqn-2}
&\big(1+\frac{\gam-1}{2}M_{\bm\nu}^2\big)|M_{\bm\nu}|^{\frac{-2(\gam-1)}{\gam+1}}
=\big(1+\frac{\gam-1}{2}M_{\infty,{\bm\nu}}^2\big)|M_{\infty, {\bm\nu}}|^{\frac{-2(\gam-1)}{\gam+1}}.
\end{align}
This is because \eqref{leftrho-expression}--\eqref{1-1} are all derived from
the Rankine-Hugoniot conditions stated in  Definition \ref{def-regular-sol}(ii-4).
By the result obtained in statement (a) and the Rankine-Hugoniot condition:
$\rho \der_{\bm\nu}\vphi=\der_{\bm\nu}\ivphi$ on $\ol{\shock}$,
\eqref{M-eqn-1} implies that $\frac{M_{\infty, {\bm\nu}}}{M_{\bm\nu}}>1$ on $\ol{\shock}$.
Since $(M_{\bm\nu}, M_{\infty, {\bm\nu}})$ satisfy \eqref{M-eqn-2} and $M_{\infty, {\bm\nu}}\neq M_{\bm\nu}$
on $\ol{\shock}$, it follows from the observation right after \eqref{1-1} that
\begin{equation*}
0<M_{{\bm\nu}}<1<M_{\infty, {\bm\nu}}\qquad\tx{on $\ol{\shock}$}.
\end{equation*}
This completes the proof of statement (b).
\end{proof}

In \eqref{fo}--\eqref{definition-f-w} and Definition \ref{def-regular-sol}, the values of
$\xi_1^{(\beta)}$, $\xi_2^{(\beta)}$, $\theta_{\infty}$, $\leftc$, and $\leftu$
depend continuously on  $\beta\in(0, \frac{\pi}{2})$ with
\begin{equation*}
\lim_{\beta\to 0+} (\xi_1^{(\beta)}, \xi_2^{(\beta)}, \theta_{\infty}, \leftc, \leftu)
=(-\infty, \neta, \pi, \rightc, 0).
\end{equation*}
As a result, we obtain
\begin{equation*}
\begin{split}
&\lim_{\beta\to 0+}|\lefttop-(-\xi_1^{\righttop}, \neta)|=0=\lim_{\beta\to 0+}|\leftbottom-(-\rightc,0)|,\\
&\lim_{\beta\to 0+}\|\leftvphi-\rightvphi\|_{C^{3}(B_R({\bf 0}))}=0\qquad\tx{for any $R>0$}.
\end{split}
\end{equation*}
For $\beta=0$, we define $\lefttop$, $\leftbottom$, $\Lbeta|_{\beta=0}$, and $\leftshockseg|_{\beta=0}$ by
\begin{equation}
\label{definition-domain-norshock-new}
\begin{split}
&\lefttop=(-\xi_1^{\righttop}, \neta),\quad \leftbottom=(-\rightc,0),\\
&\Lbeta|_{\beta=0}:=\R\times \R_+,\quad \leftshockseg|_{\beta=0}=\{(\xi_1,\neta)\,:\,\xi_1\le -\xi_1^{\righttop}\}.
\end{split}
\end{equation}
Then two points $\lefttop$ and $\leftbottom$ depend continuously on $\beta\in[0, \frac{\pi}{2})$,
so that $\Lbeta$ and $\leftshockseg$ depend continuously on $\beta\in[0, \frac{\pi}{2})$.
Using this, we extend Definition \ref{def-regular-sol} up to $\beta=0$.

\begin{definition}[Admissible solutions when $\beta=0$]
\label{def-regular-sol-normal}
Given $\gam\ge 1$ and $\iv>0$, a function $\vphi\in C^{0,1}(\R\times \R_+)$ is called \emph{an admissible solution}
corresponding to $(\iv, 0)$ if $\vphi$ satisfies the following properties{\rm :}

\smallskip
\begin{enumerate}
\item [\rm (i)] There exists a shock $\shock$ with endpoints $\lefttop=(-\nxi,\neta)$ and $\righttop=(\nxi,\neta)$
such that

\smallskip
\begin{itemize}
\item[(i-1)] Curve $\shock$ satisfies
\begin{equation}
\label{definition-1-s}
\shock\subset(\R\times \R_+)\setminus\ol{B_{1}(\Oi)},
\end{equation}
where $\der B_{1}(\Oi)$ is the sonic circle of state $\Oi=(0, -\iv)$
in $\iOm:=(\R\times \R_+)\setminus \ol{\oOm\cup\Om\cup \nOm}${\rm ;}

\smallskip
\item[(i-2)] Curve $\shock$ is $C^{2}$ in its relative interior; that is, for any $P\in \shock\setminus\{\lefttop,\righttop\}$,
 there exist $r>0$, $f\in C^{2}$, and an orthogonal coordinate system $(S,T)$ in $\R^2$
  so that $\shock\cap B_r(P)=\{S=f(T)\}\cap B_r(P)${\rm ;}

\smallskip
\item[(i-3)] Curve $\ol{\leftshockseg\cup\shock\cup\rightshockseg}$ is $C^1$, including at points $\lefttop$ and $\righttop${\rm ;}

\smallskip
\item[(i-4)] $\shock, \rightsonic, \leftsonic$, and $\Wedge:=\{(\xi_1, 0)\,:\,-\rightc<\xi_1 <\rightc\}$ do not
have common points, and  ${\shock}\cup {\rightsonic}\cup {\leftsonic}\cup{\Wedge}$ is a closed curve without self-intersection.
Denote by $\Om$ the bounded domain enclosed by this closed curve.
\end{itemize}

\smallskip
\item[\rm (ii)] $\vphi$ satisfies the following properties{\rm :}

\smallskip
\begin{itemize}
\item[(ii-1)] $\vphi\in C^{0,1}(\R\times \R_+)\cap C^1\bigl((\R\times \R_+)\setminus \ol{\leftshockseg\cup\shock\cup \rightshockseg}\bigr)${\rm ;}

\smallskip
\item[(ii-2)] $\vphi\in C^3(\Omega)\cap C^2(\ol{\Om}\setminus(\overline{\leftsonic}\cup\overline{\rightsonic}))\cap C^1(\ol{\Om})${\rm ;}

\smallskip
\item[(ii-3)]
\begin{equation*}
\vphi=
\begin{cases}
\ivphi\;\;&\text{in $(\R\times \R_+)\setminus \ol{\oOm\cup\Om\cup\nOm}$},\\
\rightvphi\;\;&\text{in $\oOm\cup\nOm$};
\end{cases}
\end{equation*}

\item[(ii-4)] $\vphi$ satisfies

\smallskip
\begin{itemize}
\item[-] Eq. \eqref{2-1} in $\Om$ with $\rho(|D\vphi|^2, \vphi)$ defined by \eqref{new-density},

\smallskip
\item[-] the slip boundary condition: $\vphi_{\etan}=0$ on $\Wedge$,

\smallskip
\item[-] the Rankine-Hugoniot conditions{\rm :} $[\vphi]_{\shock}=[\rho(|D\vphi|^2, \vphi)D\vphi\cdot{\bf n}_{\rm sh}]_{\shock}=0$
for the unit normal vector ${\bf n}_{\rm sh}$ to $\shock$ towards the interior of $\Om$.
\end{itemize}
\end{itemize}

\smallskip
\item[\rm (iii)] Eq. \eqref{2-1} is strictly elliptic in $\ol{\Om}\setminus(\overline{\leftsonic}\cup\overline{\rightsonic})$.

\smallskip
\item[\rm (iv)] $ \rightvphi\le \vphi \le \ivphi$ in $\Om$.

\smallskip
\item[\rm (v)] $\der_{{\bf e}}(\ivphi-\vphi)\le  0 $ in $\Om$ for all ${\bf e}\in \R\times \R^+$.
\end{enumerate}
\end{definition}

\begin{remark}
\label{remark-condition-v}
Condition {\rm (v)} of Definition {\rm \ref{def-regular-sol-normal}} is a continuous extension
of condition {\rm (v)} of Definition {\rm \ref{def-regular-sol}} in the sense that

\smallskip
\begin{itemize}
\item[(i)] $\coneclosure$ for $\beta>0$ defined by \eqref{definition-coneclosure-1} monotonically increases
as $\beta>0$ decreases in the sense that, if $0<\beta_1<\beta_2<\frac{\pi}{2}$, then
\begin{equation*}
\coneclosure|_{\beta_2}\subset \coneclosure|_{\beta_1};
\end{equation*}

\item[(ii)] $\displaystyle{\overline{\cup_{0<\beta<\frac{\pi}{2}}\coneclosure|_{\beta} }=\R\times \R^+}$.
\end{itemize}
\end{remark}

\begin{remark} \label{remark-admsol-wksol}
Similarly to Definition {\rm \ref{definition-weak-sol}}, it can directly be checked that
any admissible solution corresponding to $(\iv, \beta)\in\mathfrak{R}_{\rm{weak}}\cup\{(\iv, 0)\,:\,\iv>0\}$
in the sense of Definition {\rm \ref{def-regular-sol}} or {\rm \ref{def-regular-sol-normal}} satisfies
the following properties{\rm :}

\smallskip
\begin{itemize}
    \item [(i)] $\vphi\in W^{1,1}_{\rm{loc}}(\Lbeta)${\rm ;}

  \smallskip
    \item [(ii)] $\rho(|D\vphi|^2,\vphi)>0$ in $\Lbeta$ for $\rho(|D\vphi|^2,\vphi)$ defined by \eqref{new-density}{\rm ;}

  \smallskip
    \item [(iii)] $(\rho(|D\vphi|^2, \vphi),\, \rho(|D\vphi|^2, \vphi)|D\vphi|) \in L^1_{\rm{loc}}(\Lbeta) ${\rm ;}


  \smallskip
    \item [(iv)] For every $\zeta\in C^{\infty}_0(\R^2)$,
    \begin{equation*}
      \int_{\Lbeta} \left(\rho(|D\vphi|^2, \vphi)D\vphi\cdot D\zeta-2\rho(|D\vphi|^2, \vphi)\zeta\right)\,{\rm d}\bmxi =0.
    \end{equation*}
\end{itemize}

\smallskip
Specifically, property {\rm (iv)} here is obtained by condition {\rm (ii)} of
Definitions {\rm \ref{def-regular-sol}} and {\rm \ref{def-regular-sol-normal}},
and via integration by parts. Property {\rm (iv)} indicates that any admissible solution $\vphi$ is
a {\emph{weak solution}} of the boundary value problem rewritten from Problem {\rm \ref{problem-2}}
with respect to parameters $(\iv, \beta)$.
In particular, a function $\varphi$ satisfying {\rm (i)}--{\rm (iv)} is a weak solution of the boundary value problem consisting
of equation \eqref{2-1} in $\Lambda_\beta$ and the slip boundary condition $\der_{\bm \nu}\varphi=0$ on $\partial \Lambda_\beta$,
where we note that  $\partial \Lambda_\beta\subset \{(\xi_1, \xi_2)\; : \;\xi_2=0\}\cup L_{\rm w}$.

\end{remark}

\begin{lemma}\label{lemma-uniq-normal-sol}
For any given $\gam\ge 1$ and $\iv>0$, there exists at least one admissible solution corresponding to $(\iv, 0)$
in the sense of Definition {\rm \ref{def-regular-sol-normal}}.

\begin{proof}
The conditions stated in (ii-4) and (v) of Definition \ref{def-regular-sol-normal} imply that
$$
\shock=\{(\xin, \neta)\,:\,-\nxi<\xin<\nxi\};
$$
that is, $\ol{\leftshockseg\cup\shock\cup\rightshockseg}$ is a normal shock.
Therefore, the pseudo-subsonic region $\Om$ is enclosed by $\leftsonic, \rightsonic, \Wedge$,
and the line segment $(-\nxi, \nxi)\times\{\neta\}$.
It can directly be checked that a function $\vphi_{\rm norm}\in C^{0,1}(\bdomain|_{\beta=0})$ defined by
\begin{equation*}
\vphi_{\rm norm}=
\begin{cases}
\ivphi\;\;&\text{in $(\R\times \R_+)\setminus \ol{\oOm\cup\Om\cup\nOm}$},\\
\rightvphi\;\;&\text{in $\oOm\cup\Om\cup\nOm$}
\end{cases}
\end{equation*}
is an admissible solution corresponding to $(\iv, 0)$ in the sense of Definition \ref{def-regular-sol-normal}.
\end{proof}
\end{lemma}

For a fixed $(\iv, \beta)\in \mathfrak{R}_{\rm{weak}}$,
let $\vphi$ be an admissible solution corresponding to $(\iv, \beta)$ in the sense of Definition \ref{def-regular-sol}.
Let $(\iu, \leftuop)$ be given by $(\iu, \leftuop)=\mcl{T}^{-1}(\iv, \beta)\in \mathfrak{P}_{\rm{weak}}$ for map $\mcl{T}$
from Lemma \ref{lemma-parameters}. Let $\displaystyle{\theta_{\rm w}}$ be given by \eqref{expression-tw-op}. For each $\xxi'=(\xi_1',\xi_2')\in\Lbeta$, let $\xxi=(\xi_1,\xi_2)$ be given by
\begin{equation*}
\xxi^\top=\begin{pmatrix}
\phantom{-}\cos\theta_{\rm w}&\sin \theta_{\rm w}\\
-\sin \theta_{\rm w}&\cos \theta_{\rm w}
\end{pmatrix}^{-1}
\left((\xxi')^\top+\begin{pmatrix}
\iu\cos\theta_{\rm w}\\0
\end{pmatrix}
\right).
\end{equation*}
This is the inverse transformation of \eqref{2-a}. Finally, let a function $\til{\vphi}$ be given by
\begin{equation}
\label{solution-transformation}
\til{\vphi}(\xxi)=\vphi(\xxi')+\frac 12(\iu\cos \theta_{\rm w})^2\qquad\tx{for $\xxi'\in\Lbeta$.}
\end{equation}
Then $\til{\vphi}$ is an admissible solution corresponding to $(\iu,u_0)\in \mathfrak{P}_{\rm{weak}}$ in the sense
of Definition \ref{def-regular-sol-op}. From this perspective, Theorem \ref{theorem-1-op} is equivalent to the following theorem:

\begin{theorem}[Existence of admissible solutions]
\label{theorem-0}
For any given $\gam\ge 1$ and $(\iv,\beta)\in\mathfrak{R}_{\rm weak}$,
there exists an admissible solution in the sense of Definition {\rm \ref{def-regular-sol}}.
\end{theorem}

\begin{remark}[Non-existence of self-similar strong shock solutions]
\label{remark-nonexistence}
Fix $\gam\ge 1$. For $(\iv,\beta)\in \mathfrak{R}_{\rm detach}\cup \mathfrak{R}_{\rm strong}$,
let $(\Lambda_{\beta}, \ivphi, \leftvphi, \rightvphi)$ be defined as in Definition {\rm \ref{definition-domains-np}}.
We call $\vphi\in C^{0,1}(\bdomain)$
an admissible solution corresponding to $(\iv, \beta)\in  \mathfrak{R}_{\rm detach}\cup \mathfrak{R}_{\rm strong}$
if it satisfies conditions {\rm (i)}--{\rm (v)} stated in Definition {\rm \ref{def-regular-sol}}
for \textbf{Case {\rm II}}.
By the convexity of the shock polar for steady potential flow,
which is shown in \emph{Appendix {\rm \ref{section-appendix}}},
and condition {\rm (iv)} of Definition {\rm \ref{def-regular-sol}},
it follows from the non-existence result as proved in \emph{Appendix {\rm B}} (see also \cite{Elling2})
that there exists no admissible solution corresponding to
$(\iv,\beta)\in \mathfrak{R}_{\rm strong}$
in the sense of Definition {\rm \ref{def-regular-sol}}.

The existence of admissible solutions corresponding to $(\iv, \beta_{\rm d}^{(\iv)})$
is still an open question.
\end{remark}

\begin{theorem}[Regularity of admissible solutions]
\label{theorem-3}
Given $\gam\ge 1$ and $(\iv, \beta)\in \mathfrak{R}_{\rm weak}$, let $\vphi$ be a corresponding admissible solution
with the curved shock $\shock$ in the sense of Definition {\rm \ref{def-regular-sol}}.
Then the following properties hold{\rm :}
\medskip

\noindent
\textbf{Case {\rm I}.} $\beta\in(0,\betac^{(\iv)})${\rm :}

\smallskip
\begin{itemize}
\item[(a)] $\shock$ is $C^{\infty}$ in its relative interior,
and $\vphi\in C^{\infty}(\ol{\Om}\setminus({\overline{\leftsonic}}\cup{\overline{\rightsonic}}))\cap C^{1,1}(\ol{\Om})${\rm ;}

\smallskip
\item[(b)]  Define a set $\mcl{D}$ by
\begin{equation}
\label{def-D-domain}
\mcl{D}=\bdomain\cap\{\bmxi\,:\,\max\{\leftvphi(\bmxi), \rightvphi(\bmxi)\}<\ivphi(\bmxi)\}.
\end{equation}
For a constant $\sigma>0$, define $\oD_{\sigma}$ and $\nD_{\sigma}$ by
\begin{equation}\label{2-a5}
\begin{split}
&\oD_{\sigma}=\mcl{D}\cap \{\bmxi\,:\,{\rm dist}\{\bmxi,\leftsonic\}<\sigma\}\cap B_{\leftc}(\Oo),\\
&\nD_{\sigma}=\mcl{D}\cap \{\bmxi\,:\,{\rm dist}\{\bmxi,\rightsonic\}<\sigma\}\cap B_{\rightc}(O_{\mcl{N}})\\
\end{split}
\end{equation}
for $\rightc=\rightrho^{(\gam-1)/{2}}$, $\leftc=\leftrho^{(\gam-1)/{2}}$, $\Oo=(\leftu,0)$, and $O_{\mcl{N}}:=(0,0)$.
Fix any point  $\bmxi_0 \in (\overline{\leftsonic}\cup \overline{\rightsonic})\setminus\{\lefttop, \righttop\}$,
and denote $d:={\rm dist}\{\bmxi_0,\shock\}$.
Then, for any $\alp\in(0,1)$,
there exists a constant $K<\infty$ depending on $(\iv, \gam, \eps_0, \alp, d)$ and
$\|\vphi\|_{C^{1,1}(\Om\cap(\oD_{\eps_0}\cup\nD_{\eps_0}))}$ such that
\begin{equation}
    \label{Op1}
    \|\vphi\|_{2,\alp,\ol{\Om\cap B_{d/2}(\bmxi_0)\cap(\oD_{\eps_0/2}\cup\nD_{\eps_0/2})}}\le K;
    \end{equation}

\item[(c)] For any $\bmxi_0 \in (\overline{\leftsonic}\cup\overline{\rightsonic})\setminus\{\lefttop, \righttop\}$,
\begin{equation}
\label{Op2}
\lim_{\bmxi\to\bmxi_0 \atop
\bmxi \in\Om}\bigl(D_{rr}\vphi-D_{rr}\max\{\leftvphi,\rightvphi\}\bigr)(\bmxi)=\frac{1}{\gam+1},
\end{equation}
where $r=|\bmxi|$ near $\rightsonic$ and $r=|\bmxi-(\leftu,0)|$ near $\leftsonic${\rm ;}

\smallskip
\item[(d)] Limits $\displaystyle\lim_{\bmxi \to \lefttop\atop \bmxi \in\Om}D^2\vphi$ and
$\displaystyle\lim_{\bmxi \to \righttop\atop \bmxi \in\Om}D^2\vphi$ do not exist{\rm ;}

\smallskip
\item[(e)]  $\overline{\leftshockseg\cup\shock\cup\rightshockseg}$
is a $C^{2,\alp}$--curve for any $\alp\in(0, 1)$,
including at points $\lefttop$ and $\righttop$.
\end{itemize}

\smallskip
\noindent
\textbf{Case {\rm II}.} $\beta\in[\betac^{(\iv)},\betadet)${\rm :}

\smallskip
\begin{itemize}
\item[(a)] $\shock$ is $C^{\infty}$ in its relative interior,
and
$\vphi\in C^{\infty}(\ol{\Om}\setminus(\{\wedgetip\}\cup{\ol{\rightsonic}}))\cap C^{1,1}(\ol{\Om}\setminus\{P_{\beta}\})\cap C^{1,\bar{\alp}}(\ol{\Om})$
for some $\bar{\alp}\in(0,1)${\rm ;}

\smallskip
\item[(b)]  For a constant $\sigma>0$, let $\nD_{\sigma}$ be defined by \eqref{2-a5}.
Fix any point $\bmxi_0 \in \ol{\rightsonic}\setminus\{\righttop\}$, and denote $d:={\rm dist}\{\bmxi_0,\shock\}$.
Then, for any $\alp\in(0,1)$,
there exists a constant $K<\infty$ depending on $(\iv, \gam, \eps_0, \alp,d)$ and $\|\vphi\|_{C^{1,1}(\Om\cap\nD_{\eps_0})}$
such that
    \begin{equation}
    \label{Op1-a}
    \|\vphi\|_{2,\alp,\ol{\Om\cap B_{d/2}(\bmxi_0)\cap\nD_{\eps_0/2}}}\le K;
    \end{equation}

\item[(c)] For any $\bmxi_0 \in \overline{\rightsonic}\setminus\{\righttop\}$,
\begin{equation}
\label{Op2-a}
\lim_{\bmxi \to \bmxi_0 \atop
\bmxi \in\Om}\bigl(D_{rr}\vphi-D_{rr}\rightvphi\bigr)(\bmxi)=\frac{1}{\gam+1},
\end{equation}
where $r=|\bmxi|${\rm ;}

\smallskip
\item[(d)] Limit $\displaystyle\lim_{\bmxi \to \righttop\atop \bmxi \in\Om}D^2\vphi$ does not exist{\rm ;}

\smallskip
\item[(e)]
$\overline{\shock\cup\rightshockseg}$
is a $C^{1,\bar{\alp}}$--curve for the same $\bar{\alp}$ as in statement {\rm (a)}.
Furthermore, curve $\overline{\shock\cup\rightshockseg}\setminus \{P_{\beta}\}$ is $C^{2,\alp}$
for any $\alp\in(0, 1)$,
including at point $\righttop$.
\end{itemize}
\end{theorem}

Since Theorems \ref{theorem-1-op}--\ref{theorem-2-op} follow directly from Theorems \ref{theorem-0} and \ref{theorem-3}
through \eqref{solution-transformation}, the rest of the monograph is devoted to establishing Theorems \ref{theorem-0} and \ref{theorem-3}.

We will prove Theorem \ref{theorem-0} by solving the following free boundary problem:

\begin{problemL}[Free boundary problem] \label{fbp}
Given $\gam\ge 1$ and $(\iv,\beta)\in\mathfrak{R}_{\rm weak}$, define $\vphi_{\beta}$ and $\sonic$ by
\begin{equation}
\label{def-vphi-beta}
\vphi_{\beta}:=\max\{\leftvphi, \rightvphi\},\qquad
\sonic:=\leftsonic\cup\rightsonic.
\end{equation}
Find a curved shock $\shock$ and a function
$\varphi\in C^3(\Omega)\cap C^2(\overline{\Omega}\setminus (\overline{\leftsonic}\cup\overline{\rightsonic}))\cap C^1(\overline{\Omega})$
satisfying the following{\rm :}
\begin{eqnarray}
\label{3-a0}
\text{Eq. \eqref{2-1}} && \text{in $\Om$},\\
\label{3-a1}
\vphi=\vphi_{\beta},\;\; D\vphi=D\vphi_{\beta} && \text{on $\sonic$},\\
\label{3-a2}
\der_{\xi_2}\vphi=0 &&\text{on $\Wedge$},\\
\label{3-a3}
\vphi=\ivphi,\;\; \rho D\vphi\cdot{\bm \nu}_{\rm sh}=D\ivphi\cdot{\bm \nu}_{\rm sh} && \text{on $\shock$},
\end{eqnarray}
where ${\bm \nu}_{\rm sh}$ is the unit normal vector to $\shock$ towards the interior of $\Om$, and $\rho$ is defined by \eqref{new-density}.
Note that $\leftsonic$ is a closed portion of a circle,
which becomes one point for $\beta\ge \betasonic$.
Therefore, the boundary condition \eqref{3-a1} on $\leftsonic$ becomes a one-point boundary condition for $\beta\ge \betasonic$.
\end{problemL}

\begin{remark}
\label{remark-vphibeta}
It can be checked from the definitions of $(\leftvphi, \rightvphi)$ given in \eqref{def-uniform-ptnl-new}
that,
for each $\beta\in(0, \frac{\pi}{2})$, there exists a unique $\xi_1^*$ such that
\begin{equation*}
\vphi_{\beta}(\xi_1,\xi_2)=\begin{cases}
\leftvphi&\tx{for $\xi_1< \xi_1^*$},\\
\leftvphi=\rightvphi &\tx{at $\xin=\xi_1^*$},\\
\rightvphi&\tx{for $\xi_1>\xi_1^*$}.
\end{cases}
\end{equation*}
Moreover, $\xi_1^*$ satisfies that $f_{\mcl{O}}(\xi_1^*)=\xi_2^{\mcl{N}}$ and $\xi_1^{P_{\beta}}<\xi_1^*<0$.
In particular, $\varphi_{\beta}=\leftvphi$ on $\leftsonic$ and $\varphi_{\beta}=\rightvphi$ on $\rightsonic$.
\end{remark}

\section{Further Features of Problem 2.34}

\begin{figure}[htp]
\centering
\begin{psfrags}
\psfrag{ti}[cc][][0.7][0]{$\itheta$}
\psfrag{to}[cc][][0.7][0]{$\otheta$}
\psfrag{n}[cc][][0.7][0]{$(\nxi,\neta)$}
\psfrag{o}[cc][][0.7][0]{$(\oxi,\oeta)$}
\psfrag{sn}[cc][][0.7][0]{$\rightshock$}
\psfrag{so}[cc][][0.7][0]{$\leftshock$}
\psfrag{iv}[cc][][0.7][0]{$\;\;-\iv$}
\psfrag{d0}[cc][][0.7][0]{$d=1\;\;\;$}
\includegraphics[scale=.6]{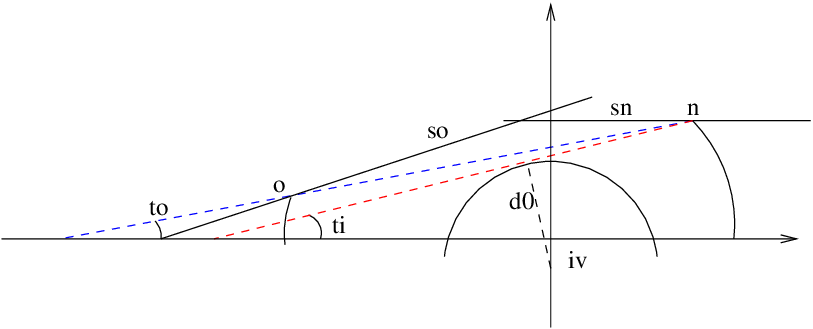}
\includegraphics[scale=.6]{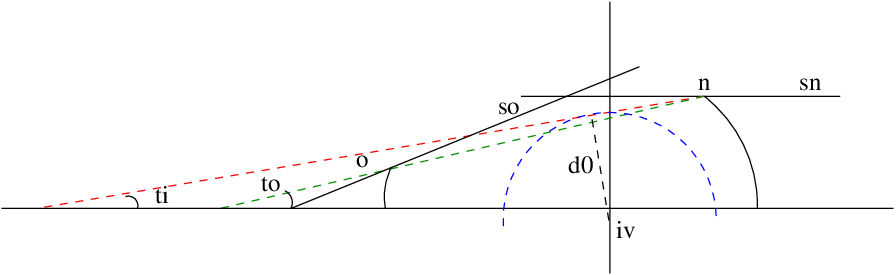}
\caption{Top: $\otheta<\itheta$; $\,$  Bottom: $\otheta>\itheta$ }
\label{fig:el}
\end{psfrags}
\end{figure}

Fix $\gam\ge 1$. For $(\iv, \beta)\in \mathfrak{R}_{\rm weak}$ with $\beta<\betac^{(\iv)}$,
let $\lefttop$ and $\righttop$ be the points as defined in Definition \ref{definition-domains-np}.
Let $\oL$ be the line segment connecting $\lefttop$
with $\righttop$.
For $0<\iv<1$, there exists a unique line $\iL$ that passes through $\righttop$ and is tangential to $\der B_1(\Oi)$
so that the intersection point of $\iL$ with $\der B_1(\Oi)$ has a negative $\xi_1$--coordinate; see Fig. \ref{fig:el}.
Let $\tan\otheta$ and $\tan\itheta$ be the slopes of $\oL$ and $\iL$, respectively. Then
\begin{equation*}
{\rm dist}(\oL, \Oi))
\begin{cases}
>1&\quad\tx{iff $\otheta<\itheta$},\\
<1&\quad \tx{iff $\otheta>\itheta$}.
\end{cases}
\end{equation*}
Note that  $\tan\itheta$ is independent of $\beta\in(0, \betac^{(\iv)})$, and $\Oi=(0,-\iv)$.

\begin{proposition}
\label{proposition:el-condition}
For any given $\gam \ge 1$, there exists a constant $v_*\in(0,1)$ so that,
if $0<\iv<v_*$, there is $\hat{\beta}^{(\iv)}\in(0,\betac^{(\iv)})$
such that
\begin{align}
\label{1-16-nversion}
&{\rm dist}(\oL, \Oi)>1\;\;\qquad \text{for $\beta\in(0,\hat\beta^{(\iv)})$},\\
\label{1-16nv-ox}
&{\rm dist}(\oL, \Oi)<1\;\;\qquad \text{for $\beta\in(\hat\beta^{(\iv)},\betac^{(\iv)})$}.
\end{align}

\begin{proof}
In this proof, we consider only case $\gam>1$. Case $\gam=1$ can be handled similarly.
The proof is divided into seven steps.

\medskip
{\bf 1.} {\emph{Claim{\rm :} For each $\gam>1$,  $\righttop=(\nxi, \neta)$ and $\rightrho$ depend continuously on $\iv>0$ and}}
\begin{equation}
\label{limits-1}
\lim_{\iv\to 0+}\nxi=0, \quad \lim_{\iv\to 0+}\nrho=\lim_{\iv\to 0+}\neta=1.
\end{equation}
Substituting $\leftrho=\rightrho$ into \eqref{2-4-a3}, we have
\begin{equation}
\label{1-7}
F_1(\nrho,\iv):=\frac{\nrho^{\gam-1}-1}{\gam-1}(\nrho-1)-\frac{1}{2}\iv^2(\nrho-1)-\iv^2=0.
\end{equation}
We differentiate $F_1$ with respect to $\rightrho$ to obtain
\begin{equation}
\label{F1-deriv}
\der_{\nrho}F_1=\nrho^{\gam-2}(\nrho-1)+\frac{\nrho^{\gam-1}-1}{\gam-1}-\frac{1}{2}\iv^2.
\end{equation}
Using \eqref{1-7} to obtain
that $\displaystyle{\frac{\rightrho^{\gam-1}-1}{\gam-1}=\frac 12 \iv^2+\frac{\iv^2}{\rightrho-1}}$,
substituting this expression into \eqref{F1-deriv}, and then applying $\rightrho>1$, we have
\begin{equation*}
\der_{\nrho}F_1=\nrho^{\gam-2}(\nrho-1)+\frac{\iv^2}{\rightrho-1}>0.
\end{equation*}
Then the implicit function theorem implies that $\rightrho$ is of the $C^1$--dependence on $\iv>0$.

The $C^1$--dependence
of $\righttop$ on $\iv$ follows directly from \eqref{1-23} and $\nxi=\sqrt{\rightc^2-(\neta)^2}$.

By the $C^1$--dependence of $\rightrho$ on $\iv$, we have
$$
0=\frac{\dd F_1(\nrho(\iv),\iv)}{\dd\iv}=\der_{\nrho}F_1\,\frac{\dd\nrho}{\dd\iv}-(\nrho-1)\iv.
$$
Since $\der_{\nrho}F_1>0$ is shown above, then
$\frac{\dd\nrho}{\dd\iv}>0$.
This implies that $\nrho(\iv)$ is bounded above by a finite constant for $\iv>0$ sufficiently
small so that it follows directly from \eqref{1-7} that
\begin{equation}
\label{limit-rightrho}
\lim_{\iv\to 0+}\rightrho=1.
\end{equation}

By \eqref{leftrho-expression} and \eqref{limit-rightrho}, we find that
$\displaystyle{\lim_{\iv \to 0+} \iM(\righttop)=1}$.
We combine this limit with \eqref{2-4-a6} to obtain
\begin{equation}\label{2.6.7a}
\lim_{\iv\to 0+}\neta=1.
\end{equation}
Finally, $\displaystyle{\lim_{\iv\to 0+}\nxi}=0$
is obtained from  $\nxi=\sqrt{\rightc^2-(\neta)^2}$,
and the limit of $\neta$ is given in \eqref{2.6.7a}.
The claim is verified.

\medskip
{\bf 2.}
For each $\gam>1$, there exists a small constant $\sigma>0$
so that $\nxi<1$ whenever $0<\iv\le \sigma$.
Fix $\gam > 1$. For  $0<\iv\le \sigma$, define a function $F: (0,\betac^{(\iv)})\rightarrow \R$ by
\begin{equation}
\label{1-17}
F(\beta):=\tan\otheta-\tan\itheta.
\end{equation}

{\emph{Claim{\rm :}}} {\it For any given $\gam >1$, there exists a constant $v_*\in(0,\sigma]$ so that, if $0<\iv<v_*$,
there is $\hat{\beta}^{(\iv)}\in(0,\betac^{(\iv)})$ such that
\begin{equation}
\label{1-16}
\begin{split}
&F(\beta)<0\;\;\qquad \text{for all $\beta\in(0,\hat\beta^{(\iv)})$},\\
&F(\beta)>  0\;\;\qquad \text{for all $\beta\in(\hat\beta^{(\iv)},\betac^{(\iv)})$}.
\end{split}
\end{equation}
}
Once the claim is verified, then \eqref{1-16-nversion} directly follows.

\medskip
{\bf 3.} We first show that, for each $\iv\in(0,\sigma]$, $F'(\beta)\ge 0$ holds for all $\beta\in(0, \betac^{(\iv)})$.
Fix $\iv\in(0,\sigma]$.

We use the equation of line $\iL$:
$$
(\xi_1-\nxi)\tan\itheta-(\xi_2-\neta)=0
$$
to see
\begin{equation*}
{\rm dist}(\iL,(0,-\iv))=\frac{|(\tan\itheta,-1,-\nxi\tan\itheta+\neta)\cdot(0,-\iv,1)|}{\sqrt{1+\tan^2\itheta}}=1,
\end{equation*}
and then solve it for $\tan \itheta$ to obtain
\begin{equation*}
\tan\itheta=\frac{\sqrt{((\iv+\neta)^2-1)+(\nxi)^2}-(\iv+\neta)\nxi}{1-(\nxi)^2}.
\end{equation*}

Let $(\iq, \leftq)$ be given by \eqref{definition-qs}.
By \eqref{1-52}--\eqref{1-36} in the proof of Lemma \ref{lemma:interval-existence}, we have shown that
$\oeta=-\iv+\iq\cos\beta-\sin\beta\sqrt{\vartriangle}$ with $\vartriangle:=\leftc^2-\leftq^2$.
Substituting this expression into $\oxi=\leftu-\sqrt{\leftc^2-(\oeta)^2}$ and then using \eqref{2-4-a4} and \eqref{2-m}, we have
\begin{equation*}
\oxi=-\iv\tan\beta-\big(\cos\beta\sqrt{\leftc^2-\leftq^2}+\leftq\sin\beta\big),
\end{equation*}
so that
\begin{equation*}
\tan\otheta=\frac{\neta-\oeta}{\nxi-\oxi}=
\frac{\iv-\iq\cos\beta+\sin\beta\sqrt{\vartriangle}+\neta}{\cos\beta\sqrt{\vartriangle}+\iq\sin\beta+\nxi}.
\end{equation*}

Since $\tan\itheta$ is independent of $\beta$, we have
\begin{equation*}
F'(\beta)=\frac{G(\beta)}{(\nxi+\iq\sin\beta+\cos\beta\sqrt{\vartriangle})^2},
\end{equation*}
where
\begin{align*}
G(\beta)
=&\,\big(\iq+\frac{1}{2\sqrt{\vartriangle}}\frac{\dd\vartriangle}{\dd\beta}\big)\big(\iq+\nxi\sin\beta-(\iv+\neta)\cos\beta\big)\\
 &\,+\big(\sqrt{\vartriangle}-\frac{\dd\iq}{\dd\beta}\big)\big(\sqrt{\vartriangle}+\nxi\cos\beta+(\iv+\neta)\sin\beta\big).
\end{align*}
By \eqref{definition-qs}, \eqref{density-mont-ox}, and \eqref{oM-monotonicity}, we obtain
\begin{equation*}
\frac{\dd\vartriangle}{\dd\beta}=\frac{\dd \leftc^2(1-\oM^2)}{\dd\beta}>0\qquad \tx{for all $\beta\in(0, \frac{\pi}{2})$}.
\end{equation*}
A direct computation yields that
$$
\iq+\nxi\sin\beta-(\iv+\neta)\cos\beta=(\righttop-\lefttop)\cdot {\bf n}_{\leftshock}>0
$$
for the unit normal vector ${\bf{n}}_{\leftshock}$ to $\leftshock$ pointing towards $\Oi=(0, -\iv)$ for all $\beta\in(0,\betac^{(\iv)})$.
Combining the two previous inequalities, we have
$$
G(\beta)>\big(\sqrt{\vartriangle}-\frac{\dd\iq}{\dd\beta}\big)\big(\sqrt{\vartriangle}+\nxi\cos\beta+(\iv+\neta)\sin\beta\big).
$$
Therefore, we can conclude that $F'(\beta)>0$, provided that  $\sqrt{\vartriangle}-\frac{\dd\iq}{\dd\beta}>0$ for $0<\beta<\betac^{(\iv)}$
can be proved.

A straightforward computation by using \eqref{1-1}, \eqref{diff-qi-beta}, and \eqref{diff-lq-qi} yields that
\begin{equation*}
\frac{\dd\iq}{\dd\beta}=\frac{(\iq^{\gam-1}-\leftq^{\gam+1})\tan\beta}{\iq^{\gam-2}+\leftq^{\gam}}.
\end{equation*}

Using \eqref{leftrho-expression} and \eqref{definition-qs}, we obtain that $\displaystyle{\leftc^2=\big(\frac{\iq}{\leftq}\big)^{\gam-1}}$.
Then
\begin{equation*}
\vartriangle-\big(\frac{\dd\iq}{\dd\beta}\big)^2
=\frac{\iq^{\gam-1}-\leftq^{\gam+1}}{\leftq^{\gam-1}}\Bigl(1-\underset{(=:\hat{\sigma})}
{\underbrace{\frac{ \leftq^{\gam-1}(\iq^{\gam-1}-\leftq^{\gam+1})\tan^2\beta}{(\iq^{\gam-2}+\leftq^{\gam})^2}}}\Bigr).
\end{equation*}
It can be checked directly that $\displaystyle{\frac{\dd\iq}{\dd\beta}>0}$,
by differentiating \eqref{2-4-a5} with respect to $\beta$ and applying \eqref{D}.
Then we have
\begin{equation}
\label{dqdb-ineq}
\iq^{\gam-1}-\leftq^{\gam+1}=\frac{\iq^{\gam-2}+\leftq^{\gam}}{\tan\beta}\frac{\dd\iq}{\dd\beta}>0.
\end{equation}
Since $\oeta=-\iv+\iq\cos\beta-\sin\beta\sqrt{\vartriangle}>0$ for $\beta<\betac^{(\iv)}$,
\eqref{2-4-a4} implies that $\leftq^2>\leftc^2\sin^2\beta$.
Substituting $\leftc^2=\leftrho^{\gam-1}=\bigl(\frac{\qi}{\leftq}\bigr)^{\gam-1}$ into this inequality,
we find that
$
\iq^{\gam-1}<\frac{\leftq^{\gam+1}}{\sin^2\beta},
$
which implies that
\begin{align*}
\hat{\sigma}=\frac{ (\frac{1}{K}-1)\tan^2\beta }{(\frac{\leftq}{K\iq}+1)^2}<1 \qquad\tx{for all $\beta\in(0, \betac^{(\iv)})$},
\end{align*}
where $K=\frac{\leftq^{\gam+1}}{\iq^{\gam-1}}$.
This implies that $\sqrt{\vartriangle}-\frac{\dd\iq}{\dd\beta}>0$ for $0<\beta<\betac^{(\iv)}$.

Therefore, $F'(\beta)> 0$ for all $\beta\in(0,\betac^{(\iv)})$.

\medskip
{\bf 4.}  At $\beta=0$, $\oeta=\neta$. This directly yields that $F(0)=-\tan\itheta<0$.

\medskip
{\bf 5.}
Fix $\iv\in(0, \sigma]$. At $\beta=\betac^{(\iv)}$, $\oeta=0$.
Let $\oxi^*$ denote the $\xi_1$--coordinate of point $\lefttop$ at $\beta=\betac^{(\iv)}$.
Then we have
\begin{equation*}
F(\betac^{(\iv)})=\frac{a-b}{(1-(\nxi)^2)(\nxi-\oxi^*)},
\end{equation*}
where
\begin{equation}\label{2.6.10b}
\begin{split}
&a:=\big(\nxi-\oxi^*\big)\big(\iv+\neta\big)\nxi-\neta\big((\nxi)^2-1\big),\\
&b:=\big(\nxi-\oxi^*\big)\sqrt{(\iv+\neta)^2+((\nxi)^2-1)}.
\end{split}
\end{equation}

{\emph{Claim}}: $\oxi^*$ depends continuously on $\iv\in (0, \sigma]$.

This can be seen as follows:  Fix $\beta=\betac^{(\iv)}$. Then

\smallskip
\begin{itemize}
\item[(5-1)] Since $\oeta=0$ at $\beta=\betac^{(\iv)}$, we derive from \eqref{1-52}--\eqref{1-36} that
$$
 \leftc\sqrt{1-\oM^2}\sin \betac^{(\iv)}=-\iv+\iq\cos \betac^{(\iv)}.
$$
We combine this equation with
\eqref{2-4-a4}
to yield that $\oM=\sin\betac^{(\iv)}$
and substitute this into \eqref{leftrho-expression} to obtain
\begin{equation}
\label{2-s}
\frac{\leftq^{\gam+1}}{\iq^{\gam-1}}=\sin^2\betac^{(\iv)}.
\end{equation}

\smallskip
\item[(5-2)] By \eqref{Bernoulli-RH-new} and the Rankine-Hugoniot jump condition: $\leftrho\leftq=\iq$, we have
$$
F_2(\leftrho, \iq):=\frac{\leftrho^{\gam-1}-1}{\gam-1}+\frac 12 \left(\frac{\iq}{\leftrho}\right)^2-\frac 12 \iq^2=0.
$$
The fact that $\der_{\rho}F_2(\orho, \iq)=\frac{1}{\leftrho}(\leftc^2-\leftq^2)>0$
implies that $\leftrho$ is of the $C^1$--dependence on $\qi$, so that $\leftq=\frac{\iq}{\leftrho}$
is of the $C^1$--dependence on $\qi$.

\smallskip
\item[(5-3)] It can be derived directly from \eqref{2-4-a4} and \eqref{2-s} that
\begin{equation}
\label{F3-new}
F_3(\iq,\iv):=(\iq-\leftq)^2\Big(1-\frac{\leftq^{\gam+1}}{\iq^{\gam-1}}\Big)-\iv^2=0,
\end{equation}
where $\leftq$ is regarded as a $C^1$--function of $\iq$ by (5-2).
A direct computation by using \eqref{2-4-a4}, \eqref{diff-lq-qi}, and \eqref{2-s}
shows that $\der_{\qi}F_2(\qi, \iv)\ge \frac{4\iv \cos\betac^{(\iv)}}{\gam+1}>0$.
This implies that $\qi$ is of the $C^1$--dependence on $\iv$.

\smallskip
\item[(5-4)] $\oxi^*$ is the $\xi_1$--intercept of $\leftshock$ so that
$\oxi^*=-\iv\tan\betac^{(\iv)}-\leftq \csc\betac^{(\iv)}$.
By the $C^1$--dependence of $\betac^{(\iv)}$ and $\leftq$ on $\iv$,
we conclude that $\oxi^*$ is of the $C^1$--dependence on $\iv$.
The claim is verified.
\end{itemize}

\smallskip
{\bf 6.}  {\emph{Claim}}: For $a$ and $b$ defined in \eqref{2.6.10b},
$\displaystyle{\lim_{\iv\to 0+}(a^2-b^2)=1}$.

It suffices to show that $\sup_{\iv\in(0,\sigma]}|\oxi^*|$ is bounded, due to \eqref{limits-1}.
From \eqref{F3-new}, we have two cases:
$\displaystyle{\lim_{\iv\to 0+}\frac{\leftq}{\qi}= 1}$ and
$\displaystyle{\lim_{\iv\to 0+}\frac{\leftq^{\gam+1}}{\iq^{\gam-1}}=1}$.

\smallskip
For the case that $\displaystyle{\lim_{\iv\to 0+}\frac{\leftq}{\iq}= 1}$,
\eqref{2-s} implies that $\displaystyle{\sup_{(0,\sigma]}\iq}$ is finite.
Then it follows from
$\qi=D\ivphi(\lefttop)\cdot{\bf n}_{\leftshock}=-\oxi^*\sin\betac^{(\iv)}+\iv\cos\betac^{(\iv)}$
that $\sup_{\iv\in(0,\sigma]}|\oxi^*\sin\betac^{(\iv)}|$ is finite.
We multiply \eqref{2-s} by $(\oxi^*)^2$ to obtain
\begin{equation*}
\sup_{\iv\in(0,\sigma]}\big(\oxi^*\big)^2\le \sup_{\iv\in(0,\sigma]}\big(\oxi^*\sin\betac^{(\iv)}\big)^2\frac{\iq^{\gam+1}}{\leftq^{\gam+1}}<\infty,
\end{equation*}
where we have used the fact that $\qi>1$ for each $\iv>0$.

For the case that $\displaystyle{\lim_{\iv\to 0+}\frac{\leftq^{\gam+1}}{\iq^{\gam-1}}=1}$,
we substitute $\orho=\frac{\iq}{\leftq}$ into $F_2(\leftrho, \qi)=0$ to obtain
\begin{equation*}
\frac{1}{2}\left(\frac{\iq}{\leftq}\right)^2
=\frac{1}{\gam-1}\left(\frac{\iq^{\gam-1}}{\leftq^{\gam+1}}-\frac{1}{\leftq^2}\right)
+\frac{1}{2}
\le \frac{1}{\gam-1}\frac{\iq^{\gam-1}}{\leftq^{\gam+1}}+\frac{1}{2}.
\end{equation*}
From this, it
follows that $\sup_{\iv\in(0, \sigma]}|\frac{\iq}{\leftq}|$ is finite. Then we use
\eqref{2-s} to see that $\sup_{\iv\in(0, \sigma]}\iq$ is finite.
Finally, we repeat the argument for the case that
$\displaystyle{\lim_{\iv\to 0+}\frac{\leftq}{\iq}= 1}$
to conclude that $\displaystyle{\sup_{\iv\in(0,\sigma]}\big(\oxi^*\big)^2}$ is finite,
which implies the claim.

\medskip
{\bf 7.}  By the result obtained from Step {\rm 6},
there exists a constant $v_*\in(0, \sigma]$ such that $F(0)<0<F(\betac^{(\iv)})$ for all $\iv\in (0,v_*]$.
Finally, the monotonicity of $F(\beta)$, proved in Step {\rm 3},
yields Proposition \ref{proposition:el-condition}.
\end{proof}
\end{proposition}

When \eqref{1-16-nversion} holds, the existence of a solution of Problem \ref{fbp} has been proved in \cite{EL}.
This implies the global existence of a weak solution of Problem \ref{problem-2}
with the structure of Fig. \ref{fig:globala}, provided that \eqref{1-16-nversion} holds.
In this monograph, we establish the global existence of admissible solutions for all $(\iv, \beta)\in \mathfrak{R}_{\rm weak}$
({\it i.e.}, the global existence of weak solutions to {\it Problem} \ref{problem-2} for all $(\iu, u_0)\in \mathfrak{P}_{\rm weak}$),
which includes the case that \eqref{1-16nv-ox} holds, or the case that $\beta\ge \betac^{(\iv)}$.




\chapter{Uniform Estimates of Admissible Solutions }\label{section-unif-est-1}

As in \cite{CF2}, we employ the Leray-Schauder degree to prove Theorem \ref{theorem-0}.
In order to construct an iteration set (as a subset of a properly defined Banach space)
and an iteration map,
we first establish uniform estimates of admissible solutions corresponding to $(\iv, \beta)$
with respect to $\beta\in[0, \betadet-\eps]$ in the sense of Definitions \ref{def-regular-sol} and \ref{def-regular-sol-normal} for each $\iv>0$
and small $\eps>0$. In particular, it is crucial to establish the uniform estimates of the size of pseudo-subsonic
region $\Om$ and the pseudo-potential function $\vphi$ restricted
to $\Om$ in properly chosen norms.
Following the approach of \cite{CF2},
we establish various uniform estimates of admissible solutions in the following order:

\smallskip
\begin{itemize}
\item Strict directional monotonicity properties of $\ivphi-\vphi$,

\smallskip
\item Strict directional monotonicity properties of $\vphi-\rightvphi$ and $\vphi-\leftvphi$,

\smallskip
\item Uniform positive lower bound of the distance between $\shock$ and $\Wedge$ away from $P_{\beta}$,

\smallskip
\item Uniform positive lower bound of ${\rm dist}(\shock, \der B_1(\Oi))$,

\smallskip
\item Uniform estimates of the ellipticity of Eq. \eqref{2-1} in $\Om$,

\smallskip
\item Uniform weighted $C^{2,\alp}$ estimates of admissible solutions in $\Om$.
\end{itemize}

\smallskip
Fix $\gam\ge 1$ and $\iv>0$. For each $\beta\in[0,\frac{\pi}{2})$,
let $(\ivphi, \rightvphi, \leftvphi)$ and $(\Oi, \Oo, \Onormal)$ be defined by Definition \ref{definition-domains-np}.
We also follow Definition \ref{definition-domains-np} for the notations of
$(\rightsonic, \leftsonic)$ and $(\lefttop, \righttop, \rightbottom,\leftbottom)$.

Note that the definitions of $(\leftsonic, \lefttop, \leftbottom)$ are different
for the respective cases $\beta\in [0, \betac^{(\iv)})$ and $\beta\in [\betac^{(\iv)}, \frac{\pi}{2})$,
but they depend continuously on $\beta\in(0, \frac{\pi}{2})$.

\numberwithin{equation}{section}
\section{Directional Monotonicity Properties of Admissible Solutions}
\label{subsection-directional-mono}
In this section,
we establish directional monotonicity properties of
$\ivphi-\vphi$, $\vphi-\rightvphi$, and $\vphi-\leftvphi$ for admissible solutions $\vphi$ in the sense of Definition \ref{def-regular-sol}.

\subsection{Strict directional monotonicity of $\ivphi-\vphi$}\label{subsec-step1}

For an admissible solution $\vphi$ in the sense of Definition \ref{def-regular-sol} for $(\iv,\beta)\in \mathfrak{R}_{\rm weak}$,
define
\begin{equation}
\label{definition-of-phi}
\phi:=\vphi-\rightvphi\qquad\text{in $\Omega$}.
\end{equation}
Then $\phi$ satisfies the equation:
\begin{equation}
\label{2-3}
(c^2-\vphi_{\xin}^2)\phi_{\xin\xin}-2\vphi_{\xin}\vphi_{\etan}\phi_{\xin\etan}
+(c^2-\vphi_{\etan}^2)\phi_{\etan\etan}=0
\end{equation}
in the pseudo-subsonic region $\Om$ for
$c^2=c^2(|D\vphi|^2, \vphi, \bmxi)$ given by
\begin{equation}
\label{2-5}
c^2(|\bm p|^2,z,\bmxi):=\rho^{\gam-1}(|\bm p|^2, z, \bmxi),
\end{equation}
where $\rho(|\bm p|^2, z, \bmxi)$ is defined by \eqref{new-density}.

\begin{lemma}
\label{lemma-7-0}
Fix $\gam\ge 1$ and $\iv>0$.
Let $\vphi$ be an admissible solution in the sense of Definition {\rm \ref{def-regular-sol}}
for $(\iv,\beta)\in\mathfrak{R}_{\rm weak}$ with $\beta>0$,
and let $\phi$ be given by \eqref{definition-of-phi}.
Then, for any given unit vector $\ee\in \R^2$, $\der_{\ee}\phi$ is not a constant in $\Om$.
\end{lemma}

\begin{proof}
By condition (ii) of Definition \ref{def-regular-sol}, $\phi$ satisfies
\begin{align}
&\der_{\ee}\phi=0 &&\text{on $\rightsonic$}, \label{6-a1}\\
&\der_{\ee}\phi=\der_{\ee}(\leftvphi-\rightvphi)={\ee}\cdot (\leftu,0)&&\text{on $\leftsonic$}\label{6-a1-b}
\end{align}
for each unit vector $\ee$ in $\R^2$.

Suppose that $\der_{\ee}\phi$ is a constant in $\Om$.
Then \eqref{6-a1}--\eqref{6-a1-b} imply that $\ee$ must be parallel to ${\ee}_{2}=(0,1)$,
because $\leftu \neq 0$, by Definition \ref{definition-domains-np}.
Then
$\der_{\etan}\phi\equiv 0$ in $\Om$, which implies that $\der_{\xin\etan}\phi=\der_{\etan\etan}\phi\equiv 0$ in $\Om$.
Since Eq. \eqref{2-3} is strictly elliptic in $\Om$, it follows that $\der_{\xin\xin}\phi\equiv 0$ in $\Om$.
Thus, there exist constants $(u,v,k)$ such that $\phi(\xin,\etan)=u\xin+v\etan+k$ in $\Om$.
Since the length of $\rightsonic$ is nonzero, we obtain from the boundary condition $\phi\equiv 0$ on $\rightsonic$
that $D\phi\equiv 0$ in $\Om$, so that $\phi\equiv 0$ in $\Om$.
However, this contradicts the boundary condition \eqref{3-a1} on $\leftsonic$,
because $\phi=\leftvphi-\rightvphi=\leftu\xin-\iv \etan^{(\beta)}+\iv\neta$ on $\leftsonic$, by Remark \ref{remark-vphibeta}.
\end{proof}

\begin{lemma}
\label{lemma-step1-1}
Fix $\gam \ge 1$ and $\iv>0$.
Let $\vphi$ be an admissible solution in the sense of Definition {\rm \ref{def-regular-sol}}
for $(\iv,\beta)\in\mathfrak{R}_{\rm weak}$ with $\beta>0$.
For vectors $\leftvec$ and $\rightvec$ given by Definition {\rm \ref{definition-domains-np}}, $\vphi$ satisfies
\begin{align}
\label{3-c2}
&\der_{\leftvec}(\ivphi-\vphi)< 0\;\;\qquad \text{in $\ol{\Om}\setminus\ol{\leftsonic}$},\\
\label{3-c3}
&\der_{\rightvec}(\ivphi-\vphi)<0\;\;\qquad \text{in $\ol{\Om}\setminus \ol{\rightsonic}$}.
\end{align}

\begin{proof}
By Definition \ref{def-regular-sol}(v), any admissible solution $\vphi$
satisfies that $\der_{\rightvec}(\ivphi-\vphi)\le 0$ and $\der_{\leftvec}(\ivphi-\vphi)\le 0$ in $\Om$.
Therefore, it suffices to prove the strict inequalities.

For ${\ee}={\leftvec}$ or $\rightvec$, we introduce a coordinate system $(S,T)$
so that $\ee=(1,0)$ and $\ee^{\perp}=(0,1)$ in the $(S,T)$--coordinates.
We note that Eq. \eqref{2-1} is invariant under a coordinate rotation.
Also, $D^2(\ivphi-\vphi)=-D^2\phi$ for $\phi$ given by \eqref{definition-of-phi}.
Then $\ivphi-\vphi$ satisfies
\begin{equation} \label{3-2}
 (c^2-\vphi_S^2)(\ivphi-\vphi)_{SS}
 -2\vphi_S\vphi_T(\ivphi-\vphi)_{ST}+(c^2-\vphi_T^2)(\ivphi-\vphi)_{TT}=0\qquad \tx{in $\Om$}.
\end{equation}
Denote $v:=\der_S(\ivphi-\vphi)$. Then $v$ satisfies the following properties:

\smallskip
\emph{{{\rm (i)} $v<0$ in $\Om$}}.
We differentiate \eqref{3-2} with respect to $S$ and use the expression:
$$
(\ivphi-\vphi)_{TT}=
-\frac{ (c^2-\vphi_S^2)(\ivphi-\vphi)_{SS}
 -2\vphi_S\vphi_T(\ivphi-\vphi)_{ST}}{c^2-\vphi_T^2}
$$
to obtain the following equation
for $v$:
\begin{align}\label{diff-nondivMainEq}
&(c^2-\varphi_S^2)v_{SS}-2 \varphi_S\varphi_T v_{ST}
+(c^2-\varphi_T^2)v_{TT}
\\[1mm]
&\,
-\big((\gamma-1)\varphi_S(\phi_{SS}+\phi_{TT})+(\gamma-1)\varphi_T\phi_{ST}+
2\varphi_S(\phi_{SS}-1)\big) v_{S} \nonumber
\\[1mm]
&\,
-\big(2\varphi_T(\phi_{SS}-1)+2\phi_{ST}\varphi_S+(\gamma+1)\varphi_T\phi_{TT}\big) v_{T}
=0.\nonumber
\end{align}
We consider \eqref{diff-nondivMainEq} as a linear second-order equation for $v$.
Then this equation does not have the zero-order terms, and its coefficients are continuous in $\Omega$.

Since Eq. \eqref{3-2} is strictly elliptic in $\ol{\Om}\setminus(\ol{\leftsonic}\cup \ol{\rightsonic})$
by  Definition \ref{def-regular-sol}(iii),
the equation for $v$ is strictly elliptic in $\ol{\Om}\setminus(\ol{\leftsonic}\cup \ol{\rightsonic})$,
because the coefficients of the principal part of the equation for $v$ are the same as
those in Eq. \eqref{3-2}.
Moreover, $v$ is not a constant in $\Om$ by Lemma \ref{lemma-7-0}, so $v$ cannot attain its maximum in
$\Om$ by the strong maximum principle. Thus, $v<0$ holds in $\Om$.

\smallskip
\emph{ {{\rm (ii)} $v<0$ on $\Wedge$}}.
On $\Wedge$, the slip boundary condition \eqref{3-a2} for $\vphi$
implies that $\der_{\etan}(\ivphi-\vphi)=-\iv$, so that
$\der_{\xin\etan}(\ivphi-\vphi)=0$.
In Eq. \eqref{3-2}, we replace $(S,T)$ by $(\xin,\etan)$ to obtain
\begin{equation}
\label{eqn-on-wedge}
(c^2-\vphi_{\xin}^2)\der_{\xin\xin}(\ivphi-\vphi)+(c^2-\vphi_{\etan}^2)\der_{\etan\etan}(\ivphi-\vphi)=0
\qquad\,\tx{on $\Wedge$}.
\end{equation}
Let $\{{\bf e}_{\xin}, {\bf e}_{\etan}\}$ form an orthonormal basis for coordinates $\xxi=(\xi_1,\xi_2)$.
By setting $a_1:=\ee\cdot{\bf e}_{\xin}$ and $a_2:=\ee\cdot{\bf e}_{\etan}$,
$v$ is expressed as $v=a_1\der_{\xin}(\ivphi-\vphi)+a_2\der_{\etan}(\ivphi-\vphi)$
so that $v_{\xin}=a_1\der_{\xin\xin}(\ivphi-\vphi)$ and $v_{\etan}=a_2\der_{\etan\etan}(\ivphi-\vphi)$ on $\Wedge$.

Substituting these expressions into \eqref{eqn-on-wedge},
we obtain the following boundary condition for $v$:
 \begin{equation}
\label{3-4}
\der_{\etan}v+\frac{a_2(c^2-\vphi_{\xin}^2)}{a_1(c^2-\vphi_{\etan}^2)}\der_{\xin}v=0\qquad \text{on $\Wedge$}.
\end{equation}

Since ${\ee}\cdot {\bf e}_{\xin}\neq 0$, {\it i.e.}, $a_1\neq 0$, \eqref{3-4} is an oblique boundary condition for $v$ on $\Wedge$.
Thus, Hopf's lemma applies. Therefore, $v$ cannot attain its maximum on $\Wedge$, which implies that $v<0$ on $\Wedge$.

\smallskip
\emph{ {{\rm (iii)} $v<0$ on $\shock$}}.
Suppose that $v(\hat P)=0$ for some $\hat P\in \shock$.
Let ${\bf n}_{\rm{sh}}$ be the unit normal vector to $\shock$ towards the interior of $\Om$,
and let $\bm{\tau}_{\rm{sh}}$ be the unit tangent vector to $\shock$ with ${\bm\tau}_{\rm{sh}}\cdot{\rightvec}<0$.
Differentiating the Rankine-Hugoniot jump condition:
$\big[\rho(|D\vphi|^2, \vphi)D\vphi\cdot{\bf n}_{\rm sh}\big]_{\shock}=0$ in the direction of $\bm{\tau}_{\rm{sh}}$, we have
\begin{equation}
\label{5-c1}
D^2(\ivphi-\vphi)[\bm{\tau}_{\rm{sh}}, \bm h]:=
\bm{\tau}_{\rm{sh}}\cdot  D^2(\ivphi-\vphi)\bm h=0\;\;\qquad \text{on $\shock$},
\end{equation}
where $\bm h=h_{\rm n}{\bf n}_{\rm{sh}}+h_{\rm t}\bm{\tau}_{\rm{sh}}$
with
\begin{equation}\label{3-13}
h_{\rm n}=-\rho\vphi_{{\bf n}_{\rm{sh}}}(c^2-\vphi_{{\bf n}_{\rm{sh}}}^2),\qquad
h_{\rm t}=( c^2+\rho \vphi_{{\bf n}_{\rm{sh}}}^2)\vphi_{\bm\tau_{\rm{sh}}}.
\end{equation}
We refer to Lemma 5.11 in
\cite{CF2} for the verification of \eqref{5-c1}--\eqref{3-13}.

 It follows from Lemma \ref{lemma-entropycond-admsbsol}(a)
 and the ellipticity of \eqref{2-1} in $\overline{\Omega}\setminus (\overline{\leftsonic}\cup \overline{\rightsonic})$
that
\begin{equation}\label{negativity-of-hn}
  h_{\rm n}<0\qquad \tx{ on $\shock$.}
\end{equation}

Since it is assumed that $v=\der_{{\ee}}(\ivphi-\vphi)$ has a local extremum at $\hat{P}\in \shock$, we have
\begin{equation}
\label{5-c2}
D^2(\ivphi-\vphi)[\bm{\tau}_{\rm{sh}}, \ee]=0\;\;\qquad\text{at $\hat{P}$}.
 \end{equation}
We express ${\ee}=b_1{\bf n}_{\rm{sh}}+b_2\bm \tau_{\rm{sh}}$.
Then we rewrite \eqref{5-c1} restricted at $\hat{P}$ and  \eqref{5-c2} as a linear system for $(\ivphi-\vphi)_{{\bm\tau}_{\rm{sh}}{{\bf n}_{\rm{sh}}}}(\hat P)$
and $(\ivphi-\vphi)_{{\bm\tau}_{\rm{sh}}{\bm\tau}_{\rm{sh}}}(\hat P)$.
By this linear system and \eqref{3-2}, we find that
$D^2(\ivphi-\vphi)(\hat P)=\bm 0$, unless
\begin{equation}\label{3.1.13b}
\det \begin{pmatrix}
h_{\rm n}&h_{\rm t}\\
b_1&b_2
\end{pmatrix}=0 \qquad\tx{at $\hat{P}$.}
\end{equation}
On the other hand, $v$ is not a constant in $\Om$ by Lemma \ref{lemma-7-0},
which implies that $D^2(\ivphi-\vphi)[{\bf n}_{\rm{sh}}, \ee]=v_{{\bf n}_{\rm{sh}}}>0$ at $\hat P$ by Hopf's lemma,
so that $D^2(\ivphi-\vphi)(\hat{P})=\bm 0$ is impossible.
%
%

Therefore, \eqref{3.1.13b} must hold,
so that
${\ee}=k{\bm h}$ at $\hat{P}$ for some constant $k\neq 0$.
This yields that
$$
|v(\hat P)|=|k h_{\rm n}(\hat P) \, D(\ivphi-\vphi)(\hat P) |>0.
$$
This contradicts the fact that $v(\hat{P})=0$.
Therefore,  we conclude that $v<0$ on $\shock$.

\smallskip
\emph{ {{\rm (iv)} $v<0$ on the sonic arcs}}.
If $\ee=\leftvec$, then
$
v=\der_{\leftvec}(\ivphi-\rightvphi)=\frac{(\iv,-\leftu)\cdot(0,-\iv)}{\sqrt{\leftu^2+\iv^2}}<0
$ on ${\rightsonic}$. This proves \eqref{3-c2}.

If $\ee=\rightvec$, then
$v=\der_{\rightvec}(\ivphi-\leftvphi)={-(\leftu,\iv)\cdot(-1,0)}<0$ on ${\leftsonic}$. This proves \eqref{3-c3}.
This computation holds even for the case that $\leftsonic=\{P_{\beta}\}$, {\it i.e.},
$\beta\ge \betasonic$ by the condition stated in (ii-3) for Case II.
\end{proof}
\end{lemma}

Define the following set:
\begin{equation}
\label{def-cone-cor}
{\rm{Cone}}(\leftvec, \rightvec):=\{\alp_1\leftvec+\alp_2\rightvec\,:\,\alp_1, \alp_2\ge 0\},
\end{equation}
and let $\cone$ be the interior of ${\rm{Cone}}(\leftvec, \rightvec)$.
By Lemma \ref{lemma-step1-1}, if $\vphi$ is an admissible solution corresponding to $(\iv, \beta)$,
then $\vphi$ satisfies
\begin{equation}
\label{cor-dir-mont}
\der_{\ee}(\ivphi-\vphi)<0\qquad \tx{in $\ol{\Om}\,$ for all $\ee\in \cone$}.
\end{equation}

\begin{remark}
\label{remark-unitnormalshock-cone}
By \eqref{def-monot-vecs},
$\cone$ can be represented as
\begin{equation*}
\cone=\{(r \cos\theta, r\sin\theta)\,:\, r>0,\,\, \beta<\theta<\pi\}.
\end{equation*}
Note that the unit normal vector ${\bf n}_{\rm{sh}}$ to $\shock$ is expressed
as $\displaystyle{{\bf n}_{\rm{sh}}}=\frac{D(\ivphi-\vphi)}{|D(\ivphi-\vphi)|}$.
It follows from \eqref{3-c2} that
$-{\bf n}_{\rm{sh}}(P)\in\{(\cos\theta,\sin\theta)\,:\, \beta-\frac{\pi}{2}<\theta<\beta+\frac{\pi}{2}\}$
for all $P\in \shock$. Moreover, it follows from \eqref{3-c3}
that
$-{\bf n}_{\rm{sh}}(P)\in\{(\cos\theta, \sin\theta)\,:\, \frac{\pi}{2}<\theta<\frac{3\pi}{2}\}$
%
for all $P\in \shock$.
Therefore, we have
\begin{equation}
\label{shocknormal-in-cone}
  -{\bf n}_{\rm{sh}}(P)\in \{(\cos\theta,\sin\theta)\,:\,\frac{\pi}{2}<\theta<\beta+\frac{\pi}{2}\}\subset
  \cone\qquad\tx{for all $P\in \shock$},
\end{equation}
since $\beta\in (0,\betadet)\subset (0, \frac{\pi}{2})$.
\end{remark}

\begin{proposition}
\label{proposition-3}
Given $\gam\ge 1$ and $\iv>0$, let $\vphi$ be an admissible solution in the sense of
Definition {\rm \ref{def-regular-sol}} for $(\iv,\beta)\in\mathfrak{R}_{\rm weak}$.
Then there exists a function $\etan=\fshock(\xin)$ such that

\smallskip
\begin{itemize}
\item[(i)] $\shock=\{\bmxi\,:\,\etan=\fshock(\xin),\;\; \xi_1^{\lefttop}<\xin<\xi_1^{\righttop}\}$,
where $\xi_1^{P_j}$ is the $\xin$--coordinate of point $P_j$ for $j=1,2${\rm ;}

\smallskip
\item[(ii)]
$\fshock$ satisfies
\begin{equation}
\label{3-c5a}
0=\fshock'(\xin^{\righttop})<\fshock'(\xin)<\fshock'(\xin^{\lefttop})=\tan\beta
\;\;\qquad \text{for $\xi_1^{\lefttop}<\xin<\xi_1^{\righttop}$}.
\end{equation}
\end{itemize}

\smallskip
\begin{proof}
Note that ${\ee}_{\etan}\in \cone$.
By \eqref{cor-dir-mont}, we have
\begin{equation}
\label{3-c4}
\der_{\etan}(\ivphi-\vphi)<0\;\;\qquad \text{on $\ol{\shock}$}.
\end{equation}
This, combined with Definition \ref{def-regular-sol}(i),
implies that there exists a unique $C^1$--function $\fshock$ satisfying statement (i) above.

Since $\ivphi-\vphi=0$ holds on $\shock$, $\fshock$ satisfies that $(\ivphi-\vphi)(\xi_1, \fshock(\xi_1))=0$
for $\xi_1^{\lefttop}<\xi_1<\xi_1^{\righttop}$. We differentiate this expression with respect to $\xi_1$ to obtain
\begin{equation*}
\fshock'(\xin)=-\frac{\der_{\xin}(\ivphi-\vphi)(\xin,\fshock(\xin))}{\der_{\etan}(\ivphi-\vphi)(\xin,\fshock(\xin))}.
\end{equation*}
By condition (i-3) of Definition \ref{def-regular-sol}, we have
\begin{equation}\label{shock-endpoints}
  \fshock'(\xi_1^{\lefttop})=\tan \beta,\qquad \fshock'(\xi_1^{\righttop})=0.
\end{equation}
By conditions (ii-3) and (iv) of Definition \ref{def-regular-sol},
the unit normal vector ${\bf n}_{\rm{sh}}$ to $\shock$ towards the interior of $\Om$ can be expressed as
\begin{equation*}
{\bf n}_{\rm{sh}}(P)=\frac{D(\ivphi-\vphi)(P)}{|D(\ivphi-\vphi)(P)|}
=\frac{(\fshock'(\xin),-1)}{\sqrt{1+(\fshock'(\xin))^2}}\qquad \tx{at $P=(\xin, \fshock(\xin))$}.
\end{equation*}
By Lemma \ref{lemma-step1-1} and the definition of $(\leftvec, \rightvec)$ given in Definition \ref{definition-domains-np}, we have
\begin{equation}
\label{shocknormal-vec-mult}
\begin{split}
&a_1\cos \beta(-\fshock'(\xin)+\tan\beta)-a_2\fshock'(\xin)\\
&=\sqrt{1+(\fshock'(\xin))^2}{\bf n}_{\rm{sh}}(P)\cdot (a_1\leftvec+a_2\rightvec)\\
&=\sqrt{1+(\fshock'(\xin))^2}\frac{D(\ivphi-\vphi)(P)\cdot (a_1\leftvec+a_2\rightvec)
}{|D(\ivphi-\vphi)(P)|}<0\qquad\tx{for $\xin^{\lefttop}<\xin<\xin^{\righttop}$}
\end{split}
\end{equation}
for any constants $a_1\ge 0$ and $a_2\ge 0$ with $a_1+a_2>0$.

If we choose $(a_1,a_2)=(1,0)$, then \eqref{shocknormal-vec-mult} yields
\begin{equation*}
 \fshock'(\xin)<\tan\beta\qquad \tx{for $\xin^{\lefttop}<\xin<\xin^{\righttop}$}.
 \end{equation*}
Choosing $(a_1, a_2)=(0,1)$, then we have
\begin{equation*}
\fshock'(\xin)>0\qquad \tx{for $\xin^{\lefttop}<\xin<\xin^{\righttop}$}.
\end{equation*}
Finally, \eqref{3-c5a} is obtained by combining the previous two inequalities with \eqref{shock-endpoints}.
\end{proof}
\end{proposition}

Given $\gam\ge 1$ and $\iv>0$, if $\beta_*\in\bigl(0,\betac^{(\iv)}\bigr)$ is fixed,
then Proposition \ref{proposition-3} directly implies
that
\begin{equation}
 \label{new-dec-1}
 \inf_{\beta\in(0,\beta_*]}{\rm dist}\{\shock, \Wedge\}\ge \inf_{(0,\beta_*]}\xi_2^{\lefttop}>0.
 \end{equation}

\begin{lemma}
\label{lemma-step3-1}
Fix $\gam\ge 1$ and $\iv>0$.
Let $\vphi$ be an admissible solution corresponding to $(\iv, \beta)\in \mathfrak{R}_{\rm weak}$
in the sense of Definition {\rm \ref{def-regular-sol}},
and let $\Om$ be its pseudo-subsonic region.
Then there exists a constant $C>0$ depending only on $(\iv, \gam)$ such that the following properties hold{\rm :}
\begin{align}
\label{3-c5}
& \Om\subset B_C(\bm 0),\\
\label{3-c6}
&\max_{\ol{\Om}}|\vphi|\le C,\quad \|\vphi\|_{C^{0,1}(\ol{\Om})}\le C,\\
\label{3-c7}
&\rho^*(\gam)\le \rho \le C\;\;\text{in $\Om$},\qquad 1 < \rho\le C\;\;\text{on $\shock$},
\end{align}
where
\begin{equation*}
\rho^*(\gam)=
\begin{cases}
\big(\frac{2}{\gam+1}\big)^{\frac{1}{\gam-1}}\qquad\mbox{for $\gam>1$},\\
e^{-\frac 12}=\lim_{\gam \to 1+}\big(\frac{2}{\gam+1}\big)^{\frac{1}{\gam-1}}
\qquad\mbox{for $\gam=1$}.
\end{cases}
\end{equation*}

\begin{proof}
To prove this lemma, we follow the ideas in the proofs for \cite[Proposition 9.1.2, Corollary 9.1.3, Lemma 9.1.4]{CF2}.

\medskip
{\textbf{1.}} {\emph{Proof of \eqref{3-c5}}.}
For an admissible solution $\vphi$, let $\fshock$ be as in Proposition \ref{proposition-3}.
From \eqref{3-c5a}, it follows that
$
0 \le  \xi_2^{{\lefttop}}\le \fshock(\xin)\le \xi_2^{\righttop}\;\;\text{on $[\xi_1^{\lefttop}, \xi_1^{\righttop}]$}.
$
Then
\begin{equation*}
\Om\subset\{\bmxi=(\xin,\etan)\,:\, \leftu-\leftc<\xin<\rightc, 0<\etan<\xi_2^{\righttop}\}.
\end{equation*}
For any given $\iv>0$, $\leftc$ and $\leftu$ depend continuously on $\beta\in[0, \frac{\pi}{2})$,
and $\beta_{\rm d}^{(\iv)}$ depends continuously on $\iv>0$.
Therefore, there exists a constant $C_1>0$ depending only on $(\iv, \gam)$ such that
$$
\displaystyle{\sup_{\beta\in[0,\beta_{\rm d}^{(\iv)}]}\big(|\leftu|+|\leftc|\big)\le C_1.}
$$
This proves \eqref{3-c5}.

\medskip
{\textbf{2.}} {\emph{Proof of \eqref{3-c6}}.}
By Definition \ref{def-regular-sol}(iv), we have
\begin{equation*}
\inf_{\Om} \max\{\leftvphi, \rightvphi\}\le \vphi \le
\sup_{\Om}\ivphi.
\end{equation*}
By \eqref{3-c5} and the definition of $(\ivphi, \leftvphi, \rightvphi)$ given in Definition \ref{definition-domains-np},
there exists a constant $C_2>0$ depending only on $(\iv, \gam)$ such that
$\displaystyle{
-C_2\le \min_{\ol{\Om}}\max\{\leftvphi, \rightvphi\}<\max_{\ol{\Om}} \ivphi \le C_2.}
$
Then condition (iv) of Definition \ref{def-regular-sol} implies that
\begin{equation}
\label{estimate-C0-C2}
\max_{\ol{\Om}}|\vphi|\le C_2.
\end{equation}

By conditions (ii)--(iii) of Definition \ref{def-regular-sol}, \eqref{1-f},
and \eqref{estimate-C0-C2},  we can choose a constant $\hat{C}_2>0$ depending only
on $(\iv, \gam)$ such that
$\displaystyle{\max_{\ol{\Om}}|D\vphi|\le \hat{C}_2}$ holds for each admissible solution corresponding
to $(\iv, \beta)\in \mathfrak{R}_{\rm weak}$. This, combined with \eqref{estimate-C0-C2}, yields \eqref{3-c6}.

\medskip
{\textbf{3.}} {\emph{Proof of \eqref{3-c7}.}}
A uniform upper bound of $\rho$ in \eqref{3-c7} is obtained directly from \eqref{3-c6} and \eqref{new-density}.

By condition (iii) of Definition \ref{def-regular-sol}, any admissible solution $\vphi$ satisfies
\begin{equation*}
h(\rho)+\frac{c^2}{2}\ge h(\rho)+\frac 12|D\vphi|^2\qquad \tx{in $\ol{\Om}$.}
\end{equation*}
Moreover, by  \eqref{1-p1} and condition (iv) of Definition \ref{def-regular-sol},
\begin{equation*}
h(\rho)+\frac 12 |D\vphi|^2\ge \underset{(=0)}{\underbrace{h(1)}}+\frac 12|D\ivphi|^2\ge 0\qquad\tx{in $\ol{\Om}$.}
\end{equation*}
Then we have
\begin{equation*}
h(\rho)+\frac{c^2}{2}\ge 0\qquad \tx{in $\ol{\Om}$},
\end{equation*}
so that the first inequality in \eqref{3-c7} is proved.

By Definition \ref{def-shocksol-self-similar-new2015} and Definition \ref{def-regular-sol}(iv),
any admissible solution satisfies that $\der_{{\bm \nu}}\ivphi> \der_{{\bm \nu}}\vphi$ on $\shock$ for the unit normal vector
${{\bm \nu}}$ to $\shock$ towards the interior of $\Om$.
Then the Rankine-Hugoniot jump condition stated in Definition \ref{def-regular-sol}(ii-4)
implies that $\rho>1$ holds on $\shock$, because $\irho=1$ is the density of the incoming state
corresponding to $\ivphi$. This verifies the second inequality in \eqref{3-c7}.
\end{proof}
\end{lemma}

\subsection{Directional monotonicity of $\vphi-\rightvphi$ and $\vphi-\leftvphi$}
\label{subsec-dir-mon-2-super}

Let $\vphi$ be an admissible solution, and let $\bm\nu$ be the unit normal vector to $\shock$ towards the interior of $\Om$.
For each point $P\in \shock$, define
\begin{equation*}
d(P):=\der_{\bm \nu}\ivphi(P),\quad \omega(P):=\der_{\bm \nu}(\ivphi-\vphi)(P)
\end{equation*}
so that
$$
\der_{\bm \nu}\vphi(P)=d(P)-\om(P).
$$
By Lemma \ref{lemma-entropycond-admsbsol}, $d(P)>1$ and $\omega(P)<d(P)$ on $\shock$.
By the Rankine-Hugoniot conditions stated in  Definition \ref{def-regular-sol}(ii-4),
$\rho(|D\vphi|^2, \vphi)=\frac{d}{d-\om}$ on $\shock$. Then it can be derived from \eqref{new-density} and $\ivphi-\vphi=0$ on $\shock$ that
\begin{equation*}
G(\omega, d):=h(\frac{d}{d-\om})+\frac{1}{2}\left((d-\omega)^2-d^2\right)=0\qquad \tx{on}\,\,\shock,
\end{equation*}
where $h(\rho)$ is defined by \eqref{def-h-and-c}.
For a fixed constant $d>0$, it is direct to see that
\begin{equation*}
\begin{split}
&G(0,d)=0,\quad \lim_{\om\to d-}G(\om, d)=\infty,\\
&G_{\om}(\om,d)=\frac{d^{\gam-1}}{(d-\om)^{\gam}}-(d-\om)\begin{cases}
\le 0\qquad\mbox{for $0\le \om \le d(1-d^{-\frac{2}{\gam+1}})$},\\
>0 \qquad\mbox{for $\om>  d(1-d^{-\frac{2}{\gam+1}})$}.
\end{cases}
\end{split}
\end{equation*}
Therefore, for each $d>0$, there exists a unique $\om_d\in(0,d)$ satisfying that $G(\om_d,d)=0$.
Define a function $H: (1,\infty)\rightarrow \R^+$ by
\begin{equation}
\label{def-H-ox}
H(d):=\om_d.
\end{equation}

By continuation, $H$ can be defined up to $d=1$ with $\displaystyle{H(1)=\lim_{d\to 1+} H(d)=0}$.
It is shown in {\cite[Lemma 6.1.3]{CF2}} that
\begin{equation}
\label{prelim2-1}
H\in C([1, \infty))\cap C^{\infty}((1,\infty)), \qquad \,
H'(d)>0 \,\,\, \text{for all $d\in (1,\infty)$}.
\end{equation}
Therefore, we have
\begin{equation}
\label{H-property-new}
H(1)=0, \qquad \,   H(d)>0\,\,\,\tx{if and only if\,\, $d>1$}.
\end{equation}

For each $P\in \shock$, we have
\begin{equation}
\label{H-relation-ox}
\der_{\bm \nu}(\ivphi-\vphi)(P)=H(\der_{\bm \nu}\ivphi(P)).
\end{equation}
The function, $H$, is useful in proving several properties of admissible solutions, which include the lemma stated below.
The lemma is essential to obtain uniform {\it a priori} estimates of admissible solutions near $\leftsonic\cup\rightsonic$.

\begin{lemma}
\label{lemma-step-2}
Fix $\gam\ge  1$ and $\iv>0$.
For vectors $(\leftvec, \rightvec)$ given by Definition {\rm \ref{definition-domains-np}},
any admissible solution $\vphi$ corresponding to $(\iv, \beta)\in\mathfrak{R}_{\rm weak}$ with $\beta>0$ satisfies
\begin{align}
\label{10-jja}
&\der_{\rightvec}(\vphi-\rightvphi),\;\;\der_{\leftvec}(\vphi-\leftvphi)\ge 0\;\;\qquad\,\,\,\, \text{in $\ol{\Om}$},\\
\label{10-iia}
&-\der_{\etan}(\vphi-\rightvphi),\;\;-\der_{\etan}(\vphi-\leftvphi)\ge 0\;\;\qquad\; \text{in $\ol{\Om}$}.
\end{align}

\begin{proof}
Since $\ivphi-\rightvphi$ is a linear function that vanishes on $\rightshock$,
$\der_{\rightvec}(\vphi-\rightvphi)=\der_{\rightvec}(\vphi-\ivphi)$ in $\ol{\Om}$.
Then \eqref{cond-ad} yields that $\der_{\rightvec}(\vphi-\rightvphi)\ge 0$ in $\ol{\Om}$.
Similarly, \eqref{cond-ad} also implies that $\der_{\leftvec}(\vphi-\leftvphi)\ge 0$ in $\ol{\Om}$. This proves \eqref{10-jja}.

Define
\begin{equation*}
w:=\der_{\xi_2}(\vphi-\rightvphi).
\end{equation*}
We first differentiate Eq. \eqref{2-3} for $\phi=\vphi-\rightvphi$ with respect to $\etan$ to obtain
\begin{equation}
\label{equation-w1}
\begin{split}
&(c^2-\vphi_{\xin}^2)w_{\xin\xin}
-2\vphi_{\xin}\vphi_{\etan}w_{\xin\etan}
+(c^2-\vphi_{\etan}^2)w_{\etan\etan}\\
&+(c^2-\vphi_{\xin}^2)_{\etan}\phi_{\xin\xin}
-2(\vphi_{\xin}\vphi_{\etan})_{\etan}w_{\xin}
+(c^2-\vphi_{\etan}^2)_{\etan}w_{\etan}=0\qquad\tx{in $\Om$}.
\end{split}
\end{equation}
 Since $c^2-\vphi_{\xin}^2>0$
from condition (iii) of Definition \ref{def-regular-sol}, we use Eq. \eqref{2-3} to express $\phi_{\xin\xin}$ as
\begin{equation*}
\phi_{\xin\xin}=
\frac{2\vphi_{\xin}\vphi_{\etan}w_{\xin}-(c^2-\vphi_{\etan}^2)w_{\etan}}
{c^2-\vphi_{\xin}^2}.
\end{equation*}
A direct computation by using \eqref{new-density} yields that $\displaystyle{c^2_{\etan}=-(\gam-1)(\vphi_{\xin}w_{\xin}+\vphi_{\etan}w_{\etan})}$.
Finally,
$(\vphi_{\xi_i}\vphi_{\xi_j})_{\etan}$, $i,j=1,2$, can be expressed
in terms of $(\vphi_{\xin}, \vphi_{\etan}, w, w_{\xin}, w_{\etan})$.
Therefore, Eq. \eqref{equation-w1} can be rewritten as
\begin{equation*}
  (c^2-\vphi_{\xin}^2)w_{\xin\xin}
-2\vphi_{\xin}\vphi_{\etan}w_{\xin\etan}
+(c^2-\vphi_{\etan}^2)w_{\etan\etan}+\sum_{j=1}^2a_j(\vphi_{\xin}, \vphi_{\etan}, w, w_{\xin}, w_{\etan})w_{\xi_j}=0
\qquad\tx{in $\Om$.}
\end{equation*}
This equation is strictly elliptic in $\Om$, and $w$ is not a constant
whenever $\beta>0$, due to Lemma \ref{lemma-7-0}.
Then the maximum principle implies that $\displaystyle{\max_{\ol{\Om}}w=\max_{\der\Om} w}$.

On $\leftsonic\cup\rightsonic$, it follows from the definition of $(\leftvphi, \rightvphi)$
given in Definition \ref{definition-domains-np} and conditions (ii-1) and (ii-3) of Definition \ref{def-regular-sol} that
\begin{equation}
\label{w-on-sonicbd}
  w=
  \begin{cases}
  \der_{\etan}(\leftvphi-\rightvphi)=0\qquad\, &\mbox{on $\leftsonic$},\\
  \der_{\etan}(\rightvphi-\rightvphi)=0\qquad\, &\mbox{on $\rightsonic$}.
  \end{cases}
\end{equation}
Using the slip boundary condition: $\der_{\etan}\vphi=0$ on $\Gam_{\rm wedge}$, stated in
Definition \ref{def-regular-sol}(ii-4), we have
\begin{equation*}
w=0\qquad\tx{ on $\Gam_{\rm wedge}$},
\end{equation*}
since $\der_{\etan}\rightvphi=0$ holds on $\Gam_{\rm wedge}$.

Suppose that there exists a point $\hat P\in \shock$ such that
\begin{equation*}
w(\hat P)=\underset{\ol{\Om}}{\max}\;w,\qquad\,\,  w(\hat{P})>0.
\end{equation*}
Let ${\bm\nu}$ be the unit normal vector to $\shock$ towards the interior of $\Om$, and let ${\bm\tau}$ be a tangent vector to $\shock$.
Since $D^2\ivphi=D^2\rightvphi=-{\bf{I}}_2$, we can rewrite \eqref{5-c1} as
\begin{equation}
\label{5-c4}
D^2(\vphi-\rightvphi)[\bm \tau, \bm h]=0\qquad\;\;\;\text{on $\shock$},
\end{equation}
with ${\bm h}=h_{\nu}{\bm\nu}+h_{\tau}{\bm\tau}$ for $(h_{\nu}, h_{\tau})$ given by \eqref{3-13}.

From the assumption that $w(\hat P)=\underset{\ol{\Om}}{\max}\;w$, it follows that
$\der_{\bm\tau}w(\hat{P})=D^2(\vphi-\rightvphi)[{\bm\tau}, {\ee}_{\xi_2}]=0$ at $\hat{P}$.
Also, by Hopf's lemma, $w$ satisfies
\begin{equation}
\label{5-c5}
\der_{\bm\nu}w(\hat{P})=D^2(\vphi-\rightvphi)[{\bm\nu}, {\ee}_{\xi_2}]<0\qquad\tx{at $\hat{P}$}.
\end{equation}
Then we can use similar arguments as to those for the proof of Lemma \ref{lemma-step1-1} to obtain
\begin{equation}
\label{vector-relation-new2015}
{\ee}_{\xi_2}=k{\bm h}(\hat{P})
\end{equation}
with some constant $k\neq 0$.
By Remark \ref{remark-normalshock}, ${\ee}_{\xi_2}\in\cone$,
so that \eqref{cor-dir-mont} implies that ${\ee}_{\xi_2}\cdot {\bm\nu}<0$ on $\shock$.
Then, at point $\hat{P}$, it follows from \eqref{3-13} and \eqref{vector-relation-new2015} that
\begin{equation*}
  kh_{\nu}(\hat{P})=k{\bm h}(\hat{P})\cdot {\bm\nu}(\hat{P})={\ee}_{\xi_2}\cdot{\bm\nu}(\hat{P})<0.
\end{equation*}
Then we obtain from \eqref{negativity-of-hn} that $k>0$.

By the invariance of Eq. \eqref{2-3} under a coordinate rotation and condition (ii) of Definition \ref{def-regular-sol},
$\phi=\vphi-\rightvphi$ satisfies
\begin{equation}
\label{equation-at-Phat}
(c^2-\vphi_{\bm\nu}^2)\phi_{{\bm\nu}{\bm \nu}}-2\vphi_{\bm\tau}\vphi_{\bm\nu}\phi_{{\bm\nu}{\bm\tau}}+(c^2-\vphi_{\bm\tau}^2)\phi_{{\bm\tau}{\bm\tau}}=0
\qquad\,\tx{at $\hat{P}$}.
\end{equation}
Here and hereafter, we denote $\vphi_{\bm\nu}=\partial_{\bm\nu}\vphi=D\vphi\cdot {\bm\nu}$ and $\vphi_{\bm\tau}=\partial_{\bm\tau}\vphi=D\vphi\cdot {\bm\tau}$
for any function $\vphi$,

Using \eqref{5-c4}, \eqref{equation-at-Phat}, and  Definition \ref{def-regular-sol}(iii), we have
\begin{equation}
\label{sec-ord-expr1}
(\phi_{{\bm\nu}{\bm\tau}},\phi_{{\bm\nu}{\bm\nu}})
=-(\frac{h_{\tau}}{h_{\nu}}, \frac{2\vphi_{\bm\nu}\vphi_{\bm\tau}\frac{h_{\tau}}{h_{\nu}}+(c^2-\vphi_{\bm\tau}^2)}{c^2-\vphi_{\bm\nu}^2}) \phi_{\bm\tau \bm\tau}
\qquad\, \tx{at $\hat P$}.
\end{equation}

Substituting ${\ee}_{\xi_2}=k{\bm h}(\hat{P})$ into \eqref{5-c5}, we obtain
\begin{equation}
\label{phi-2nd-deriv-phat}
D^2\phi[\bm\nu, \bm h]<0 \qquad \tx{at $\hat{P}$}.
\end{equation}
Using \eqref{sec-ord-expr1}, we rewrite \eqref{phi-2nd-deriv-phat} as
$$
A \phi_{{\bm\tau}{\bm\tau}}(\hat P)<0 \qquad \tx{ for $A=\frac{c^4\vphi_{\bm\tau}^2+\rho^2c^2\vphi_{\bm\nu}^2(c^2-|D\vphi|^2)}{\rho\vphi_{\bm\nu}}\,$ at $\hat{P}$}.
$$
Then it follows from Definition \ref{def-regular-sol}(iii) and Lemma \ref{lemma-entropycond-admsbsol} that $A>0$.
Thus, we conclude that $\phi_{{\bm\tau}{\bm\tau}}(\hat P)<0$. This implies that
\begin{equation*}
(\vphi-\ivphi)_{{\bm\tau}{\bm\tau}}(\hat P)<0.
\end{equation*}

Let $f:=\fshock$ be from Proposition \ref{proposition-3}. Then, using
$(\vphi-\ivphi)_{{\bm\tau}{\bm\tau}}(\hat P)<0$ and \eqref{cor-dir-mont}, we have
\begin{equation}
\label{converxity-f-Phat}
f''(\xi_1^{\hat{P}})=\frac{(\vphi-\ivphi)_{{\bm\tau}{\bm\tau}}\big(1+(f')^2\big)}{\der_{\etan}(\ivphi-\vphi)}>0\qquad \tx{at $\hat{P}$},
\end{equation}
since ${\ee}_{\etan}\in \cone$ implies that $\der_{\etan}(\ivphi-\vphi)<0$ at $\hat{P}\in \shock$, due to \eqref{cor-dir-mont}.

Let $\xi_2=L(\xi_1)$ be the equation of the tangent line to $\shock$ at $\hat P$.
Denote $F(\xi_1):=f(\xi_1)-L(\xi_1)$.
Then there exists a point $P_*\neq \hat P$ on ${\rm{int}}\,\shock$ such that
$F(\xi_1^{P_*})=\max_{ [\xi_1^{\lefttop}, \xi_1^{\righttop}]}F(\xi_1)$, due to \eqref{converxity-f-Phat}.

Note that $P_*\not\in \{\lefttop, \righttop\}$, due to \eqref{3-c5a} in Proposition \ref{proposition-3}.
If $P_*=\lefttop$, then $F'(\xin^{\lefttop})\le 0$ must hold,
but this is impossible because $f'(\xin^{P_*})=\tan\beta>f'(\xin^{\hat{P}})=L'(\xin^{P_*})$.
Similarly, if $P_*=\righttop$, then $F'(\xin^{\righttop})\ge 0$ must hold,
but this is also impossible because $f'(\xin^{P_*})=0<f'(\xin^{\hat{P}})=L'(\xin^{P_*})$.
Therefore, we conclude that
$f'(\xi_1^{P_*})=L'(\xi_1^{P_*})=f'(\xi_1^{\hat{P}})$.
This implies that $\bm\nu(P_*)=\bm \nu(\hat P)$.
Denoting $\bm\nu:=\bm\nu(P_*)=\bm \nu(\hat P)$ by $\bm\nu$, we use the definition of $\ivphi$
given in Definition \ref{definition-domains-np} to obtain
\begin{equation}
\label{expression-ivphi-normal}
\der_{\bm\nu}\ivphi(P_*)=\der_{\bm\nu}\ivphi(\hat P)-\big(\der_{\bm\nu}\ivphi(\hat P)-\der_{\bm\nu}\ivphi(P_*)\big)
=\der_{\bm\nu}\ivphi(\hat P)-(P_*-\hat P)\cdot {\bm\nu}.
\end{equation}
For each point $P\in \shock$, we represent $P$ as $(\xin, \fshock(\xin))$ and rewrite the expression as
\begin{equation*}
P=(\xin, \fshock(\xin))=(\xin, F(\xin)+L(\xin))=(\xin, L(\xin))+(0, F(\xin)).
\end{equation*}
By using this expression, $P_*-\hat{P}$ is represented as
\begin{equation*}
P_*-\hat{P}=(\xin^{P_*}-\xin^{\hat{P}})(1,L'(\xin^{\hat{P}}))+\big(F(T_{P_*})-F(T_{\hat P})\big)\bm \ee_{\etan}.
\end{equation*}
Since $L'(\xin^{\hat{P}})=f'(\xin^{\hat{P}})$, $(1, L'(\xi_1^{\hat P}))\cdot{\bm \nu}=(1, f'(\xi_1^{\hat P}))\cdot {\bm \nu}(\hat P)=0$.
This yields that
\begin{equation*}
(P_*-\hat{P})\cdot{\bm\nu}=\big(F(T_{P_*})-F(T_{\hat P})\big)\bm \ee_{\etan}\cdot{\bm\nu}.
\end{equation*}
By substituting this expression into \eqref{expression-ivphi-normal},  $\der_{\bm\nu}\ivphi(P_*)$ is represented as
\begin{equation*}
\der_{\bm\nu}\ivphi(P_*)
=\der_{\bm\nu}\ivphi(\hat P)
 -\big(F(T_{P_*})-F(T_{\hat P})\big)\bm \ee_{\etan}\cdot \bm \nu(P_*).
\end{equation*}
By \eqref{cor-dir-mont} and the definition of $P_*$, $\big(F(T_{P_*})-F(T_{\hat P})\big)\bm \ee_{\etan}\cdot \bm \nu(P_*)< 0$,
which implies that
\begin{equation*}
\der_{\bm\nu}\ivphi(P_*)>\der_{\bm\nu}\ivphi(\hat{P}).
\end{equation*}
This, combined with \eqref{prelim2-1} and \eqref{H-relation-ox}, leads to
\begin{equation}
\label{5-c6-n2017}
\der_{\bm\nu}(\ivphi-\vphi)(P_*)>\der_{\bm\nu}(\ivphi-\vphi)(\hat P)\geq 0.
\end{equation}

We rewrite $w(P_*)$ as
\begin{equation*}
w(P_*)
=\der_{\etan}(\vphi-\ivphi)(P_*)+\underset{(\equiv -\iv)}{\underbrace{\der_{\etan}(\ivphi-\rightvphi)(P_*)}},
\end{equation*}
and further express
$\der_{\etan}(\vphi-\ivphi)(P_*)=(\bm\nu(P_*)\cdot\ee_{ \etan})\der_{\bm\nu}(\vphi-\ivphi)(P_*)$,
where we have used that $\der_{\bm\tau}(\vphi-\ivphi)=0$ holds on $\shock$.
Note that $\bm\nu(P_*)\cdot\ee_{ \etan}=\bm\nu(\hat{P})\cdot\ee_{ \etan}<0$, by \eqref{cor-dir-mont}.
Then it follows from \eqref{5-c6-n2017} that
\begin{equation*}
\begin{split}
w(P_*)
&=\bigl(\bm\nu(P_*)\cdot\ee_{ \etan}\bigr)\der_{\bm \nu}(\vphi-\ivphi)(P_*)+\der_{\etan}(\ivphi-\rightvphi)(P_*)\\
&>\bigl(\bm\nu(\hat P)\cdot\ee _{\etan}\bigr)\der_{\bm\nu}(\vphi-\ivphi)(\hat P)+\der_{\etan}(\ivphi-\rightvphi)(\hat P)=w(\hat P).
\end{split}
\end{equation*}
However, this contradicts the assumption that $w(\hat P)=\underset{\ol{\Om}}{\max}\;w$.

Therefore, we conclude that
$$
\der_{\etan}(\vphi-\rightvphi)\le 0 \qquad \tx{in $\ol{\Om}$}.
$$
Since $\der_{\xi_2}(\rightvphi-\leftvphi)\equiv 0$, we also obtain that $\der_{\xi_2}(\vphi-\leftvphi)\le 0$ in $\ol{\Om}$.
This proves \eqref{10-iia}.
\end{proof}
\end{lemma}

\section{Uniform Positive Lower Bound of ${\rm dist}(\shock, \der B_1(\Oi))$}
\label{subsec-shock-sep-super}
In order to obtain a uniform estimate of the ellipticity of Eq. \eqref{2-1} in the
pseudo-subsonic regions of admissible solutions, it is essential to
make a uniform estimate of positive lower bound of ${\rm dist}(\shock, \der B_1(\Oi))$
for admissible solutions.
Once the estimate of ${\rm dist}(\shock, \der B_1(\Oi))$ is achieved,
the ellipticity of Eq. \eqref{2-1} at each point ${\bm\xi}\in \Om$ is uniformly controlled
by $\rm{dist}({\bm\xi}, \leftsonic\cup\rightsonic)$.
\begin{proposition}
\label{proposition-distance}
Fix $\gam \ge 1$ and $\iv>0$.
Then there exists a constant $C>0$ depending only on $(\iv, \gam)$ such that any admissible solution
corresponding to $(\iv,\beta)$
satisfies
\begin{align}\label{new-dec-2}
{\rm dist}(\shock, \der B_1(\Oi))\ge \frac 1{C}.
\end{align}
\end{proposition}

To prove Proposition \ref{proposition-distance}, some preliminary properties are first required, as shown in Lemmas
\ref{lemma-extcoeff-appc2-new}--\ref{lemma-w-lwrbd-new2015} below.

\smallskip
We rewrite Eq. $\eqref{2-1}$ as
\begin{equation}
\label{8-48-ps}
{\rm div} \mcl{A}(D\vphi,\vphi) +\mcl{B}(D\vphi,\vphi)=0,
\end{equation}
with ${\bf{p}}=(p_1, p_2)\in \R^2$ and $z\in \R$, where
\begin{equation}
\label{8-49-ps}
\begin{split}
&\mcl{A}({\bf p},z):=\rho(|{\bf p}|^2,z){\bf p},\qquad
\mcl{B}({\bf p},z):=2\rho(|{\bf p}|^2,z)
\end{split}
\end{equation}
for $\rho(|{\bf p}|^2,z)$, given by
\begin{equation}
\label{def-densi-ps}
\rho(|{\bf p}|^2,z)=\Bigl(1+(\gam-1)(\frac{\iv^2}{2}-\frac 12|{\bf p}|^2-z)\Bigr)^{\frac{1}{\gam-1}}.
\end{equation}
We also need the definition of $c(|{\bf p}|^2, z)$:
\begin{equation}
\label{definition-soundspeed}
c(|{\bf p}|^2, z):=\rho^{\frac{\gam-1}{2}}(|{\bf p}|^2,z).
\end{equation}
For a constant $R > 1$, define
\begin{equation}
\label{definition-kmn}
\mcl{K}_R=\left\{({\bf p},z)\in\R^2\times \R\,:\,
|{\bf p}|+|z|\le R,\; \rho(|{\bf p}|^2, z)\ge R^{-1}, \;\frac{|{\bf p}|^2}{c^2(|{\bf p}|^2, z)}\le 1-R^{-1} \right\}.
\end{equation}
For each $R>1$, there exists a constant $\lambda_R>0$ depending only on $(\iv, \gam, R)$ such that
\begin{equation*}
\sum_{i,j=1}^2 \der_{p_j}\mcl{A}_i({\bf p},z)\kappa_i\kappa_j\ge \lambda_R|{\bm \kappa}|^2\qquad
\tx{for any}\,\, ({\bf p}, z)\in \mcl{K}_R \,\,\,\tx{and}\,\,\,{\bm \kappa}=(\kappa_1,\kappa_2)\in \R^2.
\end{equation*}

\begin{lemma}[{\cite[Lemma 9.2.1]{CF2}}]
\label{lemma-extcoeff-appc2-new}
For $R>2$, let $\mcl{K}_R$ be given by \eqref{definition-kmn}.
Then there exist functions $(\til{\mcl{A}}, \til{\mcl{B}})({\bf p},z)$ in $\R^2\times \R$
satisfying the following properties{\rm :}

\smallskip
\begin{itemize}
\item[(i)] If $|({\bf p},z)-(\til{\bf p},\til z)|< {\eps}$ for some $(\til{\bf p},\til z)\in \mcl{K}_R$, then
\begin{equation}
\label{extA-1}
(\til{\mcl{A}}, \til{\mcl{B}})({\bf p},z)=(\mcl{A}, \mcl{B})({\bf p},z);
\end{equation}

\item[(ii)] For any $({\bf p},z)\in\R^2\times \R$ and ${\bm\kappa}=(\kappa_1,\kappa_2)\in \R^2$,
\begin{equation}
\label{extA-2}
\sum_{i,j=1}^2 \der_{p_j}\til{\mcl{A}}_i({\bf p},z) \kappa_i\kappa_j\ge {\lambda}|{\bm\kappa}|^2;
 \end{equation}

\item[(iii)]  For each $k=1,2,\cdots$,
\begin{equation}
\label{extA-3}
|\til{\mcl{B}}({\bf p},z)|\le C_0,\quad
|D^k_{({\bf p},z)}(\til{\mcl A}, \til{\mcl{B}})({\bf p},z)| \le C_k\qquad\,\,\tx{in}\,\,\,\R^2\times \R,
\end{equation}
\end{itemize}
where the positive constants $\eps$, $\lambda$, and $C_k$ with $k=0,1,2,\cdots$, depend only on $(\iv, \gam, R)$.
\end{lemma}

For
$\alp\in(0,1)$ and $m\in \mathbb{Z}^+$, we now define the standard H\"{o}lder norms by
\begin{equation}\label{3.2.10a}
\|u\|_{m,0,U}:=\sum_{0\le|\betaa|\le m}\sup_{\rx\in U}|D^{\betaa}u(\rx)|,\quad
[u]_{m,\alp,U}:=\sum_{|\betaa|=m}\sup_{\rx, \ry\in U,\rx\neq  \ry}\frac{|D^{\betaa}u(\rx)-D^{\betaa}u(\ry)|}{|\rx-\ry|^{\alp}},
\end{equation}
where $\betaa=(\beta_1,\beta_2)$ with $\beta_j\ge 0$ for $j=1,2$,
$D^\betaa=\partial_{x_1}^{\beta_1}\partial_{x_2}^{\beta_2}$, and $|\betaa|=\beta_1+\beta_2$.

\begin{lemma}\label{lemma-unif-est1}
Fix $\gam\ge 1$ and $\iv>0$. For any given constants  $\alp\in(0,1)$, $k\in \mathbb{N}$, and $r>0$,
there exist constants $C, C_k>0$ depending only on $(\iv, \gam, \alp, r)$ with $C_k$ depending additionally on $k$
such that any admissible solution $\vphi$ corresponding to $(\iv, \beta)\in \mathfrak{R}_{\rm weak}$
satisfies the following estimates{\rm :}

\smallskip
\begin{itemize}
\item[(i)] For any $B_{4r}(P)\subset \Om$,
\begin{align}
\label{8-40}
&\|\vphi\|_{2,\alp,\ol{B_{2r}(P)}}\le C,\\
\label{8-40-k}
&\|\vphi\|_{k,\ol{B_{r}(P)}}\le C_k.
\end{align}

\item[(ii)] If $P\in \Wedge$, and $B_{4r}(P)\cap \Om$ is the half-ball $B_{4r}^+(P)=B_{4r}(P)\cap\{\xi_2>0\}$, then
 \begin{align}
\label{8-41}
&\|\vphi\|_{2,\alp,\ol{B_{2r}(P)}\cap \Om}\le C,\\
\label{8-41-k}
&\|\vphi\|_{k,\ol{B_{r}(P)}\cap \Om}\le C_k.
\end{align}
\end{itemize}

\smallskip
\begin{proof}
Fix $\beta\in(0,\betadet)$, and let $\vphi$ be an admissible solution corresponding
to $(\iv,\beta)$ with the pseudo-subsonic region $\Om$.
Using Definition \ref{def-regular-sol}(iii) and Lemma \ref{lemma-step3-1},
we can apply Lemmas \ref{lemma6.8}--\ref{lemma6.8-slipbc} to estimate the ellipticity of Eq. \eqref{2-1}.

Suppose that $B_{4r}(P)\subset \Om$ for some constant $r\in(0,1)$.
By \eqref{3-c7}, there exists a constant $\hat c>0$ depending only on $(\iv, \gam)$ such that
any admissible solution $\vphi$ corresponding to $(\iv,\beta)\in \mathfrak{R}_{\rm weak}$
in the sense of Definition \ref{def-regular-sol} satisfies
\begin{equation*}
0<\sup_{\Om}c(|D\vphi|^2,\vphi)\le \hat c.
\end{equation*}

One can choose a smooth function $\til b(\bmxi)$ satisfying the following properties:
\begin{equation*}
\til b= 1\quad\tx{in $B_{3r}(P)$},\qquad
\til b=0 \quad\tx{on $\der B_{4r}(P)$},\qquad
|D^k \til b|\le \frac{C_k}{r^k}\quad\tx{in $B_{4r}(P)$},
\end{equation*}
for constants $C_k>0$ depending only on $k$ for each $k=1,2,\cdots$.
For a constant $\delta_r>0$ to be determined later, we define $b(\bmxi):=\delta_r\til b(\bmxi)$.
Then $b$ satisfies
\begin{equation}
\label{estimate-b-function}
|Db|+\hat c|D^2b|\le \frac{C_*}{r^2}\delta_r\qquad\tx{in $B_{4r}(P)$}
\end{equation}
for some constant $C_*$.

Since ${\rm diam} (\Om)\le \bar d$ for some constant $\bar d>0$ depending only on $(\iv, \gam)$
due to Lemma \ref{lemma-step3-1}, it follows from Lemma \ref{lemma6.8}(b) that
there exists a constant $C_0>0$ depending on $(\iv, \gam)$ such that,
for any given $\delta\in(0,1)$, if $|Db|+\hat c|D^2b|\le \frac{\delta}{\hat{c}}$ in $B_{4r}(P)$,
then either the pseudo-Mach number $M=\frac{|D\vphi|}{c(|D\vphi|^2,\vphi)}$ satisfies that $M^2\le C_0\delta$
in $B_{4r}(P)$  or $M^2+b$ does not attain its maximum in $B_{4r}(P)$.

Now we fix $\delta_r$ in the definition of $b$ as $\delta_r=\frac{ r^2}{8(C_0+1)(C_*+1)\hat c}$.
Then \eqref{estimate-b-function} leads to
\begin{equation*}
|Db|+\hat c|D^2b|\le \frac{1}{8(C_0+1)\hat{c}},
\end{equation*}
which implies that $M=\frac{|D\vphi|}{c(|D\vphi|^2,\vphi)}$ satisfies
\begin{equation*}
\tx{either $M^2\le \frac 18\,$ in $B_{4r}(P)$\qquad  or \qquad  $\displaystyle{\max_{\ol{B_{4r}(P)}} M^2+b=\max_{\der B_{4r}(P)}M^2<1}$}.
\end{equation*}
Therefore, there exists a constant $\sigma_r\in(0,1)$
depending on $(\iv, \gam, r)$ such that $\vphi$ satisfies
\begin{equation}
\label{ellip-local}
\frac{|D\vphi|^2}{c^2(|D\vphi|^2,\vphi)}\le 1-\sigma_r\quad\quad
\tx{in}\,\,B_{3r}(P).
\end{equation}

\smallskip
For a $C^1$--function $\phi$ defined in $U\subset \R^2$, denote $\mcl{E}(\phi, U)$ as
\begin{equation}
\label{8-50}
\mcl{E}(\phi, U):=\{({\bf p}, z)\,:\,z=\phi(\bmxi), {\bf p}=D\phi(\bmxi),\bmxi\in U\}.
\end{equation}

By \eqref{ellip-local} and Lemma \ref{lemma-step3-1}, there exists a constant $R_r>2$
depending only on $(\iv, \gam, r)$ so that $\mcl{E}(\vphi,B_{3r}(P))\subset \mcl{K}_{R_r}$.
Let $(\til{\mcl{A}}, \til{\mcl{B}})({\bf p},z)$ be the extensions given by Lemma \ref{lemma-extcoeff-appc2-new} for $R=R_r$.

In order to prove \eqref{8-40} by applying Theorem \ref{elliptic-t1-CF2}, we rewrite Eq. \eqref{2-1} as
\begin{equation*}
\sum_{i,j=1}^2 \underset{\big(=:A_{ij}(D\vphi,\vphi)\big)}{\underbrace{\der_{p_j}\til{\mcl{A}}_i(D\vphi,\vphi)}}\der_{ij}\vphi
+
\underset{\big(=:A(D\vphi, \vphi)\big)}{\underbrace{\sum_{i=1}^2\der_z\til{\mcl{A}}_i(D\vphi,\vphi)\der_i\vphi
+2\left( \til{\mcl{B}}(D\vphi, \vphi)- \til{\mcl{B}}({\bf 0},0)\right)}}
=-2 \til{\mcl{B}}({\bf 0},0).
\end{equation*}
By Lemma  \ref{lemma-extcoeff-appc2-new} , $(A_{ij},A)(D\vphi,\vphi)$ satisfy \eqref{app-c-prop0}--\eqref{app-c-prop3}.
Then \eqref{8-40} is obtained from Lemma \ref{lemma-step3-1} and Corollary \ref{corollary-t1-CF2}.

\smallskip
Also, \eqref{8-41} is similarly obtained from Lemma \ref{lemma6.8-slipbc} and Theorem \ref{elliptic-t2-CF2}.

\smallskip
Once we have \eqref{8-40} and \eqref{8-41},
estimates \eqref{8-40-k} and \eqref{8-41-k} can be obtained by a bootstrap argument and \cite[Theorem 6.2 and Lemma 6.29]{GT}.
\end{proof}
\end{lemma}
For an admissible solution $\vphi$ corresponding to $(\iv, \beta)\in \mathfrak{R}_{\rm{weak}}$, we define an extension $\vphi^{\rm ext}$ into $\R^2_+$ by
\begin{equation}
\label{definition-admext-new}
  \vphi^{\rm ext}(\bmxi):=\begin{cases}
                               \vphi(\bmxi) & \mbox{if } \bmxi\in \Lbeta, \\
                               \ivphi(\bmxi) & \mbox{otherwise}.
                             \end{cases}
\end{equation}
For $\leftshockseg$ and $\rightshockseg$ defined by Definition \ref{def-regular-sol}, denote $\shock^{\rm{ext}}$ as
\begin{equation*}
\shock^{\rm{ext}}=\begin{cases}
                    \leftshockseg\cup\shock\cup\rightshockseg & \mbox{if } \beta<\betasonic, \\
                    \shock\cup\rightshockseg & \mbox{otherwise}.
                  \end{cases}
\end{equation*}
By
the Rankine-Hugoniot condition: $\vphi=\ivphi$ on $\ol{\shock^{\rm ext}}$, the extension function
$\vphi^{\rm{ext}}$ satisfies the following:

\smallskip
\begin{itemize}
  \item [(i)] $\vphi^{\rm{ext}}\in C^{0,1}_{\rm{loc}}(\R^2_+)\cap C^{1}_{\rm{loc}}(\R^2_+\setminus \ol{\shock^{\rm{ext}}})$;

 \smallskip
  \item [(ii)] $\phi^{\rm{ext}}(\bmxi)=\vphi^{\rm{ext}}(\bmxi)+\frac 12 |\bmxi|^2$ satisfies
       $\|D\phi^{\rm{ext}}\|_{L^{\infty}(\R^2_+)}
        =\|D\phi\|_{L^{\infty}(\Lbeta)}$
      for $\phi(\bmxi):=\vphi(\bmxi)+\frac 12|\bmxi|^2$.
\end{itemize}

\smallskip
In the following corollary, we regard each admissible solution $\vphi$ as its extension $\vphi^{\rm{ext}}$ given by \eqref{definition-admext-new}:

\begin{corollary}
\label{corollary-10-1}
Let $\{\vphi^{(k)}\}$ be a sequence of admissible solutions corresponding to $(\iv, \beta^{(k)})\in \mathfrak{R}_{\rm{weak}}$
in the sense of Definition {\rm \ref{def-regular-sol}} with
\begin{equation*}
\lim_{k\to \infty} \beta^{(k)}=\beta^*\qquad\tx{for some $\beta^*\in[0,\betadet]$}.
\end{equation*}
Then there exists a subsequence $\{\vphi^{(k_j)}\}$ converging to a function $\vphi^*\in C^{0,1}_{\rm{loc}}(\ol{\Lambda_{\beta^*}})$
uniformly in any compact subset of $\ol{\Lambda_{\beta^*}}$,
where $\Lambda_{\beta^{*}}$ is defined by Definition {\rm \ref{definition-domains-np}} for $\beta^*>0$
and by {\rm \eqref{definition-domain-norshock-new}} for $\beta^*=0$.
Moreover, $\vphi^*$
is a weak solution of the boundary value problem consisting of equation \eqref{2-1} in $\Lambda_\beta$ and the slip boundary
condition $\der_{{\bm \nu}}\varphi=0$ on $\partial \Lambda_{\beta^*}$ in the sense of Remark {\rm \ref{remark-admsol-wksol}}.
%
%
For the rest of the statement, let superscripts $(k)$ and $*$ indicate that each object is related
to $\beta^{(k)}$ and $\beta^*$, respectively. Then we have the following properties{\rm :}

\smallskip
\begin{itemize}
\item[(a)] For $P_l$, $l=1,2,3,4$, defined by Definition {\rm \ref{definition-domains-np}},
    \begin{equation*}
      \lim_{j\to \infty} P_l^{(k_j)}=P_l^*\qquad\tx{for $l=1,4$}.
    \end{equation*}
Note that $\righttop$ and $\rightbottom$ are fixed to be the same for all $\beta\in[0,\betadet]$.

\smallskip
\item[(b)] Let $\fshock^{(k_j)}$ be the functions from Proposition {\rm \ref{proposition-3}}.
Extend $\fshock^{(k_j)}$ by
\begin{equation*}
\fshock^{(k_j)}(\xin)=\begin{cases}
f^{(k_j)}_{\mcl{O}}(\xin) &\,\,\,\tx{for}\,\,\, \xin \le \xi_1^{\lefttop^{(k_j)}},\\
\neta&\,\,\,\tx{for}\,\,\, \xin \ge \xi_1^{\righttop},
\end{cases}
\end{equation*}
where $f^{(k_j)}_{\mcl{O}}(\xin)$ is given by \eqref{fo} with $\beta=\beta^{(k_j)}$.
Then sequence $\{f_{\rm{sh}}^{(k_j)}\}$ is uniformly bounded in $C^{0,1}([\xi_1^{P_{\beta^*}}, \xi_1^{\righttop}])$
and converges uniformly on $[\xi_1^{P_{\beta^{*}}}, \xi_1^{\righttop}]$,
where $P_{\beta}$ denotes the $\xi_1$--intercept of the straight oblique shock $\leftshock$
of angle $\beta$ with the $\xi_1$--axis.
Denoting the limit function by $\fshock^{*}$, we see that $\fshock^{*}\in C^{0,1}([\xi_1^{P_{\beta^*}}, \xi_1^{\righttop}])$.

\smallskip
\item[(c)] For each $k_j$, the sonic arcs $\Gam_{\rm{sonic}}^{\mcl{O},(k_j)}$ and $\rightsonic$,
defined by Definition {\rm \ref{definition-domains-np}} corresponding to $(\iv, \beta^{(k_j)})\in \mathfrak{R}_{\rm{weak}}$, can be represented as
    \begin{equation*}
      \begin{split}
         \rightsonic &=\{(\xin, g_{\mcl{N},\rm{so}}(\xin))\,:\,\xin^{\righttop}\le \xin \le \xin^{\rightbottom}\}, \\
         \Gam_{\rm{sonic}}^{\mcl{O},(k_j)}  &=\{(\xin, g_{\mcl{O}, \rm{so}}^{(k_j)}(\xin))\,:\,\xin^{\leftbottom^{(k_j)}}\le \xin \le \xin^{\lefttop^{(k_j)}}\},
      \end{split}
    \end{equation*}
    for smooth functions $g_{\mcl{N}, \rm{so}}$ and $g_{\mcl{O}, \rm{so}}^{(k_j)}$.
Note that $g_{\mcl{N}, \rm{so}}$ is fixed to be the same for all $\beta\in[0,\betadet]$
and that $g_{\mcl{O}, \rm{so}}^{(k_j)}$ depends  continuously on $\beta\in[0,\betadet]$.
Therefore, $g_{\mcl{O}, \rm{so}}^{(k_j)}$ converges to $g_{\mcl{O}, \rm{so}}^{*}$
on $(\xin^{\leftbottom^*}, \xin^{\lefttop^*})$ as $k_j\to \infty$.
If $\beta^*\ge \betasonic$, then it follows from \eqref{10-a3} that $\Gam_{\rm{sonic}}^{\mcl{O},*}$ is a point set.

Define
    \begin{equation*}
      \widehat{\Om^*}:=\{(\xin,\etan)\in [\xin^{\leftbottom^*}, \xin^{\rightbottom}]\times \R^+\,:\, 0\le \etan\le f^*_{\rm{bd}}(\xin)\}
    \end{equation*}
    for a function $f^*_{\rm{bd}}$ given by
    \begin{equation*}
     f^*_{\rm{bd}}(\xin)=\begin{cases}
                            g_{\mcl{O}, \rm{so}}^{*}(\xin) \qquad& \mbox{for } \xin^{\leftbottom^*}\le \xin\le \xin^{\lefttop^*}, \\
                            \fshock^*(\xin) & \mbox{for } \xin^{\lefttop^*}<\xin\le \xin^{\righttop}, \\
                            g_{\mcl{N}, \rm{so}}(\xin)  & \mbox{for } \xin^{\righttop}<\xin\le \xin^{\rightbottom}.
                          \end{cases}
    \end{equation*}
Denote by $\Om^*$ the interior of $\widehat{\Om^*}$. Define
$\shock^{*}:=\{\etan=\fshock^{*}(\xin)\,:\,\xin \in(\xi_1^{\lefttop^{*}}, \xi_1^{\righttop})\}$ and
$\Wedge^{*}:=\{(\xin, 0)\,:\,\xin \in(\xi_1^{\leftbottom^{*}}, \xi_1^{\rightbottom})\}.$ Denote by $\Wedge^{*, 0}$ the relative interior of $\Wedge^*\setminus \shock^*$.
Then $\vphi^*$ satisfies the following properties{\rm :}

\smallskip
\begin{itemize}
\item[\rm (c-1)] $\, \vphi^{*}=\ivphi$ on $\shock^{*}$,

\smallskip
\item[\rm (c-2)] $\, \vphi^{*}\in C^{\infty}(\Om^{*}\cup\Wedge^{*,0})$,

\smallskip
\item[\rm (c-3)] $\, \vphi^{(k_j)}\rightarrow \vphi^{*}$ in $C^2$ on any compact subset of $\Om^{*}\cup\Wedge^{*,0}$,

\smallskip
\item[\rm (c-4)] $\, \der_{{\bf e}}(\ivphi-\vphi^{*})\le 0$ in $\Om^{*}$ for all ${\bf e}\in {\rm Cone}^0({\bf e}_{\leftshock^{*}}, {\bf e}_{\rightshock})$,

\smallskip
\item[\rm (c-5)] $\,$ Eq. \eqref{2-1} is strictly elliptic in $\Om^{*}\cup\Wedge^{*,0}$,

\smallskip
\item[\rm (c-6)] $\,$ The slip boundary condition $\der_{\xi_2}\varphi^*=0$ holds on $\Gamma_{\rm wedge}^{*,0}$,
\end{itemize}
\smallskip
\noindent
where we have followed Definition {\rm \ref{definition-domains-np}} for $(\Oi, \leftsonic, \leftvec)$.
If $\beta^*=0$, ${\rm Cone}^0({\bf e}_{\leftshock^{*}}, {\bf e}_{\rightshock})$ is understood
in the sense of Remark {\rm \ref{remark-condition-v}}.

\smallskip
\item[(d)] In $\Lambda_{\beta^{*}}\setminus {\Om^{*}}$, $\vphi^{*}$ is equal to the constant density
states $\leftvphi^{*}$, $\rightvphi$, and $\ivphi$ in their respective
domains as in \eqref{1-24} if $\beta^{*}\in (0,\betac^{(\iv)})$ and as in \eqref{1-24ab} if $\beta^{*}\in [\betac^{(\iv)},\betadet)$,
where $\leftvphi^*$ is defined by \eqref{2-4-a0} corresponding to $\beta^*$.

\smallskip
\item[(e)] $\fshock^*(\xin)>0$ for all $\xin\in(\xin^{\lefttop^*}, \xin^{\righttop})$.
\end{itemize}

\smallskip
\begin{proof}
By \eqref{3-c6}, the solution structure \eqref{1-24} and \eqref{1-24ab}
in Cases I and II of Definition \ref{def-regular-sol}, and \eqref{definition-admext-new},
it follows that, for any compact $K\subset \overline{\R^2_+}$,
there exists $C(K)<\infty$ such that, for any admissible solution $\vphi$,
$$
  \|\vphi^{\rm ext}\|_{C^{0, 1}(K)}\le C(K).
$$
It follows that there exists a subsequence
$\{\vphi^{(k_j)}\}$ such that the extensions of these functions by  \eqref{definition-admext-new}
converge to a function $\vphi^*\in C^{0,1}_{\rm{loc}}(\ol{\Lambda_{\beta^*}})$
uniformly in any compact subset of $\ol{\Lambda_{\beta^*}}$.

We divide the rest of the proof into four steps.

\smallskip
{\bf 1}. Statement (a) directly follows from Definition \ref{definition-domains-np} and the continuous dependence of $(\Oo, \leftc)$ on $(\iv, \beta)$.
Statement (b) directly follows from Proposition \ref{proposition-3} by selecting a further subsequence of $\{\vphi^{(k_j)}\}$ (without changing notations).

\smallskip
{\bf 2}. Statement (c-1) directly follows from Definition \ref{def-regular-sol}(ii-4), Corollary \ref{corollary-10-1}(a),
and the uniform convergence of $(\vphi^{(k_j)}, \fshock^{(k_j)})$ to $(\vphi^*, \fshock^*)$.
For a point $P\in \Om^*$, there are constants $r>0$ and $N\in\mathbb{N}$ such that $B_{3r}(P)\subset \Om^{(k_j)}$ for all $k_j\ge N$.
Then it follows from Lemma \ref{lemma-unif-est1}(i) and the Arzel\`{a}-Ascoli theorem
that $\vphi^*\in C^{\infty}(B_{3r}(P))$, which implies that $\vphi^*\in C^{\infty}(\Om^*)$.
We can similarly check from Lemma \ref{lemma-unif-est1}(ii) that $\vphi^*\in C^{\infty}(\Om^*\cup \Wedge^{*,0})$, which proves (c-2).

For a fixed compact set $K\subset \Om^{*}\cup\Wedge^{*,0}$,
there exists a constant $N_K\in\mathbb{N}$ so that $K$ is contained in $\Om^{(k_j)}\cap \Wedge^{(k_j)}$ for any $k_j\ge N_K$.
By Lemma \ref{lemma-unif-est1} and the compactness of $K$, $\{\vphi^{(k_j)}\}_{k_j\ge N_K}$ is sequentially compact in $C^2(K)$.
Then the uniform convergence of $\{\vphi^{(k_j)}\}$ to $\vphi^*$ in $K$ implies that the subsequence converges
to $\vphi^*$ in $C^2(K)$. This proves (c-3).

For any ${\bf e}\in {\rm Cone}^0({\bf e}_{\leftshock^{*}}, {\bf e}_{\rightshock})$,
there exists $N_{\bf e}\in \mathbb{N}$ such that
${\bf e}\in {\rm Cone}^0({\bf e}_{\leftshock^{(k_j)}}, {\bf e}_{\rightshock})$ for any $k_j\ge N_{\bf e}$.
Then (c-4) follows from Lemma \ref{lemma-step1-1} and (c-3).

For a point $P\in \Om^*$, we choose $r_P>0$ small so that $B_{r_P}(P)\subset \Om^*$.
Then we fix $N_P\in \mathbb{N}$ sufficiently large so that $B_{r_P}(P)\subset \Om^{(k_j)}$ for all $k_j\ge N_P$.
Since  $\sigma_r\in(0,1)$ in \eqref{ellip-local} is a given constant independent of admissible solutions
corresponding to $\beta\in(0,\betadet)$, we can fix a constant $\sigma_P\in(0,1)$ such that
\begin{equation*}
\frac{|D\vphi^{(k_j)}|^2}{c^2(|D\vphi^{(k_j)}|^2, \vphi^{(k_j)})}\le 1-\sigma_P\qquad\tx{in $B_{r_P}(P)$ for all $k_j\ge N_P$.}
\end{equation*}
This estimate, combined with statement (c-3), implies that Eq. \eqref{2-1} for $\vphi=\vphi^*$ is strictly elliptic in $\Om^*$.
We can use similar arguments by using Lemma \ref{lemma6.8-slipbc}
to conclude that Eq. \eqref{2-1} for $\vphi=\vphi^*$ is strictly elliptic on $\Wedge^{*,0}$, which implies (c-5).
Finally, (c-6) directly follows from (c-3) because every $\vphi^{(k)}$ satisfies the slip boundary condition $\der_{\xi_2}\vphi^{(k)}=0$
on $\Gam_{\rm wedge}^{(k)}$.

\smallskip
Statement (d) follows directly from statements (a)--(c) and Definition \ref{definition-domains-np}.

\smallskip
{\bf 3}. Observe that

\smallskip
\begin{itemize}
\item $f_{\rm w}$ given by \eqref{definition-f-w}, $\lefttop$, $\leftbottom$, $\leftsonic$, and $\leftshockseg$
    depend continuously on $\beta\in [0, \frac{\pi}{2})$;

\smallskip
\item $\righttop$, $\rightbottom$, $\rightshockseg$, and $\rightsonic$ are fixed to be the same for all $\beta\in[0, \frac{\pi}{2})$.
\end{itemize}
Combining this observation with statements (b), (c-3), and (d) implies that, for any compact set $K\subset \R^2$,

\smallskip
\begin{itemize}
\item[(i)] $K\cap \Lambda_{\beta^{(k_j)}}$ converges to $K\cap \Lambda_{\beta^*}$ in the Hausdorff metric;

\smallskip
\item[(ii)]
$D\vphi^{(k_j)}$ converges to $D\vphi^*$ almost everywhere in $K\cap \Lambda_{\beta^*}$.
\end{itemize}

\smallskip
\noindent
Then it follows from Definition \ref{def-regular-sol} that
\begin{equation*}
  \int_{\Lambda_{\beta^*}} \big(\rho(|D\vphi^*|^2,\vphi^*)D\vphi^*\cdot D\zeta-2\rho(|D\vphi^*|^2,\vphi^*)\zeta\big)\,\dd\bmxi=0
  \qquad\tx{for all $\zeta\in C^{\infty}_0(\R^2)$.}
\end{equation*}
In other words, $\vphi^*$ is a weak solution of \eqref{2-1} in $\Lambda_{\beta^*}$
in the sense of  Remark \ref{remark-admsol-wksol}(iv).

\smallskip
{\bf 4}. To prove statement (e), we consider two cases separately: $\beta<\betasonic$ and $\beta\ge \betasonic$.

By Proposition \ref{proposition-3} and statement (b), $\fshock^*$ increases monotonically on $[\xin^{\lefttop^*}, \xin^{\righttop}]$.

If $\beta^*<\betasonic$, then it follows from statement (a) and the monotonicity of $\fshock^*$ that
\begin{equation*}
  \fshock^*(\xin)\ge \fshock^*(\xin^{\lefttop^*})\ge \etan^{\lefttop^*}>0\qquad\tx{for all $\xin\in[\xin^{\lefttop^*}, \xin^{\righttop}]$}.
\end{equation*}

If $\beta^*\ge \betasonic$, it follows from statement (a) and Definition \ref{definition-domains-np}
that $\fshock^*(\xin^{\lefttop^*})=0$.
Suppose that $\fshock^*(\xin)=0$ for some $\xin\in(\xin^{\lefttop^*}, \xin^{\righttop})$. Define
\begin{equation*}
\xin^*:=\sup \{\xin\in(\xin^{\lefttop^*}, \xin^{\righttop})\,:\,\fshock^*(\xin)=0\}.
\end{equation*}
Since $\fshock^*(\xin^{\righttop})=\etan>0$, then $\xin^*\in (\xin^{P_{\beta^*}}, \xin^{\righttop})$.
Note that $\xin^{P_{\beta^*}}=\xin^{\lefttop^*}=\xin^{\leftbottom^*}$ for $\beta^*\ge \betasonic$.
By the monotonicity of $\fshock^*$ with respect to $\xin$, we have
\begin{equation}
\label{f-degeneracy-new}
\fshock^*(\xin)=0\qquad \tx{for all $\xin\in[\xin^{P_{\beta^*}}, \xin^*]$}.
\end{equation}
Let $Q$ be the midpoint of $P_{\beta^*}$ and $(\xin^*,0)$.
Then $Q$ lies on $\Wedge$. Denote $d_*:=\frac{\xin^{P_{\beta^*}}+ \xin^*}{4}$. Then it follows from \eqref{f-degeneracy-new} that
\begin{equation*}
\vphi^*=\ivphi\qquad\tx{in $B_{d_*}(Q)\cap\Lambda_{\beta^*}=B_{d_*}(Q)\cap\{\etan\ge 0\}$}.
\end{equation*}
However, this contradicts the fact that $\vphi^*$
satisfies property (iv) in Remark \ref{remark-admsol-wksol},
%
%
because a direct computation by using the definition of $\ivphi$ given by Definition \ref{definition-domains-np}
shows that a test function $\zeta\in C_0^{\infty}(B_{d_*}(Q))$ can be chosen so that
\begin{equation*}
\begin{split}
&\int_{B_{d_*}(Q)\cap\{\etan\ge 0\}}\big(\rho(|D\ivphi|^2, \ivphi)D\ivphi\cdot D\zeta-2\rho(|D\ivphi|^2,\ivphi)\zeta\big)\,\dd\bmxi\\
&=\iv \int_{\Wedge^*\cap B_{d_*}(Q)}\zeta \,\dd\xin\neq 0.
\end{split}
\end{equation*}
Therefore, we conclude that $\fshock^*(\xin)>0$ holds for any $\xin\in(\xin^{\lefttop^*}, \xin^{\righttop})$, which implies
statement (e). This completes the proof.
\end{proof}
\end{corollary}

Define
\begin{equation}
\label{r1notation}
r_1:=\min_{\beta\in[0,\betadet]}|P_{\beta}|.
\end{equation}
For each $\beta\in[0,\betac^{(\iv)}]$, we know that $|P_{\beta}|\ge\leftc\ge \rightc$, by \eqref{2-4-c2}.
For $\beta\in[\betac^{(\iv)}, \betadet]$, \eqref{12-25} implies that $|P_{\beta}|>\iv\tan\beta\ge \iv\tan\betac^{(\iv)}$.
Therefore, we have
\begin{equation*}
r_1\ge \min\{\rightc, \iv\tan\betac^{(\iv)}\} >0.
\end{equation*}

\begin{proposition}
\label{proposition-sub3}
For every $r\in(0,\frac{r_1}{2})$, there exists a constant $\mcl{C}_r>0$ depending only on
$(\iv, \gam, r)$ such that any admissible solution corresponding to $(\iv, \beta)\in \mathfrak{R}_{\rm weak}$ satisfies
\begin{equation}
\label{10-a5}
{\rm{dist}}(\shock\setminus B_r(P_{\beta}), \Wedge)>\mcl{C}_r^{-1}.
\end{equation}

\begin{proof}
This proposition is proved for two cases separately:
(i) $\leftbottom\not\in B_{\frac{r}{2}}(P_{\beta})$,
and (ii) $\leftbottom \in B_{\frac{r}{2}}(P_{\beta})$ for $\leftbottom$ defined
by Definition \ref{definition-domains-np} depending on $\beta\in[0, \frac{\pi}{2})$.
Fix $r\in(0, \frac{r_1}{2})$.

\medskip
{\textbf{1.}}  We first consider the case that $\leftbottom\not\in B_{\frac r2}(P_{\beta})$.

Define
$$
I_r:=\{\beta\in(0, \betadet)\,:\,|\leftbottom-P_{\beta}|\ge \frac r2\}.
$$
Then $I_r\subset (0, \betasonic)$.
Since $P_{\beta}$ and $\leftbottom$ depend continuously on $\beta\in(0, \betasonic)$,
$I_r$ is relatively closed in $(0, \betasonic)$.
Then there exists $\delta_0>0$ such that,
for any $\beta\in I_r$, $\leftvphi$ given by \eqref{2-4-a0} satisfies that
$\frac{|D\leftvphi(P_{\beta})|}{\leftc(\beta)}\ge 1+\delta_0$.
By Lemma \ref{lemma:interval-existence}, there exists a constant $\sigma_r\in (0, \frac{\betac^{(\iv)}}{2})$
satisfying that $I_r\subset [0, \betac^{(\iv)}-\sigma_r]$.
Then Proposition \ref{proposition-3} implies that
\begin{equation}
\label{estimate-shockwedge-case1}
\inf_{\beta\in I_r} {\rm dist}(\shock, \Wedge)\ge \inf_{\beta\in[0, \betac^{(\iv)}-\sigma_r]} \etan^{\lefttop}>0.
\end{equation}

\smallskip
{\textbf{2.}} Now consider the case that $\leftbottom\in B_{\frac r2}(P_{\beta})$.

For an admissible solution $\vphi$, define
\begin{equation*}
  J_d^{\vphi}:=\{P\in \shock\,:\,|\xin^{P}-\xin^{\leftbottom}|<d\}.
\end{equation*}

{\emph{Claim{\rm :} For any $r\in(0, \frac{r_1}{2})$,
there exists a constant $C_r>0$ such that any admissible solution corresponding to $(\iv, \beta)\in \mathfrak{R}_{\rm weak}$ satisfies
\begin{equation}
\label{10-a7}
\sup_{P\in J_{r/2}^{\vphi}}{\rm{dist}}(P, \Wedge)>C_r^{-1}.
\end{equation}
}}

This claim is proved by deriving a contradiction.
On the contrary, the claim is false.
Then there exists a sequence $\{\beta^{(k)}\}\subset (0, \betadet)$ such that, for each $k\in\mathbb{N}$,
there exists an admissible solution $\vphi^{(k)}$ corresponding to $(\iv, \beta^{(k)})$ in the sense of
Definition \ref{def-regular-sol} with
\begin{equation}
\label{sequence-shockwedge}
\sup_{P\in J_{r/2}^{\vphi^{(k)}}}{\rm dist}(P, \Wedge^{(k)})\le \frac 1k.
\end{equation}
By Corollary \ref{corollary-10-1}, such a sequence $\{\beta^{(k)}\}$ can be chosen
so that it converges to $\beta^*\in[0, \betad^{(\iv)}]$ and the corresponding solution sequence $\vphi^{(k)}$ uniformly converges in any compact subset
of $\ol{\Lambda_{\beta^*}}$ to a function $\vphi^*\in C^{0,1}_{\rm loc}(\ol{\Lambda_{\beta^*}})$
satisfying all the properties described in Corollary \ref{corollary-10-1}.
Furthermore, \eqref{sequence-shockwedge} implies that
\begin{equation*}
  \max_{P\in J_{r/4}^{\vphi^*}}{\rm dist}(P, \Wedge^{*})=0.
\end{equation*}
This contradicts Corollary \ref{corollary-10-1}(e). Thus, the claim is verified.

For each admissible solution $\vphi$, let $\fshock$ be given as an extension defined by Corollary \ref{corollary-10-1}(b).
Then
\begin{equation*}
{\rm dist}(\shock\setminus B_r(P_{\beta}), \Wedge)\ge \fshock(\xi_1^{P_{\beta}}+r)\ge
\sup_{P\in J^{\vphi}_{r/2}} {\rm dist}(P, \Wedge),
\end{equation*}
where we have used the assumption that $|\leftbottom-P_{\beta}|< \frac{r}{2}$ in the second inequality.
Finally, \eqref{10-a5} is directly obtained from this inequality, combined with \eqref{10-a7}.
\end{proof}
\end{proposition}

For $0<\iv\le 1$, define $B_1^+(\Oi):=B_1(\Oi)\cap \{\etan\ge 0\}$.
Following Definition \ref{definition-domains-np}, for each $\beta\in(0,\betadet)$,
$\leftrho>\rightrho>1$ by \eqref{density-mont-ox}.
 Moreover, the entropy condition yields that $|D\ivphi(P_{\beta})|>1$.
 By combining these properties with condition (i-1) of Definition \ref{def-regular-sol},
 any admissible solution corresponding to $(\iv, \beta)\in \mathfrak{R}_{\rm weak}$ satisfies
 \begin{equation}
\label{inclusion}
B_1^+(\Oi)\subset \ol{\Om}\setminus \ol{\leftsonic\cup\shock\cup\rightsonic}.
\end{equation}
For $\iv>1$, \eqref{inclusion} still holds, because $B_1^+(\Oi)=\emptyset$.
Therefore, any compact set
$K\subset B_1^+(\Oi)$ is contained in the pseudo-subsonic region $\Om$.

\begin{lemma}\label{lemma-separ-1}
Fix $\gam\ge 1$ and $\iv\in(0,1)$.
For every compact set $K\subset B_1^+(\Oi)$,
there  exists a constant $C_K>0$ depending only on $(\iv, \gam, K)$
such that any admissible solution $\vphi$ corresponding
to $(\iv,\beta)\in\mathfrak{R}_{\rm weak}$ satisfies
\begin{equation}
    \label{6-a7}
    \inf_{K}(\ivphi-\vphi)\ge C_K^{-1}.
    \end{equation}

\begin{proof}
Suppose that this lemma is false.
By Definition \ref{def-regular-sol}(iv), there exist
a compact set $K\subset B_1^+(\Oi)$, a sequence $\{\beta_j\}\subset (0, \betadet)$,
and a sequence of points $\{Q_j\}\subset K$ so that
\begin{equation*}
(\ivphi-\vphi^{(j)})(Q_j)\rightarrow 0\qquad\,\,\tx{as}\,\, j\to \infty,
\end{equation*}
where $\vphi^{(j)}$ is an admissible solution for each $\beta_j$ in the sense of Definition \ref{def-regular-sol}.
By passing to a subsequence (without changing index notation),
there exist $\beta_{\flat}\in [0,\betadet]$ and $Q_{\flat}\in K$ so that
\begin{equation*}
\beta_j\rightarrow \beta_{\flat}, \quad Q_j\rightarrow Q_{\flat}\qquad\,\, \tx{as $j\to \infty$}.
\end{equation*}
By \eqref{1-24} and \eqref{3-c6}, for any compact set $L\subset \R^2_+:=\{\xxi\in \R^2\,:\,\etan\ge 0\}$,
each $\vphi^{(j)}$ satisfies that
$
\|\vphi^{(j)}\|_{C^{0,1}(L\cap \ol{\Lambda_{\beta_j}})}\le C_L
$
for a positive constant $C_L$ depending only on $(\iv, \gam, L)$.
Therefore, passing to a further subsequence, we conclude that
$\vphi^{(j)}$ converges uniformly to a function
$\vphi_{\flat}\in C^{0,1}(L\cap \ol{\Lambda_{\beta_{\flat}}})$ in $L\cap \Lambda_{\beta_{\flat}}$
for a continuous function $\vphi_{\flat}$ defined in $\Lambda_{\beta_{\flat}}$,
where $\Lambda_{\beta_{\flat}}$ is given by Definition \ref{definition-domains-np}. This yields that
$
(\ivphi- \vphi_{\flat})(Q_{\flat})=0.
$

Since $K$ is compact, there exists a small constant $\epsilon\in(0, \frac{1}{10})$ such that
$
K\subset B^+_{1-2\epsilon}(\Oi).
$
By Lemma \ref{lemma-unif-est1},
sequence $\{\vphi^{(j)}\}$ of admissible solutions
is uniformly bounded in $C^3(\ol{B_{1-\epsilon/2}^+(\Oi)})$.
By the Arzel\'{a}-Ascoli theorem, there exists a subsequence (still denoted by) $\{\vphi^{(j)}\}$
that converges to a function $\vphi_{\flat}\in C^3({\ol{B^+_{1-\epsilon/2}(\Oi)}})$.
Then $\vphi_{\flat}$ satisfies Eq. \eqref{2-1} in $B^+_{1-\eps/2}(\Oi)$,
where the equation is strictly elliptic by  Definition \ref{def-regular-sol}(iii).
Moreover, $\vphi_{\flat}$ satisfies the boundary condition $\der_{\etan}(\ivphi-\vphi)=-\iv<0$
on $B_{1-\epsilon/2}^+(\Oi)\cap \{\etan=0\}$.
Note that condition (iv) of Definition \ref{def-regular-sol} implies that
$\ivphi-\vphi_{\flat}\ge 0$ in $B^+_{1-\eps/2}(\Oi)$.
By Hopf's lemma, $Q_{\flat}$ cannot lie on $B^+_{1-\eps/2}(\Oi)\cap\{\xi_2=0\}$.
Thus, $Q_{\flat}$ must lie in $B^+_{1-\eps/2}(\Oi)$. However, by the strong
maximum principle, this is impossible since $\ivphi-\vphi_{\flat}$ cannot be a constant
in $B^+_{1-\eps/2}(\Oi)$, owing to $\der_{\xi_2}(\ivphi-\vphi_{\flat})=-\iv$
on $B^+_{1-\eps/2}(\Oi) \cap\{\xi_2=0\}$. This completes the proof.
\end{proof}
\end{lemma}

Let $(r,\theta)$ be the polar coordinates centered at $\Oi$:
\begin{equation}
\label{def-polar-oi}
r(\cos\theta, \sin\theta)=(\xin, \etan)-\Oi.
\end{equation}
In $\R^2_+\setminus \{\Oi\}$, define the $(x,y)$--coordinates by
\begin{equation}
\label{prelim5-3}
(x,y)=(\cinfty-r,\theta)\qquad\tx{with $\cinfty=1$.}
\end{equation}
Suppose that a $C^2$--function $\vphi$ satisfies Eq. \eqref{2-1}. We define $w:=\ivphi-\vphi$.
Then Eq. \eqref{2-1} can be written as an equation for $w$ in the $(x,y)$--coordinates:
\begin{equation*}
\Npolar(w):=
\big(2x+(\gam+1)w_x+O_1^-\big)w_{xx}+O_2^-w_{xy}+\big(\frac{1}{\cinfty}+O_3^-\big)w_{yy}
 -(1+O_4^-)w_x+O_5^-w_y\\
=0,
\end{equation*}
with $O_j^-(Dw, w, x)=O_j(-Dw, -w, x,\cinfty)$ for $j=1,\cdots, 5$,
where $O_j(\pb, z, x,c)$ for $j=1,\cdots,5$, with  $\pb=(p_1,p_2)$, are given by
\begin{equation}
\label{prelim5-5}
\begin{split}
&\quad O_1(\pb, z, x,c)=-\frac{x^2}{c}+\frac{\gam+1}{2c}\big(2x-p_1\big)p_1
-\frac{\gam-1}{c}\Big(z+\frac{p_2^2}{2(c-x)^2}\Big),\\
&\quad O_2(\pb, z,x,c)=-\frac{2(p_1+c-x)p_2}{c(c-x)^2},\\
%
&\quad O_3(\pb, z,x,c)=\frac{1}{c(c-x)^2}\Bigl(x(2c-x)-(\gam-1)\big(z+(c-x)p_1
+\frac 12p_1^2\big)-\frac{(\gam+1)p_2^2}{2(c-x)^2}\Bigr),\\
&\quad O_4(\pb, z,x,c)=\frac{1}{c-x}\Bigl(x-\frac{\gam-1}{c}\bigl(z+(c-x)p_1
+\frac 12p_1^2+\frac{(\gam+1)p_2^2}{2(\gam-1)(c-x)^2}\bigr)\Bigr),\\
&\quad O_5(\pb, z,x,c)=-\frac{2(p_1+c-x)p_2}{c(c-x)^3}.
\end{split}
\end{equation}

\begin{lemma}\label{lemma-w-lwrbd-new2015}
For constants $\delta, \eps\ge 0$,  define
\begin{equation*}
\mcl{D}_{-\delta}^{\eps}:=B_{1+\delta}^+(\Oi)\setminus \ol{B_{1-\eps}(\Oi)}.
\end{equation*}
Suppose that $\iv\in(0,1)$ so that $\mcl{D}_{-\delta}^{\eps}\neq \emptyset$ for $\eps>0$.
Then, for any $\alp\in(\frac 12, 1)$, there exist constants $A,\eps_0>0$ depending only on $(\iv, \gam, \alp)$
such that, if $\vphi$ is an admissible solution corresponding to $(\iv, \beta)\in \mathfrak{R}_{\rm weak}$ with $\iv\in(0,1)$,
then $w:=\ivphi-\vphi$ satisfies
\begin{equation*}
w(x,y)\ge Ax^{1+\alp} \quad \quad\text{in}\;\;\mcl{D}_0^{\eps_0}.
\end{equation*}

\begin{proof}
The proof is divided into three steps.

\smallskip
{\textbf{1.}}
Define
$\widehat{O_1^-}(Dw,x):=O_1^-(Dw,w,x)-(\gam-1)w$
and
\begin{equation}
\label{def-N1-new}
\begin{split}
N_1(v):
=&\, \bigl(2x+(\gam+1)v_x+\widehat{O_1^-}+(\gam-1)w\bigr)v_{xx}+O_2^-v_{xy}+(1+O_3^-)v_{yy}\\
&\, -(1+O_4^-)v_x+O_5^-v_y,
\end{split}
\end{equation}
with $\widehat{O_1^-}=\widehat{O_1^-}(Dv,x)$ and $O_j^-=O_j^-(Dv,v,x)$ for $j=2,\cdots,5$.

Fix $\alp\in(\frac 12, 1)$, and define a function
\begin{equation*}
U(x):=Ax^{1+\alp}
\end{equation*}
for a constant $A\in(0,1)$ to be determined later. For each $\eps_0>0$, $U$ satisfies
\begin{align*}
N_1(U)
&\ge\bigl(2x+(\gam+1)U_x+\widehat{O^-_1}(DU,x)\bigr)U_{xx}
-\bigl(1+O_4^-(DU,U,x)\bigr)U_x\\
&\ge(1+\alp)Ax^{\alp}\Big(2\alp-1+\frac{\widehat{O^-_1}}{x}-O^-_4\Big)
\qquad\, \mbox{in $\D_0^{\eps_0}$},
\end{align*}
where we have applied the fact that $w\ge 0$ in $\Om$ by Definition \ref{def-regular-sol}(iv).
Using the definitions of $\widehat{O}_1$ and $O_4$, we can choose $\eps_0>0$ sufficiently small
depending only on $(\iv, \gam, \alp)$ such that
\begin{equation*}
\frac{|\widehat{O_1^-}(DU,x)|}{x}\le\frac{2\alp-1}{4},\quad
|O_4^-(DU,U,x)|\le \frac{2\alp-1}{4} \qquad\,\, \mbox{in $\mcl{D}_0^{\eps_0}$}.
\end{equation*}
Under the choice of $\eps_0$ above,
\begin{equation}
\label{N1-Uw-estimate}
N_1(U)-N_1(w)>0\quad\quad \tx{in}\,\,\mcl{D}_0^{\eps_0}.
\end{equation}

\smallskip
{\textbf{2.}} {\emph{Claim{\rm :} There exists a constant $A>0$
depending only on $(\iv, \gam, \alp)$ such that $U-w$ cannot attain its nonnegative
maximum on $\der \mcl{D}_0^{\eps_0}$.}}

\smallskip
On $\der \mcl{D}^{\eps_0}_0 \cap\{x=0\}$, condition (iv) of Definition \ref{def-regular-sol}
implies that $U-w=-w\le 0$.
By Lemma \ref{lemma-separ-1}, there exists a constant $C_{\eps_0}$ depending only
on $(\iv, \gam, \alp)$ such that
$$
U-w\le A\eps_0^{1+\alp}-C_{\eps_0}\qquad \tx{on $\der\mcl{D}^{\eps_0}_0\cap\{x=\eps_0\}$}.
$$
Thus, a constant $A\in(0,1)$ can be chosen sufficiently small to satisfy that
$A\eps_0^{1+\alp}\le \frac 12 C_{\eps_0}$.
Then we have
$$
U-w\le 0 \qquad \tx{on $\der{\mcl{D}^{\eps_0}_0}\cap\{x=\eps_0\}$}.
$$
Since $\vphi$ satisfies the slip boundary condition on $\Wedge$, $w$ satisfies that
$w_{\etan}=-\iv$ on $\der{\mcl{D}^{\eps_0}_0}\cap \Wedge$ so that
$$
\der_{\etan}(U-w)=A(1+\alp)x^{\alp}\frac{\der x}{\der\etan}+\iv
\qquad\tx{on $\der{\mcl{D}^{\eps_0}_0}\cap \Wedge$}.
$$
Therefore, we can reduce $A>0$ depending only on $(\iv, \gam, \alp)$ so that
$$
\der_{\etan}(U-w)\ge \frac{\iv}{2}
\qquad \tx{on $\der{\mcl{D}^{\eps_0}_0}\cap \Wedge$},
$$
which implies the claim.

\smallskip
{\textbf{3.}} Suppose that $\displaystyle{\max_{\ol{\mcl{D}_0^{\eps_0}}}(U-w)>0}$.
Then there exists a point $P_0\in {\rm int}\, \mcl{D}^{\eps_0}_0$ such that
$$
(U-w)(P_0)=\max_{\ol{\mcl{D}_0^{\eps_0}}}(U-w).
$$
At $P_0$, we have
\begin{equation}
\label{critical-P0}
\begin{split}
&(U-w)_x(P_0)=(U-w)_y(P_0)=0,\\
&(U-w)_{xx}(P_0)\le 0,\quad (U-w)_{yy}(P_0)\le 0,\\
&U_y(P_0)=w_y(P_0)=0,\quad  -w_{yy}(P_0)=(U-w)_{yy}(P_0)\le 0.
\end{split}
\end{equation}
A direct computation by using \eqref{prelim5-5}--\eqref{def-N1-new} and  \eqref{critical-P0} gives that
\begin{align} \label{U-W-1}
&N_1(U)-N_1(w)\\
&= \big(2x+(\gam+1)U_x+\widehat{O^-_1}(DU,x)+(\gam-1)w\big)(U-w)_{xx}\nonumber\\
&\quad\,-\frac{\gam-1}{1-x}(U-w)U_x
 -\left(1+O_3^-(DU,w)\right)w_{yy}\qquad\tx{at $P_0$}.\nonumber
\end{align}
Note that $w(P_0)>0$, by Definition \ref{def-regular-sol}(iv).
Since $|\widehat{O^-_1}(DU,x)|\le C_{O_1}A\eps_0^{2\alp}$ for some constant $C_{O_1}>0$ depending only on $\gam$,
and constant $A$ depends only on $(\gam, \iv, \alp)$,
we can choose $\eps_0>0$ sufficiently small depending on $(\gam, \iv, \alp)$
such that $2x+(\gam+1)U_x+\widehat{O^-_1}(DU,x)+(\gam-1)w>0$ at $P_0$.
Moreover, $(U-w)U_x>0$ at $P_0$. Therefore, we obtain from \eqref{U-W-1} that
\begin{equation*}
N_1(U)-N_1(w)\le -\left(1+O_3^-(DU,w)\right)w_{yy}\qquad\,\,\, {\tx{at}}\,\,P_0.
\end{equation*}
By Definition \ref{def-regular-sol}(iv) and \eqref{prelim5-5},
there exists a constant $C_*>0$ depending only on $\gam$
such that $1+O_3^-(DU, w)\ge 1-C_*\eps_0^{\alp}$ at $P_0$.
Reducing $\eps_0$ further, depending only on $(\gam,\alp)$,
to satisfy that $1-C_*\eps_0^{\alp}\ge \frac 12$,
we obtain that $N_1(U)-N_1(w)\le 0$ at $P_0$.
This contradicts \eqref{N1-Uw-estimate}.
Therefore, we conclude that there exist constants $(A, \eps_0)$
depending on $(\gam, \iv, \alp)$ such that $w\ge Ax^{1+\alp}$ in $\mcl{D}_0^{\eps_0}$.
\end{proof}
\end{lemma}

Now we are ready to prove Proposition \ref{proposition-distance}.
\begin{proofdist}
Let $\vphi$ be an admissible solution corresponding to $(\iv, \beta)\in \mathfrak{R}_{\rm weak}$.
Define
\begin{equation*}
d_{\vphi}:={\rm dist}\{B_1(\Oi), \shock\}.
\end{equation*}
We consider two separate cases: $\iv\ge 1$ and $0<\iv<1$.

\smallskip
{\textbf{1.}} We first consider the case that $\iv \ge 1$.
Then $B_1(\Oi)\subset \R\times \R^-$.
By \eqref{iq-monotonicity} and Lemma \ref{lemma-entropycond-admsbsol},
there exists a constant $d_0>0$ depending only on $(\iv, \gam)$
such that, for any $\beta\in(0, \betadet)$,
\begin{equation*}
{\rm dist}(P_{\beta}, B_1(\Oi))=|P_{\beta}\Oi|-1=|D\ivphi(P_{\beta})|-1\ge M_{\infty, \nu}(P_{\beta})-1\ge  d_0.
\end{equation*}

Denote $\bar r:=\frac 14 \min\{r_1, d_0\}$ for $r_1$ from \eqref{r1notation}.
By Proposition \ref{proposition-sub3}, there exists a constant $\mcl{C}_{\bar r}>0$ depending only
on $(\iv, \gam)$ such that
any admissible solution corresponding to $(\iv, \beta)\in\mathfrak{R}_{\rm weak}$ satisfies
\begin{equation*}
{\rm dist} (\shock \setminus \ol{B_{\frac{\bar r}{2}}(P_{\beta})}, B_1(\Oi))\ge
{\rm dist} (\shock \setminus \ol{B_{\frac{\bar r}{2}}(P_{\beta})}, \Wedge) \ge \mcl{C}_{\bar r}^{-1}>0.
\end{equation*}

By the definition of $\bar r$ above,
$
{\rm dist} (\shock \cap \ol{B_{{\bar r}}(P_{\beta})}, B_1(\Oi))\ge \frac{d_0}{4}>0.
$
Then
\begin{equation*}
d_{\vphi}\ge \min\{\mcl{C}_{\bar r}^{-1}, \frac{d_0}{4}\}>0
\end{equation*}
for any admissible solution $\vphi$ corresponding to $(\iv, \beta)\in \mathfrak{R}_{\rm{weak}}$ with $\iv\ge 1$.

\smallskip
{\textbf{2.}} Now we consider the second case that $0<\iv< 1$.
Let $P_*\in \shock$ be a point
such that
\begin{equation*}
d_{\vphi}={\rm dist}(P_*,B_1(\Oi)).
\end{equation*}
At point $P_*$,  we have
\begin{equation}
\label{d-ivphi-relation}
d_{\vphi}=\der_{\bm\nu}\ivphi(P_*)-1
\end{equation}
for  the unit normal vector ${\bm\nu}$ to $\shock$ at $P_*$ towards the interior of $\Om$.
Denote
\begin{equation}
\label{definition-om-vphi}
\omega_{\vphi}:=\der_{\bm\nu}(\ivphi-\vphi)(P_*).
\end{equation}

\smallskip
{\emph{Claim{\rm :} There exist two positive constants $d_0$ and $d_1$  depending only
on $(\iv, \gam)$ such that, if $d_{\vphi}>d_0$ does not hold,
then $\omega_{\vphi}\ge d_1$ holds.}}

\smallskip
Fix an admissible solution $\vphi$ corresponding to $(\iv, \beta)\in \mathfrak{R}_{\rm weak}$.
For the $(x,y)$--coordinates defined by \eqref{prelim5-3},
let $\eps_0>0$ be the constant from Lemma \ref{lemma-w-lwrbd-new2015} with $\alp=\frac 34$.
In other words, $w:=\ivphi-\vphi$ satisfies
\begin{equation*}
w(x,y)\ge Ax^{\frac{7}{4}}\qquad\tx{in $\mcl{D}_0^{\eps_0}$}
\end{equation*}
for some constant $A>0$ chosen depending only on $(\iv, \gam)$.
For constants $k$ and $\eps\in(0,\eps_0)$, to be determined later,
define a function $V$ in $\mcl{D}_{-d_{\vphi}}^{\eps}$ by
\begin{equation}
\label{definition-Vftn-new}
V:=(x+d_{\vphi})^2+k(x+d_{\vphi}).
\end{equation}
For a constant $d_0>0$ to be specified later, assume that $d_{\vphi}\le d_0$.
Then a direct computation by using \eqref{prelim5-3}--\eqref{prelim5-5}
and Definition \ref{def-regular-sol}(iv)
shows that $V$ satisfies
\begin{equation}\label{V-properties-new}
\begin{split}
N_1(V)\ge 3k-4d_0-C(\eps+d_0+k)^2\qquad &\tx{in $\mcl{D}^{\eps}_{-d_{\vphi}}$},\\
V=0 \qquad &\tx{on $\der \mcl{D}^{\eps}_{-d_{\vphi}}\cap\{x=-d_{\vphi}\}$},\\
V\le (\eps+d_0)^2+k(\eps+d_0)\qquad &\tx{on $\der \mcl{D}^{\eps}_{-d_{\vphi}}\cap\{x=\eps\}$}, \\
V_{\etan}\ge \frac{-\iv}{1-\eps}\left(2(\eps+d_0)+k\right)
\qquad &\tx{on $\der \mcl{D}^{\eps}_{-d_{\vphi}}\cap \Wedge$},
\end{split}
\end{equation}
for a constant $C>0$ chosen depending only on $(\gam, \iv)$.
Choosing
\begin{equation*}
k=2\eps, \qquad d_0=\eps,
\end{equation*}
we obtain from \eqref{V-properties-new}, $w\ge 0$ in $\ol{\Om}$, and \eqref{2-4-b6} that
\begin{equation}
\label{VW-properties-new}
\begin{split}
N_1(V)-N_1(w)\ge 2\eps-16C\eps^2\qquad &\tx{in $\mcl{D}^{\eps}_{-d_{\vphi}}$},\\
V-w\le 0 \qquad &\tx{on $\der \mcl{D}^{\eps}_{-d_{\vphi}}\cap\{x=-d_{\vphi}\}$},\\
V-w\le 10\eps^2-A\eps^{\frac 74}\qquad &\tx{on $\der \mcl{D}^{\eps}_{-d_{\vphi}}\cap\{x=\eps\}$}, \\
(V-w)_{\etan}\ge \iv-\frac{6\iv \eps}{1-\eps}
\qquad &\tx{on $\der \mcl{D}^{\eps}_{-d_{\vphi}}\cap \Wedge$}.
\end{split}
\end{equation}
Then we can fix a small constant ${\eps}\in(0,\eps_0)$ depending only on $(\iv, \gam)$ such that,
by \eqref{VW-properties-new}, $N_1(V)-N_1(w)\ge 0$ in $\mcl{D}^{\eps}_{-d_{\vphi}}$,
$V-w\le 0$ on $\der \mcl{D}^{\eps}_{-d_{\vphi}}\cap\{x=-d_{\vphi}\,\,\tx{or}\,\,\eps\}$, and
$(V-w)_{\etan}\ge 0$ on $\der \mcl{D}^{\eps}_{-d_{\vphi}}\cap \Wedge$.
Thus, the maximum principle yields that
\begin{equation}
\label{estimate-VW-new}
V-w\le 0\qquad \tx{in $\ol{\mcl{D}^{\eps}_{-d_{\vphi}}}$.}
\end{equation}
Since $P_*\in \der \mcl{D}^{\eps}_{-d_{\vphi}}\cap\{x=-d_{\vphi}\}$,
$\displaystyle{(V-w)(P_*)=\max_{\ol{\mcl{D}^{\eps}_{-d_{\vphi}}}}(V-w)=0}$.
Note that $\shock$ is tangential to $\der \mcl{D}^{\eps}_{-d_{\vphi}}\cap \{x=-d_{\vphi}\}$ at $P_*$
so that $(V-w)_x(P_*)=\der_{\bm \nu}(V-w)(P_*)$.
Then \eqref{estimate-VW-new} implies that $(V-w)_x(P_*)=\der_{\bm \nu}(V-w)(P_*)\le 0$.
Combining this with \eqref{definition-om-vphi}--\eqref{definition-Vftn-new} implies that
\begin{equation*}
\om_{\vphi}\ge V_x(P_*)=2\eps.
\end{equation*}
Therefore, the claim is verified by choosing
$(d_0, d_1):=(\eps, 2\eps)$.

\smallskip
According to the claim, either $d_{\vphi}$ is bounded below by $\eps$ or $\omega_{\vphi}$ is bounded below by $2\eps$.
By \eqref{H-relation-ox} and \eqref{d-ivphi-relation},  $\omega_{\vphi}=H(d_{\vphi}+1)$ for $H$ defined by \eqref{def-H-ox}.
Then it follows from  \eqref{H-property-new} that $d_{\vphi}$ is uniformly bounded below by a positive constant if and only if $\omega_{\vphi}$ is
uniformly bounded below by a positive constant.
Therefore, the claim implies that there exists a constant $\delta>0$ depending only on $(\iv, \gam)$ such that
\begin{equation*}
d_{\vphi}\ge \min\{\eps, \delta\}>0
\end{equation*}
for any admissible solution $\vphi$
corresponding to $(\iv, \gam)\in \mathfrak{R}_{\rm weak}$ with $0<\iv<1$.

\smallskip
The proof of Proposition \ref{proposition-distance} is now completed. \hfill $\Box$
\end{proofdist}

\section{Uniform Estimates for the Ellipticity of Eq. \eqref{2-1}}
\label{subsec-unif-est-ellip-super}
Given $\gam\ge 1$ and $\iv>0$,
let $\vphi$ be an admissible solution corresponding to $(\iv, \beta)\in \mathfrak{R}_{\rm weak}$.
A direct computation by using \eqref{8-49-ps} shows that Eq. \eqref{8-48-ps} (the same as Eq. \eqref{2-1}) satisfies
\begin{equation}
\label{ellip-1}
\rho(1-\frac{|D\vphi|^2}{c^2})|\bm\kappa|^2\le \sum_{i,j=1}^2\der_{\bm p_i}\mcl{A}_j(D\vphi,\vphi)\kappa_i\kappa_j
\le 2\rho |\bm\kappa|^2\qquad\,\,\text{in $\Om\,$ for any $\bm\kappa=(\kappa_1,\kappa_2)\in \R^2$}.
\end{equation}

Fix a function $h\in C^{\infty}(\R_+)$ such that
\begin{equation}
\label{8-12}
h(s)=
\begin{cases}
s&\text{if $s\in[0,\frac{1}{2}]$},\\
1&\text{if $s\ge 1$},
\end{cases}
\qquad\tx{and}\qquad
0\le h'\le 2\quad \text{on $\R_+$}.
\end{equation}
For each $\beta\in(0, \frac{\pi}{2})$, let $\Oo$ be defined by Definition \ref{definition-domains-np}, and denote
\begin{equation*}
r_{\beta}:=\min\{\leftc, |\Oo P_{\beta}|\}=
\begin{cases}
\leftc\qquad &\tx{if $\beta<\betac^{(\iv)}$},\\
|\Oo P_{\beta}|\qquad &\tx{if $\beta\ge \betac^{(\iv)}$}.
\end{cases}
\end{equation*}
Let $Q_{\mcl{O}}\in\leftshock\cap\{\etan\ge 0\}$  be the midpoint of the two intersections of
circle $|{\bm\xi}-\Oo|=r_{\beta}$ and $\leftshock\cap \{\etan\ge 0\}$, and let
\begin{equation*}
\hat{r}_{\beta}:=|\Oo Q_{\mcl{O}}|
=\begin{cases}
r_{\beta}\oM\qquad &\mbox{for $\beta<\betasonic$},\\
r_{\beta}\sin\beta\qquad &\mbox{for $\beta\ge \betasonic$},
\end{cases}
\end{equation*}
for $\oM$ defined by \eqref{1-25}.
Note that $r_{\beta}$ and $\hat{r}_{\beta}$ depend continuously on $\beta\in(0, \frac{\pi}{2})$.
It follows from \eqref{oM-monotonicity} and the definitions of $(r_{\beta}, \hat{r}_{\beta})$ stated above
that $r_{\beta}-\hat{r}_{\beta}>0$ for all $\beta\in[0, \frac{\pi}{2})$.
Therefore, there exists a constant  $\delta_0>0$ depending only on $(\iv, \gam)$
so that $r_{\beta}-\hat{r}_{\beta}\ge \delta_0$ for all $\beta\in[0, \betadet]$.

We define $(\go, \gn, \Qn)$ by
\begin{equation*}
\begin{split}
&{\go}(\bmxi):=\frac 12(r_{\beta}-\hat r_{\beta})\,h(\frac{{\rm dist}(\bmxi, \der B_{r_{\beta}}(\Oo))}{r_{\beta}-\hat r_{\beta}}),\\
&{\gn}(\bmxi):=\lim_{\beta\to 0+}\go(\bmxi)=\frac 12(\rightc-\neta)\,h(\frac{{\rm dist}(\bmxi, \der B_{\rightc}(O_{\mcl{N}}))}{\rightc-\neta}),\\
&\Qn:=\lim_{\beta\to 0+}\Qo.
\end{split}
\end{equation*}
Let $Q^*=(\xi_1^*, \neta)$ be the midpoint of $\Qn$ and $\righttop$ for point $\righttop$ given by Definition \ref{definition-domains-np}.
Moreover, we fix a function $\chi=\chi(\xin)\in C^{\infty}(\R)$
such that
\begin{equation*}
\chi(\xin)
=
\begin{cases}
1&\text{for $\xin\le \frac{\xin^*}{10}$},\\[1mm]
0&\text{for $\xin\ge \frac{\xin^*}{2}$},
\end{cases}\qquad\,\,\,
-\frac{5}{\xin^*}\le \chi'(\xin)\le 0\;\;\,\text{for all $\xin\in \R$}.
\end{equation*}
Finally, we define a function $g_{\beta}:\R^2\rightarrow \R_+$ by
\begin{equation}
\label{8-13}
g_{\beta}(\bmxi):=\chi(\xin)\big({\go}(\bmxi)+\max\big\{1-\frac{|D\leftvphi(P_{\beta})|^2}{\leftc^2},0\big\}\big)+\big(1-\chi(\xi_1)\big){\gn}(\bmxi).
\end{equation}

\begin{remark}
\label{remark-ellipticity-new2015}
By Definition {\rm \ref{def-regular-sol}} and Lemma {\rm \ref{lemma-step3-1}}, there exist constants $d>0$ and $C>1$ depending only
on $(\iv, \gam)$ such that, if $\vphi$ is an admissible solution corresponding to $(\iv, \beta)\in\mathfrak{R}_{\rm weak}$,
and if $\Om$ is its pseudo-subsonic region, then $g_{\beta}$ satisfies the following properties{\rm :}

\smallskip
\begin{itemize}
\item[(i)] For $\bmxi\in \Om$ satisfying ${\rm dist}(\bmxi, \rightsonic)<d$,
\begin{equation*}
C^{-1} {\rm dist}(\bmxi, \rightsonic)\le g_{\beta}(\bmxi) \le C {\rm dist}(\bmxi, \rightsonic);
\end{equation*}

\item[(ii)] For $\bmxi\in \Om$ satisfying ${\rm dist}(\bmxi, \leftsonic)<d$,
\begin{equation*}
C^{-1} {\rm dist}_{\beta}(\bmxi, \leftsonic)
\le g_{\beta}(\bmxi)
\le C {\rm dist}_{\beta}(\bmxi, \leftsonic),
\end{equation*}
where ${\rm dist}_{\beta}(\bmxi, \leftsonic)$ is given by
\begin{equation}
\label{def-distbe}
{\rm dist}_{\beta}(\bmxi, \leftsonic):={\rm dist}(\bmxi, \leftsonic)+\left(\leftc-|D\leftvphi(\lefttop)|\right);
\end{equation}

\item[(iii)] Furthermore, for each $\eps>0$, there exists a constant $C_{\eps}>1$ depending only
on $(\iv, \gam, \eps)$ such that,
if a point $\bmxi\in \ol{\Om}$ satisfies ${\rm dist}(\bmxi, \leftsonic\cup\rightsonic)>\eps$, then $g_{\beta}$ satisfies
\begin{equation*}
C_{\eps}^{-1}\le g_{\beta}(\bmxi)\le C_{\eps}.
\end{equation*}
\end{itemize}
In {\rm (i)}--{\rm (iii)}, $\rightsonic$, $\leftsonic$, and $\leftvphi$ are defined by Definition {\rm \ref{definition-domains-np}}.
\end{remark}

For a constant $\hat{\zeta}>0$, let us define
\begin{equation}
\label{def-distfl}
{\rm dist}^{\flat}(\bmxi, \leftsonic\cup\rightsonic):=
\min\left\{\hat{\zeta},\, {\rm dist}(\bmxi, \rightsonic),\, \rm dist_{\beta}(\bmxi, \leftsonic)\right\}.
\end{equation}
Using properties {\rm (i)}--{\rm (iii)} stated in Remark \ref{remark-ellipticity-new2015},
we can find constants $\mcl{C}>1$ and $\hat{\zeta}\in(0,1)$ depending only on $(\iv, \gam)$
such that each $g_{\beta}$ for $\beta\in(0, \betadet)$ satisfies
\begin{equation*}
\mcl{C}^{-1} {\rm dist}^{\flat}(\bmxi, \leftsonic\cup\rightsonic)
\le g_{\beta}(\bmxi) \le \mcl{C} {\rm dist}^{\flat}(\bmxi, \leftsonic\cup\rightsonic)\,\,\qquad\tx{for all}\,\,\bmxi\in \ol{\Om},
\end{equation*}
where $\Om$ is the pseudo-subsonic region of an admissible solution $\vphi$ corresponding to $(\iv, \beta)$.

\smallskip
Let $\mcl{A}({\bf p}, z)$ be given by \eqref{8-49-ps}.
The following proposition is essential to establish {\it a priori} weighted $C^{2,\alp}$ estimates of admissible solutions:
\begin{proposition}
\label{corollary-ellip}
There exists a constant $\mu>0$ such that,
if $\vphi$ is an admissible solution corresponding to $(\iv, \beta)\in \mathfrak{R}_{\rm weak}$
and $\Om$ is its pseudo-subsonic region, then the pseudo-Mach number given by
\begin{equation}
\label{definition-pseudomach}
M(\bmxi):=\frac{|D\vphi(\bmxi)|}{c(|D\vphi|^2(\bmxi),\vphi(\bmxi))}
\end{equation}
satisfies
\begin{equation}
\label{estimate-M-2015nov}
M^2(\bmxi)\le 1-\mu g_{\beta}(\bmxi)\qquad \tx{in}\,\,\ol{\Om},
\end{equation}
and there exists a constant $C>1$ such that
\begin{equation}\label{9-32}
C^{-1} {\rm dist}^{\flat}(\bmxi,\leftsonic\cup\rightsonic)|\bm\kappa|^2
\le \sum_{i,j=1}^2\mcl{A}_{p_j}^i(D\vphi(\bmxi), \vphi(\bmxi))\kappa_i\kappa_j\le C|\bm\kappa|^2
\end{equation}
for all $\bmxi\in\ol{\Om}$ and $\bm{\kappa}=(\kappa_1, \kappa_2)\in \R^2$,
where constants $\mu$ and $C$ are chosen depending only on $(\iv, \gam)$.
On the left-hand side of \eqref{9-32},
${\rm dist}^{\flat}(\cdot,\cdot)$ is given by \eqref{def-distfl}.

\begin{proof}
Once \eqref{estimate-M-2015nov} is proved, \eqref{9-32} is obtained directly from \eqref{estimate-M-2015nov},
Lemma \ref{lemma-step3-1}, \eqref{ellip-1}, and Remark \ref{remark-ellipticity-new2015}.
Therefore, it now suffices to prove \eqref{estimate-M-2015nov}.

In this proof, $\vphi$ represents any admissible solution corresponding to $(\iv, \beta)\in \mathfrak{R}_{\rm weak}$
with $\Om$ and $\shock$ being its pseudo-subsonic region and the curved transonic shock, respectively.
Unless otherwise specified, all the constants appearing in the proof are chosen depending only on $(\iv, \gam)$.
The proof is divided into four steps.

\smallskip
{\textbf{1.}} By Lemma \ref{lemma-step3-1}, there exist constants $R>1$ and $\hat c>1$ such that
\begin{equation*}
\Om\subset B_{R/2}({\bf 0}),\quad\, \|c(|D\vphi|^2,\vphi)\|_{C^0(\ol{\Om})}\le \hat c,\quad\,
\|g_{\beta}\|_{C^2(\ol{\Om})}\le \hat c
\end{equation*}
for
$g_{\beta}$ given by \eqref{8-13}.
Since $\Oo\in \{\etan=0\}$,
$\der_{\etan}g_{\beta}=0$ on $\{\etan=0\}$.
By Lemmas \ref{lemma6.8}--\ref{lemma6.8-slipbc}, we can choose constants $C_0>0$, $\delta\in(0, \frac{3}{4}C_0)$,
and $\mu_1\in(0,1)$ so that, whenever $\mu\in(0, \mu_1]$,
either the inequality: $M^2+\mu g_{\beta}\le C_0\delta <1$ holds in $\Om$,
or the maximum of $M^2+\mu g_{\beta}$ over $\ol{\Om}$ cannot be
attained in $\Om\cup \Wedge$.

Since $M^2+\mu g_{\beta}=1$ on $\rightsonic$,
the maximum of $M^2+\mu g_{\beta}$ must be attained on $\der \Om \setminus \Wedge$.

\smallskip
{\textbf{2.}} Let $\bm\nu$ be the unit normal vector to $\shock$ towards the interior of $\Om$,
and let $\bm\tau$ be a unit tangent vector to $\shock$.

\smallskip
{\emph{Claim{\rm :}
There exist constants $\alp\in(0, \frac 12)$ and $\zeta\in(0,1)$ such that
$M^2(P)\le 1-\zeta$ when
$|\vphi_{\bm\tau}|^2\le \alp |\vphi_{\bm\nu}|^2$
at $P\in \shock$.
}}

\smallskip
This claim is verified by adjusting the proof of \cite[Lemma 9.6.2]{CF2}.
For a constant $\alp\in(0, \frac 12)$ to be specified later,
assume that $|\vphi_{\bm\tau}|^2\le \alp |\vphi_{\bm\nu}|^2$
holds at $P\in \shock$.
Since $\rho\vphi_{\bm\nu}=\der_{\bm\nu}\ivphi$ and $\vphi_{\bm\tau}=\der_{\bm\tau}\ivphi$ hold along $\shock$, we have
\begin{equation*}
|D\ivphi|^2-|\der_{\bm\nu}\ivphi|^2=|\vphi_{\bm\tau}|^2\le
\alp |\vphi_{\bm\nu}|^2\le \alp \Big(\frac{\der_{\bm\nu}\ivphi}{\rho}\Big)^2,
\end{equation*}
which yields that
\begin{equation*}
|D\ivphi|^2\le \big(1+\frac{\alp}{\rho^2}\big)|\der_{\bm\nu}\ivphi|^2\qquad \tx{at}\,\,P\in\shock.
\end{equation*}
We combine this inequality with Lemma \ref{lemma-step3-1} and Proposition \ref{proposition-distance} to obtain
\begin{equation*}
|\der_{\bm\nu}\ivphi(P)|^2\ge \frac{1+d_0}{1+ \alp/C}
\end{equation*}
for some constants $d_0>0$ and $C>1$.
Therefore, we can fix constants $\bar{\alp}\in(0,\frac 12)$ and $d_1>0$ small so
that  $|\der_{\bm\nu}\ivphi(P)|\ge 1+d_1$ when $\alp\in[0, \bar{\alp}]$.

Define $M_{\infty,\bm\nu}:=|\der_{\bm\nu}\ivphi(P)|$ and  $M_{\bm\nu}:=\frac{|\vphi_{\bm{\nu}}(P)|}{c(|D\vphi|^2(P), \vphi(P))}$.
Then it follows from \eqref{1-1} that
\begin{equation*}
\big(1+\frac{\gam-1}{2}M_{\bm\nu}^2\big) M_{\bm\nu}^{-\frac{2(\gam-1)}{\gam+1}}
=\big(1+\frac{\gam-1}{2}(M_{\infty,\bm\nu})^2\big)|M_{\infty,\bm\nu}|^{-\frac{2(\gam-1)}{\gam+1}}.
\end{equation*}
Owing to $M_{\infty,\bm\nu}=|\der_{\bm\nu}\ivphi(P)|\ge 1+d_1$,
there exists a constant $\zeta_*\in(0,1)$ satisfying that $M_{\bm\nu}^2\le 1-\zeta_*$ at $P\in \shock$.
By the assumption that $|\vphi_{\bm\tau}|^2\le \alp |\vphi_{\bm\nu}|^2$ at $P\in \shock$, we have
\begin{equation*}
M^2\le (1+\alp)M_{\bm\nu}^2\le (1+\alp)(1-\zeta_*)\qquad\tx{at $P\in \shock$}.
\end{equation*}
Therefore, we can further reduce  $\alp\in(0,\bar{\alp}]$ so that the inequality right above implies that
\begin{equation*}
M^2\le 1-\frac{\zeta_*}{2}=:1-\zeta\qquad\tx{at $P\in \shock$}.
\end{equation*}
The claim is verified.

\smallskip
{\textbf{3.}}
Let $\mu_1$ be the constant from  Step 1.
In this step, we follow the approach of \cite[Steps 2--3 in the proof of Proposition 9.6.3]{CF2} to find a constant $\mu\in(0,\mu_1]$
so that $M^2+\mu g_{\beta}$ cannot attain its maximum  on $\shock$.
 Here, we give an outline to see how such a constant $\mu$ is chosen.
We refer to \cite[Proposition 9.6.3]{CF2} for further details.

\smallskip
{\textbf{3-1}}. Suppose that the maximum of $M^2+\mu g_{\beta}$ over $\ol{\Om}$ is attained at $\Pmax\in \shock$.
Then $(M^2+\mu g_{\beta})(\Pmax)\ge 1$, which implies that
\begin{equation}
\label{M-at-pmax}
M^2(\Pmax)\ge 1-C_*\mu
\end{equation}
for some constant $C_*>0$.
Moreover, we have
\begin{align}
\label{deriv-comp-1}
&\der_{\bm\tau}(M^2+\mu g_{\beta})(\Pmax)=0,\\
\label{deriv-comp-2}
&\der_{\bm\nu}(M^2+\mu g_{\beta})(\Pmax)\le 0.
\end{align}
For simplicity of notation, denote
\begin{equation}
\label{definition-kfunctiin-new}
  k(\bmxi):=\mu g_{\beta}(\bmxi)\qquad \tx{for $\bmxi\in \R^2$}.
\end{equation}
By using \eqref{new-density} and \eqref{definition-localsoundsp},
a direct computation yields that, for each unit vector ${\bf w}$,
\begin{equation}
\label{expression-Mach-deriv}
(M^2)_{\bf{w}}=\frac{\big(2+(\gam-1)M^2\big)D^2\vphi[{\bf w}, D\vphi]+(\gam-1)M^2 \vphi_{\bf{w}}}{c^2},
\end{equation}
where we have defined
\begin{equation*}
  D^2\vphi[{\bf q}_1, {\bf q}_2]:=(D^2\vphi\, {\bf q}_1)\cdot {\bf q}_2
  \qquad\tx{for ${\bf q}_1, {\bf q}_2\in \R^2$}.
\end{equation*}
By \eqref{expression-Mach-deriv}, we obtain from \eqref{deriv-comp-1} that
\begin{equation}
\label{from-pmaxassump-shock}
  D^2\vphi[{\bm\tau}, D\vphi]=-\frac{(\gam-1)M^2\vphi_{\bm\tau}+c^2 k_{\bm\tau}}{2+(\gam-1)M^2}=:B_1\qquad\tx{at $\Pmax$.}
\end{equation}

\smallskip
{\textbf{3-2}}. Next, we differentiate the Rankine-Hugoniot condition:
\begin{equation}
\label{RH-derivcond-new}
(\rho D\vphi-D\ivphi)\cdot D(\ivphi-\vphi)=0\qquad \tx{on $\shock$}
\end{equation}
in the tangential direction $\bm\tau$ of $\shock$, and then use  \eqref{2-4-b6}--\eqref{new-density}
and $(\ivphi-\vphi)_{\bm\tau}=0$ on $\shock$ to obtain
\begin{equation}
\label{tangentderiv1-RH-new}
\begin{split}
&\big(\rho D^2\vphi\, \bm\tau-\frac{\rho}{c^2}(D\vphi\cdot (D^2\vphi\, \bm\tau)+\vphi_{\bm\tau})D\vphi\big)\cdot (D\ivphi-D\vphi)\\
&-(\rho D\vphi-D\ivphi)\cdot (D^2\vphi \, \bm\tau+\bm\tau)=0 \qquad\,\,\,\,\,\, \tx{on $\shock$.}
\end{split}
\end{equation}

Using the Rankine-Hugoniot conditions \eqref{RH-derivcond-new} and $(\ivphi-\vphi)_{\bm\tau}=0$ on $\shock$,
we see that $\displaystyle{D(\ivphi-\vphi)=\der_{\bm\nu}(\ivphi-\vphi){\bm\nu}}=(\rho-1)\vphi_{\bm\nu}{\bm\nu}$.
Then we obtain
\begin{equation*}
  D\vphi\cdot D(\ivphi-\vphi)=(\rho-1)\vphi_{\bm\nu}^2 \qquad\,\,\tx{on $\shock$.}
\end{equation*}
Owing to the condition that $(\ivphi-\vphi)_{\bm\tau}=0$ on $\shock$ again, we have
\begin{equation*}
  (\rho D\vphi-D\ivphi)\cdot \bm\tau=(\rho-1)\vphi_{\bm\tau}\qquad\,\,\tx{on $\shock$.}
\end{equation*}
We substitute the expressions of $D\vphi\cdot D(\ivphi-\vphi)$ and $ (\rho D\vphi-D\ivphi)\cdot \bm\tau$
given above into \eqref{tangentderiv1-RH-new} to obtain
\begin{equation}
\label{tangentderiv2-RH-new}
\begin{split}
  &D^2\vphi[\bm\tau, \rho D(\ivphi-\vphi)+D\ivphi]\\
  &=\rho(1+\frac{\rho-1}{c^2}\vphi_{\bm\nu}^2) D^2\vphi[\bm\tau, D\vphi]+\frac{\rho}{c^2}(\rho-1)\vphi_{\bm\nu}^2\vphi_{\bm\tau}
  +(\rho-1)\vphi_{\bm\tau}\qquad\tx{on $\shock$}.
  \end{split}
\end{equation}

\medskip
{\textbf{3-3}}. Define
\begin{equation*}
M_{1}:=\frac{|\vphi_{\bm\nu}|}{c(|D\vphi|^2, \vphi)},\quad
M_{2}:=\frac{|\vphi_{\bm\tau}|}{c(|D\vphi|^2, \vphi)}.
\end{equation*}
We substitute the expression of $D^2\vphi[{\bm\tau}, D\vphi]$ given by \eqref{from-pmaxassump-shock} into the
right-hand side of \eqref{tangentderiv2-RH-new} to obtain
\begin{equation}\label{pmaxassump2}
\begin{split}
  &D^2\vphi[\bm\tau, \rho D(\ivphi-\vphi)+D\ivphi]\\
  &=\rho(1+(\rho-1)M_{1}^2) B_1
  +\rho(\rho-1)M_{1}^2\vphi_{\bm\tau}
  +(\rho-1)\vphi_{\bm\tau}=:B_2\qquad\tx{at $\Pmax$}.
  \end{split}
\end{equation}
A direct computation shows that
\begin{equation}
\label{expression-B2-new}
B_2=\frac{\big(2(\rho-1)(1+\rho M_{1}^2)-(\gam-1)M^2\big)\vphi_{\bm\tau}
   -c^2\rho\big(1+(\rho-1)M_{1}^2\big)k_{\bm\tau}}{2+(\gam-1)M^2}.
\end{equation}

Apply $\alp$ and $\zeta$ from {Step 2}, and assume that
\begin{equation}
\label{assumption-mu}
  0<\mu\le \min\big\{\mu_1, \frac{\zeta}{2C_*}\big\}.
\end{equation}
Then it follows from \eqref{M-at-pmax} and Step 2 that
\begin{equation}
\label{pos-tan-deriv}
0<\alp |\vphi_{\bm\nu}(\Pmax)|^2< |\vphi_{\bm\tau}(\Pmax)|^2,\quad\tx{or equivalently},\quad
0<\alp M_1^2(\Pmax)<M_2^2(\Pmax).
\end{equation}
Using \eqref{M-at-pmax}, \eqref{pos-tan-deriv}, and $\alp\in (0, \frac 12)$, we have
\begin{equation}
\label{inequality-estimate1}
  M_2^2(\Pmax)>\frac{\alp}{2}(1-C_*\mu).
\end{equation}

We rewrite \eqref{from-pmaxassump-shock} and \eqref{pmaxassump2}
as the following linear system
for $(\vphi_{{\bm\nu}{\bm\tau}},\vphi_{{\bm\tau}{\bm\tau}})$:
\begin{equation*}
  {\bf A}\begin{pmatrix}
     \vphi_{{\bm\nu}{\bm\tau}}\\
    \vphi_{{\bm\nu}{\bm\tau}}
  \end{pmatrix}=\begin{pmatrix}
     B_1\\
    B_2
  \end{pmatrix}
  \qquad\tx{at $\Pmax\,$ for ${\bf A}=\begin{pmatrix}
                                                  \vphi_{\bm\nu} & \vphi_{\bm\tau} \\
                                                  \rho^2 \vphi_{\bm\nu} & \vphi_{\bm\tau}
                                                \end{pmatrix}$}.
\end{equation*}
By \eqref{3-c7} and \eqref{pos-tan-deriv},
$|\det{\bf A}|=|(\rho^2-1)\vphi_{\bm\nu}\vphi_{\bm\tau}|>0$ at $\Pmax$.
Thus, $(\vphi_{{\bm\nu}{\bm\tau}},\vphi_{{\bm\nu}{\bm\tau}})$ can be written as
\begin{equation}
\label{lwrord-express-2deriv1}
\vphi_{\bm\nu\bm\tau}=\frac{B_1-B_2}{(1-\rho^2)\vphi_{\bm\nu}},\quad
\vphi_{\bm\tau\bm\tau}=\frac{\rho^2 B_1-B_2}{(\rho^2-1)\vphi_{\bm\tau}}\qquad\,\,\tx{at $\Pmax$.}
\end{equation}
Note that Eq. \eqref{2-1} is invariant under a coordinate rotation.  We rewrite Eq. \eqref{2-1} as
\begin{equation}
\label{eqn-ortheg-trans}
(c^2-\vphi_{\bm\nu}^2)\vphi_{\bm\nu\bm\nu}-2\vphi_{\bm\nu}\vphi_{\bm\tau} \vphi_{\bm\nu\bm\tau}
+(c^2-\vphi^2_{\bm\tau})\vphi_{\bm\tau\bm\tau}=|D\vphi|^2-2c^2
\qquad\tx{in $\overline{\Omega}\setminus (\overline{\leftsonic}\cup\overline{\rightsonic})$},
\end{equation}
and use this to express $\vphi_{{\bm\nu}{\bm\nu}}$ in terms of $(M_{\bm\nu}, M_{\bm\tau}, M, \rho,\vphi_{\bm\nu \bm\tau}, \vphi_{\bm\tau\bm\tau})$.
Then we use \eqref{lwrord-express-2deriv1} to obtain
\begin{equation}
\label{lwrord-express-2deriv2}
  \begin{split}
     \vphi_{\bm\nu\bm\nu}=&\, \frac{M^2-2}{1-M_{1}^2}
     -\frac{1}{1-M_{1}^2}
     \left(\frac{2M_{1}M_{2}}{(\rho^2-1)\vphi_{\bm\nu}}
        +\frac{\rho^2(1-M_{2}^2)}{(\rho^2-1)\vphi_{\bm\tau}}\right)B_1\\
       &\, +\frac{1}{1-M_{1}^2}\left(
       \frac{1-M_{2}^2}{(\rho^2-1)\vphi_{\bm\tau}}+
       \frac{2M_{1}M_{2}}{(\rho^2-1)\vphi_{\bm\nu}}
       \right)B_2\quad\quad\quad\qquad \tx{at $\Pmax$}.
  \end{split}
\end{equation}
Using \eqref{from-pmaxassump-shock}, \eqref{expression-B2-new}--\eqref{lwrord-express-2deriv1},
and \eqref{lwrord-express-2deriv2}, we can also express $(\vphi_{\bm\nu\bm\tau}, \vphi_{\bm\nu\bm\nu})$
in terms of $M, M_1$, $M_2$, $\rho$, $\vphi_{\bm\tau}, \vphi_{\bm\nu}, c$, and $k_{\bm\tau}$ at $\Pmax\in \shock$.

\smallskip
{\textbf{3-4}}. Now we choose a constant $\mu\in(0, \mu_1]$ sufficiently small so that a contradiction is derived.

By \eqref{expression-Mach-deriv}, \eqref{deriv-comp-2} can be written as
\begin{equation*}
  \big(2+(\gam-1)M^2\big)
  (\vphi_{\bm\tau}\vphi_{\bm\nu\bm\tau}
  +\vphi_{\bm\nu}\vphi_{\bm\nu\bm\nu})+(\gam-1)M^2\vphi_{\bm\nu}+c^2k_{\bm\nu}
  \le 0\qquad\tx{at $\Pmax$}.
\end{equation*}
Using \eqref{definition-kfunctiin-new} and the expressions of $(\vphi_{\bm\nu\bm\tau}, \vphi_{\bm\nu\bm\nu})$
in terms of $M, M_1, M_2, \rho, \vphi_{\bm\tau}, \vphi_{\bm\nu}, c$, and $k_{\bm\tau}$,
we can further rewrite the inequality stated above as
\begin{equation}
\label{ineq-max}
\Delta:= 2M_2^2+2(2\rho+1)M_1^2(M^2-1)
+ \mu \left( l_1 \der_{\bm\nu}g_{\beta}-
l_2 \der_{\bm\tau}g_{\beta} \right)\le 0\qquad\,\tx{at $\Pmax$}
\end{equation}
for
\begin{equation*}
l_1=c^2\vphi_{\bm\tau}(\rho+1)(1-M_1^2),\qquad
l_2=c^2\frac{(1-\rho^2)M_1^2M_2^2+\rho^2M_1^2+M_2^2}{(\rho+1)\vphi_{\bm\tau}}.
\end{equation*}

By \eqref{M-at-pmax}, Lemma \ref{lemma-step3-1},  and the definition of $g_{\beta}$ given in \eqref{8-13},
there exists a constant $C>0$ such that
\begin{equation}\label{inequality-estimate2}
\begin{split}
  2(2\rho+1)M_1^2(M^2-1)\ge -C\mu \qquad&\tx{at $\Pmax$},\\
  |l_1|\le C\qquad&\tx{on $\shock$},\\
  \|Dg_{\beta}\|_{C^{0}(\R^2)}\le C \qquad&\tx{for all $\beta\in[0,\betadet]$}.
\end{split}
\end{equation}
Moreover, by Lemma \ref{lemma-step3-1} and \eqref{inequality-estimate1}, we have
\begin{equation}\label{inequality-estimate3}
  |l_2|\le \frac{1}{\sqrt{\alp(1-C_*\mu)}}\qquad\tx{at $\Pmax$}.
\end{equation}
From \eqref{inequality-estimate1}--\eqref{inequality-estimate3}, we obtain
\begin{equation*}
  \Delta\ge \alp(1-C_*\mu)-C\mu\Big(1+\frac{1}{\sqrt{\alp(1-C_*\mu)}}\Big)
  \qquad\tx{at $\Pmax$}
\end{equation*}
for some constant $C>0$, provided that $\mu$ satisfies \eqref{assumption-mu}.
Therefore, there exists a constant $\mu_2\in(0, \mu_1^*]$ for $\displaystyle{\mu_1^*=\min\{\mu_1, \frac{\zeta}{2C_*}\}}$
such that, if $0<\mu\le \mu_2$, then $\Delta> \frac{\alp}{8}>0$ holds at $\Pmax$,
which contradicts \eqref{ineq-max}.
Therefore, we conclude that the maximum of $M^2_{\vphi}+\mu g_{\beta}$ over $\ol{\Om}$ must be attained
on $\der \Om\setminus (\Wedge\cup\shock)$,
provided that $\mu>0$ is chosen sufficiently small, depending only on $(\iv,\gam)$.

\smallskip
{\textbf{4.}}
For constant $\mu_2$ given in Step 3, we fix a constant $\mu\in(0, \mu_2]$. Then
$M^2_{\vphi}+\mu g_{\beta}$ satisfies
\begin{equation*}
\sup_{\ol{\Om}}\left(M_{\vphi}^2+\mu g_{\beta}\right)=\sup_{\leftsonic\cup\rightsonic}\left(M_{\vphi}^2+\mu g_{\beta}\right)=1.
\end{equation*}
This proves \eqref{estimate-M-2015nov}.
\end{proof}
\end{proposition}

\begin{remark}\label{remark-ellipticity-newly-added}
By Remark {\rm \ref{remark-ellipticity-new2015}} and \eqref{estimate-M-2015nov} in Proposition {\rm \ref{corollary-ellip}},
there exists a constant $\mu_{\rm{el}}>0$ depending only on $(\iv, \gam)$ such that,
if $\vphi$ is an admissible solution corresponding to $(\iv, \beta)\in\mathfrak{R}_{\rm{weak}}$,
\begin{equation}\label{distfl-ineq}
M_{\vphi}^2(\bmxi)\le 1-\mu_{\rm{el}} {\rm dist}^{\flat}(\bmxi,\leftsonic\cup\rightsonic)\qquad\tx{in $\ol{\Om}$.}
\end{equation}
\end{remark}

\section{Uniform Weighted $C^{2,\alp}$--Estimates Away From $\leftsonic$}
\label{subsec-unif-c2est-super}

According to Proposition \ref{corollary-ellip}, the ellipticity of Eq. \eqref{8-48-ps} (or equivalently, Eq. \eqref{2-1})
depends on ${\rm dist}(\bm\xi, \leftsonic\cup\rightsonic)$. In particular,
\eqref{def-distfl} indicates that the ellipticity of \eqref{8-48-ps} depends  continuously on $\beta\in(0,\betadet)$,
even across $\betac^{(\iv)}$ up to $\betadet$. For this reason, we can establish uniform weighted
$C^{2, \alp}$--estimates of admissible solutions.

We first estimate (weighted) $C^{2,\alp}$--norms of admissible solutions away from $\leftsonic$.
We will discuss the uniform  (weighted) $C^{2,\alp}$--estimates of admissible solutions
near $\leftsonic$ in \S \ref{subsec-unif-c2est-nr-leftsonic}.

\subsection{\bf $C^{2,\alp}$--estimates away from  $\leftsonic\cup\rightsonic$}
\label{subsubsec-est-away}
Fix $\gam\ge 1$ and $\iv>0$.

For a set $U\subset \R^2$ and a constant $\eps>0$, define $\mcl{N}_{\eps}(U):=\{{\bmxi}\in \R^2:{\rm dist}(\bmxi, U)<\eps\}$.

Let $C>0$ be the constant from Proposition \ref{proposition-distance}.
Then there exists a constant $d_0>0$ depending only on $(\iv, \gam)$ such that
\begin{equation}
\label{ivphi-d0-estimate}
|D\ivphi|^2\ge 1+d_0\qquad\tx{on}\,\,\mcl{N}_{\frac{1}{2C}}(\shock).
\end{equation}

(i) If $\gam=1$, then it follows directly from Definition \ref{def-regular-sol}
that any admissible solution $\vphi$ satisfies that $|D\vphi|\le 1$ in $\ol{\Om}$.
Thus, it follows from \eqref{ivphi-d0-estimate} that
\begin{equation}
\label{imp-obs1}
|D\ivphi|^2-|D\vphi|^2\ge d_0\qquad\tx{on}\,\, \mcl{N}_{\frac{1}{2C}}(\shock)\cap \overline{\Om}.
\end{equation}

(ii) If $\gam>1$, then we can rewrite the Bernoulli law \eqref{new-density} as
\begin{equation}
\label{Br-law}
\rho^{\gam-1}+\frac{\gam-1}{2}\big(|D\vphi|^2+2\vphi\big)=1+\frac{\gam-1}{2}\big(|D\ivphi|^2+2\ivphi\big).
\end{equation}
Let $\vphi$ be an admissible solution corresponding to $(\iv, \beta)\in \mathfrak{R}_{\rm weak}$.
Since $|D\vphi|^2\le \rho^{\gam-1}$ and $\ivphi-\vphi \ge 0$ hold in $\ol{\Om}$,
we obtain from \eqref{ivphi-d0-estimate} and \eqref{Br-law} that
\begin{equation*}
\frac{\gam+1}{2}\rho^{\gam-1}\ge  \rho^{\gam-1}+\frac{\gam-1}{2}|D\vphi|^2\ge 1+\frac{\gam-1}{2}(1+d_0) \qquad\tx{on $\mcl{N}_{\frac{1}{2C}}(\shock)\cap \ol{\Om}.$}
\end{equation*}
This implies that $\rho^{\gam-1}-1\ge \delta_0$ for some constant $\delta_0>0$ depending only on $(\iv, \gam)$.
Then
\begin{equation*}
|D\ivphi|^2-|D\vphi|^2=\frac{2(\rho^{\gam-1}-1)}{\gam-1}+(\ivphi-\vphi)
\ge \frac{2\delta_0}{\gam-1}+(\ivphi-\vphi)
\qquad\tx{on $\mcl{N}_{\frac{1}{2C}}(\shock)\cap \ol{\Om}.$}
\end{equation*}
Since $\ivphi-\vphi=0$ on $\shock$, it follows from \eqref{3-c6} in Lemma \ref{lemma-step3-1} that
there exist small constants $\eps\in(0, \frac{1}{4C})$ and $\delta_0'>0$ depending only on $(\gam, \iv)$  such that
\begin{equation}
\label{imp-obs2}
|D\ivphi|-|D\vphi|\ge \delta_0' \qquad\tx{on $\mcl{N}_{\eps}(\shock)\cap \ol{\Om}.$}
\end{equation}

Let $(r,\theta)$  be the polar coordinates defined by \eqref{def-polar-oi}.
Note that $|D\ivphi|=-\der_r\ivphi$.
Then \eqref{imp-obs1} and \eqref{imp-obs2} imply that there exists a constant $d_1>0$ depending only on $(\iv, \gam)$ such that
\begin{equation}
\label{r-derivative}
\der_r(\ivphi-\vphi)\le -( |D\ivphi|-|D\vphi|)\le -d_1\qquad{\text{on}\,\,\, \mcl{N}_{\eps}(\shock)\cap\ol{\Om}}.
\end{equation}

Therefore, by the implicit function theorem, there exists a unique function $\fshpolar(\theta)$ such that
\begin{equation}
\label{graph-shock}
\shock=\{r=\fshpolar(\theta), \,\,\theta_{\righttop}<\theta<\theta_{\lefttop}\},
\end{equation}
where $(\fshpolar(\theta_{P_j}), \theta_{P_j})$ represent the $(r,\theta)$--coordinates of
points $P_j$ for $j=1,2$, given by Definition \ref{definition-domains-np}.
By Lemma \ref{lemma-step3-1} and \eqref{r-derivative}, there exists a constant $C_1$ depending only on $(\iv, \gam)$
such that
\begin{equation}
\label{8-15}
\|\fshpolar\|_{C^{0,1}([\theta_{\righttop},\theta_{\lefttop}])}\le C_1.
\end{equation}

\begin{lemma}
\label{corollary-comp2}
Fix $\gam\ge  1$ and $\iv>0$.
There exists a constant $\delta_1>0$ depending only on $(\iv, \gam)$ such that, if $\vphi$ is an admissible solution
corresponding to $(\iv, \beta)\in \mathfrak{R}_{\rm weak}$, then
\begin{align}
\label{6-b2}
&\der_{\bm\nu}(\ivphi-\vphi)>\delta_1 \;\;\qquad\; \tx{on $\ol{\shock}$},\\
\label{6-b6}
&\der_{\bm\nu}\ivphi> \der_{\bm\nu}\vphi\ge \delta_1\;\qquad\; \tx{on $\ol{\shock}$}
\end{align}
for the unit normal vector ${\bm\nu}=\frac{D(\ivphi-\vphi)}{|D(\ivphi-\vphi)|}$ to $\shock$ towards the interior of $\Om$.

\begin{proof}
If $\vphi$ is an admissible solution corresponding to $(\iv, \beta)$,
then it follows from  \eqref{r-derivative} and $\ivphi-\vphi=0$ on $\shock$ that
\begin{equation}
\label{definition-d1-new2015}
\der_{\bm\nu}(\ivphi-\vphi)=|D(\ivphi-\vphi)|\ge |D\ivphi|-|D\vphi|\ge d_1\qquad\,\, \tx{on}\,\, \ol{\shock}.
\end{equation}
Since $\der_{\bm\nu}\vphi=\frac{\der_{\bm\nu} \ivphi}{\rho(|D\vphi|^2,\vphi)}$,
$\der_{\bm\nu}\ivphi>1$, and $\rho(|D\vphi|^2,\vphi)>1$ on $\ol{\shock}$,
Lemma \ref{lemma-step3-1} yields that $\der_{\bm\nu}\ivphi>\der_{\bm\nu}\vphi \ge C^{-1}$ for
a constant $C>0$ depending only on $(\iv, \gam)$. The proof is completed by choosing $\delta_1$ as
\begin{equation*}
  \delta_1=\min\{d_1, C^{-1}\}.
\end{equation*}
\end{proof}
\end{lemma}

\begin{lemma}\label{lemma-unif-est2}
Fix $\gam\ge 1$ and $\iv>0$.
Let $\vphi$ be an admissible solution corresponding to $(\iv, \beta)\in \mathfrak{R}_{\rm weak}$.
Then, for each $d>0$ and $k=2,3,\cdots$, there exist constants $s, C_k>0$ depending only on $(\iv, \gam, d)$
such that, if $P=(r_P,\theta_P)\in \shock$ in the $(r,\theta)$--coordinates, defined by \eqref{def-polar-oi},
satisfies that ${\rm dist}(P, \leftsonic\cup\rightsonic)\ge d$, then
\begin{equation}
\label{est-near-shock}
|D^k\fshpolar(\theta_P)|\le C_k,\qquad\,\,\,\,
|D^k \vphi|\le C_k\,\,\,\,\,\, \tx{in}\,\,B_s(P)\cap \Om.
\end{equation}

\begin{proof} The proof is divided into three steps.

\smallskip
{\textbf{1.}}
Let $\vphi$ be an admissible solution corresponding to $(\iv, \beta)\in \mathfrak{R}_{\rm weak}$,
and let $\Om$ be its pseudo-subsonic region.
For a constant $d>0$, define
\begin{equation*}
\ol{\Om}_{d}:=\{\bmxi\in\ol{\Om}\,:\,{\rm dist}(\bmxi,\leftsonic\cup\rightsonic)> \frac d2\}.
\end{equation*}

Let $\mcl{E}(\vphi, \ol{\Om}_{d})$ be defined by \eqref{8-50}.
Moreover, for a constant $R$, let $\mcl{K}_R$ be given by \eqref{definition-kmn}.
By Lemma \ref{lemma-step3-1} and Proposition \ref{corollary-ellip},
there exists a constant $M_{d}>0$ depending only on $(\iv, \gam, d)$ such that
$\mcl{E}(\vphi, \ol{\Om}_{d})$ is contained in $\mcl{K}_{M_{d}}$.

Let $\mcl{A}({\bf p},z)=(\mcl{A}_1,\mcl{A}_2)({\bf p},z)$ and $\mcl{B}({\bf p}, z)$ be defined
by \eqref{8-49-ps}, and let $(\til{\mcl{A}}, \til{\mcl{B}})({\bf p},z)$ be the extensions
of $({\mcl{A}}, \mcl{B})({\bf p},z)$ onto $\R^2\times \R$ described in Lemma \ref{lemma-extcoeff-appc2-new} with $M=M_d$.

\smallskip
{\textbf{2.}} We  express the Rankine-Hugoniot jump condition: $\rho D\vphi\cdot{\bm \nu}=D\ivphi\cdot{\bm\nu}$ as
\begin{equation}
\label{rhbc-gsh}
g^{\rm sh}(D\vphi, \vphi, \bmxi)=0\qquad\,\,\tx{on $\shock$}
\end{equation}
for $g^{\rm sh}({\bf p}, z, \bmxi)$ defined by
\begin{equation}
\label{def-gsh}
g^{\rm sh}({\bf p}, z, \bmxi)=\big(\mcl{A}({\bf p}, z)-D\ivphi(\bmxi)\big)\cdot \frac{D\ivphi(\bmxi)-{\bf p}}{|D\ivphi(\bmxi)-{\bf p}|}.
\end{equation}

For $\delta_1>0$ from Lemma \ref{corollary-comp2}, define a smooth function $\zeta\in C^{\infty}(\R)$ by
\begin{equation*}
\zeta(t)=\begin{cases}
t& \tx{on $t\ge \frac{3}{4}\delta_1$},\\
\frac{\delta_1}{2}& \tx{for $t<\frac{\delta_1}{2}$},
\end{cases}\qquad \quad
\zeta'(t)\ge 0\,\,\,\, \tx{on $\R$}.
\end{equation*}
Also, we define an extension of $g^{\rm sh}_{\rm mod}({\bf p}, z, \bmxi)$ onto $\R^2\times \R\times \ol{\Om}_{d}$ by
\begin{equation}
\label{gsh-ext}
g^{\rm sh}_{\rm mod}({\bf p}, z, \bmxi)
=\big(\til{\mcl{A}}({\bf p},z)-D\ivphi(\bmxi)\big)\cdot \frac{D\ivphi(\bmxi)-{\bf p}}{\zeta(|D\ivphi(\bmxi)-{\bf p}|)}.
\end{equation}
Fix a point $P\in \shock$ with ${\rm dist}(P,\leftsonic\cup\rightsonic)>2d$ for $d>0$. Then $\vphi$ satisfies
\begin{equation}
\label{new-mod-nlbvp}
\begin{split}
&{\rm div} \til{\mcl{A}}(D\vphi, \vphi)+\til{\mcl{B}}(D\vphi, \vphi)=0\qquad\tx{in}\,\, B_{d/2}(P)\cap \Om,\\
&g^{\rm sh}_{\rm mod}(D\vphi, \vphi, \bmxi)=0
\qquad\qquad\qquad\,\,\,\tx{on}\,\, B_{d/2}(P)\cap \shock.
\end{split}
\end{equation}

For $\eps>0$ from \eqref{r-derivative}, define
\begin{equation*}
  R:=\min\{\frac d2, \eps\}.
\end{equation*}
Note that such a constant $R>0$ is given depending only on $(\iv, \gam, d)$, but independent of $\vphi$ and $P$.
By \eqref{r-derivative}, we can write
$D_{{\bf p}}g^{\rm sh}_{\rm mod}(D\vphi,\vphi, \bmxi)$ as
\begin{equation*}
  D_{{\bf p}}g^{\rm sh}_{\rm mod}(D\vphi,\vphi, \bmxi)
  =D_{{\bf p}}\big((\mcl{A}({\bf p}, z, \bmxi)-D\ivphi(\bmxi))\cdot \hat{\bf n}({\bf p},\bmxi)\big)
  \qquad\tx{in $B_{R}(P)\cap \ol{\Om}$}
\end{equation*}
for
\begin{equation*}
  \hat{\bf n}({\bf p},\bmxi)=\frac{D\ivphi(\bmxi)-{\bf p}}{|D\ivphi(\bmxi)-{\bf p}|}.
\end{equation*}
Since
\begin{equation*}
\hat{\bf n}({\bf p},\bmxi)\cdot \big((\mcl{A}({\bf p}, z, \bmxi)-D\ivphi(\bmxi))D_{{\bf p}}\hat{\bf n}({\bf p},\bmxi)\big)
=\frac{1}{2}(\mcl{A}({\bf p}, z, \bmxi)-D\ivphi(\bmxi))\cdot D_{\bf p}(|\hat{\bf n}({\bf p},\bmxi)|^2)=0,
\end{equation*}
a direct computation yields that
\begin{equation*}
   D_{{\bf p}}g^{\rm sh}_{\rm mod}(D\vphi,\vphi, \bmxi)\cdot \hat{\bf n}(D\vphi,\bmxi)=\sum_{i,j=1}^2\der_{p_i}\til{\mcl{A}}_j(D\vphi,\vphi,\bm\xi)
   \hat{n}_i\hat{n}_j=:\hat{\lambda}(D\vphi, \bmxi)\qquad\tx{in $B_{R}(P)\cap \ol{\Om}$}
\end{equation*}
for $\hat{n}_i=\hat{\bf e}_i\cdot \hat{\bf n}(D\vphi,\bmxi)$.

By Lemma \ref{lemma-extcoeff-appc2-new}(ii), there exists
a constant $\lambda_d>0$ depending only on $(\iv, \gam, d)$ such that
\begin{equation*}
  D_{{\bf p}}g^{\rm sh}_{\rm mod}(D\vphi,\vphi, \bmxi)\cdot \hat{\bf n}(D\vphi,\bmxi)\ge \lambda_d>0
  \qquad\tx{in $B_{R}(P)\cap \ol{\Om}$}.
\end{equation*}
This implies that
\begin{equation}
\label{nondegeneracy-deriv-gshockmod}
|D_{\bf p} g^{\rm sh}_{\rm mod}(D\vphi,\vphi,\bmxi)|\ge
D_{{\bf p}}g^{\rm sh}_{\rm mod}(D\vphi,\vphi, \bmxi)\cdot \hat{\bf n}(D\vphi,\bmxi)\ge \lambda_d>0 \qquad\tx{in $B_{R}(P)\cap \ol{\Om}$}.
\end{equation}

\smallskip
\textbf{3.} By estimate \eqref{3-c6} of Lemma \ref{lemma-step3-1}, \eqref{8-15}, Lemma \ref{lemma-extcoeff-appc2-new},
and \eqref{nondegeneracy-deriv-gshockmod},  the boundary value problem \eqref{new-mod-nlbvp} satisfies all the conditions necessary
to apply Theorem \ref{elliptic-t3-CF2}. Therefore, there exist $\beta\in(0,1)$ and $C>0$ depending only on $(\iv, \gam, d)$ such that
\begin{equation*}
\|\vphi\|_{1,\beta, B_{d/4}(P)\cap \Om}\le C\qquad\,\, \tx{for all $P\in \shock\cap \ol{\Om}_d$}.
\end{equation*}

Combining the $C^{1,\beta}$--estimate of $\vphi$ with \eqref{r-derivative} implies that $\fshpolar$ is $C^{1,\beta}$
away from $\theta=\theta_{\lefttop}, \theta_{\righttop}$.
Then we apply Theorem \ref{elliptic-t4-CF2} to the boundary value problem \eqref{new-mod-nlbvp} to obtain the estimate:
\begin{equation*}
\|\vphi\|_{2,\beta, B_{d/8}(P)\cap \Om}\le C\qquad\tx{for all}\,\,P\in \shock\cap \ol{\Om}_d
\end{equation*}
for some constant $C>0$ depending only on $(\iv, \gam, d)$.
This implies that $\fshpolar$ is $C^{1,\alp}$ for any $\alp\in(0,1)$ away from
$\theta=\theta_{\lefttop}, \theta_{\righttop}$,
so that $\vphi$ is $C^{2,\alp}$ for any $\alp\in(0,1)$ on $\shock$ away
from $\ol{\leftsonic}\cup\ol{\rightsonic}$ by Theorem \ref{elliptic-t4-CF2}.

Finally, the $C^k$--estimates, $k=2,3,\cdots$, are obtained
by a bootstrap argument
via application of Theorem \ref{elliptic-t4-CF2} and Corollary \ref{elliptic-t5-CF2}.
\end{proof}
\end{lemma}

As a result, directly following from Lemmas \ref{lemma-unif-est1} and \ref{lemma-unif-est2},
we conclude the following uniform $C^k$--estimates of admissible solutions:
\begin{corollary}\label{corollary-unif-est-away-sn}
Fix $\gam\ge 1$ and $\iv>0$. For each $d>0$ and $k=2,3,\cdots$, there exists a constant $C_{k,d}>0$ depending only on $(\iv, \gam, k, d)$
such that any admissible solution $\vphi$ corresponding to $(\iv, \beta)\in \mathfrak{R}_{\rm weak}$ satisfies
\begin{equation*}
\|\vphi\|_{k, \ol{\Om\cap\{{\rm dist}(\bmxi, \leftsonic\cup\rightsonic)>d\}}}\le C_{k,d}.
\end{equation*}
\end{corollary}

\subsection{\bf $C^{2,\alp}$--estimates near $\rightsonic$}\label{subsubsec-sonic-est}
For fixed $\gam\ge 1$ and $\iv>0$, the sonic arc $\rightsonic$, defined by Definition \ref{definition-domains-np}
corresponding to the normal shock part of each admissible solution,
is fixed to be the same for all $\beta\in(0, \frac{\pi}{2})$.
By Definition \ref{def-regular-sol}(ii) and Proposition \ref{corollary-ellip}, the ellipticity of
Eq. \eqref{8-48-ps} (or equivalently, Eq. \eqref{2-1}) degenerates near $\rightsonic$.
In order to establish a uniform weighted $C^{2,\alp}$--estimate of admissible solutions
up to $\rightsonic$, the method of parabolic scaling is employed.
We keep following Definition \ref{definition-domains-np} for the notations used hereafter.

Define
 \begin{equation}
 \label{definition-rightch-new}
 \rightch:=\frac{\rightc+\neta}{2},
\end{equation}
which
is the same for all $\beta\in[0, \frac{\pi}{2})$.
In $U_{\mcl{N}}:=\bigl(B_{\frac{3\rightc}{2}}(\Onormal)\setminus B_{\rightch}(\Onormal)\bigr)\cap \{\xxi\,:\,\xin>0\}$,
let $(r,\theta)$ be the polar coordinates with respect to $\Onormal=(0,0)$. Define
    \begin{equation}
    \label{coord-n}
    (x,y):=(\rightc-r,\theta).
    \end{equation}
Let $\vphi$ be an admissible solution corresponding to $(\iv, \beta)\in \mathfrak{R}_{\rm weak}$,
and let $\Om$ be its pseudo-subsonic region.  Define
\begin{equation}\label{definition-nOm-new}
  \nOm:=\bigl(\Om\cap\{\xin>0\}\bigr)\setminus B_{\rightch}(\Onormal).
\end{equation}
Then $\nOm\subset B_{\rightc}(\Onormal)$
and $\nOm\subset\{(x,y)\,:\,x>0\}$.

In $\nOm$, we define a function $\psi$ by
\begin{equation}
\label{psi-def-right}
\psi:=\vphi-\rightvphi\qquad\quad\text{in}\;\;\nOm.
\end{equation}
We rewrite Eq. \eqref{2-1} and the boundary conditions \eqref{3-a1}--\eqref{3-a3} in the $(x,y)$--coordinates as follows:
\smallskip

\emph{{\rm (i)} Equation for $\psi$ in $\nOm$}:
For each $j=1,\cdots, 5$, define $\mfrakO_j^{\mcl{N}}(\pb, z, x)$ by
\begin{equation*}
\mfrakO_j^{\mcl{N}}(\pb, z, x):=O_j(\pb, z, x, \rightc)
\end{equation*}
for $O_j(\pb, z, x, c)$ given by \eqref{prelim5-5}. Then Eq. \eqref{2-1} is written as
\begin{equation}
\label{eqn-xy-right}
\big(2x-(\gam+1)\psi_x+\mfrakO_1^{\mcl{N}}\big)\psi_{xx}+\mfrakO_2^{\mcl{N}}\psi_{xy}
+\big(\frac{1}{\rightc}+\mfrakO_3^{\mcl{N}}\big)\psi_{yy}-\big(1+\mfrakO_4^{\mcl{N}}\big)\psi_x+\mfrakO_5^{\mcl{N}}\psi_y=0,
\end{equation}
with $O_j^{\mcl{N}}=O_j^{\mcl{N}}(D\psi, \psi, x)$ for $j=1,\cdots, 5$.

\smallskip
\emph{{\rm (ii)} Boundary condition for $\psi$ on $\shock \cap \der \nOm$}:
By the definitions of $(\ivphi, \rightvphi)$ given in Definition \ref{definition-domains-np},
we rewrite the condition that $\ivphi-\vphi=0$ on $\shock \cap \partial \Omega^{\mcl{N}}$
as
\begin{equation*}
\etan=\neta-\frac{\psi}{\iv}\qquad\tx{on $\shock \cap \partial \Omega^{\mcl{N}}$}.
\end{equation*}
For $g^{\rm sh}_{\rm mod}({\bf p}, z, \bmxi)$ given by \eqref{gsh-ext}, we define
\begin{equation}
\label{7-d7-pre}
M({\bf p}, z, \xin):=g^{\rm sh}_{\rm mod}({\bf p}+D\rightvphi,z+\rightvphi, \xin, \neta-\frac{z}{\iv})
\end{equation}
with $(D\rightvphi, \rightvphi)$ evaluated at $(\xin, \neta-\frac{z}{\iv})$.
Then the boundary condition \eqref{3-a3} is written as $M(D\psi, \psi, \xin)=0$ on $\shock$.
Denote
\begin{equation*}
\iphi^{\mcl{N}}:=\ivphi-\rightvphi.
\end{equation*}
Then
$|D(\iphi^{\mcl{N}}-\psi)|=|\der_{\bm\nu}(\ivphi-\vphi)|>0$ on $\shock$.
Rewriting the boundary condition
$|D(\iphi^{\mcl{N}}-\psi)|M(D\psi, \psi, \xin)=0$ on $\shock\cap \der\nOm$ in the $(x,y)$--coordinates,
we obtain
\begin{equation}
\label{bc-psi-B1-right}
\mcl{B}^{\mcl{N}}_1(\psi_x, \psi_y, \psi, x, y)=0\qquad\,\,\tx{on $\shock\cap\der{\nOm}$}
\end{equation}
for $\mcl{B}^{\mcl{N}}_1(p_x, p_y, z, x, y)$ defined by
\begin{equation}
\label{def-B1-new-2015-N}
\mcl{B}^{\mcl{N}}_1(p_x, p_y, z, x, y)
:={|D\iphi^{\mcl{N}}-\pb|}M(\pb, z, \xin)
\end{equation}
with
\begin{equation}
\label{cov-right}
\xin=(\rightc-x)\cos y,\qquad
\pb=
\begin{pmatrix}
-\cos y&- \sin y\\
-\sin y& \cos y
\end{pmatrix}
\begin{pmatrix}
p_x\\
\frac{p_y}{\rightc-x}
\end{pmatrix}.
\end{equation}

\smallskip
\emph{{\rm (iii)} Other properties of $\psi$}:
By \eqref{1-l} and Definition \ref{def-regular-sol}(ii)--(iv),
$\psi$ satisfies
\begin{equation}
\label{5-a4}
\begin{split}
&\psi \ge 0\qquad\,\,\text{in}\,\,\,\nOm,\\
&\psi=0\qquad\,\,\text{on}\,\,\rightsonic,\\
&\psi_y=0\qquad\text{on}\,\,\,\Wedge\cap\der{\nOm}.
\end{split}
\end{equation}

For each $\beta\in[0, \frac{\pi}{2})$, let $\mcl{D}$ be defined by \eqref{def-D-domain}, and define
 \begin{equation*}
 \label{definition-rightch}
\nLambda:=\mcl{D}\cap\bigl(B_{\frac{3\rightc}{2}}(\Onormal)\setminus B_{\rightch}(\Onormal)\bigr)\cap\{\xin>0\}.
\end{equation*}
Note that $\nLambda$ is the same for all $\beta\in[0, \frac{\pi}{2})$, and $\nLambda \subset \{\etan<\neta\}$.

By using the definitions of $(\rightsonic, \ivphi, \rightvphi)$  given in Definition \ref{definition-domains-np},
the following lemma can  be directly verified:

\begin{lemma}\label{lemma1-sonic-N}
Fix $\gam\ge 1$ and $\iv>0$.
There exist positive constants $\eps_1$, $\eps_0$, $\delta_0$, $\omega_0$, $C$, and $\mathfrak{M}$
depending only on $(\iv, \gam)$ with $\eps_1>\eps_0$ and $\mathfrak{M}\ge 2$
so that the following properties hold{\rm :}

\smallskip
\begin{itemize}
\item[(a)]
$\{\rightvphi<\ivphi\}\cap\nLambda\cap\mcl{N}_{\eps_1}(\rightsonic)
\subset \{0<y<\frac{\pi}{2}-\delta_0\}$,
where $\mcl{N}_{\eps}(\Gam)$ denotes the $\eps$--neighborhood of a set $\Gam$ in the $\xxi$--coordinates{\rm ;}

\smallskip
\item[(b)]
    $\{\rightvphi<\ivphi\}\cap\mcl{N}_{\eps_1}(\rightsonic)\cap\{y>y_{\righttop}\}\,
    \subset\{x>0\}${\rm ;}

\smallskip
\item[(c)] In $\{(x,y)\,:\,|x|<\eps_1,\,\, 0<y<\frac{\pi}{2}-\delta_0\}$,
$\iphi^{\mcl{N}}=\ivphi-\rightvphi$ satisfies
\begin{equation}
\label{6-a5-right}
\frac {2}{\mathfrak{M}}y \le \der_x\iphi^{\mcl{N}}(x,y)\le \frac {\mathfrak{M}}{2},
\qquad \frac {2}{\mathfrak{M}}\le -\der_y\iphi^{\mcl{N}} \le \frac {\mathfrak{M}}{2};
\end{equation}

\item[(d)]
$|(D^2_{(x,y)}, D^3_{(x,y)})\iphi^{\mcl{N}}|\le C$ in $\{|x|<\eps_1\}${\rm ;}

\smallskip
\item[(e)] There exists a unique function $\hat f_{\mcl{N},0}\in C^{\infty}([-\eps_0,\eps_0])$ such that
    \begin{equation}
    \label{est-sonic-1-right}
   \quad\!\!\!
    \begin{cases}
    \{\rightvphi<\ivphi\}\cap \nLambda\cap\mcl{N}_{\eps_1}(\rightsonic)\cap\{|x|<\eps_0\}=\{(x,y)\,:\,|x|<\eps_0, 0<y<\hat f_{\mcl{N},0}(x)\},\\[1mm]
    \rightshock\cap\mcl{N}_{\eps_1}(\rightsonic)\cap\{|x|<\eps_0\}=\{(x,y)\,:\,x\in(-\eps_0,\eps_0), y=\hat f_{\mcl{N},0}(x)\};
    \end{cases}
    \end{equation}

\item[(f)]  $\hat f_{\mcl{N},0}$ in {\rm (e)} satisfies
    \begin{equation*}
    2\omega_0\le \hat f'_{\mcl{N},0}\le C \qquad\text{on $(-\eps_0,\eps_0)$}.
    \end{equation*}
\end{itemize}
\end{lemma}

Let $\Om$ be the pseudo-subsonic region of an admissible solution $\vphi$ corresponding to $(\iv, \beta)\in \mathfrak{R}_{\rm weak}$.
For $\eps\in (0,\eps_1]$, define a set $\nOm_{\eps}$ by
\begin{equation}
\label{6-h1-right}
\nOm_{\eps}:=\Om\cap\mcl{N}_{\hat{\eps}}(\rightsonic)\cap\{x<\eps\}
\end{equation}
for some $\hat{\eps}=\hat{\eps}(\eps,\omega_0)>\eps$.

Note that $\nOm_{\eps}\subset \{0<x<\eps\}$.

\begin{lemma}
\label{lemma-est-nrsonic-N}
Let $\eps_0$, $\omega_0$, and $\mathfrak{M}$ be from Lemma {\rm \ref{lemma1-sonic-N}}.
Then there exist constants $\bar{\eps}\in(0, \eps_0]$, $L\ge 1$, $\delta\in(0,\frac 12)$,
and $\omega\in(0,\omega_0]\cap (0,1)$ depending only on $(\iv, \gam)$ such that,
whenever $\eps\in(0, \bar{\eps}]$, any admissible solution $\vphi=\psi+\rightvphi$
satisfies the following properties in $\nOm_{\eps}${\rm :}

\smallskip
\begin{itemize}
\item[(a)]
$
\psi_x(x,y)\le \frac{2-\delta}{1+\gam} x\le Lx;
$

\smallskip
\item[(b)]
$
\psi_x\ge 0
$
 and
 $
|\psi_y(x,y)|\le Lx;
$

\smallskip
\item[(c)]
$\frac{2}{\mathfrak{M}}y-\frac{2-\delta}{1+\gam}x \le \der_x(\ivphi-\vphi)(x,y)\le {\mathfrak{M}}$
and $\frac{1}{\mathfrak{M}} \le -\der_y(\ivphi-\vphi)\le {\mathfrak{M}};$

\smallskip
\item[(d)]  there exists a unique function $\fshockn\in C^1([0,\eps])$ such that
    \begin{equation*}
    \begin{split}
    &\nOm_{\eps}=\{(x,y)\,:\,x\in(0,\eps), 0<y<\fshockn(x)\},\\
    &\shock\cap\der \nOm_{\eps}=\{(x,y)\,:\,x\in(0,\eps),\,\, y=\fshockn(x)\},\\
    &\omega\le \fshockn'(x) \le L \qquad\text{for}\,\,\, 0<x<\eps;
    \end{split}
    \end{equation*}

 \item[(e)] $0\le \psi(x,y)\le Lx^2$.
\end{itemize}

\begin{proof}
We divide the proof into four steps.

\smallskip
{\textbf{1.}}
By \eqref{9-32}
and \eqref{eqn-xy-right},
there exists a constant $\bar{\delta}\in(0, \frac 14)$ depending only on $(\iv, \gam)$
such that
\begin{equation}
\label{eqn-xy-coeff11}
2x-(\gam+1)\psi_x+O_1^{\mcl{N}}(D\psi(x,y), \psi(x,y),x)\ge 2\bar{\delta} x\,\,\qquad\tx{in}\,\,\nOm
\end{equation}
for $\nOm$ defined by \eqref{definition-nOm-new}.
Since $O_1^{\mcl{N}}(D\psi(x,y), \psi(x,y),x)\le \frac{(\gam+1)}{\rightc}x\psi_x$
by \eqref{prelim5-5} and \eqref{5-a4}, we obtain from \eqref{eqn-xy-coeff11} that
\begin{equation*}
\psi_x(x,y)\le \frac{2-2\bar{\delta}}{(1+\gam)(1-\frac{\bar{\eps}_0}{\rightc})}x
\qquad\,\tx{in $\nOm_{\bar{\eps}_0}$}
\end{equation*}
for
\begin{equation*}
 \bar{\eps}_0=\min\{\rightc-\rightch, \eps_0\},
\end{equation*}
where $\rightch$ is given by \eqref{definition-rightch-new}.
Then $\bar{\eps}\in(0, \eps_0]$ can be chosen, depending only on $(\iv, \gam)$,
so that $\psi$ satisfies
\begin{equation*}
\psi_x(x,y)\le \frac{2-\bar{\delta}}{1+\gam}x\qquad\tx{in}\,\,\nOm_{\bar{\eps}}.
\end{equation*}
This proves statement (a).

By Lemma \ref{lemma-step-2}, \eqref{coord-n}, and \eqref{cov-right}, we have
\begin{equation}
\label{deriv-psi-xy-right}
\psi_x\cos y+\frac{\psi_y}{\rightc-x}\sin y\ge 0,\quad
\psi_x \sin y-\frac{\psi_y}{\rightc-x}\cos y\ge 0\,\qquad\,\,\tx{in $\nOm$}.
\end{equation}
By property (f) of Lemma \ref{lemma1-sonic-N}, there exists a constant $\delta_1\in(0, \frac{\pi}{10})$
depending only on $(\iv, \gam)$ such that
\begin{equation}
\label{y-range-nOm}
\nOm\subset\{0<y<\frac{\pi}{2}-\delta_1\}.
\end{equation}
Then \eqref{deriv-psi-xy-right}, combined with statement (a), yields that
\begin{equation}
\label{estimate-psix-xy-right}
0\le \psi_x(x,y)\le \frac{2-\bar{\delta}}{1+\gam}x\,\qquad \tx{in $\nOm_{\bar{\eps}}$}.
\end{equation}

Owing to \eqref{y-range-nOm}, the second inequality in \eqref{deriv-psi-xy-right} is equivalent to
$$
\psi_y(x,y)\le (\rightc-x)\psi_x(x,y)\tan y \qquad \mbox{in $\nOm$}.
$$
Then it follows directly from \eqref{estimate-psix-xy-right} that
\begin{equation}
\label{upper-bd-psi-y}
\psi_y\le Cx\qquad\tx{in}\,\,\nOm_{\bar{\eps}}
\end{equation}
for a constant $C>0$ chosen depending only on $(\iv, \gam)$.

\smallskip
{\textbf{2.}}
In order to obtain a lower bound of $\psi_y$ by a linear function of $x$ near $\rightsonic$,
a different approach is used.

By Proposition \ref{proposition-sub3} and \eqref{y-range-nOm},
there exists $\delta_1'\in (0, \frac{\pi}{10})$ depending only on $(\iv, \gam)$ such that
\begin{equation}
\label{shock-xy-right}
\der \nOm\cap \shock\, \subset \{\delta_1' \le y\le \frac{\pi}{2}-\delta_1'\},
\end{equation}
where $\shock$ denotes the curved pseudo-transonic shock of $\vphi$.
Thus, the first inequality in \eqref{deriv-psi-xy-right} is equivalent to
$\psi_y(x,y)\ge -(\rightc-x) \psi_x(x,y) \cot y$ on $\partial \Omega^{\mcl{N}}\cap \shock$.
Then \eqref{estimate-psix-xy-right}
implies that there exists a constant $C_{\rm sh}>0$ depending only
on $(\iv, \gam)$ such that
\begin{equation}
\label{Csh}
\psi_y \ge -C_{\rm sh} \,x\qquad \tx{on $\shock\cap \der\nOm_{\bar{\eps}}$}.
\end{equation}
By \eqref{5-a4}, we have
\begin{equation}
\label{slip-bc-right}
\psi_y=0\qquad\tx{on $\rightsonic\cup(\Wedge\cap \der \nOm_{\bar{\eps}})$}.
\end{equation}
By \eqref{3-c6} in Lemma \ref{lemma-step3-1}, there exists a constant $C_{\rm in}>0$ depending only on $(\iv, \gam)$
such that $\psi$ satisfies
\begin{equation}
\label{Cin}
\psi_y\ge -C_{\rm in}   \qquad\tx{on $\nOm$}.
\end{equation}

\smallskip
{\textbf{3.}}
By adjusting Step 3 in the proof of \cite[Lemma 11.2.6]{CF2}, the following lemma holds:
\begin{lemma}
\label{lemma-app-eqnxy-yderiv-new}
Fix constants $\gam\ge 1$, $c>0$, and $r_0\in(0, \frac{c}{2}]$. Given an open
set
\begin{equation*}
U\subset \{(x,y)\in \R^2\,:\,0<x<r_0\},
\end{equation*}
assume that a function $\psi\in C^3(U)$ satisfies the equation{\rm :}
\begin{equation*}
\mcl{N}_{pl}(\psi):=\big(2x-(\gam+1)\psi_x+O_1\big)\psi_{xx}+O_2\psi_{xy}+\big(\frac 1c+O_3\big)\psi_{yy}
 -(1+O_4)\psi_x+O_5\psi_y=0\quad\,\tx{in $U$},
\end{equation*}
with $O_j=O_j(D\psi(x,y), \psi(x,y),x,c)$ for $j=1\cdots, 5$, where each $O_j(p_x, p_y,z, x, c)$ is defined by \eqref{prelim5-5}.
Moreover, let $\psi$ satisfy the following inequalities{\rm :}
\begin{equation*}
\psi\ge 0,\quad\,\,
0\le \psi_x\le \frac{2-\delta_0}{1+\gam}x \qquad\,\,\tx{in $U$},
\end{equation*}
for some constant $\delta_0\in(0,1)$.
Then there exists a constant $\eps\in(0,r_0)$ depending only on $(\gam, c, \delta_0)$
so that
$\der_y \mcl{N}_{pl}(\psi)=0$ is rewritten as a linear equation
for $w:=\psi_y$ in the following form{\rm :}
\begin{equation}
\label{equation-polar-w-rsonic}
\begin{split}
\mcl{L}_{\psi}(w):=&\, \big(2x-(\gam+1)O_1\big)w_{xx}+O_2 w_{xy}+\big(\frac 1c+O_3\big)w_{yy}\\
&\, +b^{(\psi)}_1 w_x+b^{(\psi)}_2 w_y+b^{(\psi)}_0w
=0\qquad\,\,\tx{in $U\cap\{x<\eps\}$},
\end{split}
\end{equation}
with
\begin{equation}
\label{sign-bterms-rsonic}
b^{(\psi)}_1\le 0,\quad
b^{(\psi)}_0\le 0\qquad\,\, \tx{in $U\cap \{x<\eps\}$}.
\end{equation}
\end{lemma}

By Definition \ref{def-regular-sol}(iv) and \eqref{estimate-psix-xy-right},
we can apply Lemma \ref{lemma-app-eqnxy-yderiv-new} to $\psi=\vphi-\rightvphi$.
Therefore, we can further reduce constant $\bar{\eps}\in (0, \eps_0]$ depending only on $(\iv, \gam)$
so that $\psi_y$ satisfies the elliptic
equation:
\begin{equation*}
\mcl{L}_{\psi}(\psi_y)=0\qquad\,\, \tx{in $\nOm_{\bar{\eps}}$.}
\end{equation*}

For constants $C_{\rm sh}$ and $C_{\rm in}$ from \eqref{Csh} and \eqref{Cin}, respectively, we choose
$M:=\max\{C_{\rm sh},  \frac{C_{\rm in}}{\bar{\eps}}\}$.
Then
$w=\psi_y$ satisfies
\begin{equation*}
  \begin{split}
  &w+Mx\ge 0\qquad\,\,\tx{on $\der \nOm_{\bar{\eps}}$},\\
  &\mcl{L}_{\psi}(w+Mx)=\mcl{L}_{\psi}(Mx)
  =M\big(b_1^{(\psi)}+b_0^{(\psi)}x\big)\le 0
  \qquad\,\,\tx{in $\nOm_{\bar{\eps}}$}.
  \end{split}
\end{equation*}
The second inequality stated above is obtained from \eqref{sign-bterms-rsonic}.
Note that constant $M$ is chosen depending only on $(\iv, \gam)$.
By the maximum principle, we obtain
\begin{equation*}
  w(x,y)\ge -Mx\qquad\,\,\tx{in $\nOm_{\bar{\eps}}$}.
\end{equation*}
Combining this with \eqref{estimate-psix-xy-right}--\eqref{upper-bd-psi-y}
yields statement (b) of Lemma \ref{lemma-est-nrsonic-N}.

\smallskip
{\textbf{4.}}
By Lemma \ref{lemma1-sonic-N}(c) and Lemma \ref{lemma-est-nrsonic-N}(b), we have
\begin{equation*}
  \der_x(\ivphi-\vphi)\le \der_x\iphi^{\mcl{N}}\le \frac{\mathfrak{M}}{2}\qquad\tx{in $\nOm_{\bar{\eps}}$}.
\end{equation*}
By Lemma \ref{lemma1-sonic-N}(c) and Lemma \ref{lemma-est-nrsonic-N}(a), we obtain
\begin{equation*}
  \der_x(\ivphi-\vphi)(x,y)= \der_x\iphi^{\mcl{N}}(x,y)-\psi_x\ge \frac{2y}{\mathfrak{M}}-\frac{2-\delta}{1+\gam}x
  \qquad\tx{in $\nOm_{\bar{\eps}}$}.
\end{equation*}
The estimate of $\der_y(\ivphi-\vphi)$ stated in statement (c) of Lemma \ref{lemma-est-nrsonic-N} is similarly obtained.

The existence of a function $\fshockn:[0, \bar{\eps}]\rightarrow \R^+$ satisfying statement (d) directly follows
from $\ivphi-\vphi=0$ on $\shock$, Lemma \ref{lemma-est-nrsonic-N}(c),
and the implicit function theorem.

Finally, statement (e) directly follows from statements (a)--(b) and (d) of Lemma \ref{lemma-est-nrsonic-N}, and  Definition \ref{def-regular-sol}(iv).
\end{proof}
\end{lemma}

\begin{lemma}\label{lemma-xybvp-N}
Write Eq. \eqref{eqn-xy-right} in $\nOm$ as
\begin{equation*}
\sum_{i,j=1}^2\hat A_{ij}^{\mcl{N}}(D\psi,\psi,x)D^2_{ij}\psi+\sum_{i=1}^2\hat A_i^{\mcl{N}}(D\psi,\psi,x)D_i\psi=0,
\end{equation*}
with $(D_1,D_2)=(D_x, D_y)$ and $\hat A^{\mcl{N}}_{21}=\hat A^{\mcl{N}}_{12}$.
Then there exist $\eps_{\mcl{N}}\in(0, \frac{\bar{\eps}}{4}]$ and $\lambda_{\mcl{N}}>0$ depending only on $(\iv, \gam)$
such that, for any admissible solution $\vphi=\psi+\rightvphi$
corresponding to $(\iv, \beta)\in \mathfrak{R}_{\rm weak}$, if $(x,y)\in \ol{\nOm_{4\eps_{\mcl{N}}}}$,
then
\begin{equation}
\label{5-a2}
\frac{\lambda_{\mcl{N}}}{2}|{\bm\kappa}|^2\le \sum_{i,j=1}^2 \hat A^{\mcl{N}}_{ij}(D\psi(x,y), \psi(x,y), x) \frac{\kappa_i\kappa_j}{x^{2-\frac{i+j}{2}}}
\le \frac{2}{\lambda_{\mcl{N}}}|{\bm\kappa}|^2\qquad\tx{for all}\,\,{\bm\kappa}=(\kappa_1, \kappa_2)\in\R^2.
\end{equation}
Moreover, $\mcl{B}_1^{\mcl{N}}$ defined by \eqref{def-B1-new-2015-N} satisfies the following properties{\rm :}

\smallskip
\begin{itemize}
\item[(a)] $\mcl{B}_{1}^{\mcl{N}}({\bf 0},0,x,y)=0$ for all $(x,y)\in \R^2${\rm ;}

\smallskip
\item[(b)] For each $k=2,3,\cdots$, there exist constants $\delta_{\rm bc}>0$ and $C>1$ depending only
on $(\iv, \gam, k)$ such that, whenever $|(p_x, p_y, z, x)|\le \delta_{\rm bc}$ and $|y-y_{\righttop}|\le \delta_{\rm bc}$,
\begin{equation*}
D_j B_{1}^{\mcl{N}}(p_x, p_y, z, x, y)\le -C^{-1}   \qquad \mbox{for $j=1,2,3$,}
\end{equation*}
where $(D_1, D_2, D_3):=(D_{p_x}, D_{p_y}, D_z)$.

\item[(c)] There exist constants $\hat{\delta}_{\rm bc}>0$ and $C>1$ depending only on $(\iv, \gam)$ such that,
whenever $|(p_x, p_y, z, x)|\le \hat{\delta}_{\rm bc}$ and $|y-y_{\righttop}|\le \hat{\delta}_{\rm bc}$,
\begin{equation*}
D_{(p_x,p_y,z)}B_{1}^{\mcl{N}}(p_x, p_y, z, x, y)\le -C^{-1}.
\end{equation*}
\end{itemize}
In {\rm (b)} and {\rm (c)} above, $y_{\righttop}$ represents the $y$--coordinate
of point $\righttop$, defined by Definition {\rm \ref{definition-domains-np}}.

\begin{proof}
\eqref{5-a2} can be checked directly from \eqref{prelim5-5}.
Properties (a)--(b) of $\mcl{B}_1^{\mcl{N}}$ are the results directly following
from the definition of $\rightvphi$, \eqref{def-gsh}, \eqref{7-d7-pre}, and \eqref{def-B1-new-2015-N}.

A direct calculation by using the definition of $\rightvphi$ in Definition \ref{definition-domains-np},
\eqref{8-49-ps}--\eqref{def-densi-ps}, \eqref{def-gsh},
\eqref{7-d7-pre}, and \eqref{def-B1-new-2015-N} yields that
\begin{align*}
&\der_z\mcl{B}_1^{\mcl{N}}({\bf 0}, 0, 0, y_{\righttop})=-\frac{\rightrho\iv \neta}{\rightc^2},\\
&\der_{p_x}\mcl{B}_1^{\mcl{N}}({\bf 0}, 0, 0, y_{\righttop})=-\frac{\rightrho-1}{\rightc}(\xi_1^{\mcl{N}})^2,\\
&\der_{p_y}\mcl{B}_1^{\mcl{N}}({\bf 0}, 0, 0, y_{\righttop})
=-\frac{\xi_1^{\mcl{N}}}{\rightc^2}\big(\rightrho \iv+(\rightrho-1)\neta\big).
\end{align*}
Then property (c) is obtained by combining the results stated immediately above with property (b).
\end{proof}
\end{lemma}

\begin{lemma}\label{lemma-coeff-ext-xy-N}
Let $\eps_0>0$ and $L\ge 1$ be the constants from Lemma {\rm \ref{lemma1-sonic-N}} and Lemma {\rm \ref{lemma-est-nrsonic-N}}, respectively.
Then there exist constants $\eps\in(0,\frac{\eps_0}{2}]$ and $C>0$ depending only on $(\iv, \gam)$
such that any admissible solution $\vphi=\rightvphi+\psi$ corresponding to $(\iv, \beta)\in \mathfrak{R}_{\rm weak}$
satisfies the following equation{\rm :}
\begin{equation*}
\sum_{i,j=1}^2 {\hat A}^{\rm (mod)}_{ij}(D\psi, \psi, x)D_{ij}\psi
+\sum_{i=1}^2 {\hat A}^{\rm (mod)}_{i}(D\psi, \psi, x)D_{i}\psi=0 \qquad\,\, \tx{in $\nOm_{\eps}$},
\end{equation*}
with coefficients $(\hat{A}_{ij}^{\rm (mod)}, \hat{A}_{i}^{\rm (mod)})$
satisfying the following properties{\rm :}

\smallskip
\begin{itemize}
\item[(a)] $(\hat{A}_{ij}^{\rm (mod)}, \hat{A}_{i}^{\rm (mod)})=(\hat{A}^{\mcl{N}}_{ij}, \hat{A}^{\mcl{N}}_{i})\,$
in $\{(p_x, p_y, z, x)\,:\,|(p_x, p_y)|\le Lx, \,\, |z|\le Lx^2, \,\, x\in(0,\eps)\}$,

\smallskip
\item[(b)] $|(\hat{A}_{11}^{\rm (mod)}, \hat{A}_{12}^{\rm (mod)}, \hat{A}_2^{\rm (mod)})(p_x, p_y, z, x)|\le Cx\, $
in $\R^2\times \R\times (0,\eps)$,

\smallskip
\item[(c)] $\|(\hat{A}_{22}^{\rm (mod)}, \hat{A}_1^{\rm (mod)})\|_{0, \R^2\times \R\times (0,\eps)}\le C$,

\smallskip
\item[(d)] $\|D_{(p_x, p_y, z, x)}(\hat{A}_{ij}^{\rm (mod)}, \hat A_i^{\rm (mod)})\|_{0, \R^2\times \R\times (0,\eps)} \le C$.
\end{itemize}
\end{lemma}

\begin{proof}
This lemma can be proved by adjusting the proof of \cite[Corollary 11.2.12]{CF2}.

Choose a function $\eta\in C^{\infty}(\R)$ such that $0\le \eta \le 1$ with
$\eta(t)=1$ for $|t|\le L$ and $\eta(t)=0$ for $|t|\ge 2L$.
For such a function $\eta$, we define $(\hat{A}_{ij}^{\rm (mod)}, \hat A_{i}^{\rm (mod)})$ by
\begin{equation}
\label{definition-coeffxy-mod}
(\hat{A}_{ij}^{\rm (mod)}, \hat A_{i}^{\rm (mod)})(p_x, p_y, z, x)
=(\hat{A}_{ij}^{\mcl{N}}, \hat{A}_{i}^{\mcl{N}} )(x\eta(\frac{p_x}{x}), x\eta(\frac{p_y}{x}), x^2 \eta(\frac{z}{x^2}), x).
\end{equation}
Then Lemma \ref{lemma-coeff-ext-xy-N} directly follows from \eqref{eqn-xy-right} and Lemma \ref{lemma-est-nrsonic-N}.
\end{proof}

For the uniform weighted $C^{2,\alp}$--estimates of admissible solutions
near $\rightsonic$, we recall the definition of the norm introduced in \cite{ChenFeldman2b}.

\begin{definition} [Parabolic norms]
\label{definition-parabolic-norm}
Fix a constant $\alp\in(0,1)$.

\smallskip
\begin{itemize}
\item[(i)]
For $z=(x,y), \til z=(\til x, \til y)\in \R^2\cap\{x>0\}$, define
\begin{equation*}
\delta_{\alp}^{\rm (par)}(z,\til z):=\left(|x-\til x|^2+\max\{x,\til x\}|y-\til y|^2\right)^{\frac{\alp}{2}}.
\end{equation*}

\item[(ii)]
Let $\mcl{D}$ be an open set in $\R^2\cap\{x>0\}$. For a function $u\in C^2(\mcl{D})$ in the $(x,y)$--coordinates, define
\begin{equation*}
\begin{split}
&\|u\|_{2,0,\mcl{D}}^{\rm (par)}:=\sum_{0\le k+l\le 2}\sup_{z\in \mcl{D}}\bigl(x^{k+\frac{l}{2}-2}|\der_x^k\der_y^lu(z)|\bigr),\\
&[u]^{\rm (par)}_{2,\alp,\mcl{D}}:=
\sum_{k+l=2}\sup_{z,\til z\in \mcl{D}, z\neq \til z}\Bigl(\min\big\{x^{\alp+k+\frac l2-2}, \til{x}^{\alp+k+\frac l2-2}\big\}
\frac{|\der_x^k\der_y^lu(z)-\der_x^k\der_y^lu(\til z)|}{\delta_{\alp}^{\rm (par)}(z,\til z)}\Bigr),\\
&\|u\|^{\rm (par)}_{2,\alp,\mcl{D}}:=\|u\|^{\rm (par)}_{2,0,\mcl{D}}+[u]_{2,\alp,\mcl{D}}^{\rm (par)}.
\end{split}
\end{equation*}

\item[(iii)]
Fix an open interval $I:=(0,a)$. For a function $f\in C^2(I)$, define
\begin{equation*}
\begin{split}
&\|f\|_{2,0,I}^{\rm (par)}:=\sum_{k=0}^2\sup_{x\in I}\bigl(x^{k-2}|\der_x^kf(x)|\bigr),\\
&[f]_{2,\alp,I}^{\rm (par)}:=\sup_{x,\til x\in I, x\neq\til x}\Bigl(\min\{x^{\alp}, \til{x}^{\alp}\}\frac{|\der_x^2f(x)-\der_x^2f(\til x)|}{|x-\til x|^{\alp}}\Bigr),\\
&\|f\|^{\rm (par)}_{2,\alp,I}:=\|f\|_{2,0,I}^{\rm (par)}+[f]_{2,\alp,I}^{\rm (par)}.
\end{split}
\end{equation*}

\item[(iv)]
Given constants $\sigma>0$, $\alp\in(0,1)$, and $m\in \mathbb{Z}_+$, define
\begin{equation*}
\begin{split}
&\|u\|_{m,0,\mcl{D}}^{(\sigma),{\rm (par)}}
 :=\sum_{0\le k+l\le m}\sup_{z\in \mcl{D}}\bigl(x^{k+\frac{l}{2}-\sigma}|\der_x^k\der_y^lu(z)|\bigr),\\
&[u]^{(\sigma), {\rm (par)}}_{m,\alp,\mcl{D}}:=
\sum_{k+l=m}\sup_{z,\til z\in \mcl{D}, z\neq \til z}\Bigl(\min\big\{x^{\alp+k+\frac{l}{2}-\sigma},{\til x}^{\alp+k+\frac{l}{2}-\sigma}\big\}
  \frac{|\der_x^k\der_y^lu(z)-\der_x^k\der_y^lu(\til z)|}{\delta_{\alp}^{\rm (par)}(z,\til z)}\Bigr),\\
&\|f\|_{m,0,I}^{(\sigma),{\rm (par)}}
:=\sum_{k=0}^m\sup_{x\in I}\bigl(x^{k-\sigma}|\der_x^kf(x)|\bigr),\\
&[f]_{m,\alp,I}^{(\sigma), {\rm (par)}}
:=\sup_{x,\til x\in I, x\neq\til x}\Bigl(\min\left\{x^{\alp+m-\sigma},{\til x}^{\alp+m-\sigma}\right\}\frac{|\der_x^mf(x)-\der_x^mf(\til x)|}{|x-\til x|^{\alp}}\Bigr),\\
&\|u\|^{(\sigma), {\rm (par)}}_{m,\alp,\mcl{D}}
:=\|u\|^{(\sigma),{\rm (par)}}_{m,0,\mcl{D}}+[u]_{m,\alp,\mcl{D}}^{(\sigma),{\rm (par)}},\qquad\,\,
\|f\|^{(\sigma), {\rm (par)}}_{m,\alp,I}:=\|f\|_{m,0,I}^{(\sigma),{\rm (par)}}+[f]_{m,\alp,I}^{(\sigma), {\rm (par)}}.
\end{split}
\end{equation*}
Note that  norm $\|\cdot\|_{2,\alp,\mcl{D}}^{\rm (par)}$ in {\rm (ii)} is
norm $\|\cdot\|_{2,\alp,\mcl{D}}^{(2),{\rm (par)}}$ above here.

\smallskip
\item[(v)]
Denote $C^{m,\alp}_{(\sigma),{\rm (par)}}(\mcl{D}):= \{u\in C^{m}(\mcl{D})\,:\,\|u\|_{m,\alp,\mcl{D}}^{(\sigma),{\rm (par)}}<\infty\}$.
%
%
\end{itemize}
\end{definition}

\begin{proposition}
\label{lemma-est-sonic-general-N}
Let $\eps_{\mcl{N}}>0$ be from Lemma {\rm \ref{lemma-xybvp-N}}.
For each $\alp\in(0,1)$, there exists $C>0$ depending only on $(\iv, \gam, \alp)$
such that any admissible solution $\vphi$ corresponding to $(\iv, \beta)\in \mathfrak{R}_{\rm weak}$ satisfies
\begin{equation}
\label{par-est1-N}
\|\vphi-\rightvphi\|^{\rm (par)}_{2,\alp,\nOm_{\eps_{\mcl{N}}}}+\|\fshockn-\hat{f}_{\mcl{N},0}\|^{\rm (par)}_{2,\alp,(0,\eps_{\mcl{N}})}\le C.
\end{equation}
\end{proposition}

\begin{proof}
The proof is divided into six steps.

\smallskip

{\bf 1}. {\emph{Re-scaling coordinates}}. Fix $\eps\in(0, \frac{\eps_{\mcl{N}}}{2}]$.
For $z_0:=(x_0,y_0)\in \ol{\nOm_{\eps}}\setminus \overline{\rightsonic}$ and $r\in (0, 1]$,
define
\begin{equation*}
\til{R}_{z_0,r}:=\{(x,y)\,:\,|x-x_0|<\frac{r}{4} x_0,\,\, |y-y_0|<\frac{r}{4}\sqrt{x_0}\},\qquad
R_{z_0,r}:=\til{R}_{z_0,r}\cap \nOm_{2\eps}.
\end{equation*}
If $\eps\le y_{\righttop}^2$ and $z_0\in \shock \cap \ol{\nOm_{\eps}}$, then it follows from Lemma \ref{lemma-est-nrsonic-N}(d) that
\begin{equation}
\label{Rz-property-new}
R_{z_0,1} \subset \{(x,y): \frac{3}{4}x_0<x<\frac{5}{4}x_0, \frac{3}{4}y_0<y<\frac{5}{4}y_0,\}.
\end{equation}
For $r>0$, define the sets:
\begin{equation*}
Q_{r}:=(-r, r)^2,\qquad
Q_r^{(z_0)}:=\{(S,T)\in Q_r\,:\,z_0+\frac 14(x_0 S, \sqrt{x_0}\,T)\in R_{z_0, r}\}.
\end{equation*}

\smallskip
{\bf 2}. {\emph{Re-scaled function $\psi^{(z_0)}$}}. Let $\psi$ be given by \eqref{psi-def-right}.
For  $z_0\in \partial \Omega^{\mcl{N}}\cap \shock$, define a function $\psi^{(z_0)}(S,T)$ by
\begin{equation*}
\psi^{(z_0)}(S,T) =\frac{1}{x_0^2}\psi(x_0+\frac{x_0}{4}S, y_0+\frac{\sqrt{x_0}}{4} T)\qquad\tx{for}\,\, (S,T)\in Q_{1}^{(z_0)}.
\end{equation*}
By Lemma \ref{lemma-est-nrsonic-N} and \eqref{Rz-property-new},
we have
\begin{equation*}
|\psi^{(z_0)}|\le L,\quad |\psi^{(z_0)}_S|\le L, \quad |\psi^{(z_0)}_T|\le Lx_0^{-1/2} \qquad\,\,
\tx{in}\,\,Q^{(z_0)}_{1/2}.
\end{equation*}
Moreover, Lemma \ref{lemma-coeff-ext-xy-N} implies that $\psi^{(z_0)}$ satisfies the equation:
\begin{equation*}
\sum_{i,j=1}^2A_{ij}^{(z_0)}(D\psi^{(z_0)}, \psi^{(z_0)}, S) D_{ij}\psi^{(z_0)}
+ \sum_{i=1}^2 A_{i}^{(z_0)}(D\psi^{(z_0)}, \psi^{(z_0)}, S) D_{i}\psi^{(z_0)}=0 \qquad\tx{in}\,\,Q_{1/2}^{(z_0)},
\end{equation*}
where $(D_1, D_2)=(D_S, D_T)$, $D_{ij}=D_iD_j$, and
\begin{equation*}
\begin{split}
&A_{ij}^{(z_0)}(p_1, p_2, z, S):=x_0^{\frac{i+j}{2}-2}A_{ij}^{\rm (mod)}(4x_0p_1, 4x_0^{3/2}p_2, x_0^2 z, x_0(1+\frac S4)),\\
&A_{i}^{(z_0)}(p_1, p_2, z, S):=\frac 14 x_0^{\frac{i-1}{2}}A_{i}^{\rm (mod)}(4x_0p_1, 4x_0^{3/2}p_2, x_0^2 z, x_0(1+\frac S4)).
\end{split}
\end{equation*}

For $\fshockn$ given in Lemma \ref{lemma-est-nrsonic-N}(d), we define
\begin{equation}
\label{definition-scaled-shock-right}
F^{(z_0)}(S):=\frac{4}{\sqrt{x_0}} \left(\fshockn(x_0+\frac{x_0}{4}S)-\fshockn(x_0)\right)\qquad\tx{for $-1<S<1$}.
\end{equation}
It follows directly from  Lemma \ref{lemma-est-nrsonic-N}(d) and \eqref{definition-scaled-shock-right} that $F^{(z_0)}$ satisfies
\begin{equation}
\label{estimate-scaled-F}
F^{(z_0)}(0)=0,\qquad\,\, \|F^{(z_0)}\|_{C^{1}([-1,1])}\le C\sqrt{x_0}
\end{equation}
for some constant $C>0$ depending only on $(\iv, \gam)$.
Therefore, there exists $\eps_*\in (0, \frac{\bar{\eps}}{2}]$ depending only on $(\iv, \gam)$
such that $F^{(z_0)}(S)>-\frac{r}{2}$ for $S\in (-r,r)$, whenever $r\in(0,1)$
and  $z_0\in \overline{\nOm_{\eps_*}}\cap \shock$.

For $z_0\in \overline{\nOm_{\eps_*}}\cap \shock$, define
\begin{equation*}
\shock^{(z_0)}:=\{(S,T)\,:\, S\in(-1,1), T=F^{(z_0)}(S)\} \subset \der Q_1^{(z_0)}.
\end{equation*}
Then ${\rm dist}(\shock^{(z_0)}, \der Q_1^{(z_0)}\cap \{T=-1\})\ge \frac 12$.

By Lemma \ref{lemma-est-nrsonic-N}(a)--(b) and (e), we can fix a small constant $\eps_*\in (0, \frac{\bar{\eps}}{2}]$
depending only on $(\iv, \gam)$ so that any admissible solution satisfies
\begin{equation*}
|(\psi_x, \psi_y, \psi, y-y_{\righttop})|\le \frac 14\min\{\delta_{\rm bc}, \hat{\delta}_{\rm bc}\} \qquad\tx{in $\nOm_{2\eps_*}$}
\end{equation*}
for constants $(\delta_{\rm bc}, \hat{\delta}_{\rm bc})$ from Lemma \ref{lemma-xybvp-N}.
Then we apply Lemma \ref{lemma-xybvp-N}(c) and the implicit function theorem to rewrite
the boundary condition \eqref{bc-psi-B1-right} as
\begin{equation}
\label{bc-almost-tangential-N}
\psi_x=b_{\mcl{N}}(\psi_y, \psi, x, y)\qquad\tx{on $\shock\cap \ol{\nOm_{2\eps_*}}$}.
\end{equation}
By Lemma \ref{lemma-xybvp-N}(a)--(b), we have
\begin{equation}
\label{regularity-b-new}
\begin{split}
b_{\mcl{N}}(0,0,x,y)=0\qquad  &\tx{in $\ol{\nOm_{2\eps_*}}$},\\
|D^k b_{\mcl{N}}(p_y, z, x, y)|\le C_k\qquad &\tx{in}\,\,\R\times \R\times \ol{\nOm_{2\eps_*}}\,\,\tx{for $k=1,2,3,\cdots$},
\end{split}
\end{equation}
where constants $C_k>0$ depend only on $(\iv, \gam, k)$.

For each $z_0\in \shock\cap \overline{\Omega^{\mcl{N}}_{\eps_*}}$,
denote
\begin{equation}
\label{definition-scaledB-new}
B^{(z_0)}_{\mcl{N}}(p_T, z, S, T):=\frac{1}{4x_0} b_{\mcl{N}}(4x_0^{3/2}p_T, x_0^2z, x,y)
\qquad\,\,\mbox{for $(x,y)=z_0+(\frac{x_0}{4}S, \frac{\sqrt{x_0}}{4}T)$}.
\end{equation}
It follows directly from \eqref{regularity-b-new}
that there exists a constant $m_1>0$ depending only on $(\iv, \gam)$ such that
\begin{equation}
\label{regularity-Bn-new}
\begin{split}
B_{\mcl{N}}^{(z_0)}(0,0,S,T)=0\qquad&\tx{in $\ol{Q_1^{(z_0)}}$},\\
\|\der_{p_T}B_{\mcl{N}}^{(z_0)}(p_T, z,\cdot)\|_{0, \ol{Q_1^{(z_0)}}}\le m_1\sqrt{x_0}\qquad&\tx{for all $(p_T, z)\in \R\times \R$},\\
\|D_{(p_T, z)}B_{\mcl{N}}^{(z_0)}(p_T, z,\cdot)\|_{1, \ol{Q_1^{(z_0)}}}\le m_1\sqrt{x_0}\qquad&\tx{for all $(p_T, z)\in \R\times \R$}.
\end{split}
\end{equation}

By \eqref{bc-almost-tangential-N}, $\psi^{(z_0)}$ satisfies
\begin{equation}
\label{scaled-shockbc-new}
\psi_S^{(z_0)}=B^{(z_0)}_{\mcl{N}}(\psi_T^{(z_0)}, \psi^{(z_0)}, S, T)\qquad\tx{on $\shock^{(z_0)}$}.
\end{equation}

\smallskip
{\bf 3}. {\emph{Uniform estimates of $\psi^{(z_0)}$  for $z_0 \in \shock$}}.
By  \eqref{estimate-scaled-F} and \eqref{regularity-Bn-new}, we can apply Theorem \ref{elliptic-t8-CF2} to
find constants $(\eps,\delta, C )\in(0, \eps_*]\times (0,1)\times (0, \infty)$ depending only on $(\iv, \gam)$
so that, for any $z_0\in \ol{\nOm_{\eps}} \cap \shock$, we have
\begin{equation}
\label{expression-scaled-shock}
\|\psi^{(z_0)}\|_{1,\delta, \ol{Q^{(z_0)}_{3/4}}}\le C.
\end{equation}
By \eqref{definition-scaled-shock-right}, for each $z_0\in \ol{\nOm_{\eps}}\cap \shock$,
$\iphi^{\mcl{N}}=\ivphi-\rightvphi$ satisfies
\begin{equation}
\label{expression-F-right}
\iphi^{\mcl{N}}(x_0+\frac{x_0}{4}S, \fshockn(x_0)+\frac{\sqrt{x_0}}{4}F^{(z_0)}(S))
-x_0^2\psi^{(z_0)}(S, F^{(z_0)}(S))=0\qquad\tx{for $-1<S<1$}.
\end{equation}
Differentiating \eqref{expression-F-right} with respect to $S$, we have
\begin{equation}
\label{expression-Fdiff-new}
(F^{(z_0)})'=-\frac{\sqrt{x_0}(\der_x\iphi^{\mcl{N}}-4x_0\der_S\psi^{(z_0)})}{\der_y\iphi^{\mcl{N}}-4 x_0^{3/2}\der_T\psi^{(z_0)}}.
\end{equation}
By combining this expression with Lemma \ref{lemma1-sonic-N}(c) and \eqref{expression-scaled-shock},
a direct computation shows that there exists a small constant $\eps\in(0, \eps_*]$
depending on $(\iv, \gam)$ such that $F^{(z_0)}$ satisfies the estimate:
\begin{equation}
\label{estimate-scaledF-new}
\|F^{(z_0)}\|_{1,\delta, [-3/4, 3/4]}\le C\sqrt{x_0}\qquad\,\,\tx{for all $z_0=(x_0, y_0)\in \shock\cap \overline{\Omega^{\mcl{N}}_{\eps}}$}
\end{equation}
for some constant $C>0$ depending only on $(\iv,\gam)$.

This result, combined with Lemma \ref{lemma-unif-est2}, yields
that $\shock$ is $C^{1,\delta}$ up to $\rightsonic$ away from $\leftsonic$.

Next, it follows directly from \eqref{estimate-scaledF-new} and a direct computation
by using \eqref{regularity-b-new}--\eqref{definition-scaledB-new} that
the boundary condition \eqref{scaled-shockbc-new} satisfies all the conditions stated
in Theorem \ref{elliptic-t9-CF2} with $(\alp, \Phi, W)=(\delta, \frac{1}{\sqrt{x_0}}F^{(z_0)}, B_{\mcl{N}}^{(z_0)})$
for all $z_0\in \shock\cap\partial \Omega^{\mcl{N}}_{\eps}$,
where $\eps>0$ is the constant in \eqref{estimate-scaledF-new}.
Therefore, we can further reduce $\eps\in(0, \eps_*]$ depending on $(\iv, \gam)$ so that,
for each $z_0\in\shock\cap \der \nOm_{\eps}$,
the re-scaled function $\psi^{(z_0)}$ satisfies the estimate:
\begin{equation}
\label{estimate-scaled-right-new2}
\|\psi^{(z_0)}\|_{2,\delta, Q^{(z_0)}_{1/2}}\le C,
\end{equation}
where $C$ depends only on $(\iv, \gam)$.

We combine estimate \eqref{estimate-scaled-right-new2} with \eqref{expression-Fdiff-new}
to see that $F^{(z_0)}\in C^{1,\alp}([-\frac{1}{2}, \frac{1}{2}])$ for any $\alp\in(0,1)$.
Furthermore, we have
\begin{equation*}
\sup_{z_0\in \shock\cap \der \nOm_{\eps}}\frac{1}{\sqrt{x_0}}\|F^{(z_0)}\|_{1,\alp, [-\frac{1}{2}, \frac{1}{2}]}\le C,
\end{equation*}
where $C>0$ depends only on $(\iv, \gam)$.
Then we can repeat the previous argument by applying Theorem \ref{elliptic-t9-CF2} to conclude that,
for each $\alp\in(0,1)$, the small constant $\eps\in(0,\eps_*]$  can be further reduced so that
\begin{equation*}
\sup_{z_0\in \shock \cap \der \nOm_{\eps}}
\|\psi^{(z_0)}\|_{2,\alp,\ol{Q_{1/4}^{(z_0)}}}+\frac{1}{\sqrt{x_0}}\|F^{(z_0)}\|_{2,\alp, [-\frac 14, \frac 14]}\le C,
\end{equation*}
where $C>0$ is a constant depending only on $(\iv, \gam, \alp)$.

\smallskip
{\bf 4}. \emph{Uniform estimates of $\psi^{(z_0)}$  for $z_0 \not\in \shock$}.
If $Q_1^{(z_0)}=Q_1$, we apply Theorem \ref{elliptic-t1-CF2} to obtain that,
for each $\alp\in(0,1)$, $\|\psi^{(z_0)}\|_{2,\alp, \ol{Q_{1/2}^{(z_0)}}}$ is uniformly bounded
above by a constant depending only on $(\iv, \gam,  \alp)$.
If $z_0\in \Wedge\cap \der \nOm_{\eps}$, then $Q_1^{(z_0)}=Q_1\cap \{T>0\}$,
and $\psi^{(z_0)}$ satisfies that $\psi_T^{(z_0)}(S,0)=0$ for all $-1<S<1$.
This is owing to the slip boundary condition \eqref{slip-bc-right}.
In this case, we apply Theorem \ref{elliptic-t2-CF2} to obtain a uniform estimate
of $\|\psi^{(z_0)}\|_{2,\alp, \ol{Q_{1/2}^{(z_0)}}}$ for all $z_0\in \Wedge\cap \der\nOm_{\eps}$.

\smallskip
{\bf 5}.  {\emph{Estimate for $\|\vphi-\rightvphi\|^{\rm (par)}_{2,\alp,\nOm_{\eps_{\mcl{N}}}}$}}.
Since the estimates of $\|\psi^{(z_0)}\|_{2,\alp, \ol{Q_{1/8}^{(z_0)}}}$ are given independently of
$z_0\in \ol{\nOm_{\eps}}\setminus \rightsonic$ and $\beta\in [0,\betadet)$,
the estimate of $\|\vphi-\rightvphi\|^{\rm (par)}_{2,\alp,\nOm_{\eps_{\mcl{N}}}}$ in \eqref{par-est1-N}
is finally obtained by combining the uniform $C^k$--estimate of admissible solutions
given in Corollary \ref{corollary-unif-est-away-sn} and all the estimates of
$\|\psi^{(z_0)}\|_{2,\alp, \ol{Q_{1/8}^{(z_0)}}}$ from Steps 3--4, and by scaling back to the $(x,y)$--coordinates.
For the details, we refer to \cite[Steps 3--4 in the proof of Theorem 3.1]{BCF} or \cite[Lemma 4.6.1]{CF2}.

\smallskip

{\bf 6.}  {\emph{Estimate for $\|\fshockn-\hat{f}_{\mcl{N},0}\|^{\rm (par)}_{2,\alp,(0, \eps_{\mcl{N}})}$}}.
By Lemma \ref{lemma1-sonic-N}(e) and Lemma \ref{lemma-est-nrsonic-N}(d), we have
\begin{equation*}
\iphi^{\mcl{N}}(x, \hat{f}_{\mcl{N},0}(x))=0,\,\,\,\,
(\iphi^{\mcl{N}}-\psi)(x,\fshockn(x))=(\ivphi-\vphi)(x,\fshockn(x))=0\qquad\tx{for all $x\in[0,\eps_{\mcl{N}}]$}.
\end{equation*}
This yields that
\begin{equation}
\label{shock-expression-nrls-new2}
\iphi^{\mcl{N}}(x,\fshockn(x))-\iphi^{\mcl{N}}(x, \hat{f}_{\mcl{N},0}(x))=\psi(x,\fshockn(x))\qquad\tx{for all $x\in[0,\eps_{\mcl{N}}]$}.
\end{equation}
Since $|\der_y\iphi^{\mcl{N}}|>0$ from Lemma \ref{lemma1-sonic-N}(c), we can rewrite \eqref{shock-expression-nrls-new2} as
\begin{equation*}
\fshockn(x)-\hat{f}_{\mcl{N},0}(x)=
\frac{\psi(x,\fshockn(x))}{\int_0^1 \der_y\iphi^{\mcl{N}}(x, t\fshockn(x)+(1-t)\hat{f}_{\mcl{N},0}(x))\,\dd t}.
\end{equation*}
Then a direct computation by using Lemma \ref{lemma1-sonic-N} and the estimate
of $\|\psi\|^{\rm (par)}_{2,\alp,\nOm_{\eps_{\mcl{N}}}}\le C$ achieved in Step 5 implies that
\begin{equation*}
\|\fshockn-\hat{f}_{\mcl{N},0}\|^{\rm (par)}_{2,\alp,(0, \eps_{\mcl{N}})}\le C,
\end{equation*}
where $C>0$ is a constant depending only on $(\iv, \gam, \alp)$.
This completes the proof.
\end{proof}

\section{Weighted $C^{2,\alp}$--Estimates Near $\leftsonic$}
\label{subsec-unif-c2est-nr-leftsonic}

According to Definition \ref{definition-domains-np}, $\leftsonic$ depends continuously
on $\beta\in[0, \frac{\pi}{2})$.
In particular, the sonic arc $\leftsonic$ shrinks to a point when $\beta$ increases up
to $\betac^{(\iv)}$, and becomes a point $P_{\beta}$ for all $\beta\ge \betac^{(\iv)}$,
although the location of $P_{\beta}$ changes continuously
on $\beta\in[\betac^{(\iv)}, \frac{\pi}{2})$.
Furthermore, the ellipticity of Eq. \eqref{8-48-ps} on $\leftsonic$ also changes.
According to Proposition \ref{corollary-ellip}, the ellipticity of \eqref{8-48-ps}
degenerates on $\leftsonic$ for $\beta\le \betac^{(\iv)}$.
On the other hand, for $\beta>\betac^{(\iv)}$,
Eq. \eqref{8-48-ps} (or equivalently Eq. \eqref{2-1}) is uniformly elliptic up to $\leftsonic$ away from $\rightsonic$.
For that reason, the weighted $C^{2,\alp}$--estimates of admissible solutions near $\leftsonic$
are given for the following four cases separately:

\vspace{1mm}
1. $\beta<\betac^{(\iv)}$ away from $\betac^{(\iv)}$,

\vspace{1mm}
2. $\beta<\betac^{(\iv)}$ close to $\betac^{(\iv)}$,

\vspace{1mm}
3. $\beta\ge \betac^{(\iv)}$ close to $\betac^{(\iv)}$,

\vspace{1mm}
4. $\beta\in (\betac^{(\iv)}, \betadet)$ away from $\betac^{(\iv)}$.

\subsection{Case {\rm 1}{\rm :} Admissible solutions for $\beta<\betac^{(\iv)}$ {away from $\betac^{(\iv)}$}}
\label{subsubsec-apriori-est-leftsonic-case1}
For
$$
(\iv, \beta)\in \mathfrak{R}_{\rm weak}\cap \{\beta\,:\,0<\beta<\betac^{(\iv)}\},
$$
let $\Oo$ and $\lefttop$ be given by Definition
\ref{definition-domains-np}.
     For each $\beta>0$, let $\oM$ be defined by \eqref{1-25}. Define
     \begin{equation}
     \label{definition-leftchnew}
     \leftchnew:=\frac{|\lefttop\Oo|+\leftc\oM}{2}
     \begin{cases}
       =\frac{\leftc(1+\oM)}{2} &\quad \mbox{for $\beta\le \betasonic$}, \\[1.5mm]
       <\frac{\leftc(1+\oM)}{2} &\quad \mbox{for $\beta\ge \betasonic$}.
     \end{cases}
     \end{equation}
In $U_{\mcl{O}}:=\bigl(B_{\frac{3\leftc}{2}}(\Oo)\setminus B_{\leftchnew}(\Oo)\bigr)\cap\{\xxi\,:\,\xin<\leftu\}$,
use $(r,\theta)$ as the polar coordinates with respect to $\Oo=(\leftu, 0)$ and define
    \begin{equation}
    \label{coord-o}
    (x,y):=(\leftc-r,\pi-\theta).
    \end{equation}
Also, define a set $\oOm$ by
     \begin{equation*}
     \oOm:=\bigl(\Om\cap\{\xin<\leftu\}\bigr)\setminus B_{\leftchnew}(\Oo).
     \end{equation*}

Since $\oOm\subset B_{\leftc}(\Oo)$, $\oOm\subset\{(x,y)\,:\,x>0\}$.
In the $(x,y)$--coordinates defined by \eqref{coord-o},
$\leftvphi$ given by Definition \ref{definition-domains-np}  is written as\begin{equation}
\label{6-a4}
\begin{split}
&\leftvphi=-\frac 12(\leftc-x)^2+\frac 12 \leftu^2\leftk\qquad\text{in $U_{\mcl{O}}$}.
\end{split}
\end{equation}
For an admissible solution $\vphi$ corresponding to $(\iv, \beta)$, let $\psi$ be given by
\begin{equation}
\label{psi-def}
\psi=\vphi-\leftvphi\;\qquad\text{in $\oOm$}.
\end{equation}

\emph{{\rm (i)} Equation for $\psi$ in $\oOm$}: Similarly to \eqref{eqn-xy-right},
we rewrite Eq. \eqref{8-48-ps} for $\psi$ in the $(x,y)$--coordinates given by \eqref{coord-o}.
For each $j=1,\cdots, 5$, let $O^{\mcl{O}}_j(\pb, z, x)$ be given by
\begin{equation*}
O^{\mcl{O}}_j(\pb, z, x)=O_j(\pb, z, x,\leftc)
\end{equation*}
for $O_j(\pb, z, x,c)$ given by \eqref{prelim5-5}.
Then Eq. \eqref{2-1} is written as
\begin{equation}
\label{eqn-xy}
\big(2x-(\gam+1)\psi_x+O^{\mcl{O}}_1\big)\psi_{xx}+O^{\mcl{O}}_2\psi_{xy}
+\big(\frac{1}{\leftc}+O^{\mcl{O}}_3\big)\psi_{yy}-\big(1+O^{\mcl{O}}_4\big)\psi_x+O^{\mcl{O}}_5\psi_y=0,
\end{equation}
with $O^{\mcl{O}}_j=O^{\mcl{O}}_j(D\psi, \psi,x)$ for $j=1,\cdots, 5$.

\smallskip
\emph{{\rm (ii)} Boundary condition for $\psi$ on $\shock\cap \partial\oOm$}: Similarly to \eqref{7-d7-pre}, we define
\begin{equation}
\label{7-d7-pre-left}
M_{\beta}({\bf p}, z, \xin)=g^{\rm sh}_{\rm mod}({\bf p}+D\leftvphi,z+\leftvphi, \xin, \xi_2^{(\beta)}-\frac{\leftu\xin+z}{\iv})
\end{equation}
for $g^{\rm sh}_{\rm mod}$ given by \eqref{gsh-ext},
where $(D\leftvphi, \leftvphi)$ are evaluated at $(\xin, \xi_2^{(\beta)}-\frac{\leftu\xin+z}{\iv})$.
Note that $(\leftu, \etan^{(\beta)})$ depend continuously on $\beta\in(0,\frac{\pi}{2})$ and that
$$
\lim_{\beta\to 0+} (\leftu, \etan^{(\beta)})=(0, \neta).
$$

Define
\begin{equation}
\label{definition-iphi-O}
\iphi^{\mcl{O}}=\varphi_\infty-\varphi_{\mathcal{O}}.
\end{equation}
Rewriting the boundary condition: $|D(\iphi^{\mcl{O}}-\psi)|M_{\beta}(D\psi, \psi, \xin)=0$ on $\shock\cap \der \oOm$ in the $(x,y)$--coordinates
given by \eqref{coord-o}, we have
\begin{equation}
\label{bc-psi-B1-left}
\mcl{B}^{\mcl{O}}_1(\psi_x, \psi_y, \psi, x, y)=0\qquad\tx{on $\shock\cap\der{\oOm}$}
\end{equation}
for $\mcl{B}^{\mcl{O}}_1(p_x, p_y, z, x, y)$ given by
\begin{equation}
\label{def-B1-new-left}
\mcl{B}^{\mcl{O}}_1(p_x, p_y, z, x, y)
=|D\iphi^{\mcl{O}}-({p}_1, {p}_2)|M_{\beta}( p_1, p_2, z, {\xin})
\end{equation}
with
\begin{equation}\label{cov-left}
{\xin}=
\leftu-(\leftc-x)\cos y,\qquad
\begin{pmatrix}
p_1\\
p_2
\end{pmatrix}
=\begin{pmatrix}
\cos y& \sin y\\
-\sin y& \cos y
\end{pmatrix}
\begin{pmatrix}
p_x\\
\frac{p_y}{\leftc-x}
\end{pmatrix}.
\end{equation}

\smallskip
\emph{{\rm (iii)} Other properties of $\psi$}:
By \eqref{1-l} and conditions (ii) and (iv) of Definition \ref{def-regular-sol},
$\psi$ satisfies
\begin{equation}
\label{5-a4-left}
\begin{split}
&\psi \ge 0\qquad\,\,\text{in $\Om$},\\
&\psi=0\qquad\,\,\text{on $\leftsonic$},\\
&\psi_y=0\qquad\text{on $\Wedge\cap\der{\oOm}$}.
\end{split}
\end{equation}

For set $\mcl{D}$ defined by \eqref{def-D-domain}, let an open subset $\oLambdab$ of $\mcl{D}$ be given by
\begin{equation}
\label{def-Lambda-sets}
\begin{split}
&\oLambdab:=\mcl{D}\cap\bigl(B_{\frac{3\leftc}{2}}(\Oo)\setminus B_{\leftchnew}(\Oo)\bigr)\cap\{\xin<\leftu\}
\end{split}
\end{equation}
for $\leftchnew$ defined by \eqref{definition-leftchnew}.

\begin{lemma}\label{lemma-str-nr-sonic}
Fix $\gam\ge 1$ and $\iv>0$.
There exist positive constants $\eps_1$, $\eps_0$, $\delta_0$, $\omega_0$, $C$, and $\mathfrak{M}$
depending only on $(\iv, \gam)$ with $\eps_1>\eps_0$ and $\mathfrak{M}\ge 2$ such that,
for each $\beta\in(0,\betac^{(\iv)}]$, the following properties hold{\rm :}

\smallskip
\begin{itemize}
\item[(a)]
$\{\leftvphi<\ivphi\}\cap\oLambdab\cap\mcl{N}_{\eps_1}(\leftsonic)
\subset \{0<y<\frac{\pi}{2}-\beta-\delta_0\}${\rm ;}

\smallskip
\item[(b)]
  $\{\leftvphi<\ivphi\}\cap\mcl{N}_{\eps_1}(\leftsonic)\cap\{y>y_{\lefttop}\}
    \subset\{x>0\};$

\smallskip
\item[(c)] In $\{(x,y)\,:\,|x|<\eps_1, 0<y<\frac{\pi}{2}-\beta-\delta_0\}$,
$\iphi^{\mcl{O}}$ given by \eqref{definition-iphi-O} satisfies
\begin{equation}
\label{6-a5}
\frac {2}{\mathfrak{M}}(y+\tan \beta)\le \der_x\iphi^{\mcl{O}}\le \frac {\mathfrak{M}}{2},
\qquad \frac {2}{\mathfrak{M}}\le -\der_y\iphi^{\mcl{O}} \le \frac {\mathfrak{M}}{2};
\end{equation}

\item[(d)]
$|(D^2_{(x,y)}, D^3_{(x,y)})\iphi^{\mcl{O}}|\le C$ $\,\,$ in $\{|x|<\eps_1\}${\rm ;}

\smallskip
\item[(e)] There exists a unique function $\hat f_{\mcl{O},0}\in C^{\infty}([-\eps_0,\eps_0])$ such that
    \begin{equation}
    \label{est-sonic-1}
   \begin{cases}
   \{\leftvphi<\ivphi\}\cap \oLambdab\cap\mcl{N}_{\eps_1}(\leftsonic)\cap\{|x|<\eps_0\}=\{(x,y)\,:\,|x|<\eps_0, 0<y<\hat f_{\mcl{O},0}(x)\},\\[2mm]
    \leftshock\cap\mcl{N}_{\eps_1}(\leftsonic)\cap\{|x|<\eps_0\}=\{(x,y)\,:\,x\in(-\eps_0,\eps_0), y=\hat f_{\mcl{O},0}(x)\};
    \end{cases}
    \end{equation}

\item[(f)]
$\hat f_{\mcl{O},0}$ given in {\rm (e)} satisfies
    \begin{equation*}
    2\omega_0\le \hat f'_{\mcl{O},0}\le C \qquad\text{on $(-\eps_0,\eps_0)$}.
    \end{equation*}
\end{itemize}

\begin{proof} Note that line $\leftshock$ intersects with circle $\der B_{\leftc}(\Oo)$
at two different points, due to \eqref{oM-monotonicity} for any $(\iv, \beta)\in \mathfrak{R}_{\rm weak}$.
Point $\lefttop$ is an intersection point of $\leftshock=\{\xxi\,:\,\ivphi=\leftvphi\}$ with $\der B_{\leftc}(\Oo)$.
Let $\lefttop'$ be the other intersection point of $\leftshock$  and $\der B_{\leftc}(\Oo)$,
and let $\Qo$ be the midpoint of the line segment $\overline{\lefttop\lefttop'}$.
Then $\angle \Qo\Oo\leftbottom=\frac{\pi}{2}-\beta$.
Since $|\overline{\lefttop\Qo}|$ depends continuously on $\beta\in[0,\frac{\pi}{2})$,
there exists $\eps_1>0$ depending only on $(\iv, \gam)$
such that $|\overline{\lefttop\Qo}|\ge 2\eps_1$ for all $\beta\in[0,\betac^{(\iv)}]$.
Let $\Qo'$ be the midpoint of $\overline{\lefttop\Qo}$,
and let $(x_{\Qo'},y_{\Qo'})$ denote the $(x,y)$--coordinates of $\Qo'$.
Then there exists a constant  $\delta_0>0$ depending only on $(\iv, \gam)$
such that
\begin{equation}
\label{yQo-choice-new}
y_{\Qo'}<\frac{\pi}{2}-\beta-\delta_0.
\end{equation}

Moreover, it follows directly from \eqref{definition-iphi-O} that
\begin{equation*}
\der_x\iphi^{\mcl{O}}=\iv(\sin y+\cos y\tan\beta),
\,\,\,
\der_y\iphi^{\mcl{O}}=-\iv(\leftc-x)(\cos y-\sin y\tan\beta)
\qquad\, \text{in $\oLambdab$}.
\end{equation*}
Then statements (a)--(e) can be verified by performing a direct computation and using the observation obtained above.

Since $\iphi^{\mcl{O}}=0$ on $\leftshock$, we have
\begin{equation*}
\iphi^{\mcl{O}}(x, \hat f_{\mcl{O},0}(x))=0\qquad \tx{for $|x|<\eps_0$},
\end{equation*}
so that
$
\hat f'_{\mcl{O},0}(x)=-\frac{\der_x\iphi^{\mcl{O}}}{\der_y \iphi^{\mcl{O}}}(x, \hat f'_{\mcl{O},0}(x))
$
holds. This expression, combined with \eqref{6-a5}, yields statement (f).
\end{proof}
\end{lemma}

Similarly to \eqref{6-h1-right}, for an admissible solution $\vphi$ corresponding
to $(\iv, \beta)\in \mathfrak{R}_{\rm weak}\cap\{\beta\le \betac^{(\iv)}\}$,
let $\Om$ be its pseudo-subsonic region. Let $\eps_1$ be the constant given in Lemma \ref{lemma-str-nr-sonic}.
For $\eps\in(0,\eps_1]$, define
\begin{equation}
\label{6-h1}
\oOm_{\eps}:=\Om\cap\mcl{N}_{\eps_1}(\leftsonic)\cap\{x<\eps\}.
\end{equation}
Then $\oOm_{\eps}=\oOm_{\eps}\cap\{x>0\}$.

\smallskip
Adjusting the proof of Lemma \ref{lemma-est-nrsonic-N} by using Lemma \ref{lemma-str-nr-sonic} instead of Lemma \ref{lemma1-sonic-N},
we have the following lemma:

\begin{lemma}
\label{lemma-est-nrsonic}
Let $\eps_0$, $\omega_0$, and $\mathfrak{M}$ be three constants in Lemma {\rm \ref{lemma-str-nr-sonic}}.
Then there exist $\bar{\eps}\in(0, \eps_0]$, $L\ge 1$, $\delta\in(0,\frac 12)$, and $\omega\in(0,\omega_0]\cap (0,1)$
depending only on $(\iv, \gam)$ such that any admissible solution $\vphi=\psi+\leftvphi$ corresponding
to $(\iv, \beta)\in \mathfrak{R}_{\rm weak}\cap\{\beta\le \betac^{(\iv)}\}$
satisfies the following properties in $\oOm_{\bar{\eps}}${\rm :}

\smallskip
\begin{itemize}
\item[(a)]
$
\psi_x(x,y)\le \frac{2-\delta}{1+\gam} x\le Lx;
$

\vspace{1mm}
\item[(b)]
$
\psi_x\ge 0
$
 and
 $
|\psi_y(x,y)|\le Lx;
$

\vspace{1mm}
\item[(c)]
$\frac{2}{\mathfrak{M}}(y+\tan{\beta})-\frac{2-\delta}{1+\gamma}x
  \le \der_x(\ivphi-\vphi)(x,y)\le {\mathfrak{M}}$ and $
    \frac{1}{\mathfrak{M}} \le -\der_y(\ivphi-\vphi)\le {\mathfrak{M}};$

\vspace{1mm}
\item[(d)]  There exists a function $\fshocko\in C^1([0,\bar{\eps}])$ such that
    \begin{equation*}
    \begin{split}
    &\oOm_{\bar{\eps}}=\{(x,y)\,:\,x\in(0,\bar{\eps}), 0<y<\fshocko(x)\},\\
    &\shock\cap\der \oOm_{\bar{\eps}}=\{(x,y)\,:\,x\in(0,\bar{\eps}),\,\,y=\fshocko(x)\},\\
    &\omega\le \fshocko'(x) \le L \qquad\text{for}\,\,\, 0<x<\bar{\eps};
    \end{split}
    \end{equation*}

 \item[(e)] $0\le \psi(x,y)\le Lx^2$.
\end{itemize}
\end{lemma}

\begin{lemma}\label{lemma-xybvp-O}
Let $\vphi$ be an admissible solution corresponding to $(\iv, \beta)\in \mathfrak{R}_{\rm weak}\cap\{\beta\le \betac^{(\iv)}\}$.
Let Eq. \eqref{eqn-xy} in $\oOm$ be expressed as
\begin{equation}
\label{eqn-xy-left}
\sum_{i,j=1}^2\hat A_{ij}^{\mcl{O}}(D\psi,\psi,x)D_{ij}\psi+\sum_{i=1}^2\hat A_i^{\mcl{O}}(D\psi,\psi,x)D_i\psi=0,
\end{equation}
with $(D_1, D_2)=(D_x, D_y)$, $D_{ij}=D_iD_j$, and $\hat A^{\mcl{O}}_{21}=\hat A^{\mcl{O}}_{12}$.
Then there exist $\eps_{\mcl{O}}\in(0, \frac{\bar{\eps}}{4}]$ and $\lambda_{\mcl{O}}>0$ depending only on $(\iv, \gam)$
such that, if $(x,y)\in \ol{\oOm_{4\eps_{\mcl{O}}}}$,
\begin{equation}
\label{5-a2-left}
\frac{\lambda_{\mcl{O}}}{2}|{\bm\kappa}|^2\le \sum_{i,j=1}^2 \hat A^{\mcl{O}}_{ij}(D\psi(x,y), \psi(x,y), x) \frac{\kappa_i\kappa_j}{x^{2-\frac{i+j}{2}}}
\le \frac{2}{\lambda_{\mcl{O}}}|{\bm\kappa}|^2\qquad\,\tx{for all}\,\,{\bm\kappa}\in\R^2.
\end{equation}
Moreover, $\mcl{B}_1^{\mcl{O}}$ defined by \eqref{def-B1-new-left} satisfies the following properties{\rm :}

\smallskip
\begin{itemize}
\item[(a)] $\mcl{B}_{1}^{\mcl{O}}({\bf 0},0,x,y)=0$ holds for all $(x,y)\in \R^2${\rm ;}

\smallskip
\item[(b)] For each $k=2,3,\cdots$, there exist constants $\delta_{\rm bc}>0$ and $C>1$ depending only on $(\iv, \gam, k)$
such that, whenever $|(p_x, p_y, z, x)|\le \delta_{\rm bc}$ and $|y-y_{\lefttop}|\le \delta_{\rm bc}$,
\begin{equation*}
|D^k_{(p_x, p_y, z,x,y)}\mcl{B}_{1}^{\mcl{O}}(p_x, p_y, z, x, y)|\le C;
\end{equation*}

\item[(c)] There exist constants $\hat{\delta}_{\rm bc}>0$ and $C>1$ depending only on
$(\iv, \gam)$ such that, whenever $|(p_x, p_y, z, x)|\le \hat{\delta}_{\rm bc}$
and $|y-y_{\lefttop}|\le \hat{\delta}_{\rm bc}$,
\begin{align*}
&D_{(p_x, p_y, z)}\mcl B_{1}^{\mcl{O}}(p_x, p_y, z, x, y)\le -C^{-1} ;
\end{align*}

\item[(d)] There exists a constant $\eps'>0$ depending only on $(\iv, \gam)$,
and constants $\hat{\delta}_{\rm bc}>0$ and $C>1$ in property {\rm (c)}
can be further reduced depending
only on $(\iv, \gam)$ such that, whenever $|(p_x, p_y, z, x)|\le \hat{\delta}_{\rm bc}$
and $|y-y_{\lefttop}|\le \hat{\delta}_{\rm bc}$,
\begin{equation*}
D_{(p_x, p_y)}\mcl B_1^{\mcl{O}}(p_x, p_y,z, x, y)\cdot{\bm \nu}_{\rm sh}^{(x,y)}\ge C^{-1}
\qquad\tx{on $\shock\cap \der{\oOm_{\eps'}}$},
\end{equation*}
where $\shock$ represents the curved shock of the admissible solution,
and ${\bm \nu}_{\rm sh}^{(x,y)}$ is the unit normal vector to $\shock$.
The vector field ${\bm \nu}_{\rm sh}^{(x,y)}$ is expressed in the $(x,y)$--coordinates
and oriented towards the interior of $\Om$.
\end{itemize}
In properties {\rm (b)}--{\rm (d)},
$y_{\lefttop}$ represents the $y$--coordinate of point $\lefttop$,
defined by Definition {\rm \ref{definition-domains-np}}.
\end{lemma}

Even though this lemma is similar to Lemma \ref{lemma-xybvp-N},
the proof is more complicated, because $\leftu, \leftc, \leftvphi$, and $\leftshock$
depend on $\beta\in(0, \betac^{(\iv)}]$.

\begin{proof} We divide the proof into three steps.

\smallskip
{\textbf{1.}} As just mentioned above, $(\leftu, \leftc)$ depend continuously on $\beta\in(0,\frac{\pi}{2})$.
In particular, $|\leftu|$ and $\leftc$  increase with respect to $\beta$.
Therefore, there exists a constant $\bar{c}>1$ depending only on $(\iv, \gam)$ such that
\begin{equation*}
|\leftu|\le \bar c,\,\,\,\, 1 \le \leftc \le \bar c \qquad\,\, \tx{for all $\beta\in[0,\betadet]$}.
\end{equation*}
Then inequality \eqref{5-a2-left} and properties (a)--(b)
can be directly
checked from  \eqref{2-4-a0}, \eqref{prelim5-5}, \eqref{gsh-ext}, \eqref{7-d7-pre-left}, \eqref{def-B1-new-left}, and Lemma \ref{lemma-est-nrsonic}.

\smallskip
{\textbf{2.}}
A direct computation by using \eqref{12-25}--\eqref{2-4-a0}, \eqref{def-densi-ps},
\eqref{def-gsh}, \eqref{7-d7-pre-left}, and \eqref{def-B1-new-left} yields that
\begin{align*}
&\der_z\mcl{B}_1^{\mcl{O}}({\bf 0}, 0, 0, y_{\lefttop})=-\frac{\leftc \iv \sec\beta}{\leftrho^{\gam-2}}\sin (y_{\lefttop}+\beta),\\
&\der_{p_x}\mcl{B}_1^{\mcl{O}}({\bf 0}, 0, 0, y_{\lefttop})=-\leftc(\leftrho-1)\cos^2(y_{\lefttop}+\beta),\\
&\der_{p_y}\mcl{B}_1^{\mcl{N}}({\bf 0}, 0, 0, y_{\lefttop})=
-\Big((\leftrho-1)\sin(y_{\lefttop}+\beta)+\frac{\leftc \iv \sec\beta}{\leftc}\Big)\cos(y_{\lefttop}+\beta).
\end{align*}
For $\beta\le\betac^{(\iv)}$, we have
\begin{equation*}
\cos(\frac{\pi}{2}-\beta-y_{\lefttop})=M_{\mcl{O}}(\beta),
\end{equation*}
where $M_{\mcl{O}}$ is defined by \eqref{1-25}, which is a continuous function
of $\beta\in[0, \frac{\pi}{2})$ that satisfies $M_{\mcl{O}}<1$.
Then there exists a constant $\delta_0\in(0, \frac{\pi}{2})$ depending only on $(\iv, \gam)$
such that $y_{\lefttop}+\beta\le \frac{\pi}{2}-\delta_0$ for all $\beta\in[0, \betac^{(\iv)}]$.
This implies that there exists a constant $m_0>0$ depending only on $(\iv, \gam)$ such that
\begin{equation*}
D_{(p_x, p_y, z)}\mcl{B}_1^{\mcl{O}}({\bf 0}, 0, 0, y_{\lefttop})\le -m_0^{-1}
\qquad \tx{for all $\beta\in(0, \betac^{(\iv)}]$}.
\end{equation*}
We combine this inequality with property (b) to obtain property (c).

\smallskip
{\textbf{3.}}
By \eqref{1-25} and \eqref{app-1}, we have
$$
D_{\bf p}g^{\rm sh}_{\rm mod}(D\leftvphi(\lefttop), \leftvphi(\lefttop),\lefttop)\cdot  {\bm\nu}_{\mcl{O}}
=\leftrho(1-\oM^2)
$$
for the unit normal vector ${\bm\nu}_{\mcl{O}}$ to the straight oblique
shock $\leftshock$ pointing towards the $\xi_1$--axis.
It is shown in the proof of Lemma \ref{lemma:interval-existence} that
\begin{equation*}
\frac{\dd\oM}{\dd\beta}<0\qquad\tx{
for all $\beta\in(0, \frac{\pi}{2})$}.
\end{equation*}
Therefore, there exists a constant $m_1>0$ depending only on $(\iv, \gam)$ such that
\begin{equation*}
\displaystyle{{\bm\nu}_{\mcl{O}}\cdot D_{\bf p}g^{\rm sh}_{\rm mod}(D\leftvphi(\lefttop), \leftvphi(\lefttop),\lefttop)}\ge m_1^{-1}
\qquad\tx{for all $\beta\in (0, \betadet]$.}
\end{equation*}
A direct computation by using \eqref{coord-o}, \eqref{7-d7-pre-left}, and \eqref{cov-left} leads to
\begin{equation}
\label{estimate-B1O-pderiv-lt}
D_{(p_x, p_y)}\mcl B_1^{\mcl{O}}(0,0,0, 0, y_{\lefttop})\cdot{\bm \nu}_{\rm sh}^{(x,y)}(0,y_{\lefttop})
=\displaystyle{{\bm\nu}_{\mcl{O}}\cdot D_{\bf p}g^{\rm sh}_{\rm mod}(D\leftvphi(\lefttop), \leftvphi(\lefttop),\lefttop)}\ge m_1^{-1}.
\end{equation}
Owing to \eqref{estimate-B1O-pderiv-lt} and property (b),
there exist small constants $\hat{\delta}_{\rm bc}>0$ and $\hat{\delta}_{\bm\nu}>0$ depending only on $(\iv, \gam)$ such that, whenever
\begin{equation*}
|(p_x, p_y, z, x)|\le \hat{\delta}_{\rm bc},\quad
 |y-y_{\lefttop}|\le \hat{\delta}_{\rm bc},\quad
 |{\bm\nu}_{\rm sh}^{(x,y)}-{\bm \nu}_{\rm sh}^{(x,y)}(0,y_{\lefttop})|\le \hat{\delta}_{\bm\nu},
\end{equation*}
we have
\begin{equation*}
 D_{(p_x, p_y)}\mcl B_1^{\mcl{O}}(p_x, p_y,z, x, y)\cdot{\bm \nu}_{\rm sh}^{(x,y)}\ge
 \frac{1}{4m_1}.
 \end{equation*}

Note that
%
$\displaystyle{
{\bm\nu}_{\rm sh}^{(x,y)}=\frac{D_{(x,y)}(\ivphi-\leftvphi-\psi)}{|D_{(x,y)}(\ivphi-\leftvphi-\psi)|}
}$ on $\shock\cap \der\oOm$.
Therefore, we can choose a small constant $\eps'>0$ depending only on $(\iv, \gam)$ so that,
by properties (a)--(b) of Lemma \ref{lemma-est-nrsonic},
$\displaystyle{|{\bm\nu}_{\rm sh}^{(x,y)}-{\bm \nu}_{\rm sh}^{(x,y)}(0,y_{\lefttop})|\le \hat{\delta}_{\bm\nu}}$
on $\shock\cap \der \oOm_{\eps'}$. This completes
the proof of property (d) of  Lemma \ref{lemma-xybvp-O}.
\end{proof}

\begin{proposition}\label{lemma-est-sonic-general}
Let $\bar{\eps}>0$ be the constant introduced in Lemma {\rm \ref{lemma-est-nrsonic}}.
Fix $\sigma\in(0, \betac^{(\iv)})$. For each $\alp\in(0,1)$,
there exist $\eps\in(0, \bar{\eps}]$ depending only on $(\iv, \gam,  \sigma)$,
and $C>0$ depending only on  $(\iv, \gam,  \alp)$ such that
any admissible solution $\vphi$ corresponding to
$(\iv, \beta)\in \mathfrak{R}_{\rm weak}\cap\{\beta\le \betac^{(\iv)}-\sigma\}$ satisfies
\begin{equation}
\label{par-est1}
\|\vphi-\leftvphi\|^{\rm (par)}_{2,\alp,\oOm_{\eps}}+\|\fshocko-\hat{f}_{\mcl{O},0}\|^{\rm (par)}_{2,\alp,(0,\eps)}
\le C.
\end{equation}

\begin{proof}
For each $\beta\in(0,\betac^{(\iv)}]$, point $\lefttop$ defined by Definition \ref{definition-domains-np} satisfies
\begin{equation}
\label{expression-ylt}
\sin y_{\lefttop}=\frac{\xi_2^{\lefttop}}{\leftc}.
\end{equation}
In the proof of Lemma \ref{lemma:interval-existence}, it is shown that $\xi_2^{\lefttop}$ is a decreasing function of $\beta\in(0, \betac^{(\iv)}]$
with $\xi_2^{\lefttop}=0$ at $\beta=\betac^{(\iv)}$, and $\leftc$ is an increasing function of $\beta$.
Therefore, for each $\sigma\in(0, \betac^{(\iv)})$, there exists  a constant $\delta_1>0$ depending only
on $(\gam,\ic, \sigma)$ such that
$y_{\lefttop}\ge \delta_1$ for all $\beta\in(0, \betac^{(\iv)}-\sigma]$.
By combining this estimate with Lemma \ref{lemma-est-nrsonic}(d),
we obtain a constant $\lsonic>0$ depending only on $(\iv, \gam,  \sigma)$ such
that any admissible solution $\vphi$
corresponding to $(\iv, \beta)\in \mathfrak{R}_{\rm weak}\cap\{\beta\le \betac^{(\iv)}-\sigma\}$ satisfies
\begin{equation}
\label{fsonico-lwbd-lsonic}
\fshocko\ge \lsonic\qquad\tx{on}\,\,[0,\bar{\eps}].
\end{equation}

We choose
\begin{equation*}
\eps_*=\min\{\frac{\bar{\eps}}{2}, \lsonic^2\}.
\end{equation*}
Then we repeat the proof of Proposition \ref{lemma-est-sonic-general-N} to find a constant $\eps\in[0,\eps_*]$
depending only on $(\iv, \gam)$ such that any admissible solution $\vphi$ corresponding
to $(\iv, \beta)\in \mathfrak{R}_{\rm weak}\cap\{\beta\le \betac^{(\iv)}-\sigma\}$
satisfies estimate \eqref{par-est1} for a constant $C>0$ depending only on $(\iv, \gam,  \alp)$.

The main difference from the proof of Proposition \ref{lemma-est-sonic-general-N} is that the uniform positive lower bound
of $\fshocko$ for admissible solutions corresponding to  $(\iv, \beta)\in \mathfrak{R}_{\rm weak}\cap\{\beta\le \betac^{(\iv)}-\sigma\}$
depends on $\sigma\in(0,\betac^{(\iv)})$ so that the choice of $\eps$ to satisfy
estimate \eqref{par-est1} becomes dependent on  $\sigma$ as well, due to Theorem \ref{elliptic-t7-CF2}.
\end{proof}
\end{proposition}

\begin{remark}
\label{remark-lsonic}
Note that $\xi_2^{\lefttop}$ depends on $\beta\in[0, \frac{\pi}{2})$ continuously.
Furthermore, $\xi_2^{\lefttop}>0$ for $\beta<\betasonic$, and $\xi_2^{\lefttop}=0$ for $\beta\ge \betasonic$.
Since
\begin{equation}
\label{limit-xi2-lt}
\displaystyle{\lim_{\beta\to\betac^{(\iv)}}\xi_2^{\lefttop}=0},
\end{equation}
we have
\begin{equation*}
\lsonic=0\qquad\tx{ at $\beta=\betac^{(\iv)}$}
\end{equation*}
for constant $\lsonic$ from \eqref{fsonico-lwbd-lsonic}.
\end{remark}

\subsection{Case {\rm 2}{\rm :} Admissible solutions for $\beta<\betac^{(\iv)}$ close to $\betac^{(\iv)}$}
\label{subsec10-3}
Now we extend Proposition \ref{lemma-est-sonic-general} up to $\beta=\betac^{(\iv)}$.

\begin{proposition}\label{proposition-sub8}
Let $\bar{\eps}>0$ be the constant introduced in Lemma {\rm \ref{lemma-est-nrsonic}}.
For each $\alp\in(0,1)$, there exist $\eps\in(0,\bar{\eps}]$
and $\sigma_1\in(0,1)$ depending only on $(\iv, \gam)$,
and $C>0$ depending only on $(\iv, \gam,  \alp)$ such that any admissible solution $\vphi=\psi+\leftvphi$
corresponding to $(\iv, \beta)\in \mathfrak{R}_{\rm weak}\cap\{\betac^{(\iv)}-\sigma_1\le \beta<\betac^{(\iv)}\}$
satisfies estimate \eqref{par-est1}.
\end{proposition}

\begin{proof}  We divide the proof into five steps.

\smallskip
{\textbf{1.}}
Owing to Remark \ref{remark-lsonic}, we cannot apply Theorem \ref{elliptic-t7-CF2} directly
to establish estimate \eqref{par-est1} up to $\beta=\betac^{(\iv)}$.
We first observe that there exists a constant $k>1$ depending only on $(\iv, \gam)$ such that,
for any admissible solution corresponding to $(\iv, \beta)\in \mathfrak{R}_{\rm weak}\cap \beta<\betac^{(\iv)}\}$,
\begin{equation}
\label{10-c2}
\{0<x<2\bar{\eps},\;\; 0<y<y_{\lefttop}+\frac{x}{k}\}\subset\oOm_{2\bar{\eps}}\subset\{0<x<2\bar{\eps},\;\; 0<y<y_{\lefttop}+k x\}.
\end{equation}
Using \eqref{10-c2} and Lemmas \ref{lemma:interval-existence} and \ref{lemma-est-nrsonic},
we can adjust the proof of Proposition \ref{lemma-est-sonic-general-N}
to conclude that, for each $\alp\in(0,1)$,
there exist a small constant $\sigma^*>0$ depending on $(\iv, \gamma)$ and a constant $C>0$ depending on
$(\iv, \gam,  \alp)$ such that
any admissible solution $\vphi$ corresponding to $(\iv, \beta)\in \mathfrak{R}_{\rm weak}\cap\{\betasonic-\sigma^*\le \beta<\betac^{(\iv)}\}$
satisfies
\begin{equation}
\label{par-est-left-nr-ns}
\|\vphi-\leftvphi\|^{\rm (par)}_{2,\alp,\oOm_{{y_{\lefttop}^2}}}\le C.
\end{equation}

\smallskip
{\textbf{2.}}
\emph{Claim{\rm :} There exist $\hat{\eps}\in(0, \frac{\bar{\eps}}{2}]$, $\sigma'\in(0,\sigma^*]$,
and $C^*>0$
depending only on $(\iv, \gam)$ such that any admissible solution $\vphi=\psi+\leftvphi$
corresponding to $\mathfrak{R}_{\rm weak}\cap\{\betac^{(\iv)}-\sigma'\le \beta<\betac^{(\iv)}\}$ satisfies}
\begin{equation}
\label{10-b6}
0\le \psi(x,y)\le C^*x^4\qquad \tx{in $\oOm_{2\hat{\eps}}\cap\{x>\frac{y_{\lefttop}^2}{10}\}$}.
\end{equation}

In what follows, unless otherwise specified, the universal constant $C$ represents a positive constant
depending only on $(\iv, \gam)$, which may be different at each occurrence.

\smallskip
For an admissible solution $\vphi$ corresponding to $(\iv, \beta)\in \mathfrak{R}_{\rm weak}\cap \{\beta<\betac^{(\iv)}\}$,
let $\psi$ be given by \eqref{psi-def}. We regard Eq. \eqref{eqn-xy-left} (or equivalently, \eqref{eqn-xy-left}) as a linear
equation for $\psi$ in $\oOm_{\bar{\eps}}$, and represent it as
\begin{equation}
\label{10-c3}
\mcl{L}\psi:=\sum_{i,j=1}^2 a_{ij}(x,y)D_{ij}\psi +\sum_{i=1}^2 a_i(x,y) D_i\psi=0,
\end{equation}
with
$(a_{ij}, a_i)(x,y)=(\hat{A}^{\mcl{O}}_{ij}, \hat{A}^{\mcl{O}}_i)(D\psi(x,y),\psi(x,y),x)$ for $i,j=1,2$,
where $\hat{A}^{\mcl{O}}_{ij}$ and $\hat{A}^{\mcl{O}}_i$ are from Lemma \ref{lemma-xybvp-O}.
By  \eqref{prelim5-5} and Lemma \ref{lemma-est-nrsonic}, there exists a constant $C>0$ depending only
on $(\iv, \gam)$ such that $a_{ij}, i,j=1,2$, satisfy
\begin{align}
\label{10-c6}
&x\le a_{11}(x,y)\le 3x,\quad
C^{-1}\le a_{22}(x,y)\le C,\quad |(a_{12}, a_{21})(x,y)|\le Cx \qquad\,\,\tx{in $\oOm_{\bar{\eps}}$},\\
\label{10-c6-lot-new}
&a_1(x,y)\le 0,\quad |a_2(x,y)|\le Cx\qquad\,\,\tx{in $\oOm_{\bar{\eps}}$}.
\end{align}

By properties (a)--(b) and (e) of Lemma \ref{lemma-est-nrsonic},
there exists $\eps_1\in(0,\bar{\eps}]$ such that $\psi$
satisfies the estimates:
\begin{equation*}
|(\psi_x, \psi_y, \psi, x )|\le \frac 12 \min\{\delta_{\rm bc}, \hat{\delta}_{\rm bc}\},\quad
|y-y_{\lefttop}|\le \frac 12 \min\{\delta_{\rm bc}, \hat{\delta}_{\rm bc}\}\qquad\,\,\tx{in $\ol{\oOm_{\eps_1}}$}
\end{equation*}
for constants $(\delta_{\rm bc}, \hat{\delta}_{\rm bc})$ determined in Lemma \ref{lemma-xybvp-O}.
Then the boundary condition \eqref{bc-psi-B1-left} can be written as a linear boundary condition:
\begin{equation}
\label{10-c5}
\mcl{B}_1^L \psi:=b_1(x,y)\psi_x+b_2(x,y)\psi_y+b_3(x,y)\psi=0\qquad\,\tx{on $\shock\cap \der \oOm_{\eps_1}$},
\end{equation}
and Lemma \ref{lemma-xybvp-O} implies
\begin{equation}
\label{10-c7}
-C\le b_j\le -C^{-1} \mbox{\,\, for $j=1,2$},\,\,\quad
(b_1,b_2)\cdot{\bm\nu}_{\rm{sh}}^{(x,y)}\ge C^{-1}
\qquad\,\,\, \tx{on $\shock\cap \der\oOm_{\eps_1}$}.
\end{equation}

By \eqref{10-c2}, we have
\begin{equation}\label{kx-expl}
\oOm_{\bar{\eps}}\subset \{(x,y)\,:\ 0<x<\bar{\eps},\,\, 0<y<y_{\lefttop}+kx\}.
\end{equation}
For constants $m, \mu>1$ to be determined, define a function $v$ by
\[
v(x,y):=(x+m\mu y_{\lefttop}^2)^4-m(x+m\mu y_{\lefttop}^2)^3y^2.
\]
Suppose that
\begin{equation}
\label{condition-sig-eps}
y_{\lefttop}\le \frac{1}{(m\mu)^2},\qquad \hat{\eps} \le \frac 12 \min\{\eps_1, \eps_{\mcl{O}}, \frac{1}{m\mu}\}
\end{equation}
for $\eps_{\mcl{O}}$ from Lemma \ref{lemma-xybvp-O}.
Then a lengthy computation by using \eqref{10-c6} and \eqref{kx-expl} shows that
constants $(m, \mu)$ can be fixed sufficiently large depending only on $(\iv, \gam)$ such that
\begin{equation}\label{comp1}
\begin{split}
v(x,y)\ge \frac 12(x+m\mu y_{\lefttop})^4\qquad&\tx{in}\,\,\oOm_{2\hat{\eps}},\\
\mcl{L} v<0\qquad&\tx{in}\,\,\oOm_{2\hat{\eps}},\\
\mcl{B}_1^L v<0 \qquad&\tx{on}\,\, \shock\cap \der\oOm_{2\hat{\eps}}.
\end{split}
\end{equation}
For detailed calculations  to obtain \eqref{comp1}, we refer to \cite[Lemma 16.4.1]{CF2}.

For $\hat{\eps}:=\frac 12\min\{\eps_1, \eps_{\mcl{O}}, \frac{1}{m\mu}\}$,
we define
\begin{equation*}
a:=\frac{1}{2{\hat{\eps}}^2}\max_{\der\oOm_{2\hat{\eps}}\cap \{x=2\hat{\eps}\}}\psi.
\end{equation*}
Note that, by the strong maximum principle, $a$ is a positive constant.
By Lemma \ref{lemma-est-nrsonic}(e),
$a$ is uniformly bounded above
depending only on $(\iv, \gam)$.

Note that $\psi$ satisfies the boundary conditions \eqref{5-a4-left} on $\der \oOm_{2\hat{\eps}}\setminus (\{x=2\hat{\eps}\}\cup \shock)$.
Since $|y|\le y_{\lefttop}$ on $\leftsonic$ and $\mu>1$, we have
\begin{equation*}
av\ge 0=\psi\qquad\tx{ on $\leftsonic$}.
\end{equation*}
On $\Wedge\cap \der \oOm_{2\hat{\eps}}$, $v_y=0=\psi_y$.

By the maximum principle, we have
\begin{equation*}
\psi\le av\qquad\tx{in $\oOm_{2\hat{\eps}}$},
\end{equation*}
provided that $y_{\lefttop}$ satisfies the inequality that $y_{\lefttop}\le (m\mu)^{-2}$.

By \eqref{expression-ylt} and \eqref{limit-xi2-lt}, there exists $\sigma'\in(0,\sigma^*]$
such that each $y_{\lefttop}$
corresponding to $\beta\in [\betac^{(\iv)}-\sigma', \betac^{(\iv)})$
satisfies  the inequality that $y_{\lefttop}\le (m\mu)^{-2}$.
This verifies the claim.

\smallskip
{\textbf{3.}} Let $\vphi=\psi+\leftvphi$ be an admissible solution corresponding
to $(\iv, \beta)\in \mathfrak{R}_{\rm weak}\cap\{\betac^{(\iv)}-\sigma'\le \beta <\betac^{(\iv)}\}$.
For $z_0=(x_0, y_0)\in \ol{\oOm_{\hat{\eps}}}\cap \{x>\frac{y_{\lefttop}^2}{5}\}$ and $r\in(0,1]$,
define the sets:
\begin{equation*}
\begin{split}
&\til{R}_{z_0,r}:=\{(x,y)\,:\,|x-x_0|<\frac{x_0^{3/2}}{10k} r,\,\, |y-y_0|<\frac{{x_0}}{10k}r\},\\
&R_{z_0,r}:=\til{R}_{z_0,r}\cap \oOm_{2\hat{\eps}}.
\end{split}
\end{equation*}
Here, $R_{z_0,1}$ may intersect with $\shock\cup\Wedge$. However,
if $R_{z_0,1}\cap \shock\neq \emptyset$, then $R_{z_0,1}\cap \Wedge = \emptyset$, and vice versa.
Note that the dimensions of rectangle $\til{R}_{z_0,r}$ are given such that

\smallskip
\begin{itemize}
\item[(i)] the re-scaled function $\psi^{(z_0)}$ defined below satisfies a uniformly elliptic equation,
due to \eqref{5-a2-left}  stated in Lemma \ref{lemma-xybvp-O};

\smallskip
\item[(ii)] $\ol{R_{z_0,1}}$ does not intersect with $\shock$ and $\Wedge$ simultaneously.
\end{itemize}

For $r>0$, define the sets:
\begin{equation*}
\begin{split}
&Q_{r}:=(-r, r)^2,\\
&Q_r^{(z_0)}:=\{(S,T)\in Q_r\,:\,z_0+\frac{\sqrt{x_0}}{10k}(x_0 S, \sqrt{x_0}\,T)\in R_{z_0, r}\}.
\end{split}
\end{equation*}

For $z_0\in \ol{\oOm_{\hat{\eps}}}\cap \{x>\frac{y_{\lefttop}^2}{5}\}$, define
\begin{equation*}
\psi^{(z_0)}(S,T) =\frac{1}{x_0^4}\psi(x_0+\frac{x_0^{3/2}}{10k}S, y_0+\frac{x_0}{10k} T)\qquad\tx{for}\,\, (S,T)\in Q_{1}^{(z_0)}.
\end{equation*}

For constant $L$ from Lemma \ref{lemma-est-nrsonic}, choose a function $\eta\in C^{\infty}(\R)$ such that $0\le \eta \le 1$ with
$\eta(t)=1$ for $|t|\le L$ and $\eta(t)=0$ for $|t|\ge 2L$.
For such a function $\eta$, we define
\begin{equation}
\label{definition-coeffxy-mod-left}
(\hat{A}_{ij}^{\mcl{O},\rm (mod)}, \hat A_{i}^{\mcl{O},\rm (mod)})(p_x, p_y, z, x)
:=(\hat{A}_{ij}^{\mcl{O}}, \hat{A}_{i}^{\mcl{O}})(x\eta(\frac{p_x}{x}), x\eta(\frac{p_y}{x}), x^2 \eta(\frac{z}{x^2}), x).
\end{equation}
Then $(\hat{A}_{ij}^{\mcl{O},\rm (mod)}, \hat A_{i}^{\mcl{O},\rm (mod)}), i,j=1,2$,
satisfy the following lemma, which is a generalization of Lemma  \ref{lemma-coeff-ext-xy-N}:

\begin{lemma}
\label{lemma-coeff-ext-xy-left}
Let $\eps_0>0$ and $L\ge 1$ be the constants from Lemmas {\rm \ref{lemma-str-nr-sonic}}--{\rm \ref{lemma-est-nrsonic}}, respectively.
Then there exist constants $\eps\in(0,\frac{\eps_0}{2}]$ and $C>0$ depending only on $(\iv, \gam)$
such that any admissible solution $\vphi:=\leftvphi+\psi$ corresponding to $(\iv, \beta)\in \mathfrak{R}_{\rm weak}$
satisfies the following equation{\rm :}
\begin{equation}
\label{eqn-xy-left-extension}
\sum_{i,j=1}^2 {\hat A}^{\mcl{O},\rm (mod)}_{ij}(D\psi, \psi, x)D_{ij}\psi
+\sum_{i=1}^2 {\hat A}^{\mcl{O},\rm (mod)}_{i}(D\psi, \psi, x)D_{i}\psi=0 \qquad \tx{in}\,\,\oOm_{\eps},
\end{equation}
with coefficients $(\hat{A}_{ij}^{\mcl{O},\rm (mod)}, \hat{A}_{i}^{\mcl{O},\rm (mod)})$
satisfying the following properties{\rm :}

\smallskip
\begin{itemize}
\item[(a)] $(\hat{A}_{ij}^{\mcl{O},\rm (mod)}, \hat{A}_{i}^{\mcl{O},\rm (mod)})=(\hat{A}^{\mcl{O}}_{ij}, \hat{A}^{\mcl{O}}_{i})$ \\
in $\{(p_x, p_y, z, x)\,:\,|(p_x, p_y)|\le Lx, \,\, |z|\le Lx^2, \,\, x\in(0,\eps)\}$,

\smallskip
\item[(b)] $|(\hat{A}_{11}^{\mcl{O},\rm (mod)}, \hat{A}_{12}^{\mcl{O},\rm (mod)}, \hat{A}_2^{\rm (mod)})(p_x, p_y, z, x)|
\le Cx$ in $\R^2\times \R\times (0,\eps)$,

\smallskip
\item[(c)] $\|(\hat{A}_{22}^{\mcl{O},\rm (mod)}, \hat{A}_1^{\mcl{O},\rm (mod)})\|_{0, \R^2\times \R\times (0,\eps)}\le C$,

\smallskip
\item[(d)] $\|D_{(p_x, p_y, z, x)}(\hat{A}_{ij}^{\mcl{O},\rm (mod)}, \hat A_i^{\mcl{O},\rm (mod)})\|_{0, \R^2\times \R\times (0,\eps)} \le C$.
\end{itemize}
\end{lemma}

Substituting the definition of $\psi^{(z_0)}$ into Eq. \eqref{eqn-xy-left-extension}, we have
\begin{equation}
\label{scaled-eq-left}
\sum_{i,j=1}^2A_{ij}^{(z_0)}(D\psi^{(z_0)}, \psi^{(z_0)}, S, T) D_{ij}\psi^{(z_0)}
+ \sum_{i=1}^2 A_{i}^{(z_0)}(D\psi^{(z_0)}, \psi^{(z_0)}, S, T) D_{i}\psi^{(z_0)}=0\qquad\tx{in $Q_{1}^{(z_0)}$},
\end{equation}
with
\begin{equation*}
\begin{split}
&A_{ij}^{(z_0)}(\pb,z, S)
=x_0^{\frac{i+j}{2}-2}\hat A^{\mcl{O}, {\rm{mod}}}_{ij}(10kx_0^{4-\frac 32}p_1, 10k x_0^3p_2, x_0^4z, x_0+\frac{x_0^{3/2}}{10k}S),\\
&A_i^{(z_0)}(\pb,z, S)
=\frac{x_0^{\frac{i-1}{2}-1}}{10k}\hat A^{\mcl{O}, {\rm{mod}}}_{i}(10kx_0^{4-\frac 32}p_1, 10k x_0^3p_2, x_0^4z, x_0+\frac{x_0^{3/2}}{10k}S).
\end{split}
\end{equation*}

By \eqref{10-b6}, there exists a constant $C>0$ depending only on $(\iv, \gam)$ such that
\begin{equation}
\label{linfty-est-scaled-left}
|\psi^{(z_0)}|\le C\qquad\tx{in $Q_1^{(z_0)}$}
\end{equation}
for all $z_0\in \ol{\oOm_{\hat{\eps}}}\cap \{x>\frac{y_{\lefttop}^2}{5}\}$.

For $\fshocko$ from Lemma \ref{lemma-est-nrsonic}, define
\begin{equation}
\label{definition-scaled-shock-left-nrsonic}
F^{(z_0)}(S):=\frac{10k}{{x_0}}
\Big(\fshocko(x_0+\frac{x_0^{3/2}}{10k}S)-\fshocko(x_0)\Big)\qquad\,\,\tx{for $-1<S<1$}.
\end{equation}
Similarly to \eqref{estimate-scaled-F}, a direct computation by using \eqref{definition-scaled-shock-left-nrsonic}
and Lemma \ref{lemma-est-nrsonic}(d) shows that there exists a constant $C>0$ depending only on $(\iv, \gam)$
so that, for each $z_0=(x_0, \fshocko(x_0))\in  \shock\cap \der{\oOm_{\hat{\eps}}}$,
$F^{(z_0)}$ satisfies
\begin{equation}
\label{estimate-F-left}
F^{(z_0)}(0)=0,\qquad \|F^{(z_0)}\|_{C^1([-1,1])}\le C\sqrt{x_0}.
\end{equation}
However, it follows from $\ivphi-\vphi=0$ on $\shock$ that
\begin{equation}
\label{expression-F-left}
\iphi^{\mcl{O}}(x_0+\frac{x_0^{3/2}}{10k}S, \fshocko(x_0)+\frac{x_0}{10k}F^{(z_0)}(S))
-x_0^4\psi^{(z_0)}(S, F^{(z_0)}(S))=0
\end{equation}
for
$\iphi^{\mcl{O}}$ given by \eqref{definition-iphi-O}.

Similarly to \eqref{bc-almost-tangential-N}, by using Lemmas \ref{lemma-est-nrsonic}--\ref{lemma-xybvp-O},
we can further reduce $\hat{\eps}\in(0,\frac{\bar{\eps}}{2}]$ depending only on $(\iv, \gam)$
so that the boundary condition \eqref{bc-psi-B1-left} can be rewritten as
\begin{equation}
\label{bc-near-leftsonic-xy-new}
\psi_x=b_{\mcl{O}}(\psi_y, \psi, x, y)\qquad\tx{on}\,\,\shock\cap\der\oOm_{2\hat{\eps}},
\end{equation}
where $b_{\mcl{O}}$ satisfies the following properties:
\begin{equation}
\label{regularity-bo-new}
\begin{split}
&b_{\mcl{O}}(0,0,x,y)=0\qquad\qquad\,\,\,\, \tx{in}\,\,\ol{\oOm_{2\hat{\eps}}},\\
&|D^l b_{\mcl{O}}(p_y, z, x, y)|\le C_l\qquad \tx{in $\R\times \R\times \ol{\oOm_{2\hat{\eps}}}$, for $l=1,2,3,\cdots$},
\end{split}
\end{equation}
for $C_l>0$ chosen depending only on $(\iv, \gam,  l)$.

For each $z_0\in \shock\cap\der\oOm_{\hat{\eps}}$, we substitute $\psi^{(z_0)}$ into \eqref{bc-near-leftsonic-xy-new}
to obtain the following boundary condition on $\shock^{(z_0)}=\{T=F^{(z_0)}(S)\,:\,-1<S<1\}$:
\begin{equation}
\label{scaled-bc-nrls-new}
\psi_S^{(z_0)}=
B^{(z_0)}_{\mcl{O}}(\psi_T^{(z_0)}, \psi^{(z_0)}, S, T),
\end{equation}
for $B^{(z_0)}_{\mcl{O}}(\psi_T^{(z_0)}, \psi^{(z_0)}, S, T)$ given by
\begin{equation*}
B^{(z_0)}_{\mcl{O}}(\psi_T^{(z_0)}, \psi^{(z_0)}, S, T):=
\frac{x_0^{-4+3/2}}{10k}b_{\mcl{O}}(10kx_0^3\psi_T^{(z_0)}, x_0^4\psi^{(z_0)}, x_0+\frac{x_0^{3/2}}{10k}S, y_0+\frac{x_0}{10k}T).
\end{equation*}

It can be checked directly from  \eqref{regularity-bo-new} that,
for each  $z_0\in \shock\cap\der\oOm_{\hat{\eps}}$, $B_{\mcl{O}}^{(z_0)}$  satisfies
\begin{equation}
\label{regularity-Bo-new}
\begin{split}
B_{\mcl{O}}^{(z_0)}(0,0,S,T)=0\qquad&\tx{in $\ol{Q_1^{(z_0)}}$},\\
\|\der_{p_T}B_{\mcl{O}}^{(z_0)}(p_T, z,\cdot)\|_{0, \ol{Q_1^{(z_0)}}}\le m_2\sqrt{x_0}\qquad&\tx{for all $(p_T, z)\in \R\times \R$},\\
\|D_{(p_T, z)}B_{\mcl{O}}^{(z_0)}(p_T, z,\cdot)\|_{1, \ol{Q_1^{(z_0)}}}\le m_2\sqrt{x_0}\qquad&\tx{for all $(p_T, z)\in \R\times \R$},
\end{split}
\end{equation}
where $m_2>0$ is a constant depending only on $(\iv, \gam)$.

\smallskip
{\textbf{4.}}
Using \eqref{5-a2-left},
Lemma \ref{lemma-coeff-ext-xy-left},
\eqref{estimate-F-left}, and \eqref{regularity-Bo-new},
we see that Eq. \eqref{scaled-eq-left} and the boundary condition \eqref{scaled-bc-nrls-new} satisfy all the conditions required
to apply Theorem \ref{elliptic-t8-CF2}. Therefore, by \eqref{linfty-est-scaled-left} and Theorem \ref{elliptic-t8-CF2},
there exist $\eps\in(0, \hat{\eps}]$, $\hat{\alp}\in(0,1)$, $C$,
and $\sigma_1\in(0,\sigma']$ depending only on $(\iv, \gam)$ such that
any admissible solution $\vphi=\psi+\leftvphi$ corresponding to
$(\iv,\beta)\in \mathfrak{R}_{\rm weak}\cap\{\betac^{(\iv)}-\sigma_1\le \beta<\betac^{(\iv)}\}$ satisfies
\begin{equation}
\label{estimate-scaled-nrls-new1}
\|\psi^{(z_0)}\|_{1,\hat{\alp}, \ol{Q^{(z_0)}_{3/4}}}
\le C
\qquad\,\,\tx{for all}\,\,z_0\in \shock\cap \der \oOm_{\eps}\cap\{x>\frac{y_{\lefttop}^2}{5}\}.
\end{equation}

To obtain the $C^{1,\hat{\alp}}$--estimate of $F^{(z_0)}$,
we follow the approach given in the latter part of Step 3 in the proof of
Proposition \ref{lemma-est-sonic-general-N}.
Namely, we differentiate \eqref{expression-F-left} with respect to $S$ to obtain
\begin{equation}
\label{expression-Fdiff-new-ls}
(F^{(z_0)})'=-\frac{\sqrt{x_0}\left(\der_x\iphi^{\mcl{O}}(x_S, y_S)-10kx_0^{5/2}\der_S\psi^{(z_0)}(S,T)\right)}
{ \der_y\iphi^{\mcl{O}}(x_S, y_S)-10kx_0^3\der_T\psi^{(z_0)}(S,T)}
\end{equation}
for $\displaystyle{(x_S, y_S):=(x_0+\frac{x_0^{3/2}}{10k}S, \fshocko(x_0)+\frac{x_0}{10k}F^{(z_0)}(S))}$.

Then a direct computation by using Lemma \ref{lemma-str-nr-sonic}(c), \eqref{estimate-scaled-nrls-new1}--\eqref{expression-Fdiff-new-ls},
and the smoothness of $\iphi^{\mcl{O}}$ yields that there exists a constant $C>0$ depending only on $(\iv, \gam)$ such that
\begin{equation}
\label{estimate-F-nrls-scaled-new}
\frac{1}{\sqrt{x_0}}\|F^{(z_0)}\|_{1,\hat{\alp},[-3/4, 3/4]}\le C
\qquad\,\,\tx{for all}\,\,z_0\in \shock\cap \der \oOm_{\eps}\cap\big\{x>\frac{y_{\lefttop}^2}{5}\big\}.
\end{equation}

For higher order derivative estimates of $\psi^{(z_0)}$ and $F^{(z_0)}$,
we follow the bootstrap argument given in the latter part of Step 3 in the proof of
Proposition \ref{lemma-est-sonic-general-N} by using \eqref{estimate-scaled-nrls-new1},
\eqref{estimate-F-nrls-scaled-new}, and Theorem \ref{elliptic-t9-CF2}.
As a result, we find constants $\eps\in(0,\hat{\eps}]$ and
$\sigma_1\in(0,\sigma']$ depending only on $(\iv, \gam)$ such that,
for each $\alp\in(0,1)$, any admissible solution corresponding
to $(\iv, \beta)\in \mathfrak{R}_{\rm weak}\cap\{\betac^{(\iv)}-\sigma_1\le \beta < \betac^{(\iv)}\}$
satisfies
$$
\|\psi^{(z_0)}\|_{2,\alp, \ol{Q^{(z_0)}_{1/2}}}+\frac{1}{\sqrt{x_0}}\|F^{(z_0)}\|_{2,\alp,[-1/2, 1/2]}\le C
\qquad\, \tx{for all $z_0\in \shock\cap \der \oOm_{\eps}\cap\big\{x>\frac{y_{\lefttop}^2}{5}\big\}$},
$$
where the estimate constant $C$ depends only on $(\iv, \gam, \alp)$.

Furthermore, by repeating the argument of Step 4
in the proof of Proposition \ref{lemma-est-sonic-general-N},
it can be shown that, for each $\alp\in(0,1)$,
there exists a constant $C>0$ depending only on $(\iv, \gam, \alp)$ such that
any admissible solution $\vphi=\psi+\leftvphi$ corresponding
to $(\iv,\beta)\in \mathfrak{R}_{\rm weak}\cap\{\betac^{(\iv)}-\sigma_1\le \beta<\betac^{(\iv)}\}$
satisfies
\begin{equation*}
\|\psi^{(z_0)}\|_{2,\alp, \ol{Q^{(z_0)}_{1/2}}}+\frac{1}{\sqrt{x_0}}\|F^{(z_0)}\|_{2,\alp, [-1/2, 1/2]}\le C
\qquad\,\tx{for all $z_0\in \ol{\oOm_{\eps}}\cap \{x>\frac{y_{\lefttop}^2}{5}\}$.}
\end{equation*}

Denote $\mcl{U}_{\eps}:=\oOm_{\eps}\cap\{x>\frac{y_{\lefttop}^2}{5}\}$.
Collecting all the estimates of $\psi^{(z_0)}$ established above, scaling back to the $(x,y)$--coordinates,
and following the argument of Step 3 in the proof of \cite[Proposition 16.4.6]{CF2}, we have
\begin{equation*}
\begin{split}
&\sum_{0\le k+l\le 2} \sup_{z\in \mcl{U}_{\eps}}
\left(x^{\frac{3k}{2}+l-4}|\der_x^k\der_y^l\psi(z)|\right)\\
&\,\,+\sum_{k+l=2} \sup_{{z,\til z\in \mcl{U}_{\eps},}\atop{z\neq \til z}}
\Big(
\min\{x^{\frac 32(\alp+k)+l-4},{\til x}^{\frac 32(\alp+k)+l-4}\}
\frac{|\der_x^k\der_y^l\psi(z)-\der_x^k\der_y^l\psi(\til z)|}{\delta_{\alp}^{\rm{par}}(z,\til z)}
\Big)\le C,
\end{split}
\end{equation*}
where $k$ and $l$ are nonnegative integers, $C$ is a constant
depending only on $(\iv, \gam, \alp)$,
and we have used the notation that $z=(x,y)$ and $\tilde{z}=(\tilde{x}, \tilde{y})$.
This implies that
\begin{equation}
\label{par-est-left-away-ls}
\|\psi\|_{2,\alp,\oOm_{\eps}\cap \{x>{y_{\lefttop}^2}/{5}\}}^{{\rm (par)}}\le C.
\end{equation}

\smallskip
{\textbf{5.}} Combining estimates \eqref{par-est-left-nr-ns} and \eqref{par-est-left-away-ls} together,
we obtain
\begin{equation*}
\|\vphi-\leftvphi\|^{{\rm (par)}}_{2,0,\oOm_{\eps}}\le C,
\end{equation*}
where constant $C>0$
depends only on $(\iv, \gam,  \alp)$.

In order to estimate $[\vphi-\leftvphi]_{2,\alp,\oOm_{\eps}}^{(2), {\rm (par)}}$,
we consider two cases:
(i) either $z=(x,y), \til{z}=(\til x, \til y)\in \oOm_{y_{\lefttop}^2}$, or $z,\til z\in \oOm_{\eps}\cap\{x>\frac{y_{\lefttop}^2}{5}\}$,
and (ii) $x>y_{\lefttop}^2>\frac{y_{\lefttop}^2}{5}>\til x$.

\smallskip
For $k+l=2$, define
\begin{equation*}
q_{k,l}(z,\til z):=\min\{x^{\alp+k+\frac l2-2},\til x^{\alp+k+\frac l2-2}\}
\frac{|\der_x^k\der_y^l\psi(z)-\der_x^k\der_y^l\psi(\til z)|}{\delta^{\rm (par)}_{\alp}(z,\til z)}.
\end{equation*}
For case (i), $q_{k,l}(z,\til z)$ satisfies
\[
\sum_{k+l=2}q_{k,l}(z,\til z)
\le 4 \Big(\|\psi\|^{{\rm (par)}}_{2,\alp,\oOm_{y_{\lefttop}^2}}+
\|\psi\|^{{\rm (par)}}_{2,\alp,\oOm_{\eps}\cap\{x>{y_{\lefttop}^2}/{5}\}}
\Big).
\]
For case (ii), since $\delta_{\alp}^{\rm (par)}(z,\til z)\ge \frac{x^{\alp}}{2^{\alp}}\ge \frac{\til x^{\alp}}{2^{\alp}}$,
we have
\[
\sum_{k+l=2}q_{k,l}(z,\til z)
\le 2^{\alp+2}\Big(\|\psi\|^{{\rm (par)}}_{2,0,\oOm_{y_{\lefttop}^2}}+
\|\psi\|^{{\rm (par)}}_{2,0,\oOm_{\eps}\cap\{x>{y_{\lefttop}^2}/{5}\}}
\Big).
\]

Therefore, we conclude that there exists a constant $C>0$ depending only on $(\iv, \gam, \alp)$ such that
\begin{equation*}
\|\vphi-\leftvphi\|_{2,\alp, \oOm_{\eps}}^{\rm{(par)}}\le C.
\end{equation*}
In order to estimate $\|\fshocko-\hat{f}_{\mcl{O},0}\|_{2,\alp, (0,\eps)}^{(\rm{par})}$,
we adjust the argument of Step 6 in the proof of Proposition \ref{lemma-est-sonic-general-N}
by using Lemma \ref{lemma-str-nr-sonic}, instead of  Lemma \ref{lemma1-sonic-N}.
\end{proof}

\smallskip
\subsection{Case {\rm 3:} Admissible solutions for $\beta\ge \betac^{(\iv)}$ close to $\betac^{(\iv)}$}
\label{subsubsec-apriori-est-case2}
\begin{lemma}[Extension of Lemma \ref{lemma-str-nr-sonic} for all $\beta\in(0, \betadet)$]
\label{lemma-10-2}
For the $(x,y)$--coordinates given by \eqref{coord-o}, define
\begin{equation}
\label{definition-x-hat}
\hat x:=
x-x_{\lefttop}.
\end{equation}
Then there exist positive constants $\eps_1, \eps_0$, $\delta_0$, $\omega_0$, $C$, and $\mathfrak{M}$
depending only on $(\iv, \gam)$ with $\eps_1>\eps_0$ and $\mathfrak{M}\ge 2$ such that
Lemma {\rm \ref{lemma-str-nr-sonic}} holds for any admissible solution corresponding
to $(\iv, \beta)\in \mathfrak{R}_{\rm weak}$,
where $x$ is replaced by $\hat x$ in all the properties stated in Lemma {\rm \ref{lemma-str-nr-sonic}}.

\begin{proof}
By the definition of $\lefttop$ given in Definition {\rm \ref{definition-domains-np}},
$x_{\lefttop}=0$ for $\beta\le \betac^{(\iv)}$, which implies that $\hat x=x$ for $\beta\le \betac^{(\iv)}$.
Therefore, Lemma {\rm \ref{lemma-10-2}} coincides with Lemma {\rm \ref{lemma-str-nr-sonic}} for $\beta\le \betac^{(\iv)}$.

For $\beta>\betac^{(\iv)}$,  $\hat x<x$, since $x_{\lefttop}>0$.

For $\beta>\betasonic$, we repeat the proof of Lemma \ref{lemma-str-nr-sonic}, except for
replacing $\leftc$ by $|\lefttop\Oo|=\leftc \oM \csc \beta$ for $\oM$ defined by \eqref{1-25}.
Note that $\frac{|\lefttop\Oo|}{\leftc}=\oM\csc\beta=1$ at $\beta=\betasonic$.
Since $\oM$ is decreasing with respect to $\beta$ by \eqref{oM-monotonicity},
we see that $\frac{\dd \oM\csc\beta}{\dd\beta}\le 0$ for $\beta\in(0, \frac{\pi}{2})$ as well.
Then we conclude that
$0<\oM\csc\beta|_{\beta=\betadet}\le \oM\csc\beta<1$ for $\beta>\betasonic$ with $\oM\csc\beta=1$ at $\beta=\betasonic$,
and $|\overline{\lefttop\Qo}|>0$ depends continuously on
$\beta\in(0, \betadet]$.
Therefore, there exists a constant $\eps_1>0$ depending only on $(\iv, \gam)$ such that
$$
|\overline{\lefttop\Qo}|\ge 2\eps_1
\qquad \mbox{for all $\beta\in(0,\betadet]$}.
$$
Then we can also choose a constant $\delta_0>0$, depending only on $(\iv, \gam)$,
to satisfy \eqref{yQo-choice-new} for all $\beta\in (0, \betadet)$.
The rest of the proof is the same as for the case $\beta\le \betasonic.$
\end{proof}
\end{lemma}

\begin{lemma}
\label{lemma-oOm-nonempty}
Let $\eps_1$ be the constant introduced in Lemma {\rm \ref{lemma-10-2}}.
For $\eps\in(0,\eps_1)$, let $\oOm_{\eps}$ be given by \eqref{6-h1}.
For each $\sigma\in(0,\betadet-\betac^{(\iv)})$, define a half-open interval $I(\sigma)$ by
\begin{equation}
\label{definition-isigma-new}
I(\sigma):=(0, \betac^{(\iv)}+\sigma].
\end{equation}
Then, for any given $\eps\in(0,\eps_1)$, there exists $\sigma>0$ depending only on $(\iv, \gam, \eps)$
such that, for any admissible solution $\vphi$ corresponding
to $(\iv, \beta)\in \mathfrak{R}_{\rm weak}\cap\{\beta\in I(\sigma)\}$,
$\oOm_{\eps}$ is nonempty.
\end{lemma}

\begin{proof}
For $\beta\le \betac^{(\iv)}$, $\oOm_{\eps}$ is always nonempty, owing to Proposition \ref{proposition-sub3}.

Suppose that $\beta>\betac^{(\iv)}$.
It follows from Definition \ref{def-regular-sol}(i-4) of Case II,
Proposition \ref{proposition-sub3}, and the definition of the $(x,y)$--coordinates given by \eqref{coord-o}
that $\oOm_{\eps}$ is nonempty if $x_{P_{\beta}}<\eps$.
From this perspective, we need to find a small constant $\sigma>0$
so that $x_{P_{\beta}}<\eps$ holds for all $\beta\in I(\sigma)$.

For each admissible solution $\vphi$, define $M(P):=\frac{|D\vphi(P)|}{c(|D\vphi(P)|^2, \vphi(P))}$;
that is, $M(P)$ is the {\emph{pseudo-Mach number}} of $\vphi$ at point $P$.
For each $\beta\in(0,\frac{\pi}{2})$,
let $P_{\beta}$ be the  $\xin$--intercept $P_{\beta}$ of the straight oblique shock $\leftshock$.
By Definition \ref{def-regular-sol}(ii-3), we have
\begin{equation*}
M(P_{\beta})=\frac{|D\leftvphi(P_{\beta})|}{\leftc}=\oM\csc\beta
\end{equation*}
for $\oM$ given by \eqref{1-25}.
According to the proof of Lemma \ref{lemma:interval-existence},
$\oM$ is a decreasing function of $\beta\in(0, \frac{\pi}{2})$. This implies that
\begin{equation}
\label{10-b4}
\frac{\dd M(P_{\beta})}{\dd\beta}\le 0 \qquad\, \tx{for all}\,\, \beta\in(0, \frac{\pi}{2}),
\end{equation}
so that
\begin{equation}
\label{10-b5}
\begin{split}
&\inf_{\beta\in I(\sigma)}M(P_{\beta})=M(P_{\betac^{(\iv)}+\sigma})<1,\qquad
\lim_{\sigma\to 0+}\inf_{\beta\in I(\sigma)}M(P_{\beta})=1.
\end{split}
\end{equation}
By \eqref{coord-o}, $x_{P_{\beta}}$ can be expressed as
\begin{equation}
\label{definition-xpbeta-new}
x_{P_{\beta}}=\leftc-|D\leftvphi(P_{\beta})|=\leftc\big(1-M(P_{\beta})\big).
\end{equation}
Moreover, we obtain from \eqref{density-mont-ox} and \eqref{10-b4} that
\begin{equation}
\label{10-b9}
\frac{\dd x_{P_{\beta}}}{\dd\beta}>0\qquad\, \tx{for $\beta\in(0, \frac{\pi}{2})$}.
\end{equation}
Furthermore, \eqref{10-b5} yields that
\begin{equation}
\label{10-c1}
\sup_{\beta\in I(\sigma)} x_{P_{\beta}}=x_{P_\beta}|_{\beta=\betac^{(\iv)}+\sigma}>0,\qquad
\lim_{\sigma\to 0+}{\sup_{\beta\in I(\sigma)}}x_{P_{\beta}}=0.
\end{equation}
Therefore, for any given $\eps>0$, we can choose $\sigma>0$ depending only on $(\iv, \gam,  \eps)$
so that $x_{P_{\beta}}<\eps$ for all $\beta\in I(\sigma)$.
\end{proof}

\begin{lemma}[Extension of Lemma \ref{lemma-est-nrsonic} for $\beta>\betasonic$]
\label{proposition-sub7}
Let $\eps_0$, $\om_0$, and $\mathfrak{M}$ be from Lemma {\rm \ref{lemma-10-2}}.
Then there exist constants $\bar{\eps}\in (0, \eps_0]$, $\sigma_2\in (0,1)$, $L\ge 1$,
$\delta\in(0, \frac 12)$, and $\om\in(0,\om_0]\cap (0,1)$ depending only on $(\iv, \gam)$
such that any admissible solution $\vphi=\psi+\leftvphi$ corresponding
to $(\iv, \beta)\in \mathfrak{R}_{\rm weak}\cap\{\beta\in I(\sigma_2)\}$
satisfies properties {\rm (a)}--{\rm (e)} of Lemma {\rm \ref{lemma-est-nrsonic}}
with the following changes{\rm :}

\smallskip
\begin{itemize}
\item[(i)] The definition of $\oOm_{\bar{\eps}}$ in \eqref{6-h1} is replaced by
\begin{equation}
\label{6-h1-mod2}
\oOm_{\bar{\eps}}=\Om\cap\mcl{N}_{\eps}(\ol{\leftsonic})\cap\{x_{\lefttop}<x<x_{\lefttop}+\bar{\eps}\},
\end{equation}

\item[(ii)] $\oOm_{\bar{\eps}}=\{(x,y)\,:\,x\in(x_{\lefttop},x_{\lefttop}+\bar{\eps}), 0<y<\fshocko(x)\}$,

\smallskip
\item[(iii)] $\shock\cap \der\oOm_{\bar{\eps}}=\{(x,y)\,:\,x\in(x_{\lefttop},x_{\lefttop}+\bar{\eps}), y=\fshocko(x)\}$,

\smallskip
\item[(iv)] $\omega\le \fshocko'(x)\le L$ for $x_{\lefttop}<x< x_{\lefttop}+\bar{\eps}$,
\end{itemize}
where $I(\sigma_2)$ is given by \eqref{definition-isigma-new}.

\begin{proof}
As in Lemma \ref{lemma-est-nrsonic}, this lemma is proved by adjusting the proof of Lemma \ref{lemma-est-nrsonic-N}.

Let $\hat{x}$ be given by \eqref{definition-x-hat}.
Since $\hat x=x$ holds for $\beta\le \betasonic$ so that Lemma \ref{proposition-sub7} is the same as Lemma \ref{lemma-est-nrsonic},
it suffices to consider the case that $\beta>\betasonic$.
\smallskip

By Definition \ref{definition-domains-np}, Remark \ref{remark-ellipticity-new2015}, and Proposition \ref{corollary-ellip},
combined with \eqref{def-distbe}--\eqref{def-distfl}, \eqref{eqn-xy}, and \eqref{definition-xpbeta-new},
there exist constants $\sigma'\in (0,1)$, $\eps'\in(0,\eps_0)$, and  $\delta'\in(0,\frac 12)$
depending only on $(\iv, \gam)$ so that any admissible solution corresponding
to $(\iv, \gam)\in \mathfrak{R}_{\rm{weak}}\cap\{\beta\in I(\sigma')\cap[\betasonic, \frac{\pi}{2})\}$
satisfies
\begin{equation}
\label{estimate-principalcof-lsonic-new}
  \begin{split}
     2x-(\gam+1)\psi_x+O_1^{\mcl{O}}(D\psi, \psi, x)
     &\ge \delta'\Big({\rm dist}(\bmxi, \leftsonic)+\leftc\big(1-\frac{|D\leftvphi(\lefttop)|}{\leftc}\big)\Big)\\
     &=\delta'\big((x-x_{\lefttop})+x_{P_{\beta}}\big)
     =\delta' x\qquad\,\,\,\, \tx{in $\oOm_{\eps'}$},
  \end{split}
\end{equation}
where we have used $\lefttop=P_{\beta}$ for $\beta\ge \betasonic$, and \eqref{9-32} in Proposition \ref{corollary-ellip}.

Since $\psi\ge 0$ holds in $\oOm_{\eps_0}$ by Definition \ref{def-regular-sol}(iv), we use \eqref{prelim5-5} to obtain
 \begin{equation*}
 O_1^{\mcl{O}}(D\psi, \psi,x)\le \frac{\gam+1}{\leftc} x\psi_x\qquad \tx{ in $\oOm_{\eps_0}$}.
 \end{equation*}
 Then we can choose $\bar{\eps}\in (0, \eps']$ and $\delta\in(0,\frac 12)$ depending only on $(\iv, \beta)$ so that,
 for any admissible solution $\vphi=\psi+\leftvphi$ corresponding
 to $(\iv, \gam)\in \mathfrak{R}_{\rm{weak}}\cap\{\beta\in I(\sigma')\cap[\betasonic, \frac{\pi}{2})\}$,
 \eqref{estimate-principalcof-lsonic-new} implies that
 \begin{equation*}
  \psi_x(x,y)\le \frac{2-\delta}{1+\gam}x
\end{equation*}
in domain $\oOm_{\bar{\eps}}$ given by \eqref{6-h1-mod2}.

By Lemma \ref{lemma-est-nrsonic}, we can adjust $\delta\in(0, \delta']$ and $\bar{\eps}\in (0, \eps']$
depending only on $(\iv, \gam)$ so that property (a) of Lemma \ref{proposition-sub7} holds
for any admissible solution corresponding to $(\iv, \gam)\in \mathfrak{R}_{\rm{weak}}\cap \{\beta\in I(\sigma')\}$.

Next, we choose a constant $\sigma_2\in(0,\sigma']$ depending only on $(\iv, \gam)$
so that $\oOm_{\bar{\eps}}$ is nonempty for any admissible solution corresponding to
$(\iv, \gam)\in \mathfrak{R}_{\rm{weak}}\cap \{\beta\in I(\sigma_2)\}$.
Such a constant $\sigma_2$ can be chosen due to Lemma \ref{lemma-oOm-nonempty}.
Then  property (a) of Lemma \ref{proposition-sub7} is verified.

The proofs of properties (b)--(e) of Lemma \ref{proposition-sub7} for $\beta>\betasonic$
are the same as for the case that $\beta\le \betasonic$, except that $x$ is replaced by $\hat{x}$
for the range of variables for which the lemma holds,
and Lemma \ref{lemma-10-2} is applied instead of Lemma \ref{lemma-str-nr-sonic}.
More details for proving (b)--(e) of this lemma can be given by adjusting
the proof of Lemma \ref{lemma-est-nrsonic-N}.
\end{proof}
\end{lemma}

\begin{lemma} \label{lemma-10-3}
For each $\sigma\in(0, \betadet - \betasonic)$,
there exists a constant $\mu_0>0$ depending only on $(\iv, \gam, \sigma)$
such that, for any $\beta\in [\betasonic, \betadet-\sigma]$, $g^{\rm sh}_{\rm mod}$
defined by \eqref{gsh-ext} satisfies the following properties{\rm :}
\begin{align*}
&\der_{p_j}g^{\rm sh}_{\rm mod}(D\leftvphi(\lefttop), \leftvphi(\lefttop),\lefttop)\le -\mu_0\qquad\,\tx{for}\,\,j=1,2,\\
&\der_z g^{\rm sh}_{\rm mod}(D\leftvphi(\lefttop), \leftvphi(\lefttop),\lefttop)\le -\mu_0.
\end{align*}

\begin{proof}
Since $P_{\beta}=\lefttop$ for $\beta\ge \betasonic$ due to \eqref{10-a3} in Definition \ref{definition-domains-np},
we apply Lemma \ref{lemma2-appendix} to obtain
$$
\der_{p_1}g^{\rm sh}_{\rm mod}(D\leftvphi(\lefttop), \leftvphi(\lefttop), \lefttop)
\le -C^{-1}\qquad\tx{ for any $\beta\in[\betasonic, \betadet-\sigma]$},
$$
with a constant $C>1$ depending only on $(\iv, \gam,  \sigma)$.

A direct computation by using $\der_{\etan}\leftvphi(\lefttop)=\der_{\etan}\leftvphi(P_{\beta})=0$,  \eqref{12-25},
Definition \ref{definition-domains-np}, and \eqref{app-1} yields that
\begin{equation*}
  \der_{p_2}g^{\rm sh}_{\rm mod}(D\leftvphi(\lefttop), \leftvphi(\lefttop),\lefttop)=-(\leftrho+1)\cos\beta.
\end{equation*}

By using \eqref{new-density}, it can be directly checked that
\begin{equation*}
\der_{z}g^{\rm sh}_{\rm mod}(D\leftvphi(\lefttop), \leftvphi(\lefttop),\lefttop)=-\frac{\leftc \oM}{\leftrho^{\gam-2}}
\end{equation*}
for $\oM>0$ given by \eqref{1-25}.

Since $(\leftrho, \leftc, \oM)$ depend continuously on $\beta\in[0, \frac{\pi}{2})$,
we conclude that there exists a constant $C>1$ depending only on $(\iv, \gam)$ such that
\begin{equation*}
 (\der_{p_2},\der_z)g^{\rm sh}_{\rm mod}(D\leftvphi(\lefttop), \leftvphi(\lefttop),\lefttop)
 \le - C^{-1}\qquad\tx{for all $\beta\in[\betasonic, \betadet]$}.
 \end{equation*}
\end{proof}
\end{lemma}

\begin{corollary}\label{corollary-10-3}
Let $\bar{\eps}$ and $\sigma_2$ be the constants in Lemma {\rm \ref{proposition-sub7}}.
Then Lemma {\rm \ref{lemma-xybvp-O}} holds for all $(\iv, \beta)\in\mathfrak{R}_{\rm weak}$
with $\beta\in[\betasonic, \betasonic+\sigma_2]$.
\end{corollary}

\begin{proof}
It suffices to check property (c) of Lemma \ref{lemma-xybvp-O} for $\beta\ge \betasonic$,
as the rest of the properties of Lemma \ref{lemma-xybvp-O}
can be verified for $\beta\ge \betasonic$ in the same way
as for the case that $\beta<\betasonic$.
Since $\lefttop=P_{\beta}$ for $\beta\ge \betasonic$, $y_{\lefttop}=0$.
From \eqref{12-25} and \eqref{def-B1-new-left}--\eqref{cov-left}, we have
\begin{equation*}
(D_{p_x}, D_{p_y})\mcl{B}_1^{\mcl{O}}(0,0,0,x_{\lefttop}, y_{\lefttop})
=\iv\sec\beta (D_{p_1}, \frac{1}{\leftc} D_{p_2})g^{\rm sh}_{\rm mod}(D\leftvphi(\lefttop), \leftvphi(\lefttop), \lefttop).
\end{equation*}
Then property (c) of Lemma \ref{lemma-xybvp-O} is obtained for the case
that $\betasonic\le \beta\le \betasonic+\sigma_2$ from Lemma \ref{lemma-10-3}
and the smoothness of $\mcl{B}_1^{\mcl{O}}$.
\end{proof}

We now establish the uniform $C^{2,\alp}$--estimate of the admissible solution $\varphi=\psi+\leftvphi$
corresponding to $(\iv, \beta)\in\mathfrak{R}_{\rm weak}$ for $\beta\ge \betasonic$
close to $\betasonic$.

\begin{proposition}\label{proposition-sub9}
Let $\bar{\eps}$ and $\sigma_2$ be the constants in Lemma {\rm \ref{proposition-sub7}}.
Then, for each $\alp\in(0,1)$, there exist constants $\eps\in(0,\bar{\eps}]$ and $\sigma_3\in(0, \sigma_2]$
depending only on $(\iv, \gam)$, and a constant $C>0$ depending only on $(\iv, \gam,  \alp)$ such that
any admissible solution $\vphi=\psi+\leftvphi$ corresponding
to $(\iv, \beta)\in\mathfrak{R}_{\rm weak}\cap\{\betasonic\le \beta\le \betasonic+\sigma_3\}$ satisfies
\begin{equation}
\label{10-d1}
\begin{split}
&\|\psi\|_{C^{2,\alp}(\ol{\oOm_{\eps}})}\le C,\\
&|D^m_{\bmxi}\psi(\lefttop)|=|D^m_{\bmxi}\psi(P_{\beta})|=0\qquad\, \tx{for $m=0,1,2$}.
\end{split}
\end{equation}
Moreover, $\fshocko$ from Lemma {\rm \ref{proposition-sub7}} satisfies
\begin{equation}
\label{10-d2}
\begin{split}
&\|\fshocko-\hat{f}_{\mcl{O},0}\|_{2,\alp,[x_{\lefttop},\eps]}\le C,\\
&\frac{\dd^m}{\dd x^m}(\fshocko-\hat{f}_{\mcl{O},0})(x_{\lefttop})
=\frac{\dd^m}{\dd x^m}(\fshocko-\hat{f}_{\mcl{O},0})(x_{P_{\beta}})=0 \qquad \tx{for $m=0,1,2$}.
\end{split}
\end{equation}

\begin{proof}
In this proof, all the constants are chosen depending only
on $(\iv, \gam)$, unless otherwise specified.

\smallskip
{\textbf{1.}}
For a fixed $\beta\in[\betasonic, \betasonic+\sigma_2]$, define
\begin{equation*}
d_{\rm so}(x):=x-x_{\lefttop}.
\end{equation*}
If $\beta>\betasonic$, then $d_{\rm so}(x)<x$ for all $x\in \oOm_{\bar{\eps}}$.

\medskip
\emph{Claim{\rm :} There exist constants $\eps\in(0, \frac{\bar{\eps}}{2}]$,  $\sigma_3\in(0, \sigma_2]$, and $m>1$
such that any admissible solution $\vphi=\psi+\leftvphi$ corresponding
to $(\iv, \beta)\in\mathfrak{R}_{\rm weak}$ with $\beta\in[\betasonic, \betasonic+\sigma_3]$ satisfies
}
\begin{equation}
\label{10-c8}
\begin{split}
&x_{\lefttop}\le \frac{\eps}{10},\\
&0\le \psi(x,y)\le m \left(d_{\rm so}(x)\right)^5\qquad\tx{in}\,\,\oOm_{2\eps}.
\end{split}
\end{equation}

\smallskip
A more general version of the claim stated immediately above can be found from \cite[Lemma 16.5.1]{CF2}.

\smallskip
Note that $\psi\ge 0$ holds in $\Om$, due to Definition \ref{def-regular-sol}(iv).

For a large constant $M>1$ to be determined later, define
\begin{equation*}
v(x,y):=(x-x_{\lefttop})^5-\frac 1M(x-x_{\lefttop})^3y^2.
\end{equation*}
By Lemma \ref{proposition-sub7}, there exists a constant $k>1$ such that
\begin{equation}
\label{k-subsonic-nr-sonic}
\{(x,y):x_{\lefttop}<x<\bar{\eps},\,\,0<y<\frac 1k(x-x_{\lefttop})\}
\subset \oOm_{\bar{\eps}}
\subset
\{(x,y):x_{\lefttop}<x<\bar{\eps},\,\,0<y< k(x-x_{\lefttop})\}.
\end{equation}

As in the proof of Proposition \ref{proposition-sub8}, we regard $\psi$ as a solution of the linear boundary value problem:
\begin{align*}
\mcl{L}\psi=0\qquad &\tx{in $\oOm_{\bar{\eps}}$},\\
\mcl{B}_1^L\psi=0\qquad &\tx{on $\shock\cap \der\oOm_{\bar{\eps}}$},\\
\psi_y=0\qquad  &\tx{on $\Wedge\cap \der \oOm_{\bar{\eps}}$},
\end{align*}
where the linear operators $\mcl{L}$ and $\mcl{B}_1^L$ are given by \eqref{10-c3} and \eqref{10-c5}, respectively.

It follows from \eqref{prelim5-5} and Lemma \ref{proposition-sub7}
that there exist constants $\hat{\eps}_1\in(0,\bar{\eps}]$ and $C$ depending only on $(\iv, \gam)$
so that the linear operator $\mcl{L}$ satisfies properties \eqref{10-c6}--\eqref{10-c6-lot-new} in $\oOm_{\hat{\eps}_1}$
for any admissible solution corresponding to $(\iv, \beta)\in \mathfrak{R}_{\rm{weak}}$
with $\betasonic\le \beta\le \betasonic+\sigma_2$.

From Corollary \ref{corollary-10-3}, there also exist constants $\hat{\eps}_2\in(0,\hat{\eps}_1]$ and $C$ depending
only on $(\iv, \gam)$ so that the boundary operator $\mcl{B}_1^L$ satisfies \eqref{10-c7} in $\shock\cap\der\oOm_{\hat{\eps}_2}$
for any admissible solution corresponding to $(\iv, \beta)\in \mathfrak{R}_{\rm{weak}}$
with $\betasonic\le \beta\le \betasonic+\sigma_2$.

Similarly to Step 2 in the proof of Proposition \ref{proposition-sub8},
a lengthy computation by using \eqref{10-c6}--\eqref{10-c6-lot-new}
and \eqref{10-c7} shows that
there exist a sufficiently large constant $M>1$,
a sufficiently small constant $\eps\in(0, \frac{\hat{\eps}_2}{2}]$,
and a small constant $\sigma_3\in(0, \sigma_2]$ such that,
for any admissible solution $\vphi=\psi+\leftvphi$ corresponding
to $\beta\in\mathfrak{R}_{\rm weak}$ with $\beta\in[\betasonic, \betasonic+\sigma_3]$, we have
\begin{equation*}
\begin{split}
&x_{\lefttop}\le \frac{\eps}{10},\\
&\mcl{L}v<0\qquad\,\,\,\,\tx{on}\,\,\oOm_{2\eps},\\
&\mcl{B}_1^Lv<0 \qquad\tx{on}\,\,\shock\cap \der\oOm_{2\eps},\\
&v_y=0\qquad\,\,\,\,\,\,\tx{on}\,\,\Wedge\cap \der \oOm_{2\eps},\\
&v(x,y)\ge \frac 12 (x-x_{\lefttop})^5\qquad\tx{in}\,\,\oOm_{2\eps}.
\end{split}
\end{equation*}
Detailed calculations for the results stated above can be obtained by following the arguments in the proof
of \cite[Lemma 16.5.1]{CF2}.

Note that $\sigma_3:=\sigma_3(v_\infty, \gamma,\eps)\in(0, \sigma_2]$ can be chosen sufficiently small so that $\oOm_{2\eps}$ is nonempty
for any admissible solution $\vphi=\psi+\leftvphi$ corresponding
to $\beta\in\mathfrak{R}_{\rm weak}$ with $\beta\in[\betasonic, \betasonic+\sigma_3]$.

For $\eps\in(0, \frac{\hat{\eps}_2}{2}]$ fixed above, define $m_{\psi}$ for \eqref{10-c8} as
\begin{equation*}
m_{\psi}:= \frac{2}{\eps^5}\max_{\der \oOm_{2\eps}\cap \{x=2\eps\}}\psi(x,y).
\end{equation*}
By \eqref{3-c6} stated in Lemma \ref{lemma-step3-1}, there exists a constant $m>0$ depending only on $(\gam, \iv)$ such that
\begin{equation*}
m_{\psi}\le m
\end{equation*}
for any admissible solution $\vphi=\psi+\leftvphi$ corresponding
to $\beta\in\mathfrak{R}_{\rm weak}$ with $\beta\in[\betasonic, \betasonic+\sigma_3]$.
Moreover, we have
\begin{equation*}
\psi(x,y)\le mv(x,y)\qquad\,\tx{on $\der \oOm_{2\eps}\cap\{x=2\eps\}$}.
\end{equation*}
Then the maximum principle implies that
\begin{equation*}
\psi(x,y)\le \frac m2 (x-x_{\lefttop})^5 \qquad\tx{in  $\oOm_{2\eps}$}.
\end{equation*}
The claim is verified.

\smallskip
{\textbf{2.}} Take $\eps>0$ and $\sigma_3>0$ from Step 1.
Let $\vphi=\psi+\leftvphi$ be an admissible solution corresponding
to $(\iv, \beta)\in \mathfrak{R}_{\rm weak}$ with $\beta\in[\betasonic, \betasonic+\sigma_3]$.
For each $r\in(0,1)$ and $z_0=(x_0,y_0)\in\ol{\oOm_{\eps}}\setminus\{\lefttop\}$,
we define $Q_r$ and $Q_r^{(z_0)}$ by
\begin{equation*}
Q_r:=(-r,r)^2,\quad
Q_r^{(z_0)}:=\{(S,T)\in Q_r\,:\,z_0+\frac{\dsonic(x_0)}{10k}(\sqrt{x_0}S, T)\in \oOm_{2\eps}\},
\end{equation*}
and a re-scaled function $\psi^{(z_0)}$ by
\begin{equation*}
\psi^{(z_0)}(S,T):=
\frac{1}{(\dsonic(x_0))^5}
\psi(x_0+\frac{\dsonic(x_0)}{10k}\sqrt{x_0}S, y_0+\frac{\dsonic(x_0)}{10k}T)\qquad\mbox{for $(S,T)\in Q_1^{(z_0)}$},
\end{equation*}
where $k>1$ is the constant from \eqref{k-subsonic-nr-sonic}.

We repeat the arguments used in Steps 3--4 in the proof of Proposition \ref{proposition-sub8}
with some adjustments to obtain that,
for each $\alp\in(0,1)$, there exists a constant $\mathfrak{C}>0$ depending only on $(\iv, \gam,  \alp)$
such that any admissible solution corresponding
to $(\iv,\beta)\in \mathfrak{R}_{\rm weak}$ with $\beta\in[\betasonic, \betasonic+\sigma_3]$ satisfies
\begin{equation}
 \label{estimate-parabolic-subsonic-nrsonic}
\|\psi^{(z_0)}\|_{C^{2,\alp}(\ol{Q^{(z_0)}_{1/10}})}\le \mathfrak{C}\,\,\qquad\tx{for all $z_0\in\ol{\oOm_{\eps}}\setminus \{\lefttop\}$}.
\end{equation}

Following the argument of Step 2 in the proof of \cite[Proposition 16.5.3]{CF2} and using estimate \eqref{estimate-parabolic-subsonic-nrsonic},
we obtain
\begin{equation}
\label{estimate-psi-subsonic-nr-sonic-pre2018}
\begin{split}
&\sum_{0\le k+l\le 2}\sup_{z\in \oOm_{\eps}}\bigl((x-x_{\lefttop})^{k+l-5}x^{\frac k2}|\der_x^k\der_y^l\psi(z)|\bigr)\\
&+\sum_{k+l=2}\sup_{z,\til z\in \oOm_{\eps}, z\neq \til z}
\left(\big(\max\{x,\til x\}-x_{\lefttop}\big)^{k+l+\alp-5}\big(\max\{x,\til x\}\big)^{\frac{k+\alp}{2}}
  \frac{|\der_x^k\der_y^l\psi(z)-\der_x^k\der_y^l\psi(\til z)|}{\delta_{\alp}^{\rm (par)}(z,\til z)}\right)\\
  &\le C\mathfrak{C}
\end{split}
\end{equation}
for $\delta_{\alp}^{\rm(par)}(z, \til z)$ given by Definition \ref{definition-parabolic-norm},
where we have used the notation that $z=(x,y)$ and $\tilde{z}=(\tilde{x}, \tilde{y})$.

We further follow the proof of  \cite[Proposition 16.5.3]{CF2} to obtain that, for all $x,\til x\in (x_{\lefttop}, \eps)$,
\begin{equation}
\label{x-xtil-1}
\begin{split}
&(x-x_{\lefttop})^{k+l-5}x^{\frac k2} \ge x^{\frac 32 k+l-5}\quad\,\,\tx{for $0\le k+l\le 2$},\\[1mm]
&\left(\max\{x,\til x\}-x_{\lefttop}\right)^{k+l+\alp-5} \left(\max\{x,\til x\}\right)^{\frac{k+\alp}{2}}
\ge \left(\max\{x,\til x\}\right)^{\frac 32(k+\alp)+l-5}\quad\,\,\tx{for $k+l=2$}.
\end{split}
\end{equation}

This is because $k+l+\alp-5<0$ holds for $k,l\in \mathbb{Z}^+$ with $0\le k+l\le 2$ and $\alp\in(0,1)$.
Since $\frac 32(k+\alp)+l-5<0$ holds  for $k,l\in \mathbb{Z}^+$ with $0\le k+l\le 2$ and $\alp\in(0,1)$,
it follows from \eqref{x-xtil-1} that
\begin{equation}\label{x-xtil-2}
\begin{split}
&(x-x_{\lefttop})^{k+l-5}x^{\frac k2} \ge \eps ^{\frac 32 k+l-5}\qquad\,\,\tx{for $0\le k+l\le 2$},\\
&\left(\max\{x,\til x\}-x_{\lefttop}\right)^{k+l+\alp-5} \left(\max\{x,\til x\}\right)^{\frac{k+\alp}{2}}
\ge \eps^{\frac 32(k+\alp)+l-5}\qquad\,\,\tx{for $k+l=2$}.
\end{split}
\end{equation}

Assuming that $\eps\le 1$ without loss of generality, we also have
\begin{equation}
\label{delta-alp-zs}
\delta_{\alp}^{\rm (par)}(z,\til z)\le |z-\til z|^{\alp}\qquad\tx{for $z, \til z\in \oOm_{\eps}$}.
\end{equation}

Using \eqref{estimate-psi-subsonic-nr-sonic-pre2018} and \eqref{x-xtil-2}--\eqref{delta-alp-zs}, we obtain
\begin{equation*}
\|\psi\|_{C^{2,\alp}(\ol{\oOm_{\eps}})}\le C
\end{equation*}
for some constant $C>0$ depending only on $(\iv, \gam, \alp)$,
because the choice of $\eps$ given in Step 1 depends only on $(\iv, \gam)$.

Furthermore, it follows directly from \eqref{estimate-psi-subsonic-nr-sonic-pre2018} that
\begin{equation*}
|D^2_{(x,y)}\psi(x,y)|\le C\mathfrak{C} (x-x_{\lefttop})^2\qquad\tx{in $\oOm_{\eps}$},
\end{equation*}
which implies
that
\begin{equation*}
|D^2_{\bm\xi}\psi(\lefttop)|=0.
\end{equation*}
Note that $\psi(\lefttop)=|D_{\bm\xi}\psi(\lefttop)|=0$, due to Definition \ref{def-regular-sol}(ii-3)
for Case 2. Therefore, \eqref{10-d1} is proved.

Finally, \eqref{10-d2} can be proved by adjusting Step 6 in the proof of Proposition \ref{lemma-est-sonic-general-N} and using \eqref{10-d1}.
\end{proof}
\end{proposition}

\subsection{Case {\rm 4:} Admissible solutions for $\beta>\betasonic$ away from $\betasonic$}
\label{subsec10-4}

We first introduce a weighted H\"{o}lder space.

For a bounded connected open set $U\subset\R^2$, let $\Gam$ be a closed portion of $\der U$.
For $\rx,\ry\in U$, define
\begin{equation*}
\delta_{\rx}:=\rm{dist}(\rx,\Gam), \qquad  \delta_{\rx,\ry}:=\min\{\delta_{\rx},\delta_{\ry}\}.
\end{equation*}
For $k\in\R$, $\alp\in(0,1)$, and $m\in \mathbb{Z}^+$, define the standard H\"{o}lder norms by
\begin{align*}
&\|u\|_{m,0,U}:=\sum_{0\le|\betaa|\le m}\sup_{\rx\in U}|D^{\betaa}u(\rx)|,\qquad
[u]_{m,\alp,U}:=\sum_{|\betaa|=m}\sup_{\rx, \ry\in U,\rx\neq  \ry}\frac{|D^{\betaa}u(\rx)-D^{\betaa}u(\ry)|}{|\rx-\ry|^{\alp}},
\end{align*}
and the weighted H\"{o}lder norms by
\begin{align*}
&\|u\|_{m,0, U}^{(k),\Gam}:=\sum_{0\le|\betaa|\le m}\sup_{\rx\in U}\big(\delta_{\rx}^{\max(|\betaa|+k,0)}|D^{\betaa}u(\rx)|\big),\\
&[u]_{m,\alp, U}^{(k),\Gam}:=\sum_{|\betaa|=m}\sup_{\rx,\ry\in U, \rx\neq \ry}\Big(\delta_{\rx,\ry}^{\max\{m+\alp+k,0\}}
    \frac{|D^{\betaa}u(\rx)-D^{\betaa}u(\ry)|}{|\rx-\ry|^{\alp}}\Big),\\
&\|u\|_{m,\alp,U}:=\|u\|_{m,0,U}+[u]_{m,\alp, U},\qquad
\|u\|_{m,\alp,U}^{(k),\Gam}:=\|u\|_{m,0, U}^{(k),\Gam}+[u]_{m,\alp, U}^{(k),\Gam},
\end{align*}
where $D^{\betaa}:=\der_{x_1}^{\beta_1}\der_{x_2}^{\beta_2}$ for $\betaa=(\beta_1,\beta_2)$
with $\beta_j\in\mathbb{Z}_+$ and $|\betaa|=\beta_1+\beta_2$.
Denote
$C^{m,\alp}_{(k),\Gam}(U):=\{u\in C^m(U)\;|\;\|u\|_{m,\alp,U}^{(k),\Gam}<\infty\}$.
%

Let $\sigma_3$ be from Proposition \ref{proposition-sub9}.
Then, by Proposition \ref{corollary-ellip},
there exists $\delta\in(0,1)$ depending only on $(\iv, \gam)$ such that
any admissible solution $\vphi$ corresponding to $(\iv, \beta)\in \mathfrak{R}_{\rm weak}$
with $\betasonic+\frac{\sigma_3}{2} \le \beta<\betadet$ satisfies
\begin{equation}
\label{unif-ellip-subs-awa-sn-new}
\frac{|D\vphi|}{c(|D\vphi|^2, \vphi)}\le 1-\delta \,\qquad \tx{in}\,\, \ol{\Om}\cap \{\xin\le 0\}
\end{equation}
for $c(|{\bf p}|^2,z)$ defined by \eqref{definition-soundspeed}.
By \eqref{unif-ellip-subs-awa-sn-new} and Lemma \ref{lemma-step3-1},
there exists $M_*\ge 2$ depending only on $(\iv, \gam)$ such that
$(D\vphi(\bmxi), \vphi(\bmxi))\in \mcl{K}_{M_*}$
for $\mcl{K}_{M_*}$ defined by \eqref{definition-kmn}.
In particular, there exist $\lambda_*>0$ and $R_*>0$
depending only on $(\iv, \gam)$ such that any admissible solution $\vphi$
corresponding to $(\iv,\beta)\in\mathfrak{R}_{\rm weak}$ with $\betasonic+\frac{\sigma_3}{2} \le \beta<\betadet$ satisfies
\begin{equation*}
\sum_{i,j=1}^2 \der_{p_j}{\mcl{A}}_i(D\vphi(\bmxi), \vphi(\bmxi))\kappa_i\kappa_j\ge \lambda_*|\bm{\kappa}|^2
\end{equation*}
for any $\bmxi\in\ol{\Om}\cap {B_{R_*}(P_{\beta}
)}$ and any ${\bm\kappa}=(\kappa_1, \kappa_2)\in \R^2$.

According to Definition \ref{definition-domains-np}, $P_{\beta}=\lefttop$ for $\beta\ge \betasonic$.
In this chapter, we use $P_{\beta}$, instead of $\lefttop$, to emphasize that $P_{\beta}$ is
the $\xin$--intercept of the straight oblique shock $\leftshock$.
In order to achieve the {\it a priori} estimates of admissible solutions
for $\beta>\betasonic$ away from $\betasonic$, the convexity of the shock polar curves
is heavily used, particularly in establishing the functional independence property
of the boundary conditions for admissible solutions near $P_{\beta}$.

\begin{lemma}\label{lemma-FI}
For each small $\bar{\sigma}\in(0,\frac{\betadet}{10})$,
there exist positive constants $r$ and $M$ depending
only on $(\iv, \gam,  \bar{\sigma})$ such that any admissible solution $\vphi$ corresponding
to $(\iv, \beta)\in\mathfrak{R}_{\rm weak}\cap\{\betasonic\le \beta\le \betadet-\bar{\sigma}\}$ satisfies
\begin{equation*}
\der_{p_1}g^{\rm sh}_{\rm mod}(D\vphi(\bmxi), \vphi(\bmxi), \bmxi)\le -\frac 1M\qquad\,\tx{for all}\,\,{\bmxi}\in \shock\cap B_r(P_{\beta}),
\end{equation*}
where $g^{\rm sh}_{\rm mod}$ is given by \eqref{gsh-ext}.

{\begin{proof} In this proof, all the constants are chosen depending only on
$(\iv, \gam)$, unless otherwise specified. The proof is divided into six steps.

\smallskip
{\textbf{1.}} For $\bmxi\in \R\setminus B_1(\Oi)$, denote $\iu^{(\bmxi)}:=|D\ivphi(\bmxi)|$,
and denote $f_{\rm{polar}}^{(\bmxi)}$ as $f_{\rm{polar}}$ defined by Lemma \ref{lemma-app2}
corresponding to $(\irho, \iu)=(1, |D\ivphi(\bmxi)|)$.
Denote $(\hat u_0^{(\bmxi)}, u_{\rm d}^{(\bmxi)}, u_{\rm s}^{(\bmxi)})$
as $(\hat u_0, u_{\rm d}, u_{\rm s})$ corresponding to $(\irho, \iu)=(1, \iu^{(\bmxi)})$.

Fix $\bar{\sigma}\in (0, \frac{\betadet}{10})$. Let $\vphi$ be an admissible solution corresponding
to $(\iv, \beta)\in\mathfrak{R}_{\rm{weak}}\cap\{\betac^{(\iv)}\le \beta\le \betadet-\bar{\sigma}\}$,
and let $\shock$ be its curved pseudo-transonic shock. By Proposition \ref{proposition-distance},
$f_{\rm{polar}}^{(\bmxi)}$ is well defined for each $\bmxi\in\ol{\shock}$.
For $\bmxi\in \R^2$, denote
\begin{equation}
\label{e-vec}
{\bf e}(\bmxi):=\frac{D\ivphi(\bmxi)}{|D\ivphi(\bmxi)|},
\end{equation}
and let ${\bf e}^{\perp}(\bmxi)$ be the unit vector obtained from rotating ${\bf e}(\bmxi)$ by $\frac{\pi}{2}$ counterclockwise.
More generally, for each ${\bf e}\in \R^2\setminus\{{\bf 0}\}$, let ${\bf e}^{\perp}$ denote the vector obtained from rotating ${\bf e}$ by $\frac{\pi}{2}$ counterclockwise.

The Rankine-Hugoniot condition \eqref{3-a3} implies that $D\vphi(\bmxi)$ can be expressed as
\begin{equation}
\label{d-vphi-rep-shock}
D\vphi(\bmxi)=u{\bf e}(\bmxi)+f_{\rm{polar}}^{(\bmxi)}(u) {\bf e}^{\perp}(\bmxi)\qquad\tx{for each $\bmxi\in \ol{\shock}$},
\end{equation}
with $u=u(D\vphi, \bmxi)$ given by
\begin{equation}
\label{u-vphi-xi-def}
u(D\vphi, \bmxi):=D\vphi(\bmxi)\cdot{\bf e}(\bmxi).
\end{equation}
 By Proposition \ref{corollary-ellip}, we have
\begin{equation}
\label{upper-bound-uvphi}
u(D\vphi,\bmxi)\le u_{\rm s}^{(\bmxi)}\qquad\tx{for each $\bmxi\in \ol{\shock}$}.
\end{equation}

\medskip
{\textbf{2.}} By \eqref{1-24ab} and Lemma \ref{lemma2-appendix}, there exists a constant $M_0>1$ depending only
on $(\iv, \gam, \bar{\sigma})$ such that any admissible solution $\vphi$ corresponding to
$(\iv, \beta)\in \mathfrak{R}_{\rm{weak}}\cap\{\betac^{(\iv)}\le \beta\le \betadet-\bar{\sigma}\}$ satisfies
\begin{equation}
\label{inequaltiy-1}
\der_{p_1}g^{\rm sh}_{\rm mod}(D\vphi(P_{\beta}), \vphi(P_{\beta}), P_{\beta})
=\der_{p_1}g^{\rm sh}(D\leftvphi(P_{\beta}), \ivphi(P_{\beta}), P_{\beta})
\le -M_0^{-1}.
\end{equation}

\begin{figure}[htp]
\centering
\begin{psfrags}
\psfrag{e2}[cc][][0.8][0]{$t_2$}
\psfrag{e1}[cc][][0.8][0]{$t_1$}
\psfrag{sp}[cc][][0.8][0]{$g^{(P_{\beta})}({\bf u})=0$}
\psfrag{ls}[cc][][0.8][0]{$\leftshock$}
\psfrag{gu}[cc][][0.8][0]{$g^{(P_{\beta})}_{\bf u}(D\leftvphi(P_{\beta}))$}
\psfrag{vb}[cc][][0.8][0]{$D\leftvphi(P_{\beta})$}
\psfrag{vi}[cc][][0.8][0]{$D\ivphi(P_{\beta})\phantom{aaaa}$}
\psfrag{x1}[cc][][0.8][0]{$\xi_1$}
\psfrag{x2}[cc][][0.8][0]{$\xi_2$}
\psfrag{pb}[cc][][0.8][0]{$P_{\beta}$}
\psfrag{n}[cc][][0.8][0]{${\bf n}$}
\psfrag{T}[cc][][0.8][0]{$L_{u_*}$}
\includegraphics[scale=1.2]{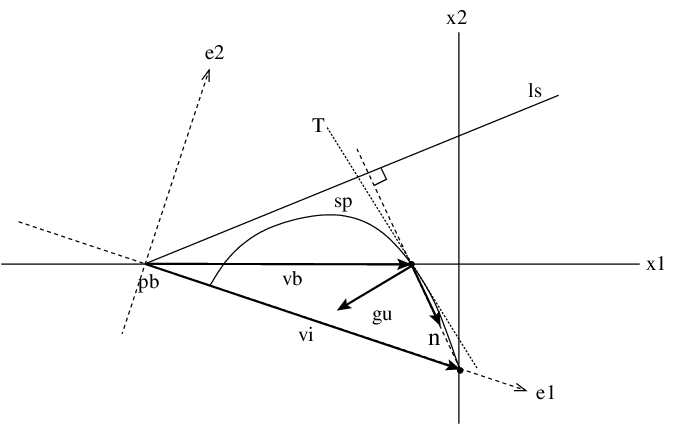}
\caption{The graph of curve $g^{(P_{\beta})}({\bf u})=0$}
\label{fig:relpolar1}
\end{psfrags}
\end{figure}
Let $(t_1,t_2)$--coordinates be given so that $(1,0)_{(t_1,t_2)}={\bf e}(P_{\beta})$ and $(0,1)_{(t_1,t_2)}={\bf e}^{\perp}(P_{\beta})$.
For $\bmxi\in \R^2\setminus B_1(\Oi)$, we define a function $g^{(\bmxi)}({\bf u})$ by
\begin{equation}
\label{definition-g-xi}
g^{(\bmxi)}({\bf u})=g({\bf u})
\end{equation}
for $g({\bf u})$ given by \eqref{g} with ${\bf u}_{\infty}=(|D\ivphi(\bmxi)|, 0)$ (see Fig. \ref{fig:relpolar1}).
If we denote
\begin{equation*}
u_*:={\bf e}(P_{\beta})\cdot D\leftvphi(P_{\beta}),
\end{equation*}
then
\begin{equation*}
D\leftvphi(P_{\beta})=(u_*, f_{\rm polar}^{(P_{\beta})}(u_*)),
\qquad g^{(P_{\beta})}(D\leftvphi(P_{\beta}))=0.
\end{equation*}
Since $D\leftvphi(P_{\beta})\cdot {\bf e}_{\xi_2}=0$,
it can be checked directly from the definitions of $g^{\rm sh}$ and $g$ given in \eqref{def-gsh} and  \eqref{g}, respectively,
that
\begin{equation}
\label{gu-deriv-expre}
g^{(P_{\beta})}_{{\bf u}}(D\leftvphi(P_{\beta}))\cdot {\bf e}_1=\der_{p_1}g^{\rm sh}(D\leftvphi(P_{\beta}), \ivphi(P_{\beta}), P_{\beta}).
\end{equation}
Moreover, from \eqref{inequaltiy-1}, we obtain
\begin{equation}
\label{inequality-gu-Pb}
g_{{\bf u}}^{(P_{\beta})}(D\leftvphi(P_{\beta}))\cdot {\bf e}_1 \le -M_0^{-1}.
\end{equation}

Note that $g^{(P_{\beta})}_{\bf u}(D\leftvphi(P_{\beta}))$ is a normal vector of
curve $(u, f_{\rm polar}^{(P_{\beta})}(u))$ at $u=u_*$.
Let $L_{u_*}$ be the tangent line of curve $(u, f_{\rm polar}^{(P_{\beta})}(u))$ at $u=u_*$.
Then $g^{(P_{\beta})}_{\bf u}(D\leftvphi(P_{\beta}))$ is perpendicular to $L_{u_*}$.
Let ${\bf n}_*$ be the unit normal vector to $L_{u_*}$ with ${\bf n}_*\cdot{\bf e}^{\perp}(P_{\beta})>0$.
Then ${\bf n}_*\cdot {\bf n}<0$ for ${\bf n}=\frac{D\ivphi-D\leftvphi}{|D\ivphi-D\leftvphi|}$,
owing to the convexity of curve $(u, f_{\rm polar}^{(P_{\beta})}(u))$.

It follows from \eqref{entropy-condition-gn} that
$g^{(P_{\beta})}_{\bf u}(D\leftvphi(P_{\beta}))\cdot {\bf n}_*=-|g^{(P_{\beta})}_{\bf u}(D\leftvphi(P_{\beta}))|<0$ (see Fig. \ref{fig:relpolar1}).
This implies that
\begin{equation*}
\frac{g^{(P_{\beta})}_{\bf u}(u, f_{\rm polar}^{(P_{\beta})}(u))}{|g^{(P_{\beta})}_{\bf u}(u, f_{\rm polar}^{(P_{\beta})}(u))|}
=\frac{(\frac{\dd}{\dd u}f_{\rm polar}^{(P_{\beta})}(u),-1)}{\sqrt{1+ \big(\frac{\dd}{\dd u}f_{\rm polar}^{(P_{\beta})}(u)\big)^2 }},
\end{equation*}
and
\begin{equation}\label{sgn-comp}
{\rm sgn}\, \big(g^{(P_{\beta})}_{\bf u}(u, f_{\rm polar}^{(P_{\beta})}(u))\cdot {\bf e}_1\big)
={\rm sgn}\,\big(\frac{\dd}{\dd u}f_{\rm polar}^{(P_{\beta})}(u)\big)
\qquad\tx{for $\hat{u}_0^{(P_{\beta})}<u< u_{\infty}^{(P_{\beta})}$},
\end{equation}
where we have continued to work in the $(t_1, t_2)$--coordinates with basis $\{{\bf e}(P_\beta), {\bf e}^\perp(P_{\beta})\}$.

By the convexity of curve $(u, f_{\rm polar}^{(P_{\beta})}(u))$, we have
\begin{equation*}
\frac{\dd^2}{\dd u^2}f_{\rm polar}^{(P_{\beta})}(u)\le 0\qquad\tx{for $\hat{u}_0^{(P_{\beta})}<u<u_{\infty}^{(P_{\beta})}$}.
\end{equation*}
Then, from \eqref{gu-deriv-expre}--\eqref{sgn-comp}, we obtain
\begin{equation*}
g^{(P_{\beta})}_{\bf u}(u, f_{\rm polar}^{(P_{\beta})}(u))\cdot {\bf e_1}\le - M_0^{-1}
\qquad\tx{for $\der_{{\bf e}(P_{\beta})}\leftvphi(P_{\beta})\le u_{\infty}^{(P_{\beta})}$}.
\end{equation*}

Note that $(P_{\beta}, D\leftvphi(P_{\beta}), {\bf e}(P_{\beta}))$ and the shock polar curve $(u, f_{\rm polar}^{(P_{\beta})}(u))$
depend smoothly on $\beta\in[\betac^{(\iv)}, \betadet]$
(for further details, see Lemma \ref{lemma-app2} or \cite[Claim 16.6.7]{CF2}).
Therefore, there exists a small constant $\eps_1>0$ depending only on $(\gam, \iv, \bar{\sigma})$ so that
\begin{equation}\label{g-deriv-Pb-new1}
g^{(P_{\beta})}_{\bf u}(u, f_{\rm polar}^{(P_{\beta})}(u))\cdot {\bf e_1}\le -\frac{1}{2M_0}
\qquad\tx{for $\der_{{\bf e}(P_{\beta})}\leftvphi(P_{\beta})-\eps_1\le u< u_{\infty}^{(P_{\beta})}$},
\end{equation}
where $\beta\in[\betasonic, \betadet-\bar{\sigma}]$.

\medskip
{\textbf{3.}} For $u(D\vphi, \bmxi)$ given by \eqref{u-vphi-xi-def}, we define
\begin{equation}
\label{q-beta-vec-def}
{\bf q}_{\beta}(u(D\vphi,\bmxi)):= u(D\vphi,\bmxi){\bf e}(P_{\beta})+f_{\rm polar}^{(P_{\beta})}(u(D\vphi,\bmxi)){\bf e}^{\perp}(P_{\beta}),
\end{equation}
provided that $\hat u_0^{(P_{\beta})}<u(D\vphi,\bmxi)<u_{\infty}^{(P_{\beta})}$ holds.

By the definitions of $g^{\rm sh}$ and $g^{(P_{\beta})}$ given in \eqref{def-gsh} and  \eqref{definition-g-xi}, respectively, we have
\begin{equation}
\label{d-p1-gsh-qb-rep}
\der_{p_1}g^{\rm sh}({\bf q}_{\beta}(u), \ivphi(P_{\beta}), P_{\beta})
=g_{\bf u}^{(P_{\beta})}({\bf q}_{\beta}(u(D\vphi,\bmxi)))\cdot {\bf{e}}_1.
\end{equation}

Since $\vphi-\ivphi=0$ holds on $\ol{\shock}$, we have
\begin{equation}
\label{inequality-5}
\begin{split}
\der_{p_1}g^{\rm sh}(D\vphi(\bmxi), \vphi(\bmxi), \bmxi)
\le&\, \der_{p_1}g^{\rm sh}({\bf q}_{\beta}(u), \ivphi(P_{\beta}), P_{\beta})\\
&+ |\der_{p_1} g_{\rm mod}^{\rm sh}(D\vphi(\bmxi), \ivphi(\bmxi), \bmxi)
  -\der_{p_1}g_{\rm mod}^{\rm sh}({\bf q}_{\beta}(u), \ivphi(P_{\beta}), P_{\beta}))|,
\end{split}
\end{equation}
where $u=u(D\vphi, \bmxi)$ for $\bm\xi\in \shock$.

\medskip
{\textbf{4.}}  {\emph{Claim{\rm :} There exist a small constant $r_1>0$ and a constant $C>0$ so that,
if $r\in(0, r_1]$ and $\vphi$ is an admissible solution corresponding to $\beta\in[\betasonic, \betadet-\bar{\sigma}]$,
then
\begin{equation}\label{inequality-new-app-a}
\der_{{\bf e}(\bmxi)}\vphi(\bmxi)\ge \der_{{\bf e}(P_{\beta})}\leftvphi(P_{\beta})-\eps_1
\qquad\tx{on}\,\,\ol{\shock}\cap B_{r_1}(P_{\beta})
\end{equation}
for constant $\eps_1>0$ from \eqref{g-deriv-Pb-new1}.
}}

\medskip
Similarly to \eqref{def-cone-cor}, define a cone generated by vectors ${\bf{u, v}}\in \R^2$ by
\[
{\rm cone}({\bf u}, {\bf v}):=\{\alp_1 {\bf u}+\alp_2{\bf v}\,:\, \alp_1, \alp_2\ge 0\}.
\]
For each $\beta\in[\betasonic, \betadet)$, it is clear that
\begin{equation}\label{mon-sub}
{\bf e}(P_{\beta})\in {\rm cone}(\leftvec, -{\bf e}_2)
\end{equation}
for $\leftvec=(\cos\beta, \sin \beta)$ and ${\bf e}_2=(0,1)$.
We also find from \eqref{2-4-b6} that
\begin{equation*}
\leftvec\cdot {\bf e}(P_{\beta})=\frac{\leftc M_{\mcl{O}}\cot\beta}{|D\ivphi(P_{\beta})|}>0,\qquad
-{\bf e}_2\cdot {\bf e}(P_{\beta})=\frac{\iv}{|D\ivphi(P_{\beta})|}>0
\end{equation*}
for $M_{\mcl{O}}$ defined by \eqref{1-25}.
Moreover, $\leftvec\cdot {\bf e}(P_{\beta})$ and $-{\bf e}_2\cdot {\bf e}(P_{\beta})$ depend continuously
on $\beta$.
Thus, there exists a constant $\kappa_0>0$ such that
\begin{equation*}
\underset{\beta\in[\betasonic, \betadet]}{\min}\{\leftvec\cdot {\bf e}(P_{\beta}), \,-{\bf e}_2\cdot {\bf e}(P_{\beta}) \}\ge \kappa_0.
\end{equation*}
Therefore, we can fix a small constant $r_1>0$ so that
\begin{equation}\label{cone-r2}
\min_{\bmxi\in \ol{B_{r_1}(P_{\beta})}}\min\{\leftvec\cdot {\bf e}(\bmxi), \, -{\bf e}_2\cdot {\bf e}(\bmxi) \}\ge \frac{\kappa_0}{2}
\qquad\tx{for all $\beta\in[\betasonic, \betadet]$}.
\end{equation}
By \eqref{mon-sub} and Lemmas \ref{lemma-step3-1}--\ref{lemma-step-2}, there exists a constant $C_{\sharp}>0$ such that
any admissible solution $\vphi$ corresponding to $(\iv, \beta)\in \mathfrak{R}_{\rm weak}\cap\{\betasonic\le \beta\le \betadet-\bar{\sigma}\}$ satisfies
\begin{equation}
\label{inequality-2}
\begin{split}
\der_{{\bf e}(\bmxi)}\vphi(\bmxi)&=\der_{{\bf e}(\bmxi)}(\vphi-\leftvphi)(\bmxi)+\der_{{\bf e}(\bmxi)}\leftvphi(\bmxi)\\
&\ge \der_{{\bf e}(P_{\beta})}\leftvphi(P_{\beta})-C_{\sharp}|\bmxi-P_{\beta}|\qquad \quad\quad \tx{for }\,\,\bmxi\in \ol{\shock}.
\end{split}
\end{equation}
We choose a constant $r_1>0$ depending only on $(\iv, \gam, \bar{\sigma})$ to satisfy
$C_{\sharp}r_1\le \frac{\eps_1}{2}$ so that \eqref{inequality-new-app-a} follows directly
from \eqref{inequality-2}.
The claim is verified.

\medskip
{\textbf{5.}} {\emph{Claim{\rm :} There exists a small constant $r_2\in(0,r_1]$ depending only on $(\gam, \iv, \bar{\sigma})$
so that, if $\vphi$ is an admissible solution corresponding to $\beta\in[\betasonic, \betadet-\bar{\sigma}]$, then
\begin{equation}
\label{inequality-4}
|D\vphi(\bmxi)-{\bf q}_{\beta}(u(D\vphi,\bmxi))|\le C|\bmxi-P_{\beta}|\qquad\,\tx{for all}\,\,\bmxi\in \ol{\shock}\cap B_{r_2}(P_{\beta}).
\end{equation}
}}

Define
\begin{equation*}
\mu_2:=\min_{\beta\in[\betac^{(\iv)}, \betadet]} \big( \iu^{(P_{\beta})}-u_{\rm s}^{(P_{\beta})}\big).
\end{equation*}
Such a constant $\mu_2$ is positive, depending only on $(\iv, \gam)$.
Choose a small constant $\hat{r}_2\in(0, r_1]$ so that $|u_{\rm s}^{(\bmxi)}-u_{\rm s}^{(P_{\beta})}|\le \frac{\mu_2}{4}$
for all $\bmxi\in B_{\hat{r}_2}(P_{\beta})$.
Then we obtain from \eqref{upper-bound-uvphi} and \eqref{inequality-new-app-a} that
\begin{equation}
\label{inequality-new-app-c}
 \der_{{\bf e}(P_{\beta})}\leftvphi(P_{\beta})-\eps_1  \le u(D\vphi,\bmxi)\le \iu^{(P_{\beta})}-\frac{\mu_2}{2}
 \qquad\tx{on}\,\,\ol{\shock}\cap B_{\hat{r}_2}(P_{\beta}).
\end{equation}

By  Lemma \ref{lemma-step3-1}, \eqref{d-vphi-rep-shock}, and \eqref{q-beta-vec-def}, we have
\begin{equation}
\label{estimate-fpolar-related}
|D\vphi(\bmxi)-{\bf q}_{\beta}(u(D\vphi,\bmxi))|\le C \big(|\bmxi-P_{\beta}|+|(f^{(P_{\beta})}_{\rm polar}-f^{(\bmxi)}_{\rm polar})(u(D\vphi, \bmxi))|\big)
\end{equation}
on $\ol{\shock}\cap B_{\hat{r}_2}(P_{\beta})$.

By the continuous dependence of $(\hat u_0^{(\bm\xi)}, u_{\infty}^{(\bm\xi)})$
and the smooth dependence of $f^{(\bmxi)}_{\rm{polar}}(u)$ on $\bmxi\in \R^2\setminus B_1(\Oi)$
for $u\in(\hat u_0^{(\bmxi)}, \iu^{(\bmxi)})$ due to Lemma \ref{lemma-app2},
and by \eqref{inequality-new-app-c} and the
continuous dependence of $P_{\beta}$ on $\beta\in[\beta_{\rm s}^{(\iv)}, \beta_{\rm d}^{(\iv)}]$,
there exist $C>0$ and $r_2\in(0, \hat{r}_2]$
depending only on $(\iv, \gam, \bar{\sigma})$ such that
\begin{equation}
\label{estimate-fpolar-xi}
|(f^{(P_{\beta})}_{\rm polar}-f^{(\bmxi)}_{\rm polar})(u(D\vphi, \bmxi))|
\le C|\bmxi-P_{\beta}|\qquad\, \mbox{on $\ol{\shock}\cap B_{r_2}(P_{\beta})$}.
\end{equation}
Then \eqref{inequality-4} follows directly from \eqref{estimate-fpolar-related}--\eqref{estimate-fpolar-xi}.

\medskip
{\textbf{6.}}
By \eqref{g-deriv-Pb-new1}, \eqref{d-p1-gsh-qb-rep}, and \eqref{inequality-new-app-c}, we have
\begin{equation}\label{spolar-estimate-gpb}
\der_{p_1}g^{\rm sh}({\bf q}_{\beta}(D\vphi,\bmxi), \ivphi(P_{\beta}), P_{\beta})\le -\frac{1}{2M_0}
\qquad\tx{for $\bmxi\in \ol{\shock}\cap B_{r_2}(P_{\beta})$}
\end{equation}
for any admissible solution $\vphi$ corresponding to $(\iv, \beta)\in \mathfrak{R}_{\rm{weak}}\cap\{\betac^{(\iv)}\le \beta \le \betadet-\bar{\sigma}\}$.

By Lemma \ref{lemma-extcoeff-appc2-new}, \eqref{gsh-ext}, and \eqref{inequality-4}, there exists a constant $C_{\rm{polar}}>0$
such that
\begin{equation}
\label{spolar-estimate-gpert}
 \big|\der_{p_1} g_{\rm mod}^{\rm sh}(D\vphi(\bmxi), \ivphi(\bmxi), \bmxi)
  -\der_{p_1}g_{\rm mod}^{\rm sh}({\bf q}_{\beta}(u), \ivphi(P_{\beta}), P_{\beta}))|\le C_{\rm{polar}}|\bmxi-P_{\beta}\big|
\end{equation}
for $\bmxi\in \ol{\shock}\cap B_{r_2}(P_{\beta})$.

Choosing
\begin{equation*}
r:=\min\{r_2, \frac{1}{4M_0 C_{\rm{polar}}}\},
\end{equation*}
we conclude from \eqref{inequality-5} and \eqref{spolar-estimate-gpb}--\eqref{spolar-estimate-gpert} that
\begin{equation*}
\der_{p_1}g^{\rm sh}(D\vphi(\bmxi), \vphi(\bmxi), \bmxi)
\le -\frac{1}{4M_0}\qquad\tx{on $\ol{\shock}\cap B_r(P_{\beta})$}
\end{equation*}
for any admissible solution $\vphi$ corresponding to
$(\iv, \beta)\in \mathfrak{R}_{\rm{weak}}\cap\{\betac^{(\iv)}\le \beta \le \betadet-\bar{\sigma}\}$.
This completes the proof.
\end{proof}
}
\end{lemma}

To simplify notations, let ${\bf e}_{\beta}$ denote ${\bf e}(P_{\beta})$ for each $\beta\in[\betasonic, \betadet)$,
and let ${\bf e}^{\perp}_{\beta}$ be the unit vector obtained from rotating ${\bf e}_{\beta}$ by $\frac{\pi}{2}$ counterclockwise.
By \eqref{def-polar-oi}, \eqref{r-derivative}, and \eqref{e-vec}, we have
\begin{equation*}
\der_{{\bf e}_{\beta}}(\ivphi-\vphi)(\bmxi)\ge d_1+({\bf e}_{\beta}-{\bf e}(\bmxi))\cdot D(\ivphi-\vphi)(\bmxi)
\qquad\tx{for all $\bmxi\in \mcl{N}_{\eps}(\shock)\cap \Om$},
\end{equation*}
where constants $d_1$ and $\eps$ are from \eqref{r-derivative}.
Therefore, we can apply Lemma \ref{lemma-step3-1} to choose a constant $s_*>0$
depending only on $(\iv, \gam)$
such that any admissible solution $\vphi$ corresponding
to $(\iv,\beta)\in \mathfrak{R}_{\rm weak}\cap\{\betasonic\le \beta< \betadet\}$ satisfies
\begin{equation}\label{inequality-6}
\der_{{\bf e}_{\beta}}(\ivphi-\vphi)\ge \frac{d_1}{8}\qquad\tx{in $B_{2s^*}(P_{\beta})\cap \Om$}.
\end{equation}

\begin{definition}
\label{definition-STcoordn}
Introduce the $(S,T)$--coordinates so that

\smallskip
\begin{itemize}
\item[(i)] $P_{\beta}$ becomes the origin in the $(S,T)$--coordinates,

\smallskip
\item[(ii)] ${\bf e}_{\beta}=(1,0)_{(S,T)}$ and ${\bf e}_{\beta}^{\perp}=(0,1)_{(S,T)}$.
\end{itemize}
In fact, the $(S,T)$--coordinates are the same as the $(t_1,t_2)$--coordinates
in Fig. {\rm \ref{fig:relpolar1}}.
\end{definition}

In the $(S,T)$--coordinates given by Definition \ref{definition-STcoordn},
$\leftshock$, $\shock$, $\Wedge$, and $\Om$ near $P_{\beta}$ can be represented as
\begin{equation*}\begin{split}
&\leftshock\cap B_{s^*}(P_{\beta})=\{S= a_{\leftshock}(\beta)T\,:\,T>0 \}\cap B_{s^*}(P_{\beta}),\\[1mm]
& \shock\cap B_{s^*}(P_{\beta})=\{S=f_{\bf e}(T)\,:\,T>0\}\cap B_{s^*}(P_{\beta}),\\[1mm]
&\Wedge\cap B_{s^*}(P_{\beta})=\{S=a_{\rm w}(\beta)T\,:\,T>0\}\cap B_{s^*}(P_{\beta}),\\[1mm]
& \Om\cap B_{s^*}(P_{\beta})=\{(S,T)\,:\,a_{\leftvec}(\beta)T\le f_{\bf e}(T)<S<a_{\rm w}(\beta)T, \;T>0 \}\cap B_{s^*}(P_{\beta}),
\end{split}
\end{equation*}
where $a_{\rm w}(\beta)$ depends continuously on $\beta\in(0, \frac{\pi}{2})$,
and $a_{\leftshock}(\beta)=\tan \theta_{\beta}$ with
$\theta_{\beta}:=\tan^{-1}a_{\rm w}(\beta)-\beta>0$ for each $\beta\in(0, \frac{\pi}{2})$.
Note that there is a constant $C>0$ depending only on $(\iv, \gam)$ such that
$C^{-1}\le a_{\rm w}(\beta)\le C$ for all $\beta\in[\betasonic, \betadet)$.
The representation of $\shock\cap B_{s^*}(P_{\beta})$ as a graph of $S=f_{\bf e}(T)$
is obtained by the implicit function theorem, combined with \eqref{inequality-6}.

\begin{proposition}
\label{lemma-gradient-est}
Let positive constants $\sigma_3$ and $r$ be from Proposition {\rm \ref{proposition-sub9}} and Lemma {\rm \ref{lemma-FI}}, respectively.
For small constants $\sigma_{\rm s}\in(0, \frac{\sigma_3}{2}]$
and $\sigma_{\rm d}\in(0, \frac{\beta_{\rm d}^{(\iv)}}{10})$,
there exist constants $s\in(0,r)$, $\alp\in(0,1)$, and $C>0$ depending only on
$(\iv, \gam,  \sigma_{\rm s}, \sigma_{\rm d})$ such that
any admissible solution $\vphi$ corresponding
to $(\iv, \beta)\in \mathfrak{R}_{\rm weak}\cap\{\betasonic+\sigma_{\rm s}\le \beta\le \betadet-\sigma_{\rm d}\}$
satisfies the estimates{\rm :}
\begin{equation*}
\begin{split}
&\|\vphi\|_{2,\alp, \Om\cap B_s(P_{\beta})}^{(-1-\alp), \{P_{\beta}\}}+
\|f_{\bf e}\|_{2,\alp, (0,s)}^{(-1-\alp),\{0\}}\le C,\\[2mm]
&|D_{\bmxi}^m(\vphi-\leftvphi)(P_{\beta})|=0\qquad\,\,\tx{for $m=0,1$}.
\end{split}
\end{equation*}

\begin{proof}
In this proof, all the estimate constants are chosen depending only on
$(\iv, \gam, \sigma_{\rm s}, \sigma_{\rm d})$, unless otherwise specified.
For fixed $\sigma_{\rm s}\in(0, \frac{\sigma_3}{2}]$ and
$\sigma_{\rm d}\in(0, \frac{\beta_{\rm d}^{(\iv)}}{10})$,
let $\vphi$ be an admissible solution for $\beta\in [\betasonic+\sigma_{\rm s}, \betadet-\sigma_{\rm d}]$.

\smallskip
{\textbf{1.}}
Denote $\bar{\phi}:=\ivphi-\vphi$, and rewrite Eq. \eqref{2-1} and the derivative boundary conditions \eqref{rhbc-gsh} and \eqref{3-a2}
in terms of $\bar{\phi}$ as follows:
\begin{equation}\label{bvp-barphi}
\begin{split}
\sum_{i,j=1}^2 A_{ij}(D\bar{\phi}, \bar{\phi}, {\bm \xi})D_{ij}\bar{\phi}=0\qquad&\tx{in $B_{s^*}(P_{\beta})\cap \Om$},\\
\hat{g}^{\rm sh}(D \bar{\phi}, \bar{\phi}, {\bm\xi})=0\qquad&\tx{on $\shock$},\\
\hat{g}^{\rm w}(D \bar{\phi}, \bar{\phi}, {\bm\xi})=0\qquad&\tx{on $\Wedge$},
\end{split}
\end{equation}
where
\begin{equation}
\label{bvp-coeff}
\begin{split}
&A_{ij}({\bf p},z, {\bm\xi})=\hat{c}^2({\bf p},z,{\bm \xi})\delta_{ij}-(\der_i\ivphi-p_i)(\der_j\ivphi-p_j) \qquad\tx{for $i,j=1,2,$}\\
&\hat{c}^2({\bf p},z,{\bm \xi})=1-(\gam-1)\big(\frac 12|D\ivphi-{\bf p}|^2+\ivphi-z\big),\\
&\hat{g}^{\rm sh}({\bf p},z,{\bm\xi})=-g^{\rm sh}(D\ivphi({\bm\xi})-{\bf p}, \ivphi(\bmxi)-z, {\bm\xi}),\\
&\hat{g}^{\rm w}({\bf p}, z, {\bm\xi})=p_2+(\etan+\iv),
\end{split}
\end{equation}
where $g^{\rm sh}$ is given by \eqref{def-gsh} and $s^*\in(0,r]$ is from \eqref{inequality-6}.

Next, we apply a partial hodograph transform to $\bar{\phi}$ in $B_{s^*}(P_{\beta})\cap \Om$
in the direction of ${\bf e}_{\beta}$.
For each $(S,T)\in B_{s^*}(P_{\beta})\cap \Om$,
define ${\bf y}=(y_1,y_2)=(\bar{\phi}(S,T),T)$.
By \eqref{inequality-6}, there exists a unique function $v({\bf y})$ such that
\begin{equation}
\label{definition-htn}
v(y_1,y_2)=S\qquad\tx{if and only if}\qquad {\bar\phi}(S, y_2)=y_1
\end{equation}
for ${\bf y}\in \mcl{D}_{s^*}^{\beta}:=\{{\bf y}=(\bar{\phi}(S,T),T)\,:\,(S,T)\in B_{s^*}(P_{\beta})\cap \Om\}$.
By taking derivatives of $v(\bar{\phi}(S, y_2), y_2)=S$, it can be directly checked that
\begin{equation}
\label{deriv-hdgraph}
\der_{y_1}v=\frac{1}{\der_S\bar{\phi}},\quad
\der_{y_2}v=-\frac{\der_T\bar{\phi}}{\der_S\bar{\phi}}.
\end{equation}

By Lemma \ref{lemma-step3-1}, \eqref{inequality-6}, and \eqref{definition-htn}--\eqref{deriv-hdgraph},
there exists a constant $K>1$ depending only on $(\gam, \iv)$ such that
\begin{equation}
\label{lip-bd-v}
\frac 1K \le \der_{y_1}v\le \frac{8}{d_1},\quad |v|+|Dv|< 2K\qquad\,\,\,\tx{in $\ol{\mcl{D}_{s^*}^{\beta}}$}.
\end{equation}

Using the definition of $v$, \eqref{bvp-barphi} can be written in terms of $v$:
\begin{equation}
\label{bvp-v}
\begin{split}
\sum_{i,j=1}^2a_{ij}(Dv,v,{\bf y})\der_{y_iy_j}v=0\qquad&\tx{in}\,\, \mcl{D}_{s^*}^{\beta},\\
g^{\rm sh}_{\rm h}(Dv,v,{\bf y})=0\qquad&\tx{on}\;\; \shock^{({\rm h})}= \{{\bf y}=(0,T)\,:\,(S,T)\in B_{s^*}(P_{\beta})\cap\shock\},\\
g^{\rm w}_{\rm h}(Dv,v,{\bf y})=0\qquad&\tx{on}\;\;\Wedge^{({\rm h})}= \{{\bf y}=({\bar\phi}(S,T),T)\,:\,(S,T)\in B_{s^*}(P_{\beta})\cap\Wedge\},
\end{split}
\end{equation}
where $(a_{ij}, g^{\rm sh}_{\rm h}, g^{\rm w}_{\rm h})({\bf p}, z, {\bf y})$ are directly computed by using \eqref{bvp-coeff}
and the definition of $v$.
More precisely, $(a_{ij}, g^{\rm sh}_{\rm h}, g^{\rm w}_{\rm h})({\bf p}, z, {\bf y})$ are given by
\begin{equation*}
\begin{split}
&a_{11}({\bf p}, z, {\bf y})=\frac{1}{p_1^3}(A_{11}-2p_2A_{12}+p_2^2A_{22}),\\
&a_{12}({\bf p}, z, {\bf y})=a_{21}({\bf p}, z, {\bf y})=\frac{1}{p_1^2}(A_{12}-p_2A_{22}),\\
&a_{22}({\bf p}, z, {\bf y})=\frac{1}{p_1}A_{22},\\
&(g_{\rm h}^{\rm sh}, g_{\rm h}^{\rm w})({\bf p}, z, {\bf y})=-(\hat g^{\rm sh}, \hat g^{\rm w}),
\end{split}
\end{equation*}
with
\begin{equation*}
(A_{11}, A_{12}, A_{22}, \hat g^{\rm sh}, \hat g^{\rm w})
=(A_{11}, A_{12}, A_{22}, \hat g^{\rm sh}, \hat g^{\rm w})((\frac{1}{p_1}, -\frac{p_2}{p_1}), y_1, (z, y_2)).
\end{equation*}

From the definition of $a_{ij}$, we find that, for $({\bf p}, z, {\bf y})$ satisfying $p_1\ne 0$,
$$
\sum_{i,j=1}^2
a_{ij}({\bf p}, z, {\bf y})\kappa_i\kappa_j=\frac 1{p_1^3}\sum_{i,j=1}^2
A_{ij}\eta_i\eta_j
$$
for $(\eta_1, \eta_2)=(\kappa_1,p_1\kappa_2-p_2\kappa_1)$, so that
\begin{equation*}
\sum_{i,j=1}^2
a_{ij}(D v, v, {\bf y})\kappa_i\kappa_j=\frac 1{v_{y_1}^3}\sum_{i,j=1}^2
A_{ij}(D\bar{\phi}, \bar{\phi}, S,T)\eta_i\eta_j,
\end{equation*}
where ${\bf y}=(\bar{\phi}(S,T), T)$ and $(\eta_1, \eta_2)=(\kappa_1,v_{y_1}\kappa_2-v_{y_2}\kappa_1)$.  This implies that there is a constant $C>0$ such that
\begin{equation*}
\frac 1C |{\boldsymbol{\kappa}}|^2
\le \sum_{i,j=1}^2
a_{ij}(D v, v, {\bf y})\kappa_i\kappa_j\le C|{\boldsymbol{\kappa}}|^2
\qquad \mbox{for all ${\bf y}\in\mcl{D}_{s^*}^{\beta}, \;{\boldsymbol{\kappa}}\in \R^2$}.
\end{equation*}

Define a set
\begin{equation*}
U:=\{({\bf p}, z, {\bf y})\in \R^2\times \R\times \mcl{D}_{s^*}^{\beta} \}.
\end{equation*}

We fix a cut-off function $\zeta\in C^{\infty}(\R)$ satisfying that
$\zeta(t)\equiv 0$ on $(-\infty, \frac{1}{10K})$ and $\zeta(t)\equiv 1$ on $(\frac{1}{4K}, \infty)$.
Furthermore, we define
\begin{equation*}
(a_{ij}^{\rm mod}, g_{\rm h}^{\rm sh, mod}, g_{\rm h}^{\rm w,\rm mod})({\bf p}, z, {\bf y})
=\zeta(p_1)(a_{ij}, g_{\rm h}^{\rm sh}, g_{\rm h}^{\rm w})({\bf p}, z, {\bf y})\qquad\,\,
\tx{for $i,j=1,2$}.
\end{equation*}
Then \eqref{bvp-v} can be rewritten as
\begin{equation}
\label{bvp-v-mod}
\begin{split}
\sum_{i,j=1}^2a_{ij}^{\rm mod}(Dv,v,{\bf y})\der_{y_iy_j}v=0\qquad&\tx{in}\,\, \mcl{D}_{s^*}^{\beta},\\
g^{\rm sh, mod}_{\rm h}(Dv,v,{\bf y})=0\qquad&\tx{on $\shock^{({\rm h})}$},\\
g^{\rm w, mod}_{\rm h}(Dv,v,{\bf y})=0\qquad&\tx{on $\Wedge^{({\rm h})}$}.
\end{split}
\end{equation}

Furthermore, for any $l=0,1,2,\cdots$, there exists a constant $C_l>0$ depending only on $(\gam, \iv, l)$ such that
\begin{equation}
\label{smooth-mod-coeff-hdgraph}
|D^l_{({\bf p}, z, {\bf y})}(a_{ij}^{\rm mod}, g_{\rm h}^{\rm sh, mod}, g_{\rm h}^{\rm w,\rm mod})|\le C_l\qquad\,\,\tx{on $U$}.
\end{equation}

\smallskip
{\textbf{2.}}
In this step, we apply Proposition \ref{appendix4-proposition1} to obtain
\begin{equation}
\label{inequality-7}
|g_{\rm h}^{\rm w}(Dv({\bf y}), v({\bf y}), {\bf y})-g_{\rm h}^{\rm w}(Dv({\bf 0}), v({\bf 0}), {\bf 0})|
\le C|{\bf y}|^{\alp_1}\qquad\tx{for ${\bf y}\in \ol{\mcl{D}_{s^*}^{\beta}\cap B_{l^*}({\bf 0})}$}
\end{equation}
for some $\alp_1\in(0,1)$, $C>0$, and $l^*>0$.

$\shock^{({\rm h})}$ is flat so that it is $C^2$ up to its endpoints, and $\Wedge^{({\rm h})}$ is Lipschitz continuous
up to its endpoints.
Then we regard $\Wedge^{({\rm h})}$ and $\shock^{({\rm h})}$ as $\Gam^1$ and $\Gam^2$, respectively,
in Proposition \ref{appendix4-proposition1}.
Then $(g_{\rm h}^{\rm w, mod}, g_{\rm h}^{\rm sh, mod}, 0)$ in \eqref{bvp-v-mod} become $(b^{(1)}$, $b^{(2)},h)$
in Proposition \ref{appendix4-proposition1}. It follows directly from \eqref{smooth-mod-coeff-hdgraph}
that \eqref{bvp-v-mod} satisfies conditions \eqref{appendix4-prop1-cond1}--\eqref{appendix4-prop1-cond4}.

Also, \eqref{lip-bd-v} implies that $v$ satisfies condition \eqref{app-D-u-Lip-est} stated in Proposition \ref{appendix4-proposition1}.

A direct computation by using the definition of $v$ in \eqref{definition-htn} yields that
\begin{equation*}
|D_{\bf p} g_{\rm h}^{\rm w}(Dv({\bf y}), v({\bf y}), {\bf y})|
=|\frac{1}{v_{y_1}^2}(v_{y_2}, -v_{y_1})|\ge \frac{1}{|v_{y_1}|}=|\bar{\phi}_S|
\qquad\,\tx{for all}\,\,{\bf y}\in \ol{\mcl{D}_{s^*}^{\beta}}.
\end{equation*}
Thus, \eqref{inequality-6} implies that
\begin{equation*}
|D_{\bf p} g_{\rm h}^{\rm w}(Dv({\bf y}), v({\bf y}), {\bf y})|
\ge \frac{d_1}{8}\qquad\,\tx{for all}\,\,{\bf y}\in \ol{\mcl{D}_{s^*}^{\beta}}.
\end{equation*}
This shows that $b^{(1)}=g_{\rm h}^{\rm w}$ satisfies condition (ii) of Proposition \ref{appendix4-proposition1}.
By \eqref{def-gsh}, \eqref{app-1}, Lemma \ref{lemma-step3-1}, Remark \ref{remark-ellipticity-new2015},
and Proposition \ref{corollary-ellip}, there exists a constant $\lambda_1>0$ depending only
on $(\iv, \gam,  \sigma_{\rm s})$ such that any admissible solution $\vphi$
for $\beta\in[\betasonic+\sigma_{\rm s}, \betadet)$
satisfies
\begin{equation*}
D_{\bf p}g^{\rm sh}_{\rm mod}(D\vphi(\bmxi), \vphi(\bmxi), \bmxi)\cdot{\bm\nu}_{\rm s}(\bmxi)
\ge \lambda_1\,\qquad\tx{for all}\,\,\bmxi\in \ol{\shock}\cap B_{s^*}(P_{\beta}),
\end{equation*}
where $\bm\nu_{\rm s}$ is the unit normal vector to $\shock$ towards the interior of $\Om$.
Then a direct computation by using \eqref{inequality-6} and \eqref{bvp-coeff}--\eqref{definition-htn}
shows that
\begin{equation*}
\der_{p_1}\hat{g}^{\rm sh}(Dv({\bf y}), v({\bf y}), {\bf y})
=|D\bar{\phi}|\bar{\phi}_SD_{\bf p}g_{\rm mod}^{\rm sh}(D\vphi(\bmxi), \vphi(\bmxi), \bmxi)\cdot{\bm \nu}_s(\bmxi)
\ge \lambda_1\big(\frac{d_1}{8}\big)^2
\qquad\tx{on $\shock^{({\rm h})}$}.
\end{equation*}
This implies that $b^{(2)}=g_{{\rm h}}^{\rm sh}$ satisfies condition (iii) of Proposition \ref{appendix4-proposition1}.
In order to apply Proposition \ref{appendix4-proposition1},
we also need to show that $(b^{(1)}, b^{(2)})=(g_{\rm h}^{\rm w}, g_{\rm h}^{\rm sh})$ satisfies condition (iv).
A direct computation by using Lemma \ref{lemma-FI}, \eqref{inequality-6}, and \eqref{bvp-coeff}--\eqref{definition-htn}
yields that
\begin{equation}
\label{FI-v}
\left|{\rm det} \begin{pmatrix}
D_{\bf p}g^{\rm sh}_{\rm h}(Dv({\bf y}), v({\bf y}), {\bf y})\\[2mm]
D_{\bf p}g^{\rm w}_{\rm h}(Dv({\bf y}), v({\bf y}), {\bf y})
\end{pmatrix}\right| =\bar{\phi}_S^3|\der_{p_1}g^{\rm sh}(D\vphi(\bmxi), \vphi(\bmxi), \bmxi)|
\ge \frac 1M \big(\frac{d_1}{3}\big)^3
\qquad\tx{for}\,\,{\bf y}\in \shock^{({\rm h})}
\end{equation}
for constant $M$ from Lemma \ref{lemma-FI}. We have shown that condition (iv) of Proposition \ref{appendix4-proposition1} holds.

Then we apply Proposition \ref{appendix4-proposition1} to conclude that there exist constants $\alp_1\in(0,1), C>0$, and $l^*>0$ depending only on
$(\iv, \gam,  \sigma_{\rm s}, \sigma_{\rm d})$ such that \eqref{inequality-7} holds.

\smallskip
{\textbf{3.}}
We know from \eqref{bvp-v} that
$v$ satisfies that
$|g_{\rm h}^{\rm sh}(Dv({\bf y}), v({\bf y}), {\bf y})-g_{\rm h}^{\rm sh}(Dv({\bf 0}), v({\bf 0}), {\bf 0})|\equiv 0$ on $\Gam^{({\rm h})}_{\rm shock}$.
This, combined with \eqref{inequality-7}, implies that condition \eqref{app4-prop2-cond3}  stated in Proposition \ref{appendix4-proposition2}
is satisfied with $\alp=\alp_1$.
It follows from \eqref{smooth-mod-coeff-hdgraph} that condition \eqref{app4-prop2-cond1} holds.
Also, \eqref{FI-v} implies that $v$ satisfies condition \eqref{app4-prop2-cond2} with ${\bf y}_0={\bf 0}$.
Moreover, condition \eqref{app4-prop2-cond4} holds for the line segment $\Gam_{\rm shock}^{({\rm h})}$.
Therefore, we obtain from Proposition \ref{appendix4-proposition2} that
\begin{equation}
\label{gradient-holder-shock-v}
|Dv({\bf y})-Dv({\bf 0})|\le C|{\bf y}|^{\alp_1}\qquad\tx{for ${\bf y}\in \Gam^{({\rm h})}_{\rm shock}\cap B_{l^*({\bf 0})}$}
\end{equation}
for a constant $C>0$ depending only on
$(\iv, \gam,  \sigma_{\rm s}, \sigma_{\rm d})$.

Since $\bar{\phi}({\bf 0})=0$ in the $(S,T)$--coordinates,
then $|{\bf y}|\le |\bar{\phi}(S,T)-\bar{\phi}({\bf 0})|+|T|$
for each ${\bf y}=(\bar{\phi}(S,T),T)\in \mcl{D}_{s^*}^{\beta}$.
We apply Lemma \ref{lemma-step3-1} to obtain
\begin{equation}
\label{dist-equiv-1}
|{\bf y}|\le C|(S,T)|=C|{\bmxi}-P_{\beta}|
\end{equation}
for a constant $C>0$ depending only on $(\gam, \iv)$.

By \eqref{definition-htn},
$|{\bmxi}-P_{\beta}|=|(S,T)|\le |v({\bf y})-v({\bf 0})|+|y_2|$ for each $(S,T)\in B_{s^*}(P_{\beta})\cap \Om$.
Then we apply \eqref{lip-bd-v} to obtain
\begin{equation}
\label{dist-equiv-2}
|{\bmxi}-P_{\beta}|=|(S,T)|\le (2K+1)|{\bf y}|
\end{equation}
for constant $K$ from \eqref{lip-bd-v}.

We write \eqref{inequality-7} and \eqref{gradient-holder-shock-v} back in the $\xxi$--coordinates
and apply \eqref{dist-equiv-1}--\eqref{dist-equiv-2} to obtain
\begin{equation}\label{inequality-9}
\begin{split}
&|\vphi_{\xi_2}(\bmxi)-\vphi_{\xi_2}(P_{\beta})|
\le C|\bmxi-P_{\beta}|^{\alp_1}\qquad \tx{in $\ol{\Om\cap B_{s_1}(P_{\beta})}$},\\[1mm]
&|D\vphi(\bmxi)-D\vphi(P_{\beta})|\le C|\bmxi-P_{\beta}|^{\alp_1}\qquad\tx{on $\shock\cap B_{s_1}(P_{\beta})$},
\end{split}
\end{equation}
where $C>0$ and $s_1\in(0,s^*]$ depend only on
$(\iv, \gam,  \sigma_{\rm s}, \sigma_{\rm d})$.

For the rest of proof, each estimate constant is chosen depending only on
$(\iv, \gam,  \sigma_{\rm s}, \sigma_{\rm d})$, unless otherwise specified.
For $\bmxi\in\ol{\Om}$, define $\mathfrak{f}(\bmxi):={\bm\tau}_{\rm w}\cdot(D\bar{\phi}(\bmxi)-D\bar{\phi}(P_{\beta}))$ for the unit
tangent vector ${\bm\tau}_{\rm w}=(1,0)$ to $\Wedge$. Then \eqref{inequality-9} implies that
\begin{equation}\label{inequality-10}
|\mathfrak{f}(\bmxi)-\mathfrak{f}(P_{\beta})|\le \hat{C}|\bmxi-P_{\beta}|^{\alp_1}\qquad\tx{for $\bmxi\in \shock\cap B_{s_1}(P_{\beta})$}.
\end{equation}

Denote $g^{\rm sh}_*({\bf p}):={\bm\tau}_{\rm w}\cdot({\bf p}-D\bar{\phi}(P_{\beta}))$
and regard $g^{\rm sh}_*(D\bar{\phi})=\mathfrak{f}$ as a boundary condition for $\vphi$ on $\shock$.
Since $\Wedge$ is flat in the $\xxi$--coordinates,
we can apply Proposition \ref{appendix4-proposition1} by setting
$(\Gam^1, \Gam^2):=(\shock,\Wedge)$ and $(b^{(1)}, b^{(2)}):=(g_*^{\rm sh}, \hat{g}^{\rm w})$
for $\Gam^j, b^{(j)}, j=1,2$, from Proposition \ref{appendix4-proposition1}.
In particular, condition \eqref{appendix4-prop1-cond4} holds with $\beta=\alp_1$, owing to
\eqref{inequality-10}. Then we obtain constants $\alp\in(0,\alp_1], C>0$, and $s_2\in(0, s_1]$
such that
\begin{equation*}
|g_*^{\rm sh}(D\vphi(\bmxi))-g_*^{\rm sh}(D\vphi(P_{\beta}))|\le C|\bmxi-P_{\beta}|^{\alp}
\qquad\, \tx{for $\bmxi\in \ol{\Om\cap B_{s_2}(P_{\beta})}$}.
\end{equation*}
Combining this with \eqref{inequality-9} and noting that both boundary
conditions $\hat g_{\rm w}$ and $g_*^{\rm sh}$ are linear with constant coefficients
and are linearly independent of each other,
we finally have
\begin{equation}
\label{inequaltiy-11}
|D\vphi(\bmxi)-D\vphi(P_{\beta})|\le C^*|\bmxi-P_{\beta}|^{\alp}\qquad\, \tx{for $\bmxi \in \ol{\Om\cap B_{s_2}(P_{\beta})}$}.
\end{equation}

\smallskip
\textbf{4.} For each $\bmxi\in\shock$, define $d(\bmxi):=|\bmxi-P_{\beta}|$.

\smallskip
{\emph{Claim{\rm :} There exist constants $\omega_0>0$ and $s_3\in(0, s_2]$ such that, for all $\bmxi\in \shock\cap B_{s_3}(P_{\beta})$,}}
\begin{equation*}
{\rm dist}(\bmxi, \Wedge)\ge \omega_0\; d(\bmxi).
\end{equation*}

\smallskip
If this claim holds, then $\Om_{s_3}=\Om\cap B_{s_3}(P_{\beta})$ satisfies
condition (ii) of Proposition \ref{appendix4-proposition3} so that Proposition \ref{lemma-gradient-est}
follows from \eqref{inequaltiy-11} and Proposition \ref{appendix4-proposition3},
where we use \eqref{inequaltiy-11} to satisfy condition \eqref{condition-propc14-grad-est}
stated in Proposition \ref{appendix4-proposition3}.

Now we show the claim.
For a fixed point $P\in\shock$, let $P'$ be the point on $\leftshock$ so that
$PP'\perp \Wedge$. Then
\begin{equation}\label{inequaltiy-12}
{\rm dist}(P, \Wedge)=d(P')\sin\beta-|P'-P|\ge d(P)\sin\beta-|P'-P|.
\end{equation}
Denote $P=(\xi_1^{P}, \xi_2^{P})$ and $P'=(\xi_1^{P'}, \xi_2^{P'})$ in the $\xxi$--coordinates.
Then we see that $P'-P=(0,\xi_2^{P'}-\xi_2^{P})$.
Since $P'\in \leftshock$ and $P\in \shock$, $(\ivphi-\leftvphi)(P')=(\ivphi-\vphi)(P)=0$ so that
\begin{equation*}
\iv|\xi_2^{P'}-\xi_2^{P}|=|(\ivphi-\leftvphi)(P')-(\ivphi-\leftvphi)(P)|=|(\leftvphi-\vphi)(P)|.
\end{equation*}
Since $(\leftvphi-\vphi)(P_{\beta})=0$ by \eqref{1-24ab},
the equation above gives
$$
|P'-P|=\frac{1}{\iv}|(\leftvphi-\vphi)(P)-(\leftvphi-\vphi)(P_{\beta})|.
$$
Then we apply \eqref{inequaltiy-11} to obtain
\begin{equation*}
|P'-P|=\frac{1}{\iv}|(\leftvphi-\vphi)(P)|\le Cd(P)^{1+\alp}\,\qquad\tx{for $P\in\ol{\Om\cap B_{s_2}(P_{\beta})}$}
\end{equation*}
for some constant $C>0$.
Combining this estimate with \eqref{inequaltiy-12},
we can choose constants $\omega_0>0$ and $s_3\in(0, s_2]$ so that the claim holds.

Then Proposition \ref{appendix4-proposition3},
combined with  \eqref{unif-ellip-subs-awa-sn-new} and the results from Steps 3--4,
leads to Proposition \ref{lemma-gradient-est}.
\end{proof}
\end{proposition}

\section{Compactness of the Set of Admissible Solutions}
\label{subsec-compact-admis-sols}
Fix $\gam\ge 1$, $\iv>0$, and $\bar{\beta}\in(0,\betadet)$.
According to all the {\it a priori} estimates obtained in Lemma \ref{lemma-unif-est2},
Corollary \ref{corollary-unif-est-away-sn},
and Propositions \ref{lemma-est-sonic-general-N}, \ref{lemma-est-sonic-general}, \ref{proposition-sub9},
and \ref{lemma-gradient-est}, there exists $\bar{\alp}\in(0,1)$ depending only on $(\iv, \gam, \bar{\beta})$
such that the set:
\begin{equation*}
\left\{\|\vphi\|_{C^{1,\bar{\alp}}(\ol{\Om})}+\|\shock\|_{C^{1,\bar{\alp}}}\,:\,
\begin{array}{ll}
\tx{$\vphi$ is an admissible solution corresponding}\\
\tx{to $(\iv, \beta)\in\mathfrak{R}_{\rm weak}\cap \{0\le \beta \le \bar{\beta}\}$}
\end{array}
\right\}
\end{equation*}
is bounded. For each admissible solution, its pseudo-subsonic region $\Om$ is a bounded domain enclosed
by $\leftsonic$, $\rightsonic$, $\shock$, and $\Wedge$.
These four curves intersect only at $P_j$ for $j=1,2,3,4$.
According to Definition \ref{definition-domains-np}, $\rightsonic$, $\Onormal$, $\righttop$, and $\rightbottom$
are fixed so as to be the same for all admissible solutions.
Moreover, $\leftsonic$, $\Oo$, $\lefttop$, and $\leftbottom$ depend continuously on $\beta\in[0, \betadet]$.
From this observation, the following lemma is obtained:

\begin{lemma}\label{102}
Fix $\gam\ge1$, $\iv>0$, and $\bar{\beta}\in(0,\betadet)$. For each $\beta\in[0, \bar{\beta}]$,
let $\Lbeta$ be defined by Definition {\rm \ref{definition-domains-np}}.
Let $\{\vphi^{(j)}\}$ be a sequence of admissible solutions corresponding
to $(\iv,\beta)\in\mathfrak{R}_{\rm weak}\cap\{0\le \beta\le \bar{\beta}\}$,
and let
$
\lim_{j\rightarrow\infty} \beta_j=\beta_{\infty}
$
for some $\beta_{\infty}\in[0,\bar{\beta}]$.
For each $j$, let $\Om^{(j)}$ and $\shock^{(j)}$ be the pseudo-subsonic region and the curved pseudo-transonic
shock of $\vphi^{(j)}$, respectively. Then there exists a subsequence $\{\vphi^{(j_k)}\}\subset \{\vphi^{(j)}\}$
such that the following properties hold{\rm :}

\smallskip
\begin{itemize}
\item[(a)]  $\{\vphi^{(j_k)}\}$ converges uniformly on any compact subset
of $\ol{\Lambda_{{\beta_{\infty}}}}$ to a function
$\vphi^{(\infty)}\in C^{0,1}_{\rm loc}(\ol{\Lambda_{{\beta_{\infty}}}})$,
and $\vphi^{(\infty)}$ is an admissible solution corresponding to $(\iv, \beta_{\infty})${\rm ;}

\smallskip
\item[(b)] $\Om^{(j_k)}\rightarrow \Om^{(\infty)}$ in the Hausdorff metric{\rm ;}

\smallskip
\item[(c)] If ${\bm\xi}^{(j_k)}\in \ol{\Om^{(j_k)}}$, and
${\bm\xi}^{(j_k)}$ converges to ${\bm\xi}^{(\infty)}\in \ol{\Om^{(\infty)}}$,
then
$$
\vphi^{(j_k)}({\bm\xi}^{(j_k)})\rightarrow {\vphi^{(\infty)}}({\bm\xi}^{(\infty)}),
\qquad
D\vphi^{(j_k)}({\bm\xi}^{(j_k)})\rightarrow D{\vphi}^{(\infty)}({\bm\xi}^{(\infty)}),
$$
where, in the case of $\bm\xi^{(j_k)}\in \shock^{(j_k)}$,
$D\vphi^{(j_k)}(\bm\xi^{(j_k)}):=\lim_{\bm\xi\in \Om^{(j_k)}, \bm\xi\to\bm\xi^{(j_k)}}D\vphi^{(j_k)}(\bm\xi)$,
and $D\vphi^{(\infty)}(\bm\xi)$ for $\bm\xi\in\shock^{(\infty)}$ is defined similarly.
\end{itemize}
\end{lemma}


\begin{proof}
By Corollary \ref{corollary-10-1}, there exists a subsequence $\{\varphi^{(j_k)}\}$ converging uniformly
on any compact subset of $\overline{\Lambda_{\beta_\infty}}$
to a function $\varphi^{(\infty)}\in C^{0,1}_{\rm loc} \overline{\Lambda_{\beta_\infty}})$
that is a weak solution of the boundary value problem consisting of equation \eqref{2-1}
in $\Lambda_{\beta_\infty}$ with boundary condition $\partial_{\bf \nu} \varphi^{(\infty)}=0$
on $\partial \Lambda_{\beta_\infty}$, especially on $\Wedge^{(\infty)}$.
Moreover, $\varphi^{(\infty)}$ satisfies further properties given in Corollary \ref{corollary-10-1}(a)--(e).
In particular, by properties (c) and (e) of Corollary \ref{corollary-10-1},
\begin{equation}\label{shockDoesNotIntersectLim}
\mbox{$\shock^{(\infty)}$ does not intersect relative interiors of
$\Gamma_{\rm sonic}^{\mathcal{O},(\infty)}$, $\Gamma_{\rm sonic}^{\mathcal{N},(\infty)}$,  and $\Wedge^{(\infty)}$.}
\end{equation}
The rest of the
proof is divided into four steps.

\medskip
{\bf 1.} The convergence: $\Omega^{(j_k)} \to \Omega^{(\infty)}$
in the Hausdorff metric follows from
Corollary \ref{corollary-10-1}(a)--(b)
and the continuity  of the parameters of state $(2)$ in $\theta_{\rm w}$.
This implies assertion  (b).

\medskip
{\bf 2.} Next, we prove that $\varphi^{(\infty)}\in C^1(\overline{\Omega^{(\infty)}})$ and assertion (c) hold.  Below we use notation:
\begin{equation}\label{notationConSubs}
\mbox{$\leftsonic=\{P_\beta\}$, \quad  $P_1=P_4=P_\beta$ $\qquad\,\,$ if
$\beta\ge \betasonic$.}
\end{equation}

According to all the {\it a priori} estimates obtained in Lemma \ref{lemma-unif-est2},
Corollary \ref{corollary-unif-est-away-sn},
Propositions \ref{lemma-est-sonic-general-N}, \ref{lemma-est-sonic-general},
\ref{proposition-sub8}, \ref{proposition-sub9},
and \ref{lemma-gradient-est}, there exists $\bar{\alp}\in(0,1)$ depending only on $(\iv, \gam, \bar{\beta})$
such that the set
\begin{equation}\label{unifC1Alp}
\left\{\|\vphi\|_{C^{1,\bar{\alp}}(\ol{\Om})}+\|\shock\|_{C^{1,\bar{\alp}}}\,:\,
\begin{array}{ll}
\tx{$\vphi$ is an admissible solution corresponding}\\
\tx{to $(\iv, \beta)\in\mathfrak{R}_{\rm weak}\cap \{0\le \beta \le \bar{\beta}\}$}
\end{array}
\right\}
\;\;\;\mbox{is bounded,}
\end{equation}
and, for each small $\delta>0$, the set
\begin{equation}\label{unifC3delt}
\begin{split}
&\left\{\|\vphi\|_{C^{4}(\ol{\Om}\setminus\mathcal{N}_\delta(\leftsonic\cup\rightsonic))}+\|\shock\setminus\mathcal{N}_\delta(\{P_1, P_2\})\|_{C^{4}}\,:\,
\begin{array}{ll}
\tx{$\vphi$ is an admissible solution}\\
\tx{corresponding to}\\
\tx{$(\iv, \beta)\in\mathfrak{R}_{\rm weak}\cap \{0\le \beta \le \bar{\beta}\}$}
\end{array}
\right\}\\[1mm]
&\;\;\mbox{is bounded,}
\end{split}
\end{equation}
For each admissible solution, its pseudo-subsonic region $\Om$ is a bounded domain enclosed by $\leftsonic$,
$\rightsonic$, $\shock$, and $\Wedge$. These four curves intersect only at $P_j$ for $j=1,2,3,4$.
According to Definition \ref{definition-domains-np}, $\rightsonic$, $\Onormal$, $\righttop$,
and $\rightbottom$ are fixed so as to be the same for all admissible solutions.
Moreover, $\leftsonic$, $\Oo$, $\lefttop$, and $\leftbottom$ depend continuously on $\beta\in[0, \betadet]$.

Also, using the uniform $C^{1,\bar{\alpha}}$ estimate of the shock functions $f_{\rm{sh}}^{(j)}$ from Proposition \ref{proposition-3} on  interval $[\xi_1^{\lefttop},\,\xi_1^{\righttop}]$ summarized in \eqref{unifC1Alp}, and Corollary \ref{corollary-10-1} (b), we obtain that  $f_{\rm{sh}}^{(j_k)}$ converges to
$f_{\rm{sh}}^{(\infty)}$ in $C^1([\xi_1^{\lefttop^{(\infty)}},\,\xi_1^{\righttop}])$, after rescaling  functions $f_{\rm{sh}}^{(j_k)}$ to be defined on $[\xi_1^{\lefttop^{(\infty)}},\,\xi_1^{\righttop}]$. It follows that
\begin{align}
 & \mbox{$f_{\rm{sh}}^{(\infty)}\in C^{1,\bar{\alpha}}([\xi_1^{\lefttop^{(\infty)}},\,\xi_1^{\righttop}])$;} \label{shockC1AlpLim}\\
& \mbox{If $t_k\in
 [\xi_1^{\lefttop^{(j_k)}},\,\xi_1^{\righttop}]$ and $t_k\to t$, then $t\in
 [\xi_1^{\lefttop^{(\infty)}},\,\xi_1^{\righttop}]$ and
 $(f_{\rm{sh}}^{(j_k)})'(t_k)
 \to (f_{\rm{sh}}^{(\infty)})'(t)$.}
\label{shockC1AlpLim-2}
\end{align}
Let points $\xxi^{(j_k)}$ and  $\xxi^{(\infty)}$ be as in (c).
Then $\xxi^{(\infty)}\in\overline{\Omega^{(\infty)}}$ by assertion (b).

\smallskip
Consider first the case: $\xxi^{(\infty)} \in \Omega^{(\infty)}$. Then, using assertion (b) verified above,
we conclude that  there exists $R>0$ such that
$\overline{B_R(\xxi^{(\infty)})}\subset\Omega^{(\infty)}$ and
$B_R(\xxi^{(j_k)})\subset\Omega^{(j_k)}$ for all sufficiently large $k$.
Then, defining $\Psi^{(j_k)}(\xxi)=\varphi^{(j_k)}(\xxi-\xxi^{(j_k)})$, we have
$$
\|\Psi^{(j_k)}\|_{C^{1,\bar{\alpha}}(\overline{B_R(0)})}\le C.
$$
Using that $\xxi^{(j_k)}\to \xxi^{(\infty)}$, and $\varphi^{(j_k)} \to \varphi^\infty$ uniformly on compact subsets of $\Lambda^\infty$, we see that
$\Psi^{(j_k)} \to \Psi^\infty$ in $C^{1,\frac{\bar{\alpha}}{2}}(\overline{B_{R/2}(\mathbf{0})})$. Then $\Psi^{(j_k)}(\mathbf{0})\to \Psi^{(\infty)}(\mathbf{0})$ and $D\Psi^{(j_k)}(\mathbf{0})\to D\Psi^{(\infty)}(\mathbf{0})$. Thus, we conclude that
$$
 \varphi^{(\infty)}\in C^1(\overline{B_{R/2}(\xxi^{(\infty)})}), \qquad\,\,
(\varphi^{(j_k)}, D\varphi^{(j_k)})(\xxi^{(j_k)})
\to (\varphi^{(\infty)}, D\varphi^{(\infty)})(\xxi^{(\infty)}).
$$

Next, consider the case: $\xxi^{(\infty)} \in \Wedge^{(\infty)}$. Then, by Proposition \ref{proposition-sub3}, there exists $R>0$ such that
$B_{2R}(\xxi^{(j_k)})\cap \partial \Omega^{(\infty)}\subset  \Wedge^{(j_k)}$
and ${ \mbox{\rm dist}}(\xxi^{(j_k)},  \Wedge^{(j_k)}) <\frac{R}{100}$ for all  $k>N$, where $N$ is sufficiently large.
Since $\Wedge^{(j_k)}$ is a straight line, there exists $C>0$ such that
$\varphi^{(j_k)}$ can be extended from $\Omega^{(j_k)}\cap B_{R}(\xxi^{(j_k)})$
to $B_{R}(\xxi^{(j_k)})$ so that the extended function $\varphi^{(j_k)}_{\rm ext}$ satisfies
\begin{equation}\label{uniformEstimatesForExtended-limCor}
\|\varphi^{(j_k)}_{\rm ext}\|_{C^{1,\bar{\alpha}}(\overline{B_{R}(\xxi^{(j_k)})})}\le C
\|\varphi^{(j_k)}\|_{C^{1,\bar{\alpha}}(\overline{ \Omega^{(j_k)}\cap B_{R}(\xxi^{(j_k)})})}\le \hat C
\qquad\mbox{for all $k>N$},
\end{equation}
where $\hat{C}>0$ is a constant independent of $k$.
Selecting a further subsequence (if needed without change of notation), we conclude that $\varphi^{(j_k)}_{\rm ext}$ converges
in $C^{1, \frac{\bar{\alpha}}{2}}$ to $\varphi^{(\infty)}_{\rm ext}$ on any compact subsets of $B_{R}(\xxi^{(\infty)})$.
Also
$\|\varphi^{(\infty)}_{\rm ext}\|_{C^{1,\bar{\alpha}}(\overline{B_{R}(\xxi^{(\infty)})})}\le \hat C$,
by (\ref{uniformEstimatesForExtended-limCor}). Note that, from the uniform convergence $\varphi^{(j_k)}\to \varphi^{(\infty)}$
on compact subsets of $\Lambda_{\beta_\infty}$, it follows that
$\varphi^{(\infty)}_{\rm ext}=\varphi^{(\infty)}$ on
$ \Omega^{(\infty)}\cap B_{R}(\xxi^{(\infty)})$.
Then we can argue as in the previous case to obtain
\begin{equation}\label{convC1Upto}
  \varphi^{(\infty)}\in C^1(\overline{B_{R/2}(\xxi_{\infty})\cap\Omega^{(\infty)}}), \qquad\,\,
(\varphi^{(j_k)}, D\varphi^{(j_k)})(\xxi^{(j_k)})
\to (\varphi^{(\infty)}, D\varphi^{(\infty)})(\xxi^{(\infty)}).
\end{equation}

Cases $\xxi^{(\infty)} \in \Gamma_{\rm sonic}^{\mathcal{N},(\infty)}$ and $\xxi^{(\infty)} \in \Gamma_{\rm sonic}^{\mathcal{O},(\infty)}$ are treated similarly.
In the latter case, we use the fact that each $\Gamma_{\rm sonic}^{\mcl{O}, (j_k)}$ is an arc whose center and radius($=\leftc^{(j_k)}$) depend continuously
on $\beta$. Furthermore, a constant $\hat{C}>0$ is fixed so that $\leftc^{(j_k)}\ge \frac{1}{\hat{C}}$ for all $k$.
Then we may assume with out loss of generality that $R\le \frac{1}{100\hat{C}}$.

Case
$\xxi^{(\infty)} \in \Shock^{(\infty)}$ is considered similarly by employing the bound in
(\ref{unifC1Alp}) for each $\Shock^{(j_k)}$.

\medskip
{\bf 3.} It remains to consider the case that $\xxi^{(\infty)}$ is one of
the corner points $P_m, m=1,2,3,4$,  of $\partial \Om^{(\infty)}$ (see Definition 2.23).

As in the previous case, we need to extend each $\varphi^{(j_k)}$ from $\Omega^{(j_k)}\cap B_{R}(\xxi^{(j_k)})$
to $B_{R}(\xxi^{(j_k)})$ so that the extended functions $\varphi^{(j_k)}_{\rm ext}$ satisfy \eqref{uniformEstimatesForExtended-limCor}
with a uniform constant $C$. Then the rest of the argument follows the previous case to obtain \eqref{convC1Upto}.

The extension satisfying \eqref{uniformEstimatesForExtended-limCor} with uniform constant $C$ for the corner points is
obtained by using the following features of domain $\Omega$ for admissible solutions:
at each corner point $P_m$, for $m=1\dots, 4$,
with the notation convention \eqref{notationConSubs},
two $C^{1,\bar{\alpha}}$ curves (with uniform $C^{1,\bar{\alpha}}$ bounds over all admissible solutions by \eqref{unifC1Alp})
meet at an angle $\theta_m\in (0, \pi)$ which depends only on the parameters of
the uniform states $\leftvphi$ and $\rightvphi$,
where we have used that $\shock$ is tangential to $\leftshock$ (resp. $\rightshock$) at
$P_1$  (resp. $P_2$). Thus, the angles $\theta_m=\theta_m(\beta)$ depend continuously on $\beta$.

From this, if $\xxi^{(\infty)}=P_m^{(\infty)}$ for $m=1,2$,
and $\beta_\infty\ne\betasonic$, we see that, if $\beta_\infty>\betasonic$ (resp. $\beta_\infty<\betasonic$), then $\beta_{j_k}>\betasonic$ (resp. $\beta_{j_k}<\betasonic$)
for all $k\ge N$ with sufficiently large  $N$, so the structure of $\Omega^{(j_k)}$ with $k\ge N$ and of $\Omega_\infty$ is the same
in the sense that both of them are as in either Fig. 2.11 or Fig. 2.12.
From the features of domain $\Omega$ for admissible solutions discussed in the previous paragraph and from \eqref{shockC1AlpLim}, the limiting
domain $\Omega_\infty$ has the same structure as domain $\Omega$ of admissible solution, {\it i.e.}, there exists $R>0$ such that $\partial\Omega\cap B_{4R}(P_m^{(\infty)})$
is the curve consisting of two $C^{1,\bar{\alpha}}$ curve segments meeting at $P_m^{(\infty)}$ at angle $\theta_m(\beta_\infty)\in (0, \pi)$ ({\it i.e.,} the same angle as for admissible solution corresponding to $\beta_\infty$). Then, in an appropriate orthonormal coordinate system $(S,T)$ in $\R^2$ with origin at $P_m^{(\infty)}$, reducing $R$ if necessary, curve
$\partial\Omega^{(\infty)}\cap B_{4R}(P_m^{(\infty)})$ is a graph of the Lipschitz function:
\begin{equation}\label{LipFunctCornerLimLim}
\partial\Omega^{(\infty)}\cap B_{4R}(P_m^{(\infty)})=\{(S,T)\;:\;S=f_\infty(T), \; T\in
(T_1^{(\infty)}, \,T_2^{(\infty)})\},
\end{equation}
where $T_1^{(\infty)}<0<T_2^{(\infty)},
f_\infty(0)=0$, and ${\rm Lip}[f_\infty]=M<\infty$.
The coordinate system $(S, T)$ can be chosen {\it e.g.} as follows: the $S$-axis is along the bisector
of the interior  for $\Omega^{(\infty)}$ angle at $P_m^{(\infty)}$.
Moreover, by \eqref{unifC1Alp} and \eqref{shockC1AlpLim-2}, it follows that,   for all $k\ge N$ (possibly increasing $N$ if needed),
\begin{equation}\label{LipFunctCornerAdmisLim}
\partial\Omega^{(j_k)}\cap B_{2R}(P_m^{({j_k})})=\{(S,T)\;:\;S=f_{j_k}(T), \; T\in (T_1^{({j_k})}, \,T_2^{({j_k})})\},
\end{equation}
where $T_1^{({j_k})}<0<T_2^{({j_k})}$,
$\xxi^{(j_k)}\in B_{ R/{100}}(P_m^{(\infty)})$, $P_m^{(j_k)}\in B_{ R/{100}}(P_m^{(\infty)})$,
and ${\rm Lip}[f_{j_k}]\le2M$. Then we can extend functions $\varphi^{({j_k})}$ from $\Omega\cap B_{R}(\xxi^{({j_k})})$
to $B_{R}(\xxi^{({j_k})})$ so that \eqref{uniformEstimatesForExtended-limCor} holds with $C$ depending only on $M$ and $R$.
For such an extension, we can use  the extension operator introduced \cite[Definition 13.9.3]{CF2}, and then \cite[Lemma 13.9.6]{CF2} to
show the $C^{2,\alpha}$ estimates for the extension operator with constant depending on $M$ and $R$ in the present case; the corresponding
$C^{1,\alpha}$ estimates are obtained similarly (and simpler).

Suppose that  $\xxi^{(\infty)}=P_m^{(\infty)}$ for $m=1,2,3,4$, and $\beta_\infty=\betasonic$.
By passing to a further subsequence (without changing notation),
we can assume that either $\beta_{j_k}<\betasonic$ for all $k\in\mathbb{N}$ or $\beta_{j_k}\ge\betasonic$
for all $k\in\mathbb{N}$. In the later case, we argue as above. It remains consider the case
$\beta_{j_k}<\betasonic$ for all $k\in\mathbb{N}$, {\it i.e.}, when the solutions
of the structure as on Fig. \ref{fig:regularsol1} converge to a solution of the structure as on Fig. \ref{fig:regularsol2}.
For $\xxi^{(\infty)}=P_m^{(\infty)}$
with $m\in\{2,3\}$, the argument is the same as before. Then consider the case: $m\in \{1,4\}$, which means $\xxi^{(\infty)}
=P_{\beta_\infty}$ by \eqref{notationConSubs} since $\beta_\infty=\betasonic$. Choose $R>0$
and a coordinate system $(S, T)$ in which \eqref{LipFunctCornerLimLim} holds. Then, for $k>N$,
curve $\partial\Omega^{(j_k)}\cap B_{2R}(P_m^{({j_k})})$ consists of three smooth parts: $\Wedge^{(j_k)}\cap B_{2R}(P_m^{({j_k})})$, $\Gamma_{\rm sonic}^{\mathcal{O},(j_k)}$, and $\shock^{(\infty)}\cap B_{2R}(P_m^{({j_k})})$ which
meet at points $P_1^{(j_k)}$ and $P_4^{(j_k)}$ respectively, and
${\Gamma_{\rm sonic}^{\mathcal{O},(j_k)}}\to \{P_{\beta_\infty}\}$ as $k\to\infty$ in the Hausdorff metric.
Then it follows from the features of domain $\Omega$ for admissible solutions discussed above that, for any sufficiently large $k$,
curve $\partial\Om^{(j_k)}\cap B_{2R}(P_m^{(j_k)})$ is a Lipschitz graph in the $(S,T)$-coordinates so that \eqref{LipFunctCornerAdmisLim} holds.
Then the case for ${\bm\xi}^{(\infty)}=P_m$ with $m=3,4$, can be handled similarly to the case for ${\bm\xi}^{(\infty)}=P_m$ with $m=1,2$.

\smallskip
Therefore, we conclude that  $\varphi^{(\infty)}\in C^1(\overline{\Omega^{(\infty)}})$ and assertion (c) of Lemma \ref{102} hold.

\medskip
{\bf 4. }  It remains to prove assertion (a).
We first prove that $\varphi^{(\infty)}$ satisfies Definition \ref{def-regular-sol} (i) (Cases I and II).

By estimate \eqref{new-dec-2} of Proposition \ref{proposition-distance}
for each
$\shock^{({j_k})}$, sending to the limit
as $k\to\infty$ by using Corollary \ref{corollary-10-1}(b), we conclude that (i-1) holds for $\shock^{({\infty})}$.
From Corollary \ref{corollary-10-1}(b) combined with the estimates of the shock functions $\fshpolar^{(j_k)}$ in Lemma \ref{lemma-unif-est2}, sending to the limit
as $k\to\infty$,
we conclude that $\shock^{(\infty)}$ is $C^\infty$ in its relative interior, so that (i-2) holds
for $\shock^{({\infty})}$.
Property (i-3) for the limiting solution $\vphi^{(\infty)}$ is obtained from property (i3) applied to each $\varphi^{(j_k)}$,
by sending $k\to\infty$ and using \eqref{shockC1AlpLim}--\eqref{shockC1AlpLim-2} and the continuous dependence of $\leftvphi$ on $\beta$.
Finally, \eqref{shockDoesNotIntersectLim} implies (i-4).
This concludes that $\varphi^{(\infty)}$ satisfies (i) of Definition \ref{def-regular-sol} (Cases I and II).

\smallskip
Next, we show that $\varphi^{(\infty)}$ satisfies Definition \ref{def-regular-sol} (ii) (Cases I and II).
In Steps 2--4 above, it is shown that $\varphi^{(\infty)}\in C^1(\overline{\Omega^{(\infty)}})$.

We now prove that $\varphi^{(\infty)}\in C^{3}(\ol{\Om^{(\infty)}}\setminus(\leftsonic\cup\rightsonic))$.
For a constant $\delta>0$, let $K\subset
{\Om^{(\infty)}}\setminus\mathcal{N}_\delta(\leftsonic\cup\rightsonic)$ be compact.
Then, for a sufficiently large $N$,
$K\subset {\Om^{(j_k)}}\setminus\mathcal{N}_{\delta/2}(\leftsonic\cup\rightsonic)$ for all $k\ge N$.
Since $\varphi^{(j_k)}\to \varphi^{(\infty)}$ uniformly on $K$ and \eqref{unifC3delt}
holds for each $\varphi^{(j_k)}$, we obtain
$$
\|\varphi^{(\infty)}\|_{C^3(K)}\le C_1(\delta),
$$
where the estimate constant  $C_1(\delta)$ depends on $\delta$, but is independent of $K \subset
{\Om^{(\infty)}}\setminus\mathcal{N}_\delta(\leftsonic\cup\rightsonic)$. This implies that
$\varphi^{(\infty)}\in C^{3}(\ol{\Om^{(\infty)}}\setminus\mathcal{N}_\delta(\leftsonic\cup\rightsonic))$.
Since $\delta>0$ is arbitrary, we obtain that $\varphi^{(\infty)}\in C^{3}(\ol{\Om^{(\infty)}}\setminus(\leftsonic\cup\rightsonic))$.
Also, by Corollary \ref{corollary-10-1}(d),
$\varphi^{(\infty)}$ satisfies \eqref{1-24} in Case I and \eqref{2-4-a0} in Case II of Definition \ref{def-regular-sol}.

Then, in order to complete the proof of (ii-1)--(ii-3), it remains to show that

\begin{itemize}
\item $\vphi^{(\infty)}$ is $C^1$ across $\rightsonic$ in Cases I and II of Definition \ref{def-regular-sol};

\item $\vphi^{(\infty)}$ is $C^1$ across $\leftsonic$ in Case I of Definition \ref{def-regular-sol};

\item conditions at $P_\beta$ in
 (ii-3) of Definition \ref{def-regular-sol} (Case II) hold for $\vphi^{(\infty)}$.
\end{itemize}
Indeed, the first two statements imply (ii-1) in Cases I and II, while the last statement yields  (ii-3) in Case 2.

Estimate \eqref{par-est1-N} in Proposition \ref{lemma-est-sonic-general-N} holds for each $\varphi^{(j_k)}$, which implies
that $D^m\vphi^{({j_k})}=D^m\rightvphi$  on $\rightsonic$ for $m=0,1$. In the
limit: $k\to\infty$, using Lemma \ref{102}(c) (proved above),  we obtain that
$D^m\vphi^{(\infty)}=D^m\rightvphi$  on $\rightsonic$ for $m=0,1$. That is, $\vphi^{(\infty)}$ is $C^1$ across $\rightsonic$.

If $\beta_\infty<\betac^{(\iv)}$, then $\beta_{j_k}<\betac^{(\iv)}$ for all $k\ge N$ with sufficiently large $N$. Estimate \eqref{par-est1}
in Propositions \ref{lemma-est-sonic-general} and \ref{proposition-sub8} holds for each $\varphi^{(j_k)}$, which implies that
$D^m\vphi^{({j_k})}=D^m\leftvphi^{({j_k})}$ on $\Gamma_{\rm sonic}^{\mathcal{O},({j_k})}$ for $m=0,1$. In the limit: $k\to\infty$,
using the continuous dependence of parameters of the uniform state $\leftvphi$  on $\beta$,
the continuous dependence of  $\leftsonic$ on $\beta$ in the Hausdorff metric, and Lemma \ref{102}(c) (proved above),
we obtain that $D^m\vphi^{(\infty)}=D^m\leftvphi^{(\infty)}$ on $\Gamma_{\rm sonic}^{\mathcal{O},(\infty)}$ for $m=0,1$,
if $\beta_\infty<\betac^{(\iv)}$. That is,
$\vphi^{(\infty)}$ is $C^1$ across $\Gamma_{\rm sonic}^{\mathcal{O},(\infty)}$.

If $\beta_\infty\ge\betac^{(\iv)}$, we may have both cases $\beta_{j_k}<\betac^{(\iv)}$ and
$\beta_{j_k}\ge\betac^{(\iv)}$.
Then we use
estimate \eqref{par-est1} in Propositions \ref{lemma-est-sonic-general} and \ref{proposition-sub8}
and the results in Proposition \ref{lemma-gradient-est} to obtain
$D^m\vphi^{({j_k})}(P_1^{(j_k)})=D^m\leftvphi^{({j_k})})(P_1^{(j_k)})$ for $m=0,1$,
where we use the notation convention \eqref{notationConSubs}.
In the limit: $k\to\infty$, using the continuous dependence of parameters of the uniform state $\leftvphi$  on $\beta$,
the continuous dependence of $\overline{\leftsonic}$ in the Hausdorff metric (again, using notations \eqref{notationConSubs}),
and Lemma \ref{102} (c) (proved above),
we obtain that
$P_1^{(j_k)}\to P_1^{(\infty)}=P_{\beta_\infty}$ and
$D^m\vphi^{(\infty)}(P_{\beta_\infty})=D^m\leftvphi^{(\infty)}(P_{\beta_\infty})$ for $m=0,1$.
That is, conditions at $P_\beta$ in
Definition \ref{def-regular-sol}(ii-3) (Case II) hold for $\vphi^{(\infty)}$.

Now (ii-1)--(ii-3) in Cases I--II are proved.

Property (ii-4) follows from the fact that $\varphi^{(\infty)}$
is a weak solution of the boundary-value problem consisting
of equation \eqref{2-1} in $\Lambda_{\beta_{\infty}}$ with boundary condition $\partial_{\bf \nu} \varphi^{(\infty)}=0$
on $\partial \Lambda_{\beta_{(\infty)}}$, especially on $\Wedge^{(\infty)}$,
in the sense of Remark {\rm \ref{remark-admsol-wksol}},
and from the regularity of $\varphi^{(\infty)}$ in (ii-1)--(ii-3).

This completes the proof that $\varphi^{(\infty)}$ satisfies Definition \ref{def-regular-sol} (ii) (Cases I--II).

Properties (iii)--(v) of  Definition \ref{def-regular-sol} (Cases I--II)
for $\varphi^{(\infty)}$ directly follow
from the corresponding properties for $\varphi^{(j_k)}$,
Corollary \ref{corollary-10-1} (b)-(c), and the continuous dependence of the parameters of
$\leftvphi$ on $\beta$.

This completes the proof of Lemma \ref{102}(a), so does the proof of Lemma \ref{102}.
\end{proof}



\chapter{Iteration Set}
\label{section-itr-set}

\numberwithin{equation}{section}

In order to prove the existence of admissible solutions in the sense of Definition \ref{def-regular-sol}
for all $(\iv, \beta)\in\mathfrak{R}_{\rm weak}$ by employing the Leray-Schauder degree for a fixed point,
we first introduce the iteration set.

\section{Mapping the Admissible Solutions to the Functions Defined in $\iter$}\label{subsec-Q}
Fix $\gam \ge 1$ and $\iv>0$. We continue to follow Definition \ref{definition-domains-np}
for the notations: $\Oi$, $\Oo$, $\Onormal$, $\rightsonic$, $\leftsonic$, and $P_j$ for $j=1,2,3,4$, {\it etc.}.
Denote $\iter = (-1,1)\times (0,1)$.

\begin{definition}\label{definition-Qbeta}
Let $(\ivphi, \rightvphi, \leftvphi)$ be defined by \eqref{def-uniform-ptnl-new}.

{\rm (i)} \emph{Definition of $\leftch$}. For each $\beta\in [0,\betadet]$, define $\leftch$ by
\begin{equation*}
\leftch:=\rm{dist}(\leftsonic, \Oo)=
\begin{cases}
\leftc\qquad&\,\,\tx{for $\beta<\betasonic$},\\
|\Oo P_{\beta}|\qquad&\,\,\tx{for $\beta\ge \betasonic$}.
\end{cases}
\end{equation*}
Note that $\leftch<\leftc$ if $\beta>\betasonic$.
\smallskip

{\rm (ii)} \emph{Extended sonic arcs}.
Since $\leftch$ depends continuously on $\beta\in [0, \frac{\pi}{2})$,
a constant $\delta_0>0$ can be chosen depending only on $(\iv, \gam)$ such that
$\rightshock^{\delta_0}=\{\bmxi\in\R^2\,:\,(\ivphi-\rightvphi)(\bmxi)=-\delta_0\}$
and $\der B_{\rightc}(\Onormal)$ intersect at two distinct points,
and $\leftshock^{\delta_0}=\{\bmxi\in \R^2\,:\,(\ivphi-\leftvphi)(\bmxi)=-\delta_0\}$
and $\der B_{\leftch}(\Oo)$ intersect at two distinct points for each $\beta\in[0,\betadet]$.
Let $\leftsonicdelta$ be the smaller arc lying on $\der B_{\leftch}(\Oo)$
with endpoints $\leftbottom$ and $\lefttop'$,
where $\lefttop'$ is the intersection point of $\leftshock^{\delta_0}$ and $\der B_{\leftch}(\Oo)$
closer to $\lefttop$.
Similarly, let $\rightsonicdelta$ be the smaller arc lying on $\der B_{\rightc}(\Onormal)$
between $\rightshock^{\delta_0}$ and $\etan=0$ with endpoints $\righttop'$ and $\rightbottom$,
where $\righttop'$ is the intersection point of $\rightshock^{\delta_0}$
and $\der B_{\rightc}(\Onormal)$ closer to $\righttop$.
\smallskip

{\rm (iii)} {\emph{Definition of $\Qbeta$}}.
Define $\Qbeta$ as the bounded region enclosed by $\leftsonicdelta$, $\rightsonicdelta$,
$\leftshock^{\delta_0}$, $\rightshock^{\delta_0}$, and $\Wedge$.
\end{definition}

For each $\beta\in[0,\betadet]$, we first define a map $G_1:\Qbeta\to \R^2$ such that
\begin{equation}
\label{prop-G1}
G_1(\bmxi)=
\begin{cases}
(x+\leftu-\leftc,y)\qquad\,&\tx{for $\bmxi$ near $\leftsonicdelta$},\\[2mm]
(\rightc-x,y)\qquad\,&\tx{for $\bmxi$ near $\rightsonicdelta$},
\end{cases}
\end{equation}
for the $(x,y)$--coordinates defined by \eqref{coord-o} near $\leftsonicdelta$ and by \eqref{coord-n} near $\rightsonicdelta$.
We take several steps to construct $G_1$.
The definition of $G_1$ is given in \eqref{12-16}.
First, we define a map $F_1: \Qbeta \rightarrow \R^2$ such that $F_1(\bmxi)\cdot (1,0)=x+\leftu-\leftc$
for $\bmxi\,\,\tx{near}\,\,\leftsonicdelta$ and $F_1(\bmxi)\cdot (1,0)=\rightc-x$ for $\bmxi\,\,\tx{near}\,\,\rightsonicdelta$.
Then we define a map $F_2: F_1(\Qbeta)\rightarrow \R^2$ so that $(F_2\circ F_1)(\bmxi)\cdot (1,0)=F_1(\bmxi)$,
and $(F_2\circ F_1)(\bmxi)\cdot (0,1)=y$ for $\bmxi\,\,\tx{near}\,\,\leftsonicdelta\cup\rightsonicdelta$.
Finally, $G_1$ is defined by $G_1=F_2\circ F_1$ as in \eqref{12-16}.

For $\eps>0$, define two sets $\oD_{\eps}$ and $\nD_{\eps}$ by
\begin{equation}\label{definition-Dr-ext}
\begin{split}
&\oD_{\eps}:=(\Qbeta\cap\{\xin<\leftu\})\setminus \ol{B_{\leftch-\eps}(\Oo)},\\[1mm]
&\nD_{\eps}:=(\Qbeta\cap\{\xin>0\})\setminus \ol{B_{\rightc-\eps}(\Onormal)}.
\end{split}
\end{equation}
Since $\leftch, \leftshock$, and $\Oo$ depend continuously on $\beta\in[0, \frac{\pi}{2})$,
there exist constants $k>4$ and $\delta_1\in(0, \frac{\pi}{2})$
depending only on $(\iv, \gam)$ such that, for each  $\beta\in[0, \betadet]$, we have
\begin{equation}
\label{7-a1}
\begin{split}
&\oD_{\frac{4}{k}\leftch}\subset \{x_{\lefttop}<x<x_{\lefttop}+\frac{4}{k}\leftch,\,\,
 \beta<y+\beta<\frac{\pi}{2}-\delta_1\},\\
&\nD_{\frac{4}{k}\rightc}\subset \{0<x<\frac 4k \rightc,\,\,
 0<y<\frac{\pi}{2}-\delta_1\}.
\end{split}
\end{equation}

Define cut-off functions $\zetao,\zetan, \leftchi$, and $\rightchi$ as follows:

\smallskip
\begin{itemize}
\item[(i)] $\zetao, \zetan\in C^4(\R)$ satisfy
 \begin{align}
    \label{7-b4}
   &\zetao(r)=\begin{cases}
   1&\tx{for $r\ge \leftch(1-\frac{2}{k})$},\\[1mm]
    0&\tx{for $r<\leftch(1-\frac{3}{k})$},
    \end{cases}
    \qquad\;\;
    0\le \zetao'(r)\le \frac{2k}{\leftch} \;\;\;\tx{on $\R$};\\[2mm]
 \label{7-b5}
&\zetan(r)=
\begin{cases}
1&\tx{for}\;r\ge \rightc(1-\frac{2}{k}),\\[1mm]
0&\tx{for}\;r<\rightc(1-\frac 3k),
\end{cases}
\qquad\,\, 0\le \zetan'(r)\le \frac{2k}{\rightc}\;\;\;\tx{on $\R$};
\end{align}

\item[(ii)] Let $q^{\delta_0}_{\mcl{O}}$ be the distance between $\Oo=(\leftu,0)$ and $\leftshock^{\delta_0}$, and denote
    \begin{equation}
    \label{definition-leftu-delta-new}
    \leftu^{\delta_0}:=\leftu-q^{\delta_0}_{\mcl{O}}\sin\beta.
    \end{equation}
Since $\leftu=-\iv\tan\beta<0$, $\leftu^{\delta_0}<0$.
Then $\leftchi,\rightchi\in C^4(\R)$ satisfy
\begin{align}
\label{10-e9-a}
&\leftchi(\xin)=
\begin{cases}
1&\quad \text{for $\xin\le \leftu^{\delta_0}-\frac{2\leftch}{k}$},\\[1mm]
0&\quad \text{for $\xin\ge \leftu^{\delta_0}$},
\end{cases}\qquad -\frac{2k}{\leftch}\le \leftchi'(\xin)\le 0\;\;\;\;\mbox{on $\R$};\\[2mm]
\label{10-e9}
&\rightchi(\xin)=
\begin{cases}
0&\quad  \text{for $\xin\le \frac{\rightc}{k}$},\\[1mm]
1&\quad  \text{for $\xin\ge \frac{2\rightc}{k}$},
\end{cases}
\qquad\qquad\,\,\,\,  0\le \rightchi'(\xin)\le \frac{2k}{\rightc}\;\;\;\;\;\;\mbox{on $\R$}.
\end{align}
\end{itemize}
Choose constant $k>4$ sufficiently large, depending only on $(\iv, \gam)$, such that
\begin{equation}
\label{10-e6}
\oD_{\frac{3\leftch}{k}}\cap \{\xin<\leftu^{\delta_0}\}
\subset \{\xin<\leftu^{\delta_0}-\frac{3\leftch}{k}\},\quad
\nD_{\frac{3\rightc}{k}}\subset \{\xin>\frac{3\rightc}{k}\}.
\end{equation}

Next, define a variable $r$ by
\begin{equation}
\label{def-r}
r=\begin{cases}
\sqrt{(\xin-\leftu)^2 +\etan^2}\qquad\, &\mbox{for $\xin\le \leftu^{\delta_0}$},\\[1mm]
\sqrt{\xin^2+\etan^2}&\mbox{for $\xin\ge 0$}.
\end{cases}
\end{equation}
Since $\leftu^{\delta_0}<0$, $r$ is well defined by \eqref{def-r}.

For the cut-off functions $(\zetao, \zetan, \leftchi, \rightchi)$ given by \eqref{7-b4}--\eqref{10-e9}
under the choice of $k$ to satisfy \eqref{10-e6}, we define a function $h_1: \Qbeta\to \R$ as
\begin{equation}
\label{10-e8}
\begin{split}
h_1(\xin,\etan):=&\;\bigl((\leftu-r)\zetao(r)+(1-\zetao(r))\xin\bigr)\leftchi\\
&+\bigl(1-\leftchi \bigr)
\bigl(\xin(1-\rightchi)+(r \zetan(r)+(1-\zetan(r))\xin)\rightchi\bigr).
\end{split}
\end{equation}
In \eqref{10-e8}, $\leftchi$ and $\rightchi$ are evaluated at $\xin$.

Define a map $F_1: \Qbeta\to \R^2$ by
\begin{equation}
\label{10-e7}
F_1(\xin,\etan):=\bigl(h_1(\xin,\etan),\etan\bigr).
\end{equation}

\begin{lemma}\label{lemma-10-1}
There exist constants $C>0$ and $\delta_{F_1}>0$ depending only on $(\iv, \gam)$ such that,
for each $\beta\in[0,\betadet]$, $F_1$ defined by \eqref{10-e7} satisfies the following properties{\rm :}

\smallskip
\begin{itemize}
\item[(a)] $\|F_1\|_{C^4(\ol{\Qbeta})}+\|F_1^{-1}\|_{C^4(\ol{F_1(\Qbeta)})}\le C $, and $\det(DF_1)\ge \delta_{F_1}$ in $\Qbeta${\rm ;}

\smallskip
\item[(b)] Denoting $F_1(\xxi):=(s,t)$, then
    \begin{equation}
    \label{7-a5}
    F_1(\Wedge)=\{(s,0)\,:\,s\in(\leftu-\leftch,\rightc)\};
    \end{equation}

\item[(c)] For $\iphi:=\ivphi+\frac 12|\xxi|^2$,
\begin{equation*}
\der_t\iphi\bigl(F_1^{-1}(s,t)\bigr)= -\iv\qquad \text{for all $(s,t)\in \ol{F_1(\Qbeta)}$};
\end{equation*}

\item[(d)] For each $j=1,\cdots,4$, denote $P_j=(\xi_1^{P_j},\xi_2^{P_j})$ in the $\xxi$--coordinates. Then
\begin{equation*}
\begin{split}
&F_1(\lefttop)=(\leftu-\leftch,\xi_2^{\lefttop}), \quad F_1(\righttop)=(\rightc,\xi_2^{\righttop}), \\
&F_1(\rightbottom)=(\rightc, 0), \qquad\quad\quad\,  F_1(\leftbottom)=(\leftu-\leftch, 0);
\end{split}
\end{equation*}

\item[(e)] For $h_1$ defined by \eqref{10-e8},
\begin{equation*}
h_1(\bmxi)=
\begin{cases}
\leftu-\leftc+x&\quad \text{if ${\rm dist}(\bmxi, \leftsonicdelta)<\frac{\leftch}{k}$},\\[1mm]
\rightc-x&\quad \text{if ${\rm dist}(\bmxi,\rightsonicdelta)<\frac{\rightc}{k}$}
\end{cases}
\end{equation*}
for the $(x,y)$--coordinates defined by \eqref{coord-n} and \eqref{coord-o}.
\end{itemize}

\begin{proof}
By the definition of $F_1$ in \eqref{10-e7}, we have
\begin{equation}
\label{12-9}
\det(DF_1)=\der_{\xin}h_1.
\end{equation}
Choose constant $k$ large to satisfy that $\rightchi\leftchi'=0$ and $\zetan\rightchi'=\zetao\leftchi'=0$.
Then, from definition \eqref{10-e8} of $h_1$ and \eqref{7-b4}--\eqref{10-e9},
\begin{equation}
\label{h1-deriv-xi1}
\der_{\xin}h_{1}(\xxi)=\sum_{j=1}^3 a_j,
\end{equation}
where
\begin{equation*}
\begin{split}
&a_1=\Big(\frac{\leftu-\xin}{r}\zetao+(1-\zetao)+\frac{\leftu-\xin}{r}(r-(\leftu-\xin))\zetao' \Big)\leftchi,\\
&a_2=\Big( \frac{\xin}{r}\zetan+(1-\zetan)+\frac{\xin}{r}(r-\xin)\zetan'\Big)\rightchi(1-\leftchi),\\[1mm]
&a_3=(1-\rightchi)(1-\leftchi).
\end{split}
\end{equation*}
Then \eqref{7-a1} implies that
\begin{equation}
\label{12-7}
\begin{split}
\der_{\xin}h_{1}
&\ge \Bigl(\frac{\leftu-\xin}{r}\zetao+(1-\zetao) \Bigr)\leftchi
+\Bigl( \bigl(\frac{\xin}{r}\zetan+(1-\zetan)\bigr)\rightchi+(1-\rightchi)\Bigr)(1-\leftchi)\\
&\ge \cos(\frac{\pi}{2}-\delta_1)
\end{split}
\end{equation}
for $\delta_1$ from \eqref{7-a1}.

Moreover, it follows from \eqref{h1-deriv-xi1} that
\begin{equation}
\label{upper-bd-hi-deriv-xi1}
\sup_{\bmxi\in\Qbeta}\der_{\xin}h_1(\bm\xi)\le C
\end{equation}
for a constant $C>0$ depending only on $(\gam, \iv)$.

For a constant $a$, if $Q_{\beta}\cap \{\xi_2=a\}$ is nonempty,
then \eqref{12-7} implies that the one-dimensional
map $(\xi_1,a)\in Q_{\beta}\cap \{\xi_2=a\}\mapsto h_1(\xi_1,a) $ is invertible.
Then it follows directly from the definition of $F_1$ given in \eqref{10-e7} that $F_1$ is invertible.
Also, we can directly check that $F_1$ and $F_1^{-1}$ are  $C^4$ from \eqref{10-e8}, which yields (a).
Finally, (b), (d), and (e)
follow from  \eqref{10-e8}--\eqref{10-e7}.

By \eqref{2-4-b6} and \eqref{10-e7},
$\iphi\bigl(F_1^{-1}(s,t)\bigr)=-\iv t$, which gives
\begin{equation*}
\der_t\iphi\bigl(F_1^{-1}(s,t)\bigr)=-\iv  \qquad\mbox{for all $(s,t)\in F_1(Q)$}.
\end{equation*}
This proves (c).
\end{proof}
\end{lemma}

By the definition of $h_1$ in \eqref{10-e8}, we have
 \begin{equation*}
 \ol{F_1(\Qbeta)}\subset [\leftu-\leftch,\rightc]\times [0,\infty).
 \end{equation*}

\begin{lemma}\label{lemma-7-1}
Fix $\gam \ge 1$, $\iv>0$, and $\bar{\beta}\in(0, \betadet)$.
Then there exists a constant $m_0>0$ depending only on $(\iv, \gam, \bar{\beta})$
such that
any admissible solution $\vphi$ corresponding to $(\iv, \beta)\in\mathfrak{R}_{\rm weak}\cap\{0\le \beta\le \bar{\beta}\}$
satisfies
\begin{equation}\label{7-a4}
\der_t(\ivphi-\vphi)(F_1^{-1}(s,t))\le -m_0<0\,\qquad\text{in $\ol{F_1({\Om})}$}.
\end{equation}
Therefore, there exists a unique function $\tgshock: [\leftu-\leftch,\rightc]\rightarrow \R_+$ such that
\begin{equation*}
F_1(\shock)=\{(s,\tgshock(s))\,:\, \leftu-\leftch< s <\rightc\}.
\end{equation*}

\begin{proof}
For each $\beta\in[0,\betadet]$, we represent $F_1^{-1}$ as
$$
F_1^{-1}(s,t)=(\til h_1(s,t),t) \qquad \mbox{in $\ol{F_1(\Qbeta)}$}.
$$
This expression yields that
\begin{equation}\label{7-b6}
\der_t(\ivphi-\vphi)(F_1^{-1}(s,t))=D(\ivphi-\vphi)|_{F_1^{-1}(s,t)}\cdot(\der_t\til h_1(s,t),1).
\end{equation}
It follows from $(F_1\circ F_1^{-1})(s,t)=\bigl(h_1(\til h_1(s,t)),t\bigr)=(s,t)$
that $\der_t\til h_1(s,t)=-\frac{\partial_{\etan}h_1}{\partial_{\xin}h_1}$. This implies that
\begin{equation*}
(\der_t\til h_1(s,t),1)=-\frac{1}{\partial_{\xin}h_1}(\partial_{\etan}h_1, -\partial_{\xin}h_1),
\end{equation*}
where $D_{(\xin,\etan)}h_1$ is evaluated at $\xxi=F_1^{-1}(s,t)$.

Next, we compute  $\bm v:=\frac{1}{\partial_{\xin}h_1}(-\der_{\etan}h_1, \der_{\xin}h_1)$.

\smallskip
{\it Case} 1. If $\leftchi\neq 0$ so that $\rightchi=\rightchi'=0$, we use $\zetao(r)\leftchi'(\xi_1)\equiv 0$ to obtain
\begin{equation}
\label{7-a2}
\der_{\xin }h_1\bm v =k_1\bm a_1+k_2\bm a_2,
\end{equation}
where
\begin{equation*}
\bm a_1=(\sin y, \cos y),\quad \bm a_2=(0,1),\quad k_1=\big(\zetao+r(1-\cos y)\zetao'\big)\leftchi,
\quad k_2=1-\zetao\leftchi
\end{equation*}
for the $(x,y)$--coordinates defined by \eqref{coord-o}.

\smallskip
{\it Case} 2. If $\leftchi=0$ so that $\leftchi=\leftchi'=0$, we use $\zetan(r)\rightchi'(\xi_1)\equiv 0$ to obtain
\begin{equation}
\label{7-a3}
 \der_{\xin}h_1\bm v=l_1\bm b_1+l_2\bm b_2,
\end{equation}
where
\begin{equation*}
\bm b_1=(-\sin y, \cos y),\quad \bm b_2=(0,1),\quad l_1=\big(\zetan+r(1-\cos y)\zetan'\big)\rightchi,
\quad l_2=1-\zetan\rightchi
\end{equation*}
for the $(x,y)$--coordinates defined by \eqref{coord-n}.

\medskip
\noindent
{\emph{Claim{\rm :} There exists a constant $\til m>0$ depending only on $(\iv, \gam,  \bar{\beta})$
such that any admissible solution $\vphi$ corresponding to $(\iv, \beta)\in\mathfrak{R}_{\rm weak}\cap\{0\le \beta\le \bar{\beta}\}$
satisfies}}
\begin{equation*}
\sup_{P\in \ol{\Om}}\big(D(\ivphi-\vphi)\cdot \bm v\big)(P)\le -\til m.
\end{equation*}

\smallskip
Fix an admissible solution $\vphi$ for $\beta\in[0, \bar{\beta}]$.
Let the unit vectors ${\bm a_1}$, ${\bm a_2}$, ${\bm b_1}$, and ${\bm b_2}$ be from \eqref{7-a2}--\eqref{7-a3}.
Then ${\bm a_1}, \bm a_2\in \cone$ for all $y\in[0,\frac{\pi}{2}-\beta-\delta_1]$ for $\delta_1>0$ from \eqref{7-a1},
and ${\bm b_1}, \bm b_2\in \cone$ for all $y\in[0, \frac{\pi}{2}-\delta_1]$.
Moreover, $k_j$ and $l_j$, $j=1,2$, are nonnegative and satisfy that $k_1+k_2\ge1$ and $l_1+l_2\ge 1$ for all $P\in \ol{\Om}$.
Then \eqref{cor-dir-mont} yields
\begin{equation*}
\sup_{P\in \ol{\Om}}\big(\der_{\xin} h_1\, D(\ivphi-\vphi)\cdot \bm v  \big)(P)\le -m_{\vphi}<0
\end{equation*}
for a constant $m_{\vphi}>0$.
Furthermore, Lemma \ref{102} implies that there exists a constant $m_1>0$ depending only
on $(\iv, \gam,  \bar{\beta})$ such that any admissible solution $\vphi$ corresponding
to $(\iv, \beta)\in \mathfrak{R}_{\rm weak}\cap\{0\le \beta\le \bar{\beta}\}$ satisfies
\begin{equation}
\label{2014-0821}
\sup_{P\in \ol{\Om}} \big(\der_{\xin} h_1\, D(\ivphi-\vphi)\cdot \bm v\big)(P)\le -m_1.
\end{equation}
Combining \eqref{2014-0821} with \eqref{12-7}--\eqref{upper-bd-hi-deriv-xi1},
we conclude that there exists a constant $m_0>0$
depending only on $(\iv, \gam,  \bar{\beta})$ such that any admissible solution $\vphi$
for $\beta\in[0,\bar{\beta}]$ satisfies
\begin{equation}
\label{est-t}
\der_t(\ivphi-\vphi)(F_1^{-1}(s,t))=(D(\ivphi-\vphi)\cdot \bm v)(F_1^{-1}(s,t))\le -m_0
<0
\end{equation}
for all $(s,t)\in F_1(\ol{\Om})$.
\end{proof}
\end{lemma}

Next, we define a map $F_2:F_1(\ol{\Qbeta})\rightarrow \R^2$ so that
map $G_1:=F_2\circ F_1$ satisfies property \eqref{prop-G1} in $\ol{\Qbeta}$.

For each $\beta\in[0,\betadet]$, we define $F_2:F_1(\ol{\Qbeta})\rightarrow \R^2$ by
\begin{equation}
\label{12-14}
F_2(s,t):=\bigl(s, h_2(s,t)\bigr),
\end{equation}
and define a function $h_2:{F_1(\ol{\Qbeta})}\to [0, \infty)$  by
\begin{equation}
\label{12-15}
h_2(s,t):={\tleftchi}\sin^{-1}(\frac{t}{\leftu-s})+
(1-{\tleftchi})\big(t(1-{\trightchi})+{\trightchi}\sin^{-1}(\frac ts)\big)
\end{equation}
for the cut-off functions $\tleftchi, \trightchi\in C^{4}(\R)$ satisfying the following conditions:
\begin{equation*}
\label{definition-cutoff-tildexi-new2015}
\begin{split}
&{\tleftchi}(s)=\begin{cases}
1&\quad\text{for $s<\leftu-\leftch(1-\frac{1}{2k})$},\\[1mm]
0&\quad \text{for $s>\leftu-\leftch(1-\frac{1}{k})$},
\end{cases}\\[2mm]
&\trightchi(s)=
\begin{cases}
0&\quad \text{for $s<\rightc(1-\frac{1}{k})$},\\[1mm]
1&\quad \text{for $s>\rightc(1-\frac{1}{2k})$},
\end{cases}
\end{split}
\end{equation*}
\begin{equation*}
\begin{split}
&0\le \tleftchi, \trightchi\le 1, \qquad -\frac{4k}{\leftch}\le \tleftchi'\le 0 \le \trightchi'\le
\frac{4k}{\rightc},\qquad \tleftchi'\trightchi'=0,
\end{split}
\end{equation*}
where $k>4$ is the constant chosen to satisfy \eqref{10-e6} and all the properties used in the proof of Lemma \ref{lemma-10-1}.

Then $h_2$ satisfies
\begin{equation}
\label{h2-st-y}
h_2(s,t)=y\qquad\tx{ for $(s,t)$ near $F_1(\leftsonicdelta\cup\rightsonicdelta)$.}
\end{equation}

\begin{lemma}
\label{lemma-prop-F2}
There exist constants $C>0$ and $\kappa_1>0$ depending only on $(\iv, \gam)$ such that,
for each $\beta\in[0,\betadet]$, $F_2$ defined by \eqref{12-14} satisfies the following properties{\rm :}

\smallskip
\begin{itemize}
\item[(a)] $\|F_2\|_{C^4(F_1({\Qbeta}))}+ \|F_2^{-1}\|_{C^4(F_2\circ F_1({\Qbeta}))}\le C$,
and  ${\rm{det}} (DF_2)= \der_t h_2\ge \kappa_1$ in $\ol{F_1(\Qbeta)}${\rm ;}

\smallskip
\item[(b)] For $F_2(s,t):=(\til s, \til t)$,
$(F_2\circ F_1)(\Wedge)=\{(\til s, 0)\,:\, \til s\in (\leftu-\leftch, \rightc)\}$.
\end{itemize}

\smallskip
\begin{proof}
A direct computation by using \eqref{12-14} shows that
\begin{equation*}
\det(DF_2)
=\der_th_2(s,t)
=\frac{\tleftchi}{\sqrt{(\leftu-s)^2-t^2}}
+(1-\tleftchi)\Bigl((1-\trightchi)+\frac{\trightchi}{\sqrt{s^2-t^2}}\Bigr).
\end{equation*}
For $s<\leftu-\leftch(1-\frac{1}{2k})$, we can write
\begin{equation*}
\sqrt{(\leftu-s)^2-t^2}=r\cos y,
\end{equation*}
by \eqref{7-b4} and \eqref{10-e8},
where $r$ and $y$ are given by \eqref{def-r} and \eqref{coord-o} for $\xxi=F_1^{-1}(s,t)$.
Similarly, for $s>\rightc(1-\frac 1{2k})$, we can write as $\sqrt{s^2-t^2}=r\cos y$,
where $r$ and $y$ are given by \eqref{def-r} and \eqref{coord-n} for $\xxi=F_1^{-1}(s,t)$.
Then there exists a constant $\kappa_1>0$ depending only on $(\iv, \gam)$ such that
\begin{equation}
\label{h2-deiv-t}
\det(DF_2)=\der_th_2\ge \kappa_1\qquad \tx{in $\ol{F_1({\Qbeta})}$}.
\end{equation}
For a constant $a$, if
$F_1(Q_{\beta})\cap \{s=a\}$ is nonempty,
then \eqref{h2-deiv-t} implies that the one-dimensional map $(a,t)\in F_1(Q_{\beta})\cap \{s=a\}\mapsto h_2(a,t)$
is invertible. Then map $F_2$ given by \eqref{12-14} is also invertible.

The $C^4$--estimates of $F_2$ and $F_2^{-1}$ and property (b) are obtained directly from \eqref{7-a5} and \eqref{12-15}.
\end{proof}
\end{lemma}

By \eqref{12-15} and the invertibility of $F_2$, there exists a function $\til{h}_2: [\leftu-\leftch, \rightc]\to \R_+$ such that
\begin{equation*}
F_2^{-1}(s', t')=(s', \til{h}_2(s', t'))\qquad\,\mbox{for all $(s', t')\in \ol{(F_2\circ F_1)(\Qbeta)}$}.
\end{equation*}

For $F_1$ and $F_2$ given by \eqref{10-e7} and \eqref{12-14} respectively,
define a map $G_1:\ol{\Qbeta}\to [\leftu-\leftch,
\rightc]\times \R_+$ by
\begin{equation}
\label{12-16}
G_1:=F_2\circ F_1,
\end{equation}
and denote $G_1(\xxi)=(s',t')$. Map $G_1$ satisfies property \eqref{prop-G1}.

For each $\beta\in [0, \betadet]$, define
\begin{equation}
\label{def-sbeta}
\sbeta:=\leftu-\leftch.
\end{equation}
Note that $\sbeta$ varies continuously on $(\gam, \iv)$ and $\beta\in[0,\frac{\pi}{2})$.
Define a linear function $\Lb(s')$ by
\begin{equation}
\label{7-b8}
\Lb(s'):=\frac{2}{\rightc-\sbeta}(s'-\sbeta)-1.
\end{equation}
Then $\Lb$ maps $[\sbeta, \rightc]$ onto $[-1,1]$.
We define a map $\mcl{G}_1^{\beta}:\ol{\Qbeta}\to [-1,1]\times \R_+$  by
\begin{equation}
\label{12-16-mod}
\mcl{G}_1^{\beta}(\xxi)=(L_{\beta}(s'), t')\qquad\tx{for $(s',t')=G_1(\xxi)$}.
\end{equation}

\begin{lemma}
\label{lemma-12-1-mod}
There exist constants $C>0$ and $\kappa>0$ depending only on $(\iv, \gam)$ such that,
for any $\beta\in[0,\betadet]$,
$\mcl{G}_1^{\beta}$ defined by \eqref{12-16-mod} satisfies the following properties{\rm :}

\smallskip
\begin{itemize}
\item[(a)] $\|\mcl{G}_1^{\beta}\|_{C^4(\ol{\Qbeta})}+ \|(\mcl{G}_1^{\beta})^{-1}\|_{C^4({\mcl{G}_1^{\beta}(\ol{\Qbeta})})}\le C${\rm ;}

\smallskip
\item[(b)] $|\det (D\mcl{G}_1^{\beta})|\ge \kappa\quad$ in $\ol{\Qbeta}${\rm ;}

\smallskip
\item[(c)] $\mcl{G}_1^{\beta}(\Wedge)=\{(s,0)\,:\,s\in(-1, 1)\}${\rm ;}

\smallskip
\item[(d)] For $\iphi:=\ivphi+\frac 12|\xxi|^2$,
$\der_{t'}\iphi\bigl((\mcl{G}^{\beta}_1)^{-1}(s,t')\bigr)\le -\kappa<0$  for all $(s,t')\in \ol{\mcl{G}^{\beta}_1(\Qbeta)}$.
\end{itemize}
In addition, for any $\bar{\beta}\in(0,\betadet)$, there exists $m_2>0$ depending only on $(\iv, \gam,  \bar{\beta})$
such that any admissible solution $\vphi$ corresponding to $(\iv, \beta)\in\mathfrak{R}_{\rm weak}\cap\{0\le \beta\le \bar{\beta}\}$
satisfies
\begin{equation}
\label{7-a7-mod}
\der_{t'}(\ivphi-\vphi)((\mcl{G}_1^{\beta})^{-1}(s,t'))\le -m_2<0\,\qquad\text{in}\,\,\, \ol{\mcl{G}^{\beta}_1({\Om})}.
\end{equation}

\begin{proof}
Fix $\bar{\beta}\in(0,\betadet)$.
It follows
from \eqref{10-e8}, \eqref{12-15}, \eqref{12-16},
and Lemmas \ref{lemma-10-1} and \ref{lemma-prop-F2} that
there exist constants $C, \kappa_2>0$ depending only on $(\iv, \gam)$ such that,
for any $\beta\in[0,\betadet]$, map $G_1$ defined by \eqref{12-16} satisfies the following properties:

\smallskip
\begin{itemize}
\item[(a$'$)] $\|G_1\|_{C^4({\Qbeta})}+\|G_1^{-1}\|_{C^4({G_1(\Qbeta)})}\le C$;

\smallskip
\item[(b$'$)] $|\det (DG_1)|\ge \kappa_2$ in $\ol{\Qbeta}$;

\smallskip
\item[(c$'$)] $G_1(\Wedge)=\{(s',0)\,:\, s'\in(\leftu-\leftch, \rightc)\}$.
\end{itemize}
These properties, combined with \eqref{12-16-mod}, yield (a)--(c) for some $\kappa<\kappa_2$.

By \eqref{10-e7} and \eqref{12-14}--\eqref{12-16}, we find that, at $\xxi=G_1^{-1}(s',t')$,
\begin{equation*}
\der_{t'}(\ivphi-\vphi)( G_1^{-1}(s',t'))
=D_{\xxi}(\ivphi-\vphi)\cdot (\der_{t} \til h_1,1)\der_{t'}\til h_2
=\frac{D_{(\xin,\etan)}(\ivphi-\vphi)\cdot {\bm v} }{\der_t h_2}
\end{equation*}
for ${\bm v}$ given by \eqref{7-a2}--\eqref{7-a3}.
Then \eqref{7-a7-mod} follows by combining \eqref{12-7} and \eqref{est-t}
with Lemma \ref{lemma-prop-F2}(a) and \eqref{12-16-mod}.
 Assertion (d) can be verified similarly.
\end{proof}
\end{lemma}

By using \eqref{12-25} and the definitions of $(\ivphi, \leftvphi, \rightvphi)$ given in \eqref{def-uniform-ptnl-new},
it can be checked
that $\leftshock=\{\bmxi\,:\,(\ivphi-\leftvphi)(\bmxi)=0\}$ and $\rightshock=\{\bmxi\,:\,(\ivphi-\rightvphi)(\bmxi)=0\}$
intersect at a unique point:
\begin{equation}
\label{def-PI}
P_I=(\xi_1^{I}, \neta)\qquad\tx{for $\xi_1^I=-\frac{\etan^{(\beta)}-\neta}{\tan\beta}$},
\end{equation}
where $\etan^{(\beta)}$ is the $\etan$--intercept of $\leftshock$. Then
$\leftshock^{\delta_0}$ and $\rightshock^{\delta_0}$ intersect
at $(\xi_1^{I},\neta+\frac{\delta_0}{\iv})$.
It follows from \eqref{2-4-a6} and \eqref{iq-monotonicity} that $\frac{\dd\xi_2^{(\beta)}}{\dd\beta}>0$ for $\beta\in(0, \frac{\pi}{2})$
so that
\begin{equation}
\label{xi1-I-1}
\xi_1^I<0.
\end{equation}

Regarding \eqref{1-1} as an implicit function for $\oM$ with respect to $\iM(\neq 1)$, a direct computation yields
\begin{equation*}
  \frac{{\rm d}\oM}{{\rm d}\iM}=
  \left(\frac{\oM}{\iM}\right)^{\frac{4\gam}{\gam+1}}
  \frac{\iM^2-1}{\oM^2-1}\quad\tx{for $\iM\neq 1$}.
\end{equation*}
Next, we differentiate \eqref{2-4-a5} with respect to $\beta$, and use the representation of $\frac{d\oM}{d\iM}$ given in the right above to get
\begin{equation*}
  \frac{{\rm d}\iM}{{\rm d}\beta}=\frac{\iv \sec\beta\tan \beta}{1-\frac{\gam-1}{\gam+1}\left(\frac{\oM}{\iM}\right)^{\frac{2}{\gam+1}}
  -\frac{2(\gam-1)}{\gam+1}\left(\frac{\oM}{\iM}\right)^{\frac{3\gam+1}{\gam+1}}\frac{\iM^2-1}{\oM^2-1}}\quad\tx{for $\beta>0$}.
\end{equation*}
Then we apply L'{H}\^{o}pital's rule and use \eqref{2-4-a6} to derive from \eqref{def-PI} that
\begin{equation*}
  \lim_{\beta\to 0+}\xi_1^I=-\lim_{\beta\to 0+}\frac{\frac{{\rm d}\iM}{{\rm d}\beta}+\iM \tan \beta}{\sec \beta}=0.
\end{equation*}

Since point $P_I$ lies on $\leftshock$, and its $\xi_2$--coordinate is greater than
the $\xi_2$--coordinate of $\lefttop$, we have
\begin{equation}
\label{xi1-I-2}
\xi_1^{I}>\xi_1^{\lefttop}.
\end{equation}

By \eqref{12-25}, \eqref{7-a1}, and \eqref{7-a2}--\eqref{7-a3},
there exists a constant $m_3>0$ depending only on $(\iv, \gam)$ such that,
for each $\beta\in[0,\betadet]$,
\begin{equation}
\label{7-a8}
\begin{split}
&\der_{t'}\bigl((\ivphi-\leftvphi)\circ (\mcl{G}^{\beta}_1)^{-1}(s,t')\bigr)\le -m_3,\\
&\der_{t'}\bigl((\ivphi-\rightvphi)\circ (\mcl{G}_1^{\beta})^{-1}(s,t')\bigr)\le -m_3
\end{split}
\end{equation}
for all $(s,t')\in \ol{\mcl{G}_1^{\beta}(\Qbeta)}$.
By the implicit function theorem, there exists a unique function $f_{\beta}\in C^{0,1}([-1,1])$ such that
\begin{equation}
\label{7-c3}
\mcl{G}_1^{\beta}(\Qbeta)=\{(s,t')\,:\, -1<s<1, \;0<t'<f_{\beta}(s)\},
\quad \|f_{\beta}\|_{C^{0,1}([-1,1])}\le C
\end{equation}
for a constant $C$ depending only on $(\iv, \gam)$.

\begin{proposition}
\label{lemma-12-3}
Fix $\gam\ge 1$ and $\iv>0$. For each admissible solution $\vphi$ corresponding to $(\iv, \beta)\in\mathfrak{R}_{\rm weak}$,
there exists a unique function
\[
\gshock:[-1,1]\to \R_+
\]
satisfying the following properties{\rm :}

\smallskip
\begin{itemize}
\item[(a)]
$\mcl{G}^{\beta}_1(\Om)=\{(s,t')\,:\,-1<s<1, 0<t'< \gshock(s) \}$, \\[1mm]
 $\mcl{G}^{\beta}_1(\shock)=\{(s,\gshock(s))\,:\,-1<s<1\}${\rm .}

\smallskip
\item[(b)] For any constant $\hat{\eps}\in(0, \frac{1}{10}]$, there exists a constant $C_{\hat{\eps}}>0$
depending only on $(\iv, \gam)$ such that
\begin{equation*}
\|\gshock\|_{C^3([-1+\hat{\eps},1-\hat{\eps}])}\le C_{\hat{\eps}}.
\end{equation*}

\item[(c)] Let $\eps_0^*>0$ be the minimum of $\eps_0$ from Lemmas {\rm \ref{lemma1-sonic-N}}
and {\rm \ref{lemma-10-2}}.  For each $\eps\in(0, \eps_0^*]$, denote
\begin{equation}
\label{def-heps-mod}
\hat{\eps}:=\frac{2}{\rightc-\sbeta}\eps.
\end{equation}
Let $Q^{\beta}_0$ be the bounded region enclosed by $\leftsonic$, $\rightsonic$, $\leftshock$, $\rightshock$, and $\Wedge$.
Then
\begin{equation*}
\Om\subset Q^{\beta}_0 \subset \Qbeta
\end{equation*}
for $\Qbeta$ given by Definition {\rm \ref{definition-Qbeta}(iii)}.
For $\nD_{\eps}$ and $\oD_{\eps}$ defined by \eqref{definition-Dr-ext}, there exist unique functions $\gsn$ and $\gso$ so that
\begin{equation}
\label{12-21}
\begin{split}
&\mcl{G}_1^{\beta}(Q^{\beta}_0\cap \nD_{\eps})=\{(s,t')\,:\,1-\hat{\eps}<s< 1,\;\; 0<t'< \gsn(s)\},\\[1mm]
&\mcl{G}^{\beta}_1(Q^{\beta}_0\cap \oD_{\eps})=\{(s,t')\,:\,-1<s<-1+\hat{\eps},\;\; 0<t'<\gso(s)\},
\end{split}
\end{equation}
for $\hat{\eps}$ defined by \eqref{def-heps-mod}.
Moreover, there exists a constant $C>0$ depending only on $(\iv, \gam)$ such that
\begin{equation}
\label{7-b1}
\|\gsn\|_{C^3([1-\hat{\eps}^*_0, 1])}+ \|\gso\|_{C^3([-1, -1+\hat{\eps}^*_0])}\le C.
\end{equation}
For any $\alp\in(0,1)$, there exists $C_{{\rm par}}>0$ depending only on $(\iv, \gam,  \alp)$
such that, for any admissible solution corresponding to $(\iv, \beta)\in\mathfrak{R}_{\rm weak}$,
\[
\|\gsn-\gshock\|_{2,\alp,(1-\hat{\eps}^*_0,1)}^{\rm (par)}\le C_{\rm par},
\]
where the norm, $\|\cdot\|_{2,\alp, (1-\hat{\eps}^*_0,1)}^{\rm (par)}$, is defined
by Definition {\rm \ref{definition-parabolic-norm}(iii)}
with the replacement of $x$ by $1-|s|$ for the weight of the norm.

\smallskip
\item[(d)] For each $\bar{\beta}\in(0, \betadet)$, there exist $\bar{\alp}\in(0,1)$ and $C_{\bar{\beta}}>0$
depending only on $(\iv, \gam, \bar{\beta})$ such that, for any admissible solution corresponding
to $\beta\in[0, \bar{\beta}]$,
\begin{equation}
\label{7-b2}
\|\gshock\|_{2,\bar{\alp}, (-1, -1+\hat{\eps}_0^*)}^{(-1-\bar{\alp}),\{-1\}}\le C_{\bar{\beta}},\,\quad
 (\gshock-\gso)(-1)=0,\,\quad  (\gshock-\gso)'(-1)=0.
\end{equation}
Property \eqref{7-b2} is equivalent to
\begin{equation*}
\|\gshock-\gso\|_{2,\hat{\alp}, (-1,-1+\hat{\eps}_0^*)}^{(1+\hat{\alp}), (\rm par)}\le C'_{\bar{\beta}}
\end{equation*}
for a constant $C'_{\bar{\beta}}>0$ depending only on $(\iv, \gam, \bar{\beta})$, where the
norm, $\|\cdot \|_{2,\hat{\alp}, (-1,-1+\hat{\eps}_0^*)}^{(1+\hat{\alp}), (\rm par)}$,
is defined by Definition {\rm \ref{definition-parabolic-norm}(iv)} with the replacement of $x$ by $1-|s|$.

\smallskip
\item[(e)] For each $\bar{\beta}\in(0,\betadet)$, there exists a constant $\hat{k}>1$ depending
only on $(\iv, \gam,  \bar{\beta})$ such that, for any admissible solution $\vphi$ for $\beta\in[0,\bar{\beta}]$,
\[
\min\{\gshock(-1)+\frac{s+1}{\hat{k}}, \frac{1}{\hat{k}}\}
\le \gshock(s)\le \min\{f_{\beta}(s)-\frac{1}{\hat{k}}, \gshock(-1)+\hat k(s+1)\}
\]
for all $-1\le s\le 1$.
\end{itemize}

\begin{proof}
By \eqref{7-a7-mod} and the implicit function theorem, property (a) is obtained.
For an admissible solution $\vphi$, we differentiate the equation:
$
(\ivphi-\vphi)\circ (\mcl{G}_1^{\beta})^{-1}(s, \gshock(s))=0
$
with respect to $s$ to obtain
\begin{equation*}
\gshock'(s)=\frac{\der_{\rm s}\bigl((\ivphi-\vphi)\circ (\mcl{G}^{\beta}_1)^{-1}\bigr)}{\der_{t'}\bigl((\ivphi-\vphi)\circ (\mcl{G}^{\beta}_1)^{-1}\bigr)},
\end{equation*}
where the right-hand side is evaluated at $(s, \gshock(s))$.
Then property (b) is obtained from Lemma \ref{lemma-unif-est2}, Corollary \ref{corollary-unif-est-away-sn},
and Lemma \ref{lemma-12-1-mod}.
Similarly, properties (c) and (d) are obtained from \eqref{1-24}, \eqref{1-24ab},
and Propositions \ref{lemma-est-sonic-general-N}, \ref{proposition-sub8}, \ref{proposition-sub9},
and \ref{lemma-gradient-est}.

By Lemma \ref{lemma-10-2} and \eqref{prop-G1}, there exist constants $\hat{\eps}_1\in(0, \hat{\eps}_0^*]$ and $m>1$
depending only on $(\iv, \gam)$ such that, for each $\beta\in[0,\betadet]$, $\gso$ satisfies
\[
\frac 1m\le \gso'(s)\le m\,\qquad\tx{for all $-1\le s\le -1+\hat{\eps}_1$}.
\]
For each $\bar{\beta}\in(0,\betadet)$, by \eqref{7-b2}, we can choose $\hat{\eps}_2\in(0,\hat{\eps}_1]$
depending only on $(\iv, \gam,  \bar{\beta})$ such that,
for any admissible solution corresponding to $(\iv, \beta)\in\mathfrak{R}_{\rm weak}\cap\{0\le \beta\le \bar{\beta}\}$,
\[
\frac{1}{2m}\le \gshock'(s)\le 2m\,\qquad\tx{for $-1\le s\le -1+\hat{\eps}_2$}.
\]
By combining this estimate with Proposition \ref{proposition-sub3}, property (e) is obtained as a result.
\end{proof}
\end{proposition}

\begin{remark}
\label{remark-gshock-par-est-nr-leftsonic}
By Propositions {\rm \ref{lemma-est-sonic-general}} and {\rm \ref{proposition-sub8}},
for each $\alp\in(0,1)$, there exist constants $\hat{\eps}_3>0$ and $C_{\alp}>0$ depending
only on $(\iv, \gam,  \alp)$ such that,
for any admissible solution corresponding to $(\iv, \beta)$ with $0\le \beta<\betasonic$,
\[
\|\gshock-\gso\|_{2,\alp, (-1, -1+\hat{\eps}_3)}^{\rm (par)}\le C_{\alp},
\]
where the norm, $\|\cdot\|_{2,\alp, (-1,-1+\hat{\eps}_2)}^{\rm (par)}$,
is defined by Definition {\rm \ref{definition-parabolic-norm}(iii)}
with the replacement of $x$ by $1-|s|$ for the weight of the norm.

By Proposition {\rm \ref{proposition-sub9}}, for each $\alp\in(0,1)$, there exist
constants $\hat{\eps}_4>0$ and $C'_{\alp}>0$ depending only on $(\iv, \gam,  \alp)$
such that, for any admissible solution corresponding
to $(\iv, \beta)$ for $\betasonic\le \beta\le \betasonic+\sigma_3$,
\[
\|\gshock-\gso\|_{C^{2,\alp}([-1, -1+\hat{\eps}_4]}\le C'_{\alp},\qquad
\frac{\dd^m}{\dd s^m} (\gshock-\gso)(-1)=0\,\,\,\,\,\tx{for $m=0,1,2$}.
\]
\end{remark}

\smallskip
By \eqref{xi1-I-1}--\eqref{xi1-I-2},
$\xi_1^I$ given by \eqref{def-PI} satisfies that $\xi_1^{\lefttop}<\xi_1^I\le 0$ for any $\beta\in [0, \betadet]$.

\begin{definition}
\label{definition-smooth-connection}
Fix $\beta\in[0, \betadet]$. For $\xi_1^{I}$ given by \eqref{def-PI},
fix a smooth function $\chi^*_{\beta}$ such that
\begin{equation*}
\chi^*_{\beta}(\xin)=\begin{cases}
1&\tx{for $\xin \le \xi_1^{I}-\frac{\xi_1^I-\xi_1^{\lefttop}}{20}$},\\[1mm]
0&\tx{for $\xin \ge \xi_1^{I}$},
\end{cases}
\qquad  -\frac{10C}{\xi_1^I-\xi_1^{\lefttop}}\le (\chi_{\beta}^*)'\le 0,
\qquad \|\chi_{\beta}^*\|_{C^3(\R)}\le C
\end{equation*}
for some constant $C>0$ depending only on $(\iv, \gam)$.
For such a smooth cut-off function, define
\begin{equation}
\label{12-32}
\vphib^*(\xxi):=\leftvphi(\xxi)\chi^*_{\beta}(\xin)
+\rightvphi(\xxi)(1-\chi^*_{\beta}(\xin)),
\end{equation}
which is uniformly continuous when $\beta\to 0$.
\end{definition}

For later use, we list the following useful properties of $\vphi_{\beta}^*$ for $\beta\in[0,\betadet]${\rm :}

\smallskip
\begin{itemize}
\item[(i)] Define
\begin{equation}
\label{definition-vphi-beta}
\vphi_{\beta}:=\max\{\leftvphi, \rightvphi\}.
\end{equation}
By \eqref{def-uniform-ptnl-new} and the definition of $\xi_1^I$ given in \eqref{def-PI}, we have
\begin{equation*}
\vphi_{\beta}(\xi_1, \xi_2)=\begin{cases}
\leftvphi(\xi_1, \xi_2) \quad&\mbox{if $\xi_1<\xi_1^I$},\\
\leftvphi(\xi_1,\xi_2)=\rightvphi(\xi_1,\xi_2)\quad&\mbox{if $\xi_1=\xi_1^I$},\\
\rightvphi(\xi_1,\xi_2)\quad&\mbox{if $\xi_1>\xi_1^I$},
\end{cases}
\end{equation*}
so that
\begin{equation}
\label{vphi-beta-star}
\vphi_{\beta}^*\le \vphi_{\beta}\qquad \tx{in} \,\, \R^2.
\end{equation}

\item[(ii)] Let $\oD_r$ and $\leftch$ be given by \eqref{definition-Dr-ext} and Definition \ref{definition-Qbeta}, respectively.
Then there exists a sufficiently large constant $\bar{k}>1$ depending only on $(\iv, \gam)$ such that,
for any $\beta\in[0,\betadet]$, $\vphi_{\beta}^*$ satisfies
\begin{equation}
\label{7-b7}
\vphi_{\beta}^*=\vphi_{\beta}=
\begin{cases}
\leftvphi&\tx{in $\oD_{\frac{\leftch}{\bar{k}}}$},\\[1mm]
\rightvphi&\tx{in $\{\bmxi\in\R^2\,:\,\xi_1\ge 0\}$}.
\end{cases}
\end{equation}

\item[(iii)]
The set, $\{\bmxi\,:\, \xi_1^{\lefttop}< \xin < \xi_1^{\righttop},\;\; (\ivphi-\vphi_{\beta}^*)(\bmxi)=0\}$,
is contained in $\Qbeta$ and
\begin{equation}
\label{sup-inf}
\sup_{\Qbeta}(\ivphi-\vphi_{\beta}^*)-\inf_{\Qbeta}(\ivphi-\vphi_{\beta}^*)\ge \bar{\delta}>0
\end{equation}
for some constant $\bar{\delta}$ depending only on $(\iv, \gam)$.
\end{itemize}

\begin{lemma}\label{lemma-7-2}
There exists a constant $m>0$ depending only on $(\iv, \gam)$ such that
each $\vphib^*$ for $\beta\in [0,\betadet]$ satisfies
\begin{equation*}
\der_{t'}(\ivphi-\vphib^*)((\mcl{G}^{\beta}_1)^{-1}(s,t'))\le -m\,\qquad\tx{for all $(s,t')\in \ol{\mcl{G}_1^{\beta}(\Qbeta)}$}.
\end{equation*}

\begin{proof}
We have seen in the proof of Lemma \ref{lemma-12-1-mod} that
\begin{equation*}
\der_{t'}(\ivphi-\vphib^*)((\mcl{G}_1^{\beta})^{-1}(s,t'))=\frac{1}{\der_t h_2} D_{\bmxi}(\ivphi-\vphib^*)\cdot {\bm v}
\end{equation*}
for ${\bm v}$ given by \eqref{7-a2}--\eqref{7-a3},
where $D_{\bmxi}(\ivphi-\vphi_{\beta}^*)$ is evaluated at  $(\mcl{G}_1^{\beta})^{-1}(s,t')$.
By using \eqref{def-uniform-ptnl-new} and \eqref{12-32}, a direct computation yields that
\begin{equation*}
D_{\bmxi}(\ivphi-\vphib^*)=\iv \sec\beta(\sin\beta,-\cos\beta)\chi_{\beta}^*
+(0,-\iv)(1-\chi_{\beta}^*)+(\rightvphi-\leftvphi)(\chi_{\beta}^*)'(1,0).
\end{equation*}
From \eqref{7-a1} and \eqref{7-a2}--\eqref{7-a3}, there exists a constant
$m_*>0$  depending only on $(\iv, \gam)$ such that
\begin{equation}
\label{7-c1}
D_{\bmxi}(\ivphi-\leftvphi)\cdot \bm v\le -m_*,\quad
D_{\bmxi}(\ivphi-\rightvphi)\cdot \bm v\le -m_*
\qquad\tx{for all}\,\,(s,t')\in \ol{\mcl{G}_1^{\beta}(\Qbeta)}.
\end{equation}

By \eqref{10-e9-a}--\eqref{10-e9} and the definition of $\chi_{\beta}^*$,
we see that $\leftchi(\chi_{\beta}^*)'=\rightchi(\chi_{\beta}^*)'=0$ on $\R$.
This, combined with \eqref{7-a2}--\eqref{7-a3}, yields that
$
(\rightvphi-\leftvphi)(\chi_{\beta}^*)'(1,0)\cdot \bm v=0.
$
Then \eqref{7-c1} implies that
\begin{equation}
\label{2014-0822}
D_{\bmxi}(\ivphi-\vphib^*)\cdot{\bm v}\le -m_*\,\qquad\tx{for all $(s,t')\in \ol{\mcl{G}_1^{\beta}(\Qbeta)}$}.
\end{equation}
The proof is completed by \eqref{2014-0822} and Lemma \ref{lemma-prop-F2}.
\end{proof}
\end{lemma}

Each admissible solution $\vphi$ corresponding to $(\iv, \beta)\in\mathfrak{R}_{\rm weak}$ has
a unique function $\gshock:(-1,1)\rightarrow \R_+$ satisfying all the properties stated
in Proposition \ref{lemma-12-3}.
For such a function $\gshock$, define a map $G_{2,\gshock}:\mcl{G}_1^{\beta}(\Qbeta)\rightarrow \R^2$ by
\begin{equation}
\label{7-b9}
G_{2,\gshock}:\;(s,t')\mapsto \bigl(s,\frac{t'}{\gshock(s)}\bigr)=:(s,t).
\end{equation}
By Proposition \ref{lemma-12-3}(e), $G_{2,\gshock}$ is well defined and invertible with
\begin{equation*}
G_{2,\gshock}^{-1}(s,t)
=\bigl(s, t\gshock(s)\bigr).
\end{equation*}
More importantly, we have
\[
G_{2,\gshock}\circ \mcl{G}_1^{\beta}(\Om)=(-1,1)\times (0,1)=:\iter.
\]
Therefore, a function $u$ given by
\begin{equation}
\label{10-12}
u(s,t):=(\vphi-\vphib^*)\circ (\mcl{G}_1^{\beta})^{-1}\circ G_{2,\gshock}^{-1}(s,t)\qquad\tx{for $(s,t)\in\iter$}
\end{equation}
is well defined. To establish a uniform estimate of $u$ given by \eqref{10-12} for admissible solutions
corresponding to $(\iv, \beta)\in\mathfrak{R}_{\rm weak}$,
we introduce a new weighted $C^{2,\alp}$--norm in $\iter$.

\begin{definition}
\label{definition-unified-norms-u}
Fix constants $\sigma>0$, $\alp\in(0,1)$, and $m\in\mathbb{Z}_+$.

\smallskip
{\rm (i)} For ${\bf s}=(s,t), {\bf \til{s}}=(\til s, \til t)\in \iter$, define
\begin{equation*}
\delta_{\alp}^{(\rm{subs})}({\bf s}, \til{\bf s}):=\bigl((s-\til s)^2+(\max\{1-|s|, 1-|\til s|\})^2(t-\til t)^2\bigr)^{\frac{\alp}{2}}.
\end{equation*}
For an open set $U\subset \iter$, define
\begin{equation*}
\begin{split}
&\|u\|_{m,0,U}^{(\sigma),{\rm  (subs)}}
:=\sum_{0\le k+l\le m} \sup_{{\bf s}\in U}\, \big((1-|s|)^{k-\sigma} |\der_s^k\der_t^l u({\bf s})|\big),\\
&[u]_{m,\alp,U}^{(\sigma), {\rm (subs)}}
:=\sup_{{\bf s}\neq \til{\bf s}\in U} \Big(\min\bigl\{(1-|s|)^{\alp+k-\sigma}, (1-|\til s|)^{\alp+k-\sigma}\bigr\}
\frac{|\der_s^k\der_t^lu({\bf s})-\der_s^k\der_t^l u(\til{\bf s})|}{\delta_{\alp}^{\rm (subs)}({\bf s}, \til{\bf s})}\Big),\\
&\|u\|_{m,\alp, U}^{(\sigma), {\rm (subs)}}:= \|u\|_{m,0,U}^{(\sigma), {\rm (subs)}}
+ [u]_{m,\alp, U}^{(\sigma), {\rm (subs)}}.
\end{split}
\end{equation*}

{\rm (ii)} {\rm H\"{o}lder norms with parabolic scaling}:
For ${\bf s}=(s,t), {\bf \til{s}}=(\til s, \til t)\in \iter$, define
\begin{equation*}
\delta_{\alp}^{\rm (par)}({\bf s}, \til{{\bf s}}):=
\bigl((s-\til s)^2+\max\{1-|s|, 1-|\til s|\}(t-\til t)^2\bigr)^{\frac{\alp}{2}}.
\end{equation*}
For an open set $U\subset \iter$, define
\begin{equation*}
\begin{split}
&\|u\|_{m,0,U}^{(\sigma), {\rm (par)}}
:=\sum_{0\le k+l\le m}\sup_{{\bf s}\in U}\,\big( (1-|s|)^{k+\frac l2 -\sigma} |\der_s^k \der_t^l u({\bf s})|\big),\\
&[u]_{m,\alp, U}^{(\sigma), {\rm (par)}}
:= \sum_{k+l=m} \sup_{{\bf s}\neq \til{\bf s}\in U}\Big( \min\bigl\{(1-|s|)^{\alp+k+\frac l2-\sigma}, (1-|\til s|)^{\alp+k+\frac l2-\sigma}\bigr\}
\frac{|\der_s^k\der_t^l u({\bf s})-\der_s^k \der_t^l u(\til{\bf s})|}{\delta_{\alp}^{\rm (par)}({\bf s}, \til{\bf s})}\Big),\\
&\|u\|_{m,\alp, U} ^{(\sigma), {\rm (par)}}:=\|u\|_{m,0,U}^{(\sigma), {\rm (par)}}
+[u]_{m,\alp, U}^{(\sigma), {\rm (par)}}.
\end{split}
\end{equation*}
\end{definition}

\smallskip
For a constant $r\in(0,1)$, denote
\begin{equation}\label{4.1.50b}
\mcl{Q}_r^{\mcl{O}}:=\iter \cap\{-1<s<-1+r\},\,\,\,\mcl{Q}_r^{\mcl{N}}:=\iter\cap\{1-r<s<1\},\,\,\,
\mcl{Q}_r^{\rm{int}}:=\iter \cap\{|s|<1-r\}.
\end{equation}

\begin{remark}[Compact embedding properties of the norms in Definition \ref{definition-unified-norms-u}]
\label{remark-cpt-emd-prop}
For $m\in \mathbb{Z}_+$, $\alp\in[0,1)$, $\sigma>0$, and an open bounded set $U$ in $\R^2$,
let $C^{m,\alp}_{(\sigma), {\rm par}}(U)$ be the completion under the norm, $\|\cdot\|_{m,\alp, U}^{(\sigma), (\rm par)}$,
of the set of all smooth functions whose $\|\cdot\|_{m,\alp, U}^{(\sigma), (\rm par)}$--norms are finite.
Moreover, let $C^{m,\alp}_{(\sigma), (\rm subs)}(U)$ be the completion under the norm, $\|\cdot\|_{m,\alp, U}^{(\sigma), (\rm subs)}$,
of the set of all smooth functions whose $\|\cdot\|_{m,\alp, U}^{(\sigma), (\rm subs)}$--norms are finite.
Then the following compact embedding properties hold{\rm :}

\smallskip
\begin{enumerate}
\item[\rm (i)]
Let $r\in(0,1)$, $\alp, \hat{\alp}\in [0,1)$ with $\alp<\hat{\alp}$, and $m\in\{1,2\}$.
Then $C^{m,\hat{\alp}}_{(1+\hat{\alp}), (\rm sub)}(\mcl{Q}_r^{\mcl{O}})$ is compactly
embedded into $C^{m,{\alp}}_{(1+{\alp}), (\rm sub)}(\mcl{Q}_r^{\mcl{O}})$; see \cite[Corollary 17.2.7]{CF2}.

\smallskip
\item[\rm (ii)]
Let $m_1$ and $m_2$ be nonnegative integers, $\alp_1,\alp_2\in[0,1)$, and $m_1+\alp_1>m_2+\alp_2$,
and let $\sigma_1>\sigma_2>0$. Then $C^{m_1,\alp_1}_{(\sigma_1), (\rm par)}(U)$ is compactly
embedded into $C^{m_2,\alp_2}_{(\sigma_2), (\rm par)}(U)$; see \cite[Lemma 4.6.3]{CF2}.
\end{enumerate}
\end{remark}

For simplicity, let $\varepsilon_0$ denote $\varepsilon_0^*$ from Proposition \ref{lemma-12-3}. Define
\begin{equation}
\label{defi-epsilon-new}
\varepsilon_0'
:=\min_{\beta\in[0,\betadet]}\hat{\varepsilon}_0,
\end{equation}
for $\hat{\varepsilon}_0$ given by \eqref{def-heps-mod}.
%
%

\begin{proposition}
\label{proposition-unif-est-u-new}
For each $\bar{\beta}\in(0,\betadet)$, there exist constants $M>0$ and $\bar{\alp}\in(0, \frac 13]$
depending only on $(\iv, \gam,  \bar{\beta})$ such that,
for any admissible solution $\vphi$ corresponding
to $(\iv, \beta)\in\mathfrak{R}_{\rm weak}\cap\{0\le \beta\le \bar{\beta}\}$,
$u:\iter\rightarrow \R$ defined by \eqref{10-12} satisfies
\begin{equation}
\label{u-estimate}
\|u\|_{C^{2,\bar{\alp}}(\mcl{Q}^{\rm int}_{{\eps}'_0/4})}
+\|u\|^{(2),{\rm (par)}}_{2,\bar{\alp}, \mcl{Q}^{\mcl{N}}_{{\eps}'_0}}
+\|u\|^{(1+\bar{\alp}), {\rm (par)}}_{2,\bar{\alp},\mcl{Q}^{\mcl{O}}_{{\eps}'_0}}
+\|u\|^{(1+\bar{\alp}), {\rm (subs)}}_{1,\bar{\alp}, \mcl{Q}^{\mcl{O}}_{{\eps}'_0}}\le M.
\end{equation}

\begin{proof}  We divide the proof into six steps.

\smallskip
{\textbf{1.}} {\emph{Estimate of $u$ away from $s=-1$}}:
A direct computation by using Corollary \ref{corollary-unif-est-away-sn}, Proposition \ref{lemma-est-sonic-general-N},
Lemma \ref{lemma-10-1}, Proposition \ref{lemma-12-3}, \eqref{7-b7}, and \eqref{10-12} shows that,
for any $\alp\in(0,1)$, there exists a constant $M_1>0$ depending only on $(\iv, \gam,  \alp)$ such that
\begin{equation}
\label{u-estimate-away-ls-all}
\|u\|_{C^{2,{\alp}}(\mcl{Q}^{\rm int}_{{\eps}'_0/4})}
+\|u\|^{(2),{\rm (par)}}_{2,{\alp}, \mcl{Q}^{\mcl{N}}_{{\eps}'_0}}
\le M_1
\end{equation}
for any admissible solution $\vphi$ corresponding to $(\iv, \beta)\in\mathfrak{R}_{\rm weak}$.

\medskip
\textbf{2.} To obtain the {\it a priori} estimates of $u$ near $s=-1$,
the following two embedding inequalities from \cite{CF2} are applied in the next two steps:

\begin{lemma}[Lemma 17.2.10 in \cite{CF2}]
\label{lemma-appG-w2}
For a nonnegative integer $m$, $\alp\in(0,1)$, and $\sigma>0$,
let both norms $\|\cdot\|_{m,\alp, U}^{(\sigma), {\rm (subs)}}$
and $\|\cdot\|_{m,\alp, U}^{(\sigma), {\rm (par)}}$ be defined
in Definition {\rm \ref{definition-unified-norms-u}}.
For $r\in(0,1]$, there exists a constant $C>0$ independent of $(r, \alp)$ such that
\[
\|u\|_{m,\alp, \mcl{Q}^{\mcl{O}}_r}^{(\sigma), {\rm (par)}}
\le  \|u\|_{m,\alp, \mcl{Q}^{\mcl{O}}_r}^{(\sigma), {\rm (subs)}}.
\]
\end{lemma}

\begin{lemma}[Lemma 17.2.11 in \cite{CF2}]
\label{lemma-appG-w1}
For a nonnegative integer $m$, $\alp\in(0,\frac 13]$, $\sigma>0$, and $r\in(0,1)$,
there exists a constant $C>0$ independent of $(r, \alp)$ such that
\[
\|u\|_{1,\alp, \mcl{Q}^{\mcl{O}}_r}^{(1+\alp), {\rm (subs)}}
\le C \|u\|_{2,0, \mcl{Q}^{\mcl{O}}_r}^{(2),{\rm (par)}}.
\]
\end{lemma}

The estimates of $u$ near $s=-1$ for the admissible solution are given for two cases
separately: (i) $\beta\in[0, \betasonic)$ and (ii) $\beta\in[\betasonic, \bar{\beta}]$.

\medskip
{\textbf{3.}} {\emph{Estimate of $u$ near $s=-1$ for $\beta\in [0, \betac^{(\iv)})$}}:
For each $\beta\in[0, \betadet]$, by \eqref{prop-G1}, \eqref{12-16-mod},
and Definition \ref{definition-Gset-shocks-new}, we have
\begin{equation}
\label{u-vphi-relation}
u(s,t)=(\vphi-\leftvphi)(x,y)\qquad\tx{for}\,\,(s,t)\in \iter\cap\{-1<s<-1+\eps_0'\}
\end{equation}
with
\begin{equation*}
(s,t)=(L_{\beta}(x+\leftu-\leftc), \frac{y}{(\gshock\circ L_{\beta})(x+\leftu-\leftc)})
\end{equation*}
for the $(x,y)$--coordinates defined
by \eqref{coord-o}.
Differentiating \eqref{u-vphi-relation}, we have
\begin{equation}
\label{c_o_v}
\begin{split}
&u_s =\frac{\rightc-\sbeta}{2}\psi_x+t\gshock'\psi_y, \,\quad u_t = \gshock \psi_y,\\
&u_{ss} =\big(\frac{\rightc-\sbeta}{2}\big)^2\psi_{xx}+2t\gshock' \frac{\rightc-\sbeta}{2}\psi_{xy}+t\gshock''\psi_y+(t\gshock')^2\psi_{yy},\\
&u_{st} =\gshock'\psi_y+\frac{\rightc-\sbeta}{2}\gshock\psi_{xy}+t\gshock' \gshock\psi_{yy},\\
& u_{tt} =\gshock^2\psi_{yy}.
\end{split}
\end{equation}
A direct computation by using \eqref{10-12} and Propositions \ref{lemma-est-sonic-general} and \ref{proposition-sub8} shows that,
for $\beta\in [0, \betasonic)$ and  $\alp\in(0,1)$,
there exists a constant $C>0$ depending only on $(\iv, \gam, \alp)$ such that
\begin{equation}
\label{u-estimate-near-ls-supersonic1}
\|u\|^{(2),{\rm (par)}}_{2,\alp, \mcl{Q}^{\mcl{O}}_{\eps'_0}}\le C.
\end{equation}
Furthermore, \eqref{u-estimate-near-ls-supersonic1}, combined with Lemma \ref{lemma-appG-w1}, implies
that there exists a constant $M_2'>0$ depending only on $(\iv, \gam)$ such that
\begin{equation}
\label{u-estimate-near-ls-supersonic2}
\|u\|_{1,\frac 13, \mcl{Q}^{\mcl{O}}_{\eps_0'}}^{(1+\frac 13), {\rm (subs)}}\le M_2'
\end{equation}
for any admissible solution corresponding to $(\iv, \beta)\in\mathfrak{R}_{\rm weak}\cap\{0\le \beta<\betasonic\}$.
Combining the two estimates \eqref{u-estimate-near-ls-supersonic1}--\eqref{u-estimate-near-ls-supersonic2} together,
we have
\begin{equation}
\label{u-estimate-near-ls-supersonic3}
\|u\|_{2,\frac 13, \mcl{Q}^{\mcl{O}}_{\eps_0'}}^{(1+\frac 13), {\rm (par)}}
+\|u\|_{1,\frac 13, \mcl{Q}^{\mcl{O}}_{\eps_0'}}^{(1+\frac 13), {\rm (subs)}}\le M_2
\end{equation}
for a constant $M_2>0$ depending only on $(\iv, \gam)$.
\smallskip

\textbf{4}. {\emph{Estimate of $u$ near $s=-1$ for $\beta\in [\betac^{(\iv)}, \betac^{(\iv)}+\sigma_3]$}}:
Denote $\psi:=\vphi-\leftvphi$.
By Proposition \ref{proposition-sub9}, any admissible solution corresponding
to $(\iv, \beta)\in\mathfrak{R}_{\rm weak}\cap \{\betasonic\le \beta\le \betasonic+\sigma_3\}$ satisfies
\begin{equation}
\label{at-Pbeta}
\psi(P_{\beta})=|D\psi(P_{\beta})|=0.
\end{equation}

Regarding $\psi$ as a function of $(x,y)$ in $\oD_{\eps_0}$ for $\eps_0>0$ from Proposition \ref{proposition-sub7},
one can directly check by using \eqref{at-Pbeta} that $\psi$ satisfies the following estimate:
For ${\rm{\bf x}}=(x,y)$, ${\til{\rm{\bf x}}}=(\til x, \til y)\in \oD_{\eps_0}$,
\begin{equation}
\label{new-wet-est-psi-nr-ls}
\begin{split}
{\|\psi\|'}^{(-1-\alp)}_{2,\alp, \oD_{\eps_0}}:=&\;\sum_{0\le k+l \le2 } \sup_{\rm{\bf x}\in \oD_{\eps_0}} \big(|x-x_{P_{\beta}}|^{k+l-(1+\alp)}|\der_x^k\der_y^l \psi(\rm{\bf x})|\big)\\
&\;+\sum_{k=0}^2 \sup_{\rm{\bf x}, \til{\rm{\bf x}} \in \oD_{\eps_0}, \rm{\bf x}\neq \til{\rm{\bf x}}}
\Big(\min\{|x-x_{P_{\beta}}|, |\til x-x_{P_{\beta}}|\}\frac{|\der_x^k\der_y^{2-k}\psi({\rm{\bf x}})-\der_x^k\der_y^{2-k}\psi(\til{\rm{\bf x}})|}{|\rm{\bf x}- \til{\rm{\bf x}} |^{\alp}}\Big)\\
\le &\;  \kappa_1\|\psi\|_{2,\alp,\Om\cap \oD_{\eps_0}}^{(-1-\alp),\{P_{\beta}\}}
\end{split}
\end{equation}
for some constant $\kappa_1>0$ depending only on $(\iv, \gam, \alp)$.

Since $\gshock(-1)=0$ for $\beta\ge \betasonic$, Proposition \ref{lemma-12-3}(e) implies that
\begin{equation*}
\frac{1-|s|}{\hat k}\le \gshock(s) \le \hat k(1-|s|)\,\qquad \tx{for $s\in[-1, -1+\eps_0']$.}
\end{equation*}
Then, following the calculations in the proof of \cite[Lemma 17.2.5]{CF2},
we obtain from \eqref{c_o_v} and Remark \ref{remark-gshock-par-est-nr-leftsonic} that
\begin{equation*}
\|u\|_{2,\alp, \mcl{Q}^{\mcl{O}}_{\eps_0'}}^{(1+\alp), {\rm (subs)}}\le \kappa_2 {\|\psi\|'}^{(-1-\alp)}_{2,\alp, \oD_{\eps_0}}
\end{equation*}
for some constant $\kappa_2>0$ depending only on $(\iv, \gam, \alp)$.

By Corollary \ref{corollary-unif-est-away-sn} and Proposition \ref{proposition-sub9},
for each $\alp\in(0,1)$, there exists a constant $C>0$ depending only on $(\iv, \gam, \alp)$
such that any admissible solution corresponding
to $(\iv, \beta)\in\mathfrak{R}_{\rm weak}\cap \{\betasonic\le \beta\le \betasonic+\sigma_3\}$ satisfies
\begin{equation}
\label{estimate-psi-nr-Pbeta-subs-nr-sonic}
\|\psi\|_{2,\alp,\Om\cap \oD_{\eps_0}}^{(-1-\alp),\{P_{\beta}\}}
\le C
\end{equation}
for $\eps_0>0$ from Proposition \ref{proposition-sub7}.
Therefore, there exists a constant $M_3>0$ depending only on $(\iv, \gam,  \alp)$ such that
$u$ given by \eqref{10-12} associated with $\vphi$ satisfies
\begin{equation}
\label{u-estimate-near-ls-subsonic1}
\|u\|_{2, {\alp}, \mcl{Q}^{\mcl{O}}_{\eps_0'}}^{(1+{\alp}), {\rm (par)}}
\le \|u\|_{2, {\alp}, \mcl{Q}^{\mcl{O}}_{\eps_0'}}^{(1+\alp), {\rm (subs)}}\le M_3.
\end{equation}

\textbf{5}.  {\emph{Estimate of $u$ near $s=-1$ for $\beta\in [\betac^{(\iv)}+\frac{\sigma_3}{2}, \bar{\beta}]$}}:
By Propositions \ref{lemma-gradient-est} and \ref{lemma-12-3}, there exists $\hat{\alp}\in(0,1)$
depending on $(\iv, \gam,  \bar{\beta})$ so that $\psi =\vphi-\leftvphi$ still satisfies
estimate \eqref{estimate-psi-nr-Pbeta-subs-nr-sonic} for all $\beta\in [\betasonic+\frac{\sigma_3}{2}, \bar{\beta}]$
and $\alp\in (0,\hat{\alp}]$.
Then there exists $M_4>0$ depending only on $(\iv, \gam,  \bar{\beta})$ such that
any admissible solution $\vphi$ corresponding
to $(\iv, \beta)\in\mathfrak{R}_{\rm weak}\cap \{\betasonic+\frac{\sigma_2}{2}\le \beta\le \bar{\beta}\}$
satisfies estimate \eqref{u-estimate-near-ls-subsonic1} with $\alp=\hat{\alp}$ and $M_3=M_4$.

\smallskip
{\textbf{6.}} Finally, \eqref{u-estimate} is proved by choosing $\bar{\alp}=\min\{\hat{\alp}, \frac 13\}$
and $M=4\max\{M_1, M_2, M_3, M_4\}$.
\end{proof}
\end{proposition}

\section{Mapping the Functions in $\iter$ to Approximate Admissible Solutions}
\label{subsec-mapping-invert}
\quad
Fix $\gam \ge 1$ and $\iv>0$. For each $\beta\in[0, \betadet]$, let $\Qbeta$ be defined by
Definition {\rm \ref{definition-Qbeta}(iii)}.
For each $s^*\in(-1, 1)$, define
\begin{equation}
\label{definition-Qbeta-str}
\Qbeta(s^*):=\Qbeta\cap (\mcl{G}_1^{\beta})^{-1}\bigl(\{s=s^*\}\bigr).
\end{equation}
For each $\beta\in[0, \frac{\pi}{2})$, let $\vphib^*$ be defined by \eqref{12-32}. Then
\begin{equation*}
\inf_{\Qbeta(-1)} (\ivphi-\vphib^*)<0\le \sup_{\Qbeta(-1)}(\ivphi-\vphib^*).
\end{equation*}
In particular, the nonstrict inequality on the right above
becomes strict when $\beta<\betasonic$ and becomes an equality when $\beta\ge \betasonic$.

\begin{definition}
\label{definition-Gset-shocks-new}
Fix $\alp\in(0,1)$, $\bar{\beta}\in(0, \betadet)$, and $\beta\in (0, \bar{\beta}]$.
Let $u\in C^{1,\alp}(\ol{\iter})$ be a function satisfying that, for any $s\in (-1,1)$,
\begin{equation}
\label{zero-level-set}
\inf_{\Qbeta(s)}(\ivphi-\vphib^*)<u(s,1)<  \sup_{\Qbeta(s)}(\ivphi-\vphib^*).
\end{equation}
We define functions $\gshock^{(u,\beta)}$, $\mathfrak{F}^{(u,\beta)}$,
and $\vphi^{(u,\beta)}$ as follows{\rm :}

\smallskip
\begin{itemize}
\item[(i)] By Lemma {\rm \ref{lemma-7-2}}, for each $s\in(-1,1)$, there exists a unique
$\bar{t}'> 0$ such that
\[
(\ivphi-\vphib^*)\circ (\mcl{G}_1^{\beta})^{-1}(s,\bar t')=u(s,1).
\]
Define a function $\gshock^{(u,\beta)}: (-1,1)\to \R^+$ by
\begin{equation}
\label{definition-gshock-new2015}
\gshock^{(u,\beta)}(s)=\bar t'.
\end{equation}

\item[(ii)] For $\gshock^{(u,\beta)}$ from {\rm (i)},
define $G_{2,\gshock^{(u,\beta)}}$ by \eqref{7-b9}. For $\mcl{G}_1^{\beta}$ given by \eqref{12-16-mod},
define a map $\mathfrak{F}_{(u,\beta)}:\iter\to \Qbeta$ by
\begin{equation*}
\mathfrak{F}_{(u,\beta)}=(\mcl{G}_1^{\beta})^{-1}\circ G_{2,\gshock^{(u,\beta)}}^{-1}.
\end{equation*}

\item[(iii)] For $\mathfrak{F}_{(u,\beta)}$ from {\rm (ii)}, define the sets{\rm :}
\begin{equation*}
\begin{split}
&\shock(u,\beta):=\mathfrak{F}_{(u,\beta)}((-1,1)\times \{1\}),\qquad
\Om(u,\beta):=\mathfrak{F}_{(u,\beta)}(\iter).
\end{split}
\end{equation*}
Moreover, define a function $\vphi^{(u,\beta)}$ in $\Om(u,\beta)$ by
\begin{equation}
\label{def-vphi}
\vphi^{(u,\beta)}(\bmxi)=(u\circ \mathfrak{F}_{(u,\beta)}^{-1})(\bmxi)+\vphib^*(\bmxi)\;\;\qquad \tx{for all $\bmxi\in \Om(u,\beta)$}.
\end{equation}
\end{itemize}
\end{definition}

For $\alp\in(0,1)$ and $\bar{\beta}\in(0,\betadet)$, define
\begin{equation}
\label{12-74}
\mathfrak{G}_{\alp}^{\bar{\beta}}:=\left\{
(u,\beta)\in C^{1,\alp}(\ol{\iter})\times [0,\bar{\beta}]\,:\,
\begin{array}{ll}
\tx{$(u,\beta)$ satisfy \eqref{zero-level-set} for each $s\in (-1,1)$}\\[1.5mm]
\tx{and $(u, Du)(\pm 1,\cdot)=(0, {\bf 0})$}
\end{array}
\right\}.
\end{equation}

The next lemma follows from Definition \ref{definition-Gset-shocks-new}.
For details of the proof,
we refer to \cite[Lemmas 12.2.7 and 17.2.13]{CF2}.

\begin{lemma}
\label{lemma-7-4}
Fix $\alp\in(0,1)$ and $\bar{\beta}\in (0,\betadet)$.
For each $(u,\beta)\in\mathfrak{G}_{\alp}^{\bar{\beta}}$, the following properties hold{\rm :}

\smallskip
\begin{itemize}
\item[(a)] $\gshock^{(u,\beta)}\in C^{1,\alp}([-1,1])$.

\smallskip
\item[(b)] For domain $\Lbeta$ defined by Definition {\rm \ref{definition-domains-np}},
\[
\Om(u,\beta)\cup\shock(u,\beta) \subset \Qbeta\subset \Lbeta.
\]
Denote
$\lefttop=\mathfrak{F}_{(u,\beta)}(-1,1)$, $\righttop=\mathfrak{F}_{(u,\beta)}(1,1)$,
$\rightbottom=\mathfrak{F}_{(u,\beta)}(1,0)$, and $\leftbottom=\mathfrak{F}_{(u,\beta)}(-1,0)$.
Then $\shock(u,\beta)$ is a $C^{1,\alp}$--curve up to its endpoints $\lefttop$ and $\righttop$, and
is tangential to $\leftshock$ at $\lefttop$ and to $\rightshock$ at $\righttop$.
For $\hat f_{\mcl{O},0}$ and $\hat f_{\mcl{N},0}$ defined in Lemmas {\rm \ref{lemma1-sonic-N}}
and {\rm \ref{lemma-str-nr-sonic}},
\begin{equation}
\label{gshock-at-endpts}
\begin{split}
&\gshock^{(u,\beta)}(-1)=\hat f_{\mcl{O},0}(x_{\beta}),\qquad \qquad\qquad
\gshock^{(u,\beta)}(1)=\hat f_{\mcl{N},0}(0), \\
&\frac{\dd}{\dd s}\gshock^{(u,\beta)}(-1)=\frac{\rightc-\sbeta}{2}\hat f_{\mcl{O},0}'(x_{\beta}),\quad\,
\frac{\dd}{\dd s}\gshock^{(u,\beta)}(1)=-\frac{\rightc-\sbeta}{2}\hat f_{\mcl{N},0}'(0),
\end{split}
\end{equation}
where $\sbeta$ is defined by \eqref{def-sbeta} and $x_{\beta}$ is given by
\begin{equation*}
x_{\beta}=\begin{cases}
0\quad&\mbox{if $\beta<\betasonic$},\\
x_{P_{\beta}}\quad&\mbox{if $\beta\ge \betasonic$}.
\end{cases}
\end{equation*}
In the above, $P_{\beta}$ is the $\xi_1$--intercept of $\leftshock$,
and $x_{P_{\beta}}$ represents the $x$--coordinate of $P_{\beta}$ in
the $(x,y)$--coordinates defined by \eqref{coord-o}.
Note that $\frac{\dd^k}{\dd s^k} \gshock^{(u,\beta)}(\pm 1)$, $k=0,1$,
are uniquely determined depending only on $(\iv, \beta)$,
but independent of $u\in \mathfrak{G}_{\alp}^{\bar{\beta}}$.
Boundary $\der \Om(u,\beta)$ consists of $\Wedge=\mathfrak{F}_{(u,\beta)}\bigl((-1,1)\times \{0\}\bigr)$,
$\rightsonic=\mathfrak{F}_{(u,\beta)}\bigl(\{1\}\times (0,1)\bigr)$,
$\leftsonic=\mathfrak{F}_{(u,\beta)}\bigl(\{-1\}\times (0,1)\bigr)$,
and $\shock(u,\beta)=\mathfrak{F}_{(u,\beta)}\bigl((-1,1)\times \{1\}\bigr)$ which do not
intersect at the points of their relative interiors.

\smallskip
\item[(c)] Let $\delta_0>0$ be from Definition {\rm \ref{definition-Qbeta}}.
Let the $(x,y)$--coordinates be defined by \eqref{coord-o} near $\leftsonic$,
and by \eqref{coord-n} near $\rightsonic$.
For a constant $\eps>0$, define the two sets $\oOm_{\eps}$ and $\nOm_{\eps}$ by
\begin{equation*}
\begin{split}
&\oOm_{\eps}:=\mcl{N}_{\eps_0}(\Gamma_{\rm sonic}^{\mcl{O},\delta_0})\cap \{x_{\lefttop}<x<x_{\lefttop}+\eps\}\cap \Om(u,\beta),\\[1mm]
&\nOm_{\eps}:=\mcl{N}_{\eps_0}(\Gamma_{\rm sonic}^{\mcl{N},\delta_0})\cap \{0<x<\eps\}\cap \Om(u,\beta)
\end{split}
\end{equation*}
for $\eps_0>0$ to be fixed, where $\mcl{N}_{r}(\Gam)$ denotes an open $r$--neighborhood of $\Gam$.
Then there exists a constant $\eps_0>0$ depending only on $(\iv, \gam)$ such that
the following holds{\rm :}
for $L_{\beta}$ defined by \eqref{7-b8}, define the two functions $\fshocko$ and $\fshockn$ by
\[
\fshocko(x)=\gshock^{(u,\beta)}\circ L_{\beta}(x+\leftu-\leftc),\,\,\quad
\fshockn(x)=\gshock^{(u,\beta)}\circ L_{\beta}(\rightc-x).
\]
Then
\begin{equation*}
\begin{split}
&\oOm_{\eps}=\{(x,y)\,:\, x\in(x_{\lefttop},x_{\lefttop}+\eps),\, 0<y<\fshocko(x)\},\\
&\shock(u,\beta)\cap \der\oOm_{\eps}=\{(x,\fshocko(x))\,:\,x\in(x_{\lefttop}, x_{\lefttop}+\eps)\},\\
&\Wedge\cap\der\oOm_{\eps}=\{(x,0)\,:\, x\in(x_{\lefttop}, x_{\lefttop}+\eps)\}, \\
&\leftsonic=\leftsonic\cap \der \oOm_{\eps}=\{(x_{P_1}, y)\,:\,0<y<\fshocko(0)\},
\end{split}
\end{equation*}
and
\begin{equation*}
\begin{split}
&\nOm_{\eps}=\{(x,y)\,:\, x\in(0,\eps),\, 0<y<\fshockn(x)\},\\
&\shock(u,\beta)\cap \der\nOm_{\eps}=\{(x,\fshockn(x))\,:\,  x\in(0, \eps)\},\\
&\Wedge\cap\der\nOm_{\eps}=\{(x,0)\,:\,  x\in(0, \eps)\}, \\
&\rightsonic=\rightsonic\cap \der \nOm_{\eps}=\{(0, y)\,:\, 0<y<\fshockn(0)\}.
\end{split}
\end{equation*}

\item[(d)] Suppose that $(u,\beta), (\til{u}, \til{\beta})\in \mathfrak{G}_{\alp}^{\bar{\beta}}$ satisfy that
$\|(u, \til u)\|_{C^{1,\alp}(\ol{\iter})}< M$ for some constant $M>0$.
Then there exists a constant $C>0$, depending only on $(\iv, \gam, \bar{beta}, M, \alp)$, satisfying
the following estimates{\rm :}
\begin{align*}
&
\|\gshock^{(u,\beta)}\|_{C^{1,\alp}([-1,1])}+\|\mathfrak{F}_{(u,\beta)}\|_{C^{1,\alp}(\ol{\iter})}\le C,\\
&
\|\gshock^{(u,\beta)}-\gshock^{(\til u,\til{\beta})}\|_{C^{1,\alp}([-1,1])}\le
C\big(\|u-\til u\|_{C^{1,\alp}(\ol{\iter})}+|\beta-\til{\beta}|\big),\\
&
\|\mathfrak{F}_{(u,\beta)}-\mathfrak{F}_{(\til u,\til{\beta})}\|_{C^{1,\alp}(\ol{\iter})}\le
C\big(\|u-\til u\|_{C^{1,\alp}(\ol{\iter})}+|\beta-\til{\beta}|\big),\\
&
\|\vphi^{(u,\beta)}\circ \mathfrak{F}_{(u,\beta)}-{\vphi}^{(\til u, \til{\beta})}\circ \mathfrak{F}_{(\til u, \til{\beta})}\|_{C^{1,\alp}(\ol{\iter})}
\le C\big(\|u-\til u\|_{C^{1,\alp}(\ol{\iter})}+|\beta-\til{\beta}|\big),\\
&
\|(\vphi^{(u,\beta)}-\vphi_{\beta}^*)\circ \mathfrak{F}_{(u,\beta)}-({\vphi}^{(\til u, \til{\beta})}
 -\vphi_{\til{\beta}}^*)\circ \mathfrak{F}_{(\til u, \til{\beta})}\|_{C^{1,\alp}(\ol{\iter})}\\
&\quad \le C\big(\|u-\til u\|_{C^{1,\alp}(\ol{\iter})}+|\beta-\til{\beta}|\big).
\end{align*}

\item[(e)] $\psi^{(u,\beta)}:=\vphi^{(u,\beta)}-\max\{\leftvphi, \rightvphi\}=0$ holds on $\leftsonic\cup\rightsonic$.

\smallskip
\item[(f)] For $\eps>0$, let  $\varepsilon_0'$ 
be defined by \eqref{defi-epsilon-new}
Let $\eps_0>0$ be the constant from {\rm (c)}.
Assume that, for constants $\alp\in(0, 1)$, $\sigma\in(1,2]$, and $M>0$,
\begin{equation}
\label{estimate-u-new}
\|u\|_{2,\alp,\iter\cap\{|s|<1-\frac{\varepsilon_0'}{10}\}}+\|u\|^{(\sigma), {\rm (par)}}_{2,\alp,\iter\cap\{|s|>1-\varepsilon_0'\}}\le M.
\end{equation}
Then there exist $C>0$ depending only on $(\iv, \gam, \bar{\beta}, \alp, \sigma)$
and $C_0>0$ depending only on $(\iv, \gam, \bar{\beta})$ such that
\begin{equation}
\label{12-38}
\quad\, \|\gshock^{(u,\beta)}\|_{2,\alp,[-1+\frac{\varepsilon_0'}{10},1-\frac{\varepsilon_0'}{10}]}
+\|\gshock^{(u,\beta)}-\gso\|_{2,\alp,(-1, -1+\varepsilon_0')}^{(\sigma),{\rm (par)}}
+\|\gshock^{(u,\beta)}-\gsn\|_{2,\alp,(1-\varepsilon_0',1)}^{(\sigma),{\rm (par)}}\le CM,
\end{equation}
$\mathfrak{F}_{(0,\beta)}$ in $\{1-|s|<\eps_0\}\times (0, \infty)$ defined by
\begin{equation*}
\mathfrak{F}_{(0,\beta)}(s,t')=
\begin{cases}
\big(G_{2,\gso}\circ\mcl{G}_1^{\beta}\big)^{-1}(s,t')\qquad\tx{for $s\in (-1,-1+\varepsilon_0')$},\\
\big(G_{2,\gsn}\circ\mcl{G}_1^{\beta}\big)^{-1}(s,t')\qquad\tx{for $s\in (1-\varepsilon_0', 1)$}
\end{cases}
\end{equation*}
satisfies
\begin{align*}
&\|\mathfrak{F}_{(0,\beta)}\|_{C^3(\ol{\iter}\cap\{|s|\ge 1-\varepsilon_0'\})}\le C_0,\\
&\|\mathfrak{F}_{(u,\beta)}\|_{2,\alp,\iter\cap\{|s|<1-\frac{\varepsilon_0'}{10}\}}
+\|\mathfrak{F}_{(u,\beta)}-\mathfrak{F}_{(0,\beta)}\|_{2,\alp,\iter\cap\{|s|>1-\varepsilon_0'\}}^{(\sigma),{\rm (par)}}\le C.
\end{align*}

\item[(g)] Let $f_{\beta}$ be from \eqref{7-c3}.
For constants $M>0$ and $\delta_{\rm sh}>0$,
assume that $(u,\beta)\in \mathfrak{G}_{\alp}^{\bar{\beta}}$ satisfies \eqref{estimate-u-new}, 
$\gshock^{(u,\beta)}(-1)\le\delta_{\rm sh}$, 
and
\begin{equation*}
\qquad\quad \min\big\{\gshock^{(u,\beta)}(-1)+\frac{s+1}{M}, \delta_{\rm sh}\big\}\le \gshock^{(u,\beta)}(s)\le
\min\big\{\gshock^{(u,\beta)}(-1)+M(s+1), f_{\beta}(s)-\frac{1}{M}\big\}
\end{equation*}
for all $-1\le s\le 1$.
%
Then, for any $\eps\in(0, \frac 14\min\{\sbeta,\rightc\})$, there exists a constant $C_{\eps}>0$
depending only on $(\iv, \gam,  \alp, \bar{\beta}, \delta_{\rm sh}, \eps, M)$ such that
\begin{equation*}
\begin{split}
&\|\mathfrak{F}^{-1}_{(u,\beta)}\|_{2,\alp, \Om(u,\beta)\setminus (\ol{\oOm_{\eps}}\cup\ol{\nOm_{\eps}})}
+\|\mathfrak{F}^{-1}_{(u,\beta)}-\mathfrak{F}^{(-1)}_{(0,\beta)}\|_{2,\alp, \nOm_{\eps_0}}^{(\sigma), {\rm (par)}}\le C_{\eps},\\
&\|\vphi-\vphib^*\|_{2,\alp, \Om(u,\beta)\setminus (\ol{\oOm_{\eps}}\cup\ol{\nOm_{\eps}})}
+\|\vphi-\vphib^*\|_{2,\alp, \nOm_{\eps_0}}^{(\sigma), {\rm (par)}}\le C_{\eps}.
\end{split}
\end{equation*}

\item[(h)] Let $(u,\beta)$ and $(\til u, \til{\beta})$ be as in {\rm (d)}.
For any open set $K\Subset\iter$ so that $K\subset (-1+\delta,1-\delta)\times (0,1)$ for some $\delta>0$,
there exists a constant $C_{\delta}>0$ depending only on $(\iv, \gam,  \alp, \bar{\beta}, \sigma, \delta)$ such that
\begin{align*}
&\|\mathfrak{F}_{(u,\beta)}-\mathfrak{F}_{(\til u, \til{\beta})}\|_{C^{2,\alp}(\ol{K})}
\le C_{\delta}\big(\|(u-\til u)(\cdot, 1)\|_{C^{2,\alp}([-1+\delta, 1-\delta])}+|\beta-\til{\beta}|\big),\\
&\|\vphi^{(u,\beta)}\circ \mathfrak{F}_{(u,\beta)}-\vphi^{(\til u, \til{\beta})}\circ \mathfrak{F}_{(\til u, \til{\beta})}\|_{C^{2,\alp}(\ol{K})}
\le C_{\delta}\big(\|u-\til u\|_{C^{2,\alp}(\ol{K})}+|\beta-\til{\beta}|\big),\\
&\|\psi^{(u,\beta)}\circ \mathfrak{F}_{(u,\beta)}-{\psi}^{(\til u, \til{\beta})}\circ \mathfrak{F}_{(\til u, \til{\beta})}\|_{C^{2,\alp}(\ol{K})}
\le C_{\delta}\big(\|u-\til u\|_{C^{2,\alp}(\ol{K})}+|\beta-\til{\beta}|\big),
\end{align*}
where $\psi^{(u,\beta)}$ is given by $\psi^{(u,\beta)}:=\vphi^{(u,\beta)}-\vphib^*$ for each $(u,\beta)\in\mathfrak{G}_{\alp}^{\bar{\beta}}$.
\end{itemize}
\end{lemma}

\begin{remark}
By \eqref{prop-G1} and \eqref{gshock-at-endpts}, for any $(u,\beta)\in \mathfrak{G}_{\alp}^{\bar{\beta}}$, we have
\begin{equation*}
\gshock^{(u,\beta)}(1)=\sin^{-1}(\frac{\neta}{\rightc})>0
\end{equation*}
Fix $\delta\in(0,\betasonic)$, and suppose that $(u,\beta)\in\mathfrak{G}_{\alp}^{\bar{\beta}}$ and $\beta \in[0, \betasonic-\delta]$.
Then it follows from \eqref{fsonico-lwbd-lsonic}, \eqref{prop-G1}, and \eqref{gshock-at-endpts}
that there exists a constant $l_{\rm so}>0$ depending only on $(\iv, \gam, \delta)$ such that
\[
\gshock^{(u,\beta)}(-1)\ge l_{\rm so}.
\]
Therefore, there exists $b\in (0,1)$ depending only on $(\iv, \gam,  \sigma, \delta, M)$ such that,
for any $(u,\beta)\in\mathfrak{G}_{\alp}^{\bar{\beta}}$ with $\beta\in[0, \betasonic-\delta]$,
$\gshock^{(u,\beta)}$ satisfies
\begin{equation}
\label{gshock-bd-new-a}
b\le \gshock^{(u,\beta)}(s)\le b^{-1} \,\qquad\tx{for all $s\in[-1,1]$}.
\end{equation}
Then there exist $\hat{C}>0$ depending on $(\iv, \gam,  \alp,\sigma, \delta)$ and
$\hat{C}_0>0$ depending only on $(\iv, \gam,  \delta)$
such that
\begin{equation}
\label{estimate-iter-domain-2015j}
\begin{split}
&\|\mathfrak{F}^{-1}_{(0,\beta)}\|_{C^3(\ol{\Qbeta\cap \mcl{D}_{\eps_0}})}\le \hat C_0\qquad\;\;\tx{for}\,\,\,
\mcl{D}_{\eps_0}=\mcl{N}_{\eps_0}(\Gam_{\rm sonic}^{\mcl{O}, \delta_0})
\cup \mcl{N}_{\eps_0}(\Gam_{\rm sonic}^{\mcl{N}, \delta_0}) ,\\
&\|\mathfrak{F}^{-1}_{(u,\beta)}\|_{C^{2,\alp}(\ol{\Om(u,\beta)\setminus \mcl{D}_{{\eps_0}/{10}}})}
+\|\mathfrak{F}^{-1}_{(u,\beta)}-\mathfrak{F}^{-1}_{(0,\beta)}\|_{2,\alp, \Om(u,\beta)\cap \mcl{D}_{\eps_0}}^{(\sigma), {\rm (par)}}
\le \hat C M.
\end{split}
\end{equation}
Furthermore, $\vphi=\vphi^{(u,\beta)}$ defined by \eqref{def-vphi}
corresponding to $(u, \beta)$ satisfies
\begin{equation}
\label{estimate-vphi-u-new}
\|\vphi-\vphib^*\|_{C^{2,\alp}(\ol{\Om(u,\beta)\setminus \mcl{D}_{{\eps_0}/{10}}})}+
\|\vphi-\vphib^*\|_{2,\alp, \Om(u,\beta)\cap \mcl{D}_{\eps_0}}^{(\sigma), {\rm (par)}}\le \hat C M.
\end{equation}
\end{remark}

\smallskip
\section{Definition of the Iteration Set}
\label{subsection-def-iterset}
\begin{definition}
\label{pre-definition-new}
For $\eps_0>0$ from \emph{Lemma {\rm \ref{lemma-7-4}(c)}}, let $\varepsilon_0'$ be given by \eqref{defi-epsilon-new}.

\smallskip
\begin{itemize}
\item[(i)] Define $\un\in C^{3}(\ol{\iter})$ by \eqref{10-12} with $\beta=0$ and $\vphi=\rightvphi$.
Note that $\vphib^*\equiv \rightvphi$ in $Q^\beta$
by \eqref{12-32} because $\leftvphi=\rightvphi$ when $\beta=0$,
which yields that
\begin{equation*}
\un\equiv 0\qquad\tx{ in $\iter$}.
\end{equation*}

\item[(ii)]  For $\alp\in(0,1)$ and $\alp'\in(0,1]$, we introduce the norm{\rm :}
\begin{equation*}
\|u\|_{2,\alp,\iter}^{(*,\alp')}
:=\|u\|_{C^{2,\alp}(\ol{\mcl{Q}^{\rm int}_{\varepsilon'_0/4}})}
+\|u\|^{(1+\alp'), {\rm (par)}}_{2,\alp, \mcl{Q}^{\mcl{N}}_{\varepsilon_0'}}
+\|u\|^{(1+\alp), {\rm (par)}}_{2,\alp, \mcl{Q}^{\mcl{O}}_{\varepsilon_0'}}
+\|u\|^{(1+\alp), {\rm (subs)}}_{1,\alp, \mcl{Q}^{\mcl{O}}_{\varepsilon_0'}},
\end{equation*}
where $\mcl{Q}^{\rm int}_{\varepsilon_0'/4}$, $\mcl{Q}^{\mcl{N}}_{\varepsilon_0'}$, and $\mcl{Q}^{\mcl{O}}_{\varepsilon_0'}$
are defined in \eqref{4.1.50b}.
Denote by $C^{2,\alp}_{(*,\alp')}(\iter)$ the set of all $C^2(\iter)$--functions
whose $\|\cdot\|_{2,\alp,\iter}^{(*,\alp')}$--norms are finite.
Note that $C^{2,\alp}_{(*,\alp')}(\iter)$ is compactly embedded
into $C^{2,\til{\alp}}_{(*,\til{\alp}')}(\iter)$ whenever $0\le \til{\alp}<\alp<1$ and $0\le \til{\alp}'< \alp'\le 1$.
\end{itemize}
\end{definition}

For fixed $\gam\ge 1, \iv>0$, and $\beta_*\in(0, \betadet)$, we define the
iteration set $\mcl{K}\subset C^{1,\alp}(\ol{\iter})\times [0,\beta_*]$ for some appropriate $\alp\in(0,1)$.
For each
$\beta\in[0,\beta_*]$, $\mcl{K}_{\beta}:=\{u\in C^{1,\alp}(\ol{\iter})\,:\,(u,\beta)\in\mcl{K}\}$.
In the definition to come, the iteration set $\mcl{K}$ is given such that

\smallskip
\begin{itemize}
\item $\mcl{K}_0$ contains $\un$;

\smallskip
\item If $\beta$ is sufficiently close to $0$, then $u\in \mcl{K}_{\beta}$ is also close to $\un$ in an appropriate norm;

\smallskip
\item If $\beta$ is away from $0$, then any $\vphi^{(u,\beta)}$ given by \eqref{def-vphi} for $u\in \mcl{K}_{\beta}$
satisfies the strict directional monotonicity properties \eqref{3-c2}--\eqref{3-c3};

\smallskip
\item $\mcl{K}_{\beta}$ varies continuously on $\beta\in[0,\beta_*]$.
\end{itemize}

For $\gam\ge1$ and $\iv>0$, fix $\beta_*\in (0, \betadet)$. For $\bar{\alp}\in(0,\frac 13]$ 
from Proposition \ref{proposition-unif-est-u-new} with $\bar{\beta}$ replaced by $\beta_*$,
define
\begin{equation}
\label{alpha-choice}
\alp_*:=\frac{\bar{\alp}}{2}.
\end{equation}
Let $\eps_0>0$ be from  Lemma \ref{lemma-7-4}.
For constants $\alpha\in(0,\alpha_*]$, $\alpha_1\in(0,1)$,
$\delta_1$, $\delta_2$, $\delta_3$, $\eps\in(0, \frac{\eps_0}{2})$,
and $N_1>1$ to be specified later, we now define the iteration
set $\mcl{K}\subset C^{2,\alp}_{(*,\alp_1)}(\iter)\times [0, \beta_*]$.

\begin{definition}
\label{definition-10-3}
For fixed $\beta_*\in(0 , \betadet)$, the iteration set $\mcl{K}\subset C^{2,\alp}_{(*,\alp_1)}(\iter)\times [0, \beta_*]$
is the set of all $(u, \beta)$ satisfying the following properties{\rm :}

\smallskip
\begin{itemize}
\item[(i)] Fix
$\alp_1=\frac 78$.
Then $(u,\beta)$ satisfies
\begin{equation*}
\|u-\un\|_{2,\alp,\Q}^{(*,\alp_1)}<\mathscr{K}_1(\beta)
\end{equation*}
for $\mathscr{K}_1\in C^{0,1}(\R)$ given by
\begin{equation*}
\mathscr{K}_1(\beta)=\begin{cases}
\delta_1&\tx{if $\beta\le \frac{\delta_1}{N_1}$},\\[1mm]
N_0&\tx{if $\beta\ge \frac{2\delta_1}{N_1}$},\\[1mm]
\tx{linear}&\tx{if $\beta\in(\frac{\delta_1}{N_1}, \frac{2\delta_1}{N_1})$},
\end{cases}
\end{equation*}
with $N_0=\max\{10M, 1\}$ for constant $M$ from Proposition {\rm \ref{proposition-unif-est-u-new}}.

\smallskip
\item[(ii)] For  set $\mathfrak{G}_{\alp}^{\beta_*}$ defined by \eqref{12-74}, $(u,\beta)$ is
contained in $\mathfrak{G}_{\alp}^{\beta_*}$. Moreover,
let $\gshock=\gshock^{(u,\beta)}$, $\shock=\shock(u,\beta)$, $\Om=\Om(u,\beta)$,
and $\vphi=\vphi^{(u,\beta)}$ be defined by Definition {\rm \ref{definition-Gset-shocks-new}}.

\smallskip
\item[(iii)]
$\shock$ and $\gshock$ satisfy
\begin{equation*}
{\rm dist} (\shock, B_1(\Oi))>N_2^{-1},
\end{equation*}
\begin{equation}
\label{gshock-positivity-uterset-def2018}
\quad \min\{\gshock(-1)+ N_3^{-1}(s+1), N_3^{-1}\}   <\gshock(s)<\min\{\gshock(-1)+N_3(s+1),  f_{\beta}(s)-N_3^{-1}\}
\end{equation}
for all $-1< s < 1$
with $N_2=2C$ for $C$ from Proposition {\rm \ref{proposition-distance}},
and $N_3=2\hat{k}$ for $\hat{k}$ from Proposition {\rm \ref{lemma-12-3}(e)} with $\gshock(-1)\ge 0$,
where $f_{\beta}$ is defined by \eqref{7-c3}.

\smallskip
\item[(iv)] Let the $(x,y)$--coordinates be defined by \eqref{coord-n} near $\rightsonic$
and by \eqref{coord-o} near $\leftsonic$.
For $\vphib=\max\{\leftvphi, \rightvphi\}$, denote $\psi:=\vphi-\vphib$.
For $r>0$, let $\mcl{D}_r^{\mcl{O}}$ and $\mcl{D}_r^{\mcl{N}}$ be defined by \eqref{definition-Dr-ext}.
Let $\vphi$ and $\psi$ satisfy the following:
\begin{align}
\label{definiiton-iterset-ineq1}
\psi>\mathscr{K}_2(\beta)\qquad \text{in $\ol{\Om}\setminus (\mcl{D}_{\varepsilon_0'/10}^{\mcl{O}}\cup \mcl{D}_{\varepsilon_0'/10}^{\mcl{N}})$},&\\
\label{definiiton-iterset-ineq2}
\der_{\leftvec}(\ivphi-\vphi)<-\mathscr{K}_2(\beta)\;\qquad\qquad &\text{in $\ol{\Om}\setminus \mcl{D}^{\mcl{O}}_{\varepsilon_0'/10}$},\\
\label{definiiton-iterset-ineq3}
-\der_{\xi_1}(\ivphi-\vphi)<-\mathscr{K}_2(\beta)\;\qquad\qquad&\text{in $\ol{\Om}\setminus \mcl{D}^{\mcl{N}}_{\varepsilon_0'/10}$},\\
\label{definiiton-iterset-ineq4-rsonic}
|\der_x\psi(x,y)| <\frac{2-\mu_0}{1+\gam}x\quad\qquad\quad &\tx{in $\ol{\Om}\cap (\nD_{\eps_0}\setminus \nD_{\varepsilon_0'/10})$},\\
\label{definiiton-iterset-ineq4-lsonic}
|\der_x\psi(x,y)| <\mcl{K}_3(\beta)x\quad\qquad\quad &\tx{in $\ol{\Om}\cap (\oD_{\eps_0}\setminus \oD_{\varepsilon_0'/10} )$},\\
\label{definiiton-iterset-ineq5}
|\der_y \psi(x,y)|<N_4 x\qquad\qquad &\tx{in $\ol{\Om}\cap\big( (\oD_{\eps_0}\setminus\oD_{\varepsilon_0'/10})  \cup
 (\nD_{\eps_0}\setminus \nD_{\varepsilon_0'/10})\big)$},\\
\label{definiiton-iterset-ineq6}
|(\der_x\psi, \der_y\psi)|<N_4\eps \qquad\qquad&\tx{in $\ol{\Om}\cap (\ol{\mcl{D}^{\mcl{O}}_{\varepsilon_0'}}\cup \ol{\mcl{D}^{\mcl{N}}_{\varepsilon_0'}})$},\\
\label{definiiton-iterset-ineq7}
\|\vphi-\rightvphi\|_{C^{0,1}(\ol{\Om})}+\|\vphi-\leftvphi\|_{C^{0,1}(\ol{\Om})}&<N_5,\\
\label{definiiton-iterset-ineq8}
\der_{\bm\nu}(\ivphi-\vphi)>\mu_1, \quad \der_{\bm\nu}\vphi>\mu_1\qquad\;&\text{on $\ol{\shock}$},
\end{align}
for the unit normal vector $\bm{\nu}$ to $\shock$ towards the interior of $\Om$.
In the above conditions, functions $\mcl{K}_2, \mcl{K}_3\in C(\R)$ are defined by
\begin{equation}
\label{definition-K2b-2018}
 \mathscr{K}_2(\beta)=\delta_2\min\bigl\{\beta-\frac{\delta_1}{N_1^2}, \frac{\delta_1}{N_1^2}\bigr\},
\end{equation}
\begin{equation*}
\mcl{K}_3(\beta)=\begin{cases}
\frac{2-\mu_0}{1+\gam}&\tx{if $0\le \beta\le \betasonic+\frac{\sigma_2}{2}$},\\
\tx{linear}&\tx{if $\betasonic+\frac{\sigma_2}{2}<\beta<\betasonic+\sigma_2$},\\
N_4&\tx{if $\betasonic+\sigma_2 \le \beta$},
\end{cases}
\end{equation*}
for constants $\eps_0$, $\sigma_2, \mu_0,\mu_1$, $N_4$, and $N_5$ chosen as follows{\rm :}

\smallskip
\begin{itemize}
\item[(iv-1)] $\eps_0$ is from Lemma {\rm \ref{lemma-7-4}}.

\smallskip
\item[(iv-2)]$\sigma_2>0$ is from Lemma {\rm \ref{proposition-sub7}}, and
$ \mu_0=\frac{\delta}{2}$ for $\delta>0$ from Lemmas {\rm \ref{lemma-est-nrsonic}} and {\rm\ref{proposition-sub7}}.

\smallskip
\item[(iv-3)] $ \mu_1=\frac{\delta_1}{2}$ for ${\delta}_1>0$ from Corollary {\rm \ref{corollary-comp2}}.

\smallskip
\item[(iv-4)] Choice of $N_4${\rm :}
By \eqref{10-b9}--\eqref{10-c1}, for each $\sigma\in(0, \betadet-\betasonic)$,
\begin{equation}
\label{x-Pbeta}
\inf_{\betasonic+\sigma\le \beta<\betadet }x_{P_{\beta}}
= x_{P_{\beta}}|_{\beta=\betasonic+\sigma}=:x_{\sigma}>0.
\end{equation}
By Propositions  {\rm \ref{lemma-est-sonic-general}}, {\rm \ref{proposition-sub8}}, and {\rm \ref{proposition-sub9}},
there exists $C_1>0$ depending only on $(\iv, \gam)$ such that any admissible
solution $\vphi=\psi+\vphib$ for $\beta\in(0, \betasonic+\sigma_3]$
satisfies that $|(\der_x, \der_y)\psi(x,y)|\le C_1x$ in $\ol{\Om}\cap \mcl{D}^{\mcl{O}}_{\eps_0}$.

Let $\bar{\alp}\in(0,1)$ be from Proposition {\rm \ref{proposition-unif-est-u-new}}.
By Proposition {\rm \ref{lemma-gradient-est}} and \eqref{x-Pbeta},
any admissible solution $\varphi=\psi+\varphi_{\beta}$ for $\beta\ge \beta_s^{(\iv)}+\frac{\sigma_3}{2}$
satisfies
\begin{equation*}
|(\der_x, \der_y)\psi(x,y)|\le C_2x^{\bar{\alp}}\le C_2 \big(x_{P_{\beta}}|_{\beta=\betasonic+\frac{\sigma_3}{2}}\big)^{\bar{\alp}-1}x
\qquad\tx{in $\ol{\Om}\cap \mcl{D}^{\mcl{O}}_{\eps_0}$}
\end{equation*}
for a constant $C_2>0$ depending only on $(\iv, \gam, \beta_*)$.
Then there exists a constant $C^*_1>0$ depending only
on $(\iv, \gam,  \beta_*, \sigma)$
such that any admissible solution $\vphi=\psi+\vphib$ for $\beta\in(0, \beta_*]$ satisfies
\begin{equation*}
|(\der_x, \der_y)\psi(x,y)|
\le C^*_1x\qquad\tx{in $\ol{\Om}\cap \mcl{D}^{\mcl{O}}_{\eps_0}$}.
\end{equation*}
By combining this inequality with Proposition {\rm \ref{lemma-est-sonic-general-N}},
there exists a constant $C^*>0$ depending only on $(\iv, \gam,  \beta_*)$ such
that any admissible solution $\vphi$ for $\beta\in[0,\beta_*]$ satisfies
\begin{equation*}
|(\der_x, \der_y)\psi(x,y)|\le C^* x \qquad\tx{in $\ol{\Om}\cap (\mcl{D}^{\mcl{O}}_{\eps_0}\cup \mcl{D}^{\mcl{N}}_{\eps_0})$}.
\end{equation*}
We choose $N_4:=10C^*$.

\smallskip
\item[(iv-5)] By Lemma {\rm \ref{lemma-step3-1}} and the continuous dependence
of $\leftu$ and $\leftc$ on $\beta\in[0, \frac{\pi}{2})$, there exists a constant $\hat{C}>0$
depending only on $(\iv, \gam)$ such that any admissible solution $\vphi$ for $\beta\in[0,\betadet)$ satisfies
\begin{equation*}
\|\vphi-\rightvphi\|_{C^{0,1}(\ol{\Om})}+\|\vphi-\leftvphi\|_{C^{0,1}(\ol{\Om})}\le \hat C.
\end{equation*}
For such $\hat C>0$, we choose $N_5:=10\hat C$.
\end{itemize}

\smallskip
\item[(v)] Let $c(|D\vphi|^2,\vphi)$ be defined by
\begin{equation}
\label{def-sound}
c(|D\vphi|^2,\vphi):=\rho^{\frac{\gam-1}{2}}(|D\vphi|^2,\vphi)
\end{equation}
for $\rho(|{\bf p}|^2,z)$ given by \eqref{new-density}.
Then $\vphi$ satisfies
 \begin{equation}  \label{12-c7}
\frac{|D\vphi(\bmxi)|^2}{c^2(|D\vphi(\bmxi)|^2,\vphi(\bmxi))}
< 1-\tilde{\mu}\,{\rm dist}^{\flat}({\bm\xi}, \leftsonic\cup\rightsonic)
\qquad\tx{for $\bmxi\in\ol{\Om}\setminus ({\mcl D}^{\mcl{N}}_{{\eps}/{10}}\cup {\mcl D}^{\mcl{O}}_{{\eps}/{10}})$}.
\end{equation}
In \eqref{12-c7}, $\til{\mu}=\frac{\mu_{\rm el}}{2}$ for $\mu_{\rm el}>0$ from Remark {\rm \ref{remark-ellipticity-newly-added}}.

\smallskip
\item[(vi)] $\rho(|D\vphi|^2, \vphi)$ given by \eqref{new-density} satisfies
\begin{equation*}
\frac{a_*}{2}<\rho(|D\vphi|^2,\vphi)<2C
\qquad \text{in $\ol{\Om}\setminus (\mcl{D}^{\mcl{N}}_{\varepsilon_0'/10}\cup {\mcl D}^{\mcl{O}}_{\varepsilon_0'/{10}})$},
\end{equation*}
for $a_*=(\frac{2}{\gam+1})^{\frac{1}{\gam-1}}$ and $C$ from \eqref{3-c7} in Lemma {\rm \ref{lemma-step3-1}}.
For such constants, denote
\begin{equation*}
\rho_{\min}:=\frac{a_*}{2},\qquad \rho_{\rm max}=2C.
\end{equation*}

\item[(vii)] The boundary value problem
\begin{equation}\label{12-50}
\begin{cases}
\mcl{N}_{(u,\beta)}(\hat{\phi}):=A_{11}\hat{\phi}_{\xin\xin}+2A_{12}\hat{\phi}_{\xin\etan}
+A_{22}\hat{\phi}_{\etan\etan}=0 &\mbox{in $\Om$},\\[1mm]
\M_{(u,\beta)}(D\hat{\phi},\hat{\phi},\bmxi)=0 &\text{on $\shock$},\\[1mm]
\hat{\phi}=\max\{\rightvphi, \leftvphi\}-\rightvphi &\mbox{on $\leftsonic\cup\rightsonic$},\\[1mm]
\hat{\phi}_{\etan}=0 &\text{on $\Wedge$}
\end{cases}
\end{equation}
has a unique solution $\hat{\phi}\in C^2(\Om)\cap C^1(\ol{\Om})$,
where $\mcl{N}_{(u,\beta)}$ and $\M_{(u,\beta)}$ are determined
by $(u,\beta)$ in \S {\rm \ref{section-eqn-bdry-itr}}.
Moreover, this solution satisfies
that
$\hat u(s,t)$, defined by
\begin{equation}
\label{definition-sol-ubvp}
\hat u(s,t):=(\hat{\phi}+\rightvphi-\vphib^*)\circ \mathfrak{F}_{(u,\beta)}(s,t)\qquad\text{in $\iter$},
\end{equation}
satisfies
\begin{equation}
\label{12-65}
\|\hat u-u\|^{(*, \alpha_1)}_{2,\frac{\alpha}{2}, \iter}
<\delta_3.
\end{equation}
\end{itemize}
\end{definition}

\begin{remark}
\label{remark-vphi-beta}
By \eqref{7-b7}, the boundary condition
$\hat{\phi}=\max\{\rightvphi, \leftvphi\}-\rightvphi$ on $\leftsonic\cup\rightsonic$ given in \eqref{12-50} is equivalent to
\begin{equation*}
\hat{\phi}
=\begin{cases}
\leftvphi-\rightvphi\quad\,\,&\mbox{on $\leftsonic$},\\[1mm]
0\quad&\mbox{on $\rightsonic$}.
\end{cases}
\end{equation*}
\end{remark}

\begin{remark}
\label{remark-a}
For a fixed $\beta_*\in(0, \betadet)$,
let the iteration set $\mcl{K}$ be defined by Definition {\rm \ref{definition-10-3}}.
For each $(u,\beta)\in \mcl{K}$, let $\gshock=\gshock^{(u,\beta)}$, $\Om=\Om(u,\beta)$, $\shock=\shock(u,\beta)$,
and $\vphi=\vphi^{(u,\beta)}$ be defined by {\rm{Definition \ref{definition-Gset-shocks-new}}}.
Then there exist constants $M_{\rm dom}>0$ depending only on $(\iv, \gam)$,
$C>0$ depending only on $(\iv, \gam,  \alp)$, and $C_{\beta_*}>0$ depending only
on $(\iv, \gam,  \beta_*,\alp)$ such that the following properties hold{\rm :}

\smallskip
\begin{itemize}
\item[(i)] Let $\gso$ and $\gsn$ be from \eqref{12-21}.
For $N_0$ from {\rm{Definition \ref{definition-10-3} (i),
$\gshock$ satisfies
\begin{equation}
\label{iterset-gshocko-est}
\begin{split}
&\|\gshock\|_{2,\alp, (-1,1)}^{(-1-\alp), \{\pm 1\}}\le CN_0,\\
&\frac{\dd^k}{\dd s^k}(\gshock-\gso)(-1)=\frac{\dd^k}{\dd s^k}(\gshock-\gsn)(1)=0 \qquad \tx{for}\,\,k=0,1.
\end{split}
\end{equation}
}}

\item[(ii)] $\shock$ is a $C^{1,\alp}$--curve up to its endpoints.
Furthermore, $\shock\cap \oD_{\eps_0}$ and $\shock\cap \nD_{\eps_0}$ are
graphs $y=\fshocko(x)$ and $y=\fshockn(x)$ for
\begin{equation}
\label{shock-ftns-nr-sonic-2016sep}
\fshocko(x)=(\gshock\circ L_{\beta}^{-1})(\sbeta+x),\qquad
\fshockn(x)=(\gshock\circ L_{\beta}^{-1})(\rightc-x),
\end{equation}
with $\fshockn$ and $\fshocko$ satisfying that
\begin{equation*}
    \|\fshockn-\hat{f}_{\mcl{N},0}\|_{2,\alp,(0,\eps_0)}^{(1+\alp_1),{\rm (par)}}
    +\|\fshocko-\hat{f}_{\mcl{O},0}\|_{2,\alp,(0,\eps_0)}^{(1+\alp),{\rm (par)}}
    < C\mathscr{K}_1(\beta)
    \end{equation*}
   for $\hat{f}_{\mcl{N},0}$ and $\hat{f}_{\mcl{O},0}$
   from Lemmas {\rm \ref{lemma1-sonic-N}(e)} and {\rm \ref{lemma-str-nr-sonic}(e)}, respectively.

\smallskip
\item[(iii)] $\Om\subset B_{M_{\rm dom}}(\bm 0)$.

\smallskip
\item[(iv)] $\psi=\vphi-\vphib^*$ satisfies
\begin{equation*}
\begin{split}
&D^k\psi=0\qquad\tx{on}\,\,\,\ol{\leftsonic}\cup\ol{\rightsonic}\qquad\tx{for}\,\,k=0,1,\\[1mm]
&\|\psi\|_{C^{1,\alp}(\ol{\Om})}< C\mathscr{K}_1(\beta).
\end{split}
\end{equation*}
By Lemma {\rm \ref{lemma-str-nr-sonic}(e)} and \eqref{iterset-gshocko-est}, we can adjust $\eps_0$ depending on $(\iv, \gam)$ to satisfy
\begin{equation*}
0<\frac 12 \gso'(-1)\le \gshock'(s)\le 4\gso'(-1)\qquad\tx{for all $s\in[-1, -1+\hat{\eps}_0]$}.
\end{equation*}
Then, for each $\beta <\betasonic$,
\begin{equation*}
\qquad\quad\, |\der_y\psi(x,y)|=\frac{|u_t(s,t)|}{\gshock(s)}\le
\|u\|^{(1+\alp), (\rm{subs})}_{1,\alp, \mcl{Q}^{\mcl{O}}_{\hat{\eps}_0}}\frac{(1-|s|)^{1+\alp}}{\gshock(-1)}\le
 Cx^{\frac 12+\alp}\qquad \tx{for $(x,y)\in \Om\cap \oD_{\hat{r}}$},
\end{equation*}
where
$\hat{r}=\min\{\gshock^2(-1), \eps_0\}$
{\rm (}note that $\gshock(-1)>0$ for
each $(u,\beta)\in \mcl{K}\cap \{\beta<\betac^{(\iv)}\}${\rm )}.
For each $\sigma\in(0, \betasonic)$, there exists a constant $N_0^*(\sigma)$ depending
only on $(\iv, \gam,  \beta_*, \sigma)$ such that,
if $(u,\beta)\in\mcl{K}\cap \{\beta<\betac^{(\iv)}-\sigma\}$, then
\begin{equation*}
\|\psi\|_{2,\alp, \oD_{\eps_0}}^{(1+\alp), {\rm (par)}}< N_0^*(\sigma).
\end{equation*}

\item[(v)] For each $r\in(0, \eps_0)$, there exists a constant $C_{\beta_*,r}>0$
depending only on $(\iv, \gam,  \beta_*, r, \alp)$ such that
\begin{equation*}
\|\vphi\|_{C^{2,\alp}(\ol{\Om\setminus (\oD_{r}\cup\nD_{r})})}<C_{\beta_*,r}.
\end{equation*}
\end{itemize}
\end{remark}

\begin{definition}
\label{definition-K-ext}
Define the following sets{\rm :}

\smallskip
\begin{itemize}
\item[(i)] Denote $\Kext$ as
\begin{equation}
\label{def-kext-subs-new}
\Kext:=\{(u,\beta)\in C^{2,\alp}_{(*,\alpha_1)}(\Q)\,:\,\mbox{$(u,\beta)$ satisfy Definition {\rm \ref{definition-10-3}(i)--(vi)}}\};
\end{equation}

\item[(ii)] $\ol{\mcl{K}}$ and
$\ol{\Kext}$ are the closures of $\mcl{K}$ and $\Kext$ in $C^{2,\alp}_{(*, \alpha_1)}(\iter)\times [0, \beta_*]$, respectively;

\smallskip
\item[(iii)] For each $\mcl{C}\in \{\mcl{K}, \Kext, \ol{\mcl{K}}, \ol{\Kext}\}$ and each $\beta\in[0, \beta_*]$, denote
\begin{equation*}
\mcl{C}_{\beta}:=\{u\,:\,(u,\beta)\in \mcl{C}\}.
\end{equation*}
Note that $\mcl{C}_{\beta}\subset C^{2,\alp}_{(*, \alpha_1)}(\iter)$.
\end{itemize}
\end{definition}

\begin{remark}
\label{lemma-closure}
Each $(u,\beta)\in \ol{\Kext}$ satisfies property {\rm (ii)} of {\rm{Definition \ref{definition-10-3}}},
as well as properties {\rm (i)} and {\rm (iii)}--{\rm (vi)} of {\rm{Definition \ref{definition-10-3}}},
and all the properties stated in {\rm{Remark \ref{remark-a}}} with nonstrict inequalities in the estimates.
\end{remark}

\section{Boundary Value Problem \eqref{12-50}}
\label{section-eqn-bdry-itr}
In order to complete Definition \ref{definition-10-3}, it remains to define the nonlinear differential
operators $\mcl{N}_{(u,\beta)}$ and $\M_{(u,\beta)}$ in \eqref{12-50} for each $(u,\beta)\in \mcl{K}$.

For each $(u,\beta)\in \Kext$, let $\gshock=\gshock^{(u,\beta)}$, $\mathfrak{F}=\mathfrak{F}_{(u,\beta)}$,
$\Om=\Om(u,\beta)$, and $\shock=\shock(u,\beta)$,
and let $\vphi=\vphi^{(u,\beta)}$ be defined
by \eqref{def-vphi}.

\subsection{Definition of $\mcl{N}_{(u,\beta)}$ in \eqref{12-50}}
\label{subsubsec-defA}

For $\rightvphi$ defined by \eqref{def-uniform-ptnl-new}, denote
\[
\phi:=\vphi-\rightvphi.
\]
For a $C^2$--function $\hat{\phi}$ in $\Om$, we define $\mcl{N}_{(u,\beta)}(\hat{\phi})$ by
\begin{equation}
\label{definition-nub-hphi}
\mcl{N}_{(u,\beta)}(\hat{\phi}):=\sum_{i,j=1}^2 A_{ij}(D\hat{\phi}, \bmxi)\der_{\xi_i\xi_j}\hat{\phi}
\end{equation}
so that the following properties hold:

\smallskip
\begin{itemize}
\item Equation $\mcl{N}_{(u,\beta)}(\hat{\phi})=0$  is strictly elliptic in $\ol{\Om}\setminus (\ol{\leftsonic}\cup\ol{\rightsonic})$;

\smallskip
\item If $\phi$ is a solution of \eqref{12-50}, then equation $\mcl{N}_{(u,\beta)}({\phi})=0$  coincides with \eqref{2-3}.
\end{itemize}

\smallskip
The coefficient functions $A_{ij}({\bf p}, \bmxi), i,j=1,2$, of the nonlinear operator $\mcl{N}_{(u,\beta)}$
are defined in the following six steps:

\medskip
{\textbf{1.}} For a constant $r>0$,
let $\oD_r$ and $\nD_r$ be defined by \eqref{definition-Dr-ext},
and let $\mcl{D}_r:=\oD_r\cup\nD_r$.
Let $\eps_0>0$ be from Lemma \ref{lemma-7-4}.
For a constant $\eps_{\rm eq}\in (0, \frac{\eps_0}{2})$ to be chosen later, we define $A^{(1)}_{ij}(\bmxi)$
for $\bmxi\in\Om\setminus {\mcl{D}}_{{\eps_{\rm eq}}/{10}}$ by
\begin{equation}
\label{def-A1}
A^{(1)}_{ij}(\bmxi):=A^{\rm potn}_{ij}(D\phi(\bmxi), \phi(\bmxi),\bmxi),
\end{equation}
where
\begin{equation}\label{def-A-potn-new}
\begin{split}
&A^{\rm potn}_{11}({\bf p}, z, \bmxi)=c^2-(p_1+\der_{\xin}\rightvphi)^2,\\[1mm]
&A^{\rm potn}_{12}({\bf p}, z, \bmxi)=A^{\rm potn}_{21}({\bf p}, z, \bmxi)
=-\big(p_1+\der_{\xin}\rightvphi(\bmxi)\big)\big(p_2+\der_{\etan}\rightvphi(\bmxi)\big),\\[1mm]
&A^{\rm potn}_{22}({\bf p}, z, \bmxi)=c^2-(p_2+\der_{\etan}\rightvphi(\bmxi))^2
\end{split}
\end{equation}
for $c^2=c^2(|{\bf p}+D\rightvphi|^2,z+\rightvphi)$ given by \eqref{def-sound}.

\smallskip
{\textbf{2.}} For $\mu_0>0$ from Definition \ref{definition-10-3}(iv-1), fix a function $\zeta_1\in C^{3}(\R)$ such that
\begin{align}
\label{4-7}
&\zeta_1(s):=\begin{cases}
s&\quad\text{if $|s|\le\frac{2-\frac{\mu_0}{5}}{1+\gam}$},\\[1mm]
\frac{(2-\frac{\mu_0}{10}){\rm sgn}(s)}{1+\gam}&\quad\text{if $|s|> \frac{2}{1+\gam}$},
\end{cases}\\[1.5mm]
&0\le \zeta_{1}'(s)\le 10,\quad \zeta_{1}(-s)=-\zeta_{1}(s)\,\,\qquad\tx{for all $s\in\R$}, \label{7-c5}\\[1.5mm]
&-\frac{20(1+\gam)}{\mu_0}\le \zeta_1''(s)\le 0\,\, \qquad\tx{for all $s\ge 0$}.\label{7-c5-a}
\end{align}
Define $\cbeta$, $u_{\beta}$, $r$, and $\phi_{\beta}$ by
\begin{align}
\label{def-cbeta}
&(\cbeta, u_{\beta}):=
\begin{cases}
(\leftc,\leftu)&\quad \tx{in $\oD_{2\eps_{\rm eq}}$},\\[1.5mm]
(\rightc, 0)&\quad \tx{in $\nD_{2\eps_{\rm eq}}$},
\end{cases}\\[1mm]
&r=\sqrt{(\xin-u_{\beta})^2+\etan^2},\\[1mm]
\label{def-phibeta-star}
&\phi_{\beta}:=\vphi_{\beta}^*-\rightvphi
\end{align}
for $\vphi_{\beta}^*$ given by \eqref{12-32}.

Denote
$\psi:=\phi-\phi_{\beta}=\vphi-\vphi_{\beta}^*$.
Suppose that $\hat{\phi}$ is a solution of \eqref{12-50}. We denote
\begin{equation}
\label{def-hat-Psi}
\hat{\psi}:=\hat{\phi}-\phi_{\beta}.
\end{equation}

Let the $(x,y)$--coordinates be defined by \eqref{coord-n} and \eqref{coord-o} in $\nD_{2\eps_{\rm eq}}$ and $\oD_{2\eps_{\rm eq}}$, respectively.
For ${\bf p}\in \R^2$, denote
\begin{equation*}
{\bf p}':={\bf p}-D_{(x,y)}\phi_{\beta}.
\end{equation*}
Note that ${\bf p}'=\bf p$ in $\nD_{2\eps_{\rm eq}}$ and
${\bf p}'={\bf p}-D_{(x,y)}(\leftvphi-\rightvphi)$ in $\oD_{2\eps_{\rm eq}}$.
Let $N_4$ be the constant from  Definition \ref{definition-10-3}(iv-4).
In $\D_{2\eps_{\rm eq}}=\nD_{2\eps_{\rm eq}}\cup \oD_{2\eps_{\rm eq}}$,
define $O_j^{\rm mod}({\bf p}, x, y)$ by
\begin{equation}
\label{def-O-mod}
O_j^{\rm mod}(p_1, p_2, x, y)
=O_j(x^{3/4}\zeta_1(\frac{p_1'}{x^{3/4}}),\, (\gam+1)N_4 x\zeta_1(\frac{p_2'}{(\gam+1)N_4 x}), \,\psi(x,y), \, x, \, \cbeta)
\end{equation}
for $j=1,\cdots, 5$, where each $O_j({\bf p}, z, x)$ is given by \eqref{prelim5-5}.
With $O_j^{\rm mod}=O_j^{\rm mod}(\hat{\phi}_x, \hat{\phi}_y, x, y)$ for $j=1,\cdots, 5$,
define a nonlinear differential operator $\mcl{N}_{(u,\beta)}^{{\rm polar}}$ by
\begin{equation}
\label{definition-nlop-N2}
\begin{split}
\mcl{N}_{(u,\beta)}^{{\rm polar}}(\hat{\phi}):=&
\Big(2x-(\gam+1)x\zeta_1(\frac{\hat{\psi}_x}{x})+O_1^{\rm mod}\Big)\hat{\psi}_{xx}+
O_2^{\rm mod}\hat{\psi}_{xy}+\Big(\frac{1}{\cbeta}+O_3^{\rm mod}\Big)\hat{\psi}_{yy}\\
&-(1+O_4^{\rm mod})\hat{\psi}_x+O_5^{\rm mod}\hat{\psi}_y\\
=:&\,\, a_{11}(D_{(x,y)}\hat{\phi}, x, y)\hat{\psi}_{xx}+2 a_{12}(D_{(x,y)}\hat{\phi}, x, y)\hat{\psi}_{xy}+{a}_{22}(D_{(x,y)}\hat{\phi}, x, y)\hat{\psi}_{yy}\\
&+{a}_1(D_{(x,y)}\hat{\phi}, x, y)\hat{\psi}_x+{a}_2(D_{(x,y)}\hat{\phi}, x, y)\hat{\psi}_y.
\end{split}
\end{equation}

\smallskip
{\textbf{3.}}
For a $C^2$--function $\hat{\phi}=\hat{\psi}+\phi_{\beta}$,  the expression of $\cbeta\mcl{N}_{(u,\beta)}^{\rm{polar}}(\hat{\phi})$
in the $\xxi$--coordinates is given in the form:
\begin{equation}
\label{definition-A-nr-sonic}
\cbeta\mcl{N}_{(u,\beta)}^{\rm{polar}}(\hat{\phi})=\sum_{i,j=1}^2A_{ij}^{(2)}(D_{\bmxi}\hat{\phi},\bmxi)\der_{\xi_i\xi_j}\hat{\phi}
+\sum_{i=1}^2 A_{i}^{(2)}(D_{\bmxi}\hat{\phi},\bmxi)\der_{\xi_i}\hat{\phi}
\qquad\tx{in $\Om\cap \D_{2\eps_{\rm eq}}$},
\end{equation}
where we have used that $D^2_{\bmxi}\hat{\psi}\equiv D^2_{\bmxi}\hat{\phi}$ holds in $\Om\cap \D_{2\eps_{\rm eq}}$.
In the expression above, $c_{\beta}$ is multiplied to $\mcl{N}_{(u,\beta)}^{\rm{polar}}$ because
the expression of $\cbeta \mcl{N}_{(u,\beta)}^{\rm{polar}}$ without cutoffs in the $\xxi$--coordinates coincides with the left-hand side of Eq. \eqref{2-3}.

In $\Om\cap \oD_{2\eps_{\rm eq}}$, a direct computation shows that
\begin{equation*}
\begin{split}
&A^{(2)}_1=\big((\leftc-x)O_5^{\rm mod}-O_2^{\rm mod}\big)\sin y+
\Big((\leftc-x)\big(\frac{1}{\leftc}+O_3^{\rm mod}\big)-(1+O_4^{\rm mod})\Big)\cos y,\\
&A^{(2)}_2=\big((\leftc-x)O_5^{\rm mod}-O_2^{\rm mod}\big)\cos y-
\Big((\leftc-x)\big(\frac{1}{\leftc}+O_3^{\rm mod}\big)-(1+O_4^{\rm mod})\Big)\sin y.
\end{split}
\end{equation*}
From this, combined with \eqref{prelim5-5} and \eqref{def-O-mod}, we see that
$A^{(2)}_1=A^{(2)}_2=0$ in $\Om\cap \oD_{2\eps_{\rm eq}}.$
Similarly, it can be checked that
$A^{(2)}_1=A^{(2)}_2=0$ in $\Om\cap \nD_{2\eps_{\rm eq}}$. Therefore, we have
\begin{equation*}
A^{(2)}_1=A^{(2)}_2=0\qquad\tx{in $\Om\cap \D_{2\eps_{\rm eq}}$}.
\end{equation*}
For ${\bm\xi}\in\Om\cap \nD_{2\eps_{\rm eq}}$, define
$A_{ij}^{\mcl{N}}$ as
\begin{equation}
\label{definition-coeff-itereqn-rsonic}
A_{ij}^{\mcl{N}}({\bf p}, {\bm \xi}):=A_{ij}^{(2)}({\bf p},{\bm \xi}).
\end{equation}
For ${\bm\xi}\in\Om\cap \oD_{2\eps_{\rm eq}}$, define
$A_{ij}^{\mcl{O}}$ as
\begin{equation}
\label{definition-coeff-itereqn-lsonic}
A_{ij}^{\mcl{O}}({\bf p}, {\bm \xi}):=A_{ij}^{(2)}({\bf p},{\bm \xi}).
\end{equation}
By using Definition \ref{definition-10-3}, the next two lemmas can be  directly derived.
We first discuss the properties of coefficients $(a_{ij}, a_i)$ near $\rightsonic$.

\begin{lemma}[Coefficients $(a_{ij}, a_i)({\bf p},x,y)$ in $\Om\cap \nD_{2\eps_{\rm eq}}$]
\label{lemma-8-3}
There exist constants $\lambda_1\in(0,1), \eps_{\rm eq}\in (0, \frac{\eps_0}{2})$, and $N_{\rm eq}\ge 1$ depending only
on $(\iv, \gam, \beta_*)$ such that,
for any $(u,\beta)\in\ol{\Kext}\cap\{0\le \beta< \betadet\}$,
coefficients $(a_{ij}, a_i)({\bf p},x,y)$ defined by \eqref{definition-nlop-N2}
satisfy the following properties{\rm :}

\smallskip
\begin{itemize}
\item[(a)] For any $(x,y)\in \Om\cap \nD_{2 \eps_{\rm eq}}$ and ${\bf p}, \bm{\kappa}=(\kappa_1,\kappa_2)\in \R^2$,
    \begin{equation*}
    \lambda_1|\bm\kappa|^2\le \sum_{i,j=1}^2 a_{ij}(\bm p, x, y)\frac{\kappa_i\kappa_j}{x^{2-\frac{i+j}{2}}}\le \lambda_1^{-1}|\bm\kappa|^2.
    \end{equation*}

\item[(b)] ${a}_{ij}, {a}_i\in C^{1,\alp}(\R^2\times (\ol{\Om\cap \nD_{\eps_{\rm eq}}}\setminus \ol{\rightsonic}))$ for $j=1,2$, and
\begin{equation*}
\begin{split}
&\|({a}_{11}, {a}_{12}, a_2)\|_{C^{0,1}(\R^2\times \ol{\Om\cap \nD_{\eps_{\rm eq}}})}\le N_{\rm eq},\\[1mm]
&\|({a}_{22},  a_1)\|_{L^{\infty}(\R^2\times \ol{\Om\cap \nD_{\eps_{\rm eq}}})}+
\|D_{({\bf p}, y)}(a_{22},  a_1)\|_{L^{\infty}(\R^2\times \ol{\Om\cap \nD_{\eps_{\rm eq}}})}\le N_{\rm eq},\\[1mm]
&\sup_{({\bf p},x,y)\in \R^2\times \ol{\Om\cap \nD_{\eps_{\rm eq}}}} |x^{1/4}D_x(a_{22}, a_2)({\bf p}, x, y)|\le N_{\rm eq},\\[1mm]
&\sup_{{\bf p}\in \R^2}\|(a_{ij}, {a}_i)({\bf p},\cdot, \cdot)\|_{C^{3/4}(\ol{\Om\cap \nD_{\eps_{\rm eq}}})}\le N_{\rm eq}\qquad\tx{for $i,j=1,2$}.
\end{split}
\end{equation*}

\item[(c)] For each $k=1,2$, $D_{\bf p}^k(a_{ij}, {a}_i)\in C^{1,\alp}(\R^2\times (\ol{\Om\cap \nD_{\eps_{\rm eq}}}\setminus \ol{\rightsonic}))$
and
\begin{equation*}
\qquad\sup_{{\bf p}\in \R^2} \|D_{\bf p}^k(a_{ij}, a_i)({\bf p}, \cdot, \cdot)\|_{C^{1,\alp}(\R^2\times (\ol{\Om\cap \nD_{\eps_{\rm eq}}}\setminus \mcl{N}_r(\rightsonic)))}
\le N_{\rm eq}r^{-5}\quad \mbox{for each $r\in (0, \frac{\eps_{\rm eq}}{2})$}.
\end{equation*}

\item[(d)] There exists a constant $\hat{C}>0$ depending only on $(\iv, \gam,  \beta_*)$
such that
\begin{equation*}
|\der_y(a_{11}, a_{12})({\bf p}, x, y)|\le \hat C x^{1/2}
\qquad\,\tx{for all ${\bf p}\in \R^2$ and $(x,y)\in \Om \cap \nD_{\eps_{\rm eq}}$}.
\end{equation*}

\item[(e)] For every $({\bf p}, x, y)\in \R^2\times \ol{\Om\cap \nD_{\eps_{\rm eq}}}$,
\begin{equation*}
\begin{split}
&({a}_{11}, a_{22}, a_2)((p_1, -p_2), x, y)=({a}_{11},  a_{22}, a_2)((p_1, p_2), x, y),\\
&|a_{ii}({\bf p},x,y)-a_{ii}({\bf 0}, 0, y)|\le N_{\rm eq} x^{3/4}\qquad\tx{for $i=1,2$},\\
&|{a}_{12}({\bf p}, x, y)|\le N_{\rm eq}x,\\
&{a}_1({\bf p}, x, y)\le -\frac 12.
\end{split}
\end{equation*}

\item[(f)] For any ${\bf p}\in \R^2$, the values of $({a}_{ij}, {a}_i)({\bf p}, \cdot, \cdot)$ are given
on $\rightsonic=\{x=0\}\cap \der(\Om\cap \nD_{\eps_{\rm eq}})$
by fixing ${\bf p}$ and taking a limit in $(x,y)$ from $\Om\cap \D_{\eps_{\rm eq}}\subset \{x>0\}$.
More explicitly, for any ${\bf p}\in \R^2$ and $(0, y)\in \rightsonic$,
\begin{equation*}
\begin{split}
&{a}_{ij}({\bf p}, 0, y)=0\qquad\tx{for all}\,\,\,\, (i,j)\neq (2,2),\\
&{a}_{22}({\bf p}, 0, y)=\rightc^{-1},\quad {a}_1({\bf p}, 0, y)=-1,\quad {a}_2({\bf p}, 0, y)=0.
\end{split}
\end{equation*}

\item[(g)]
$\phi=\psi+\phi_{\beta}$ satisfies
\begin{equation*}
\quad O_j^{\rm mod}(\phi_x, \phi_y,x,y)=O_j(\psi_x, \psi_y, \psi, x, y,\cbeta)\qquad\tx{in $\Om\cap \nD_{\eps_{\rm eq}}$\, for $j=1,\cdots, 5$}.
\end{equation*}
In addition, if $\psi$ satisfies
\begin{equation*}
|\psi_x|\le \frac{2-\frac{\mu_0}{5}}{1+\gam}x\qquad\tx{in $\Om\cap \nD_{\eps/10}$}
\end{equation*}
for $\eps\in(0, \frac{\eps_{\rm eq}}{2}]$ from Definition {\rm \ref{definition-10-3}(iv)},
then, in $\Om\cap \nD_{\eps_{\rm eq}}$,
\begin{equation*}
\qquad\,\,\mcl{N}_{(u,\beta)}^{\rm{polar}}(\phi)
=(2x-(\gam+1)\psi_x+O_1){\psi}_{xx}+
O_2{\psi}_{xy}+(\frac{1}{\cbeta}+O_3){\psi}_{yy}
-(1+O_4){\psi}_x+O_5{\psi}_y
\end{equation*}
for $O_j=O_j(\psi_x, \psi_y, \psi, x, y,\rightc)$.
Therefore, equation $\mcl{N}_{(u,\beta)}^{\rm{polar}}(\phi)=0$ coincides with Eq. \eqref{2-3} in $\Om\cap \nD_{\eps_{\rm eq}}$.
\end{itemize}
\end{lemma}

Let $\sigma_3$ be from Proposition \ref{proposition-sub9}.
Coefficients $A_{ij}^{\mcl{O}}, i,j=1,2$,
are used only for $(u,\beta)\in\ol{\Kext}\cap\{\beta\,:\,\beta\in[0,\betasonic+\sigma_3]\}$ to define $\mcl{N}_{(u,\beta)}$.

In the next lemma, we discuss the properties of coefficients $(a_{ij}, a_i)$ near $\leftsonic$
for $\beta\le \betasonic+\sigma_3$. While $\rightsonic$ is fixed to be the same for all $\beta\in[0, \frac{\pi}{2})$,
$\leftsonic$ changes as $\beta$ varies. As $\beta\in[0,\betasonic)$ tends to $\betasonic$, $\leftsonic$ shrinks
to a point set $\{\lefttop\}$ for $\lefttop$ given in Definition \ref{definition-domains-np},
and it remains to be the point set $\{\lefttop\}$ for $\beta>\betasonic$.
For that reason, the properties of $(a_{ij}, a_i)$ near $\leftsonic$ are different from Lemma \ref{lemma-8-3}.

\begin{lemma}[Coefficients $(a_{ij}, a_i)({\bf p},x,y)$ in $\Om\cap \oD_{2\eps_{\rm eq}}$]
\label{lemma1-coeff-iter-eqn-new}
For each $(u,\beta)\in \ol{\Kext}\cap\{\beta\,:\,\beta\in[0,\betasonic+\sigma_3]\}$,
let $(a_{ij}, a_i)$
be
defined by \eqref{definition-nlop-N2}.
Then there exists a constant $\eps_{\rm eq}\in(0, \frac{\eps_0}{2})$
depending only on $(\iv, \gam, \beta_*)$ satisfying the following properties{\rm :}

\smallskip
\begin{itemize}
\item[(a)]
There exist constants $\lambda_1\in(0,1)$ and $N_{\rm eq}\ge 1$ depending only
on $(\iv, \gam, \beta_*)$ such that,
for each $(u,\beta)\in\ol{\Kext}$ with $\beta\in[0, \betac^{(\iv)}+\sigma_3]$,
coefficients $(a_{ij}, a_i)$
satisfy all the assertions
of Lemma {\rm \ref{lemma-8-3}} except for assertions {\rm (d)} and {\rm (g)} of {Lemma {\rm \ref{lemma-8-3}}}
by replacing $(\nD_{\eps_{\rm eq}}, \rightsonic)$ with $(\oD_{\eps_{\rm eq}}, \leftsonic)$.

\smallskip
\item[(b)] Assertion {\rm (d)} of Lemma {\rm \ref{lemma-8-3}}
now takes the following form{\rm :}

\smallskip
\begin{itemize}
\item[(b-1)]
There exists a constant $\hat{C}>0$ depending only on $(\iv, \gam,  \beta_*,\alp)$ such that,
for each $(u,\beta)\in\ol{\Kext}$ with $\beta\in[0, \betasonic)$,
\begin{equation*}
|D_y( a_{11}, {a}_{12})({\bf p},x, y)|\le \hat C x^{1/2}
\qquad\tx{for}\,\,({\bf p}, x, y)\in \R^2\times (\Om\cap \oD_r),
\end{equation*}
where $r=\min\{\gshock^2(-1), \eps_{\rm eq}\}${\rm ;}

\smallskip
\item[(b-2)] Let $\sigma_1>0$ be from Proposition {\rm \ref{proposition-sub8}}.
For any $\delta\in(0, \frac{\sigma_1}{2})$,
there exists a constant $\hat C_{\delta}>0$ depending on $(\iv, \gam,  \beta_*, \delta)$
such that, for each $(u,\beta)\in\ol{\Kext}\cap\{\beta\in(0, \beta_{\rm s}^{(\iv)}-\delta]\}$,
\begin{equation*}
|D_y( a_{11}, {a}_{12})({\bf p},x, y)|\le \hat C_{\delta} x^{1/2}
\qquad\tx{for}\,\,({\bf p}, x, y)\in \R^2\times (\Om\cap \oD_{\eps_{\rm eq}}).
\end{equation*}
\end{itemize}

\item[(c)] Assertion {\rm (g)} of Lemma {\rm \ref{lemma-8-3}}
now takes the following form{\rm :}
suppose that $\psi$ satisfies
\begin{equation}
\label{property-g}
|\psi_x|\le C' x,\quad |\psi_y|\le C' x^{3/2}\qquad \tx{in}\,\,\,\, \Om\cap \oD_{\eps_{\rm eq}}
\end{equation}
for some constant $C'>0$; then there exists a small constant $\eps^{(1)}\in(0, \frac{\eps_{\rm eq}}{2})$ depending
on $(\iv, \gam,  C')$ such that, whenever $\eps$ from Definition {\rm \ref{definition-10-3}(iv)}
with $\eps\le \eps^{(1)}$, $\phi=\psi+\phi_{\beta}$ satisfies
\begin{equation*}
\qquad O_j^{\rm mod}(\phi_x, \phi_y,x,y)=O_j(\psi_x, \psi_y, \psi, x, y,\cbeta)\qquad\tx{in $\Om\cap \oD_{\eps_{\rm eq}}\,$ for $j=1,\cdots, 5$}.
\end{equation*}

\begin{itemize}
\item[(c-1)]
For $P_{\beta}$ given by \eqref{def-Pbeta}, suppose that
\begin{equation*}
x_{P_{\beta}}<\frac{\eps}{10},\quad {i.e.}, \,\,\,\, \Om\cap\oD_{\eps/10}\neq \emptyset.
\end{equation*}
If $\psi$ satisfies
\begin{equation*}
|\psi_x|\le \frac{2-\frac{\mu_0}{5}}{1+\gam}x\qquad\tx{in}\,\, \Om\cap \oD_{\eps/10},
\end{equation*}
then, in $\Om\cap \oD_{\eps_{\rm eq}}$,
\begin{equation*}
\qquad\,\,\, \mcl{N}_{(u,\beta)}^{\rm{polar}}(\phi)
=(2x-(\gam+1)\psi_x+O_1){\psi}_{xx}+
O_2{\psi}_{xy}+(\frac{1}{\cbeta}+O_3){\psi}_{yy}
-(1+O_4){\psi}_x+O_5{\psi}_y
\end{equation*}
for $O_j=O_j(\psi_x, \psi_y, \psi, x, y,\cbeta)$.
Therefore, if $\mcl{N}_{(u,\beta)}^{\rm{polar}}(\phi)=0$ holds in $\Om\cap \oD_{\eps_{\rm eq}}$,
then $\vphi$ satisfies Eq. \eqref{2-3} in $\Om\cap \oD_{\eps_{\rm eq}}$.

\smallskip
\item[(c-2)]
For $\beta\in(\betac^{(\iv)}, \betac^{(\iv)}+\sigma_3]$, suppose that
\begin{equation*}
  x_{P_{\beta}}\ge \frac{\eps}{10},
\end{equation*}
which is equivalent to the case that $\Om\cap \oD_{\eps/10}=\emptyset$.
Then equation $\mcl{N}_{(u,\beta)}^{\rm{polar}}(\phi)=0$ coincides with Eq. \eqref{2-3} in $\Om\cap \oD_{\eps_{\rm eq}}$.
\end{itemize}

\smallskip
\item[(d)] For all $(u,\beta)\in\ol{\Kext}$ with $\beta>\betac^{(\iv)}$, $({a}_{ij}, {a}_i)({\bf p},\cdot, \cdot)$
and $D_{\bf p}^k({a}_{ij}, {a}_i)({\bf p},\cdot, \cdot)$, $k=1,2$,
 are in $C^{1,\alp}(\ol{\Om\cap \oD_{\eps_{\rm eq}}})$.
In particular, for each $\delta\in(0, \frac{\sigma_3}{2})$, there exists a constant $C_{\delta}>0$ depending
only on $(\iv, \gam, \beta_*,  \delta)$ such that,
if $(u,\beta)\in\ol{\Kext}$ with $\beta\in [\betasonic+\delta, \betasonic+ \frac{\sigma_3}{2})$, then
\begin{equation*}
\begin{split}
&\sup_{{\bf p}\in \R^2}\|({a}_{ij}, a_i)({\bf p}, \cdot,\cdot)\|_{C^{1,\alp}(\ol{\Om\cap \oD_{\eps_{\rm eq}}})}\le C_{\delta},\\
&\sup_{{\bf p}\in \R^2} \|D^k_{\bf p}(a_{ij}, a_i)({\bf p}, \cdot,\cdot)\|_{C^{1,\alp}(\ol{\Om\cap \oD_{\eps_{\rm eq}}})}\le C_{\delta}
\qquad\tx{for $k=1,2$}.
\end{split}
\end{equation*}
\end{itemize}
\end{lemma}
\smallskip
{\textbf{4.}}
In this step, we define $\mcl{N}_{(u,\beta)}$ near $\leftsonic$ for $(u,\beta)\in \ol{\Kext}$ with $\beta\ge \betasonic+ \frac{\sigma_3}{4}$.
\begin{lemma}
\label{lemma-c4-approximation-new}
For each $(u,\beta)\in\ol{\Kext}$,
let $\gshock=\gshock^{(u,\beta)}$, $\mathfrak{F}=\mathfrak{F}_{(u,\beta)}$, $\vphi=\vphi^{(u,\beta)}$, and
$\Om=\Om(u,\beta)$ be defined by Definition {\rm \ref{definition-Gset-shocks-new}}, and let
\begin{equation}
\label{definition-phi-ub}
\phi:=\vphi^{(u,\beta)}-\rightvphi
\end{equation}
for $\rightvphi$ given by \eqref{def-uniform-ptnl-new}. For any given $\sigma\in(0,1)$, there exists a constant $C_{\sigma}>0$
depending only on $(\iv, \gam,  \beta_*, \sigma)$ such that,
for each $(u,\beta)\in\ol{\Kext}$, there exists a function $v_{\sigma}^{(u,\beta)}\in C^4(\ol{\Om})$
satisfying the following two properties{\rm :}

\smallskip
\begin{itemize}
\item[(a)]
$\|v_{\sigma}^{(u,\beta)}-\phi\|_{C^1(\ol{\Om})}\le \sigma^2$
and $\|v_{\sigma}^{(u,\beta)}\|_{C^4(\ol{\Om})}\le C_{\sigma}${\rm ;}

\smallskip
\item[(b)] $v_{\sigma}^{(u,\beta)}$ depends continuously on $(u, \beta)\in\ol{\Kext}$
in the sense that, if $\{(u_k, \beta_k)\}\subset \ol{\Kext}$  converges
to $(u,\beta)$ in  $C^{1,\alp}(\ol{\iter})\times[0, \beta_*]$ for some $(u,\beta)\in\ol{\Kext}$, then
\begin{equation*}
v_{\sigma}^{(u_k,\beta_k)}\circ \mathfrak{F}_{(u_k,\beta_k)}\rightarrow v_{\sigma}^{(u,\beta)}\circ \mathfrak{F}_{(u,\beta)}
\qquad\, \tx{in $C^{1,\alp}(\ol{\iter})$}.
\end{equation*}
\end{itemize}

\begin{proof}
For $\mcl{G}_1^{\beta}$ defined by \eqref{12-16-mod}, denote
\begin{equation*}
w(s,t'):=\phi\circ (\mcl{G}_1^{\beta})^{-1}(s,t')
\end{equation*}
for $(s,t')\in \mcl{G}_1^{\beta}(\Om)=\{(s,t')\,:\,-1<s<1, 0<t'<\gshock^{(u,\beta)}(s)\}$.
For each small constant $\eps>0$, define a function $\til{w}_{\eps}(s,t')$ by
\begin{equation*}
\til{w}_{\eps}(s,t'):=w(\frac{s}{1+\frac{\eps}{M_1}}, \frac{t'+\frac{\eps}{2M_2}}{1+\eps})
\end{equation*}
for constants $M_1>1$ and $M_2>1$ to be determined later.
Then $\til w_{\eps}$ is well defined in the set:
\begin{equation*}
\mcl{A}_{\eps}:=\Big\{(s,t')\,:\,|s|<1+\frac{\eps}{M_1},\,\,
-\frac{\eps}{2M_2}<t'<(1+\eps)\gshock(\frac{s}{1+\eps/M_1})-\frac{\eps}{2M_2}\Big\}.
\end{equation*}
Using (i) and (iii) of Definition \ref{definition-10-3}, and Remark \ref{remark-a}(i),
we choose constants $M_1, M_2, M_3>1$
depending only on $(\iv, \gam, \beta_*)$ such that the $\frac{\eps}{M_3}$--neighborhood
$\mcl{N}_{\frac{\eps}{M_3}}(\mcl{G}_1^{\beta}(\Om))$
of $\mcl{G}_1^{\beta}(\Om)$
is contained in $\mcl{A}_{\eps}$.

Define
\begin{equation*}
  w_{\eps}(s,t'):=(\til{w}_{\eps} * \chi_{\frac{\eps}{2M_3}})(s,t')\qquad\,\tx{in $\mcl{G}_1^{\beta}(\Om)$}
\end{equation*}
with $\chi_{\delta}({\bmxi}):=\frac{1}{\delta^2}
\chi(\frac{\bmxi}{\delta})$, where $\chi(\cdot)$ is a standard mollifier:
$\chi\in C_0^{\infty}(\R^2)$ is a nonnegative function with ${\rm supp}(\chi)\subset B_1({\bf 0})$
and $\int_{\R^2} \chi({\bmxi})\,{\rm d}{\bmxi}=1$.
Then we define
\begin{equation*}
V_{\eps}^{(u,\beta)}(\bmxi):=w_{\eps}\circ \mcl{G}_1^{\beta}(\bmxi)\,\qquad\mbox{in $\Om$}.
\end{equation*}
For each $\sigma\in(0,1)$, there exists a small
constant $\eps_*(\sigma)>0$ depending on $(\iv, \gam,  \beta_*, \sigma)$
such that $v_{\sigma}^{(u,\beta)}:=V^{(u,\beta)}_{\eps_*(\sigma)}$ satisfies
properties (a)--(b).
\end{proof}
\end{lemma}

Let $\varsigma\in C^{\infty}(\R)$ be a cut-off function satisfying that
\[
\varsigma(t)=\begin{cases}
1\quad\tx{for}\,\, t<1,\\
0\quad\tx{for}\,\,t>2,
\end{cases}\qquad
0\le \varsigma\le 1\quad\tx{on $\R$}.
\]
For a constant $\sigma>0$, denote
\begin{equation}
\label{definition-varsigma}
\varsigma_{\sigma}(t):=\varsigma(\frac{t}{\sigma}).
\end{equation}

Let $\sigma_{\rm{cf}}\in(0,1)$ be a constant to be specified later.
For each $(u,\beta)\in\ol{\Kext}$, let $v^{(u,\beta)}_{\sigma_{\rm{cf}}}$ be the function given by Lemma \ref{lemma-c4-approximation-new}.
For each $i,j=1,2$, we define
\begin{align}
A_{ij}^{\mcl{O}, {\rm subs}}({\bf p}, {\bm\xi})
=&\,\varsigma_{{\sigma}_{\rm{cf}}}(|{\bf p}
-Dv_{\sigma_{\rm{cf}}}^{(u,\beta)}(\bm\xi)|)A_{ij}^{\rm potn}({\bf p}, \phi(\bm\xi), {\bm\xi})\notag \\
&\,+\big(1- \varsigma_{{\sigma}_{\rm{cf}}}(|{\bf p}
-Dv_{\sigma_{\rm{cf}}}^{(u,\beta)}(\bm\xi)|)\big)A_{ij}^{\rm potn}(Dv_{\sigma_{\rm{cf}}}^{(u,\beta)}(\bm\xi),\phi(\bm\xi),{\bm\xi})
\label{definition-coeff-itereqn-sub-new}
\end{align}
for $A_{ij}^{\rm potn}({\bf p}, z, \bmxi)$ defined by \eqref{def-A-potn-new}.

\begin{lemma}
\label{lemma2-coeff-iter-eqn-new}
There exist two small constants $\eps^{(2)}>0$ and $\delta_1^{(1)}>0$ depending only on $(\iv, \gam)$ such that,
whenever $\eps$ and $\delta_1$ from Definition {\rm \ref{definition-10-3}} satisfy
\begin{equation*}
\eps\le \eps^{(2)},\quad \delta_1\le \delta_1^{(1)},
\end{equation*}
there exist $C>0$ depending only on $(\iv, \gam,  \beta_*)$
and $\lambda\in(0,1)$ depending only on $(\iv, \gam)$ so that,
for each $(u,\beta)\in\ol{\Kext}\cap\{\beta\ge \betasonic+\frac{\sigma_3}{4}\}$,
the associated coefficients $A_{ij}^{\mcl{O},{\rm subs}}$ defined
by \eqref{definition-coeff-itereqn-sub-new}
with $\sigma_{\rm{cf}}=\sqrt{\delta_1}$ satisfy the following properties{\rm :}

\smallskip
\begin{itemize}
\item[(a)] For all $({\bf p}, \bmxi)\in\R^2\times \ol{\Om\cap \oD_{\eps_{\rm eq}}}$
satisfying that $|{\bf p}-D\phi(\bmxi)|<\frac{\sqrt{\delta_1}}{2}$,
\begin{equation*}
A_{ij}^{\mcl{O}, {\rm subs}}({\bf p}, \bmxi)=A_{ij}^{\rm potn}({\bf p}, \phi(\bmxi), \bmxi),
\end{equation*}
so that
\begin{equation*}
A_{ij}^{\mcl{O}, {\rm subs}}(D\phi(\bmxi), \bmxi)=A_{ij}^{\rm potn}(D\phi(\bmxi), \phi(\bmxi), \bmxi)\qquad\tx{in $\Om${\rm ;}}
\end{equation*}

\item[(b)] For all $({\bf p}, \bmxi)\in\R^2\times \ol{\Om\cap \oD_{\eps_{\rm eq}}}$,
\begin{equation*}
|A_{ij}^{\mcl{O}, {\rm subs}}({\bf p}, \bmxi)-A_{ij}^{\mcl{O}, {\rm subs}}(D\phi(\bmxi), \bmxi)| \le C\sqrt{\delta_1};
\end{equation*}

\item[(c)] For each ${\bf p}\in \R^2$, $D^k_{\bf p}A_{ij}^{\mcl{O}, {\rm subs}}({\bf p},\cdot)$
are in $C^{1,\alp}(\ol{\Om\cap \oD_{\eps_{\rm eq}}})$ for $k=0,1,2$, with
\begin{equation*}
\sum_{k=0}^2\|D^k_{\bf p}A_{ij}^{\mcl{O}, {\rm subs}}({\bf p},\cdot)\|_{C^{1,\alp}(\ol{\Om\cap \oD_{\eps_{\rm eq}}})} \le C;
\end{equation*}

\item[(d)] For all $\bmxi\in \Om\cap \oD_{\eps_{\rm eq}}$ and ${\bf p}$, ${\bm\kappa}=(\kappa_1, \kappa_2)\in \R^2$,
\begin{equation*}
\lambda|{\bm\kappa}|^2 \le \sum_{i,j=1}^2 A_{ij}^{\mcl{O}, {\rm subs}}({\bf p}, \bmxi) \kappa_i \kappa_j \le \lambda^{-1} |{\bm\kappa}|^2.
\end{equation*}
\end{itemize}
\end{lemma}

{\textbf{5.}}
Let $\chi_{\rm eq}\in C^{\infty}(\R)$ be a function satisfying that
\begin{equation*}
\chi_{\rm eq}(\beta)=\begin{cases}
1&\quad
\tx{if}\,\,\beta\le \betasonic+\frac{\sigma_3}{4},
\\[1mm]
0&\quad
\tx{if}\,\,\beta \ge \betasonic+\frac{\sigma_3}{2},
\end{cases} \qquad\,\,
\chi_{\rm eq}'(\beta)\le 0\,\,\,\,\tx{on $\R$}.
\end{equation*}
For such a cut-off function $\chi_{\rm eq}$, we define
\begin{equation}
\label{definition-coeff-itereqn-near-pt-new}
A_{ij}^{(3)}({\bf p}, {\bm\xi})
=
\begin{cases}
 \chi_{\rm eq}(\beta) A_{ij}^{\mcl{O}}({\bf p}, {\bm\xi})+(1-\chi_{\rm eq}(\beta))A_{ij}^{\mcl{O}, {\rm subs}}({\bf p}, {\bm\xi})=:A_{ij}^{(3,\mcl{O})}({\bf p}, \bmxi)
 &\,\tx{for $\xin<0$},\\[2mm]
A_{ij}^{\mcl{N}}({\bf p}, {\bm\xi})&\,\tx{for $\xin>0$}
 \end{cases}
\end{equation}
for $A_{ij}^{\mathcal{N}}$ and $A_{ij}^{\mathcal{O}}$ given by
\eqref{definition-coeff-itereqn-rsonic} and
\eqref{definition-coeff-itereqn-lsonic}, respectively.

\smallskip
{\textbf{6.}} Finally, we combine \eqref{def-A1} with
\eqref{definition-coeff-itereqn-near-pt-new}
to complete the definition of $\mcl{N}_{(u,\beta)}(\hat{\phi})$ in \eqref{definition-nub-hphi}.

\begin{definition}\label{definition-cutoff-2016-1}
We define the following:
\begin{enumerate}
\item[\rm (i)] For a parameter $\tau\in(0, \frac 12]$, introduce a family of functions $\bar\zeta_2(s,t;\tau)$ so that

\smallskip
\begin{itemize}
\item $\bar\zeta_2(\cdot,\cdot;\tau)\in C^4(\R^2)$ for each $\tau\in(0, \frac 12]${\rm ;}

\smallskip
\item $\der_{t}\bar\zeta_2(s,t;\tau)= 0$ for each $\tau\in(0, \frac 12]$ and $(s,t)\in \R^2${\rm ;}

\smallskip
\item For each $\tau\in(0, \frac 12]$,
$\bar\zeta_2(s,t;\tau)=\begin{cases}
1&\quad \tx{for $|s|<1-\tau$},\\
0&\quad \tx{for $|s|\ge 1-\frac{\tau}{2}$};
\end{cases}
$

\smallskip
\item  $\bar\zeta_2(-s,t;\tau)=\bar\zeta(s,t;\tau)$ for all $s\in\R$ and $\tau\in(0, \frac 12]${\rm ;}

\smallskip
\item $-\frac{10}{\tau}\le \der_s\bar\zeta_2(s,t;\tau)\le 0$  for all $s\ge 0$ and $\tau\in(0, \frac 12]${\rm ;}

\smallskip
\item  $\|\bar\zeta_2(\cdot,\cdot; \tau)\|_{C^4(\R^2)}$ is a continuous function of $\tau\in(0, \frac 12]$.
\end{itemize}

\item[\rm (ii)] For $\beta_*\in(0,\betadet)$, define a set $Q_{\beta_*}^{\cup} \subset \R^2_+\times [0, \frac{\pi}{2})$ as
\begin{equation*}
Q_{\beta_*}^{\cup}:=\cup_{\beta\in[0, \beta_*]} \Qbeta\times \{\beta\}
\end{equation*}
for $\Qbeta$ defined by Definition {\rm \ref{definition-Qbeta}(iii)}.

For $\eps>0$ and $\beta\in[0, \beta_*]$, let $\hat{\eps}$ be given by \eqref{def-heps-mod}.
For $(\bmxi,\beta)\in Q_{\beta_*}^{\cup}$, define a function $\zeta_2^{(\eps,\beta)}: Q_{\beta_*}^{\cup}\rightarrow \R$ by
\begin{equation}
\label{definition-zeta2-subs-new}
{\zeta}_2^{(\eps,\beta)}(\bmxi):=\bar\zeta_2(\mcl{G}_1^{\beta}(\bmxi); \hat{\eps}).
\end{equation}
\end{enumerate}
\end{definition}

The $C^1$--dependence of $(\sbeta, \cbeta, \leftu)$ on $\beta\in[0, \frac{\pi}{2})$ yields the following lemma:

\begin{lemma}
\label{lemma-properties-zeta2}
Let $\eps_0>0$ be from Lemma {\rm \ref{lemma-7-4}(c)}. For each $\eps\in(0, \frac{\eps_0}{2})$, ${\zeta}_2^{(\eps,\beta)}$
satisfies the following properties{\rm :}

\smallskip
\begin{itemize}
\item[(a)] ${\zeta}_2^{(\eps,\beta)}: Q_{\beta_*}^{\cup}\rightarrow \R$ is $C^4$ with respect
to $\bmxi\in \Qbeta$ for $\beta\in[0,\beta_*]$, and is continuous with respect to $\beta\in[0, \beta_*]${\rm ;}

\smallskip
\item[(b)] There exists a constant $C_{\eps}>0$ depending only on $(\iv, \gam,  \eps)$ such that
\begin{equation*}
\|\zeta_2^{(\eps,\beta)}\|_{C^4(\ol{\Qbeta})}\le C_{\eps};
\end{equation*}

\item[(c)]
$
\zeta_2^{(\eps,\beta)}=
\begin{cases}
1&\quad\tx{in $\Om{(u,\beta)}\setminus \D_{\eps}$},\\
0&\quad \tx{in $\Om{(u,\beta)}\cap \D_{{\eps}/{2}}$}.
\end{cases}
$
\end{itemize}
\end{lemma}
Finally, we define coefficients $A_{ij}({\bf p}, {\bm \xi})$ for the nonlinear differential operator $\mcl{N}_{(u,\beta)}$ given by \eqref{definition-nub-hphi} as follows:
\begin{equation}
\label{12-57}
A_{ij}({\bf p}, {\bm \xi}):=\zeta_2^{(\eps_{\rm eq},\beta)}({\bm\xi})A_{ij}^{(1)}(\bm\xi)
+\big(1-\zeta_2^{(\eps_{\rm eq},\beta)}({\bm\xi})\big)A_{ij}^{(3)}({\bf p}, {\bm\xi}), \quad i,j=1,2.
\end{equation}
Hereafter, we continue to adjust $\eps_{\rm eq}>0$ depending only on $(\iv, \gam)$.

\begin{lemma}
\label{lemma-6-1}
For each $(u,\beta)\in\ol{\Kext}$, let coefficients $A_{ij}({\bf p},\bmxi), i,j=1,2$, of $\mcl{N}_{(u,\beta)}$
in \eqref{definition-nub-hphi} be given by \eqref{12-57}.
Then there exist constants $\eps_{\rm eq}\in(0, \frac{\eps_0}{2})$, $\lambda_0\in(0,1)$, $N_{\rm eq}\ge 1$,
and $C>0$
with $\lambda_0$ depending only on $(\iv, \gam)$, $(N_{\rm eq}, \eps_{\rm eq})$ depending on
$(\iv, \gam,  \beta_*)$,
and  $C>0$ depending only on $(\iv, \gam,  \beta_*,\alp)$
such that the following properties hold{\rm :}

\smallskip
\begin{itemize}
\item[(a)] For all $\bmxi\in \Om$ with $\Om=\Om(u,\beta)$
and all ${\bf p}, \bm{\kappa}=(\kappa_1,\kappa_2)\in \R^2$,
\begin{equation*}
\lambda_0 \;{\rm dist}(\bmxi,\leftsonic\cup \rightsonic)|\bm{\kappa}|^2\le \sum_{i,j=1}^2A_{ij}(\bm p, \bmxi)\kappa_i\kappa_j\le \lambda_0^{-1}
|\bm{\kappa}|^2;
\end{equation*}

\item[(b)] $A_{12}({\bf p},\bmxi)=A_{12}({\bf p}, \bmxi)$ holds in $\R^2\times \Om$, and each $A_{ij}$ satisfies
\begin{equation*}
\|A_{ij}\|_{L^{\infty}(\R^2\times \Om)}\le N_{\rm eq};
\end{equation*}

\item[(c)] For $\bmxi=(\xin,\etan)\in\Om\setminus \D_{\eps_{\rm eq}}$,  $A_{ij}({\bf p}, \bmxi)=A_{ij}^{(1)}(\bmxi)$ and
\begin{equation*}
    \|A_{ij}\|_{C^{1,\alp}(\ol{\Om\setminus \D_{\eps_{\rm eq}})}}\le C;
    \end{equation*}

\item[(d)] For each ${\bf p}\in\R^2$,
\begin{equation*}
\|A_{ij}({\bf p},\cdot,\cdot)\|_{C^{3/4}(\ol{\Om})}
+\|D_{{\bf p}}A_{ij}({\bf p}, \cdot, \cdot)\|_{L^{\infty}(\Om)}\le N_{\rm eq};
\end{equation*}

\item[(e)] For each $k=0,1,2$,
$D_{\bf p}^kA_{ij}\in C^{1,\alp}(\R^2\times (\ol{\Om}\setminus \ol{\leftsonic}\cup\ol{\rightsonic})$.
Furthermore, for each $s\in (0, \frac{\eps_0}{2})$, $D_{\bf p}^kA_{ij}$ satisfies
\begin{equation*}
\|D_{\bf p}^kA_{ij}\|_{C^{1,\alp}(\R^2\times (\ol{\Om}\setminus \mcl{N}_s(\ol{\leftsonic}\cup\ol{\rightsonic})))}\le C s^{-5};
\end{equation*}

\item[(f)] For each $i,j=1,2$, $A_{ij}({\bf p}, \bmxi)=A_{ij}^{\mcl{N}}({\bf p}, \bmxi)$
holds for all $({\bf p}, \bmxi)\in\R^2\times (\Om\cap \nD_{\eps_{\rm eq}/2})${\rm ;}

\smallskip
\item[(g)]
If $\beta\le \betasonic+\frac{\sigma_3}{4}$,
then $A_{ij}({\bf p}, \bmxi)=A_{ij}^{\mcl{O}}({\bf p}, \bmxi)$ holds for all
$({\bf p},\bmxi)\in \R^2\times \oD_{\eps_{\rm eq}/2}${\rm ;}

\smallskip
\item[(h)] If $\beta\in[\betasonic+\delta, \beta_*]$
for $\delta\in(0, \frac{\sigma_3}{2})$,
then $A_{ij}({\bf p}, \bmxi)=A_{ij}^{(3)}({\bf p}, \bmxi)$ holds for all $({\bf p}, \bmxi)\in\R^2\times (\Om \cap \oD_{\eps_{\rm eq}/2})$, and
\begin{equation*}
\begin{split}
&\qquad\,\,\lambda_0\big({\rm dist}(\bmxi, \ol{\leftsonic})+\delta\big)|\bm{\kappa}|^2
 \le \sum_{i,j=1}^2 A_{ij}({\bf p}, \bmxi)\kappa_i\kappa_j\le \lambda_0^{-1}|\bm\kappa|^2
\qquad
\tx{for all}\,\,{\bm\kappa}=(\kappa_1, \kappa_2)\in\R^2,\\
&\qquad\,\,\sup_{{\bf p}\in \R^2}\|D^k_{\bf p} A_{ij}({\bf p}, \cdot, \cdot)\|_{C^{1,\alp}(\ol{\Om\cap \oD_{\eps_{\rm eq}/2}})} \le C\qquad\tx{for}\,\,k=0,1,2;
\end{split}
\end{equation*}

\item[(i)] For each $(u,\beta)\in \ol{\Kext}$, let $\phi=\phi^{(u,\beta)}$ be defined by \eqref{definition-phi-ub}.
Suppose that $\eps$ from Definition {\rm \ref{definition-10-3}} satisfies that $0<\eps<\frac{\eps_{\rm eq}}{2}$.
Then equation $\mcl{N}_{(u,\beta)}(\phi)=0$ coincides with \eqref{2-3} in $\Om\setminus (\oD_{\eps/10}\cup \nD_{\eps/10})$.
In addition, if $x_{P_{\beta}}\ge \frac{\eps}{10}$
or $\beta\ge \betasonic+\frac{\sigma_3}{2}$ holds,
then equation $\mcl{N}_{(u,\beta)}(\phi)=0$ coincides with \eqref{2-3} in $\Om\setminus \nD_{\eps/10}$.
\end{itemize}
\end{lemma}

\subsection{\it Definition of $\M_{(u,\beta)}({\bf p},z,\bmxi)$ in \eqref{12-50}}
\label{subsubsec-BC-shock}
The definition 
is given in the following five steps:

\smallskip
{\textbf{1.}}
For $\rightvphi$ and $g^{\rm sh}$ given by \eqref{def-uniform-ptnl-new} and \eqref{def-gsh}, respectively, define
\begin{equation}
\label{12-67}
\M_0({\bf p},z,\bmxi):
=g^{\rm sh}({\bf p}+D\rightvphi(\bmxi), z+\rightvphi(\bmxi), \bmxi) \qquad \text{for ${\bf p},\bmxi\in \R^2$ and $z\in \R$}.
\end{equation}
The nonlinear function $\M_0({\bm p},z,{\bm \xi})$ is well defined on the set:
\begin{equation*}
\mcl{A}_{\M_0}:=
\left\{
\begin{array}{ll}
({\bm p}, z, {\bm\xi})\in B_{4N_5}(\bm 0)\times (-4N_5,4N_5)\times B_{4M_{\rm dom}}(\bm 0)\\
\,: \,\,\,\, 2\rho_{\rm max}^{\gam-1}>\rightrho^{\gam-1}+(\gam-1)\big({\bm \xi}\cdot{\bf p}-\frac{|\bm p|^2}{2}-z\big)>\frac{\rho_{\rm min}^{\gam-1}}{2},\\
\,\,\,\,\,\,\,\,\,\, |{\bf p}-(0,-\iv)|>\frac{\mu_1}{2}
\end{array}
\right\}
\end{equation*}
for constants $(\mu_1, N_5, \rho_{\rm min}, \rho_{\rm max})$ from properties (iv) and (vi) of Definition \ref{definition-10-3},
and $M_{\rm dom}$ from Remark \ref{remark-a}.
Since these constants are chosen depending only on $(\iv, \gam)$, for each $k=1,2,\cdots$, there exists a constant $C_k>0$ depending only on $(\iv, \gam,  k)$
to satisfy
\begin{equation}
\label{estimate-M0-iter}
\|\M_0\|_{C^k(\ol{\mcl{A}_{\M_0}})}\le C_k.
\end{equation}

\smallskip
{\textbf{2.}}
Similarly to \eqref{7-d7-pre}, we define a function $\M_1({\bf p},z,\xin)$ by
\begin{equation}
\label{7-d7}
\M_1({\bf p}, z, \xin)=\M_0({\bf p},z, \xin,\neta-\frac{z}{\iv}).
\end{equation}
$\M_1$ is well defined in the set:
\begin{equation*}
\mcl{A}_{\M_1}:=\left\{
\begin{array}{ll}
({\bf p}, z, {\bm\xi})\in B_{3 N_5}(\bm 0)\times (-3N_5, 3N_5)\times B_{3M_{\rm dom}}(\bm 0)\\
\,: \,\,\,\,
2\rho_{\rm max}^{\gam-1}>\rightrho^{\gam-1}+(\gam-1)\big(p_1\neta+p_2(\xin-\frac{z}{\iv})
-\frac{|{\bf p}|^2}{2}-z\big)>\frac{\rho_{\rm min}^{\gam-1}}{2},
\\
\,\,\,\,\,\,\,\,\,\, |{\bf p}-(0,-\iv)|>\frac{\mu_1}{2}
\end{array}
\right\}.
\end{equation*}

For each $k=1,2,\cdots$, there exists a constant $C_k>0$ depending only on $(\iv, \gam,  k)$
such that
\begin{equation}
\label{estimate-M1-iter}
\|\M_1\|_{C^k(\ol{\mcl{A}_{\M_1}})}\le C_k.
\end{equation}
In particular, $\M_1$ is homogeneous in the sense of
\begin{equation}
\label{hom-M1}
\M_1({\bf 0}, 0,\xi_1)=0,\quad \M_1(D(\leftvphi-\rightvphi),\leftvphi-\rightvphi,\xi_1)=0
\qquad\,\, \tx{for all $\xi_1\in\R$}.
\end{equation}

\smallskip
{\textbf{3.}}
For $(\leftvphi, \rightvphi)$ given by \eqref{def-uniform-ptnl-new}, denote
\begin{equation}
\label{definition-left-psi}
\leftphi:=\leftvphi-\rightvphi.
\end{equation}
For a constant $\sigma>0$, let function $\varsigma_{\sigma}$ be given by \eqref{definition-varsigma}.
For a constant ${\sigma}_{\rm{bc}}>0$ to be determined later, we define
\begin{equation}
\label{7-d8}
\begin{split}
\M(\bm p, z, \bm \xi)
=&\, \varsigma_{\sigma_{\rm{bc}}}(|({\bf p},z)|)\M_1({\bf p}, z, \xin)\\
&\,+ \bigl(1-\varsigma_{\sigma_{\rm{bc}}}(|({\bf p},z)|)\bigr)
\Bigl(
\varsigma_{\sigma_{\rm{bc}}}(|({\bf p},z)-(D\leftphi,\leftphi(\bm\xi))|)\M_1({\bf p}, z, \xin)\\
&\qquad\quad\qquad\qquad\quad\,\,\,\,\,\,\,\,\,\,\,+\big(1-\varsigma_{\sigma_{\rm{bc}}}(|({\bf p},z)-(D\leftphi,\leftphi(\bm\xi))|)\big) \M_0({\bf p}, z, \bm\xi)
\Bigr)
\end{split}
\end{equation}
for $({\bf p},z,\bm\xi)\in  \mcl{A}_{\M}:=\mcl{A}_{\M_0}\cap \mcl{A}_{\M_1}$.
\smallskip

For each $(u,\beta)\in \ol{\Kext}$, let $\gshock=\gshock^{(u,\beta)}$, $\mathfrak{F}=\mathfrak{F}_{(u,\beta)}$,
$\Om=\Om(u,\beta)$, $\shock=\shock(u,\beta)$,
and $\vphi=\vphi^{(u,\beta)}$ be defined by Definition \ref{definition-Gset-shocks-new}.
Denote $\phi:=\vphi-\rightvphi$.

\smallskip
For a constant $\sigma>0$, we define
\begin{equation*}
\mcl{E}(\phi,\shock)=\{({\bf p}, z, {\bm \xi})\in \R^2\times \R\times \R^2\,:\;
{\bm p}=D\phi(\bm\xi),\; z=\phi(\bm\xi),\;{\bm\xi}\in \shock\}
\end{equation*}
and
\begin{equation*}
\mcl{E}_{\sigma}(\phi,\shock)=\{({\bf p}, z, {\bm \xi})
\in \R^2\times \R\times \R^2:\;
{\rm dist}({\bm \xi},\shock)<\sigma,
|{\bf p}-D\phi({\bm\xi})|<\sigma,
|z-\phi({\bm \xi})|<\sigma
\}.
\end{equation*}

\begin{lemma}
\label{lemma-Mproperty}
There exists a constant $\bar{\sigma}_{\rm{bc}}>0$ depending only on $(\iv, \gam)$ such that,
whenever $\sigma_{\rm{bc}}\in(0,\bar{\sigma}_{\rm{bc}}]$, there exists
a constant $C_{\sigma_{\rm{bc}}}>0$ depending only on $(\iv, \gam,  \sigma_{\rm{bc}})$ so that
\begin{equation*}
\|\M\|_{C^4(\ol{\mcl{A}_{\M}})}\le C_{\sigma_{\rm{bc}}}.
\end{equation*}
Furthermore,
for each $(u,\beta)\in \ol{\Kext}$, the following properties hold{\rm :}

\smallskip
\begin{itemize}
\item[(a)] $\mcl{E}_{\sigma_{\rm{bc}}}(\phi,\shock)\subset \mcl{A}_M${\rm ;}

\smallskip
\item[(b)] The map{\rm :} $\beta\mapsto \M$ is in $C([0,\beta_*];C^4(\ol{\mcl{A}_{\M}}))${\rm ;}
\item[(c)] On $\shock$, ${\M}(D\phi,\phi,{\bm\xi})={\M}_0(D\phi,\phi,{\bm\xi})$ and
$\der_{\bf p}\M(D\phi,\phi,{\bm\xi})=\der_{\bf p}\M_0(D\phi,\phi,{\bm\xi})${\rm ;}

\smallskip
\item[(d)] $\phi$ satisfies
\begin{equation}
\label{bc-M}
\M(D{\phi},{\phi},\bm\xi)=0\qquad\mbox{on $\shock$}
\end{equation}
if and only if $\vphi$ satisfies \eqref{rhbc-gsh}{\rm ;}

\smallskip
\item[(e)] $\M$ is homogeneous
in the sense that
\begin{equation}
\label{hom}
\M({\bf 0}, 0,\bm \xi)=0,\quad \M(D(\leftvphi-\rightvphi),\leftvphi-\rightvphi,\bm\xi)=0
\qquad\,\,\mbox{for all ${\bm \xi}\in B_{2M_{\rm dom}}({\bf 0})$}.
\end{equation}
\end{itemize}
\end{lemma}

\begin{lemma}
\label{lemma-7-7}
For constant $\bar{\sigma}_{\rm{bc}}$ from Lemma {\rm \ref{lemma-Mproperty}},
there exist constants $\sigma_{\rm{bc}}\in(0,\bar{\sigma}_{\rm{bc}}]$, $\bar{\eps}_{\rm{bc}}>0$,
and $\delta_{bc}>0$
depending only on $(\iv,\gam)$
such that, if $\eps$ from Definition {\rm \ref{definition-10-3}} satisfies that $0<\eps\le \bar{\eps}_{\rm{bc}}$,
then, for each $(u,\beta)\in\ol{\Kext}$, $\M({\bf p},z, {\bm\xi})$ satisfies
that,
for all $\bmxi\in\shock$,
\begin{align}
\label{estimate-Dp-M}
&\delta_{\rm bc}\le D_{{\bf p}}\M(D\phi(\bm\xi),\phi(\bm\xi),\bm\xi)\cdot\bm\nu_{\rm sh}(\bm\xi)\le \delta_{\rm bc}^{-1},\\[1mm]
\label{estimate-Dz-M}
&D_z\M(D\phi(\bm\xi),\phi(\bm\xi),\bm\xi)\le -\delta_{\rm bc},
\end{align}
where $\bm\nu_{\rm sh}$ is the unit normal vector to $\shock$ towards the interior of $\Om$.

\begin{proof}
By Lemma \ref{lemma-Mproperty}(c),
it suffices to estimate $D_{\bf p}\M_0(D\phi,\phi,\bmxi)\cdot {\bm\nu}_{\rm sh}$
to prove \eqref{estimate-Dp-M}.
Following Definition \ref{definition-domains-np}, let $\bm \xi^{\lefttop}$ and $\bm \xi^{\righttop}$ be the $\xxi$--coordinates
of points $\lefttop$ and $\righttop$, respectively.
By
Definition \ref{definition-10-3}(i),
$Du(\pm 1,1)=0$, which implies that
$D\phi=D\phi_{\beta}-D\rightvphi$ at ${\bm \xi}^{\lefttop}$ and $\bmxi^{\righttop}$,
for $\phi_{\beta}$ given by \eqref{def-phibeta-star}. By \eqref{7-b7}, we have
 \begin{equation*}
D_{\bm p}\M_0(D\phi(\bmxi^{P_j}),\phi({\bm\xi}^{P_j}), {\bm\xi}^{P_j})\cdot\bm\nu_{\rm sh}(\bm\xi^{P_j})=
\begin{cases}
\leftrho(1-\oM^2)\quad&\tx{for $j=1$},\\[1.5mm]
\rightrho\big(1-(\frac{\neta}{\rightc})^2\big)\quad&\tx{for $j=2$},
\end{cases}
\end{equation*}
for $\oM$ given by \eqref{1-25}.
For each $\beta\in[0, \frac{\pi}{2})$,  $\oM<1\le \leftrho$.
Furthermore, it is shown in \eqref{density-mont-ox}--\eqref{oM-monotonicity}
that $\frac{\dd\leftrho}{\dd\beta}>0$ and $\frac{\dd\oM}{\dd\beta}<0$ for all $\beta\in(0, \frac{\pi}{2})$.
Then there exists a constant $\delta_{\rm bc}^{(1)}\in(0,1)$ depending only on $(\iv, \gam)$ such that
\begin{equation*}
\delta_{\rm bc}^{(1)}
\le \inf_{\beta\in[0,\betadet] }
D_{\bm p}\M_0(D\phi(\bmxi^{P_j}),\phi({\bm\xi}^{P_j}), {\bm\xi}^{P_j})\cdot\bm\nu_{\rm sh}(\bm\xi^{P_j})
\le \frac{1}{\delta_{\rm bc}^{(1)}}\qquad\tx{for $j=1,2$}.
\end{equation*}
By \eqref{estimate-M0-iter}, there exists a constant $\bar{\eps}_{\rm{bc}}\in (0, \eps_0)$
depending only on $(\iv, \gam)$ such that, for each $(u,\beta)\in\ol{\Kext}$,
\begin{equation*}
\frac{\delta_{\rm bc}^{(1)}}{2}\le
D_{\bm p}\M_0(D\phi, \phi, \bm\xi)\cdot\bm\nu_{\rm sh}({\bmxi}) \le \frac{2}{\delta_{\rm bc}^{(1)}}
\qquad
\tx{for all ${\bm\xi}\in \shock\cap \mcl{D}_{\bar{\eps}_{\rm{bc}}}$}.
\end{equation*}

By Definition \ref{definition-10-3}(v)--(vi),
if $\eps$ from Definition \ref{definition-10-3} satisfies that $0<\eps<\bar{\eps}_{\rm{bc}}$,
then there exists a constant $\delta_{\rm bc}^{(2)}>0$ depending only on $(\iv, \gam)$ such that
\begin{equation*}
D_{\bm p}\M_0(D\phi, \phi, \bm\xi)\cdot\bm\nu_{\rm sh}({\bm\xi})
= \rho\Bigl(1-\frac{|D\vphi(\bmxi)|^2}{c^{2}(|D\vphi(\bmxi)|^2, \vphi(\bmxi))}\Bigr)
\ge \delta_{\rm bc}^{(2)}\qquad \tx{for all ${\bm\xi}\in\shock\setminus \D_{\bar{\eps}_{\rm{bc}}/4}$}.
\end{equation*}
Then \eqref{estimate-Dp-M} is obtained from the previous two inequalities.

A direct computation by using \eqref{7-d7} yields that,
for all $\bm\xi=(\xin, \etan)\in B_{M_{\rm dom}}({\bf 0})$,
\begin{align*}
&D_z\M_1(D\leftphi(\bm\xi),\leftphi(\bm\xi),\xin)
    = -\leftrho \oM-(\leftrho-1)\frac{\cos\beta}{\iv},\\
&D_z\M_1({\bf 0},0, \xin)
    = -\rightrho^{2-\gam}\neta-\frac{\rightrho-1}{\iv}.
    \end{align*}
Then there exists a constant $\delta_{\rm bc}^{(3)}>0$ depending only on $(\iv, \gam)$ such that
\begin{equation*}
\underset{\beta\in[0,\betadet]}{\max}
\{D_z\M_1(D\leftphi(\bm\xi),\leftphi(\bm\xi),\xin),
D_z\M_1(\bm 0, 0, \xin)\}\le -\delta_{\rm bc}^{(3)}\;\;\qquad\tx{for all $\bm\xi\in \shock$}.
\end{equation*}
By \eqref{estimate-M1-iter}, there exists a constant ${\sigma}_{\rm{bc}}\in(0,\bar{\sigma}_{\rm{bc}}]$
depending on $(\iv, \gam)$ such that
\begin{equation}
\label{dz-M1-2018}
D_z\M_1({\bf p}, z, {\xin})\le -\frac{\delta_{\rm bc}^{(3)}}{2}
\end{equation}
for all $\bmxi\in B_{M_{\rm dom}}({\bf 0})$ and for all $({\bf p},z)$
satisfying that either $|({\bf p},z)|\le \sigma_{\rm{bc}}$
or $|({\bf p},z)-(D\leftphi, \leftphi(\bm\xi))|\le \sigma_{\rm{bc}}$.
By \eqref{definiiton-iterset-ineq8}, \eqref{12-67},
and  Definition \ref{definition-10-3}(vi),
there exists a constant $\delta_{\rm bc}^{(4)}>0$ depending on $(\iv, \gam)$ such that
\begin{equation*}
D_z\M_0(D\phi({\bm\xi}),\phi(\bm\xi), \bm\xi)=-\frac{1}{\rho^{\gam-2}}D\vphi\cdot {\bm\nu_{\rm sh}}(\bm\xi)\le -\delta_{\rm bc}^{(4)}
\qquad\tx{on $\shock\setminus (\nD_{\eps/10}\cup \oD_{\eps/10})$}.
\end{equation*}
By  Definition \ref{definition-10-3}(i),
$\rho(|D\vphi|^2, \vphi)=\leftrho$ on $\leftsonic$
and $\rho(|D\vphi|^2, \vphi)=\rightrho$ on $\rightsonic$.
Using Definition \ref{definition-10-3}(i),  we can further reduce
$\bar{\eps}_{\rm bc}>0$ depending only on
$(\iv, \gam, \beta_*)$ so that $\rho(|D\vphi|^2, \vphi)\ge \frac{1}{10}\min\{\leftrho, \rightrho\}>0$
on $\shock\cap  (\nD_{\bar{\eps}_{\rm bc}}\cup \oD_{\bar{\eps}_{\rm bc}})$.
Therefore, if $\eps\in (0, \bar{\eps}_{\rm bc})$, then we obtain
\begin{equation}
\label{M0-z-ineq-2018}
D_z\M_0(D\phi({\bm\xi}),\phi(\bm\xi), \bm\xi)=-\frac{1}{\rho^{\gam-2}}D\vphi\cdot {\bm\nu_{\rm sh}}(\bm\xi)\le -\delta_{\rm bc}^{(5)}
\qquad\,\,\tx{on $\shock$}
\end{equation}
for a constant $\delta_{\rm bc}^{(5)}>0$ depending on $(\iv, \gam)$.

Then \eqref{estimate-Dz-M} is obtained by combining inequalities \eqref{dz-M1-2018}--\eqref{M0-z-ineq-2018}.
\end{proof}
\end{lemma}

Hereafter, let $\sigma_{\rm{bc}}>0$ in \eqref{7-d8} be fixed as in Lemma \ref{lemma-7-7}.
This completes the definition of $\M$ in \eqref{7-d8}.

\smallskip
{\textbf{4.}} For $\phi_{\beta}$ given by \eqref{def-phibeta-star}, denote $\psi:=\phi-\phi_{\beta}=\vphi-\vphi_{\beta}^*$.

Let the $(x,y)$--coordinates be defined by \eqref{coord-n} and \eqref{coord-o}
near $\rightsonic$ and $\leftsonic$, respectively.
For  ${\bmxi}=((\rightc-x)\cos y,(\rightc-x)\sin y)$ near $\rightsonic$, and for $\M$ given by \eqref{7-d8},
we use \eqref{cov-right} to define $\hat{\M}^{\mcl{N}}$ by
\begin{equation}
\label{Mn}
\begin{split}
&\hat{\M}^{\mcl{N}}(q_1,q_2,z,x,y)\\
&:=\M(-q_1\cos y-\frac{q_2\sin y}{\rightc-x},-q_1\sin y+\frac{q_2\cos y}{\rightc-x},z,(\rightc-x)\cos y,(\rightc-x)\sin y).
\end{split}
\end{equation}

For $\bmxi=(\leftu-(\leftc-x)\cos(\pi-y),(\leftc-x)\sin(\pi-y))$ near $\leftsonic$, we first denote
\begin{equation*}
\M^{\mcl{O}}({\bf q},z,\bm\xi):=\M({\bf q}+D\leftphi, z+\leftphi,\bm\xi),
\end{equation*}
and then define $\hat{\M}^{\mcl{O}}$ by
\begin{eqnarray}
&&\hat{\M}^{\mcl{O}}(q_1,q_2,z,x,y)\label{Mo} \\
&&:=\M^{\mcl{O}}(-q_1\cos(\pi-y)+\frac{q_2\sin(\pi-y)}{\leftc-x},-q_1\sin(\pi-y)-\frac{q_2\cos(\pi-y)}{\leftc-x},z,\nonumber\\
&&\phantom{aaaaaaaaaaaaaaaaaaaaaaaaa}
\leftu-(\leftc-x)\cos(\pi-y),(\leftc-x)\sin(\pi-y)). \nonumber
\end{eqnarray}

\begin{lemma}
\label{lemma-7-9}
Let
constant $\sigma_2>0$ be
from Lemma {\rm \ref{proposition-sub7}}.
Following Definition {\rm \ref{definition-domains-np}},
let $(x_{P_j}, y_{P_j})$ be the $(x,y)$--coordinates of $P_j$ for $j=1,2$.
Let $\bar{\eps}_{\rm bc}$ be from Lemma {\rm \ref{lemma-7-7}}.
Then there exist $\eps_{\rm{bc}}\in(0,\bar{\eps}_{\rm{bc}})$, $\hat{\sigma}_{\rm{bc}}>0$,
and $C>0$ depending only on $(\iv, \gam)$ such that,
for any $\beta\in[0,\betasonic+\frac{\sigma_2}{4}]$
and all $({\bf q},z)$ satisfying that
\begin{equation}
\label{q-condition}
|({\bf q},z)|\le \hat{\sigma}_{{\rm{bc}}},
\end{equation}
the following properties hold{\rm :}
\begin{itemize}
\item[(a)] If $0<x-x_{\lefttop}\le \eps_{\rm bc}$, then
\begin{equation*}
D_{q_i}\hat \M^{\mcl{O}}({\bf q},z,x,y)\le -C^{-1} \,\,\,\,\tx{for $i=1,2$},\qquad
D_z\hat \M^{\mcl{O}}({\bf q},z,x,y)\le -C^{-1};
\end{equation*}

\item[(b)] If $0<x-x_{\righttop}\le \eps_{\rm bc}$, then
\begin{equation*}
D_{q_i}\hat \M^{\mcl{N}}({\bf q},z,x,y)\le -C^{-1}\;\;\,\,\tx{for $i=1,2$},\qquad
D_z\hat \M^{\mcl{N}}(q_1,q_2,z,x,y)\le -C^{-1}.
\end{equation*}
\end{itemize}

\begin{proof}
By \eqref{cov-right} and \eqref{cov-left}, there exists a constant $\hat{\sigma}_{\rm bc}^*$
depending
only on $(\iv, \gam)$ such that, for
each $\beta\in[0,\betasonic+\frac{\sigma_2}{4}]$, if
$|({\bf q},z)|\le \hat{\sigma}^*_{\rm{bc}}$,
then $\M$ on the right-hand side of \eqref{Mn} and \eqref{Mo} is the same as $\M_1$
given by \eqref{7-d7}.
A direct computation shows that there exists a constant $\til{C}>0$ depending only on
$(\iv, \gam)$ such that, for each $\beta\in[0, \betasonic+\frac{\sigma_2}{4}]$,
 \begin{equation*}
\begin{split}
&D_{q_i}\hat \M^{\mcl{O}}(\mathbf{0},0,x_{\lefttop},y_{\lefttop})\le -{\til C}^{-1},\quad\,\,
D_z\hat{\M}^{\mcl{O}}(\mathbf{0},0,x_{\lefttop},y_{\lefttop})\le -{\til C}^{-1},
\\
&D_{q_i}\hat \M^{\mcl{N}}(\mathbf{0},0,x_{\righttop},y_{\righttop})\le -{\til C}^{-1}, \quad\,\,
 D_z\hat{\M}^{\mcl{N}}(\mathbf{0},0,x_{\righttop},y_{\righttop})\le -{\til C}^{-1}
\end{split}
\end{equation*}
for $i=1,2$.
Then, by Lemma \ref{lemma-Mproperty},
there exist constants $\hat{\sigma}_{\rm{bc}}\in(0, \hat{\sigma}^*_{\rm{bc}}]$ and $C>0$ depending
only on $(\iv, \gam)$ such that properties (a) and (b) hold.
\end{proof}
\end{lemma}

{\textbf{5.}}
The next step is to extend the definition of $\M$ in \eqref{7-d8} to all $({\bf p}, z)\in \R^2\times \R$.

For each $(u,\beta)\in\ol{\Kext}$ and a constant $\sigma>0$, let $v^{(u,\beta)}_{\sigma}\in C^4(\ol{\Om})$ (from Lemma \ref{lemma-c4-approximation-new}) be given.
For a constant $\sigma>0$ to be fixed later, we define a linear operator:
 \begin{equation}
\label{definition-L-iter-bc}
\begin{split}
\mcl{L}^{(u,\beta)}_{\sigma}({\bf p}, z, \bm\xi):=&\,\M(Dv_{\sigma}^{(u,\beta)}(\bm\xi), v_{\sigma}^{(u,\beta)}(\bm\xi),\bm\xi)\\
&\, +D_{\bf p}\M(Dv_{\sigma}^{(u,\beta)}(\bm\xi), v_{\sigma}^{(u,\beta)}(\bm\xi),\bm\xi)\cdot{\bf p}
+D_{z}\M(Dv_{\sigma}^{(u,\beta)}(\bm\xi), v_{\sigma}^{(u,\beta)}(\bm\xi),\bm\xi) z.
\end{split}
\end{equation}
Let $\sigma_{\rm{bc}}>0$ be from Lemma \ref{lemma-7-7}.
By Lemma \ref{lemma-c4-approximation-new}(a), if $\sigma^2<\sigma_{\rm{bc}}$,
then $\mcl{L}^{(u,\beta)}_{\sigma}$ is well defined for all $({\bf p}, z, {\bm\xi})\in \R^2\times \R\times \ol{\Om}$.
For a constant $\sigma\in(0, \sigma_{\rm{bc}})$ to be determined later, depending
only on $(\iv, \gam, \beta_*)$, we define $\M_{(u,\beta)}({\bf p}, z, \bm\xi)$ by
\begin{equation}
\label{definition-final-M}
\begin{split}
\M_{(u,\beta)}({\bf p}, z, \bm\xi):=\varsigma_{\sigma} \M({\bf p}, z, \bm\xi)+
(1-\varsigma_{\sigma})\mcl{L}_{\sigma}^{(u,\beta)}({\bf p}-Dv_{\sigma}^{(u,\beta)}(\bm\xi), z-v_{\sigma}^{(u,\beta)}(\bm\xi), \bm\xi)
\end{split}
\end{equation}
for $\varsigma_{\sigma}=\varsigma_{\sigma}(|({\bf p},z)-(Dv_{\sigma}^{(u,\beta)}(\bm\xi), v_{\sigma}^{(u,\beta)}(\bm\xi))|)$,
where $\varsigma_{\sigma}$ is defined by \eqref{definition-varsigma}.
\smallskip

The following lemma is obtained by adjusting the proofs of
\cite[Lemmas 12.5.7 and 17.3.23]{CF2}
via use of Definition \ref{definition-10-3}, Lemmas \ref{lemma-Mproperty}--\ref{lemma-7-9},
and \eqref{definition-L-iter-bc}--\eqref{definition-final-M}:

\begin{lemma}
\label{lemma-ext-M-property}
Let constants $\bar{\eps}_{\rm{bc}}$ and $\eps_{\rm{bc}}$ be from Lemmas {\rm \ref{lemma-7-7}} and {\rm \ref{lemma-7-9}}, respectively.
Then there exist positive constants $\delta_1^{(1)}$, $N_1^{(1)}$, $\delta_{\rm bc}$, $C$, $C_{\beta_*}$, and
$\eps_{\M}\in(0,{\eps}_{\rm{bc}}]$ with $(\delta_1^{(1)}, N_1^{(1)},\delta_{\rm bc}, C)$
depending on $(\iv, \gam)$, $\eps_{\M}$ depending on $(\iv, \gam,  \beta_*)$,
and $C_{\beta_*}$ depending on $(\iv, \gam,  \beta_*,\alp)$ such that,
if parameters $(\eps, \delta_1, N_1)$ from Definition {\rm \ref{definition-10-3}}
satisfy that $\eps\in(0, \bar{\eps}_{\rm{bc}}]$, $\delta_1\in(0, \delta_1^{(1)}]$,
and $N_1\ge N_1^{(1)}$, then, for each $(u,\beta)\in\ol{\Kext}$,
$\M_{(u,\beta)}:\R^2\times \R\times \ol{\Om}\rightarrow \R$ given by \eqref{definition-final-M}
with $\sigma=\sqrt{\delta_1}$ satisfies the following properties{\rm :}

\smallskip
\begin{itemize}
\item[(a)]
$\M_{(u,\beta)}: \R^2\times \R\times \ol{\Om}\rightarrow \R$ is in $C^3$ and,
for all $({\bf p}, z)\in \R^2\times \R$,
\begin{equation*}
\|(\M_{(u,\beta)}({\bf 0}, 0, \cdot), D^k_{({\bf p}, z)} \M_{(u,\beta)}({\bf p}, z, \cdot))\|_{C^3(\ol{\Om})}
\le C_{\beta_*}\qquad \tx{for $k=1,2,3$};
\end{equation*}

\item[(b)] For $|{\bf p}-D\phi(\bm\xi)|+|z-\phi(\bm\xi)|\le \frac{\sqrt{\delta_1}}{2}$,
\[
\M_{(u,\beta)}({\bf p}, z, \bm\xi)=\M({\bf p}, z, \bm\xi)
\]
for $\M$ defined by \eqref{7-d8}{\rm ;}

\smallskip
\item[(c)] For all $({\bf p}, z,\bm\xi)\in \R^2\times \R\times \ol{\Om}$,
\begin{equation*}
|D_{({\bf p}, z)}\M_{(u,\beta)}({\bf p}, z, \bm\xi)-D_{({\bf p}, z)}\M(D\phi(\bm\xi), \phi(\bm\xi), \bm\xi)|\le C\sqrt{\delta_1};
\end{equation*}

\item[(d)]  For all $({\bf p}, z, \bm\xi)\in \R^2\times \R\times \ol{\shock}$,
\begin{equation*}
\delta_{\rm bc} \le D_{{\bf p}}\M_{(u,\beta)}({\bf p}, z, \bm\xi)\cdot \bm\nu_{\rm sh}\le \frac{1}{\delta_{\rm bc}},\quad
D_z\M_{(u,\beta)}({\bf p}, z, \bm\xi)\le -\delta_{\rm bc},
\end{equation*}
where $\bm\nu_{\rm sh}$ is the unit normal vector to $\shock$ towards the interior of $\Om${\rm ;}

\smallskip
\item[(e)] Representing as
$\mcl{L}_{\sigma}^{(u,\beta)}({\bf p}-Dv^{(u,\beta)}_{\sqrt{\delta_1}}(\bm\xi), z-v^{(u,\beta)}_{\sqrt{\delta_1}}(\bm\xi), \bm\xi)
=\mcl{B}^{(u,\beta)}_{\sigma,\shock}({\bf p},z, \bm\xi),$
define
\begin{equation*}
\mcl{B}^{(u,\beta)}_{\sigma,\shock}({\bf p},z, \bm\xi )=
b_1^{\rm (sh)}(\bm\xi)p_1+b_2^{\rm (sh)}(\bm\xi)p_2
+b_0^{\rm (sh)}(\bm\xi)z+h^{\rm (sh)}(\bm\xi).
\end{equation*}
Then
\begin{equation*}
\|(b_i^{\rm (sh)},h^{\rm (sh)} )\|_{C^3(\ol{\shock})}\le C_{\beta_*}\qquad \tx{for $i=0,1,2$},
\end{equation*}
and, for all $({\bf p}, z,\bm\xi)\in \R^2\times \R\times \ol{\Om}$,
\begin{equation*}
\begin{split}
&|\M_{(u,\beta)}({\bf p}, z, \bm\xi)-\mcl{B}^{(u,\beta)}_{\sqrt{\delta_1},\shock}({\bf p}, z, \bm\xi)|
\le C\sqrt{\delta_1}\big(|{\bf p}-Dv_{\sqrt{\delta_1}}^{(u,\beta)}(\bm\xi)|+|z-v_{\sqrt{\delta_1}}^{(u,\beta)}(\bm\xi)|\big),\\
&|D_{({\bf p}, z)}\M_{(u,\beta)}({\bf p}, z, \bm\xi)-D_{({\bf p}, z)}\mcl{B}^{(u,\beta)}_{\sqrt{\delta_1},\shock}(\bm\xi)|\le C\sqrt{\delta_1};
\end{split}
\end{equation*}

\item[(f)] $\M_{(u,\beta)}$ is homogeneous in the sense that
\[
\begin{cases}
\M_{(u,\beta)}({\bf 0}, 0,\bm\xi)=0,\\[1mm]
\M_{(u,\beta)}(D\leftphi(\bm\xi),\leftphi(\bm\xi),\bm\xi)=0
\end{cases}
\]
for all $\bmxi\in \shock$ when $\beta\in[0, \frac{\delta_1}{N_1}]$,
and for all $\bmxi\in\shock\cap \D_{\eps_{\M}}$ when $\beta\in ( \frac{\delta_1}{N_1},\beta_*]$.

\smallskip
\item[(g)]
Let the $(x,y)$--coordinates be defined by \eqref{coord-n} and \eqref{coord-o}
near $\rightsonic$ and  $\leftsonic$, respectively.  For $\bm\xi\in \ol{\shock\cap \nD_{\eps_{\rm bc}}}$, define
\begin{equation}
\label{Mn-ext}
\begin{split}
&\hat{\M}^{\mcl{N}}_{(u,\beta)}(q_1,q_2,z,x,y)\\
&:=\M_{(u,\beta)}(-q_1\cos y-\frac{q_2\sin y}{\rightc-x},-q_1\sin y+\frac{q_2\cos y}{\rightc-x},z,(\rightc-x)\cos y,(\rightc-x)\sin y).
\end{split}
\end{equation}
For $\bm\xi\in \ol{\shock\cap \oD_{\eps_{\rm bc}}}$, define
$
\M^{\mcl{O}}_{(u,\beta)}({\bf p}, z, \bm\xi):=\M_{(u,\beta)}({\bf p}+D\leftphi, z+\leftphi, \bm\xi),
$
and
\begin{equation}
\label{Mo-ext}
\begin{split}
&\hat{\M}_{(u,\beta)}^{\mcl{O}}(q_1,q_2,z,x,y)\\
&:=\M^{\mcl{O}}_{(u,\beta)}(-q_1\cos(\pi-y)+\frac{q_2\sin(\pi-y)}{\leftc-x},-q_1\sin(\pi-y)-\frac{q_2\cos(\pi-y)}{\leftc-x},\\
&\qquad\qquad\qquad \leftu-(\leftc-x)\cos(\pi-y),(\leftc-x)\sin(\pi-y)).
\end{split}
\end{equation}
Then $\hat{\M}^{\mcl{N}}_{(u,\beta)}$ and $\hat{\M}^{\mcl{O}}_{(u,\beta)}$ satisfy the following properties,
provided that $\shock\cap \oD_{\eps_{\rm{bc}}}$ is nonempty{\rm :}

\smallskip
\begin{enumerate}
\item[\rm (g-1)] $\|\hat{\M}^{\mcl{N}}_{(u,\beta)}\|_{C^3(\R^2\times \R\times \ol{\shock\cap \nD_{\eps_{\rm bc}}})}+
\|\hat{\M}^{\mcl{O}}_{(u,\beta)}\|_{C^3(\R^2\times \R\times \ol{\shock\cap \oD_{\eps_{\rm bc}}})}
\le C_{\beta_*}${\rm ;}

\smallskip
\item[\rm (g-2)] For all $|({\bf q}, z)|\le \frac{\delta_{\rm bc}}{C}$,
\begin{equation*}
\begin{split}
&\hat{\M}^{\mcl{N}}_{(u,\beta)}({\bf q}, z, x, y)=\hat{\M}^{\mcl{N}}({\bf q}, z, x, y)\qquad\;\tx{in $\shock\cap \nD_{\eps_{\rm bc}}$},\\
&\hat{\M}^{\mcl{O}}_{(u,\beta)}({\bf q}, z, x, y)=\hat{\M}^{\mcl{O}}({\bf q}, z, x, y)\qquad\;\tx{in $\shock\cap \oD_{\eps_{\rm bc}}$}
\end{split}
\end{equation*}
for $\hat{\M}^{\mcl{N}}$ and $\hat{\M}^{\mcl{O}}$ defined by \eqref{Mn} and \eqref{Mo}, respectively{\rm ;}

\smallskip
\item[\rm (g-3)] For each $({\bf q},z)\in \R^2\times \R$ and $i=1,2$,
\begin{equation*}
D_{q_i}\hat{\M}^{\mcl{N}}_{(u,\beta)}({\bf q}, z, x, y)
\le -\delta_{\rm bc},\,\,\,\,
D_{z}\hat{\M}^{\mcl{N}}_{(u,\beta)}({\bf q}, z, x, y)
\le -\delta_{\rm bc}\qquad\tx{in $\shock\cap \nD_{\eps_{\M}}$},
\end{equation*}
\begin{equation*}
D_{q_i}\hat{\M}^{\mcl{O}}_{(u,\beta)}({\bf q}, z, x, y)
\le -\delta_{\rm bc},\,\,\,\,
D_{z}\hat{\M}^{\mcl{O}}_{(u,\beta)}({\bf q}, z, x, y)
\le -\delta_{\rm bc}\qquad\tx{in $\shock\cap D^{\mcl{O}}_{\eps_{\M}}$},
\end{equation*}
provided that $\shock\cap \oD_{\eps_{\M}}$ is nonempty{\rm ;}
\end{enumerate}

\smallskip
\item[(h)] $\M_{(u,\beta)}(D{\phi}, {\phi}, \bm\xi)=0$ on $\shock$
if and only if $\vphi=\phi+\rightvphi$ satisfies the Rankine-Hugoniot jump condition \eqref{rhbc-gsh} on $\shock=\{\vphi=\ivphi\}$.
\end{itemize}
\end{lemma}

By \eqref{12-57} and \eqref{definition-final-M}, the definition of the nonlinear boundary value problem \eqref{12-50}
is completed.

\subsection{Well-posedness of the boundary value problem \eqref{12-50}}
\label{subsubsec--existence-hpsi-bvp}

\begin{lemma}
\label{proposition-existence-hpsi-bvp}
Fix $\gam\ge 1$, $\iv>0$, and $\beta_*\in (0,\betadet)$.
Let $\eps_0>0$ be from Lemma {\rm \ref{lemma-7-4}(c)} with $\bar{\beta}$ replaced by $\beta_*$.
Let constant $\sigma_2>0$ be from Lemma {\rm \ref{proposition-sub7}}.
Moreover, let $\bar{\alpha}\in(0,\frac{1}{3}]$
be from Proposition {\rm \ref{proposition-unif-est-u-new}} with $\bar{\beta}$ replaced by $\beta_*$, 
and let $\alpha\in(0,\frac{\bar{\alpha}}{2}]$ be from Definition {\rm 4.19},
Then there exist constants $\eps^{({\rm w})}\in(0,\eps_0]$, $\delta_1^{({\rm w})}\in (0,1)$,
$N_1^{({\rm w})}\ge 1$, and $\alp_1^*\in(0,\bar{\alp}]$
depending only on $(\iv, \gam,  \beta_*)$ such that, whenever parameters $(\eps, \delta_1, N_1)$
from Definition {\rm \ref{definition-10-3}}
satisfy that $\eps\in(0, \eps^{({\rm w})}]$, $\delta_1\in(0, \delta_1^{({\rm w})}]$, and $N_1\ge N_1^{({\rm w})}$, the following properties hold{\rm :}

\smallskip
{\rm {\bf Case} 1}.
If $\beta\le \betasonic+\sigma_2$,
then the boundary
value problem \eqref{12-50} associated
with $(u,\beta)\in \ol{\Kext}\cap\{\beta\le \betasonic+\sigma_2\}$
has a unique solution $\hat{\phi}\in C^2(\Om)\cap C^1(\ol{\Om}\setminus (\ol{\leftsonic}\cup\ol{\rightsonic}))\cap C^0(\ol{\Om})$
for $\Om=\Om(u,\beta)$. Moreover, there exists a constant $C>0$ depending
only on $(\iv, \gam,  \beta_*,\alp)$ such that solution $\hat{\phi}$ satisfies
\begin{equation}
\label{estimate1-iter-bv-sol}
\|\hat{\phi}\|_{L^{\infty}(\Om)}\le C,
\qquad\,
|\hat{\phi}(\bmxi)-\phi_{\beta}^*(\bmxi)|\le C\,{\rm dist}(\bmxi, \leftsonic\cup\rightsonic)\,\,\,\,\,\,\tx{in $\Om$}
\end{equation}
for $\phi_{\beta}^*=\max\{\leftvphi, \rightvphi\}-\rightvphi$.
Furthermore, for each $d\in (0,\eps_0)$, there exists a constant $C_d>0$ depending
only on $(\iv, \gam,  \beta_*, d, \alp)$ such that
\begin{equation}
\label{estimate2-iter-bv-sol}
 \|\hat{\phi}\|_{2,\alp_1^*, \Om\setminus {\D_d}}\le C_d.
 \end{equation}

{\rm {\bf Case} 2}.
For each $\delta\in(0, \frac{\sigma_2}{2})$,
if $\betasonic+\delta\le \beta\le \beta_*$,
then the boundary value problem \eqref{12-50}
associated with $(u,\beta)\in \ol{\Kext}$
 has a unique solution $\hat{\phi}\in C^2(\Om)\cap C^1(\ol{\Om}\setminus (\ol{\leftsonic}\cup\ol{\rightsonic}))\cap C^0(\ol{\Om})$
 for $\Om=\Om(u,\beta)$, and the solution satisfies \eqref{estimate1-iter-bv-sol}--\eqref{estimate2-iter-bv-sol}
 for constants $C>0$ depending only on $(\iv, \gam,  \beta_*, \delta)$
 and $C_d>0$ depending only on $(\iv, \gam,  \beta_*, \delta, d, \alp)$.

\begin{proof}
Fix $(u,\beta)\in \ol{\Kext}\cap\{\beta\le \betasonic+\sigma_2\}$.
Using $\mcl{G}_1^{\beta}$ defined by \eqref{12-16-mod},
we rewrite the boundary value problem \eqref{12-50} associated with fixed $(u,\beta)$
in domain $\mcl{R}=\mcl{G}_1^{\beta}(\Om(u,\beta))$.
Then we follow the argument of Step 1 in the proof of \cite[Proposition 17.4.2]{CF2},
by using
Lemmas \ref{lemma-10-1}, \ref{lemma-12-1-mod}, \ref{lemma-8-3}--\ref{lemma-c4-approximation-new},
and \ref{lemma-ext-M-property},
to choose
constants $\eps^{\rm (w)}\in(0,\eps_0]$, $\delta_1^{\rm (w)}\in (0,1)$, and $N_1^{\rm (w)}\ge 1$
such that, whenever parameters $(\eps, \delta_1, N_1)$ from Definition {\rm \ref{definition-10-3}}
satisfy that $\eps\in(0, \eps^{\rm (w)}]$, $\delta_1\in(0, \delta_1^{\rm (w)}]$, and $N_1\ge N_1^{\rm (w)}$,
the newly written boundary value problem in $\mcl{R}$ satisfies all the conditions
of Proposition \ref{proposition-app-C-wp}.
Then the existence and uniqueness of solution $\hat{\phi}$ of problem \eqref{12-50}
satisfying \eqref{estimate1-iter-bv-sol}--\eqref{estimate2-iter-bv-sol} directly follows from Proposition \ref{proposition-app-C-wp}.

In the case of $\betasonic+\delta\le \beta\le \beta_*$ for $\delta\in(0, \frac{\sigma_2}{2})$,
we follow the argument of Step 2 in the proof of \cite[Proposition 17.4.2]{CF2} by
using
Lemma \ref{lemma2-coeff-iter-eqn-new} and Proposition \ref{proposition-app-C-wp2}
to prove that the boundary value problem \eqref{12-50} associated with $(u,\beta)\in \ol{\Kext}$
has a unique solution $\hat{\phi}$ that satisfies \eqref{estimate1-iter-bv-sol}--\eqref{estimate2-iter-bv-sol}.
\end{proof}
\end{lemma}

For each $(u,\beta)\in\ol{\Kext}$, the corresponding pseudo-subsonic region $\Om=\Om(u,\beta)$ depends
continuously on $(u,\beta)$.
For later discussions, it is useful to rewrite \eqref{12-50} as a boundary value problem for
\begin{equation}
\label{definition-uhat-new}
\hat u(s,t)=(\hat{\phi}+\rightvphi-\vphib^*)\circ \mathfrak{F}_{(u,\beta)}(s,t) \qquad\,\,\tx{in $\iter$}
\end{equation}
for map $\mathfrak{F}=\mathfrak{F}_{(u,\beta)}$ defined by Definition \ref{definition-Gset-shocks-new}(ii),
where $\vphi_{\beta}^*$ is given by \eqref{12-32}.

Substitute expression $\hat{\phi}=\hat u \circ (\mathfrak{F}_{(u,\beta)})^{-1}-(\rightvphi- \vphib^*)$ into \eqref{12-50}
and then rewrite \eqref{12-50} in terms of $\hat u$ to obtain
\begin{equation}
\label{iter-nlbvp-new}
\begin{split}
\sum_{i,j=1}^2 \mcl{A}^{(u,\beta)}_{ij}(D\hat u, s, t)\der_{ij}\hat u+\sum_{i=1}^2\mcl{A}^{(u,\beta)}_i(D\hat u, s, t)\der_i \hat u
  =f^{(u,\beta)} \qquad&\mbox{in $\iter=(-1,1)\times (0,1)$},\\
\mathscr{M}_{(u,\beta)}(D\hat u, \hat u, s)=0 \qquad&\mbox{on $\dershock:=(-1,1)\times \{1\}$},\\
\hat u=0 \qquad &\mbox{on $\dersonic:=\{-1,1\}\times (0,1)$},\\
\mathscr{B}^{\rm (w)}_{(u,\beta)}(D\hat u, s):=b_{1}^{\rm (w)}(s)\der_1\hat u+ b_{2}^{\rm (w)}(s) \der_2\hat u=0
\qquad &\mbox{on $\derwedge:=(-1,1)\times\{0\}$},
\end{split}
\end{equation}
where $(\der_1, \der_2)=(\der_s, \der_t)$.

Since $\rightvphi-\vphi_{\beta}^*=0$ when $\beta=0$, we have
\begin{align}
\label{iter-bvp-f-b0}
&f^{(u,\beta)}\equiv 0\qquad\qquad\qquad \tx{if $\beta=0$},\\
\label{M-for-u-homeg}
&\mathscr{M}_{(u,0)}({\bf 0}, 0, s)=0\qquad\,\,\tx{on $\dershock$},
\end{align}
where \eqref{M-for-u-homeg} follows from Lemma \ref{lemma-ext-M-property}(f).

From Lemmas \ref{lemma-7-4}, \ref{lemma-6-1}, and \ref{lemma-ext-M-property}--\ref{proposition-existence-hpsi-bvp},
the following lemma is obtained:

\begin{lemma}
\label{remark-prop-coeff-iter-bvp}
For each $(u,\beta)\in \ol{\Kext}$, let $\mcl{A}^{(u,\beta)}_{ij}$, $\mcl{A}^{(u,\beta)}_i$, $f^{(u,\beta)}$, $\mathscr{M}_{(u,\beta)}$,
$\mathscr{B}^{\rm (w)}_{(u,\beta)}$, and $b_{j,(u,\beta)}^{\rm (w)}$ be as those in \eqref{iter-nlbvp-new}.
Then the following properties hold{\rm :}

\smallskip
\begin{itemize}
\item[(a)]
$\mcl{A}^{(u,\beta)}_{ij}, \mcl{A}^{(u,\beta)}_i\in C(\R^2\times \iter)$, $f^{(u,\beta)}\in C(\iter)$,
$\mathscr{M}_{(u,\beta)}\in C(\R^2\times \R\times \dershock)$, and\\
$\mathscr{B}^{\rm (w)}_{(u,\beta)}\in C(\R^2\times \R\times \derwedge)${\rm ;}

\smallskip
\item[(b)] Suppose that a sequence $\{(u_k,\beta_k)\}_{k=1}^{\infty}\subset \ol{\Kext}$ converges
to $(u,\beta)\in \ol{\Kext}$ in $C^{2,\alp}_{(*, \alpha_1)}(\iter)\times [0, \beta_*]$ as $k\to \infty$.
Then the following properties hold{\rm :}

\smallskip
\begin{itemize}
\item $(\mcl{A}^{(u_k,\beta_k)}_{ij}, \mcl{A}^{(u_k,\beta_k)}_i)\rightarrow (\mcl{A}_{ij}^{(u,\beta)}, \mcl{A}_i^{(u,\beta)})$
   uniformly on compact subsets of $\R^2\times \iter${\rm ;}

\smallskip
\item $f^{(u_k,\beta_k)}\rightarrow f^{(u,\beta)}$ uniformly on compact subsets of $\iter${\rm ;}

\smallskip
\item $\mathscr{M}_{(u_k, \beta_k)}\rightarrow \mathscr{M}_{(u,\beta)}$ uniformly on compact  subsets of $\R^2\times \R\times \dershock${\rm ;}

\smallskip
\item $\mathscr{B}^{\rm (w)}_{(u_k, \beta_k)}\rightarrow \mathscr{B}^{\rm (w)}_{(u,\beta)}$ uniformly on compact subsets of $\R^2\times \derwedge$.
\end{itemize}
\end{itemize}
\end{lemma}

From Lemmas \ref{lemma-7-4} and \ref{proposition-existence-hpsi-bvp}--\ref{remark-prop-coeff-iter-bvp},
we obtain the following corollary:
\begin{corollary}
\label{corollary-u-convergence}
Let constants $\eps^{\rm (w)}$, $\delta_1^{\rm (w)}$, and $N_1^{\rm (w)}$
be from {\rm{Lemma \ref{proposition-existence-hpsi-bvp}}}.
Let parameters $\eps, \delta_1$, and $N_1$ from Definition {\rm \ref{definition-10-3}}
satisfy that $\eps\in(0, \eps^{\rm (w)}]$, $\delta_1\in(0, \delta_1^{\rm (w)}]$, and $N_1\ge N_1^{\rm (w)}$.

\smallskip
\begin{itemize}
\item[(a)] For each $(u,\beta)\in \ol{\Kext}$, $\hat{\phi}$ solves the boundary value problem \eqref{12-50}
if and only if $\hat u$ given by \eqref{definition-uhat-new} solves the boundary value problem \eqref{iter-nlbvp-new}.
Thus, \eqref{iter-nlbvp-new} has a unique
solution $\hat u\in C^2(\iter)\cap C^1(\ol{\iter}\setminus \ol{\dersonic})\cap C(\ol{\iter})$.
Furthermore, there exists a constant $C\ge 1$ depending on $(\iv, \gam,  \beta_*,\alp)$ such that
\begin{equation*}
|\hat u(s,t)|\le C(1-|s|)\qquad\mbox{in $\iter$}.
\end{equation*}
For each $\hat d\in(0, \frac 12)$, there exists $C_{\hat d}$ depending on $(\iv, \gam,  \beta_*, \hat{d},\alp)$ such that
\begin{equation*}
\|\hat u\|_{2,\alp_1^*,\iter\cap\{1-|s|>\hat d\}}\le C_{\hat d},
\end{equation*}
where constant $\alp_1^*\in(0,\bar{\alp}]$ is from Lemma {\rm \ref{proposition-existence-hpsi-bvp}}.

\smallskip
\item[(b)] For each $(u_k, \beta_k)\in\ol{\Kext}$, let $\hat u_k$ be the solution of the boundary value
problem \eqref{iter-nlbvp-new} associated with $(u_k,\beta_k)$.
Suppose that sequence $\{(u_k, \beta_k)\}$ converges to $(u,\beta)\in\ol{\Kext}$ in $C^{1}(\ol{\iter})\times [0,\beta_*]$.
Then there exists a unique solution $\hat u\in C^2(\iter)\cap C^1(\ol{\iter}\setminus \ol{\dersonic})\cap C(\ol{\iter})$ to
the boundary value problem \eqref{iter-nlbvp-new} associated with $(u,\beta)$. Moreover, $\hat u_k$ converges to
$\hat u$ in the following senses{\rm :}

\smallskip
\begin{itemize}
\item uniformly in $\ol{\iter}$,

\smallskip
\item in $C^{1,\alp'}(K)$ for any compact subset $K\subset \ol{\iter}\setminus \ol{\dersonic}$ and any $\alp'\in[0,\alp_1^*)$,

\smallskip
\item in $C^{2,\alp'}(K)$ for any compact subset $K\subset \iter$ and any $\alp'\in[0, \alp_1^*)$.
\end{itemize}

\smallskip
\item[(c)] If $(u,\beta)\in \ol{\mcl{K}}$, then $(u,\beta)$ satisfies property {\rm (vii)} of Definition {\rm \ref{definition-10-3}}
with nonstrict inequality in \eqref{12-65}.
\end{itemize}
\end{corollary}

\begin{remark}
\label{remark-u-convergence}
For a constant $M>0$, define a set $\mcl{K}^{E}_M$ by
\begin{equation*}
\mcl{K}^{E}_{M}:=
\big\{(u,\beta)\in C^{2,\alp}_{(*,\alp_1)}(\iter)\;:\; \|u\|_{2,\alp,\iter}^{(*,\alp_1)}\le M, \,
\tx{$(u,\beta)$ satisfy {\rm (ii)}--{\rm (vi)} of Definition {\rm \ref{definition-10-3}}}
\big\}.
\end{equation*}

Let $\ol{\mcl{K}^{{E}}_{M}}$ be the closure of $\mcl{K}^{{E}}_M$ in $C^{2,\alp}_{(*,\alp_1)}(\iter)\times [0,\beta_*]$.
Then Lemma {\rm \ref{remark-prop-coeff-iter-bvp}} and Corollary {\rm \ref{corollary-u-convergence}}
still hold  when $\ol{\Kext}$
is replaced by $\ol{\mcl{K}^{E}_{M}}$
for some constant $M>0$.
\end{remark}

\section{Properties of the Iteration Set $\mcl{K}$}
\label{subsec-prop-iterset}

\subsection{Admissible solutions}
\quad
As stated in Definition \ref{definition-10-3}, parameter $\alp$ for the iteration set $\mcl{K}$ will be
chosen in $(0, \frac{\bar{\alp}}{2}]$, where $\bar{\alp}\in(0,1)$ is the constant in
Proposition  \ref{proposition-unif-est-u-new}.

\begin{lemma}
\label{lemma:convergence-adm-sol-normal}
Given $\gam \ge 1$ and $\iv>0$, fix $\beta_*\in(0,\betadet]$.
Take a sequence $\{\beta_j\}_{j=1}^{\infty}\subset (0,\beta_*]$ such that $\beta_j$ converges to $0$ as $j\to \infty$.
For each $j\in \mathbb{N}$, let $\vphi^{(j)}$ be an admissible solution corresponding to $(\iv, \beta_j)$.
Let $u^{(j)}$ be defined by \eqref{10-12} corresponding to $(\vphi_j, \beta_j)$.
Then there exists a subsequence of $\{u^{(j)}\}$ converging in $C^{2,\alp}_{(*,\alp_1)}(\iter)$ to $\un\equiv 0$.

\begin{proof}
By Proposition \ref{proposition-unif-est-u-new} and \eqref{alpha-choice},
sequence $\{u^{(j)}\}$ is uniformly bounded in $C^{2,2\alp}_{(*,1)}(\iter)$.
Since $C^{2,2\alp}_{(*,1)}(\iter)$ is compactly embedded into $C^{2,\alp}_{(*,\alp_1)}(\iter)$,
there exists a subsequence (still denoted as) $\{u^{(j)}\}$
such that the subsequence converges in $C^{2,\alp}_{(*,\alp_1)}(\iter)$
to a function $u^{(\infty)}\in C^{2,\alp}_{(*,\alp_1)}(\iter)$.

By \eqref{iter-bvp-f-b0}, Lemma \ref{remark-prop-coeff-iter-bvp}, Corollary \ref{corollary-u-convergence},
and Remark \ref{remark-u-convergence},  we see that  $u=u^{(\infty)}$ is the solution of the nonlinear boundary value problem:
\begin{equation}
\label{nlbvp-limit}
\begin{split}
\sum_{i,j=1}^2 \mcl{A}^{(u,0)}_{ij}(Du, s, t)\der_{ij}u+\sum_{i=1}^2\mcl{A}^{(u,0)}_i(Du, s, t)\der_i u=0 \qquad&\mbox{in $\iter$},\\
\mathscr{M}_{(u,0)}(Du, u, s)=0 \qquad&\mbox{on $\dershock$},\\
u=0 \qquad &\mbox{on $\dersonic$},\\
\mathscr{B}^{\rm (w)}_{(u,0)}(Du, s):=b_{1}^{\rm (w)}(s)\der_1u+ b_{2}^{\rm (w)}(s) \der_2u=0\qquad &\mbox{on $\derwedge$}.
\end{split}
\end{equation}
Owing to \eqref{M-for-u-homeg}, $u=0$ is the solution of the boundary value problem \eqref{nlbvp-limit}.
Then $u^{(\infty)}=0$ in $\iter$ by the uniqueness of solutions.
In other words,  $u^{(\infty)}=\un$ in $\iter$.
\end{proof}
\end{lemma}

\begin{corollary}
\label{corollary-admisol-iterset}
Let constants $\eps^{\rm (w)}$, $\delta_1^{\rm (w)}$, and $N_1^{\rm (w)}$ be
from Lemma {\rm \ref{proposition-existence-hpsi-bvp}},
and let parameters $(\eps, \delta_1)$ in Definition {\rm \ref{definition-10-3}}
be fixed from $(0,\eps^{\rm (w)}]\times (0,\delta_1^{\rm (w)}]$.
For each admissible solution $\vphi$ corresponding to $(\iv, \beta)\in\mathfrak{R}_{\rm weak}\cap\{0\le \beta\le \beta_*\}$
in the sense of Definition {\rm \ref{def-regular-sol}},
let a function $u=u^{(\vphi,\beta)}$ be given by \eqref{10-12}.
Let $N_1$ be the parameter in Definition {\rm \ref{definition-10-3}}.
For each $\delta_1\in(0, \delta_1^{\rm (w)}]$, there exists a constant
$N_1^{\rm (a)}\in [N_1^{\rm (w)},\infty)$ depending only on $(\iv, \gam,  \beta_*, \delta_1)$ such that,
if $N_1\ge N_1^{\rm (a)}$, then $(u^{(\vphi,\beta)}, \beta)\in \mcl{K}$
for each admissible solution $\vphi$ corresponding to $(\iv, \beta)\in\mathfrak{R}_{\rm weak}\cap \{0\le \beta\le \beta_*\}$.

\begin{proof}
For a fixed admissible solution $\vphi$ corresponding to $(\iv, \beta)\in\mathfrak{R}_{\rm weak}\cap\{0\le \beta\le \beta_*\}$,
let $u=u^{(\vphi,\beta)}$ be given by \eqref{10-12}. For simplicity of notation,
denote $u$ as $u^{(\vphi,\beta)}$ in this proof.

By the choice of constants $N_i\, (i=2,3,4,5)$, $\mu_j\, (j=0,1)$, $\til{\mu}$, $\sigma_1$, $\hat{\zeta}$, and $C$
in Definition \ref{definition-10-3}, $(u,\beta)$ satisfy properties (ii)--(vi) of Definition \ref{definition-10-3}.

 By the choice of constant $N_0$ in Definition \ref{definition-10-3}(i), $u$ satisfies
 \begin{equation*}
 \|u-\un\|_{2,\alp,\iter}^{(*,\alp_1)}<N_0
 \end{equation*}
 for any admissible solution $\vphi$ corresponding to $(\iv, \beta)\in\mathfrak{R}_{\rm weak}\cap\{0\le \beta\le \beta_*\}$.
 Lemma \ref{lemma:convergence-adm-sol-normal} implies that,
 for any given constant $\delta_1\in(0,\delta_1^{\rm (w)}]$,
 a constant $N_1^{\rm (a)}\in[N_1^{\rm (w)}, \infty)$  can be chosen depending only
 on $(\iv, \gam,  \beta_*, \delta_1)$ such that,
 whenever $\beta\in[0, \frac{2\delta_1}{N_1^{\rm(a)}}]$, $u$ satisfies
 \[
 \|u-\un\|_{2,\alp,\iter}^{(*,\alp_1)}<\frac{\delta_1}{2}.
 \]
Therefore, if $N_1\ge N_1^{\rm(a)}$, then
any $(u,\beta)$ given by \eqref{10-12} for an admissible solution $\vphi$
corresponding to $(\iv,\beta)\in\mathfrak{R}_{\rm weak}\cap\{0\le \beta\le\beta_*\}$ satisfies
property (i) of Definition \ref{definition-10-3}.
This implies that $(u,\beta)\in \Kext$.
Therefore, Lemmas \ref{lemma-8-3}, \ref{lemma2-coeff-iter-eqn-new}, \ref{lemma-6-1},
and \ref{lemma-ext-M-property} apply to the nonlinear differential operators
$(\mcl{N}_{(u,\beta)}, \mcl{M}_{(u,\beta)})$.
Then, by Propositions \ref{lemma-est-sonic-general}, \ref{proposition-sub8}, and \ref{proposition-sub9},
and Corollary \ref{corollary-u-convergence},
we conclude that $u$
is the unique solution of the boundary value problem \eqref{iter-nlbvp-new}
associated
with $(u,\beta)$.
That is,
$\hat u=u$ in $\iter$, for $\hat u$ is given by \eqref{definition-sol-ubvp}.
Thus, $(u,\beta)$
satisfies property (vii) of Definition \ref{definition-10-3}.

Therefore, we conclude that $(u^{(\vphi,\beta)},\beta)\in \mcl{K}$ for any admissible solution $\vphi$
corresponding to $(\iv,\beta)\in\mathfrak{R}_{\rm weak}\cap\{0\le \beta\le \beta_*\}$
in the sense of Definition \ref{def-regular-sol}.
\end{proof}
\end{corollary}

\subsection{Openness of $\mcl{K}$}
\label{subsec-bvp-iter}
Let $\eps$, $\delta_1$, $\delta_2$, $\delta_3$, and $N_1$ be the parameters from Definition \ref{definition-10-3}.
In this chapter, we further adjust parameters $(\eps, \delta_1)$, then choose $\delta_3>0$ small, depending only
on $(\eps, \delta_1)$ such that Definition \ref{definition-10-3} determines a relatively open
subset of $C^{2,\alp}_{(*,\alp_1)}(\iter)\times [0,\beta_*]$.

\begin{lemma}
\label{lemma-Kext-open}
For each $\beta_*\in(0,\betadet)$, the function set $\Kext$ given by Definition {\rm \ref{definition-K-ext}}
is relatively open in $C^{2,\alp}_{(*,\alp_1)}(\iter)\times [0,\beta_*]$.

\begin{proof}
For each $j=1,2,3$, function $\mathscr{K}_j(\beta)$ of $\beta$ in Definition \ref{definition-10-3}
is continuous for $\beta\in[0,\beta_*]$. Since $\leftvphi$ defined in \eqref{2-4-b6} depends continuously
on $\beta\in [0, \frac{\pi}{2})$,
$\vphi_{\beta}=\max\{\leftvphi, \rightvphi\}$ and $\vphi_{\beta}^*$ defined in \eqref{12-32}
also depend continuously on $\beta\in[0, \frac{\pi}{2})$.
Moreover, $\sbeta$ and $L_{\beta}$ defined in \eqref{def-sbeta} and \eqref{7-b8}, respectively,
depend continuously on $\beta\in[0, \frac{\pi}{2})$.
Furthermore, for each $\beta\in[0,\beta_*]$,
\begin{equation*}
\sup_{\Qbeta(s^*)}(\ivphi-\vphi_{\beta}^*)-\inf_{\Qbeta(s^*)}(\ivphi-\vphi_{\beta}^*)>0
\qquad \mbox{for all $s^*\in[\sbeta, \rightc]$},
\end{equation*}
where  $\Qbeta(s^*)$ is defined in \eqref{definition-Qbeta-str}.

By Lemma \ref{lemma-7-4} and the observations stated above,
the set determined by conditions (i)--(vi) of Definition \ref{definition-10-3} is relatively open
in $C^{2,\alp}_{(*,\alp_1)}(\iter)\times [0,\beta_*]$,
because $C^{2,\alp}_{(*,\alp_1)}(\iter)$ is compactly
embedded in $C^1(\ol{\iter})$;
for further details, we refer to the proofs of \cite[Lemmas 12.8.1 and 17.5.1]{CF2}.
\end{proof}
\end{lemma}

\begin{lemma}
\label{lemma-positivity-hatphi-2018}
Let $\eps^{\rm (w)}, \delta_1^{\rm (w)}, N_1^{\rm (w)}$, and $\alp_1\in(0, \bar{\alp}]$ be
from {\rm{Lemma \ref{proposition-existence-hpsi-bvp}}}.
Let $\eps_0>0$ be from Lemma {\rm \ref{lemma-7-4}(c)}.
Then there exists $\eps^{\rm (lb)}\in(0,\eps^{\rm (w)}]$
depending only on $(\iv, \gam,  \beta_*)$ such that,
whenever parameters $(\eps, \delta_1, N_1)$ in Definition {\rm \ref{definition-10-3}}
are from $(0,\eps^{\rm (lb)}]\times (0, \delta_1^{\rm (w)}]\times[N_1^{\rm (w)}, \infty)$,
there is $\bar{\delta}_3>0$  depending only on $(\iv, \gam,  \beta_*, \delta_1, \delta_2, N_1)$
for $\delta_2$ from Definition {\rm \ref{definition-10-3}(iv)} so that, if parameter
$\delta_3$ in Definition {\rm \ref{definition-10-3}(vii)} satisfies that
$\delta_3\in(0, \bar{\delta}_3]$, then the following properties hold{\rm :}
For each $(u^{\sharp}, \beta^{\sharp})\in \ol{\mcl{K}}$,
a constant $\delta^{\sharp}>0$ can be chosen depending only
on $(\iv, \gam,  \beta_*, u^{\sharp}, \beta^{\sharp})$ such that
solution $\hat{\phi}$ of the boundary value problem \eqref{12-50} associated with $(u,\beta)$ satisfies
\begin{equation}
 \label{strict-positivity}
\hat{\phi}-(\vphi_{\beta}^*-\rightvphi)>  0\qquad\tx{in $\Om$}
\end{equation}
for $\Om=\Om(u,\beta)$, provided that
$(u,\beta)\in \ol{\Kext}$ satisfies
\begin{equation}
\label{7-e9-new2015}
\|u^{\sharp}-u\|_{C^1(\ol{\iter})}+|\beta^{\sharp}-\beta| \le \delta^{\sharp}.
\end{equation}

\begin{proof}
We consider two cases separately:
(i) $\beta^{\sharp}\in [\frac{2\delta_1}{N_1^2}, \beta_*]$
and (ii) $\beta^{\sharp}\in [0, \frac{2\delta_1}{N_1^2}]$.
\smallskip

{\textbf{1.}} Suppose that $\beta^{\sharp}\in[\frac{2\delta_1}{N_1^2}, \beta_*]$.
By  \eqref{definiiton-iterset-ineq1} in Definition \ref{definition-10-3}(iv),
$u^{\sharp}$ satisfies
\[
u^{\sharp}> \frac{\delta_1\delta_2}{N_1^2}\qquad\tx{in}\,\, \ol{\iter}\cap\{1-|s|\ge \frac{\til{\eps}}{10}\}
\]
for $\til{\eps}=\frac{2\eps}{\rightc-s_{\beta^\sharp}}$.
If $\delta_3>0$ satisfies
\begin{equation}
\label{delta3-choice}
\delta_3\le \frac{\delta_1}{2N_1^2}\delta_2,
\end{equation}
then it follows from \eqref{12-65} that
$\hat u^{\sharp}:=(\hat{\phi}^{\sharp}+\rightvphi-\vphi_{\beta^{\sharp}}^*)
\circ \mathfrak{F}_{(u^{\sharp}, \beta^{\sharp})}$
satisfies
\begin{equation}
\label{est-u-sharp}
\hat u^{\sharp}> \frac{\delta_1}{2N_1^2}\delta_2\qquad\mbox{in}\,\,\ol{\iter}\cap\{1-|s|\ge \frac{\til{\eps}}{10}\}
\end{equation}
for $\til{\eps}=\frac{2\eps}{\rightc-s_{\beta^\sharp}}$,
provided that $\hat{\phi}^{\sharp}$
is the solution of the boundary value problem \eqref{12-50} associated
with $(u^{\sharp}, \beta^{\sharp})$.

Note that $\hat{u}^{\sharp}$ is the solution of \eqref{iter-nlbvp-new} determined by
$(u^{\sharp}, \beta^{\sharp})$.
Then, by Corollary \ref{corollary-u-convergence},
there exists a constant $\delta^{\sharp}>0$ small, depending
on $(\iv, \gam,  \beta_*, \delta_3, u^{\sharp}, \beta^{\sharp})$,
such that, if $(u,\beta)\in\ol{\Kext}$ satisfies \eqref{7-e9-new2015},
then \eqref{est-u-sharp} implies that $\hat{u}$ given by \eqref{definition-sol-ubvp} satisfies
\begin{equation}
\label{estimate-hatu-posbd-new}
\hat u> \frac{\delta_1}{4N_1^2}\delta_2\qquad\tx{in}\,\, \ol{\iter}\cap\{1-|s|\ge \frac{\til{\eps}}{10}\}.
\end{equation}

For a constant $r>0$, denote $\D_r:=\nD_r\cup\oD_r$ for $\nD_r$ and $\oD_r$
defined by \eqref{definition-Dr-ext}.
By Proposition \ref{lemma-7-4}(c),
$\mathfrak{F}_{(u,\beta)}^{-1}({\D_{\eps/10}})={\iter}\cap \{1-|s|< \frac{\til{\eps}}{10}\}$.
Thus, \eqref{estimate-hatu-posbd-new} implies
\begin{equation}
\label{pos-hphi-away-sonic}
\hat{\phi}-(\vphi_{\beta}^*-\rightvphi)=\hat u\circ \mathfrak{F}_{(u,\beta)}^{-1}>0\qquad\tx{ in $\Om\setminus \D_{\eps/10}$}.
\end{equation}

Define
\begin{equation}
\label{definition-hpsi-2018}
\hat{\psi}:=\hat{\phi}-(\vphi_{\beta}^*-\rightvphi)\qquad \tx{in $\Om\cap \D_{\eps/2}$}.
\end{equation}
By \eqref{7-b7}, we have
\begin{equation}
\label{hpsi-hphi}
\hat{\psi}=\begin{cases}
\hat{\phi}-(\leftvphi-\rightvphi)\quad &\mbox{in $\Om\cap \oD_{\eps/2}$},\\[1mm]
\hat{\phi}\quad &\mbox{in $\Om\cap \nD_{\eps/2}$},
\end{cases}
\end{equation}
provided that the condition{\rm :}
\begin{equation}
\label{eps-condition-vphi-betastar}
\eps<\frac{2\leftch}{\bar{k}}
\end{equation}
holds for $\bar{k}>1$ from \eqref{7-b7}.

By \eqref{def-uniform-ptnl-new},
$\leftvphi-\rightvphi$ is a linear function depending only on $\xi_1$.
Since $\hat{\phi}$ is a solution of the boundary value problem \eqref{12-50}
associated with $(u,\beta)$, $\hat{\psi}$ satisfies
\begin{equation*}
  \begin{split}
  \mcl{L}_{(u,\beta)}(\hat{\psi}):=\sum_{i,j=1}^2A_{ij}(D\hat{\phi}, {\bmxi})\der_{\xi_i\xi_j}\hat{\psi}=0\qquad\,\,&\mbox{in $\Om\cap \oD_{\eps/2}$},\\
  \hat{\psi}=0\qquad&\mbox{on $\leftsonic$},\\
  \der_{\xi_2}\hat{\psi}=0\qquad&\mbox{on $\Wedge\cap \der \oD_{\eps/2}$},
  \end{split}
\end{equation*}
where $(A_{ij}(D\hat{\phi}, {\bmxi}))_{i,j=1}^2$ is given by \eqref{12-57}.
By Lemma \ref{lemma-6-1}(g)--(h),
$\mcl{L}_{(u,\beta)}(\hat{\psi})=0$ is strictly elliptic in $\oD_{\eps/2}$.
By Lemma \ref{lemma-ext-M-property}(f),
the boundary condition $\mcl{M}_{(u,\beta)}(D\hat{\phi}, \hat{\phi},\bmxi)=0$ on $\shock\cap \der\oD_{\eps/2}$ is equivalent to
\begin{equation*}
  \mcl{M}_{(u,\beta)}(D\hat{\phi}, \hat{\phi}, \bmxi)-\mcl{M}_{(u,\beta)}(D(\leftvphi-\rightvphi), \leftvphi-\rightvphi, \bmxi)=0
  \qquad\tx{on $\shock\cap \der\oD_{\eps/2}$}.
\end{equation*}
By Lemma \ref{lemma-ext-M-property}(d), the boundary condition stated immediately above can be rewritten as
\begin{equation*}
  {\bm\beta}\cdot\nabla\hat{\psi}-\mu\hat{\psi}=0\qquad\,\,\tx{on $\shock\cap \der\oD_{\eps/2}$},
\end{equation*}
where ${\bm\beta}$ and $\mu$ satisfy
\begin{equation*}
  \delta_{\rm bc}\le {\bm\beta}\cdot{\bm\nu}_{\rm sh}\le \delta_{\rm bc}^{-1},\quad \mu\ge \delta_{\rm bc}\qquad\,\,\tx{on $\shock\cap \der\oD_{\eps/2}$}
\end{equation*}
for constant $\delta_{\rm bc}>0$ from Lemma \ref{lemma-ext-M-property}(d)
and   the unit normal vector ${\bm\nu}_{\rm sh}$ to $\shock$ towards the interior of $\Om$.

By \eqref{pos-hphi-away-sonic}, the strong maximum principle, and Hopf's lemma,
we obtain that $\hat{\psi}>0$ in $\oD_{\eps/2}$, which implies that
\begin{equation}
\label{hat-u1}
  \hat u>0\qquad\tx{in $\iter\cap \{-1<s<-1+\frac{\til{\eps}}{2}\}$},
\end{equation}
provided that condition \eqref{eps-condition-vphi-betastar} holds.

By using \eqref{hpsi-hphi}, Lemma \ref{lemma-6-1}(a), and properties (d) and (f) of Lemma \ref{lemma-ext-M-property},
it can be similarly checked that
\begin{equation}
\label{hat-u2}
  \hat u>0\qquad\tx{in $\iter\cap \{1-\frac{\til{\eps}}{2}<s<1\}$}.
\end{equation}
From \eqref{estimate-hatu-posbd-new} and \eqref{hat-u1}--\eqref{hat-u2}, we obtain that
$\hat{u}>0$ in $\iter$, provided that $\delta^{\sharp}>0$ is chosen sufficiently small and
$\eps$ satisfies \eqref{eps-condition-vphi-betastar}.
This proves \eqref{strict-positivity} for $\beta^{\sharp}\in[\frac{2\delta_1}{N_1^2}, \beta_*]$.

\smallskip
{\textbf{2.}} Suppose that $\beta^{\sharp}\in [0, \frac{2\delta_1}{N_1^2}]$.
Choose $\delta^{\sharp}\in(0, \frac{2\delta_1}{N_1^2})$ so that \eqref{7-e9-new2015}
implies that $\beta\in [0, \frac{\delta_1}{N_1})$.
By Lemma \ref{lemma-ext-M-property}(d),
the maximum principle applies to solution $\hat{\phi}$ of the boundary value problem \eqref{12-50}
associated with $(u,\beta)\in \ol{\Kext}$ satisfying \eqref{7-e9-new2015} so that
\begin{equation}
\label{minimum-1}
\hat{\phi}>0\qquad\mbox{in $\Om$}.
\end{equation}
For $(\leftvphi, \rightvphi)$ given by \eqref{def-uniform-ptnl-new},
denote $\phi_{\beta}:=\leftvphi-\rightvphi$.
Since ${\phi}_{\beta}$ is a linear function of $\xxi$,
$\hat{\phi}-\phi_{\beta}$ satisfies
\begin{equation*}
\mcl{N}_{(u,\beta)}(\hat{\phi}-\phi_{\beta})=\mcl{N}_{(u,\beta)}(\hat{\phi})=0
\qquad\mbox{in $\Om$}
\end{equation*}
for the second-order differential operator \eqref{definition-nub-hphi}.
From properties (d) and (f) of Lemma \ref{lemma-ext-M-property}, it follows that
$\M_{(u,\beta)}(D\hat{\phi}, \hat{\phi}, \bm\xi)-\M_{(u,\beta)}(D{\phi_{\beta}}, {\phi_{\beta}}, \bm\xi)=0$
for all $\bm\xi\in \shock$.
This condition can be written as
\begin{equation*}
{\bf b}\cdot D_{\bmxi}(\hat\phi-{\phi}_{\beta})+b_0(\hat{\phi}-{\phi}_{\beta})=0\qquad \mbox{on $\shock$},
\end{equation*}
where ${\bf b}$ and $b_0$ satisfy that ${\bf b}\cdot {\bm\nu}_{\rm sh}>0$
and $b_0<0$ on $\shock$ for  the unit normal vector ${\bm\nu}_{\rm sh}$
to $\shock$ towards the interior of $\Om$.
Then the comparison principle implies that $\hat{\phi}\ge {\phi}_{\beta}$ in $\Om$.
Furthermore, $\hat{\phi}=0>\phi_{\beta}$ on $\rightsonic$.
By the strong maximum principle, we conclude that
\begin{equation}
\label{minimum-2}
\hat{\phi}>{\phi}_{\beta}\qquad\mbox{in $\Om$}.
\end{equation}
Then \eqref{strict-positivity} is obtained from
\eqref{minimum-1}--\eqref{minimum-2},
because
$\max\{0, \phi_{\beta}\}\ge \vphib^*-\rightvphi$ holds in $\Om$.
\end{proof}
\end{lemma}

\begin{lemma}[{\emph{Estimate of $\hat{\phi}$ away from $\leftsonic$}}]
\label{proposition-est-hpsi-new2015}
Let $\eps_0>0$ be from Lemma {\rm \ref{lemma-7-4}(c)}.
Let $\eps^{\rm (w)}$, $\delta_1^{\rm (w)}$, $N_1^{\rm (w)}$, and $\alp_1^*\in(0, \bar{\alp}]$ be
from {\rm{Lemma \ref{proposition-existence-hpsi-bvp}}}.
Let $\eps^{\rm(lb)}$ and $\bar{\delta}_3$ be from Lemma {\rm \ref{lemma-positivity-hatphi-2018}}.
For a constant $r>0$, let $\oD_r$ be defined by \eqref{definition-Dr-ext}.
Then there exist $\eps^{\rm (par)}\in(0,\eps^{\rm (lb)}]$
depending only on $(\iv, \gam,  \beta_*)$
and $C>0$ depending only on $(\iv, \gam,  \beta_*,\alp)$ such that,
whenever parameters $(\eps, \delta_1, N_1)$ in Definition {\rm \ref{definition-10-3}}
are from $(0,\eps^{\rm (lb)}]\times (0, \delta_1^{\rm (w)}]\times[N_1^{\rm (w)}, \infty)$,
and $\delta_3\in(0, \bar{\delta}_3]$, then the following properties hold{\rm :}
For each $(u^{\sharp}, \beta^{\sharp})\in \ol{\mcl{K}}$,
a constant $\delta^{\sharp}>0$ can be chosen depending only
on $(\iv, \gam,  \beta_*, u^{\sharp}, \beta^{\sharp})$ so that,
if $(u,\beta)\in \ol{\Kext}$ satisfies \eqref{7-e9-new2015},
solution $\hat{\phi}$ of the boundary value problem \eqref{12-50} associated
with $(u,\beta)$ satisfies the estimate{\rm :}
\begin{equation}\label{est-hat-psi-new2015}
\|\hat{\phi}\|_{2,{\alp_1^*},\Om\setminus \oD_{\eps_0/10}}^{(2), {\rm (par)}}\le C \end{equation}
for $\Om=\Om(u,\beta)$,
where norm $\|\cdot\|_{2,\alp_1^*, \Om\setminus \oD_{\eps_0/10}}^{(2), {\rm (par)}}$
is defined by {\rm{Definition \ref{definition-parabolic-norm}}}.

\begin{proof}
The proof is divided into two steps.

\smallskip
{\textbf{1.}}
{\emph{Claim{\rm :} There exists a constant $C>0$ depending only on
$(\iv, \gam,  \beta_*)$ such that, for each $(u,\beta)\in\ol{\Kext}$,
$\hat{\phi}$ satisfies
\begin{equation}
\label{est-quad-upper-bd}
\hat{\phi}(x,y)\le C x^2\qquad\tx{in}\,\, \Om\cap\nD_{\eps_0}
\end{equation}
in the $(x,y)$--coordinates defined by \eqref{coord-n}.
}}

For the $(x,y)$--coordinates defined by \eqref{coord-n}, denote
\[
v(x,y):=\frac A2 x^2
\]
for a constant $A\ge \frac{2-\frac{\mu_0}{10}}{\gam+1}$ to be determined later,
where $\mu_0$ is from Definition \ref{definition-10-3}(iv-1).
For the elliptic cut-off $\zeta_1$ defined by \eqref{4-7},
$\zeta_1(\frac{v_x}{x})=\frac{2-\frac{\mu_0}{10}}{\gam+1}$.
By Lemma \ref{lemma-8-3} and \eqref{12-57},
equation $\mcl{N}_{(u,\beta)}(\hat{\phi})=0$ is rewritten
in the $(x,y)$--coordinates as
\begin{equation*}
  \mcl{N}^{\rm{polar}}_{(u,\beta)}(\hat{\phi})=0\qquad \tx{in} \,\,
\Om\cap \nD_{\eps_{\rm eq}/2}
\end{equation*}
for the nonlinear differential operator $\mcl{N}^{\rm{polar}}_{(u,\beta)}$ given
by \eqref{definition-nlop-N2}, where $\eps_{\rm eq}\in(0, \frac{\eps_0}{2})$ is from Lemma \ref{lemma-8-3}.

By $\zeta_1(\frac{v_x}{x})=\frac{2-\frac{\mu_0}{10}}{\gam+1}$ and \eqref{definition-nlop-N2}, we have
\begin{equation*}
  \mcl{N}^{\rm{polar}}_{(u,\beta)}(v)
  =Ax\Big(-(1-\frac{\mu_0}{10})+\frac{O_1^{\rm mod}}{x}+O_4^{\rm mod}\Big)\qquad \tx{in} \,\,
\Om\cap \nD_{\eps_{\rm eq}/2},
\end{equation*}
with $O_j^{\rm mod}=O_j^{\rm mod}(v_x, 0, x, y)$ for $j=1,4$.
It follows from \eqref{def-O-mod} that
$|\frac{O_1^{\rm mod}}{x}|+|O_4^{\rm mod}|\le C\sqrt{x}$ for $C>0$ depending only on $(\iv, \gam)$.
Therefore, there exists $\bar{\eps}\in(0, \frac 12 \min\{\eps_0, \eps_{\rm eq}, \bar{\eps}_{\rm{bc}}\})$
depending only on $(\iv, \gam)$ such that
\begin{equation*}
\begin{split}
\mcl{N}^{\rm{polar}}_{(u,\beta)}(v)
  \le Ax\big(-(1-\frac{\mu_0}{10})+C\sqrt{\bar{\eps}}\big)
  < -\frac{Ax}{2}\big(1-\frac{\mu_0}{10}\big)
  <0= \mcl{N}^{\rm{polar}}_{(u,\beta)}(\hat{\phi})
  \qquad \mbox{in $\Om\cap \nD_{\bar{\eps}}$}.
\end{split}
\end{equation*}
Note that $0<\mu_0<1$ by Definition \ref{definition-10-3}(iv-1) and Lemma \ref{lemma-est-nrsonic}.

On $\shock(u,\beta)\cap \nD_{\bar{\eps}}$, properties (f)--(g) of Lemma \ref{lemma-ext-M-property} imply that
\begin{equation*}
\begin{split}
\M_{(u,\beta)}(Dv, v, \bm\xi )&=\M_{(u,\beta)}(Dv, v, \bm\xi )-\M_{(u,\beta)}({\bf 0}, 0, \bm\xi )\\
&\le -\delta_{\rm bc}(Ax+\frac{A}{2}x^2)<0=\M_{(u,\beta)}(D\hat{\phi}, \hat{\phi},\bm\xi)
\end{split}
\end{equation*}
for constant $\delta_{\rm bc}>0$ from Lemma \ref{lemma-ext-M-property}(g).
On $\Wedge\cap \ol{\nD_{\bar{\eps}}}$, $\der_{\etan}v=\der_y v=0=\der_{\bf n_{\rm w}} \hat{\phi}$.
On $\rightsonic$, $v=0=\hat{\phi}$.

By \eqref{estimate1-iter-bv-sol} and Remark \ref{remark-a}(ii),
there exists a constant $\hat C>0$ depending only on $(\iv, \gam,  \beta_*)$ such that
$\hat{\phi}$  satisfies
\begin{equation}
\label{upperbd-hpsi-nr-so-new2015}
\hat{\phi}(x,y)\le \hat C x \qquad\mbox{on}\,\,  \Om\cap \nD_{\eps_0}.
\end{equation}
Choose $A=\max\{\frac{2\hat C}{\bar{\eps}}, \frac{2-\frac{\mu_0}{10}}{1+\gam}\}$ so that $v$ satisfies
\[
\hat{\phi}\le v\qquad \tx{on}\,\, \Om\cap \{x=\bar{\eps}\}.
 \]
By Lemmas \ref{lemma-6-1} and \ref{lemma-ext-M-property}, and the comparison principle, we have
\begin{equation}\label{label-a}
\hat{\phi}\le v\qquad\mbox{in}\,\, \Om\cap \nD_{\bar{\eps}}.
\end{equation}
In order to extend this result onto $\Om\cap \nD_{\eps_0}$, we adjust the choice of $A$ as
\begin{equation*}
A=\max\Big\{\frac{2\hat C}{\bar{\eps}}, \frac{2-\frac{\mu_0}{10}}{1+\gam}, \frac{2\hat C \eps_0}{\bar{\eps}^2}\Big\},
\end{equation*}
so that, from \eqref{upperbd-hpsi-nr-so-new2015},
\begin{equation}\label{label-b}
\hat{\phi}(x,y)\le \hat C\eps_0\le \frac A2 \bar{\eps}^2\le v(x,y)\qquad\,\,\mbox{in}\,\,\,
\Om\cap (\nD_{\eps_0}\setminus \nD_{\bar{\eps}}).
\end{equation}
Combining \eqref{label-a} with \eqref{label-b},
we obtain \eqref{est-quad-upper-bd} with $C=A$ for $A$ given above before \eqref{label-b}.

\medskip
{\textbf{2.}} By Definition \ref{definition-10-3}(iii) and Remark \ref{remark-a}(ii),
there exists a constant $l>0$ depending only on $(\gam, \iv)$ such that
\begin{equation}
\label{fn-lwbd-newcor}
\fshockn(x)\ge l\qquad\,\,\tx{on $[0, \eps_0]$.}
\end{equation}

By Remark \ref{remark-a}(ii), $\fshockn$ satisfies the estimate:
\begin{equation}\label{Fz0-estimate-nr-rsonic}
\|\fshockn\|_{2,\alp, (0,\eps_0)}^{(-1-\alp),\{0\}}\le
\|\hat f_{\mcl{N},0}\|_{C^3([0, \eps_0])}+CN_0.
\end{equation}
By \eqref{est-quad-upper-bd}, \eqref{fn-lwbd-newcor}--\eqref{Fz0-estimate-nr-rsonic},
Lemmas \ref{lemma-8-3}
and \ref{lemma-7-9}--\ref{lemma-ext-M-property},
the boundary value problem \eqref{12-50} associated with $(u,\beta)\in\ol{\Kext}$ satisfying \eqref{7-e9-new2015}
 satisfies all the conditions of Theorem \ref{elliptic-t7-CF2}.
 Therefore, we conclude from Theorem \ref{elliptic-t7-CF2} that,
 for each $\alp'\in(0,1)$,
there exists a constant $C_{\alp'}>0$ depending only on
$(\iv, \gam,  \beta_*, \alp')$ such that $\hat{\phi}$ satisfies
\begin{equation}
\label{hphi-par-est-nr-rsonic-2016sep}
\|\hat{\phi}\|_{2,\alp', \Om\cap \nD_{\eps_0}}^{(2),{\rm (par)}}\le C_{\alp'}.
\end{equation}

Finally, \eqref{est-hat-psi-new2015} is obtained by combining estimate \eqref{hphi-par-est-nr-rsonic-2016sep} with Lemma \ref{proposition-existence-hpsi-bvp}.
\end{proof}
\end{lemma}

As pointed out earlier, $\leftsonic$ defined in Definition \ref{definition-domains-np}
depends continuously on $\beta\in[0, \frac{\pi}{2})$.
Therefore, the pseudo-subsonic region $\Om(u,\beta)$ associated with $(u,\beta)\in\ol{\Kext}$ depends continuously
on $(u,\beta)$.
In particular, $\Om(u,\beta)\cap \oD_{\eps_0}$ changes from a rectangular domain to a triangular domain as $\beta$
increases from $\beta<\betasonic$ to $\beta>\betasonic$.
Furthermore, the ellipticity of equation $\mcl{N}_{(u,\beta)}(\hat{\phi})=0$ near $\leftsonic$ changes as $\beta$ varies.
For that reason, the {\it a priori} estimate of a solution $\hat{\phi}$ of the boundary value problem \eqref{12-50}
is given for the three cases separately:

(i) $\beta< \betasonic$;

(ii) $\beta\ge\betasonic$ close to $\betasonic$;

(iii) $\beta>\betasonic$ away from $\betasonic$.

\begin{lemma}[Estimates of $\hat{\phi}$ near $\leftsonic$]
\label{lemma-est-hpsi-nrlsonic}
Let $\eps^{\rm (par)}$ be from Lemma {\rm \ref{proposition-est-hpsi-new2015}}.
There exist $\eps^{\mcl{O}}\in(0,\eps^{\rm (par)}]$ and $\delta_1^{\rm(E)}$ depending only on $(\iv, \gam, \beta_*)$
such that,
whenever parameters $(\eps, \delta_1, \delta_3, N_1)$ in Definition {\rm \ref{definition-10-3}} are chosen
as in Lemma {\rm \ref{proposition-est-hpsi-new2015}},
and $(\eps, \delta_1)$ further satisfy
\begin{equation*}
0<\eps<\eps^{\mcl{O}},\quad 0<\delta_1\le \delta_1^{(E)},
\end{equation*}
then, for each $(u^{\sharp}, \beta^{\sharp})\in\ol{\mcl{K}}$,
there is a constant $\delta^{\sharp}$ depending on $(\iv, \gam,  \beta_*, \delta_2, \delta_3, u^{\sharp}, \beta^{\sharp})$
so that, if $(u,\beta)\in\ol{\Kext}$ satisfies \eqref{7-e9-new2015}, then the following properties hold{\rm :}

\smallskip
\begin{itemize}
  \item[(i)] If $\beta\in[0, \betasonic)$, for each $\alp'\in(0,1)$, there exist constants
  $\hat{\eps}_p\in(0,\eps_0]$ and $C_{\alp'}>0$
  depending only on $(\iv, \gam,  \beta_*, \alp')$
  such that solution $\hat{\phi}\in C^2(\Om)\cap C^1(\ol{\Om})$ of
  the boundary value problem \eqref{12-50} associated with $(u,\beta)$ satisfies
      \[
      \|\hat{\phi}-(\leftvphi-\rightvphi)\|^{(2),{\rm (par)}}_{2,\alp', \Om\cap \oD_{\hat{\eps}_p}}\le C_{\alp'};
      \]

  \item[(ii)] There exists a constant $\hat{\delta}\in(0, \beta_*-\betasonic)$ depending only
  on $(\iv, \gam,  \beta_*)$ such that, if $\beta\in[\betasonic, \betasonic+\hat{\delta}]$,
  then, for each $\alp'\in(0,1)$, there exist constants $\hat{\eps}_p\in(0,\eps_0]$
  depending on $(\iv, \gam, \beta_*)$ and $C_{\alp'}>0$
  depending only on $(\iv, \gam,  \beta_*, \alp')$ so that $\hat{\phi}$ satisfies
      \begin{equation*}
      \begin{split}
      &\|\hat{\phi}-(\leftvphi-\rightvphi)\|_{C^{2,\alp'}(\Om\cap \oD_{\hat{\eps}_p})}\le C_{\alp'},\\
      &D^m(\hat{\phi}-\leftvphi+\rightvphi)(P_{\beta})=0\qquad
      \tx{for}\,\,m=0,1,2,
      \end{split}
      \end{equation*}
      where $P_{\beta}$ is defined in Definition {\rm \ref{definition-domains-np}}{\rm ;}

\smallskip
  \item[(iii)] There exist constants $\hat{\alp}\in(0, \frac 13)$ depending only on $(\iv, \gam,  \beta_*)$
  and $C>0$
  depending only on $(\iv, \gam,  \beta_*)$ so that,
  if $\beta\in[\betasonic+\frac{\hat{\delta}}{2}, \beta_*]$,
    then $\hat{\phi}$ satisfies
      \begin{align}
      \label{estimate-corner-subsonic}
      &\|\hat{\phi}-(\leftvphi-\rightvphi)\|_{2,\hat{\alp}, \Om\cap \oD_{\eps_0}}^{(-1-\hat{\alp}), \{P_{\beta}\}}
      \le C,\\
      \label{vanishing-corner-subsonic}
      &D^m(\hat{\phi}-\leftvphi+\rightvphi)(P_{\beta})=0\qquad
      \tx{for}\,\,m=0,1.
      \end{align}
\end{itemize}

\begin{proof} We divide the proof into two steps.

\smallskip
{\textbf{1.}} {\emph{Assertion {\rm (i)}}}:
Owing to Remark \ref{remark-lsonic},
we need to consider two cases separately: (i) $\beta<\betasonic$ away from $\betasonic$ and (ii) $\beta<\betasonic$ close to $\betasonic$.

By Lemma \ref{lemma-10-1}(e), \eqref{h2-st-y}, \eqref{12-16-mod}, Proposition \ref{lemma-12-3},
and Definition \ref{definition-10-3}(iii), there exist $\hat{\eps}\in(0, \eps^{(\rm par)}]$ and
$\hat{\sigma}_1\in(0, \frac{\betasonic}{10})$  so that, for any $(u,\beta)\in\ol{\Kext}$, it holds that,
if $\sigma\in(0, \hat{\sigma}_1]$, then
we can fix  $\hat{m}>1$ depending only on $(\iv, \gam, \sigma)$ and $\hat{k}>1$ depending only on $(\iv, \gam)$  such that

\smallskip
\begin{itemize}
\item[(a)] if $0\le \beta\le \betasonic-\frac{\sigma}{2}$, then
\begin{equation}
\label{domain-nr-lsonic1}
\{0<x<2\hat{\eps}, 0<y< \frac{1}{2\hat{m}}\}\subset \Om\cap \oD_{2\hat{\eps}}\subset \{0<x<2\hat{\eps}, 0<y<2\hat{m}\};
\end{equation}

\item[(b)] if $\betasonic-\sigma\le \beta < \betasonic$, then
\begin{equation}
\label{domain-nr-lsonic2}
\{0<x<2\hat{\eps}, 0<y< y_{\lefttop}+\frac{x}{2\hat{k}}\}\subset \Om\cap \oD_{2\hat{\eps}}\subset \{0<x<2\hat{\eps}, 0<y<y_{\lefttop}+2\hat{k}x\}.
\end{equation}
\end{itemize}

For a fixed $\sigma\in(0, \hat{\sigma}_1]$, suppose that $0\le \beta\le \betasonic-\frac{\sigma}{2}$.
Let $\hat{\psi}$ be given by \eqref{definition-hpsi-2018}. By Lemma \ref{lemma-positivity-hatphi-2018}, we have
\begin{equation}
\label{iter-bvp-sol-nrlsonic1}
\hat{\psi}>0 \qquad \tx{ in $\Om\cap \oD_{\eps/2}$,}
\end{equation}
provided that $(u,\beta)\in \ol{\Kext}$ satisfies (4.5.2) for $\delta^{\sharp}>0$ from Lemma 4.42.

Owing to \eqref{7-b7}, if condition \eqref{eps-condition-vphi-betastar} holds,
then we can repeat Step 1 in the proof of Lemma \ref{proposition-est-hpsi-new2015}  to obtain
\begin{equation}
\label{iter-bvp-sol-nrlsonic2}
\hat{\psi}(x,y)\le Cx^2\qquad\tx{in $\Om\cap \oD_{\hat{\eps}_0}$} \,\,\,\tx{for $\hat{\eps}_0:=\min\{\eps_0, \frac{\leftch}{\bar k}\}$}
\end{equation}
for $C>0$ depending only on $(\iv, \gam, \beta_*)$,
where the $(x,y)$--coordinates are given by \eqref{coord-o},
and $\leftch$ and $\bar{k}$ are given by Definition \ref{definition-Qbeta} and \eqref{7-b7}, respectively.
Repeating Step 2 in the proof of Lemma \ref{proposition-est-hpsi-new2015}  with \eqref{iter-bvp-sol-nrlsonic1}--\eqref{iter-bvp-sol-nrlsonic2}
and $\fshocko$ given by \eqref{shock-ftns-nr-sonic-2016sep}, and using \eqref{domain-nr-lsonic1},
we can show that, for each $\alp'\in(0,1)$, there exists a constant $C_{\alp'}>0$ depending only on
$(\iv, \gam,  \beta_*, \alp')$ such that
\begin{equation*}
 \|\hat{\phi}-(\leftvphi-\rightvphi)\|^{(2),{\rm (par)}}_{2,\alp', \Om\cap \oD_{\hat{\eps}_0}}=
 \|\hat{\psi}\|_{2,\alp', \Om\cap \oD_{\hat{\eps}_0}}^{(2),{\rm (par)}}
\le C_{\alp'}.
\end{equation*}

Next, suppose that $\betasonic-\sigma\le \beta<\betasonic$.
In this case, we need to combine two estimates: (i) in  $\Om\cap \{x<y_{\lefttop}^2\}$
and (ii) in $\Om\cap \{x>\frac{y_{\lefttop}^2}{10}\}$ near $\leftsonic$.

In $\Om\cap \{x<y_{\lefttop}^2\}$, we repeat the argument of Step 2  in the proof of Lemma \ref{proposition-est-hpsi-new2015}
to obtain
\begin{equation*}
 \|\hat{\phi}-(\leftvphi-\rightvphi)\|^{(2),{\rm (par)}}_{2,\alp', \Om\cap \oD_{y_{\lefttop}^2}}=
 \|\hat{\psi}\|_{2,\alp', \Om\cap \oD_{y_{\lefttop}^2}}^{(2),{\rm (par)}}
\le C_{\alp'}
\end{equation*}
for each $\alp'\in(0,1)$, where $C^{\alp'}>0$ is given, depending only on $(\iv, \gam, \beta_*,\alp')$.

In $\Om\cap \{x>y_{\lefttop}^2\}$ near $\leftsonic$, we adjust the argument in Step 2 in the proof of Proposition \ref{proposition-sub8}
to show that there exist sufficiently small constants $\bar{\sigma}\in(0, \sigma_1]$
and $\eps^*\in(0, \hat{\eps}_0]\cap  (0, \eps^{\rm{(par)}}]$ depending only on $(\iv, \gam, \beta_*)$ so that $\hat{\psi}$ satisfies
\begin{equation*}
\hat{\psi}(x,y)\le Cx^4\qquad\tx{in $\Om\cap \oD_{\eps^*}\cap\{x>\frac{y_{\lefttop}^2}{10}\}$}
\end{equation*}
for $C>0$ depending only on $(\iv, \gam, \beta_*)$.
For $\fshocko$ defined by (4.3.20) and $z_0=(x_0, y_0)\in \Om\cap \mcl{D}^{\mcl{O}}_{\eps^*}\cap\{x>\frac{y_{\lefttop}^2}{5}\}$,
we define $F^{(z_0)}(S)$ by (3.5.39) given in the proof of Proposition 3.32.
By Remark 4.21(i)--(ii), $F^{(z_0)}$ satisfies
\begin{equation*}
\|F^{(z_0)}\|_{C^2([-1,1])}\le CN_0\sqrt{x_0}
\end{equation*}
for $C>0$ depending only on $(\iv, \gam,\alp)$.
Then we apply Theorem C.6 and adjust the later part of Step 4 in the proof of Proposition 3.32 to
conclude that
\begin{equation*}
 \|\hat{\phi}-(\leftvphi-\rightvphi)\|^{(2),{\rm (par)}}_{2,\alp', \Om\cap \oD_{\eps^*}}=
 \|\hat{\psi}\|_{2,\alp', \Om\cap \oD_{\eps^*}}^{(2),{\rm (par)}}
\le C_{\alp'}
\end{equation*}
for each $\alp'\in(0,1)$, where $C_{\alp'}>0$ is given, depending only on $(\iv, \gam, \beta_*,\alp')$,
provided that $\sigma\in(0, \bar{\sigma}]$.

The proof of assertion (i) is completed.

\medskip
{\textbf{2.}} {\emph{Assertions {\rm (ii)} and {\rm (iii)}}}:
Assertion (ii) can be proved in a way similar to Proposition \ref{proposition-sub9}.
Estimate \eqref{estimate-corner-subsonic} in assertion (iii) directly follows
from Proposition \ref{proposition-app-C-wp2}.

For $\beta\ge \betac+\frac{\hat{\delta}}{2}$, \eqref{estimate1-iter-bv-sol} implies that
\begin{equation}
\label{van-0th-corner-subsonic}
(\hat{\phi}-\leftphi)(P_{\beta})=0
\end{equation}
for $\leftphi=\leftvphi-\rightvphi$.
By Lemma  \ref{lemma-ext-M-property}(f) and \eqref{van-0th-corner-subsonic}, $\hat{\phi}$ satisfies
\begin{equation*}
\int_0^1 \mcl{M}_{(u,\beta)}(tD\hat{\phi}+(1-t)D\phi_0, t\hat{\phi}+(1-t)\phi_0, \bmxi)\,\dd t\cdot D(\hat{\phi}-\phi_0)=0
\qquad\tx{at $\bmxi=P_{\beta}$}.
\end{equation*}
By \eqref{12-67}, \eqref{7-d7}, \eqref{7-d8}, and Lemma  \ref{lemma-ext-M-property},
\begin{equation*}
|\der_{p_1}\mcl{M}_{(u,\beta)}({\bf p}, z, P_{\beta})-\der_{p_1}g^{\rm sh}(D\leftvphi(P_{\beta}), \leftvphi(P_{\beta}), P_{\beta})|\le C\sqrt{\delta_1}
\end{equation*}
for some $C>0$ depending only on $(\iv, \gam, \beta_*)$.
This inequality, combined with Lemma \ref{lemma-10-3},
implies that, if $\delta_1>0$ is chosen small, depending only on $(\iv, \gam, \beta_*)$,
then the boundary conditions: $\mcl{M}_{(u,\beta)}(D\hat{\phi}, \hat{\phi}, \bmxi)=0$
on $\shock$ and $\hat{\phi}_{\xi_2}=0$ on $\Wedge$ are functionally independent at $P_{\beta}$ so that
\begin{equation*}
D(\hat{\phi}-\leftphi)(P_{\beta})=0.
\end{equation*}

In proving assertions (i)--(iii), all the required properties of $\mcl{N}_{(u,\beta)}$
and $\mcl{M}_{(u,\beta)}$ are provided by Lemmas \ref{lemma1-coeff-iter-eqn-new}, \ref{lemma2-coeff-iter-eqn-new}, \ref{lemma-6-1},
and \ref{lemma-7-7}--\ref{lemma-ext-M-property}.
\end{proof}
\end{lemma}

\begin{corollary}
\label{corollary-hatu-2018}
In {\rm{Definition \ref{definition-10-3}}}, choose parameters $(\alp, \eps, \delta_1, \delta_3, N_1)$ as follows{\rm :}

\smallskip
\begin{itemize}
\item[(i)] For $\bar{\alp}$, $\alp_1$, and $\hat{\alp}$ from {\rm{Lemmas \ref{proposition-existence-hpsi-bvp},
\ref{proposition-est-hpsi-new2015}, and \ref{lemma-est-hpsi-nrlsonic}}}, respectively, choose
\begin{equation*}
\alp=\frac 12 \min\{\bar{\alp}, \alp_1, \hat{\alp}\};
\end{equation*}

\item[(ii)]
Choose $(\eps, \delta_1,  N_1)$ to satisfy
\begin{equation*}
(\eps, \delta_1, N_1)\in
(0,\eps^{\mcl{O}}]\times (0, \delta_1^{\rm (w)}]\times[N_1^{\rm (a)}, \infty)
\end{equation*}
for $N_1^{\rm (a)}\in[N_1^{\rm (w)}, \infty)$ from Corollary {\rm \ref{corollary-admisol-iterset}}{\rm ;}

\smallskip
\item[(iii)] For $(\delta_1, N_1)\in (0, \delta_1^{\rm (w)}]\times [N_1^{\rm (a)}, \infty)$,
denote $\bar{\delta}:=\frac{\delta_1}{2N_1^2}\delta_2$, where
$\delta_2>0$ is a parameter to be determined later.
Choose $\delta_3$ to satisfy
\begin{equation*}
\delta_3\in(0, \bar{\delta}_3].
\end{equation*}
\end{itemize}
Under the choices of parameters $(\alp, \eps, \delta_1, \delta_3, N_1)$ above,
there exists a constant $C>0$ depending only on $(\iv, \gam,  \beta_*)$ such that,
for each $(u,\beta)\in\ol{\mcl{\Kext}}$, denoting  the unique solution of the boundary value problem \eqref{12-50}
associated with $(u,\beta)$ by $\hat{\phi}\in C^2(\Om(u,\beta))\cap C^1(\ol{\Om(u,\beta)})$
and defining  $\hat u:\ol{\iter}\rightarrow \R$ by \eqref{definition-sol-ubvp}, then
\begin{equation}
\label{7-e8}
\|\hat u\|_{2,2\alp,\iter}^{(*,1)}\le C.
\end{equation}
\end{corollary}

\begin{proof}
By the choice of parameters $\alp\in(0, \frac 16)$ and $(\eps, \delta_1, \delta_3, N_1)$,
estimate \eqref{7-e8} follows from Lemmas \ref{proposition-est-hpsi-new2015}--\ref{lemma-est-hpsi-nrlsonic} by
repeating the argument in the proof of Proposition \ref{proposition-unif-est-u-new}.
\end{proof}

\begin{proposition}
\label{lemma-7-10}
Under the choices of parameters $(\alp, \eps, \delta_1, \delta_3, N_1)$ as in Corollary {\rm \ref{corollary-hatu-2018}},
the iteration set $\mcl{K}$ defined in Definition {\rm \ref{definition-10-3}}
is relatively open in $C^{2,\alp}_{(*,\alp_1)}(\iter)\times [0,\beta_*]$.
\end{proposition}

\begin{proof}
We have shown in Lemma \ref{lemma-Kext-open}
that $\Kext$ is relatively open in $C^{2,\alp}_{(*,\alp_1)}(\iter)\times [0,\beta_*]$.
Therefore, it remains to check that property (vii) of Definition \ref{definition-10-3} defines a relatively open subset
of $C^{2,\alp}_{(*,\alp_1)}(\iter)\times [0,\beta_*]$ under the choice of $\delta_3$ given by (iii) in the statement
of Corollary \ref{corollary-hatu-2018}.

Suppose that this is not true. Then there exist $(u^{\sharp}, \beta^{\sharp}) \in \mcl{K}$
and a sequence $\{(u_n,\beta_n)\}_{n=1}^{\infty}\subset \Kext$ such that
\begin{equation*}
\lim_{n\to \infty}\|u_n-u^{\sharp}\|_{2,\alp/2, \iter}^{(*,\alp_1)}+|\beta_n-\beta^{\sharp}|=0,
\qquad
\|\hat{u}_n-u_n\|_{2,\alp/2, \iter}^{(*,\alp_1)}\ge \delta_3\qquad\tx{for all $n\in \mathbb{N}$},
\end{equation*}
where each $\hat{u}_n$ for $n\in \mathbb{N}$ is given by \eqref{12-65} for $(u,\beta)=(u_n,\beta_n)$.

Let $\hat{u}^{\sharp}$ be given by \eqref{definition-sol-ubvp} with $(u,\beta)=(u^{\sharp}, \beta^{\sharp})$,  and denote
\begin{equation*}
\delta^{\sharp}:=\frac{\delta_3-\|\hat u^{\sharp}-u^{\sharp}\|_{2,\alp/2, \iter}^{(*,\alp_1)}}{10}.
\end{equation*}
By \eqref{12-65}, it holds that $\delta^{\sharp}>0$.
Therefore, we can choose $n^{\sharp}\in \mathbb{N}$ sufficiently large
such that $\|u_n-u^{\sharp}\|_{2,\alp/2, \iter}^{(*,\alp_1)}+|\beta_n-\beta^{\sharp}|\le \delta^{\sharp}$ for all $n\ge n^{\sharp}$.
Then we have
\begin{equation*}
\|\hat u_n-\hat u^{\sharp}\|_{2,\alp/2, \iter}^{(*,\alp_1)}\ge 9\delta^{\sharp}\qquad\tx{for all $n\ge n^{\sharp}$}.
\end{equation*}

By Corollary \ref{corollary-hatu-2018}, $\{\hat{u}_n\}$ is bounded in $C^{2,2\alp}_{(*,1)}(\iter)$.
It is noted in Definition \ref{pre-definition-new} that $C^{2,2\alp}_{(*,1)}(\iter)$ is compactly
embedded into $C^{2,\alp}_{(*,\alp_1)}(\iter)$.
Therefore,  $\{\hat{u}_n\}$ has a subsequence $\{\hat{u}_{n_j}\}$ that converges in $C^{2,\alp}_{(*,\alp_1)}(\iter)$
to a function $\hat{u}^*\in C^{2,\alp}_{(*,\alp_1)}(\iter)$ so that
\begin{equation}
\label{nonzero-distance-lmt}
\|\hat u^*-\hat u^{\sharp}\|_{2,\alp/2, \iter}^{(*,\alp_1)}\ge 9\delta^{\sharp}.
\end{equation}
Define
\begin{equation*}
\hat{\phi}^*:=\hat u^*\circ \mathfrak{F}^{-1}_{(u^{\sharp}, \beta^{\sharp})}-\rightvphi+{\vphi_{\beta^\sharp}^*}
\qquad\tx{in $\mathfrak{F}^{-1}_{(u^{\sharp}, \beta^{\sharp})}(\ol{\iter})=\ol{\Om(u^{\sharp}, \beta^{\sharp})}$}.
\end{equation*}
By Lemma \ref{remark-prop-coeff-iter-bvp}, $\hat{\phi}^*$ solves the nonlinear boundary value problem \eqref{12-50}
associated with $(u^{\sharp}, \beta^{\sharp})$. Then the uniqueness of solutions of \eqref{12-50}
stated in Lemma \ref{proposition-existence-hpsi-bvp} implies that $\hat{u}^*=u^{\sharp}$,
which is in contradiction to \eqref{nonzero-distance-lmt}.
Therefore, we conclude that property (vii) of Definition \ref{definition-10-3} defines a relatively open subset
of $C^{2,\alp}_{(*,\alp_1)}(\iter)\times [0,\beta_*]$ under the choice of $\delta_3$
given by (iii) in the statement of Corollary \ref{corollary-hatu-2018}.
\end{proof}

\begin{remark}
\label{remark-parameter-chocie1}
In Proposition {\rm \ref{lemma-7-10}},
the choice of $(\alp, \eps, \delta_1, N_1)$ depends only on $(\iv, \gam,  \beta_*)$,
and the choice of $\delta_3$ depends only on $(\iv, \gam,  \beta_*,\delta_1,  \delta_2, N_1)$,
where parameter $\delta_2$ is to be determined later.
\end{remark}



\chapter{Existence of Admissible Solutions up to $\betadet$\\ -- Proof of Theorem 2.31}
\label{sec-proof-theorem5-case1}
\numberwithin{equation}{section}

Fix $\gam \ge 1$, $\iv>0$, and $\beta_*\in(0, \betadet)$.
For the iteration set $\mcl{K}$ defined in Definition \ref{definition-10-3}, define
\begin{equation*}
\mcl{K}(\beta):=\{u\in C^{2,\alp}_{(*,\alp_1)}(\iter)\,:\,(u,\beta)\in \mcl{K}\}
\qquad\tx{for each $\beta\in [0,\beta_*]$}.
\end{equation*}
In this chapter, we define an iteration map $\mcl{I}: \ol{\mcl{K}}\rightarrow C^{2,\alp}_{(*,\alp_1)}(\iter)$ with the following properties:

\smallskip
\begin{itemize}
\item[(i)] For each $\beta\in[0, \beta_*]$, there exists $u\in \mcl{K}(\beta)$ such that $\mcl{I}(u,\beta)=u$;

\smallskip
\item[(ii)] If $\mcl{I}(u,\beta)=u$, then $\vphi$ given by \eqref{10-12}
yields an admissible solution corresponding to $(\iv,\beta$).
\end{itemize}

\section{Definition of the Iteration Map}
\label{section-itr-mp}
Let parameters $(\alp, \eps, \delta_1, \delta_3, N_1)$ in Definition \ref{definition-10-3}
be fixed as  in Proposition \ref{lemma-7-10}.

In order to define an iteration map satisfying (i)--(ii) stated above,
and to employ the Leray-Schauder degree argument for proving the existence
of a fixed point of $\mcl{I}(\cdot,\beta)$ in $\mcl{K}(\beta)$ for all $\beta\in(0,\betadet)$,
we require the compactness of $\mcl{I}$.

\smallskip
For each $(u,\beta)\in \ol{\mcl{K}}$,
let $(\gshock^{(u,\beta)}, \shock(u,\beta), \Om(u,\beta), \vphi^{(u,\beta)})$ be defined
by Definition \ref{definition-Gset-shocks-new},
and denote them as $(\gshock, \shock, \Om, \vphi)$.
For such a function $\gshock$, we define $(\mcl{G}_1^{\beta}, G_{2,\gshock})$
by \eqref{12-16-mod} and \eqref{7-b9}, respectively.
Let $\hat{\phi}\in C^2(\Om)\cap C^1(\ol{\Om})$ be the unique solution of the boundary value problem \eqref{12-50}
associated with $(u,\beta)$.
Then function $\hat u:\ol{\iter}\rightarrow \R$ is given by \eqref{definition-sol-ubvp},
and function $\hat{\vphi}=\hat{\vphi}^{(u,\beta)}$ is given by
\begin{equation}
\label{definition-hat-vphi-ubeta}
\hat{\vphi}^{(u,\beta)}=\vphi_{\beta}^*+\hat u\circ \mathfrak{F}_{(u,\beta)}^{-1}
\end{equation}
for $\vphi_{\beta}^*$ given by \eqref{12-32}.

Next, we define functions $w$, $w_{\infty}$, and $\hat w$ by
\begin{equation}
\label{def-w-infty}
\begin{split}
&w(s,t'):=(\vphi-\vphi_{\beta}^*)\circ (\mcl G_1^{\beta})^{-1}(s,t'),\\
&w_{\infty}(s,t'):=(\ivphi-\vphib^*)\circ(\mcl G_1^{\beta})^{-1}(s,t'),\\
&\hat{w}(s,t'):=(\hat{\vphi}-\vphib^*)\circ (\mcl G_1^{\beta})^{-1}(s,t').
\end{split}
\end{equation}

\begin{lemma}
\label{lemma-winfty-new2015}
For each $\beta\in[0, \betad^{(\iv)}]$,
there exists a unique function $\mathfrak{g}_{\beta}:[-1, 1]\rightarrow \R_+$ such that

\smallskip
\begin{itemize}
\item[(a)] $w_{\infty}(s, \mathfrak{g}_{\beta}(s))=0$ for all $s\in[-1,1]${\rm ;}

\smallskip
\item[(b)] $\{(s, \mathfrak{g}_{\beta}(s))\,:\, s\in(-1,1)\} \subset \mcl G_1^{\beta}(Q^{\beta})$
for $\Qbeta$ defined in Definition {\rm \ref{definition-Qbeta}(iii)}{\rm ;}

\smallskip
\item[(c)] $\|\mathfrak{g}_{\beta}\|_{C^3([-1, 1])}\le C$ for $C>0$ depending only on $(\gam, \iv)$.
\end{itemize}

\begin{proof}
By property (iii) stated right after Definition \ref{definition-smooth-connection}, the set:
$$
\{(s,t')\,:\,w_{\infty}(s,t')=0\}
$$
is contained in $\mcl{G}_1^{\beta}(\Qbeta)$.
Then the existence and uniqueness of $\mathfrak{g}_{\beta}$ satisfying  statements (a)--(b) follow from Lemma \ref{lemma-7-2},
combined with the implicit function theorem.
Statement (c) is obtained from Lemma \ref{lemma-7-2} and the smoothness of $\ivphi-\vphib^*$, owing to \eqref{12-32}.
\end{proof}
\end{lemma}

For each $(u,\beta)\in\ol{\mcl{K}}$, $\gshock :[-1, 1]\rightarrow \overline{\R_+}$ 
is in $C^{0,1}([-1,1])$ and satisfies $\gshock>0$ on $(-1, 1)$.
Define
\begin{equation}
\label{definition-Rg-Sg-2018}
\begin{split}
&R_{\gshock}:=\{(s, t')\in \R^2\,:\,-1<s<1, 0<t'<\gshock(s)\},\\[1mm]
&\Sigma_{\gshock}:=\{(s, \gshock(s))\,:\,-1<s<1\}.
\end{split}
\end{equation}
Note that $w$ and $\hat w$ are defined in $R_{\gshock}$, and $w_{\infty}$ is defined in ${R}_{\infty}:=(-1,1)\times \R_+$.

In order to define an iteration map $\mcl{I}$, the first step is to introduce
an extension of $\hat w$ onto $R_{(1+\kappa)\gshock}$ for some $\kappa\in(0, \frac 13]$.

\smallskip
\begin{lemma}[Regularized distance]
\label{lemma-reg-distance}
Let $R_{\infty}:=(-1,1)\times \overline{\R_+}$.
For each $g\in C^{0,1}([-1,1])$ satisfying
\begin{equation}
\label{positivity-g-reg-dist}
g>0\qquad\tx{on $(-1, 1)$},
\end{equation}
define
\begin{equation}
\label{definition-Rg-Sg}
R_{g}:=\{(s, t')\in \R^2\,:\,-1<s<1, 0<t'<g(s)\},\quad
\Sigma_{g}:=\{(s, g(s))\,:\,-1<s<1\}.
\end{equation}
Then there exists a function $\delta_{g}\in C^{\infty}(\ol{R_{\infty}}\setminus \ol{R_g})$, {\emph{the regularized distance}},
such that

\smallskip
\begin{itemize}
\item[(i)] For all ${\bf x}=(s,t')\in \ol{R_{\infty}}\setminus \Sigma_g$,
\begin{equation*}
\frac 12 {\rm dist}({\bf x}, \Sigma_g)\le \delta_g({\bf x})\le \frac 32 {\rm dist}({\bf x},\Sigma_g).
\end{equation*}

\item[(ii)] For all ${\bf x}=(s,t')\in \ol{R_{\infty}}\setminus \Sigma_g$,
\begin{equation*}
|D^m \delta_g({\bf x})|\le C(m)\big({\rm dist}({\bf x},\Sigma_g) \big)^{1-m} \qquad\mbox{for $m=1,2,3,\cdots$},
\end{equation*}
where $C(m)$ depends only on $m$.

\smallskip
\item[(iii)] There exists $C_*>0$ depending only on ${\rm Lip}[g]$ such that
\begin{equation*}
\delta_g({\bf x})\ge C_*(t'-g({s})) \qquad \tx{for all ${\bf x}\in\ol{R_{\infty}}\setminus \ol{R_g}$}.
\end{equation*}

\item[(iv)]  Suppose that $g_i\in C^{0,1}([-1,1])$ and $g\in C^{0,1}([-1,1])$ satisfy \eqref{positivity-g-reg-dist} and
\begin{equation*}
\|g_i\|_{C^{0,1}([-1,1])}\le L\qquad\tx{for all $i\in \mathbb{N}$}
\end{equation*}
for some constant $L>0$.
If $\{g_i(s)\}_{i\in \mathbb{N}}$ converges to $g(s)$ uniformly on $[-1,1]$,
then $\{\delta_{g_i}({\bf x})\}_{i\in \mathbb{N}}$ converges to $\delta_g({\bf x})$
in $C^m(K)$ for any $m=0,1,2,\cdots$, and any compact set $K\subset \ol{R_{\infty}}\setminus \ol{R_g}$.

\smallskip
\item[(v)]
For $C_*$ from {\rm (iii)}, define
\begin{equation}
\label{definition-delg-star}
\delta^*_g({\bf x}):=\frac{2}{C_*}\delta_g({\bf x}).
\end{equation}
Then there exists $\kappa\in(0, \frac 13]$ depending only on ${\rm Lip}[g]$  such that, for each ${\bf x}=(s,t')\in R_{(1+\kappa)g}\setminus \ol{R}_g$,
\begin{equation*}
(s,t'-\lambda \delta_g^*({\bf x}))\in \{s\}\times [\frac{g(s)}{3}, g(s)-(t'-g(s))]\Subset R_g\qquad\tx{for all $\lambda\in[1,2]$.}
\end{equation*}

\item[(vi)]  There exist constants $C_*>0$ and $\kappa\in(0, \frac 13]$ depending only on $(\gam, \iv, \beta_*)$ such that,
for each $(u,\beta)\in\ol{\Kext}$, the regularized distance $\delta^{(u,\beta)}_{\gshock}$
can be given so that properties {\rm (i)}--{\rm (iii)} and {\rm (v)} stated above are satisfied.

\smallskip
\item[(vii)]
If $\{(u_j, \beta_j)\}_{j=1}^n\subset \ol{\Kext}$ converges to $(u,\beta)$ in $C^{2,\alp}_{(*, \alpha_1)}(\iter)\times [0,\beta_*]$,
then $\delta_{\gshock}^{(u_j, \beta_j)}$ converges to $\delta_{\gshock}^{(u,\beta)}$ in $C^m(K)$
for any $m=0,1,2,\cdots$, and any compact subset $K\subset \ol{R_{\infty}}\setminus \ol{R_{\gshock}^{(u,\beta)}}$.
\end{itemize}

\begin{proof}
Statements (i)--(iv) of this lemma follow directly from \cite[Lemma 13.9.1]{CF2}.
Statement (v) can be
verified by using statement (iii).
We refer to  \cite[Lemma 13.9.4]{CF2} for a proof of statement (v).
Finally, statements (i)--(v), combined with (d) and (g)--(h) of Lemma \ref{lemma-7-4}
and (i) of Remark \ref{remark-a}, lead to statements (vi) and (vii).
\end{proof}
\end{lemma}

By \cite[Lemma 13.9.2]{CF2}, there exists a function $\Psi\in C^{\infty}_c(\R)$ satisfying that
\begin{equation}
\label{definition-Psi-ext-2018}
\begin{split}
&{\rm supp} \Psi \subset [1,2],\\[1mm]
&\int_{-\infty}^{\infty} \Psi(y) \,\dd\lambda=1,\qquad \int_{-\infty}^{\infty} \lambda^m\Psi(\lambda)\,\dd\lambda=0\quad\tx{for $m=1,2$}.
\end{split}
\end{equation}
For a function $g\in C^{0,1}([-1,1])$ satisfying \eqref{positivity-g-reg-dist},
let $R_g$ and $\delta_g^*$ be given by \eqref{definition-Rg-Sg} and \eqref{definition-delg-star}, respectively.
Let $\kappa\in(0, \frac 13]$ be fixed depending on ${\rm Lip}[g]$ to satisfy Lemma \ref{lemma-reg-distance}(v).
For a function $v\in C^0(\overline{R_g})\cap C^2(R_{g}\cup \Sigma_{g})$, we define its extension $\mcl{E}_{g}(v)$ onto $R_{(1+\kappa)g}$ by
\begin{equation}
\label{definition-extension-operator-2018}
\mcl{E}_{g}(v)({\bf x})=
\begin{cases}
v({\bf x})\quad&\tx{for ${\bf x}=(s,t')\in \ol{R_g}$},\\[1mm]
\int_{-\infty}^{\infty} v\left(s,t'-\lambda \delta_g^*({\bf x})\right)\,\Psi(\lambda)\,\dd\lambda\quad
&\tx{for ${\bf x}\in R_{(1+\kappa)g}\setminus \ol{R_g}$}.
\end{cases}
\end{equation}

\begin{definition}[Extension map]
\label{definition-extension-map}

For each $(u,\beta)\in\ol{\Kext}$,  let $g$ denote $\gshock^{(u,\beta)}$,
and let $\delta_g$ be the regularized distance given in Lemma {\rm \ref{lemma-reg-distance}}.
For constant $C_*>0$ from Lemma {\rm \ref{lemma-reg-distance}(vi)},
let $\delta^*_g$ be given by \eqref{definition-delg-star}.
Let $\kappa\in(0, \frac 13]$ be from Lemma {\rm \ref{lemma-reg-distance}(vi)}.
Then, for each $v\in C^0(\ol{R}_g)\cap C^2(R_{g}\cup \Sigma_{g})$,
define its extension $\mcl{E}_{g}(v)$ onto $R_{(1+\kappa)g}$ by \eqref{definition-extension-operator-2018}
for $\Psi$ given by \eqref{definition-Psi-ext-2018}.
\end{definition}

\begin{proposition}[Properties of extension operator $\mcl{E}$]
\label{proposition-prop-ext-map}
For each $(u,\beta)\in \ol{\Kext}$, the extension operator $\mcl{E}_g$ given
by Definition {\rm \ref{definition-extension-map}}
maps $C^2(R_g\cup \Sigma_g)$ into $C^2(R_{(1+\kappa)g})$
with the following properties{\rm :} Fix $\alp\in(0,1)$.
Then

\smallskip
\begin{itemize}
\item[\rm (a)]
Fix $b_1,b_2$ with $-1< b_1<b_2 <1$.

\smallskip
\begin{itemize}
\item[(a-1)] There exists $C>0$ depending only on $(\iv, \gam, \beta_*, \alp)$ such that
\begin{equation*}
\|\mcl{E}_g(v)\|_{2,\alp, R_{(1+\kappa)g}\cap\{b_1<s<b_1\}} \le
C\|v\|_{2,\alp, R_{g}\cap\{b_1<s<b_1\}}.
\end{equation*}
More precisely,
\begin{equation*}
\begin{split}
&\|\mcl{E}_g(v)\|_{m,0, R_{(1+\kappa)g}\cap\{b_1<s<b_2\}}\le C \|v\|_{m,0, R_{g}\cap\{b_1<s<b_2\}} \qquad \mbox{for $m=0,1,2$},\\
&[D^2\mcl{E}_g(v)]_{\alp, R_{(1+\kappa)g}\cap\{b_1<s<b_2\}}\le  C [D^2v]_{\alp, R_{g}\cap\{b_1<s<b_2\}}.
\end{split}
\end{equation*}

\item[(a-2)] $\mcl{E}_g: C^{2,\alp}(\ol{R_g\cap\{b_1<s<b_2\}})\longrightarrow C^{2,\alp}(\ol{R_{(1+\kappa)g}\cap\{b_1<s<b_2\}})$
is linear and continuous.

\smallskip
\item[(a-3)] Suppose that $\{(u_j, \beta_j)\}\subset \ol{\Kext}$ converges
to $(u, \beta)$ in $C^{2,\til{\alp}}_{(*,\alp_1)}(\iter)\times [0,\beta_*]$
for some $\til{\alp}\in (0,1)$.
If $\{v_j\}$ satisfies
\begin{equation*}
\qquad\quad\qquad\,\,\,  v_j\in C^{2,\alp}(\ol{R_{\gshock}^{(u_j,\beta_j)}\cap\{b_1<s<b_2\}}),\,\, \|v_j\|_{2,\alp, R_{\gshock}^{(u_j,\beta_j)}\cap\{b_1<s<b_2\}}\le M
\,\,\,\,\,\tx{for all $j\in \mathbb{N}$}
\end{equation*}
for some constant $M>0$ and
converges uniformly to $v$ on compact subsets of $R_{\gshock^{(u,\beta)}}$
for some 
$v\in C^{2,\alp}(\ol{R_{\gshock^{(u,\beta)}}\cap\{b_1<s<b_2\}})$
%
then $\mcl{E}_{\gshock^{(u_j,\beta_j)}}(v_j)$ converges to
$\mcl{E}_{\gshock^{(u,\beta)}}(v)$ in $C^{2,\alp'}(\ol{R_{(1+\frac{\kappa}{2})g}\cap \{b_1<s<b_2\}})$ for all $\alp'\in(0,\alp)$,
where
$\mcl{E}_{\gshock^{(u_j,\beta_j)}}(v_j)$ is well defined on $\ol{R_{(1+\frac{\kappa}{2})\gshock^{(u,\beta)}}\cap\{b_1<s<b_2\}}$ for large $j$.
\end{itemize}

\smallskip
\item[\rm (b)] Fix $\sigma>0$ and $\eps\in(0,\frac 14]$.

\smallskip
\begin{itemize}
\item[(b-1)]
There exists $C_{\rm{par}}>0$ depending only on $(\iv, \gam, \beta_*, \alp, \sigma)$ such that
\begin{equation*}
\begin{split}
&\|\mcl{E}_g(v)\|^{(\sigma), (\rm{par})}_{2,\alp, R_{(1+\kappa)g}\cap\{-1<s<-1+\eps\}}\le
C_{\rm{par}} \|v\|^{(\sigma), (\rm{par})}_{2,\alp, R_{g}\cap\{-1<s<-1+\eps\}},\\[1mm]
&\|\mcl{E}_g(v)\|^{(\sigma), (\rm{par})}_{2,\alp, R_{(1+\kappa)g}\cap\{1-\eps<s<1\}}\le
C_{\rm{par}} \|v\|^{(\sigma), (\rm{par})}_{2,\alp, R_{g}\cap\{1-\eps<s<1\}}.
\end{split}
\end{equation*}

\item[(b-2)] The map
$$
\qquad\qquad\qquad\mcl{E}_g:\, C^{2,\alp}_{(\sigma), ({\rm par})}(R_g\cap \{-1<s<-1+\eps\})\to
C^{2,\alp}_{(\sigma), ({\rm par})}(R_{(1+\kappa)g}\cap \{-1<s<-1+\eps\})
$$
is linear and continuous. The same is true when we replace $\{-1<s<-1+\eps\}$ by $\{1-\eps<s<1\}$.

\smallskip
\item[(b-3)] If $\{(u_j, \beta_j)\}\subset \ol{\Kext}$ converges to $(u, \beta)$ in $C^{2,\til{\alp}}_{(*,\alp_1)}(\iter)\times [0,\beta_*]$
for some $\til{\alp}\in (0,1)$,
and if
\begin{align*}
&\{v_j\}\subset C^{2,\alp}_{(\sigma), (\rm par)}(\ol{R_{\gshock}^{(u_j,\beta_j)}\cap\{-1<s<-1+\eps\}}),\\
&v\in C^{2,\alp}_{(\sigma), (\rm par)}(\ol{R_{\gshock}^{(u,\beta)}\cap\{-1<s<-1+\eps\}}),
\end{align*}
and $v_j$ converges uniformly to $v$ on compact subsets of $R_{\gshock^{(u,\beta)}}$,
then $\mcl{E}_{\gshock^{(u_j,\beta_j)}}(v_j)$ converges to
$\mcl{E}_{\gshock^{(u,\beta)}}(v)$ in $C^{2,\alp'}_{(\sigma'), (\rm par)}(R_{(1+\frac{\kappa}{2})g}\cap \{-1<s<-1+\eps\})$
for all $\alp'\in(0,\alp)$
and all $\sigma'\in(0,\sigma)$. The same is true when we replace $\{-1<s<-1+\eps\}$ by $\{1-\eps<s<1\}$.
\end{itemize}

\smallskip
\item[\rm (c)] Consider the case that $s\in (-1, \frac{1}{2})$.

\smallskip
\begin{itemize}
\item[(c-1)]
There exists $C_{\rm{sub}}>0$ depending only on $(\iv, \gam, \beta_*, \alp)$ such that
\begin{equation*}
\|\mcl{E}_g(v)\|^{(-1-\alp), \{s=-1\}}_{2,\alp, R_{(1+\kappa)g}\cap\{-1<s<-\frac 12\}}\le
C_{\rm{sub}} \|v\|^{(-1-\alp), \{s=-1\}}_{2,\alp, R_{g}\cap\{-1<s<-\frac 12\}}.
\end{equation*}
Furthermore, if $(v, Dv)=(0, {\bf 0})$ on $\ol{R_g}\cap\{s=-1\}$, then
\begin{equation*}
(\mcl{E}_g(v), D\mcl{E}_g(v))=(0, {\bf 0})\qquad\tx{on $\ol{R_{(1+\kappa)g}}\cap\{s=-1\}$.}
\end{equation*}

\item[(c-2)]  $\mcl{E}_g: C^{2,\alp}_{(-1-\alp), \{s=-1\}}(R_{g}\cap\{-1<s<-\frac 12\})$\\
${}\quad\,\,\,$ $\longrightarrow C^{2,\alp}_{(-1-\alp), \{s=-1\}}(R_{(1+\kappa)g}\cap\{-1<s<-\frac 12\})$ is linear and continuous.

\smallskip
\item[(c-3)]  If $\{(u_j, \beta_j)\}\subset \ol{\Kext}$ converges to $(u, \beta)$ in $C^{2,\til{\alp}}_{(*,\alp_1)}(\iter)\times [0,\beta_*]$
for some $\til{\alp}\in (0,1)$,
and if
\begin{align*}
&\{v_j\}\subset C^{2,\alp}_{(-1-\alp), \{s=-1\}}(R_{\gshock^{(u_j,\beta_j)}}\cap\{-1<s<-\frac 12\}),\\
&v\in C^{2,\alp}_{(-1-\alp), \{s=-1\}}(R_{\gshock^{(u,\beta)}}\cap\{-1<s<-\frac 12\}),
\end{align*}
and $v_j$ converges uniformly to $v$ on compact subsets of $R_{\gshock^{(u,\beta)}}$,
then $\mcl{E}_{\gshock^{(u_j,\beta_j)}}(v_j)$ converges
to $\mcl{E}_{\gshock^{(u,\beta)}}(v)$ in $C^{2,\alp'}_{(-1-\alp'), \{s=-1\}}(R_{(1+\frac{\kappa}{2})g}\cap \{-1<s<-\frac 12 \})$
for all $\alp'\in(0,\alp)$.
\end{itemize}
\end{itemize}
\end{proposition}

\begin{proof}  We divide the proof into three steps.

\medskip
{\textbf{1.}} By Remark \ref{remark-a}, ${\rm Lip}[\gshock]$ is uniformly bounded by a constant $C>0$ depending
only on $(\iv, \gam, \beta_*)$ for all $(u,\beta)\in\ol{\Kext}$.
Then statements (a-1)--(a-2) follow from \cite[Lemma 13.9.6(i)--(ii)]{CF2}.
By Lemma \ref{lemma-7-4}(d), if $\{(u_j, \beta_j)\}\subset \ol{\Kext}$ converges
to $(u, \beta)$ in $C^{2,\til{\alp}}_{(*,\alp_1)}(\iter)\times [0,\beta_*]$ for some $\til{\alp}\in (0,1)$,
then $\gshock^{(u_j,\beta_j)}$ converges to $\gshock^{(u,\beta)}$ in $C^1([-1,1])$.
Thus, we apply \cite[Lemma 13.9.6 (iii)]{CF2} to obtain statement (a-3).
\smallskip

{\textbf{2.}} Statements (b-1)--(b-2) can be proved by following Steps 2--3
in the proof of \cite[Theorem 13.9.5]{CF2}.
Since ${\rm Lip}[\gshock]$ is uniformly bounded by a constant $C>0$ depending only
on $(\iv, \gam, \beta_*)$ for all $(u,\beta)\in\ol{\Kext}$,
the estimate constant $C_{\rm par}$ in (b-1) can be given uniformly, depending only
on $(\iv, \gam, \beta_*, \alp, \sigma)$, for all $(u,\beta)\in\ol{\Kext}$.
Moreover, statement (b-3) can be proved by following Step 4 in the proof
of \cite[Theorem 13.9.5]{CF2} and using the uniform convergence
of $\gshock^{(u_j,\beta_j)}$ to $\gshock^{(u,\beta)}$ on $[-1,1]$ when $\{(u_j, \beta_j)\}\subset \ol{\Kext}$
converges to $(u, \beta)$ in $C^{2,\til{\alp}}_{(*,\alp_1)}(\iter)\times [0,\beta_*]$ for some $\til{\alp}\in (0,1)$.
\smallskip

{\textbf{3.}} Finally, we follow the proof of \cite[Theorem 13.9.8]{CF2} to obtain statements (c-1)--(c-3).
Similarly to Steps 1--2, the uniform boundedness of ${\rm Lip}[\gshock]$ for all $(u,\beta)\in\ol{\Kext}$
implies that the estimate constant $C_{\rm sub}$ depends only on $(\iv, \gam, \beta_*, \alp)$
for all $(u,\beta)\in\ol{\Kext}$.
To prove (c-3), we use the uniform convergence of $\gshock^{(u_j,\beta_j)}$ to $\gshock^{(u,\beta)}$ on $[-1,1]$
when $\{(u_j, \beta_j)\}\subset \ol{\Kext}$ converges
to $(u, \beta)$ in $C^{2,\til{\alp}}_{(*,\alp_1)}(\iter)\times [0,\beta_*]$ for some $\til{\alp}\in (0,1)$.
\end{proof}

\begin{lemma}
\label{lemma-existence-hat-g}
Let parameters $(\alp, \eps, \delta_1, \delta_3, N_1)$ in Definition {\rm \ref{definition-10-3}}
be fixed as in Proposition {\rm \ref{lemma-7-10}}.
Then there exists a constant $\delta_3^{\rm{(imp)}}>0$ depending only
on $(\iv, \gam,  \beta_*, \delta_2)$ {\rm (}where parameter $\delta_2$ in Definition {\rm \ref{definition-10-3}}
is determined later{\rm )} such that, if
$\delta_3$ further satisfies $0<\delta_3 \le \delta_{3}^{\rm{(imp)}}$,
for each $(u,\beta)\in\ol{\mcl{K}}$, there exists a unique function $\hgshock:[-1,1]\rightarrow \R_+$ such that
\begin{equation}
\label{hgshock}
(w_{\infty}-\mcl{E}_{\gshock}(\hat w))(s,\hgshock(s))=0\qquad\mbox{for all}\,\, s\in[-1, 1].
\end{equation}
Furthermore, there exists a constant $C>0$ depending only on $(\iv, \gam,  \beta_*)$ such that $\hgshock$ satisfies
\begin{align}
\label{12-77}
&\begin{cases}
\|\hgshock-\mathfrak{g}_{\beta}\|_{2,2\alp,(-\frac{1}{2},1)}^{(2),(\rm par)}+
\|\hgshock-\mathfrak{g}_{\beta}\|_{2,2\alp,(-1,0)}^{(-1-2\alp),\{-1\}}\le C,\\[1.8mm]
\frac{\dd^k}{\dd s^k}(\hgshock-\mathfrak{g}_{\beta})(-1)=0\qquad \tx{for $k=0,1$},
\end{cases}\\[2mm]
\label{12-77-add-2018}
&\|\hgshock-\gshock\|_{1, \frac{\alp}{2}, (-1,1)}\le C\delta_3,\\[1.5mm]
\label{12-77-add2-2018}
& (\hgshock-\gshock)(\pm 1)=(\hgshock-\hgshock)'(\pm 1)=0,
\end{align}
where $\mathfrak{g}_{\beta}$ is from Lemma {\rm 5.1}.

\begin{proof} We divide the proof into three steps.

\smallskip
{\textbf{1.}}
By  Definition \ref{definition-Gset-shocks-new}(i), $w$ given by \eqref{def-w-infty} satisfies
\begin{equation}
\label{equation-gshock}
w_{\infty}-w=0\qquad\mbox{on $\Sigma_{\gshock}$}.
\end{equation}
By \eqref{definiiton-iterset-ineq8} in Definition \ref{definition-10-3}(iv), Lemma \ref{lemma-12-1-mod}(a),
and \eqref{def-w-infty},
there exists a constant $C'>0$ depending on $(\iv, \gam)$ such that
\begin{equation}
\label{estimate-w-new2015}
|D(w_{\infty}-w)|\ge C' \mu_1>0\qquad\, \mbox{on $\Sigma_{\gshock}$}.
\end{equation}
Therefore, we have
\begin{equation*}
\frac{D(w_{\infty}-w)}{|D(w_{\infty}-w)|}
=-\frac{(-\gshock',1)}{\sqrt{1+(\gshock')^2}}\qquad\,\mbox{on $\Sigma_{\gshock}$}.
\end{equation*}
Since ${\rm Lip}[\gshock]$ is uniformly bounded by a constant $C>0$ depending only on $(\iv, \gam, \beta_*)$
for all $(u,\beta)\in\ol{\mcl{K}}$, there exists a constant $\bar{\mu}>0$ depending only on $(\iv, \gam,  \beta_*)$ to satisfy
\begin{equation}
\label{t-deriv-w-Gshock}
\der_{t'}(w_{\infty}-w)=-\frac{|D(w_{\infty}-w)|}{\sqrt{1+(\gshock')^2}}\le  - \bar{\mu}\qquad\, \mbox{on}\,\, \Sigma_{\gshock}.
\end{equation}

For each $(u,\beta)\in\ol{\Kext}$, the corresponding function $\gshock=\gshock^{(u,\beta)}$
satisfies that $\gshock(-1)\ge 0$. Therefore, Definition \ref{definition-10-3}(iii) implies that
\begin{equation}
\label{gshock-range-2018}
\begin{split}
\frac{1}{N_3}(1+s)\le \gshock(s)\le \gshock(-1)+N_3(1+s)\qquad\,&\tx{for $-1\le s\le -1+\hat{\eps}_0$}
\end{split}
\end{equation}
for $\hat{\eps}_0=\frac{1}{5}$, where $N_3>1$ is the constant from Definition \ref{definition-10-3}(iii).
The lower bound of $\gshock(s)$ in \eqref{gshock-range-2018} is obtained from
Definition \ref{definition-10-3}(iii),
and $\gshock(-1)\ge 0$ which follows from \eqref{gshock-at-endpts}.

Let $\kappa\in(0, \frac 13]$ be fixed as in Definition \ref{definition-extension-map}.
In other words, let $\kappa$ be from  Lemma \ref{lemma-reg-distance}(vi).
By Definition \ref{definition-10-3}(i),
Remark \ref{remark-a}, \eqref{t-deriv-w-Gshock},
and Proposition \ref{proposition-prop-ext-map},
there exists a small constant $\sigma\in(0, \frac 14 \min\{1,\kappa\}]$ depending only
on $(\iv, \gam,  \beta_*)$ such that, for each $(u,\beta)\in\ol{\mcl{K}}$,
$\gshock$ satisfies
\[
0<\gshock(s)-\sigma<\gshock(s)+\sigma<(1+\kappa)\gshock(s)\qquad\tx{for $-1+\frac{\hat{\eps}_0}{2}\le s\le 1$},
\]
and the corresponding function $w$ given by \eqref{def-w-infty} satisfies
\begin{equation}
\label{t-deriv-w-ext}
  \begin{split}
   \,\, \der_{t'}(w_{\infty}-\mcl{E}_{\gshock}(w))(s,t')\le -\frac{\bar{\mu}}{2}
    \qquad&\tx{for $-1< s\le -1+\hat{\eps}_0$ and $1-\sigma \le\frac{ t'}{\gshock(s)} \le 1+\sigma$},\\
   \,\,  \der_{t'}(w_{\infty}-\mcl{E}_{\gshock}(w))(s,t')\le -\frac{\bar{\mu}}{2}
    \qquad&\tx{for $-1+\frac{\hat{\eps}_0}{2}\le s\le 1$ and $|t'-\gshock(s)| \le \sigma$}.
  \end{split}
\end{equation}

{\textbf{2.}} By \eqref{equation-gshock} and the linearity of the extension operator $\mcl{E}_{\gshock}$, we have
\begin{equation*}
 (w_{\infty}-\mcl{E}_{\gshock}(\hat w))(s,(1+\sigma)\gshock(s))
 =\mcl{A}_1+\mcl{A}_2,
 \end{equation*}
 where
 \begin{equation*}
 \begin{split}
&\mcl{A}_1=  (w_{\infty}-\mcl{E}_{\gshock}(w))(s,(1+\sigma)\gshock(s))-
 (w_{\infty}-\mcl{E}_{\gshock}(w))(s,\gshock(s)),\\
&\mcl{A}_2= \mcl{E}_{\gshock}(w-\hat{w})(s,(1+\sigma)\gshock(s)).
 \end{split}
\end{equation*}
By \eqref{gshock-range-2018}--\eqref{t-deriv-w-ext}, we have
\begin{equation}
\label{estimate-A1-gshock}
\mcl{A}_1\le -\frac{\bar{\mu}\sigma}{2N_3}(1-|s|)\qquad\tx{for $-1\le s \le -1+\hat{\eps}_0$}.
\end{equation}
By \eqref{12-65}, \eqref{def-w-infty}, and properties (b-1) and (c-1) of Proposition \ref{proposition-prop-ext-map},
there exists a constant $C>0$ depending only on $(\iv, \gam, \beta_*)$ such that
\begin{equation}
\label{estimate-A2-gshock}
|\mcl{A}_2|\le C\delta_3(1-|s|)\qquad\tx{for $-1\le s\le 1$},
\end{equation}
where $\delta_3>0$  is the constant in \eqref{12-65}.
From \eqref{estimate-A1-gshock}--\eqref{estimate-A2-gshock}, we obtain
\begin{equation*}
(w_{\infty}-\mcl{E}_{\gshock}(\hat w))(s,(1+\sigma)\gshock(s))
\le  (1-|s|)\Big(C\delta_3-\frac{\bar{\mu}\sigma}{2N_3}\Big)\qquad\tx{for $-1\le s \le -1+\hat{\eps}_0$}.
\end{equation*}
Therefore, a constant $\delta_3^{{\rm{(imp)}}}\in(0,\bar{\delta}_3]$  can be chosen depending only
on $(\iv, \gam,  \beta_*)$ such that, whenever $\delta_3\in(0, \delta_3^{{\rm{(imp)}}}]$,
the inequality above implies that, for any $(u,\beta)\in\ol{\mcl{K}}$,
\begin{equation}
\label{gshock-imp1}
\left(w_{\infty}-\mcl{E}_{\gshock}(\hat w)\right)(s,(1+\sigma)\gshock(s))<0
\qquad\tx{for $-1< s\le -1+\hat{\eps}_0$}.
\end{equation}
Under the same choice of $\delta_3$, we also have
\begin{equation}
\label{gshock-imp2}
\left(w_{\infty}-\mcl{E}_{\gshock}(\hat w)\right)(s,(1-\sigma)\gshock(s))>0
\qquad\tx{for $-1< s\le -1+\hat{\eps}_0$}.
\end{equation}

Adjusting the argument above, we can further reduce $\delta_3^{\rm{(imp)}}>0$ depending only on $(\iv, \gam,  \beta_*)$
so that, whenever $\delta_3\in(0,\delta_3^{\rm{(imp)}}]$,
\begin{equation}
\label{gshock-imp3}
(w_{\infty}-\mcl{E}_{\gshock}(\hat w))(s,\gshock(s)+\sigma)<0<
(w_{\infty}-\mcl{E}_{\gshock}(\hat w))(s,\gshock(s)-\sigma)\quad\tx{for $-1+\frac{\hat{\eps}_0}{2}\le s\le 1$.}
\end{equation}

{\textbf{3.}} Finally, by \eqref{12-65}, \eqref{t-deriv-w-ext}, and Proposition \ref{proposition-prop-ext-map},
we can reduce $\delta_3^{\rm{(imp)}}>0$
depending only on $(\iv, \gam,  \beta_*)$ so that, whenever $\delta_3\in(0,\delta_3^{\rm{(imp)}}]$,
$\hat{w}$ satisfies
\begin{equation}
\label{t-deriv-hatw-ext}
  \begin{split}
    \quad\der_{t'}(w_{\infty}-\mcl{E}_{\gshock}(\hat w))(s,t')\le -\frac{\bar{\mu}}{4}
    \quad\, &\tx{for $-1\le s\le -1+\hat{\eps}_0$ and $1-\sigma \le\frac{ t'}{\gshock(s)} \le 1+\sigma$},\\
    \quad \der_{t'}(w_{\infty}-\mcl{E}_{\gshock}(\hat w))(s,t')\le -\frac{\bar{\mu}}{4}
    \quad\, &\tx{for $-1+\frac{\hat{\eps}_0}{2}\le s\le 1$ and $|t'-\gshock(s)| \le \sigma$}.
  \end{split}
\end{equation}
Then \eqref{hgshock} follows from the implicit function theorem.
By \eqref{gshock-range-2018} and \eqref{gshock-imp1}--\eqref{gshock-imp3},
there exists a constant $C>0$ depending only on $(\iv, \gam, \beta_*)$ such that
\begin{equation*}
\|\hgshock-\gshock\|_{C^0([-1,1])}<C\sigma.
\end{equation*}

By Lemmas \ref{lemma-12-1-mod} and \ref{proposition-est-hpsi-new2015}, and definition \eqref{def-w-infty}, for any $\eps\in(0,1)$, we have
\begin{equation*}
\|\hat w\|_{2,2\alp, R_{\gshock}\cap\{s>-1+\eps\}}^{(2), (\rm par)}\le C_{\eps},
\end{equation*}
where constant $C_{\eps}>0$ depends only on $(\iv, \gam, \beta_*)$ and $\eps$.
Furthermore, by Lemmas \ref{lemma-12-1-mod} and \ref{lemma-est-hpsi-nrlsonic}, we obtain
\begin{equation*}
\|\hat w\|_{2,2\alp, R_{\gshock}\cap\{-1<s<0\}}^{(-1-2\alp),\{s=-1\}}\le C,\qquad
\hat w(-1, t')= D\hat{w}(-1,t')=0\,\,\,\,\,\tx{for $0<t'<\gshock(-1)$}
\end{equation*}
for a constant $C>0$ depending only on $(\iv, \gam, \beta_*)$.
Combining these two estimates of $\hat w$ with \eqref{hgshock}, \eqref{t-deriv-hatw-ext}, and Proposition \ref{proposition-prop-ext-map},
we obtain \eqref{12-77}.

Next, we use \eqref{definition-hat-vphi-ubeta}--\eqref{def-w-infty}, Lemma \ref{lemma-12-1-mod}, Definition \ref{definition-Gset-shocks-new}(ii),
Lemma \ref{lemma-7-4}(d), and estimate \eqref{12-65} given in Definition \ref{definition-10-3}(vii) to obtain
\begin{equation*}
\|\hat w-w\|_{1,\alp/2, \mcl{G}_1^{\beta}(\Om)}
=\|(\hat u-u)\circ \mathfrak{F}_{(u,\beta)}^{-1} \circ (\mcl{G}_1^{\beta})^{-1}\|_{1,\alp/{2}, \mcl{G}_1^{\beta}(\Om)}
\le C\delta_3
\end{equation*}
for a constant $C>0$ depending only on $(\iv, \gam, \beta_*)$.
Using this estimate and \eqref{t-deriv-w-ext}, we obtain \eqref{12-77-add-2018}.
Finally, \eqref{12-77-add2-2018} follows directly from \eqref{12-77} and
the fact that $\frac{\dd^k}{\dd x^k}(\gshock-\mathfrak{g}_{\beta})(\pm 1)=0$ for $k=0,1$.
\end{proof}
\end{lemma}

Let parameters $(\alp, \eps, \delta_1, \delta_3, N_1)$ in Definition \ref{definition-10-3}
be chosen as in Lemma \ref{lemma-existence-hat-g}.
For each $(u,\beta)\in\ol{\mcl{K}}$, let $\hgshock:[-1,1]\rightarrow \R_+$ be
given by \eqref{hgshock}.
From \eqref{12-77-add-2018}--\eqref{12-77-add2-2018}, further reducing $\delta_3$,
we obtain that $\hgshock$ satisfies
estimate \eqref{gshock-positivity-uterset-def2018} in  Definition \ref{definition-10-3}(iii) with $N_3$ replaced by $2N_3$.
We define a function $\til u:\ol{\iter}\rightarrow \R$ by
\begin{equation}
\label{12-78}
\til u=\mcl{E}_{\gshock}(\hat{w})\circ (G_{2,\hgshock})^{-1}
\end{equation}
for $G_{2,\hgshock}$ defined by \eqref{7-b9}.
By Corollary \ref{corollary-hatu-2018}, Proposition  \ref{proposition-prop-ext-map}, and Lemma \ref{lemma-existence-hat-g},
there exists a constant $C>0$ depending only on $(\iv, \gam,  \beta_*)$
such that $\til u$ satisfies
\begin{equation}
\label{estimate-tilu-new2015}
\|\til u\|_{2,2\alp,\iter}^{(*,1)}\le C.
\end{equation}

Now we define the iteration map $\mcl{I}: \ol{\mcl{K}} \rightarrow C^{2,\alp}_{(*,\alp_1)}(\iter)$.

\begin{definition}
\label{definition-iterationmap}
Let parameters $(\alp, \eps, \delta_1, \delta_3, N_1)$ in Definition {\rm \ref{definition-10-3}}
be fixed as in Proposition {\rm \ref{lemma-7-10}}.
Then we adjust $\delta_3\in(0, \delta_3^{\rm imp}]$ for $\delta_3^{\rm imp}$ from Lemma {\rm \ref{lemma-existence-hat-g}}
so that Lemma {\rm \ref{lemma-existence-hat-g}} holds for all $(u,\beta)\in \ol{\mcl{K}}$.
For each $(u,\beta)\in \ol{\mcl{K}}$, let $\til u$ be given by \eqref{12-78}.
Then define an iteration map $\Itr: \ol{\mcl{K}} \rightarrow  C^{2,\alp}_{(*,\alp_1)}(\iter)$ by
\begin{equation*}
\Itr(u,\beta)=\til u.
\end{equation*}
\end{definition}

\begin{lemma}
\label{lemma-12-6}
The iteration map $\Itr$ defined in Definition {\rm \ref{definition-iterationmap}}
satisfies the following properties{\rm :}

\smallskip
\begin{itemize}
\item[(a)] For any $\beta\in[0,\beta_*]$, define
\begin{equation*}
\ol{\mcl{K}}(\beta):=\big\{u\in C^{2,\alp}_{(*,\alp_1)}(\iter)\,:\,(u,\beta)\in \ol{\mcl{K}}\big\}.
\end{equation*}
For each $(u,\beta)\in \ol{\mcl{K}}$, define
\begin{equation*}
\Itrf^{(\beta)}(u)=\hat u,
\end{equation*}
where $\hat u$ is given by \eqref{definition-sol-ubvp}.
Then $u\in \ol{\mcl{K}}(\beta)$ satisfies $\Itr(u, \beta)=u$
if and only if $\Itrf^{(\beta)}(u)=u$.

\smallskip
\item[(b)] For $\til{\alp}=\frac{\alp}{2}$, there exists a constant $C>0$ depending only on
$(\iv, \gam,  \beta_*)$ such that, for each $(u,\beta)\in\ol{\mcl{K}}$,
    \begin{equation*}
    \|\Itr(u,\beta)\|_{2,\alp+\til{\alp},\iter}^{(*,1)}\le C.
    \end{equation*}
\end{itemize}

\begin{proof}
For a fixed $\beta\in [0, \beta_*]$, suppose that $\Itr(u,\beta)=u$ for some $u\in \ol{\mcl{K}}(\beta)$;
that is, $\til u=u$ for $\til u$ given by \eqref{12-78}.
Then, by Definition \ref{definition-Gset-shocks-new} and \eqref{12-78}, we see that, for all $s\in [-1,1]$,
\begin{equation*}
w_{\infty}(s,\gshock(s))=u(s,1)=\mcl{E}_{\gshock}(\hat w)(s, \hgshock(s))=w_{\infty}(s,\hgshock(s)).
\end{equation*}
This, combined with  Lemma \ref{lemma-7-2} and \eqref{def-w-infty},
implies that $\gshock=\hgshock$ on $[-1,1]$.
Then it follows from \eqref{12-78} that $\til u=\mcl{E}_{\gshock}(\hat{w})\circ (G_{2,\gshock})^{-1}=\hat u$,
which implies that
$u=\hat u=\Itr_1^{(\beta)}(u)$ in $\iter$.

Next, suppose that $\Itr_1^{(\beta)}(u)=u$ for some $u\in \ol{\mcl{K}}(\beta)$.
Then $\gshock=\hgshock$ on $[-1,1]$.
This, combined with \eqref{12-78}, implies that
$\til u=\Itr(u,\beta)=\mcl{E}_{\gshock}(\hat{w})\circ (G_{2,\gshock})^{-1}=\hat u$.
Therefore, we obtain that $\til  u=u$ in $\iter$.

Finally, statement (b) directly follows from \eqref{estimate-tilu-new2015}.
\end{proof}
\end{lemma}

\section{Fixed Points of $\Itr(\cdot, \beta)$ and Admissible Solutions}
For the iteration map $\Itr$ defined in Definition \ref{definition-iterationmap},
we show that, if $u\in \ol{\mcl{K}}(\beta)$ is a fixed point of $\Itr(\cdot,\beta)$
for some $\beta\in(0, \beta_*]$,
then $\vphi$ defined by \eqref{def-vphi} in Definition \ref{definition-Gset-shocks-new}
is an admissible solution corresponding to $(\iv, \beta)\in\mathfrak{R}_{\rm weak}$
in the sense of Definition \ref{def-regular-sol}.

\begin{proposition}
\label{lemma-8-2}
Let parameters $(\alp, \eps, \delta_1, \delta_3, N_1)$ in Definition {\rm \ref{definition-10-3}}
be fixed as in Definition {\rm \ref{definition-iterationmap}}.
Then parameters $(\eps, \delta_1)$ can be further reduced depending only
on $(\iv, \gam, \beta_*)$ so that, for each $\beta\in(0, \beta_*]$,
$u\in \mcl{K}(\beta)$ is a fixed point of
$\Itr(\cdot, \beta)\,:\,\ol{\mcl{K}}(\beta)\rightarrow C^{2,\alp}_{(*,\alp_1)}(\iter)$
if and only if
$\vphi$, defined by \eqref{def-vphi} in Definition  {\rm \ref{definition-Gset-shocks-new}},
yields an admissible solution corresponding to $(\iv, \beta)\in\mathfrak{R}_{\rm weak}$
in the sense of Definition {\rm \ref{def-regular-sol}} by extending $\vphi$
into $\Lambda_{\beta}$ via \eqref{1-24}
 if $\beta<\betasonic$,
 and via \eqref{1-24ab} if $\beta\ge \betasonic$.

\begin{proof}  By Corollary \ref{corollary-admisol-iterset}, it suffices to prove that,
if $u\in \mcl{K}(\beta)$ is a fixed point of $\Itr(\cdot, \beta):\ol{\mcl{K}}(\beta)\rightarrow C^{2,\alp}_{(*,\alp_1)}(\iter)$,
then $\vphi$, defined by \eqref{def-vphi} in Definition  {\rm \ref{definition-Gset-shocks-new}},
yields an admissible solution corresponding to $(\iv, \beta)\in\mathfrak{R}_{\rm weak}$ in the sense of
Definition {\rm \ref{def-regular-sol}}.
We divide the proof into six steps.

\smallskip
{\textbf{1.}} For $(u,\beta)\in\ol{\mcl{K}}$,
let $(\Om, \shock, \vphi)=(\Om(u,\beta), \shock(u,\beta), \vphi^{(u,\beta)})$ be defined
by Definition \ref{definition-Gset-shocks-new}, and denote $\phi:=\vphi-\rightvphi$.
Let $\hat{\phi}\in C^2(\Om)\cap C^1(\ol{\Om})$ be the unique solution of the boundary value
problem \eqref{12-50} determined by $(u, \beta)$.

Suppose that
\begin{equation*}
\Itr(u,\beta)=u\qquad\tx{for some $u\in \ol{\mcl{K}}(\beta)$.}
\end{equation*}
By Lemma \ref{lemma-12-6}, we have
\begin{equation}
\label{fp-equality}
\hat{\phi}=\phi\qquad\tx{in $\Om$}.
\end{equation}
Let $\vphi$ be extended onto $\Lambda_{\beta}$ by \eqref{1-24} for $\beta<\betasonic$,
and by \eqref{1-24ab}
for $\beta\ge \betasonic$.
Moreover, let $\leftsonic$, $\rightsonic$, $\leftvec$, $\rightvec$, $\ivphi$, $\leftvphi$,
and $\rightvphi$ be defined by Definition \ref{definition-domains-np}.

\smallskip
{\textbf{2.}} {\emph{Verification of properties {\rm (i-2)}--{\rm (i-4)} and {\rm (ii-1)}--{\rm (ii-3)}
of Definition {\rm \ref{def-regular-sol}}}}.

Properties (i-2)--(i-3) follows from Remark \ref{remark-a}(i).
By using Lemma \ref{lemma-7-4}(b), it can directly be checked that property (i-4) holds.

\smallskip
By Definition \ref{definition-10-3}(i) (or Corollary \ref{corollary-hatu-2018})
and the extension of $\vphi$ onto $\Lambda_{\beta}$ described in Step 1,
$\vphi$ satisfies properties (ii-1) and (ii-3).

We define
\begin{equation*}
A_{ij}({\bmxi}):=A_{ij}(D\hat{\phi},\bmxi), \qquad i,j=1,2,
\end{equation*}
for $A_{ij}(D\hat{\phi},\bmxi)$ given by \eqref{12-57}.
By Definition \ref{definition-10-3}(i) (or Corollary \ref{corollary-hatu-2018}),
coefficients $A_{ij}({\bmxi}), i,j=1,2$, of equation  $\mcl{N}_{(u,\beta)}(\hat{\phi})=0$ in \eqref{12-50}
are in $C^{1,\alp}(\ol{\Om}\setminus (\ol{\leftsonic}\cup\ol{\rightsonic}))$.
Furthermore, Lemma \ref{lemma-6-1}(a) implies that $\mcl{N}_{(u,\beta)}(\hat{\phi})=0$ is strictly elliptic in $\Om$.
Then the standard interior Schauder estimates for linear elliptic equations
imply that $\vphi\in C^{3,\alp}(\Om)$.
This, combined with Definition \ref{definition-10-3}(i) (or Corollary \ref{corollary-hatu-2018}),
implies that $\vphi$ satisfies property (ii-2).

\smallskip
{\textbf{3.}} {\emph{Verification of property {\rm (iv)} of Definition {\rm \ref{def-regular-sol}}}}.

For
$A_{ij}(\bmxi)$ defined in Step 2, we define a linear operator $\mcl{L}_{(u,\beta)}$ by
\begin{equation*}
\mcl{L}_{(u,\beta)} (v):=\sum_{i,j=1}^2A_{ij}\der_{\xi_i\xi_j}v.
\end{equation*}
Since $\ivphi-\rightvphi$ is a linear function of $\xxi$, and $\vphi-\ivphi=\hat{\phi}-(\ivphi-\rightvphi)$, we have
\begin{equation}
\label{L-equation-2018}
\mcl{L}_{(u, \beta)}(\vphi-\ivphi)=\mcl{L}_{(u, \beta)}(\hat{\phi})=0\qquad\tx{in}\,\, \Om.
\end{equation}
By Lemma \ref{lemma-6-1}(a), the equation stated above is strictly elliptic in $\Om$ so that
the maximum principle applies to $\vphi-\ivphi$ in $\Om$.
From \eqref{equation-gshock} and \eqref{fp-equality}, we obtain that $\vphi-\ivphi=0$ on $\shock$.
By Definition \ref{definition-smooth-connection}(ii),
it follows directly from the boundary condition $\hat{\phi}=\max\{\leftvphi,\rightvphi\}-\rightvphi$
on $\leftsonic\cup\rightsonic$ given in \eqref{12-50} that $\vphi-\ivphi=\leftvphi-\ivphi\le 0$ on $\leftsonic$,
and $\vphi-\ivphi=\rightvphi-\ivphi\le 0$ on $\rightsonic$.
Furthermore, the boundary condition for $\hat{\phi}_{\xi_2}=0$ on $\Wedge$ given in \eqref{12-50} implies that
\begin{equation}
\label{fp-W-bc-wedge-2018}
\der_{\etan}(\ivphi-\vphi)=-\iv<0\qquad\tx{on $\Wedge$.}
\end{equation}
Therefore, by the maximum principle and Hopf's lemma, we obtain
\begin{equation}
\label{fp-vphi-upper-bd2018}
\vphi\le \ivphi\quad\tx{in $\Om$}.
\end{equation}

When $\beta<\frac{2\delta_1}{N_1^2}$, we have shown in Step 2 in the proof of Lemma \ref{lemma-positivity-hatphi-2018} that
\begin{equation}
\label{left-inequality-new}
\max\{\leftvphi, \rightvphi\}\le \vphi\qquad \mbox{in} \,\, \Om.
\end{equation}

When $\beta\ge \frac{2\delta_1}{N_1^2}$, \eqref{definiiton-iterset-ineq1} in Definition \ref{definition-10-3}(iv)
implies that $\max\{\leftvphi, \rightvphi\}\le \vphi$ holds
in $\Om\setminus (\oD_{\eps/{10}}\cup \nD_{\eps/{10}})$.
Note that parameter $\eps$ in Definition \ref{definition-10-3} has been chosen
so that $\eps<\frac{\leftch}{\bar{k}}$ for $\frac{\leftch}{\bar{k}}$ from \eqref{7-b7}
in Definition \ref{definition-smooth-connection}.
Therefore, $\vphi_{\beta}^*=\max\{\leftvphi, \rightvphi\}$
in $\Om\cap (\oD_{\eps}\cup \nD_{\eps})$ for $\vphi_{\beta}^*$ given by \eqref{12-32}.
Then we obtain from  \eqref{strict-positivity} in Lemma \ref{proposition-est-hpsi-new2015}
that $\max\{\leftvphi, \rightvphi\}\le \vphi$ holds in $\Om\cap (\oD_{{\eps}/{10}}\cup \nD_{{\eps}/{10}})$.

Therefore, we conclude that inequality \eqref{left-inequality-new} holds for any $\beta\in (0,\beta_*]$.
Combining this inequality with \eqref{fp-vphi-upper-bd2018}, we conclude that $\vphi$ satisfies
property (iv) of Definition \ref{def-regular-sol}.

\medskip
{\textbf{4.}}
{\emph{Verification of  property {\rm (v)} of Definition {\rm \ref{def-regular-sol}}}}.
In order to show that $\vphi$ satisfies property (v) of Definition \ref{def-regular-sol},
it suffices to verify the following claim:

\smallskip
{\emph{Claim. There exist small constants $\eps_{\rm fp}>0$ and $\delta_{\rm fp}>0$ depending only
on $(\gam, \iv, \beta_*)$ so that, if parameters $(\eps, \delta_1)$ in Definition {\rm \ref{definition-10-3}}
satisfy $\eps\in(0,\eps_{\rm fp}]$ and $\delta_1\in(0, \delta_{\rm fp}]$, then $\vphi$ satisfies
\begin{equation}
\label{monotonicity-of-fp}
\der_{\leftvec}(\ivphi-\vphi)\le 0,\quad \der_{\rightvec}(\ivphi-\vphi)\le 0\qquad\tx{in $\Om$.}
\end{equation}
}}

\smallskip
Similarly to the previous step, we consider two cases: $\beta\in[\frac{\delta_1}{N_1^2}, \beta_*]$
and $\beta\in(0, \frac{\delta_1}{N_1^2})$, separately.

\smallskip
{\textbf{4-1}.}
Suppose that $\beta \in[\frac{\delta_1}{N_1^2}, \beta_*]$. Define
\begin{equation*}
W:=\ivphi-\vphi\qquad\tx{in $\Om$}.
\end{equation*}

Let $(X,Y)$ be the rectangular coordinates such that $(\leftvec, \leftvec^{\perp})=({\bf e}_X, {\bf e}_Y)$.
By \eqref{L-equation-2018}, $W$ satisfies that $\mcl{L}_{(u,\beta)}(W)=0$ in $\Om$.
Since the $(X,Y)$--coordinates are obtained from rotating the $(\xi_1,\xi_2)$--plane by $\beta$ counter-clockwise,
equation $\mcl{L}_{(u,\beta)}(W)=0$ can be rewritten in the $(X,Y)$--coordinates as follows:
\begin{equation}
\label{fp-equation-W}
\hat{A}_{11}W_{XX}+2\hat{A}_{12}W_{XY}+\hat{A}_{22}W_{YY}=0\qquad\tx{in $\Om$},
\end{equation}
with $\hat{A}_{ij}\in C^{\alp}(\ol{\Om})\cap C^{1,\alp}(\ol{\Om}\setminus (\ol{\leftsonic}\cup\ol{\rightsonic})), i,j=1,2$.

Define
\begin{equation*}
w:=W_X=\der_{\leftvec}(\ivphi-\vphi).
\end{equation*}

By \eqref{definiiton-iterset-ineq2} in Definition \ref{definition-10-3}(iv), $w$ satisfies
\begin{equation}
\label{inequality-w-left}
w<0 \qquad \tx{in} \,\,\ol{\Om}\setminus \oD_{{\eps}/{10}}.
\end{equation}
Next, we prove that $w\le 0$ in $\Om\cap \oD_{{\eps}/{10}}$.

Differentiating \eqref{fp-equation-W} with respect to $X$, we have
\begin{equation*}
\hat{A}_{11}w_{XX}+2\hat{A}_{12}w_{XY}+\hat{A}_{22}w_{YY}
+\der_X\hat{A}_{11}w_X+2\der_X\hat{A}_{12}w_Y+\der_X\hat{A}_{22}W_{YY}=0\qquad\tx{in $\Om$}.
\end{equation*}
Using the strict ellipticity of operator $\mcl{L}_{(u,\beta)}$ following from Lemma \ref{lemma-6-1}(a),
we obtain that $\hat{A}_{22}>0$ in $\Om$ such that  $W_{YY}$ can be expressed
as
$$
W_{YY}=-\frac{\hat{A}_{11}w_{X}+2\hat{A}_{12}w_{Y}}{\hat{A}_{22}} \qquad \mbox{in $\Om$}.
$$
Substituting this expression into the equation immediately above,
we obtain a strictly elliptic equation for $w$ in the following form:
\begin{equation}
\label{equation-w}
\hat{A}_{11}w_{XX}+2\hat{A}_{12}w_{XY}+\hat{A}_{22}w_{YY}+\hat{A}_1w_{X}+\hat{A}_2w_Y=0\qquad\,\,\tx{in $\Om$}.
\end{equation}
Since $\hat{A}_{ij}\in C^{\alp}(\ol{\Om})\cap C^{1,\alp}(\ol{\Om}\setminus (\ol{\leftsonic}\cup\ol{\rightsonic})), i,j=1,2$,
we see that $\hat{A}_i\in C^{\alp}(\ol{\Om}\setminus (\ol{\leftsonic}\cup\ol{\rightsonic})), i=1,2$.

By a direct computation, applying Lemma \ref{lemma-est-hpsi-nrlsonic} and the definitions of $(\leftvec, \ivphi, \leftvphi)$ given
in Definition \ref{definition-domains-np},
we have
\begin{equation}
\label{fp-w-bc-lsonic2018}
w=\der_{\leftvec}(\ivphi-\leftvphi)=0\qquad\,\,\tx{on $\leftsonic$}.
\end{equation}
On $\Wedge$, $w$ satisfies the homogeneous oblique boundary condition:
\begin{equation}
\label{bc-w-wedge}
{\bf b}_{\rm w}\cdot \nabla w=0\qquad\,\, \tx{with ${\bf b}_{\rm w}\cdot {\bf n}_{\rm w}>0$ on $\Wedge$}
\end{equation}
for the inward unit normal vector ${\bf n}_{\rm w}$ to $\Wedge$.
This can be verified as follows:
Differentiating the boundary condition \eqref{fp-W-bc-wedge-2018} along $\Wedge\subset \{\xi_2=0\}$,
we find that $W_{\xi_1\xi_2=0}$ on $\Wedge$.
Equation \eqref{L-equation-2018}, combined with $\der_{\xin\etan}W=0$ on $\Wedge$,
leads to
$$
A_{11}W_{\xin\xin}+A_{22}W_{\etan\etan}=0 \qquad\mbox{on $\Wedge$}.
$$
Note that $A_{11}>0$ and $A_{22}>0$ hold on $\Wedge$ by Lemma \ref{lemma-6-1}(a).
Then a direct computation by using the definition of $\leftvec$ shows that
$$
\frac{A_{11}}{\cos\beta}w_{\xin}+\frac{A_{22}}{\sin\beta}w_{\etan}=0 \qquad\mbox{with $\frac{A_{22}}{\sin\beta}>0\,\,$ on $\Wedge$}.
$$
This implies the strict obliqueness of the boundary condition for $w$ on $\Wedge$.

In order to obtain a boundary condition for $w$ on $\shock$, we apply \cite[Lemma 13.4.5]{CF2}.
For this purpose, we need to check that all the conditions to apply \cite[Lemma 13.4.5]{CF2} are satisfied.

Let $\oM$ and $\leftc$ be given by \eqref{1-25}, and let $\leftshock$ and $\Oo$ be given by Definition \ref{definition-domains-np}.
Then $\leftc-{\rm dist}(\leftshock, \Oo)>0$ if and only if $\oM<1$.
By Lemma \ref{remark-1-1}, $\oM<1$ for $\beta=0$.
Then \eqref{oM-monotonicity} given in the proof of Lemma \ref{lemma:interval-existence}
implies that $\oM<1$ for $\beta\in(0,\beta_*]$.
Therefore, there exists a constant $\mu_0>0$ depending only on $(\iv, \gam)$ such that
\begin{equation}
\label{definition-mu0-2018}
\leftc-{\rm dist}(\leftshock, \Oo)\ge \mu_0\qquad\,\,\tx{for all $\beta\in(0, \beta_*]$}.
\end{equation}

By Lemma \ref{lemma-ext-M-property}(h) and \eqref{fp-equality}, $\vphi$ satisfies the
Rankine-Hugoniot condition \eqref{3-a3} on $\shock$.

Let $\bm\nu$ be the unit normal vector to $\shock$ towards the interior of $\Om$,
and let $\bm\tau$ be obtained from rotating $\bm\nu$ by $\frac{\pi}{2}$
counter-clockwise ($\bm\tau$ is a unit tangent vector to $\shock$).
By Definition \ref{definition-10-3}(i) (or by Corollary \ref{corollary-hatu-2018}), we have
\begin{equation}
\label{est-bc-shock-case1}
\|\vphi-\leftvphi\|_{C^1(\ol{\Om\cap \oD_{\eps}})}
+ \|\bm\tau-{\bf e}_X\|_{C^0(\ol{\shock\cap \oD_{\eps}})}+
\|\bm\nu-(-{\bf e}_Y)\|_{C^0(\ol{\shock\cap \oD_{\eps}})}\le C\eps^{\alp}
\end{equation}
for a constant $C>0$ depending only on $(\iv, \gam,  \beta_*)$,
where point $\lefttop$  is defined in Definition \ref{definition-domains-np}.
Note that point $\lefttop$ lies on $\leftsonic$.
At $\lefttop$, ${\bm\tau}={\leftvec}={\bf e}_X$ and ${\bm\nu}=-{\bf e}_Y$.

By the definition of $A_{ij}$ given in \eqref{12-57},
Corollary \ref{corollary-hatu-2018}, and \eqref{fp-equality}, we have
\begin{equation}
\label{Aij-expression-P1}
A_{ij}=A_{ij}^{\rm{potn}}(D(\leftvphi-\rightvphi), (\leftvphi-\rightvphi)(\lefttop), \lefttop) \qquad\tx{at $\lefttop$}
\end{equation}
for $A_{ij}^{\rm{potn}}, i,j=1,2$, defined by \eqref{def-A-potn-new}.
By \eqref{def-uniform-ptnl-new}, this yields
\begin{equation*}
A_{11}=\leftc^2-(\der_{\xi_1}\leftvphi)^2,\,\,\, A_{12}=A_{21}=-\der_{\xi_1}\leftvphi\der_{\xi_2}\leftvphi=0,
\,\,\, A_{22}=\leftc^2-(\der_{\xi_2}\leftvphi)^2\qquad\tx{at $\lefttop$}.
\end{equation*}
Then we have
\begin{equation}
\label{A11-at-P1}
\sum_{i,j=1}^2A_{ij}\nu_i\nu_j=\leftc^2-(\der_{\bm\nu}\leftvphi)^2=\leftc^2-({\rm dist}(\leftshock, \Oo))^2=\leftc^2(1-\oM^2)
>\lambda_0\qquad\tx{at $\lefttop$}
\end{equation}
for some constant $\lambda_0>0$.
By \eqref{definition-mu0-2018}, constant $\lambda_0>0$ in \eqref{A11-at-P1} can be fixed, depending only on $(\iv, \gam)$.
By \eqref{est-bc-shock-case1} and \eqref{A11-at-P1},
there exists a small constant $\eps^{(1)}_{\rm fp}>0$ depending only on $(\iv, \gam, \beta_*)$ such that
\begin{equation}
\label{A11-at-nrP1}
\sum_{i,j=1}^2A_{ij}\nu_i\nu_j\ge \frac{\lambda_0}{2} \qquad\tx{in $\shock\cap \oD_{\eps^{(1)}_{\rm fp}}$}.
\end{equation}
By Lemma \ref{lemma-6-1}(a), there exists a constant $\lambda_1>0$ depending only on $(\gam, \iv, \beta_*)$ such that
\begin{equation}
\label{A11-at-awayP1}
\sum_{i,j=1}^2A_{ij}\nu_i\nu_j\ge \lambda_1 \qquad
\tx{in $(\shock\cap \oD_{{\leftch}/{10}} )\setminus \oD_{{\eps^{(1)}_{\rm fp}}/{2}}$}
\end{equation}
for $\leftch$ defined in Definition \ref{definition-Qbeta}.

Since $\vphi$ satisfies the Rankine-Hugoniot condition \eqref{3-a3} on $\shock$,
it follows from \eqref{est-bc-shock-case1} and \eqref{A11-at-nrP1}--\eqref{A11-at-awayP1}
that $\vphi$ satisfies all the conditions required to apply \cite[Lemma 13.4.5]{CF2}.
Then, by  \cite[Lemma 13.4.5]{CF2}, we obtain a boundary condition for $w$ in the form:
\begin{equation}
\label{fp-w-bc-shock-away-nsonic}
{\bf b}_{\rm sh}\cdot \nabla w=0\qquad\tx{on $\shock\cap \oD_{\eps_{\rm fp}^{(2)}}$}
\end{equation}
for some small constant $\eps_{\rm fp}^{(2)}>0$ depending on $(\gam, \iv, \beta_*)$,
where ${\bf b}_{\rm sh}$ satisfies
\begin{equation*}
{\bf b}_{\rm sh}\cdot {\bm\nu}>0\qquad\tx{on $\shock\cap \oD_{\eps_{\rm fp}^{(2)}}$}.
\end{equation*}

In conclusion, $w$ satisfies the strictly elliptic equation \eqref{equation-w} in $\Om\cap \oD_{\eps}$
for $\eps>0$ to be specified later, the boundary condition $w=0$ on $\leftsonic$,
and the oblique boundary conditions \eqref{bc-w-wedge}
on $\Wedge$ and \eqref{fp-w-bc-shock-away-nsonic} on $\shock\cap \oD_{\eps_{\rm fp}^{(2)}}$.
Therefore, if parameter $\eps>0$ in Definition \ref{definition-10-3} satisfies
\begin{equation}
\label{fp-condition-eps}
0<\eps\le \eps_{\rm fp}^{(2)},
\end{equation}
then it follows from the maximum principle, Hopf's lemma, and \eqref{inequality-w-left} that
\begin{equation*}
w\le 0\qquad\tx{in $\Om\cap \oD_{\eps_{\rm fp}^{(2)}}$}.
\end{equation*}
Finally, we combine this result with \eqref{inequality-w-left} to conclude that
\begin{equation*}
\der_{\leftvec}(\ivphi-\vphi)\le 0\qquad\tx{in $\Om$}\,\, \tx{for $\beta\in[\frac{\delta_1}{N_1^2}, \beta_*]$},
\end{equation*}
provided that $\eps$ satisfies condition \eqref{fp-condition-eps}.

\smallskip
{\textbf{4-2}.} Suppose that  $\beta\in(0, \frac{\delta_1}{N_1^2})$.
Note that $w$ satisfies \eqref{equation-w}--\eqref{bc-w-wedge}.
By the definitions of $(\leftvec, \ivphi, \rightvphi)$ given in Definition \ref{definition-domains-np}
and Corollary \ref{corollary-hatu-2018}, $w$ satisfies
\begin{equation*}
w=\der_{\leftvec}(\ivphi-\rightvphi)=-\iv\sin \beta<0\qquad\tx{on $\rightsonic$.}
\end{equation*}

By \eqref{12-25} and \eqref{def-uniform-ptnl-new},
$\leftvphi-\rightvphi=\iv(\xi_1 \tan\beta-\xi_2^{(\beta)}+\neta)$.
Note that $\neta=\xi_2^{(\beta)}|_{\beta=0}$.
Then, by \eqref{2-4-a6} and the continuous differentiability of $\iM$ with respect
to $\beta\in[0,\betadet]$, there exists a constant $C>0$ depending only on $(\gam, \iv)$ such that
\begin{equation}
\label{cont-dep-bgd-potential-beta}
\|\leftvphi-\rightvphi\|_{C^{1,\alp}(\ol{\Om})} \le C\beta \qquad\mbox{for all $\beta\in[0,\betadet]$}.
\end{equation}
By Definition \ref{definition-10-3}(i) and \eqref{cont-dep-bgd-potential-beta},
we see that, for any  $\beta\in(0, \frac{\delta_1}{N_1^2})$,
\begin{equation}
\label{est-bc-shock-case2}
\|\vphi-\leftvphi\|_{C^{1,\alp}(\ol{\Om})}\\
\le \|\vphi-\rightvphi\|_{C^{1,\alp}(\ol{\Om})}+\|\leftvphi-\rightvphi\|_{C^{1,\alp}(\ol{\Om})} \le C\delta_1
\end{equation}
for some constant $C>0$ depending only on $(\gam, \iv, \beta_*)$ so that
\begin{equation*}
[A_{ij}]_{\alp, \Om}+[{\bm\nu}]_{\alp, \shock}+[{\bm\tau}]_{\alp, \shock}\le C\delta_1
\end{equation*}
for $C>0$ depending only on $(\gam, \iv, \beta_*)$.
By \eqref{A11-at-P1} and the estimate immediately above,
there exists a small constant $\delta_{\rm fp}>0$ depending only on $(\iv, \gam, \beta_*)$ so that, if
\begin{equation}
\label{fp-condition-delta1}
\delta_1\in(0,\delta_{\rm fp}],
\end{equation}
then
\begin{equation*}
\sum_{i,j=1}^2A_{ij}\nu_i\nu_j\ge \frac{\lambda_0}{2} \qquad\tx{on $\shock$}
\end{equation*}
for $\lambda_0>0$ from \eqref{A11-at-P1}.
Then \cite[Lemma 13.4.5]{CF2} implies that $w$ satisfies a boundary condition in the form:
\begin{equation}
\label{bc-w-shock-small-beta-new2019}
{\bf b}_{\rm sh}\cdot \nabla w=0\qquad\tx{on $\shock$},
\end{equation}
with ${\bf b}_{\rm sh}$ satisfying
${\bf b}_{\rm sh}\cdot {\bm\nu}>0$ on $\shock$.

Since $w$ satisfies the strictly elliptic equation \eqref{equation-w} in $\Om$,
$w\le 0$ on $\leftsonic\cup\rightsonic$, and the strictly oblique boundary conditions \eqref{bc-w-wedge} on $\Wedge$
and \eqref{bc-w-shock-small-beta-new2019} on $\shock$, it follows from the maximum principle and Hopf's lemma that
\begin{equation*}
w\le 0\qquad\tx{in} \,\,\Om,
\end{equation*}
provided that parameter $\delta_1>0$ in Definition \ref{definition-10-3} satisfies \eqref{fp-condition-delta1}.

\smallskip
{\textbf{4-3.}} By repeating the argument in Steps 4-1 and 4-2 with $w=\der_{\leftvec}(\ivphi-\vphi)$
replaced by $w=\der_{\rightvec}(\ivphi-\vphi)$, we can also show that
\begin{equation*}
\der_{\rightvec}(\ivphi-\vphi)\le 0\qquad\tx{in $\Om$},
\end{equation*}
provided that constants $(\eps_{\rm fp}^{(2)}, \delta_{\rm fp})$ from \eqref{fp-condition-eps}
and \eqref{fp-condition-delta1} are adjusted, depending only on $(\iv, \gam, \beta_*)$.

For the rest of the proof, parameters $(\eps, \delta_1)$ in Definition \ref{definition-10-3} satisfy
\begin{equation*}
0<\delta_1< \delta_{\rm fp},\qquad 0<\eps<\min\{\eps_{\rm fp}^{(1)}, \eps_{\rm fp}^{(2)}\}.
\end{equation*}

\smallskip
{\textbf{5.}}  {\emph{Verification of property {\rm (ii-4)} of Definition {\rm \ref{def-regular-sol}}}}.
Since Eq. \eqref{2-1} is equivalent to \eqref{2-3}, it suffices to check that
equation $\mcl{N}_{(u,\beta)}(\phi)=0$ coincides with Eq. \eqref{2-3}.

\smallskip
{\textbf{5-1.}} {\emph{Equation $\mcl{N}_{(u,\beta)}(\phi)=0$ away from $\leftsonic\cup\rightsonic$}}.
In order to show that $\vphi$ satisfies property (ii-4), it suffices to show that
equation $\mcl{N}_{(u,\beta)}(\phi)=0$ from \eqref{12-50} coincides with Eq. \eqref{2-3} in $\Om$.
By Lemma \ref{lemma-6-1}(i), equation $\mcl{N}_{(u,\beta)}(\phi)=0$ coincides with Eq. \eqref{2-3}
in $\Om\setminus (\oD_{\eps/10}\cup \nD_{\eps/10})$ for parameter $\eps>0$ in Definition \ref{definition-10-3}
fixed as in Definition \ref{definition-iterationmap}.

\smallskip
{\textbf{5-2.}} {\emph{Equation $\mcl{N}_{(u,\beta)}(\phi)=0$ near $\rightsonic$}}.
In $\nOm_{\eps}:=\Om\cap \nD_{\eps}$,
let the $(x,y)$--coordinates be defined by \eqref{coord-n}.
Define $\psi:=\vphi-\rightvphi=\hat{\vphi}-\rightvphi$ in $\nOm_{\eps}$.
By Lemma \ref{lemma-8-3}(g), if it can be shown that
\begin{equation}
\label{fp-eqn-coinc-ne-rsonic}
\big|\psi_x(x,y)\big|<\frac{2-\frac{\mu_0}{5}}{1+\gam}x\qquad \tx{in $\Om\cap\nD_{\frac{\eps}{2}}$}
\end{equation}
for $\mu_0\in (0,1)$ from Definition \ref{definition-10-3}(iv-1),
then equation $\mcl{N}_{(u,\beta)}(\phi)=0$ coincides with Eq. \eqref{2-3} in $\nOm_{{\eps}/{10}}$.

Define
\begin{equation*}
v(x,y):=Ax-\psi_x(x,y)\qquad  \tx{for $A=\frac{2-\frac{\mu_0}{5}}{1+\gam}$.}
\end{equation*}
Then $v$ satisfies
\begin{equation}
\label{v-rsonic-bc1}
v=0\,\,\,\,\tx{on $\rightsonic=\{x=0\}$},\qquad\,\,\, v_y=0\,\,\,\,\tx{on $\Wedge\cap \der\nOm_{\eps}$},
\end{equation}
because $\der_{\etan}\vphi=\der_{\etan}\rightvphi=0$ on $\Wedge$.

By \eqref{fp-equality} and properties (a), (f), and (g-3) of Lemma \ref{lemma-ext-M-property},
the boundary condition
on $\shock$ in \eqref{12-50}
can be written as
\begin{equation*}
b_1\psi_x+b_2\psi_y+b_0\psi=0\qquad\tx{on $\shock\cap \nD_{\eps}$}
\end{equation*}
for $(b_0, b_1, b_2)$ satisfying that
\begin{equation*}
-\delta^{-1}\le b_j\le -\delta \qquad\tx{on $\shock\cap \nD_{\eps}$}
\end{equation*}
for a constant $\delta\in(0,1)$ depending only on $(\iv, \gam, \beta_*)$.
Then
$|\psi_x|\le C(|\psi_y|+|\psi|)$ on $\shock\cap \nD_{\eps}$ for $C>0$
depending only on $(\iv, \gam, \beta_*)$.
By combining this inequality with estimate \eqref{est-hat-psi-new2015}
given in Lemma \ref{proposition-est-hpsi-new2015}, we have
\begin{equation*}
|\psi_x|\le Cx^{3/2}\qquad\tx{on $\shock\cap \nD_{\eps}$}
\end{equation*}
for $C>0$ depending only on $(\iv, \gam, \beta_*)$.
Then we can fix a small constant  $\eps_{\rm{fp}}^{(3)}$ depending only on $(\iv, \gam,  \beta_*)$
so that, if
\begin{equation}
\label{eps-fp-2}
0\le \eps \le \eps_{\rm{fp}}^{(3)},
\end{equation}
we have
\begin{equation}
\label{v-rsonic-bc2}
v\ge 0\qquad \tx{on $\shock\cap \der \nD_{\eps}$.}
\end{equation}
By \eqref{definiiton-iterset-ineq4-rsonic} in Definition \ref{definition-10-3}(iv), we obtain
\begin{equation}
\label{v-rsonic-bc3}
v\ge \frac{4\mu_0\eps }{5(1+\gam)}>0
\qquad\tx{on $\der\nOm_{\eps}\cap\{x=\eps\}.$}
\end{equation}

By Lemma \ref{proposition-est-hpsi-new2015}, $\eps_{\rm fp}^{(3)}$ can be further reduced, depending only on $(\iv, \gam,  \beta_*)$,
so that, if \eqref{eps-fp-2} holds, then
\begin{equation*}
\zeta_1(\frac{\psi_x}{x^{3/4}})=\frac{\psi_x}{x^{3/4}},\quad
\zeta_1(\frac{\psi_y}{(\gam+1)N_4x})=\frac{\psi_y}{(\gam+1)N_4x}\qquad\,\,\tx{in $\nOm_{\eps_{\rm fp}^{(3)}}$}
\end{equation*}
for $\zeta_1$ given by \eqref{4-7}. This implies that
\begin{equation*}
O_j^{\rm mod}(\psi_x,\psi_y,x,y)=O_j(\psi_x, \psi_y, \psi, x,y)\qquad \tx{in $\nOm_{\eps_{\rm fp}^{(3)}}$ for all $j=1,\cdots, 5$},
\end{equation*}
for  $O_j^{\rm mod}$ and $O_j$ defined by \eqref{def-O-mod} and \eqref{prelim5-5}, respectively.

By \eqref{12-57} and \eqref{fp-equality}, equation $\mcl{N}_{(u,\beta)}(\hat{\phi})=0$ in $\nOm_{\eps/2}$
becomes $\mcl{N}_{(u,\beta)}^{\rm polar}(\psi)=0$ in the $(x,y)$--coordinates given by \eqref{coord-n}
for $\mcl{N}_{(u,\beta)}^{\rm polar}$ defined by \eqref{definition-nlop-N2}.
We differentiate $\mcl{N}^{\rm polar}_{(u,\beta)}(\psi)=0$ with respect to $x$ in
$\nOm_{\eps/2}$ and then rewrite the resulting equation as an equation for $v(x,y)$ in the following form:
\begin{equation}
\label{equation-v-xy}
a_{11}v_{xx}+a_{12}v_{xy}+a_{22}v_{yy}+a_1 v_x+a_0 v=-A\big((\gam+1)A-1\big)+E(x,y)\qquad
\tx{in}\; \nOm_{\eps/2},
\end{equation}
where
\begin{equation*}
\begin{split}
&a_{ij}=a_{ij}(D_{(x,y)}\psi, x, y)\qquad\tx{for $a_{ij}(D_{(x,y)}\psi, x, y)$ given by \eqref{definition-nlop-N2}},\\
&a_1=1-(\gam+1)\left(\zeta_1(A-\frac vx)+\zeta_1'(A-\frac vx)(\frac vx-v_x+A)\right), \\
&a_0=(\gam+1)\frac{A}{x}\Big(\zeta_1'(A-\frac{v}{x})-\int_0^1\zeta_1'(A-s \frac vx)\;\dd s\Big),\\
&E(x,y)=\psi_{xx}\der_x \hat O_1+\psi_{xy}\der_x\hat O_2+\psi_{yy}\der_x \hat O_3-\psi_{xx}\hat O_4-\psi_x\der_x\hat O_4+\psi_{xy}\hat O_5+\psi_y\der_x \hat O_5,\\
&\hat O_j(x,y)=O_j(\psi_x(x,y), \psi_y(x,y), \psi(x,y), x,y) \quad\tx{for $j=1,\cdots, 5$}.
\end{split}
\end{equation*}

By Lemma \ref{lemma-8-3}(a), Eq. \eqref{equation-v-xy} is strictly elliptic in $\nOm_{\eps/2}$.
Estimate \eqref{est-hat-psi-new2015} given in Lemma \ref{proposition-est-hpsi-new2015}
implies that $a_{ij}, a_1, a_0\in C(\ol{\Om}\setminus \{x=0\})$.
Since $\zeta_1''\le 0$ by \eqref{7-c5-a}, $ a_0v \ge 0$ in $\nOm_{\eps/2}$.
By \eqref{prelim5-5} and \eqref{est-hat-psi-new2015},
there exists a constant $C>0$ depending on $(\iv, \gam,  \beta_*)$ such
that $|E(x,y)|\le Cx$ in $\nOm_{\eps/2}$.
Therefore, we can fix a small constant  $\eps_{\rm{fp}}^{(4)}$ depending only
on $(\iv, \gam,  \beta_*)$ so that, if
\begin{equation}
\label{eps-fp-3}
0\le \eps \le \eps_{\rm{fp}}^{(4)},
\end{equation}
then $-A\big((\gam+1)A-1\big)+E(x,y)<0$ in $\nOm_{\eps/2}$. Thus, for such $\eps$, we have
\begin{equation}
\label{inequality-for-v-new2018}
a_{11}v_{xx}+a_{12}v_{xy}+a_{22}v_{yy}+a_1 v_x+a_0 v<0\qquad\tx{in $\nOm_{\eps/2}$.}
\end{equation}

By properties \eqref{v-rsonic-bc1}, \eqref{v-rsonic-bc2}--\eqref{v-rsonic-bc3},
and \eqref{inequality-for-v-new2018}, we can apply the maximum principle and
Hopf's lemma to conclude that, if
\begin{equation}
0<\eps<\min\{\eps_{\rm{fp}}^{(3)},\eps_{\rm{fp}}^{(4)}\},
\end{equation}
then  $v\ge 0$ in $\nOm_{\eps/2}$, which is equivalent to stating that
\begin{equation*}
\psi_x(x,y)\le \frac{2-\frac{\mu_0}{5}}{1+\gam}\qquad\tx{in $\nOm_{\eps/2}$}.
\end{equation*}

Next, we show that $\psi_x\ge  -\frac{2-\frac{\mu_0}{5}}{1+\gam} x$ in $\nOm_{\eps/2}$.
Since $\der_{\rightvec}(\ivphi-\rightvphi)=0$, we obtain from \eqref{monotonicity-of-fp} that
\begin{equation}
\label{monotonicity-fp}
\der_{\rightvec}\psi=\der_{\rightvec}(\vphi-\ivphi)\ge 0\qquad \tx{in} \,\, \Om.
\end{equation}
By \eqref{cov-right}, $\der_{\rightvec}\psi$ is represented as
\begin{equation}
\label{Psi-expression}
\der_{\rightvec}\psi=\psi_x\cos y+\frac{\sin y}{\rightc-x} \psi_y\qquad\tx{in}\,\,\nOm_{\rightc/2}.
\end{equation}
By Remark \ref{remark-a}(i)--(ii), we can fix a small constant $\eps_{\rm fp}^{(5)}>0$
depending only on $(\iv, \gam,  \beta_*)$ such that
$
\nOm_{\eps_{\rm fp}^{(5)}}\subset \{(x,y): x\in(0, \eps_{\rm fp}^{(5)}), 0<y<\frac{\pi}{2}-\sigma_0\}
$
for some constant $\sigma_0>0$ that is chosen depending only on $(\iv, \gam)$.
Then it follows from estimate \eqref{est-hat-psi-new2015} given
in Lemma \ref{proposition-est-hpsi-new2015} and \eqref{monotonicity-fp}--\eqref{Psi-expression}
that there exists a constant $C>0$ depending only on $(\iv, \gam,  \beta_*)$ such that
\begin{equation*}
\psi_x\ge -\tan(\frac{\pi}{2}-\sigma_0)\psi_y\ge -Cx^{3/2}\qquad\tx{in}\,\,\, \nOm_{\eps_{\rm fp}^{(5)}}.
\end{equation*}
Therefore, $\eps_{\rm fp}^{(5)}$ can be further reduced, depending
only on $(\iv, \gam,  \beta_*)$,
so that the inequality above implies
\begin{equation*}
\psi_x\ge -\frac{2-\frac{\mu_0}{5}}{1+\gam}x\qquad\tx{in}\,\, \nOm_{\eps_{\rm fp}^{(5)}}.
\end{equation*}

We finally conclude that $\vphi$ satisfies \eqref{fp-eqn-coinc-ne-rsonic},
provided that parameter $\eps$ in Definition \ref{definition-10-3} satisfies
\begin{equation}
\label{fp-eps-rsonic-eqn-coinc-cond}
0<\eps\le\min\{\eps_{\rm{fp}}^{(3)},\eps_{\rm{fp}}^{(4)}, \eps_{\rm{fp}}^{(5)} \}.
\end{equation}
Therefore, equation $\mcl{N}_{(u,\beta)}(\phi)=0$ coincides with Eq. \eqref{2-3} in $\nOm_{\eps/10}$,
provided that condition \eqref{fp-eps-rsonic-eqn-coinc-cond} holds.

\smallskip
{\textbf{5-3.}} {\emph{Equation $\mcl{N}_{(u,\beta)}(\phi)=0$ near $\leftsonic$}}.
In $\oOm_{\eps}:=\Om\cap \oD_{\eps}$, let the $(x,y)$--coordinates be defined by \eqref{coord-o}.

By \eqref{definition-xpbeta-new}--\eqref{10-c1},
there exists a small constant $\eps_{\rm fp}^{(6)}>0$ depending only on $(\iv, \gam)$ so that,
if $x_{P_{\beta}}< \frac{\eps_{\rm fp}^{(6)}}{10}$, then
$\beta<\betasonic+\frac 12\min\{\sigma_3, \hat{\delta}\}$ for $\hat{\delta}>0$ from Lemma \ref{lemma-est-hpsi-nrlsonic}(ii)
and $\sigma_3$ from Proposition \ref{proposition-sub9}.

Assume that parameter $\eps$ in Definition \ref{definition-10-3} satisfies
\begin{equation}
\label{fp-eps-lsonic-eqn-coinc-cond1}
0<\eps<\eps_{\rm fp}^{(6)},
\end{equation}
and suppose that $x_{P_{\beta}}< \frac{\eps}{10}$.
By \eqref{definition-coeff-itereqn-near-pt-new} and \eqref{12-57}, if we can show that
\begin{equation}
\label{fp-eqn-coinc-ne-lsonic}
\big|\psi_x(x,y)\big|<\frac{2-\frac{\mu_0}{5}}{1+\gam}x\qquad\,\, \tx{in $\Om\cap\oD_{\eps/2}$},
\end{equation}
then it follows from Lemma \ref{lemma1-coeff-iter-eqn-new}(c-1) that
equation $\mcl{N}_{(u,\beta)}(\phi)=0$ coincides with Eq. \eqref{2-3} in $\oOm_{\eps/10}$.
To prove \eqref{fp-eqn-coinc-ne-lsonic}, we can
mostly repeat the argument in Step 5-2 by using Lemma \ref{lemma-est-hpsi-nrlsonic}(i)--(ii)
and the positivity of $\der_{\leftvec}(\vphi-\ivphi)$ in $\Om$ given in \eqref{monotonicity-of-fp},
instead of Lemma \ref{proposition-est-hpsi-new2015} and the positivity of $\der_{\rightvec}(\vphi-\ivphi)$ in $\Om$.
Then there exists a small constant $\eps_{\rm fp}^{(6)}>0$ depending only on $(\iv, \gam)$ such that,
if $\eps$ satisfies condition \eqref{fp-eps-lsonic-eqn-coinc-cond1},
then equation $\mcl{N}_{(u,\beta)}(\phi)=0$ coincides with Eq. \eqref{2-3} in $\oOm_{\eps/10}$.

If parameter $\eps$ in Definition \ref{definition-10-3} satisfies condition \eqref{fp-eps-lsonic-eqn-coinc-cond1},
and if  $x_{P_{\beta}}\ge \frac{\eps}{10}$, then it follows from Lemma \ref{lemma-6-1}(i) that
equation $\mcl{N}_{(u,\beta)}(\phi)=0$ coincides with Eq. \eqref{2-3} in $\oOm_{\eps/10}$.

\smallskip
For the rest of the proof, parameters $(\eps, \delta_1)$ in Definition \ref{definition-10-3} satisfy
\begin{equation}
\label{fp-cond-eps-delta-final}
0<\delta_1< \delta_{\rm fp},\qquad 0<\eps<\min\{\eps_{\rm fp}^{(j)}\,:\, j=1,\cdots,6\},
\end{equation}
where $\delta_{\rm fp}$ is from \eqref{fp-condition-delta1}.

\smallskip
{\textbf{6.}}  It remains to check that properties (i-1) and (iii) of Definition \ref{def-regular-sol} hold.

\smallskip
{\emph{Verification of property {\rm (iii)} of Definition {\rm \ref{def-regular-sol}}}}.
In Step 5, we have shown that Eq. \eqref{2-3}  coincides with equation $\mcl{N}_{(u,\beta)}(\phi)=0$ in $\Om$.
Therefore, it directly follows from Lemma \ref{lemma-6-1}(a) and Lemmas \ref{proposition-est-hpsi-new2015}--\ref{lemma-est-hpsi-nrlsonic}
that Eq. \eqref{2-3} is strictly elliptic in $\ol{\Om}\setminus (\ol{\leftsonic}\cup\ol{\rightsonic})$.
This proves that property (iii) of Definition \ref{def-regular-sol} holds,
because Eq. \eqref{2-1} is equivalent to \eqref{2-3} in $\Om$.

\smallskip
{\emph{Verification of property {\rm (i-1)} of Definition {\rm \ref{def-regular-sol}}}}.  The strict ellipticity
of Eq. \eqref{2-3} in $\ol{\Om}\setminus (\ol{\leftsonic}\cup\ol{\rightsonic})$ implies
\begin{equation*}
 \frac{|\der_{\bm\nu}\vphi(\bmxi)|^2}{c^2(|\nabla \vphi(\bmxi)|^2, \vphi(\bmxi), \bmxi)}
 \le \frac{|\nabla\vphi(\bmxi)|^2}{c^2(|\nabla \vphi(\bmxi)|^2, \vphi(\bmxi), \bmxi)}<1
 \qquad\tx{on $\shock\setminus(\ol{\leftsonic}\cup\ol{\rightsonic})$}.
\end{equation*}
for a unit normal vector $\bm\nu$ to $\shock$.
We have shown in Step 4-1 that $\vphi$ satisfies the Rankine-Hugoniot condition \eqref{3-a3} on $\shock$.
Define $M:=\frac{|\der_{\bm\nu}\vphi(\bmxi)|}{c^2(|\nabla \vphi(\bmxi)|^2, \vphi(\bmxi), \bmxi)}$
and $M_{\infty}:=|\der_{\bm\nu}\ivphi(\bmxi)|$.
We substitute $M_{\mcl{O}}=M$ into the left-hand side of \eqref{1-1} in the proof of Lemma \ref{lemma-2-0}.
Then, by repeating the argument right after \eqref{1-1} in the proof of Lemma \ref{lemma-2-0}, we obtain that
$M_{\infty}>1$ on $\shock$, which yields that
\begin{equation}
\label{super-on-shock}
|D\ivphi(\bmxi)|>1\qquad\,\,\tx{on $\shock$}.
\end{equation}
By the definition of $\ivphi$ given in \eqref{def-uniform-ptnl-new},
\eqref{super-on-shock} implies that $\bmxi \not\in \ol{B_1(\Oi)}$ for all ${\bm\xi\in \shock}$.
Furthermore, $\{\lefttop,\righttop\}\not\subset \ol{B_1(\Oi)}$,
because $\lefttop$ and $\righttop$ lie on $\leftshock$ and $\rightshock$, respectively.

Now it remains to show that $\xi_1^{\lefttop}\le \xi_1\le \xi_1^{\righttop}$ for all $\bmxi=(\xi_1, \xi_2)\in \shock$.
Since we have shown that $\vphi$ satisfies properties (i-2), (i-4), and (ii)--(v) of Definition \ref{def-regular-sol}
in the previous steps, we can repeat the proof of Lemma \ref{lemma-step1-1} to show that $\vphi$ satisfies the directional
monotonicity properties \eqref{3-c2}--\eqref{3-c3}.
Then, by repeating the proof of Proposition \ref{proposition-3},
we obtain a function $\fshock$ satisfying
\begin{equation*}
\shock=\{{\bmxi}=(\xi_1, \xi_2)\,:\, \xi_2=\fshock(\xi_1),\,\,\xi_1^{\lefttop}<\xi_2<\xi_1^{\righttop}\}.
\end{equation*}
Therefore, property (i-1) holds.

With these, we complete the proof.
\end{proof}
\end{proposition}

\section{Existence of Admissible Solutions for All $(\iv, \beta)\in  \mathfrak{R}_{\rm{weak}}$}
\label{subsec-pf-mthm-p1}

In order to prove the existence of admissible solutions for all $(\iv, \beta)\in \mathfrak{R}_{\rm{weak}}$,
we employ the Leray-Schauder fixed point index and its generalized homotopy invariance property.

\subsection{Leray-Schauder degree theorem}
\label{subsec-LS-degree}

\begin{definition}[Compact map]
\label{definition-compact-mapping}
Let ${\rm X}$ and ${\rm Y}$ be two Banach spaces. For an open subset $G$ in ${\rm X}$, a map ${\rm{\bf{f}}}:\ol{G}\rightarrow Y$ is called {\emph{compact}} if

\smallskip
\begin{itemize}
\item[(i)] ${\rm{\bf{f}}}$ is continuous{\rm ;}

\smallskip
\item[(ii)] ${\rm{\bf{f}}}({\rm U})$ is precompact in ${\rm Y}$ for any bounded subset ${\rm U}$ of $\ol{G}$.
\end{itemize}
\end{definition}

\begin{definition}
\label{definition-V-class}
Let $G$ be an open bounded set in a Banach space ${\rm X}$.
Denote by $V(G, {\rm X})$ the set of all maps ${\bf{\rm f}}: \ol{G}\rightarrow {\rm{X}}$ satisfying the following{\rm :}

\smallskip
\begin{itemize}
\item[(i)] ${\bf{\rm f}}$ is compact in the sense of Definition {\rm \ref{definition-compact-mapping}}{\rm ;}

\smallskip
\item[(ii)] ${\bf{\rm f}}$ has no fixed points on the boundary $\der G$.
\end{itemize}
\end{definition}

\begin{definition}
Two maps ${\rm{\bf f}}, {\rm{\bf g}} \in V(G,{\rm X})$ are called {\emph{compactly homotopic}}
on $\der G$ if there exists a map ${\rm{\bf H}}$ with the following properties{\rm :}

\smallskip
\begin{itemize}
\item[(i)] ${\rm{\bf H}}: \ol{G}\times [0,1]\rightarrow {\rm X}$ is compact in the sense of
Definition {\rm \ref{definition-compact-mapping}}{\rm ;}

\smallskip
\item[(ii)] ${\rm{\bf H}}({\bf x}, \tau)\neq {\bf x}$ for all $({\bf x},\tau)\in \der G\times [0,1]${\rm ;}

\smallskip
\item[(iii)] ${\rm{\bf H}}({\bf x},0)={\rm{\bf f}}({\bf x})$ and ${\rm{\bf H}}({\bf x},1)={\rm{\bf g}}({\bf x})$ in $\ol{G}$.
\end{itemize}
We write $\der G: {\rm{\bf f}}\cong  {\rm{\bf g}}$ if ${\rm{\bf f}}$ and ${\rm{\bf g}}$ are compactly homotopic on $\der G$,
and call ${\rm{\bf H}}$
a compact homotopy.
\end{definition}

\begin{theorem}[Leray-Schauder degree theorem]
\label{theorem-LS-theory}
Let $G$ be an open bounded set in a Banach space ${\rm X}$. Then, to each map ${\rm{\bf f}}\in V(G,X)$,
a unique integer ${\rm{\bf Ind}}({\rm{\bf f}}, G)$ can be assigned with the following properties{\rm :}

\smallskip
\begin{itemize}
\item[(i)] If ${\rm{\bf f}}({\bf x})\equiv {\bf x}_0$ for all ${\bf x}\in\ol{G}$ and some fixed ${\bf x}_0\in G$, then
${\rm{\bf Ind}}({\rm{\bf f}}, G)=1${\rm ;}

\smallskip
\item[(ii)] If ${\rm{\bf Ind}}({\rm{\bf f}}, G)\neq 0$, then there exists ${\bf x}\in G$ such that ${\rm{\bf f}}({\bf x})={\bf x}${\rm ;}

\smallskip
\item[(iii)] ${\rm{\bf Ind}}({\rm{\bf f}}, G)=\sum_{j=1}^n {\rm{\bf Ind}}({\rm{\bf f}}, G_j)$,
whenever $f\in V(G, X)\cap (\cap_{j=1}^n V(G_j, X))$, where $G_i\cap G_j=\emptyset$ for $i\neq j$ and $\ol{G}=\cup_{j=1}^n \ol{G_j}${\rm ;}

\smallskip
\item[(iv)] If $\der G: {\rm{\bf f}} \cong {\rm{\bf g}}$, then ${\rm{\bf Ind}}({\rm{\bf f}}, G)={\rm{\bf Ind}}({\rm{\bf g}}, G)$.
\end{itemize}
Such a number ${\rm{\bf Ind}}({\rm{\bf f}}, G)$ is called the {\emph{fixed point index}} of ${\rm{\bf f}}$ over $G$.
\end{theorem}

A generalized homotopy invariance of the fixed point index is given in the following theorem:
\begin{theorem}[\cite{Zeidler}, \S 13.6, A4*]
\label{theorem-hom-inv-fp}
Let ${\rm{X}}$ be a Banach space, and let $t_2>t_1$.
Let ${\rm U}\subset {\rm X}\times [t_1,t_2]$, and let $U_t=\{{\bf x}\,:\,({\bf x},t)\in {\rm U}\}$. Then
\[
{\rm{\bf Ind}}({\rm{\bf h}}(\cdot, t), U_t)=const.\qquad\tx{for all $t\in[t_1, t_2]$},
\]
provided that ${\rm U}$ is bounded and open in ${\rm X}\times [t_1, t_2]$, and map ${\rm{\bf h}}:\ol{\rm{U}}\rightarrow {\rm X}$
is compact in the sense of Definition {\rm \ref{definition-compact-mapping}}
with ${\rm{\bf h}}({\bf x},t)\neq {\bf x}$ on $\der{\rm U}$.
\end{theorem}

\medskip
\subsection{Proof of Theorem {\rm \ref{theorem-0}}}
\label{subsubsection-proof-admiss-exist}
In this subsection, we complete the proof of Theorem {\rm \ref{theorem-0}}.

\smallskip
{\emph{{\textbf{Parameters $(\alp, \eps, \delta_1, \delta_3, N_1)$ in Definition {\rm \ref{definition-10-3}}}}}}:
Let parameters $(\alp, \eps, \delta_1, \delta_3, N_1)$ in {\rm{Definition \ref{definition-10-3}}} be fixed
as in Definition \ref{definition-iterationmap}.
We further reduce $(\eps, \delta_1)$ depending only on $(\iv, \gam, \beta_*)$ so that
Proposition \ref{lemma-8-2} implies that, for each $\beta\in(0, \beta_*]$,
$u\in \mcl{K}(\beta)$ is a fixed point of $\Itr(\cdot, \beta):\ol{\mcl{K}}(\beta)\rightarrow C^{2,\alp}_{(*,\alp_1)}(\iter)$
if and only if $\vphi$, defined by \eqref{def-vphi} in Definition  {\rm \ref{definition-Gset-shocks-new}},
yields an admissible solution corresponding to $(\iv, \beta)\in\mathfrak{R}_{\rm weak}$
in the sense of Definition {\rm \ref{def-regular-sol}}.
\smallskip

In the proof of Theorem  {\rm \ref{theorem-0}}, we adjust $N_1$ and choose $\delta_2$ so that
$\Itr(\cdot, \beta)$ has a fixed point in $\mcl{K}(\beta)$ for each $\beta\in(0,\beta_*]$.
Then the existence of
an admissible solution for each $(\iv, \beta)\in \mathfrak{R}_{\rm weak}\cap \{\beta\le \beta_*\}$
follows from Proposition \ref{lemma-8-2}. This proves Theorem \ref{theorem-0}, since $\beta_*$ is arbitrarily chosen
in $(0,\betadet)$.

\smallskip
{\emph{{\textbf{Further adjustment of $\delta_3$ in Definition {\rm \ref{definition-10-3}}}}}}:
Note that, if parameter $N_1$ in Definition \ref{definition-10-3} is adjusted such that
the new choice of $N_1$ is greater than the previous one, all the properties stated previously hold.
Then we  choose $N_1$ greater than the previous choice in the proof of Theorem  {\rm \ref{theorem-0}}.
Also, once parameters $(N_1, \delta_2)$ are fixed,
$\delta_3$ can be adjusted to satisfy the conditions of $\delta_3$ in
Lemmas \ref{lemma-positivity-hatphi-2018}--\ref{proposition-est-hpsi-new2015}.
As long as the new choice of $\delta_3$ is less than the previous choice,
all the properties stated previously hold.
Since $N_1$ is adjusted to be greater than the previous one, the new choice of $\delta_3$
is less than the previous one.
Since the previous choice of  $(\alpha, \eps, \delta_1, \delta_2, N_1)$ was independent
of $\delta_3$, we can reduce $\delta_3$ as described above.

\medskip
\begin{proof}
[Proof of Theorem {\rm \ref{theorem-0}}]
The proof is divided into three steps.

\smallskip
{\textbf{1.}} {\emph{Claim {\rm 1:}
The iteration map $\Itr:\ol{\mcl{K}}\rightarrow C^{2,\alp}_{(*,\alp_1)}(\iter)$
defined by Definition {\rm \ref{definition-iterationmap}}
is continuous. Moreover,
$\Itr:\ol{\mcl{K}}\rightarrow C^{2,\alp}_{(*,\alp_1)}$ is compact in the sense
of Definition {\rm \ref{definition-compact-mapping}}. }}
\smallskip

{\textbf{1-1.}} {\emph{Continuity of $\Itr:\ol{\mcl{K}}\rightarrow C^{2,\alp}_{(*,\alp_1)}$}}.
Suppose that $\{(u_j,\beta_j)\}_{j=1}^{\infty}\subset \mcl{K}$ converges
to $(u,\beta)$ in $C^{2,\alp}_{(*,\alp_1)}(\iter)\times [0, \beta_*]$.
For each $j\in \mathbb{N}$, define
$(\Om_j, \gshock^{(j)}):=(\Om(u_j, \beta_j), \gshock^{(u_j,\beta_j)})$
for $\Om(u_j, \beta_j)$ and $\gshock^{(u_j,\beta_j)}$ given
by Definition \ref{definition-Gset-shocks-new}.
By Lemma \ref{proposition-existence-hpsi-bvp}, the nonlinear boundary value problem \eqref{12-50}
associated with $(u_j, \beta_j)$ has a unique
solution $\hat{\phi}^{(j)}\in C^2(\Om_j)\cap C^1(\ol{\Om_j}\setminus(\ol{\Gam^{\mcl{O},j}_{{\rm sonic}}}\cup\ol{\rightsonic}))\cap C^0(\ol{\Om_j})$,
where $\Gam^{\mcl{O},j}_{{\rm sonic}}$ is $\leftsonic$ corresponding to $(\iv, \beta_j)$.
For such $\hat{\phi}^{(j)}$, define
\begin{equation}
\label{definition-hatw-j-2018}
\hat w^{(j)}:=(\hat{\phi}^{(j)}+\rightvphi-\vphi^*_{\beta_j})\circ(\mcl{G}_1^{\beta_j})^{-1}
\end{equation}
for $\mcl{G}_1^{\beta_j}$ and $\vphi^*_{\beta_j}$ defined by \eqref{12-16-mod} and \eqref{12-32}, respectively.

Let $\hat{u}_j$ be given by \eqref{definition-sol-ubvp}  associated with $(u_j, \beta_j, \hat{\phi}^{(j)})$.
Then Definition \ref{definition-Gset-shocks-new}(ii) implies that
\begin{equation}
\label{w-hatu-relation2018}
\hat w^{(j)}=\hat u_j\circ G_{2,\gshock^{(j)}}
\end{equation}
for $G_{2,\gshock^{(j)}}$ defined by \eqref{7-b9}.

For each $\hat{w}^{(j)}$, let $\hgshock^{(j)}$ be given from \eqref{hgshock} with $\hat w=\hat w^{(j)}$.
We also define $\Om$, $\gshock$, $\hat{\phi}$,
$\hat w$, $\hat u$, and $\hgshock$, similarly associated with $(u,\beta)\in\ol{\mcl{K}}$.

By Lemma \ref{lemma-7-4}(d), we have
\begin{equation}
\label{gshock-cov1-2018}
\gshock^{(j)} \rightarrow \gshock\qquad\tx{in $C^{1,\alp}([-1,1])$}.
\end{equation}
Fix a compact set $K \subset \mcl{G}_1^{\beta}(\Om)=\{(s,t')\,:\,-1<s<1, 0<t'<\gshock(s)\}$.
Then there exists a constant $\sigma_K\in(0,1)$ depending only on $K$ such that $K\subset \{s\ge -1+\sigma_K\}$.
Thus, by Lemma \ref{lemma-7-4}(g), there exists a constant $C_K>1$ depending only on $(\iv, \gam, \beta_*)$ and $K$
such that, for any $(u^{\sharp}, \beta^{\sharp})\in \ol{\mcl{K}}$,
\begin{equation}
\label{gshock-cov2-2018}
C_K^{-1}<\gshock^{\sharp}(s)<C_K\qquad\tx{for all $(s,t')\in K$}.
\end{equation}
By \eqref{7-b9} and \eqref{gshock-cov1-2018}--\eqref{gshock-cov2-2018}, we have
\begin{equation}
\label{G2-gshock-conv-2018}
G_{2,\gshock^{(j)}}\rightarrow G_{2,\gshock}\qquad\tx{in $C^{1,\alp}(K)$}.
\end{equation}
This implies that there exists a compact set $\mcl{Q}_K\subset \iter$ such
that $G_{2,\gshock^{(j)}}(K)\subset \mcl{Q}_K$ for all $j$,
and $G_{2,\gshock}(K)\subset \mcl{Q}_k$.
By Corollary \ref{corollary-u-convergence}(b),
$\hat u_j$ converges to $\hat{u}$ in $C^2(\mcl{Q}_K)$.
Therefore, it follows from \eqref{w-hatu-relation2018} and \eqref{G2-gshock-conv-2018} that
\begin{equation}
\label{wj-unif-conv-cpt-sets2018}
\hat{w}^{(j)}\rightarrow \hat w\qquad\tx{in $C^{1,\alp}(K)$}.
\end{equation}
Since $K$ is an arbitrary compact subset of $\mcl{G}_1^{\beta}(\Om)$,
we conclude that $\hat w_j$ converges to $\hat w$ in $C^{1,\alp}$
for any compact subset of $\mcl{G}_1^{\beta}(\Om)$.

By \eqref{definition-hatw-j-2018}, \eqref{wj-unif-conv-cpt-sets2018}, and
Lemmas \ref{lemma-12-1-mod} and \ref{proposition-est-hpsi-new2015}--\ref{lemma-est-hpsi-nrlsonic},
we can apply Proposition \ref{proposition-prop-ext-map}(a-3) to obtain the convergence of
sequence $\{\mcl E_{\gshock^{(j)}}(\hat w^{(j)})\}$ to $\mcl E_{\gshock}(\hat{w})$
in $C^{2,\alp}(R_{(1+\frac{\kappa}{2})\gshock}\cap \{b_1<s<b_2\})$ for any $b_1$ and $b_2$
with $-1<b_1<b_2<1$, where $\kappa\in(0,\frac 13]$ is from Definition \ref{definition-extension-map}.
Note that, for any $\sigma\in(0,1)$,
$\{(s, \hgshock^{(j)}(s))\,:\,-1+\sigma<s<1-\sigma\}\subset R_{(1+\frac{\kappa}{2})\gshock}$
holds for all $j$ sufficiently large depending on $\sigma$.
Therefore, by using the $C^2$--estimate of $\hgshock$ given in Lemma \ref{lemma-existence-hat-g} and \eqref{t-deriv-w-ext},
it can be directly checked that $\{\hgshock^{(j)}\}$ converges to $\hgshock$ in $C^2([-1+\sigma, 1-\sigma])$
for any $\sigma\in (0,1)$. Then we obtain from \eqref{12-77} that
\begin{equation}
\label{hgshock-conv-2018}
\hgshock^{(j)}\rightarrow \hgshock \qquad\tx{in $C^{2,\alp}_{(*,\alp_1)}((-1,1))$.}
\end{equation}
By \eqref{12-78}, \eqref{hgshock-conv-2018}, and properties (a-3), (b-3), and (c-3) of Proposition \ref{proposition-prop-ext-map},
we conclude that $\til{u}_j:=\Itr(u_j,\beta_j)$ converges to $\til u=\Itr(u,\beta)$ in $C^{2,\alp}_{(*,\alp_1)}(\iter)$.
This implies that
$\Itr:\ol{\mcl{K}}\rightarrow C^{2,\alp}_{(*,\alp_1)}$ is continuous.

\smallskip
{\textbf{1-2.}} {\emph{Compactness of $\Itr:\ol{\mcl{K}}\rightarrow C^{2,\alp}_{(*,\alp_1)}$}}.
Let $U$ be a subset of $\ol{\mcl{K}}\subset C^{2,\alp}_{(*,\alp_1)}(\iter)\times[0,\beta_*]$. Then $U$ is bounded in $C^{2,\alp}_{(*,\alp_1)}(\iter)\times [0,\beta_*]$.
Since $C^{2,2\alp}_{(*,1)}(\iter)$ is compactly embedded into $C^{2,\alp}_{(*,\alp_1)}(\iter)$,
Lemma \ref{lemma-12-6}(b) implies that $\Itr(U)$ is pre-compact
in  $C^{2,\alp}_{(*,\alp_1)}(\iter)$.
From this property, combined with the continuity of $\Itr$ proved in the previous step,
we conclude that $\Itr:\ol{\mcl{K}}\rightarrow C^{2,\alp}_{(*,\alp_1)}(\iter)$
is compact in the sense of Definition \ref{definition-compact-mapping}.
This verifies {\emph{Claim} 1}.

\smallskip
{\textbf{2.}} {\emph{Claim {\rm 2:}
In Definition {\rm \ref{definition-10-3}}, $N_1$ can be increased and $\delta_2>0$ can be fixed such that,
for any $\beta\in(0,\beta_*]$, no fixed point of $\mcl{I}(\cdot, \beta)$ lies on
boundary $\der \mcl{K}(\beta)$ of $\mcl{K}(\beta)$,
where $\der\mcl{K}(\beta)$ is considered relative to space $C^{2,\alp}_{(*,\alp_1)}(\iter)$.
Furthermore, the choices of $(N_1, \delta_2)$ depend only on $(\iv, \gam,  \beta_*)$.
}}
\smallskip

{\textbf{2-1.}}
Let $\mcl{I}(u,\beta)=u$ for some $(u,\beta)\in \ol{\mcl{K}}$,
and let $\vphi=\vphi^{(u,\beta)}$ be given by \eqref{def-vphi}.
We extend $\vphi$ onto $\Lambda_{\beta}$ by \eqref{1-24}
if $\beta<\betasonic$,
and by \eqref{1-24ab} if $\beta\ge \betasonic$.
By Proposition \ref{lemma-8-2}, $\vphi$ is an admissible solution corresponding
to $(\iv, \beta)\in\mathfrak{R}_{\rm weak}$ in the sense of Definition \ref{def-regular-sol}.

In order to verify Claim 2, we need to show the following:

\smallskip
\begin{itemize}
\item[-] $u$ satisfies the strict inequality given in condition (i) of Definition  \ref{definition-10-3};

\smallskip
\item[-] $\vphi$ satisfies all the strict inequalities given in conditions (iii)--(vi) given in Definition  \ref{definition-10-3}.
 \end{itemize}

\smallskip
{\textbf{2-2.}} {\emph{The strict inequalities in condition {\rm (i)} of Definition {\rm \ref{definition-10-3}}}}:
Note that $N_1$ satisfies that $N_1\ge N_1^{\rm (a)}$ for $N_1^{\rm (a)}$ from Corollary \ref{corollary-admisol-iterset}.
Therefore, $u$ satisfies the strict inequality given in condition (i) of Definition  \ref{definition-10-3}.

\smallskip
 {\textbf{2-3.}}  {\emph{The strict inequalities in conditions {\rm (iii)} and {\rm (v)}--{\rm (vi)} of Definition {\rm \ref{definition-10-3}}}}.
 In conditions (iii) and (v)--(vi) of Definition  \ref{definition-10-3},
 constants $(N_2,  \til{\zeta}, \til{\mu}, a_*, C)$ are fixed so that any admissible solution satisfies
 the strict inequalities in conditions (iii) and (v)--(vi) of Definition  \ref{definition-10-3}
 by Propositions \ref{proposition-distance} and \ref{lemma-12-3},
 Remark \ref{remark-ellipticity-newly-added}, and Lemma \ref{lemma-step3-1}.

\smallskip
{\textbf{2-4.}} {\emph{The strict inequalities in condition {\rm (iv)} of Definition {\rm \ref{definition-10-3}}}}.
Suppose that $0<\beta<\frac{\delta_1}{N_1^2}$.
Then $\mathscr{K}_2(\beta)$ defined by \eqref{definition-K2b-2018} satisfies that $\mathscr{K}_2(\beta)<0$ for any $\delta_2>0$.
Moreover, $\vphi$ satisfies \eqref{definiiton-iterset-ineq1} in the whole domain $\Om$ by Definition \ref{def-regular-sol}(iv),
the strong maximum principle, and Hopf's lemma.
The strict inequalities in \eqref{definiiton-iterset-ineq2}--\eqref{definiiton-iterset-ineq3}
are satisfied by Lemma \ref{lemma-step1-1}.

\smallskip
Next, suppose that $\beta\ge \frac{\delta_1}{N_1^2}$.
Then it follows directly from \eqref{def-uniform-ptnl-new}
that $\rightvphi-\leftvphi$ is a nontrivial linear function.
By Definition \ref{def-regular-sol}(iv),  $\psi=\vphi-\max\{\leftvphi, \rightvphi\}\ge 0$
in $\Om$. Since $\vphi=\leftvphi$ on $\leftsonic$,
 $\rightvphi$ on $\rightsonic$, and $\leftvphi-\rightvphi$ is a nonzero function,
 the strong maximum principle and Hopf's lemma apply to $\vphi$,
 so that $\vphi-\leftvphi>0$ and $\vphi-\rightvphi>0$ in $\Om$ hold, which yields that
\begin{equation}
\label{vphi-strict-positivity}
\psi=\vphi-\max\{\leftvphi, \rightvphi\}>0\qquad\tx{in}\,\,\ol{\Om}\setminus (\oD_{\eps/10}\cup \nD_{\eps/10})
\end{equation}
for fixed $\eps>0$ in Definition \ref{definition-10-3}.
By \eqref{vphi-strict-positivity}, Lemmas \ref{lemma-step1-1} and \ref{102},
and the continuous dependence of $(\leftsonic, \leftvphi)$ on $\beta$,
there exists a constant $\sigma>0$ depending only on $(\gam, \iv, \beta_*)$ such that
\begin{equation*}
\psi=\vphi-\max\{\leftvphi, \rightvphi\}>\sigma \qquad\tx{in}\,\,\ol{\Om}\setminus (\oD_{\eps/10}\cup \nD_{\eps/10}).
\end{equation*}

By Lemma \ref{lemma-step1-1}, we also have
\begin{equation*}
\der_{\leftvec}(\ivphi-\vphi)<0\quad\tx{in $\ol{\Om}\setminus \oD_{\eps/10}$},\qquad \,\,
-\der_{\xi_1}(\ivphi-\vphi)<0\quad\tx{in $\ol{\Om}\setminus \nD_{\eps/10}$}.
\end{equation*}
By Corollary \ref{corollary-unif-est-away-sn}, and Propositions \ref{lemma-est-sonic-general-N},
\ref{lemma-est-sonic-general}, \ref{proposition-sub8}, \ref{proposition-sub9}, and \ref{lemma-gradient-est},
the set of admissible solutions corresponding to $(\iv, \beta)\in\mathfrak{R}_{\rm weak}\cap\{\beta\le \beta_*\}$
is uniformly bounded in $C^{1,\alp}$.
Therefore,  there exists a constant $\hat{\sigma}>0$ depending only on $(\gam, \iv, \beta_*)$ such that
\begin{equation*}
\der_{\leftvec}(\ivphi-\vphi)<-\hat{\sigma} \quad\tx{in $\ol{\Om}\setminus \oD_{\eps/10}$},\qquad
-\der_{\xi_1}(\ivphi-\vphi)<-\hat{\sigma} \quad\tx{in $\ol{\Om}\setminus \nD_{\eps/10}$}.
\end{equation*}

Since $\delta_1>0$ is fixed, depending on $(\iv, \gam, \beta_*)$,
we can choose $N_1$ sufficiently large and $\delta_2>0$ sufficiently small, depending only on $(\iv, \gam,  \beta_*, \delta_1, N_1)$, such that
\begin{equation*}
\mathscr{K}_2(\beta)\le \frac{\delta_1\delta_2}{N_1^2}< \min\{\sigma,\hat{\sigma}\}\qquad\tx{ for all $\beta\in[0,\beta_*]$. }
\end{equation*}
With the choices of $(N_1, \delta_2)$, $\vphi$ satisfies \eqref{definiiton-iterset-ineq1}--\eqref{definiiton-iterset-ineq3}
in  Definition \ref{definition-10-3}(iv).

In inequalities \eqref{definiiton-iterset-ineq4-rsonic}--\eqref{definiiton-iterset-ineq8},
parameters $\mu_0$, $\mcl{K}_3(\beta)$, $N_4$, $N_5$, and $\mu_1$ are fixed
so that any admissible solution corresponding to $(\iv, \beta)\in\mathfrak{R}_{\rm weak}\cap\{\beta<\beta_*\}$
satisfies all the strict inequalities.

\smallskip
{\textbf{2-5.}}
With the choices of $(N_1, \delta_2)$ determined in Step 2-4, we conclude that
any fixed point of $\Itr(\cdot, \beta)$ for $\beta\in(0,\beta_*]$ lies in $\mcl{K}(\beta)$.
In the next step, we also show that no fixed point of $\Itr(\cdot, 0)$ lies on $\der \mcl{K}(0)$.

\smallskip
{\textbf{3.}}
Let parameters $(\alp, \eps, \delta_1, \delta_3, N_1)$ in {\rm{Definition \ref{definition-10-3}}}
be fixed as described at the beginning of \S \ref{subsubsection-proof-admiss-exist}.
Let $N_1$ be further adjusted, and let $\delta_2$ be fixed as in Step 2  so that {\emph{Claim}} 2 holds.
Finally, let $\delta_3$ be further adjusted to satisfy the conditions in Lemmas \ref{lemma-positivity-hatphi-2018} and \ref{proposition-est-hpsi-new2015}
as described at the beginning of \S \ref{subsubsection-proof-admiss-exist}.
In particular, let $\delta_3$ be adjusted to satisfy \eqref{delta3-choice} given in the proof of Lemma \ref{lemma-positivity-hatphi-2018}.
With these choices of parameters $(\alp, \eps, \delta_1, \delta_2, \delta_3, N_1)$,
the definition for the iteration set $\mcl{K}$ given in Definition \ref{definition-10-3} is now complete.
\smallskip

{\textbf{3-1.}} {\emph{Claim {\rm 3:} The iteration map $\Itr(\cdot, 0)$ has a unique fixed point $0$ with
\begin{equation*}
{\rm{\bf Ind}}(\Itr(\cdot, 0), \mcl{K}(0))=1.
\end{equation*}
}}

At $\beta=0$, it follows from \eqref{def-uniform-ptnl-new} that $\leftvphi-\rightvphi\equiv 0$,
so that the boundary condition on $\leftsonic\cup\rightsonic$ of the boundary value
problem \eqref{12-50} associated with any $u\in \ol{\mcl{K}}(0)$ becomes homogeneous.
Then it follows from Lemmas \ref{lemma-ext-M-property}(f) and \ref{proposition-existence-hpsi-bvp}
that, for any $u\in \ol{\mcl{K}}(0)$,
the associated boundary value problem \eqref{12-50} has a unique
solution $\hat{\phi}=0$ in $\Om(u,0)$. From this, we have
 \begin{equation*}
 \Itr(u,0)=0\qquad\tx{for all}\,\, u\in \ol{\mcl{K}}(0).
 \end{equation*}
It can directly be checked from Definition \ref{definition-10-3} that
the fixed point $u=0$ of $\Itr(\cdot,0)$ lies in $\mcl{K}(0)$.
Also, we have shown in Step 1 that $\Itr:\ol{\mcl{K}}\rightarrow C^{2,\alp}_{(*,\alp_1)}$ is compact in the sense
of Definition {\rm \ref{definition-compact-mapping}}.
Therefore, the fixed point index ${\rm{\bf Ind}}(\Itr(\cdot, \beta), \mcl{K}(\beta))$ satisfying properties
(i)--(iv) stated in Theorem \ref{theorem-LS-theory} is well defined.
Then Theorem \ref{theorem-LS-theory}(i) implies that
\begin{equation}
\label{fp-index-normal-2018}
{\rm{\bf Ind}}(\Itr(\cdot, 0), \mcl{K}(0))=1.
\end{equation}

\smallskip
{\textbf{3-2.}}  Combining Claim 2 in Step 2 with Claim 3 in Step 3-1, we see that
no fixed point of $\Itr(\cdot, \beta)$ lies on
the boundary $\der \mcl{K}(\beta)$ of $\mcl{K}(\beta)$ for all $\beta\in[0,\beta_*]$.
Then, using \eqref{fp-index-normal-2018} and properties (a) and (d) of Theorem \ref{theorem-hom-inv-fp}, we have
\begin{equation}
\label{index-iter-map}
{\rm{\bf Ind}}(\Itr(\cdot, \beta), \mcl{K}(\beta))={\rm{\bf Ind}}(\Itr(\cdot, 0), \mcl{K}(0))
\qquad \mbox{for all $\beta\in[0,\beta_*]$}.
\end{equation}
By Theorem \ref{theorem-LS-theory}(ii), \eqref{index-iter-map} implies
that $\Itr(\cdot, \beta)$ has a fixed point in $\mcl{K}(\beta)$ for all $\beta\in[0,\beta_*]$.
Then Proposition \ref{lemma-8-2} implies that,
for each $(\iv, \beta)\in \mathfrak{R}_{\rm weak}\cap\{0\le \beta\le \beta_*\}$,
an admissible solution corresponding to $(\iv, \beta)$ exists.
Since $\iv>0$ is arbitrary, and $\beta_*$ is also arbitrary  in $(0,\betadet)$,
we finally conclude that there exists an admissible solution for
any $(\iv, \beta)\in\mathfrak{R}_{\rm weak}$.
This completes the proof of Theorem \ref{theorem-0}.
\end{proof}



\chapter{Optimal Regularity of Admissible Solutions \\ -- Proof of Theorem 2.33}
\label{sec-optimal}
\numberwithin{equation}{chapter}

This chapter is devoted to the complete proof of Theorem \ref{theorem-3}.

\medskip
Let $\vphi$ be an admissible solution corresponding to $(\iv, \beta)\in\mathfrak{R}_{\rm weak}$
in the sense of Definition \ref{def-regular-sol}.
We now prove statements (a)--(e) of Theorem \ref{theorem-3}, respectively.

\smallskip
{\textbf{1.}} {\emph{Proof of statement {\rm (a)} of Theorem {\rm \ref{theorem-3}}}}. It follows from Lemmas \ref{lemma-unif-est1} and
\ref{lemma-unif-est2} that $\shock$ is $C^{\infty}$ in its relative interior,
and $\vphi\in C^{\infty}(\ol{\Om}\setminus \ol{\leftsonic}\cup\ol{\rightsonic})$.
By Definition \ref{definition-domains-np}, $\ol{\leftsonic}$ is a closed portion of a circle
when $\beta<\betasonic$ and becomes a point $P_{\beta}$ when $\beta\ge \betasonic$.
Near $\rightsonic$, we combine Proposition \ref{lemma-est-sonic-general-N} with
the smoothness of $\vphi$ away from $\ol{\leftsonic}\cup\ol{\rightsonic}$
to obtain $\vphi\in C^{1,1}(\ol{\Om}\setminus \ol{\leftsonic})$.

Near $\leftsonic$, we consider two cases separately: (i) $\beta<\betasonic$ and (ii) $\beta\ge \betasonic$.
If $\beta < \betasonic$, it follows from Propositions \ref{lemma-est-sonic-general} and \ref{proposition-sub8}
that $\vphi$ is $C^{1,1}$ up to $\ol{\leftsonic}$.
If $\beta\ge \betasonic$, then Propositions \ref{proposition-sub9} and \ref{lemma-gradient-est} imply
that $\vphi$ is $C^{1,\alp}$ up to $\ol{\leftsonic}=\{P_{\beta}\}$
for some $\alp\in(0,1)$.
This completes the proof of statement (a).

\medskip
{\textbf{2.}} {\emph{Proof of statements {\rm (b)}--{\rm (c)} of Theorem {\rm \ref{theorem-3}}}}.
Let the $(x,y)$--coordinates be defined by \eqref{coord-n} and \eqref{coord-o}
near $\rightsonic$ and $\leftsonic$, respectively.
Define
\begin{equation*}
\psi:=\vphi-\max\{\rightvphi, \leftvphi\}
\end{equation*}
for $\leftvphi$ and $\rightvphi$ given by \eqref{def-uniform-ptnl-new}.
Note that $\psi=\vphi-\rightvphi$ near $\rightsonic$ and $\psi=\vphi-\leftvphi$ near $\leftsonic$.

By \eqref{prelim5-5}, \eqref{eqn-xy-right}, \eqref{5-a4}, Lemma \ref{lemma-est-nrsonic-N},
and Proposition \ref{lemma-est-sonic-general-N},
we can apply the following theorem to $\psi$ near $\rightsonic$:

\begin{theorem}[Theorem 3.1 in \cite{BCF}]
\label{theorem-BCF1}
For constants $r,R>0$, define $Q^+_{r,R}$ by
$$
Q^+_{r,R}:=\{(x,y)\,:\,x\in(0,r),\;\;|y|<R\}.
$$
For positive constants $a,b,M,N$, and $\kappa\in(0,\frac 14)$,
suppose that $\psi\in C(\ol{Q^+_{r,R}})\cap C^2(Q^+_{r,R})$ satisfies
\begin{eqnarray*}
\quad (2x-a\psi_x+O_1)\psi_{xx}+O_2\psi_{xy}+(b+O_3)\psi_{xy}
-(1+O_4)\psi_x+O_5\psi_y=0 &&\text{in $Q^+_{r,R}$},\label{bcf-1}\\
\label{bcf-2}
\quad \psi>0 &&\text{in $Q^+_{r,R}$},\\
\label{bcf-3}
\quad \psi=0 &&\text{on $\der Q^+_{r,R}\cap\{x=0\}$},\\
\label{bcf-4}
\quad -Mx\le \psi_x\le \frac{2-\kappa}{a}x &&\text{in $Q^+_{r,R}$},
\end{eqnarray*}
where terms $O_i(x,y)$, $i=1,\cdots,5$, are continuously differentiable and
\begin{align}
\label{bcf-5}
\frac{|O_1(x,y)|}{x^2}+\frac{|DO_1(x,y)|}{x^2}
+\sum_{k=2}^5 \Big(\frac{|O_k(x,y)|}{x}+
|DO_k(x,y)|\Big)\le N \qquad\mbox{in $Q^+_{r,R}$}.
\end{align}
Then
\begin{equation*}
\psi\in C^{2,\alp}(\ol{Q^+_{r/2,R/2}})\qquad\text{for any $\alp\in(0,1)$},
\end{equation*}
with
\begin{equation*}
\psi_{xx}(0,y)=\frac 1a,\;\;\psi_{xy}(0,y)=\psi_{yy}(0,y)=0\qquad\text{for all $|y|<\frac{R}{2}$}.
\end{equation*}
\end{theorem}

For $\beta\in[0, \betasonic)$, it can be directly checked from the results
in \S \ref{subsubsec-apriori-est-leftsonic-case1} that Theorem \ref{theorem-BCF1} applies to $\psi$ near $\leftsonic$.
Then the admissible solution $\vphi$ satisfies statements (b)--(c) of Theorem \ref{theorem-3}.

\medskip
{\textbf{3.}} {\emph{Proof of statement {\rm (d)} in Theorem {\rm \ref{theorem-3}}}}.
By Lemma \ref{lemma-est-nrsonic-N}(d), $\shock\cap {\nD_{\bar{\eps}}}$ is represented
as the graph of $y=\fshockn(x)$ for $0\le x\le \bar{\eps}$,
where $\nD_{\bar{\eps}}$ is defined by \eqref{definition-Dr-ext}.

Let $\{y_m^{(1)}\}$ be a sequence satisfying $0<y_m^{(1)}<\fshockn(0)$ for each $m\in \mathbb{N}$,
and $\displaystyle{\lim_{m\to \infty}y_m^{(1)}=\fshockn(0)}$.
By \eqref{Op2}, \eqref{Op2-a}, and Theorem \ref{theorem-3}(c),
we can choose a sequence $\{x_m^{(1)}\}$ such that  $\{(x_m^{(1)}, y_m^{(1)})\}\subset \Om$,
$x_m^{(1)}\in(0, \frac 1m)$,
and
\begin{equation*}
\big|\psi_{xx}(x_m^{(1)}, y_{m}^{(1)})-\frac{1}{\gam+1}\big|<\frac 1m\qquad\tx{for each $m\in\mathbb{N}$}.
\end{equation*}
By Lemma \ref{lemma-est-nrsonic-N}(d), $0<y_m^{(1)}<\fshockn(0)<\fshockn(x_m^{(1)})$ for each $m\in\mathbb{N}$.
Therefore, we have
\begin{equation}
\label{limit-nr-sonic2018}
\lim_{m\to \infty}(x_m^{(1)}, y_m^{(1)})=(0, \fshockn(0)),\qquad
\lim_{m\to \infty}\psi_{xx}(x_m^{(1)}, y_m^{(1)})=\frac{1}{\gam+1}.
\end{equation}

By properties (a) and (c) of Lemma \ref{lemma-xybvp-N},
and Proposition \ref{lemma-est-sonic-general-N},
there exists $\eps\in(0,\bar{\eps}]$ such that, on $\shock\cap \nD_{\eps}$,
the boundary condition \eqref{bc-psi-B1-right} can be rewritten as
\begin{equation}
\label{bc-psi-shock-nr-rsonic2018}
\psi_x+b_1\psi_y+b_0\psi=0\qquad\tx{on $\shock \cap {\nD_{\eps}}$}
\end{equation}
for $(b_0,b_1)=(b_0,b_1)(\psi_x, \psi_y, \psi, x, \fshockn(x))$.
Let $\om>0$ be from Lemma \ref{lemma-est-nrsonic-N}(d). Then
$$
\{(x,\fshockn(x)-\frac{\om}{10}x)\,:\,0<x<\eps\}\subset \Om.
$$
Denote  $\mcl{F}(x):=\psi_x(x,\fshockn(x)-\frac{\om}{10}x)$.
By \eqref{bc-psi-shock-nr-rsonic2018}, we have
\begin{equation*}
\begin{split}
\mcl{F}(x)
&=\psi_x(x, \fshockn(x))-\frac{\om}{10}x\int_0^1 \psi_{xy}(x,\fshockn(x)-\frac{t\om}{10}x)\,\dd t\\
&=-(b_1\psi_y+b_0\psi)(x,\fshockn(x))-\frac{\om}{10}x\int_0^1 \psi_{xy}(x,\fshockn(x)-\frac{t\om}{10}x)\,\dd t
\qquad\tx{for $0<x<\eps$.}
\end{split}
\end{equation*}
From the last equality and Proposition \ref{lemma-est-sonic-general-N},
we obtain that $\mcl{F}(0)=0$, $\mcl{F}\in C([0,\eps])\cap C^1((0,\eps))$,
and $\lim_{x\to 0+}\frac{\mcl{F}(x)}{x}=0$.
Then, by the mean value theorem, there exists a sequence $\{x_m^{(2)}\}\subset (0,\eps)$ such that
\begin{equation}
\label{limit-of-F-shock2018}
\lim_{m\to \infty} x_m^{(2)}=0,\qquad \mcl{F}'(x_m^{(2)})=0.
\end{equation}
For each $m\in \mathbb{N}$, define $y_m^{(2)}:=\fshockn(x_m^{(2)})-\frac{\om}{10}x_m^{(2)}$
so that $\{(x_m^{(2)}, y_m^{(2)})\}\subset \Om$.
By the definition of $\mcl{F}$ and \eqref{limit-of-F-shock2018}, we have
\begin{equation}
\label{limit2-of-F-shock2018}
\begin{split}
\lim_{m\to \infty} \psi_{xx}(x_m^{(2)}, y_m^{(2)})
&=\lim_{m\to \infty}\mcl{F}'(x_m^{(2)})-\lim_{m\to \infty}(\fshockn'(x_m^{(2)})-\frac{\om}{10})\psi_{xy}(x_m^{(2)}, y_m^{(2)})\\
&=-\lim_{m\to \infty}(\fshockn'(x_m^{(2)})-\frac{\om}{10})\psi_{xy}(x_m^{(2)}, y_m^{(2)}).
\end{split}
\end{equation}
Since $\displaystyle{\lim_{m\to \infty}(x_m^{(2)}, y_m^{(2)})=(0,\fshockn(0))}$,
we combine \eqref{limit2-of-F-shock2018} with Proposition \ref{lemma-est-sonic-general-N} to obtain
\begin{equation}
\label{limit-near-shock-right}
\lim_{m\to \infty} \psi_{xx}(x_m^{(2)}, y_m^{(2)})=0.
\end{equation}
In \eqref{limit-nr-sonic2018} and \eqref{limit-near-shock-right}, we have shown that there are
two sequences, $\{(x_m^{(1)}, y_m^{(1)})\}$ and $\{(x_m^{(2)}, y_m^{(2)})\}$, in $\Om$ such that
the limits of both sequences are $(0, \fshockn(0))$. On the other hand,
\begin{equation*}
\lim_{m\to \infty}\psi_{xx}(x_m^{(1)}, y_m^{(2)})\neq \lim_{m\to \infty}\psi_{xx}(x_m^{(1)}, y_m^{(2)}).
\end{equation*}

For $\beta\in(0,\betasonic)$, we can repeat the argument above by using Lemma \ref{lemma-est-nrsonic}(d)
and Propositions \ref{lemma-est-sonic-general}  and \ref{proposition-sub8} to show that there are
two sequences, $\{(\til x_m^{(1)}, \til y_m^{(1)})\}$ and $\{(\til x_m^{(2)}, \til y_m^{(2)})\}$, in $\Om$
such that the limits of both sequences are $(0, \fshocko(0))$, but
it can similarly be shown that
\begin{equation*}
\lim_{m\to \infty}\psi_{xx}(\til x_m^{(1)}, \til y_m^{(1)})=\frac{1}{\gam+1}\neq 0 =\lim_{m\to \infty}\psi_{xx}(\til x_m^{(2)}, \til y_m^{(2)}),
\end{equation*}
where $\fshocko$ is from Lemma \ref{lemma-est-nrsonic}.
This proves statement (d) of Theorem \ref{theorem-3}.

\medskip
{\textbf{3.}} {\emph{Proof of statement {\rm (e)} of Theorem {\rm \ref{theorem-3}}}}.
By Lemma \ref{lemma1-sonic-N}(e), $\rightshock$ is represented as the graph
of $y=\hat{f}_{\mcl{N},0}(x)$ near point $\righttop$ in the $(x,y)$--coordinates given by \eqref{coord-n}.
We extend the definition of $\fshockn$ into $(-\bar{\eps}, \bar{\eps})$ by
\begin{equation}
\label{extension-fshockn}
\fshockn(x)=\hat{f}_{\mcl{N},0}(x)\qquad\,\,\tx{for}\,\,x\in(-\bar{\eps},0].
\end{equation}
By Proposition \ref{lemma-est-sonic-general-N}, $\fshockn$ satisfies
\begin{equation}
\label{fn-vanish}
(\fshockn-\hat{f}_{\mcl{N},0})(0)=(\fshockn-\hat{f}_{\mcl{N},0})'(0)=0,
\end{equation}
so that curve $\ol{\shock\cup \rightshockseg}$
is $C^{1,1}$, including at point $\righttop$.

Define
\begin{equation*}
\phi_{\infty}^{\mcl{N}}:=\ivphi-\rightvphi.
\end{equation*}
Since $\phi_{\infty}^{\mcl{N}}(x, \hat f_{\mcl{N},0}(x))=0$ and $(\varphi_\infty-\varphi)(x, \hfshock(x))=0$,
$\psi$ satisfies
\begin{equation}
\label{bc-fdiff-2018}
\phi_{\infty}^{\mcl{N}}(x,\hat f_{\mcl{N},0}(x))-\phi_{\infty}^{\mcl{N}}(x, \fshockn(x))=\psi(x, \hfshock(x))\qquad\tx{for $0<x<\bar{\eps}$}.
\end{equation}
A direct computation yields that
\begin{equation}
\label{fderiv-expression-2018}
\begin{split}
&\qquad\,\, \frac{\dd^2\phi_{\infty}^{\mcl{N}}(x,\hat f_{\mcl{N},0}(x))}{\dd x^2}
=\hat f_{\mcl{N},0}''(x)\der_y\phi_{\infty}^{\mcl{N}}(x,\hat f_{\mcl{N},0}(x))
   +\sum_{k=0}^2a_k(\hat f_{\mcl{N},0}'(x))^{k}\der_x^{2-k}\der_y^{k}\phi_{\infty}^{\mcl{N}}(x,\hat f_{\mcl{N},0}(x)),\\
&\qquad\,\, \frac{\dd^2\phi_{\infty}^{\mcl{N}}(x,\fshockn(x))}{\dd x^2}
=\fshockn''(x)\der_y\phi_{\infty}^{\mcl{N}}(x,\fshockn(x))
 +\sum_{k=0}^2a_k(\fshockn'(x))^{k}\der_x^{2-k}\der_y^{k}\phi_{\infty}^{\mcl{N}}(x,\fshockn(x)),
\end{split}
\end{equation}
\begin{equation}
\label{fderiv-expression-2018}
\begin{split}
& \frac{\dd^2\phi_{\infty}^{\mcl{N}}(x,\hat f_{\mcl{N},0}(x))}{\dd x^2}\\
&\,\,\,=\hat f_{\mcl{N},0}''(x)\der_y\phi_{\infty}^{\mcl{N}}(x,\hat f_{\mcl{N},0}(x))
   +\sum_{k=0}^2a_k(\hat f_{\mcl{N},0}'(x))^{k}\der_x^{2-k}\der_y^{k}\phi_{\infty}^{\mcl{N}}(x,\hat f_{\mcl{N},0}(x)),\\
&\frac{\dd^2\phi_{\infty}^{\mcl{N}}(x,\fshockn(x))}{\dd x^2}\\
&\,\,\,=\fshockn''(x)\der_y\phi_{\infty}^{\mcl{N}}(x,\fshockn(x))
 +\sum_{k=0}^2a_k(\fshockn'(x))^{k}\der_x^{2-k}\der_y^{k}\phi_{\infty}^{\mcl{N}}(x,\fshockn(x)),
\end{split}
\end{equation}
with $(a_0, a_1, a_2)=(1,2,1)$.

We differentiate \eqref{bc-fdiff-2018} with respect to $x$ twice and use \eqref{fderiv-expression-2018} to obtain
the following expression:
\begin{equation}
\label{f-2nd-deriv-right}
(\fshockn-\hat{f}_{\mcl{N},0})''(x)
=\frac{A_1(x)+A_2(x)+A_3(x)}{\der_y\phi_{\infty}^{\mcl{N}}(x,\fshockn(x))},
\end{equation}
where
\begin{equation*}
\begin{split}
&A_1(x)=\sum_{k=0}^2 a_k\left((\hat f_{\mcl{N},0}'(x))^k\der_x^{2-k}\der_y^k\phi_{\infty}^{\mcl{N}}(x,\hat{f}_{\mcl{N},0}(x))-
(\fshockn'(x))^k\der_x^{2-k}\der_y^k\phi_{\infty}^{\mcl{N}}(x,\fshockn(x))\right),\\
&A_2(x)=\left(\der_y\phi_{\infty}^{\mcl{N}}(x,\hat{f}_{\mcl{N},0}(x))-\der_y\phi_{\infty}^{\mcl{N}}(x,\fshockn(x))\right)\hat{f}_{\mcl{N},0}''(x),\\
&A_3(x)=-\Big(\fshockn''(x)\psi_y(x,\fshockn(x))+\sum_{k=0}^2a_k(\fshockn'(x))^{k}\der_x^{2-k}\der_y^{k}\psi(x,\fshockn(x))\Big).
\end{split}
\end{equation*}
By \eqref{fn-vanish}, we have
\begin{equation}
\label{A-estimate1-2018}
A_1(0)=A_2(0)=0.
\end{equation}
We differentiate the boundary condition \eqref{bc-psi-B1-right} in the tangential direction along $\shock$,
and apply Lemma \ref{lemma-xybvp-N}(a)--(c) and Proposition \ref{lemma-est-sonic-general-N} to obtain that
there exists a constant $C>0$ such that
\begin{equation*}
\begin{split}
&|\psi_{xx}(x,\fshockn(x))|\\
&\le C
\big(|\psi(x,\fshockn(x))|+|D_{(x,y)}\psi(x,\fshockn(x))|+|D_{(x,y)}\psi_y(x,\fshockn(x))|\big)
\qquad\tx{on $\shock\cap \nD_{\bar{\eps}}$}.
\end{split}
\end{equation*}
From this estimate and Proposition \ref{lemma-est-sonic-general-N}, we see that
$\displaystyle{\lim_{x\to 0+}\psi_{xx}(x,\fshockn(x))=0}$, which implies that
\begin{equation}
\label{A-estimate2-2018}
\lim_{x\to 0+}A_3(x)=0.
\end{equation}
By Lemma \ref{lemma1-sonic-N}(c), $\der_y\phi_{\infty}^{\mcl{N}}(x,\fshockn(x))\neq 0$ on $\shock\cap \nD_{\bar{\eps}}$.
Then we conclude from \eqref{f-2nd-deriv-right}--\eqref{A-estimate2-2018}
that
\begin{equation*}
(\fshockn-\hat{f}_{\mcl{N},0})''(0)=0.
\end{equation*}
This implies that the extension of $\fshockn$ given by \eqref{extension-fshockn}
is in $C^2([-\bar{\eps}, \bar{\eps}])$.
Furthermore, we conclude from \eqref{f-2nd-deriv-right} and Proposition \ref{lemma-est-sonic-general-N}
that the extension of $\fshockn$ given by \eqref{extension-fshockn}
is in $C^{2,\alp}((-\bar{\eps}, \bar{\eps}))$ for any $\alp\in(0, 1)$.
This implies that $\ol{\shock\cup\rightshockseg}$ is $C^{2,\alp}$ for any $\alp\in(0, 1)$,
including at point $\righttop=(0, \fshockn(0))$.
For $\beta\in(0, \betasonic)$, it can similarly be checked
that $\ol{\leftshockseg\cup\shock}$ is $C^{2,\alp}$ for any $\alp\in(0, 1)$, including
at point $\lefttop=(0, \fshocko(0))$ for $\fshocko$ from Lemma \ref{lemma-est-nrsonic}.
Therefore,  statement (e) of Theorem \ref{theorem-3} is proved.

\appendix

\chapter{The Shock Polar for Steady Potential Flow}
\label{section-appendix}
\numberwithin{equation}{chapter}

According to \cite{CF}, for any given uniform supersonic state,
a shock polar curve for the two-dimensional steady full Euler system
exist and is convex.
In this appendix, we show the same for the potential flow.
The convexity of the shock polar curve leads to Lemma \ref{lemma2-appendix},
which is the key ingredient for proving
the existence of admissible solutions in the sense of
Definition \ref{def-regular-sol-op} for $(\iu, u_0)\in \mathfrak{P}_{\rm weak}$
with $u_0\le u_{\rm s}^{(\irho, \iu)}$, and
the non-existence of admissible solutions for $(\iu, u_0)\in \mathfrak{P}_{\rm strong}$.
The existence of convex shock polar curves for potential flow
is proved by combining the results from \cite{EL3, K}.

The two-dimensional steady potential flow for an ideal polytropic gas
is governed by the equations:
\begin{equation*}
\begin{cases}
(\rho u)_{x_1}+(\rho v)_{x_2}=0,\\
u_{x_2}-v_{x_1}=0,\\[1mm]
\frac 12(u^2+v^2)+i(\rho)=B_0\qquad\tx{\rm{(Bernoulli's law)}}
\end{cases}
\end{equation*}
for a constant $B_0>0$, where $i(\rho)$ is given by
\begin{equation*}
h(\rho)=\begin{cases}
\frac{\rho^{\gam-1}-1}{\gam-1}&\,\,\tx{for}\;\;\gam>1,\\
\ln\rho&\,\,\tx{for}\;\;\gam=1.
\end{cases}
\end{equation*}

\begin{lemma}
\label{lemma-appendix}
Fix $\gam\ge 1$ and the incoming constant state $(\irho, {\bf u}_{\infty})=(\irho, (\iu,0))$
with $\iu>\irho^{(\gam-1)/2}>0$.
Denote $M_{\infty}:=\frac{\iu}{\irho^{(\gam-1)/2}}>1$ as the Mach number of the incoming supersonic flow.
For each $\beta\in[0,\cos^{-1}(\frac{1}{M_{\infty}}))$,
there exists a unique ${\bf u}=(\leftu, \leftv)\in(\R_+)^2\setminus \{{\bf u}_{\infty}\}$ such that
\begin{align}
\label{ap-1}
&\leftrho {\bf u}\cdot {\bf n}=\irho {\bf u}_{\infty}\cdot{\bf n},\\
\label{ap-2}
&({\bf u}_{\infty}-{\bf u})\cdot {\bf t}=0,\\
\label{ap-8}
&\frac 12({\bf u}\cdot {\bf n})^2+i(\leftrho)=\frac 12({\bf u}_{\infty}\cdot{\bf n})^2+i(\irho)
\end{align}
for  ${\bf n}=(\cos\beta,-\sin\beta)$ and ${\bf t}=(\sin\beta,\cos\beta)$,
where $\leftrho$ is given by
\begin{equation}
\label{app-bernoulli}
\leftrho=i^{-1}(i(\irho)+\frac 12(\iu^2-|{\bf u}|^2)).
\end{equation}
In other words, ${\bf u}$ becomes the downstream velocity behind a straight oblique shock $\leftshock$ of
angle $\frac{\pi}{2}-\beta$ from the horizontal axis.
Moreover, the collection of such ${\bf u}=(\leftu, \leftv)$ for $\beta\in[0,\cos^{-1}(\frac{1}{M_{\infty}}))$
forms
a concave curve
on the $(u,v)$--plane.

\begin{proof}
The existence of the curve for $(\leftu, \leftv)$ is verified by following the proof
of \cite[Proposition 2.1]{K}, and the convexity of this curve can be checked by adjusting
the proof of \cite[Theorem 1]{EL3}.
We prove the lemma for the case that $\gam>1$.
The case that $\gam=1$ can be treated in the same way. The proof is divided into two steps.

\smallskip
{{\bf 1}.} {\emph{Existence of shock polar}}.
Fix constants $\gam>1$, $\irho>0$, and $\iu$ with $\iu>\irho^{(\gam-1)/2}$.
Let $\leftshock$ be a straight oblique shock with angle $\frac{\pi}{2}-\beta$ from the horizontal ground,
and let $\leftrho$ and ${\bf u}=(\leftu,\leftv)$ be the density and the velocity behind shock $\leftshock$.
By \eqref{ap-2}, the angle between vector ${\bf u}-{\bf u}_{\infty}$ and
the horizontal axis in Fig. \ref{fig:steadypolar} is $\beta$.
\begin{figure}[htp]
\centering
\begin{psfrags}
\psfrag{u}[cc][][0.8][0]{$u$}
\psfrag{v}[cc][][0.8][0]{$v$}
\psfrag{iu}[cc][][0.8][0]{$\iu$}
\psfrag{p}[cc][][0.8][0]{$\phantom{aaa}(\leftu,\leftv)$}
\psfrag{b}[cc][][0.8][0]{$\beta$}
\psfrag{q}[cc][][0.8][0]{$\leftq$}
\includegraphics[scale=.8]{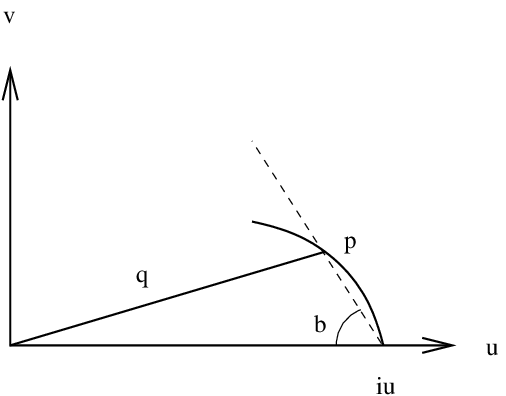}
\caption{The shock polar for potential flow}
\label{fig:steadypolar}
\end{psfrags}
\end{figure}
By the expression of $\{{\bf n}, {\bf t}\}$, we have
\begin{equation}
\label{ap-3}
\begin{split}
&{\bf u}_{\infty}\cdot{\bf n}=\iu\cos\beta,\qquad \qquad\qquad{\bf u}_{\infty}\cdot{\bf t}=\iu\sin\beta,\\
&{\bf u}\cdot{\bf n}=\leftu\cos\beta-\leftv\sin\beta,\qquad \;\;{\bf u}\cdot {\bf t}=\leftu\sin\beta+\leftv\cos\beta.
\end{split}
\end{equation}
Denote $M_{\infty,n}=\frac{{\bf u}_{\infty}\cdot{\bf n}}{\irho^{(\gam-1)/2}}$.
For each $\beta\in [0, \frac{\pi}{2})$, $M_{\infty, n}$ is fixed and $M_{\infty, n}>0$ holds.
It has been shown in the proof of Lemma \ref{lemma-2-0} that there exists a unique $M_n$
with $M_n\neq M_{\infty,n}$ as a solution of the equation:
\begin{equation}
\label{ap-6}
{\mathfrak g}(M_n)={\mathfrak g}(M_{\infty,n})
\end{equation}
for
${\mathfrak g}(M)=(M^2+\frac{2}{\gam-1})M^{-\frac{2(\gam-1)}{\gam+1}}$,
unless $M_{\infty,n}=1$.
Substitute ${\bf u}\cdot{\bf n}=M_n \leftrho^{\frac{\gam-1}{2}}$ into \eqref{ap-8}
and solve the resultant equation for $\leftrho$ to obtain
\begin{equation*}
\leftrho^{\gam-1}=\frac{({\bf u}_{\infty}\cdot {\bf n})^2+2i(\irho)+\frac{2}{\gam-1}}{M_n^2+\frac{2}{\gam-1}}.
\end{equation*}
By the entropy condition,
shock $\leftshock$ is admissible only if $\irho<\leftrho$, which
is equivalent to saying that $0<M_n<1<M_{\infty,n}$.
Since $M_{\infty,n}= M_{\infty}\cos\beta$ for
$M_{\infty}=\frac{\iu}{\irho^{(\gam-1)/2}}$,
we restrict our consideration only to
the case that $\beta\in[0,\cos^{-1}(\frac{1}{M_{\infty}}))$.
Then \eqref{ap-2} and \eqref{ap-3} yield that
\begin{equation*}
  \begin{pmatrix}
    \leftu \\
    \leftv
  \end{pmatrix}
  =\begin{pmatrix}
     \cos\beta & \sin\beta \\
     -\sin\beta & \cos\beta
   \end{pmatrix}
   \begin{pmatrix}
     M_n\leftrho^{\frac{\gam-1}{2}} \\
     \iu\sin\beta
   \end{pmatrix}.
\end{equation*}
Therefore, curve $(\leftu, \leftv)(\beta)$ is given
for $\beta\in[0,\cos^{-1}(\frac{1}{M_{\infty}}))$ in the $(u,v)$--plane;
see Fig. \ref{fig:steadypolar}.

Since $\lim_{\beta\to \cos^{-1}(\frac{1}{M_{\infty}})}M_n=1=\lim_{\beta\to \cos^{-1}(\frac{1}{M_{\infty}})}M_{\infty,n}$,
the shock polar curve is extended up to $\beta=\cos^{-1}(\frac{1}{M_{\infty}})$ by $(\leftu, \leftv)=(\iu,0)$.

This
curve $(u,v)=(\leftu, \leftv)(\beta)$ for
$\beta\in[0,\cos^{-1}(\frac{1}{M_{\infty}})]$
is called a \emph{shock polar for potential flow}.

\smallskip
{\bf 2.} {\emph{Convexity of shock polar}}.
Let ${\bf u}=(u,v)$ denote each point on the shock polar curve.
By \eqref{ap-1}--\eqref{ap-2},
each point ${\bf u}$ on the shock polar satisfies the equation:
\begin{equation}
\label{g}
g({\bf u})=\big(\rho(|{\bf u}|^2){\bf u}-\irho{\bf u}_{\infty}\big)\cdot \frac{{\bf u}_{\infty}-{\bf u}}{|{\bf u}_{\infty}-{\bf u}|}=0
\end{equation}
for ${\bf u}_{\infty}=(\iu,0)$, where $\rho(|{\bf u}|^2)$ is given by \eqref{app-bernoulli}
so that
$D_{\bf u}\rho=-\frac{{\bf u}}{c^2}\rho$ for $c^2(|{\bf u}|^2)=\rho^{\gam-1}(|{\bf u}|^2)$.
Combining this with \eqref{g} gives that
\begin{equation}
\label{ap-4}
\begin{split}
&g_{\bf u}\cdot{\bf n}=\rho\left(1-\big(\frac{{\bf u}\cdot{\bf n}}{c}\big)^2\right),\qquad
g_{\bf u} \cdot{\bf t}=-({\bf u}_{\infty}\cdot{\bf t})\left(\frac{\rho {\bf u}\cdot{\bf n}}{c^2}
+\frac{\rho-\irho}{|{\bf u}_{\infty}-{\bf u}|}\right).
\end{split}
\end{equation}
By the entropy condition, we have
\begin{equation}
\label{entropy-condition-gn}
g_{\bf u}\cdot{\bf n}>0.
\end{equation}
Define
\begin{equation*}
{\bf q}:=\frac{g_{\bf u}}{g_{\bf u}\cdot{\bf n}},
\end{equation*}
and express ${\bf q}$ as
${\bf q}
={\bf n}+\frac{g_{\bf u}\cdot{\bf t}}{g_{\bf u}\cdot{\bf n}}{\bf t}.
$

\medskip
\emph{Claim}:  $\, {\bf q}\times \frac{\dd{\bf q}}{\dd\beta}<0$ for all $\beta\in(0,\cos^{-1}(\frac{1}{M_{\infty}}))$.

\medskip
Denote $A:=-\frac{g_{\bf u}\cdot{\bf t}}{g_{\bf u}\cdot{\bf n}}$. Then
$
\frac{\dd{\bf q}}{\dd\beta}=-(1+\frac{\dd A}{\dd\beta}){\bf t}-A {\bf n},
$
which implies that
\begin{equation}
\label{A-eqn}
{\bf q}\times\frac{\dd{\bf q}}{\dd\beta}=-\big(1+A^2+\frac{\dd A}{\dd\beta}\big).
\end{equation}
By \eqref{ap-1}, \eqref{ap-3}, and \eqref{ap-4}, we can rewrite $A$ as
$A=\frac{\iu \sin\beta}{1-M_n^2}(\frac{M_n}{c}+\frac{1}{\iu\cos\beta})$ for $M_n:=\frac{{{\bf u}\cdot{\bf n}}}{c}.$
Differentiate \eqref{ap-6} with respect to $\beta$ to obtain
\begin{equation*}
\frac{\dd M_n}{\dd\beta}
=-\frac{{\mathfrak g}'(M_{\infty,n})}{{\mathfrak g}'(M_n)}
\frac{\iu\sin\beta}{\irho^{(\gam-1)/2}}>0\qquad\tx{for $\beta\in (0,\cos^{-1}(\frac{1}{M_{\infty}}))$}.
\end{equation*}
From   $\rho^{\frac{\gam+1}{2}}M_n=\irho^{\frac{\gam+1}{2}}M_{\infty,n}=\irho\iu\cos\beta$
and $\frac{\dd M_n}{\dd\beta}>0$,
we see that $\frac{\dd\rho}{\dd\beta}<0$
so that
$\frac{\dd A}{\dd\beta}\ge 0$
holds for all $\beta\in(0,\cos^{-1}(\frac{1}{M_{\infty}}))$.
Combining this with \eqref{A-eqn}, we have
$$
{\bf q}\times \frac{\dd{\bf q}}{\dd\beta}\le -1\qquad \mbox{for $\beta\in(0,\cos^{-1}(\frac{1}{M_{\infty}}))$}.
$$
The claim is verified.

\medskip
The inequality above gives the useful property:
\begin{equation}
\label{ap-7}
\frac{\bf q}{|\bf q|}\times \frac{\dd}{\dd\beta}\Big(\frac{\bf q}{|\bf q|}\Big)
=\frac{{\bf q}\times \frac{\dd{\bf q}}{\dd\beta}}{|{\bf q}|^2}\le -\frac{1}{|{\bf q}|^2}< 0
\end{equation}
at each point on the shock polar curve.

Fix a point ${\bf u}_0=(u_0,v_0)$ on the shock polar $\{{\bf u}=(u,v)\,:\,g({\bf u})=0\}$,
and define ${\bf n}_0=\frac{{\bf u}_0-{\bf u}_{\infty}}{|{\bf u}_0-{\bf u}_{\infty}|}$.
We introduce a new coordinate system $(s,t)$ so that the following properties hold
in the new $(s,t)$--coordinates:

\smallskip
\begin{itemize}
\item[(i)] ${\bf u}_0=(0,0)$, ${\bf n}_0=(0,1)$;

\smallskip
\item[(ii)] If ${\bm \tau}_0$ is the unit vector perpendicular to ${\bf n}_0$ and oriented
to satisfy ${\bf u}_{\infty}\cdot {\bm \tau}_0>0$, then ${\bm \tau}_0=(1,0)$.
\end{itemize}

\smallskip
\noindent
Define a function $\mcl{G}(s,t)$ by
\begin{equation*}
\mcl{G}(s^{(\bf u)},t^{(\bf u)})=g(\bf u),
\end{equation*}
where $(s^{(\bf u)}, t^{(\bf u)})$ is the $(s,t)$--coordinates of ${\bf u}$ on the shock polar.
Since the value of $g_{\bf u}\cdot {\bf n}$ for ${\bf n}=\frac{{\bf u}_{\infty}-{\bf u}}{|{\bf u}_{\infty}-{\bf u}|}$
is invariant under the rotation,
$
\mcl{G}_t(0,0)=-(g_{\bf u}\cdot {\bf n})({\bf u}_0)<0.
$
By the implicit function theorem,
there exists a function $f_{{\bf u}_0}:(-\eps_0,\eps_0)\to \R$ for some small constant $\eps_0>0$ so that
the shock polar curve is represented by $t=f_{{\bf u}_0}(s)$ near ${{\bf u}_0}$ in the $(s,t)$--coordinates.
Such a function $f_{{\bf u}_0}$ satisfies the relation:
\begin{equation*}
\frac{f''_{{\bf u}_0}(0)}{\sqrt{1+(f'_{{\bf u}_0}(0))^2}}
=\frac{\bf q}{|\bf q|}\times \frac{\dd}{\dd\beta}\Big(\frac{\bf q}{|\bf q|}\Big)\Big|_{{\bf u}={\bf u}_0}
\le -\frac{1}{|{\bf q}({\bf u}_0)|^2}<0.
\end{equation*}
Therefore, we conclude that the shock polar for potential flow is concave.
\end{proof}
\end{lemma}

\begin{remark}
\label{remark-appendix-de}
Fix $\gam\ge 1$ and $(\irho, \iu)$ with $\iu>\irho^{(\gam-1)/2}>0$.
Let $\Upsilon^{(\irho,\iu)}$ be the shock polar curve lying in the first quadrant
in the $(u,v)$--plane for the steady potential flow with
the incoming supersonic state $(\irho, \iu)$.
Owing to the concavity of the shock polar,
there exists a unique $\theta_{\rm d}^{(\irho, \iu)}\in(0, \frac{\pi}{2})$
such that the following properties hold{\rm :}

\smallskip
\begin{itemize}
\item[(i)] If $0\le \theta_{\rm w}<\theta_{\rm d}^{(\irho, \iu)}$,
then line $\frac vu=\tan\theta_{\rm w}$ intersects with $\Upsilon^{(\irho, \iu)}$
at two distinct points{\rm ;}

\smallskip
\item[(ii)] Line $\frac vu=\tan\theta_{\rm d}^{(\irho, \iu)}$
and $\Upsilon^{(\irho, \iu)}$ have a unique intersection point so that
$\frac vu=\tan\theta_{\rm d}^{(\irho, \iu)}$ is tangential
to $\Upsilon^{(\irho, \iu)}$ at the intersection point{\rm ;}

\smallskip
\item[(iii)] If $\theta_{\rm d}^{(\irho, \iu)}<\theta_{\rm w}<\frac{\pi}{2}$,
then line $\frac vu=\tan\theta_{\rm w}$ never
intersects
with $\Upsilon^{(\irho, \iu)}$.
\end{itemize}
\end{remark}

\begin{lemma}
\label{lemma-app2}
Fix $\gam\ge 1$.
For each $(\irho, \iu)$ with $\iu>\irho^{(\gam-1)/2}>0$,
there exist a unique constant $\hat u_0^{(\irho, \iu)}=:\hat u_0 \in(0,\iu)$ and a unique smooth function
$f_{\rm{polar}}\in C^0([\hat u_0, \iu])\cap C^{\infty}((\hat u_0, \iu))$ such that
\begin{equation}
\label{app-5}
\Upsilon^{(\irho,\iu)}=
\{(u, f_{\rm{polar}}(u))\,:\, u\in[\hat u_0,\iu]\}.
\end{equation}
Furthermore, the following properties hold{\rm :}

\smallskip
\begin{itemize}
\item[(a)] Let $\theta_{\rm s}^{(\irho, \iu)}$ be from Lemma {\rm \ref{lemma-spolar-steady-new}(c)}.
Then there exist unique $u_{\rm d}, u_{\rm s}\in(\hat u_0, \iu)$ such that
\begin{equation}
\label{definition-ud-us}
\frac{f_{\rm{polar}}(u_{\rm s})}{u_{\rm s}}=\tan \theta_{\rm s}^{(\irho, \iu)},
\qquad \frac{f_{\rm{polar}}(u_{\rm d})}{u_{\rm d}}=\tan \theta_{\rm d}^{(\irho, \iu)}.
\end{equation}
Moreover, $u_{\rm d}<u_{\rm s}$ holds, and $(u_{\rm d}, u_{\rm s})$ vary
continuously on $(\irho, \iu)$.

\smallskip
\item[(b)] Denote by $f_{\rm{polar}}(\cdot, \irho, \iu)$ the shock polar function $f_{\rm{polar}}(\cdot)$
for the incoming flow $(\irho, \iu)$.
Then $f_{\rm{polar}}$ as a function of $(u, \irho, \iu)$ is $C^{\infty}$ on the domain{\rm :}
\begin{equation*}
\{(u, \irho, \iu)\,:\,\irho>0, \,\,\iu>\irho^{(\gam-1)/2},\,\,u\in(\hat u^{(\irho, \iu)}, \iu)\}.
\end{equation*}
\end{itemize}

\begin{proof}
The proof is divided into four steps.

\smallskip
{\textbf{1.}} For each $\beta\in[0, \cos^{-1}(\frac{1}{M_{\infty}})]$,
let $(\leftrho, \leftu, \leftv)$ be from Lemma \ref{lemma-appendix}, and let $\leftq:=\sqrt{\leftu^2+\leftv^2}$.
Since $(\leftrho, \leftu, \leftv)$ is uniquely determined for $\beta\in [0, \cos^{-1}(\frac{1}{M_{\infty}})]$,
$\leftq$ is considered as a function of $\beta$.
Substituting \eqref{ap-3} into \eqref{ap-1}--\eqref{ap-2}, we obtain
$$
(\leftu,\leftv)=\iu (1-(1-\frac{\irho}{\leftrho})\cos^2\beta,\,(1-\frac{\irho}{\leftrho})\cos\beta\sin\beta),
$$
so that
\begin{equation}
\label{q-eqn}
\cos^2\beta=\frac{1-(\frac{\leftq}{\iu})^2}{1-(\frac{\irho}{\leftrho})^2}=:h(\leftq).
\end{equation}
It follows from \eqref{app-bernoulli} and \eqref{q-eqn} that
\begin{equation*}
\begin{split}
{ h}'(\leftq)=:\frac{2\leftq}{(1-\frac{\irho^2}{\leftrho^2})^2\leftrho^2\leftc^2\iu^2}I(\leftq)
\end{split}
\end{equation*}
for $I(\leftq)$ satisfying $I(\iu)=0$ and $I'(\leftq)=(\gam+1)\leftq(\leftrho^2-\irho^2)$.
Inequality $\leftrho>\irho$ holds, owing to the entropy condition for the admissible shock
so that $I'(\leftq)>0$ and $I(\leftq)<I(\iu)=0$ for $0<\leftq<\iu$,
which implies that $h'(\leftq)<0$ for $0<\leftq<\iu$. Then \eqref{q-eqn} yields that
\begin{equation}
\label{leftq-monotinicity-new}
\frac{\dd\leftq}{\dd\beta}=-\frac{2\cos\beta\sin\beta}{{h'}(\leftq)}>0\qquad\,
\tx{for all}\,\, \beta\in(0,\cos^{-1}\frac{1}{M_{\infty}}).
\end{equation}

\smallskip
{\textbf{2.}}
Let $g({\bf u}), {\bf n}$, and ${\bf t}$ be given by \eqref{g}. Then \eqref{ap-4} implies
\begin{equation*}
\der_{v}g({\bf u})=-(g_{\bf u}\cdot{\bf n})\sin\beta+(g_{\bf u}\cdot{\bf t})\cos\beta<0
\end{equation*}
for any interior point ${\bf u}=(u,v)$ in $\Upsilon^{(\irho,\iu)}$.
By the implicit function theorem, there exists a unique function $f_{\rm{polar}}:[\hat u_0,\iu]\to [0,\infty)$
so that \eqref{app-5} holds, where $\hat u_0=\leftq(\beta)|_{\beta=0}$ for $\leftq$ defined through \eqref{q-eqn}.
The smoothness of map $(u, \irho, \iu)\mapsto f_{\rm{polar}}(u ,\irho, \iu)$ follows from the implicit function
theorem and the smooth dependence of $g({\bf u})$ on $(\irho, \iu)$.

\smallskip
{\textbf{3.}} The existence and uniqueness of $u_{\rm d}\in(\hat u_0, \iu)$ result directly
from the concavity of the shock polar curve $\Upsilon^{(\irho, \iu)}$.
Since point $(\hat u_0, 0)$ on the shock polar $\Upsilon^{(\irho, \iu)}$ corresponds to a normal shock,
$(\hat u_0, 0)$ is subsonic; that is, $\leftrho^{\gam-1}-\leftq^2>0$ holds at $\beta=0$.
At $\beta=\cos^{-1}(\frac{1}{M_{\infty}})$,
$\leftrho^{\gam-1}-\leftq^2<0$, because $(\leftrho, \leftq)|_{\beta=\cos^{-1}(\frac{1}{M_{\infty}})}=(\irho, \iu)$.
From \eqref{leftq-monotinicity-new} and Bernoulli's law that $\frac 12 \leftq^2+\irho(\leftrho)=B_0$,  we have
\begin{equation*}
\frac{\dd (\leftrho^{\gam-1}-\leftq^2)}{\dd\beta}<0\qquad\tx{for all}\,\,\beta\in(0, \cos^{-1}(\frac{1}{M_{\infty}})).
\end{equation*}
Therefore, there exists a unique $u_{\rm s}\in(\hat u_0, \iu)$ such that
$\frac{f_{\rm{polar}}(u_{\rm s})}{u_{\rm s}}=\tan \theta_{\rm s}^{(\irho, \iu)}$ holds.
Furthermore, Lemma \ref{lemma-spolar-steady-new}(c) and the concavity of $\Upsilon^{(\irho, \leftu)}$
imply that $u_{\rm d}<u_{\rm s}$.

\smallskip
{\textbf{4.}} By Bernoulli's law and the concavity of $\Upsilon^{(\irho, \iu)}$,
\eqref{definition-ud-us} is equivalent to
\begin{equation}
\label{equiv-ud-us}
\begin{split}
&u_{\rm s}^2+f_{\rm{polar}}^2(u_{\rm s}, \irho, \iu)=\frac{2(\gam-1)}{\gam+1}\Big(\frac 12 \iu^2+\frac{\irho^{\gam-1}}{\gam-1}\Big),\\
&f_{\rm{polar}}(u_{\rm d}, \irho, \iu)-u_{\rm d}f'_{\rm{polar}}(u_{\rm d}, \irho, \iu)=0
\end{split}
\end{equation}
for each $(\irho, \iu)$ with $\iu>\irho^{(\gam-1)/2}>0$.

For each $k\in\mathbb{N}$, let a sequence
$\{(\irho^{(k)}, \iu^{(k)})\}$ satisfy $\iu^{(k)}>(\irho^{(k)})^{(\gam-1)/2}>0$.
Also, suppose that $\{(\irho^{(k)}, \iu^{(k)})\}$ converges to $(\irho^*, \iu^*)$
with $\iu^*>(\irho^*)^{(\gam-1)/2}>0$.
Let $(u_{\rm d}^{(k)}, u_{\rm s}^{(k)})$ and $(u_{\rm d}^*, u_{\rm s}^*)$ be
the values of $(u_{\rm d}, u_{\rm s})$ corresponding
to $(\irho^{(k)}, \iu^{(k)})$ and $(\irho^*, \iu^*)$, respectively.
Note that
$$
(u_{\rm d}^{(k)}, \irho^{(k)}, \iu^{(k)}), (u_{\rm s}^{(k)}, \irho^{(k)}, \iu^{(k)})
\in\{(u, \irho, \iu)\,:\, \irho>0, \,\iu>\irho^{(\gam-1)/2},\,u\in(\hat u^{(\irho, \iu)}, \iu)\}
$$
for each $k\in\mathbb{N}$
and that $\hat u_0$ varies continuously on $(\irho, \iu)$ so that
$\{(u_{\rm d}^{(k)},u_{\rm s}^{(k)})\}$ is bounded in $(\R^+)^2$.
Therefore, there exist a convergent subsequence $\{(u_{\rm d}^{(k_j)},u_{\rm s}^{(k_j)})\}$
and states $(u_{\rm d}^{\sharp}, u_{\rm s}^{\sharp})$ such that
$\lim_{j\to \infty}(u_{\rm d}^{(k_j)},u_{\rm s}^{(k_j)})=(u_{\rm d}^{\sharp}, u_{\rm s}^{\sharp})$.
Then assertion (b) (proved in Step 2) and \eqref{equiv-ud-us} yield
\begin{equation*}
\begin{split}
&(u_{\rm s}^{\sharp})^2+f_{\rm{polar}}^2(u_{\rm s}^{\sharp}, \irho^*, \iu^*)
  =\frac{2(\gam-1)}{\gam+1}\left(\frac 12 (\iu^*)^2+\frac{(\irho^*)^{\gam-1}}{\gam-1}\right),\\
&f_{\rm{polar}}(u_{\rm d}^{\sharp}, \irho^*, \iu^*)
  -u_{\rm d}^{\sharp}f'_{\rm{polar}}(u_{\rm d}^{\sharp}, \irho^*, \iu^*)=0.
\end{split}
\end{equation*}
This implies that $(u_{\rm d}^{\sharp},u_{\rm s}^{\sharp})=(u_{\rm d}^*, u_{\rm s}^*)$,
since it has been shown in Step 3 that $(u_{\rm d}, u_{\rm s})$ satisfying \eqref{definition-ud-us}
for $(\irho^*, \iu^*)$ uniquely exists.
Therefore, we conclude that $(u_{\rm d}, u_{\rm s})$ varies continuously on $(\irho, \iu)$.
\end{proof}
\end{lemma}

In Lemma \ref{lemma-parameters}, the one-to-one correspondence between parameter
sets $\mathfrak{P}$ and $\mathfrak{R}$ is established.
For each $(\iu, u_0)\in \mathfrak{P}$, there exists a unique $\theta_{\rm w}\in(0, \frac{\pi}{2})$
such that $\iv$ is given by \eqref{2-4-b2}, where $(\iv, \beta)\in\mathfrak{R}$ corresponds to $(\iu, u_0)$.
The convexity of the shock polar obtained
in Lemma \ref{lemma-app2} yields the following property:

\begin{lemma}
\label{lemma2-appendix}
Fix $\gam\ge 1$ and $\iv>0$. For each $\beta\in (0, \frac{\pi}{2})$, let $\ivphi, \leftvphi, \leftrho$, and $P_{\beta}$
be defined by \eqref{2-4-b6}, \eqref{2-4-a0}, \eqref{2-4-a3}, \eqref{def-Pbeta}, respectively.
Denote $G({\bf p}, z, \bmxi)=g^{\rm sh}({\bf p}, z, \bmxi)$ for $g^{\rm sh}({\bf p}, z, \bmxi)$ defined by \eqref{def-gsh}.
Then there exists $\beta_{\rm d}^{(\iv)}\in(0,\frac{\pi}{2})$ depending only on $(\iv, \gam)$ such that
$G({\bf p}, z, \bmxi)$ satisfies
\begin{equation}
\label{11-a4}
G_{p_1}(D\leftvphi, \ivphi,P_{\beta})
\begin{cases}
<0& \quad\tx{for $\beta\in(0,\beta_{\rm d}^{(\iv)})$},\\[1mm]
=0& \quad \tx{for $\beta=\beta_{\rm d}^{(\iv)}$},\\[1mm]
>0& \quad \tx{for $\beta\in (\beta_{\rm d}^{(\iv)},\frac{\pi}{2})$}.
\end{cases}
\end{equation}

\begin{proof}
The following facts are useful to compute $G_{q_1}(D\leftvphi,\ivphi,P_{\beta})$:

\smallskip
\begin{itemize}
\item[(i)] The unit normal vector ${\bm n}_{\mcl{O}}$ to $\leftshock$ towards the downstream is
$
{\bm n}_{\mcl{O}}=\frac{D\ivphi-D\leftvphi}{|D\ivphi-D\leftvphi|}=(\sin\beta,-\cos\beta)
$
so that
$
(\leftrho D\leftvphi-D\ivphi)\cdot (1,0)=(\leftrho-1)(\leftu-\xi_1)\cos^2\beta
$,
where $D\leftvphi$ and $D\ivphi$ are evaluated at $\bmxi=(\xin,\etan)\in\R^2$.

\smallskip
\item[(ii)] It is shown from a direct computation that, if $G({\bf p}, z, \bmxi)=0$, then
\begin{equation}
\label{app-1}
G_{{\bf p}}({\bf p},z,\bmxi)=\frac{1}{\rho^{\gam-2}}
\left(c^2\frac{D\ivphi-{\bf p}}{|D\ivphi-{\bf p}|}
-\Big({\bf p}\cdot\frac{D\ivphi-{\bf p}}{|D\ivphi-{\bf p}|}\Big){\bf p}\right)
-\frac{\rho{\bf p}-D\ivphi}{|D\ivphi-{\bf p}|}
\end{equation}
for $\rho=\rho({\bf p},z)$.
\smallskip
\end{itemize}
It follows from (i)--(ii) that
\begin{equation}
\label{app-4}
\begin{split}
G_{p_1}(D\leftvphi,\ivphi,P_{\beta})
&=\left(\leftc^2-(\leftu-\xib)^2\right)\frac{\sin\beta}{\leftrho^{\gam-2}}
-\frac{(\leftrho-1)(\leftu-\xib)\cos^2\beta}{\sqrt{\leftu^2+\iv^2}}
\end{split}
\end{equation}
for $\leftc=\leftrho^{(\gam-1)/2}$.
Denote $\leftq:=D\leftvphi(P_{\beta})\cdot {\bm n}_{\mcl{O}}$.
Then
$\leftu-\xib=\leftq\csc\beta$,
where $P_{\beta}$ is denoted as $P_{\beta}=(\xi_1^{P_{\beta}},0)$.
Also, $\xi_2^m$ in the proof of Lemma \ref{lemma:interval-existence}
can be written as $\xi_2^m=\leftq\cos\beta$.
Substituting these two expressions into \eqref{app-4} and using the relations that $\leftu=-\iv\tan\beta$
and $\frac{(\leftrho-1)\leftq}{\iv\sec\beta}=1$ obtained from \eqref{2-4-b6}, \eqref{12-25},
and \eqref{2-4-a8}, we have
 \begin{equation*}
G_{p_1}(D\leftvphi,\ivphi,P_{\beta})
=\leftrho(1-\oM^2)\sin\beta-\frac{(\xi_2^m)^2}{\leftrho^{\gam-2}}\csc\beta
-\frac{\cos\beta}{\tan\beta},
\end{equation*}
where $\oM$ is defined by \eqref{1-25} with $c=\leftc$.
Then it can be directly checked that
$
\frac{\dd G_{p_1}(D\leftvphi,\ivphi,P_{\beta})}{\dd\beta}>0
$
for all $0<\beta<\frac{\pi}{2}$.

It follows from $\underset{\beta\to 0+}{\lim}(\leftrho, \xi_2^m)=(\rightrho,\neta)$ that
$\lim_{\beta\to 0+}G_{p_1}(D\leftvphi,\ivphi,P_{\beta})=-\infty$.

Relations \eqref{2-4-a4} and \eqref{1-36}  yield $\xi_2^m=\leftq\cos\beta$, which gives that
$$
G_{p_1}(D\leftvphi,\ivphi,P_{\beta})
=\leftrho\big((1-\oM^2)\sin\beta-\oM^2\cos\beta^2\csc\beta\big)
-\frac{\cos\beta}{\tan\beta}.
$$
It is shown in the proof of Lemma \ref{lemma:interval-existence}
that $\lim_{\beta\to \frac{\pi}{2}-}\leftc=\infty$
and $\frac{\dd\oM}{\dd\beta}<0$ for all $0<\beta<\frac{\pi}{2}$.
This implies that
$
\lim_{\beta\to \frac{\pi}{2}-}G_{p_1}(D\leftvphi,\ivphi,P_{\beta})=\infty.
$
Therefore, there exists a unique $\beta_{\rm d}^{(\iv)}\in(0,\frac{\pi}{2})$ satisfying \eqref{11-a4}.
\end{proof}
\end{lemma}

\chapter{Non-Existence of Self-Similar Strong Shock Solutions}
\label{section-nonexistence}
\numberwithin{equation}{chapter}

For the completeness of this monograph, we include the proof of the non-existence of admissible solutions
corresponding to $(\iv, \beta)\in\mathfrak{R}_{\rm strong}$
in the sense of Definition \ref{def-regular-sol},
or equivalently, the non-existence of admissible solutions corresponding
to $(\iu, u_0)\in \mathfrak{P}_{\rm strong}$ in the sense of
Definition \ref{def-regular-sol-op}. The non-existence of self-similar strong shock
solutions was first studied in Elling \cite{Elling2}.
In this appendix, we combine the convexity of the shock polar shown in Lemma \ref{lemma-appendix}
for steady potential flow with the result from \cite{Elling2} to show the non-existence
of admissible solutions corresponding to $(\iv, \beta)\in\mathfrak{R}_{\rm strong}$.

\begin{propositionnolabel}[Non-existence of admissible solutions with a strong shock]
For each $\gam \ge 1$, there is no admissible solution corresponding to $(\iv,\beta)\in \mathfrak{R}_{\rm strong}$
in the sense of Definition {\rm \ref{def-regular-sol}};
equivalently, there is no admissible solution corresponding to $(\iu, u_0)\in \mathfrak{P}_{\rm strong}$.
\end{propositionnolabel}

\begin{proof}
The proof is divided into six steps.

\smallskip
{\bf 1.}
On the contrary, suppose that there is an admissible solution $\vphi$ for some $(\iv, \beta)\in \mathfrak{R}_{\rm strong}$ in
the sense of Definition \ref{def-regular-sol}.
Then
$\psi:=\vphi-\leftvphi\in C^3(\ol{\Om}\setminus(\ol{\leftsonic}\cup\ol{\rightsonic}))\cap C^1(\ol{\Om})$ satisfies
\begin{align}
\label{fbp-eq}
(c^2-\vphi_{\xin}^2)\psi_{\xin\xin}-2\vphi_{\xin}\vphi_{\etan}\psi_{\xin\etan}
+(c^2-\vphi_{\etan}^2)\psi_{\etan\etan}=0\qquad&\tx{in $\Om$},\\
\label{fbp-shockbc}
\psi=\ivphi-\leftvphi,\quad \mathfrak{g}(D\psi,\psi, \bmxi)=0\qquad&\tx{on $\shock$},\\
\label{fbp-sonicbc}
\psi=|D\psi|=0\quad\tx{on}\,\,\, \leftsonic,\qquad\,\, \psi=\rightvphi-\leftvphi\qquad&\tx{on $\rightsonic$},\\
\label{fbp-wedgebc}
\der_{\etan}\psi=0\qquad&\tx{on $\Wedge$}
\end{align}
for $c^2=c^2(|D\vphi|^2, \vphi)$ and $\leftsonic=\{P_{\beta}\}$ by \eqref{10-a3},
where
\begin{equation}
\label{nonex-1}
\begin{split}
&\mathfrak{g}({\bf q}, z, \bmxi):=G(D\leftvphi(\bmxi)+{\q}, \leftvphi(\bmxi)+z, \bmxi),\\
&G({\q},z,\bmxi):=\big(\rho(\q,z)\q-D\ivphi(\bmxi)\big)\cdot \frac{D\ivphi(\bmxi)-\q}{|D\ivphi(\bmxi)-\q|},\\
&\rho(\q,z):=\begin{cases}
\left(1+(\gam-1)(\frac 12 \iv^2-\frac 12|\q|^2-z)\right)^{\frac{1}{\gam-1}}&\tx{for}\;\;\gam>1,\\[1mm]
\exp\big(\frac{\iv^2}{2}-\frac 12|{\bf q}|^2-z\big)&\tx{for}\;\;\gam=1,
\end{cases}\\[1mm]
& c^2(|\q|^2, z)=\rho^{\gam-1}(|\q|^2, z),
\end{split}
\end{equation}
for ${\q}\in \R^2, z\in \R$, and $\bmxi\in\ol{\Om}$.

\medskip
{\bf 2.} {\emph{Claim}: $\psi$ attains its minimum
at $P_{\beta}$.}

\smallskip
Since \eqref{fbp-sonicbc}, combined with  Remark \ref{remark-vphibeta},
implies that $\psi$ is not a constant in $\ol{\Om}$,
then the minimum of $\psi$ over $\ol{\Om}$ is attained on $\der\Om$ by
the strong maximum principle.
Also, $\psi$ cannot attain its minimum over $\ol{\Om}$ on $\Wedge$ by Hopf's lemma.
The proof of Proposition \ref{proposition-3} applies
to $\vphi$ such that $\shock$ lies strictly below $\leftshock$, and $\psi>0$ on $\shock$.
Therefore, we conclude that
$\min_{\ol{\Om}}\psi=\psi(P_{\beta})=0$.

\smallskip
{{\bf 3.}} Divide equation \eqref{fbp-eq} by $c^2(|D\vphi|^2, \vphi)$ to rewrite \eqref{fbp-eq} as
\begin{equation*}
\mcl{L}\psi:=\big(1-\frac{|D\leftvphi(P_{\beta})|^2}{\leftc^2}+\mcl{O}_{11}(\bmxi)\big)\psi_{\xin\xin}
+2\mcl{O}_{12}(\bmxi)\psi_{\xin\etan}+\bigl(1+O_{22}(\bmxi)\bigr)\psi_{\etan\etan}=0\qquad\tx{in}\,\,\Om
\end{equation*}
for  $\bmxi=(\xin,\etan)\in\Om$, where each $\mcl{O}_{ij}=\mcl{O}_{ij}(D\vphi,\vphi)$  satisfies
that $\lim_{\bmxi\to P_{\beta}}|\mcl{O}_{ij}(\bmxi)|=0$ for $i,j=1,2$.
Define $k:=\frac{1}{\sqrt{1-{|D\leftvphi(P_{\beta})|^2}/{\leftc^2}}}$
and $\tilde{\xi}_1:=k (\xin-\xib)$.
Let $(r,\theta)$ be the polar coordinates
of $(\tilde{\xi}_1,\etan)$ centered at $P_{\beta}$.
Then
$\Om\subset\{r>0, 0<\theta< \tilde{\beta}\}$
for $\tan\tilde{\beta}=\frac{\tan\beta}{k}$.

Next, define
\begin{equation}
\label{bPsi-def}
\Psi(r,\theta):=\eps r\cos(\omega_0 \theta)
\end{equation}
for constants $\eps, \omega_0>0$ to be determined later.
As in \cite{Elling2}, choose $\eps>0$ small and $\omega_0\in(0,1)$  close to 1.
A direct computation by using the definition of $(r,\theta)$ shows that
\begin{equation}
\label{nonex-2}
\mcl{L}\Psi=\frac{\eps}{r}(1-\omega_0^2)\left(\cos(\omega_0\theta)+\mcl{O}_1^{\rm (polar)}(r,\theta)\right)\qquad\tx{in $\Om$},
\end{equation}
with
$
\lim_{r\to 0+}|\mcl{O}_1^{\rm (polar)}(r,\theta)|=0.
$

A direct computation by using  \eqref{app-1} and Lemma \ref{lemma2-appendix} leads to 
\begin{equation*}
\mathfrak{g}_{{\bf q}}({\bf 0}, 0, P_{\beta})\cdot(\cos\beta,\sin\beta)<0
< \mathfrak{g}_{\bf q}({\bf 0}, 0, P_{\beta})\cdot (1,0).
\end{equation*}
Therefore, there exists $\theta_0\in(-\frac{\pi}{2}, -\frac{\pi}{2}+\beta)$ such that
$
\frac{\mathfrak{g}_{{\bf q}}({\bf 0}, 0, P_{\beta})}{|\mathfrak{g}_{{\bf q}}({\bf 0}, 0, P_{\beta})|}
=(\cos\theta_0, \sin\theta_0).
$
Then it can directly be checked that
\begin{equation}
\label{nonex-3}
\mathfrak{g}_{{\bf q}}({\bf 0}, 0, P_{\beta})\cdot D_{\bmxi}\Psi(r,\theta)
=\eps\big(k\cos\theta_0\cos((1-\omega_0)\theta)+\mcl{O}_2^{\rm (polar)}(\theta)\big),
\end{equation}
where
$|\mcl{O}_2^{\rm (polar)}(\theta)|\le C^{\sharp}|1-\omega_0|$ for all $\theta\in[0, \tilde{\beta}]$
with a constant $C^{\sharp}>0$ chosen independently of $\eps$ and $r$.

\medskip
{\bf 4.}
{\emph{Claim{\rm :} There exist $\om_*\in(0,1)$ and $R_2>0$ such that, whenever $\om_0\in[\om_*, 1)$
in \eqref{bPsi-def} and $R\le R_2$, the minimum of $\psi-\Psi$ over $\ol{\Om\cap B_R(P_{\beta})}$
cannot be attained on $\shock\cap B_R(P_{\beta})$.
Furthermore, $\om_*$ and $R_2$ can be chosen independently of $\eps$. }}

\smallskip
Suppose that $\displaystyle{(\psi-\Psi)(P_*)=\min_{\ol{\Om\cap B_{R}(P_{\beta})}}(\psi-\Psi)}$
for $P_*\in \shock\cap \Om_R(P_{\beta})$ for some $R>0$.
Since $\psi-\Psi=0$ at $P_{\beta}$,
$\psi-\Psi\le 0$  at $P_*.$
Let ${\bm\nu}_{\rm sh}$ be the unit normal vector to $\shock$ at $P_*$ oriented towards the interior of $\Om$,
and let ${\bm\tau}_{\rm sh}$ be a unit tangent vector to $\shock$ at $P_*$.
Then $\psi-\Psi$ satisfies
\begin{equation}
\label{nonex-6}
\der_{{\bm \tau}_{\rm sh}}(\psi-\Psi)(P_*)=0,\qquad \der_{{\bm \nu}_{\rm sh}}(\psi-\Psi)(P_*)\ge 0.
\end{equation}
Let $P_{\beta}P_*'$
be the projection of
$P_{\beta}P_*$ onto $\leftshock$.
Since $(\ivphi-\leftvphi)(P_*')=0$, it follows from \eqref{2-4-b6} and \eqref{12-25}--\eqref{2-4-a0} that
\begin{equation*}
\eps |P_*-P_{\beta}|\ge \Psi(P_*)-\Psi(P_{\beta})\ge \psi(P_*)=[(\ivphi-\leftvphi)(\bmxi)]^{P_*}_{\bmxi=P_*'}\ge \iv \sec\beta\, |P_*-P_*'|,
\end{equation*}
which yields 
\begin{equation}
\label{nonex-7}
|P_*-P_*'|\le \frac{\eps}{\iv \sec\beta}\,|P_*-P_{\beta}|.
\end{equation}
From \eqref{nonex-6}, we have
\begin{equation}
\label{nonex-9}
D\psi(P_*)=D\Psi(P_*')+\big(D\Psi(P_*)-D\Psi(P_*')\big)+|D(\psi-\Psi)(P_*)|{\bm\nu_{\rm sh}}.
\end{equation}

Since $|D(\ivphi-\vphi)\cdot{\bm\nu}_{\rm sh}|>0$ on $\ol{\shock}$,
there exist constants $\hat{\eps}, \delta>0$
such that
%
$|D(\ivphi-\vphi)|\ge \delta$
on the open $\hat{\eps}$--neighborhood $\mcl{N}_{\hat{\eps}}(\shock)$ of $\shock$.
Since $\psi=\ivphi-\leftvphi$ on $\shock$,
$\mathfrak{g}(D\psi, \psi, \bmxi)=\mathfrak{g}(D\psi, \ivphi-\leftvphi, \bmxi)$ on $\shock$. Define
$
\mathfrak{g}_{\sharp}(\q, \bmxi):=\mathfrak{g}(\q, (\ivphi-\leftvphi)(\bmxi),\bmxi).
$
Choose constants $\sigma_0, R_1>0$ small so that

\smallskip
\begin{itemize}
\item[(i)] $\mathfrak{g}_{\sharp}(\q, \bmxi)$ is well defined in $U_{\sigma_0, R_1}=\{(\q,\bmxi)\,:\,|\q|\le 2\sigma_0, |\bmxi-P_{\beta}|\le 2R_1\}$;

\smallskip
\item[(ii)] There is a constant $C_{\mathfrak{g}}>0$ such that
\begin{equation}\label{nonex-8}
\begin{split}
&\|\mathfrak{g}_{\sharp}\|_{C^1(\ol{U_{\sigma_0, R_1}})}\le C_{\mathfrak{g}},\\
&\der_{\q}\mathfrak{g}_{\sharp}(\q, \bmxi)\cdot \frac{D\ivphi(\bmxi')-\q'}{|D\ivphi(\bmxi')-\q'|}\ge  C_{\mathfrak{g}}^{-1}
\qquad \quad  \tx{for $(\q, \bmxi), (\q', \bmxi')\in U_{\sigma_0, R_1}$}.
\end{split}
\end{equation}
Such a constant $C_{\mathfrak{g}}$ can be chosen independently of $(\eps,\omega_0)$.
\end{itemize}

Owing to $|D\psi(P_{\beta})|=0$,
there exists $R_1>0$ small, depending on $\sigma_0$, such
that $(D\psi(\bmxi), \bmxi)\in U_{\sigma_0, R_1}$ for all $\bmxi\in\ol{\Om\cap B_{R_1}(P_{\beta})}$.

If $P_*\in \ol{\Om\cap B_{{R_1}/{2}}(P_{\beta})}$ and $\frac{\eps}{\iv \sec\beta}\le \frac 14$,
then \eqref{nonex-7} implies that $P_*'\in B_{{3R_1}/{4}}(P_{\beta})$.
Choose $\eps_1\in(0, \frac{\iv\sec\beta}{4}]$ so that, whenever $\eps\in(0, \eps_1]$,
$(D\Psi(P_*'), P_*')\in U_{\sigma_0, R_1}$.
Note that $\eps_1$ can be chosen depending only on $\sigma_0$.
Then
\begin{equation*}
\begin{split}
0&=\mathfrak{g}_{\sharp}(D\psi(P_*), P_*)-\mathfrak{g}_{\sharp}({\bf 0}, P_*')\\
&=\big(\mathfrak{g}_{\sharp}(D\psi(P_*), P_*)-\mathfrak{g}_{\sharp}(D\psi(P_*), P_*')\big)
+\big(\mathfrak{g}_{\sharp}(D\psi(P_*), P_*')-
\mathfrak{g}_{\sharp}({\bf 0}, P_*')\big)\\
&=:J_1+J_2.
\end{split}
\end{equation*}
By \eqref{nonex-7} and \eqref{nonex-8}, $J_1$ is estimated as
\begin{equation}
\label{estimate-nonex-J1}
|J_1|\le \frac{C_{\mathfrak{g}}\eps}{\iv \sec \beta}|P_*-P_{\beta}|.
\end{equation}
$J_2$ is estimated more carefully by using \eqref{nonex-3} and \eqref{nonex-7}--\eqref{nonex-8}
as follows:
\begin{equation*}
\begin{split}
J_2
&=\big(D \Psi(P_*')+(D\Psi(P_*)-D\Psi(P_*'))+|D(\psi-\Psi)(P_*)|{\bm\nu_{\rm sh}}\big)
\cdot\int_0^1\der_{\q}\mathfrak{g}_{\sharp}(tD\psi(P_*), P_*')\,\dd t\\
&\ge
\bigl(D \Psi(P_*')+(D\Psi(P_*)-D\Psi(P_*'))\big)
\cdot\int_0^1\der_{\q}\mathfrak{g}_{\sharp}(tD\psi(P_*), P_*')\,\dd t.
\end{split}
\end{equation*}
Let $C^{\sharp}$ be from Step 3. By \eqref{nonex-3} and \eqref{nonex-8},
\begin{equation*}
\begin{split}
&D \Psi(P_*')\cdot\int_0^1\der_{\q}\mathfrak{g}_{\sharp}(tD\psi(P_*), P_*')\,\dd t\\
&\ge \eps\big(k\cos\theta_0\cos((1-\omega_0)\beta)-C^{\sharp}|1-\om_0|-C|P_*-P_{\beta}|^{\alp}\big)
\end{split}
\end{equation*}
for some $C>0$ depending on $C_{\mathfrak{g}}$ and $\|\psi\|_{C^{1,\alp}(\ol{\Om})}$.
By \eqref{bPsi-def}, \eqref{nonex-7}, and \eqref{nonex-8},
\begin{equation*}
(D\Psi(P_*)-D\Psi(P_*'))
\cdot\int_0^1\der_{\q}\mathfrak{g}_{\sharp}(tD\psi(P_*), P_*')\dd t
\ge C\eps^2|P_*-P_{\beta}|
\end{equation*}
for some $C>0$ depending on $C_{\mathfrak{g}}$. Therefore, $J_2$ is estimated as
\begin{equation*}
J_2\ge \eps\big(k\cos\theta_0\cos((1-\omega_0)\beta)-C^{\sharp}|1-\om_0|-Ch(|P-P_{\beta}|)\big)
\end{equation*}
for  a non-increasing continuous function $h(r)$ that tends to $0$ as $r$ tends to $0$,
where $C^{\sharp}$ and $C$ are chosen, independent of $P_*$ and $\om_0$.
Combine this estimate with \eqref{estimate-nonex-J1} to obtain
\begin{equation}\label{nonex-10}
\eps\left(k\cos\theta_0\cos((1-\omega_0)\beta)-C^{\sharp}|1-\om_0|-C\big(h(|P_*-P_{\beta}|)+|P_*-P_{\beta}|\big)\right)\le 0.
\end{equation}
Choose $\omega_*\in(0,1)$ close to $1$ and $R_2\in(0,R_1]$ small,
so that
$$
\varepsilon k\cos\theta_0\cos((1-\omega_*)\beta)-C^{\sharp}|1-\omega_*|-C(h(R_2)+R_2)\ge \frac{\eps}{2} k\cos\theta_0.
$$
Under such choices of $(\om_*, R_2)$,
we arrive at a contradiction whenever $\om_0\in[\om_*,1)$ and $P_*\in \shock\cap B_{R_2}(P_{\beta})$.
Thus, $\psi-\Psi$ cannot attain its minimum on $\shock\cap B_{R}(P_{\beta})$
whenever $\om_0\in[\om_*,1)$ and $R\le R_2$.

\medskip
{{\bf 5.}}
{\emph{Claim{\rm :}
Let $\om_*$ and $R_2$ be from Step {\rm 4}.
There exist $\eps>0$, $\omega_0\in[\om_*, 1)$, and $R\in(0, R_2]$ such that,
for $\Psi$ defined by \eqref{bPsi-def},
$\psi-\Psi$ attains its minimum over $\Om_R(P_{\beta}):=\Om\cap B_R(P_{\beta})$ at $P_{\beta}$.}}

\smallskip
By \eqref{nonex-2}, there exists
a small constant $R_3\in(0, R_2]$ so that $\mcl{L}$ is uniformly elliptic in $\Om_{R_3}(P_{\beta})$ and
\begin{equation*}
\mcl{L}(\psi-\Psi)\le -\frac{\eps}{2R_3}(1-\omega^2_0)\cos(\omega_0\tilde{\beta})\qquad\,\,\tx{in $\Om_{R_3}(P_{\beta})$}.
\end{equation*}
By the strong maximum principle and Hopf's lemma,
the minimum of $\psi-\Psi$ over $\ol{\Om_{R}(P_{\beta})}$ must be attained on $\der \Om_{R_3}(P_{\beta})\setminus \Wedge$.
It is shown in Step 4 that $\psi-\Psi$ cannot attain its minimum on $\shock\cap B_{R_3}(P_{\beta})$.

Denote $m:=\inf_{\Om\cap \der B_{R_3}(P_{\beta})}\psi$. The claim in Step 2 implies that $m>0$.
Choose $\eps>0$ small, depending only on $R_3$, so that $\psi-\Psi>0$ on $\Om\cap \der B_{R_3}(P_{\beta})$.
For such a choice of $\eps$, since $(\psi-\Psi)(P_{\beta})=0$, we conclude that
$$
\min_{\Om_{R_3}(P_{\beta})}(\psi-\Psi)=(\psi-\Psi)(P_{\beta})=0.
$$

\smallskip
{{\bf 6.}} In Steps 4--5, it is shown that we can choose $(\eps, \om_0)$ in \eqref{bPsi-def}
so that, if $R>0$ is sufficiently small, the minimum of $\psi-\Psi$ over $\Om_R(P_{\beta})$
must be attained at $P_{\beta}$, provided that there is an admissible solution $\vphi$ corresponding
to some $(\iv, \beta)\in \mathfrak{R}_{\rm{strong}}$, and that $\psi$ is given by $\psi=\vphi-\leftvphi$.

By the definition of $\Psi$ with $\omega_0\in(0,1)$ and \eqref{fbp-sonicbc}, and by the $C^1$--regularity of $\vphi$ up to $P_{\beta}$,
there exists a small constant $\delta>0$ so that $\der_r(\psi-\Psi)<-\frac{\eps}{2}$ in $\Om_{\delta}(P_{\beta})$.
However, this contradicts the fact that
$$
(\psi-\Psi)(P_{\beta})=\min_{\Om_R(P_{\beta})}(\psi-\Psi).
$$
Therefore, we conclude that there exists no admissible solution corresponding to $(\iv, \beta)\in \mathfrak{R}_{\rm{strong}}$
in the sense of Definition \ref{def-regular-sol}.
\end{proof}



\chapter{Quasilinear Elliptic Equations in Two Variables}
\label{Appendix-elliptic-PDEs}
\numberwithin{equation}{section}

For the completeness of this work, this appendix includes several properties of quasilinear elliptic equations,
which are used to prove Theorem \ref{theorem-0}.
We refer the reader to \cite{CF2} for the proofs of these properties as stated below.

\section{Ellipticity Principle for Self-Similar Potential Flow}
The following lemma is an extension of the ellipticity principle of Elling-Liu \cite{EL}:

\begin{lemma}[Theorem 5.2.1, \cite{CF2}]\label{lemma6.8}
Fix $\gam\ge 1$ and $\iv>0$.
In a bounded domain $\Om\subset \R^2$, let $\vphi\in C^3(\Om)$ satisfy the equation{\rm :}
\begin{equation}
\label{2-1repeat}
{\rm div}\left(\rho(|D\vphi|^2,\vphi)D\vphi\right)+2\rho(|D\vphi|^2,\vphi)=0
\end{equation}
for $\rho(|D\vphi|^2,\vphi)$ given by \eqref{new-density}.
Denote the pseudo-Mach number as $M:=\frac{|D\vphi|}{c(|D\vphi|^2,\vphi)}$ for
$c(|D\vphi|^2, \vphi)=\rho^{\frac{(\gam-1)}{2}}(|D\vphi|^2,\vphi)$.
Let $\vphi$ satisfy that $\rho>0$ and $M\le 1$ in $\Om$.
Then the following properties hold{\rm :}

\smallskip
\begin{itemize}
\item[{\rm (a)}]
Either $M\equiv 0$ holds in $\Om$ or $M$ does not attain its maximum in $\Om${\rm ;}

\smallskip
\item[{\rm (b)}]
Suppose that ${\rm diam}(\Om)\le d$ for some constant $d>0$.
Then there exists a constant $C_0>0$ depending only on $(\iv, \gam,  d)$
such that, for any given $\delta\ge 0$, $\hat c\ge 1$, and $b\in C^2(\Om)$ with
$|Db|+\hat c|D^2b|\le \frac{\delta}{\hat c}$ in $\Om$,
if $c(|D\vphi|^2, \vphi)\le \hat c$ holds in $\Om$,
then either $M^2\le C_0\delta$ holds in $\Om$ or $M^2+b$ does not attain its maximum in $\Om$.
\end{itemize}
\end{lemma}

\begin{lemma}
[Theorem 5.3.1, \cite{CF2}]\label{lemma6.8-slipbc}
In a bounded domain $\Om\subset \R^2$ with  a relatively open flat segment $\Gam\subset \der \Om$,
let $\vphi\in C^3(\Om\cup \Gam)$ satisfy \eqref{2-1repeat} in $\Om$ and
\begin{equation*}
\der_{\bm\nu}\vphi=0\qquad\tx{on $\Gam$}
\end{equation*}
for the unit normal vector ${\bm\nu}$ to $\Gam$ towards the interior of $\Om$.
Assume that $\rho>0$ and $M\le 1$ in $\Om\cup\Gam$.
Then the following properties hold{\rm :}

\smallskip
\begin{itemize}
\item[{\rm (a)}]
Either $M\equiv 0$ holds in $\Om \cup \Gam$ or $M$ does not attain its maximum in $\Om \cup \Gam${\rm ;}

\smallskip
\item[{\rm (b)}]
Let ${\rm diam}(\Om)\le d$ for some constant $d>0$.
Then there exists a constant $C_0>0$ depending only on $(\iv, \gam,  d)$ such that,
for any given $\delta\ge 0$, $\hat c\ge 1$, and $b\in C^2(\Om)$ with $|Db|+\hat c|D^2b|\le \frac{\delta}{\hat c}$
in $\Om$ and $\der_{\bm\nu}b=0$ on $\Gam$,
if $c(|D\vphi|^2, \vphi)\le \hat c$ holds in $\Om\cup \Gam$, then either $M^2\le C_0\delta$ holds in $\Om\cup \Gam$
or $M^2+b$ does not attain its maximum in $\Om\cup \Gam$.
\end{itemize}
\end{lemma}

\section{Uniformly Elliptic Equations Away From the Corners}
Consider a quasilinear elliptic equation of the form:
\begin{equation}
\label{ql-eq-app}
\mcl{N}(u)=f({\bf x})\qquad\tx{in $\Om$},
\end{equation}
with
\begin{equation*}
\mcl{N}(u):= \sum_{i,j=1}^2 A_{ij}(Du, u, {\bf x})D_{ij}u+A(Du, u, {\bf x}),
\end{equation*}
where
\begin{equation}
\label{app-c-prop0}
\quad A_{ij}({\bf p}, z, {\bf x})=A_{ji}({\bf p}, z, {\bf x}),\,\,\,
A({\bf 0}, 0, {\bf x})=0\quad \tx{for all}\,\,
({\bf p}, z, {\bf x})\in \R^2\times \R\times \Om
\,\,\tx{and}\,\,i,j=1,2.
\end{equation}
Suppose that there exist $\lambda>0$ and $\alp\in(0,1)$ such that
\begin{align}
\label{app-c-prop1}
&\lambda|{\bm\mu}|^2\le \sum_{i,j=1}^2 A_{ij}(Du({\bf x}), u({\bf x}), {\bf x})\mu_i\mu_j \le \lambda^{-1}|{\bm \mu}|^2
\quad\tx{for all ${\bf x}\in \Om$ and ${\bm\mu}=(\mu_1, \mu_2)\in \R^2$},\\
\label{app-c-prop2}
&\|(A_{ij}, A)({\bf p}, z,\cdot)\|_{0, \alp, \ol{\Om}}\le \lambda^{-1}\quad\tx{for all $({\bf p}, z)\in \R^2\times \R$},\\
\label{app-c-prop3}
&\|D_{({\bf p},z)}(A_{ij}, A)\|_{0, \R^2\times \R\times \ol{\Om}}\le \lambda^{-1}.
\end{align}

For $r>0$, let $B_r$ denote a ball of radius $r$ in $\R^2$.
\begin{theorem}[Theorem 4.2.1, \cite{CF2}] \label{elliptic-t1-CF2}
For $\Om=B_2$,  if $u\in C^{2,\alp}(B_2)$ is a solution of \eqref{ql-eq-app} with
\begin{equation*}
\|u\|_{0, B_2}+\|f\|_{0,\alp, B_2}\le M,
\end{equation*}
then there exists a constant $C>0$ depending only on $(\lambda, M, \alp)$ such that
\begin{equation*}
\|u\|_{2,\alp, B_1}\le C\big(\|u\|_{0, {B_2}}+\|f\|_{0,\alp, {B_2}}\big).
\end{equation*}
\end{theorem}

Applying Theorem \ref{elliptic-t1-CF2} to $v(x)=\frac 1r u(rx)$, we have the following corollary:
\begin{corollary}
\label{corollary-t1-CF2}
  If $u\in C^{2,\alp}(B_{2r})$ is a solution of \eqref{ql-eq-app} for $r\in(0,1]$ with
  \begin{equation*}
\|u\|_{0, B_{2r}}+\|f\|_{0,\alp,B_{2r}}\le M,
\end{equation*}
then there exists a constant $C>0$ depending only on $(\lambda, M, \alp)$ such that
\begin{equation*}
\|u\|_{2,\alp, B_r}\le \frac{C}{r^{2+\alp}}\big(\|u\|_{0, {B_{2r}}}+r^2\|f\|_{0,\alp, {B_{2r}}}\big).
\end{equation*}
\end{corollary}

\begin{theorem}[Theorem 4.2.3, \cite{CF2}] \label{elliptic-t8-CF2}
For $\lambda\in(0,1)$, let $\Phi\in C^1(\R)$ satisfy
\begin{equation*}
\|\Phi\|_{1,\R}\le \lambda^{-1},\qquad \Phi(0)=0.
\end{equation*}
For $R>0$, denote
\begin{equation*}
\Om_R:=B_R({\bf 0})\cap \{x_2>\eps \Phi(x_1)\},\qquad
\Gam_R:=B_R({\bf 0})\cap \{x_2=\eps \Phi(x_1)\}.
\end{equation*}
In addition to assumptions \eqref{app-c-prop0}--\eqref{app-c-prop3} with $\Om=\Om_{2r}$, let $W(p_2, z, x)$ satisfy
\begin{equation*}
\begin{split}
W(0,0, {\bf x})=0&\qquad\tx{on $\Gam_{2r}$},\\
|\der_{p_2} W(p_2, z, {\bf x})|\le \eps&\qquad\tx{for all $(p_2, z, {\bf x})\in \R\times \R\times \Gam_{2r}$},\\
\|D_{(p_2, z)}{W}(p_2, z, \cdot)\|_{1,\Gam_{2r}}\le \lambda^{-1}&\qquad \tx{for all $(p_2, z)\in \R\times \R$}.
\end{split}
\end{equation*}
Then there exist constants $\eps, \beta\in(0,1)$ and $C>0$ depending only on $\lambda$
such that, for $u\in C^2(\Om_{2r})\cap C^{1,\beta}(\Om_{2r}\cup \Gam_{2r})$ satisfying \eqref{ql-eq-app}
with $f=0$ in $\Om_{2r}$ and
\begin{equation}
\label{bc-almost-tangential}
u_{x_1}=W(u_{x_2}, u, {\bf x})\qquad\,\, \tx{on $\Gam_{2r}$},
\end{equation}
the following estimate holds{\rm :}
\begin{equation*}
\|u\|_{1,\beta, {{\Om_{9r/5}}}}\le \frac{C}{r^{1+\beta}}\|u\|_{0, \Om_{2r}}.
\end{equation*}
\end{theorem}

\begin{theorem}[Theorem 4.2.8, \cite{CF2}]
\label{elliptic-t9-CF2}
In addition to the assumptions of Theorem {\rm \ref{elliptic-t8-CF2}}, for $\alp\in(0,1)$, assume that
\begin{align*}
&\|\Phi\|_{1,\alp,\R}\le \lambda^{-1},\\
&\|D_{(p_2, z)}W(p_2, z, \cdot)\|_{1,\alp, \Gam_{2r}}\le \lambda^{-1}\qquad \tx{for all $(p_2, z)\in \R\times \R$},\\
&\|D^2_{(p_2, z)}W\|_{1,0, \R\times \R\times {\Gam_{2r}}}\le \lambda^{-1}.
\end{align*}
Then there exist $\eps\in(0,1)$ and $C>0$ depending only on $(\lambda, \alp, \|u\|_{0, {\Om_{2r}}})$ such that,
for $u\in C^{2,\alp}(\Om_{2r}\cup \Gam_{2r})$ satisfying \eqref{ql-eq-app} with $f=0$ in $\Om_{2r}$
and \eqref{bc-almost-tangential} on $\Gam_{2r}$,
\begin{equation*}
\|u\|_{2,\alp, {\Om_{9r/5}}}\le \frac{C}{r^{2+\alp}}\|u\|_{0, {\Om_{2r}}}.
\end{equation*}
\end{theorem}

\begin{theorem}[Theorem 4.2.10, \cite{CF2}] \label{elliptic-t2-CF2}
For $\lambda\in(0,1)$ and $\alp\in(0,1)$, let
$\Phi \in C^{2,\alp}(\R)$ satisfy
\begin{equation*}
\|\Phi\|_{2,\alp, \R}\le \lambda^{-1},\qquad \Phi(0)=\Phi'(0)=0,
\end{equation*}
and denote
\begin{equation*}
\Om_R:= B_R(0)\cap \{x_2>\Phi(x_1)\},\quad \Gam_R:=\der \Om_R\cap \{x_2=\Phi(x_1)\}\qquad\,\,{\tx{for $R\in(0,2)$}}.
\end{equation*}
Let $u \in C^{2,\alp}(\Om_R\cup \Gam_R)$ satisfy \eqref{ql-eq-app} in $\Om_R$ and
\begin{equation*}
{{\bm\omega}}\cdot Du+b_0 u=h\qquad \tx{on}\,\,\Gam_R.
\end{equation*}
Assume that ${{{\bm\omega}=(\om_1, \om_2)({\bf x})}}$ and $b_0=b_0({\bf x})$ satisfy the following conditions{\rm :}
\begin{align*}
{{\bm\omega}}\cdot{\bm\nu}\ge \lambda\,\,\,\, \tx{on}\,\,\Gam_R,\qquad
&\|({{\bm\omega}}, b_0)\|_{1,\alp, {{\Gam_R}}}\le \lambda^{-1},
\end{align*}
where ${\bm\nu}$ represents the unit normal vector to $\Gam_R$ towards the interior of $\Om_R$.
If $u$ satisfies
\begin{equation*}
\|u\|_{0, {{\Om_R}}}+\|f\|_{0,\alp, {{\Om_R}}}+\|h\|_{1,\alp, {{\Gam_R}}}\le M,
\end{equation*}
then there exists a constant $C>0$ depending only on ${{(\lambda, \alp)}}$ such that
\begin{equation*}
\|u\|_{2,\alp, {{\Om_{R/2}}}}
\le {{\frac{C}{R^{2+\alp}}}}\Big(\|u\|_{0, {{\Om_R}}}+{{R^2}}\|f\|_{0,\alp,{{\Om_R}}}+{{R}}\|h\|_{1,\alp, {{\Gam_R}}}\Big).
\end{equation*}
In addition, there exist $\beta\in(0,1)$ and $\hat C>0$ depending only on ${{\lambda}}$ such that
\begin{equation*}
\|u\|_{1,\beta,{{\Om_{R/2}}}}\le {{\frac{\hat C}{R^{1+\beta}}}}\Big(\|u\|_{0,{{\Om_R}}}+{{R^2}}\|f\|_{0,\alp,{{\Om_R}}}+{{R}}\|h\|_{0,\beta, {{\Gam_R}}}\Big).
\end{equation*}
Note that $\beta$ is independent of $\alp$.
\end{theorem}

\begin{theorem}[Theorem 4.3.2, \cite{CF2}]
\label{elliptic-t3-CF2}
Let $R>0$, $\lambda\in (0,1)$, $\gam\in(0,1)$, and $K>0$.
Let $\Phi \in C^1(\R)$ satisfy
\begin{equation*}
\|\Phi\|_{0,1, \R}\le \lambda^{-1},\qquad \Phi(0)=0.
\end{equation*}
Let $\Om_R$ and $\Gam_R$ be as in Theorem {\rm \ref{elliptic-t2-CF2}}  for $R>0$.
Define
\begin{equation*}
d({\bf x}):={\rm dist}({\bf x}, \Gam_R)\qquad\,\, \tx{for ${\bf x}\in \Om_R$}.
\end{equation*}
Assume that $u\in C^3(\Om_R)\cap C^1(\ol{\Om_R})$ is a solution of \eqref{ql-eq-app}
with $f = 0$ in $\Om_R$ and the boundary condition{\rm :}
\begin{equation*}
B(Du, u, {\bf x})=0\qquad \tx{on $\Gam_R$}.
\end{equation*}
Assume that $A_{ij}({\bf p}, z, {\bf x}), i,j=1,2$, and $A({\bf p}, z, {\bf x})$
satisfy \eqref{app-c-prop1}--\eqref{app-c-prop3} and the additional property{\rm :}
\begin{equation*}
d({\bf x})^{\gam}|D_{\bf x}(A_{ij}, A)({\bf p}, z, {\bf x})|\le \lambda^{-1}
\qquad \mbox{for all ${\bf x}\in \Om_R$ and $|{\bf p}|+|z|\le 2K$},
\end{equation*}
and  that $B({\bf{p}}, z, {\bf x})$ satisfies
\begin{align}
\label{app-c-B-cond1}
|D_{{\bf p}}B(Du({\bf x}), u({\bf x}), {\bf x})|\ge \lambda\,\,\,\,\tx{for all ${\bf x}\in \ol{\Om_R}$},
\qquad \|B\|_{2,\{|{\bf p}|+|z|\le 2K, {\bf x}\in\ol{\Om_R}\}}\le \lambda^{-1}.
\end{align}
Assume that $u$ satisfies
\begin{equation*}
|u|+|Du|\le K\qquad \,\, \tx{on $\Om_R\cup\Gam_R$}.
\end{equation*}
Then there exist both $\beta\in(0,1]$ depending only on $(\lambda, K, \gam)$ and $C>0$ depending only on $(R,\lambda, K, \gam)$ such that
\begin{equation*}
\|u\|_{1,\beta, {\Om_{R/2}}}\le C, \qquad
\|u\|_{2,\beta, \Om_{R/2}}^{(-1-\beta), \Gam_{R/2}}\le C.
\end{equation*}
\end{theorem}

\begin{theorem} [Theorem 4.3.4, \cite{CF2}]
\label{elliptic-t4-CF2}
Let the assumptions of Theorem {\rm \ref{elliptic-t3-CF2}} be satisfied with $\gam=0$.
In addition, for $\alp, \sigma\in (0,1)$, assume that
\begin{align*}
&\|\Phi\|_{C^{1,\sigma}(\R)}\le \lambda^{-1},\qquad \Phi(0)=0,\\
&\|(A_{ij}, A)\|_{C^{1,\alp}(\{|{\bf p}|+|z|\le 2K,{{\, {\bf x}\in \ol{\Om_R}}}\})}
+\|B\|_{C^{2,\alp}(\{|{\bf p}|+|z|\le 2K,{\, {\bf x}\in\ol{\Om_R}} \})}\le \lambda^{-1} \qquad \mbox{for $j=1,2$}.
\end{align*}
Then
\begin{equation*}
\|u\|_{{{2,\sigma, \Om_{R/4}}}}\le C,
\end{equation*}
where $C$ depends only on $(\lambda, K, \alp, \sigma, R)$.
\end{theorem}

\begin{corollary}[Corollary 4.3.5, \cite{CF2}]
\label{elliptic-t5-CF2}
Let the assumptions of Theorem {\rm \ref{elliptic-t3-CF2}}
be satisfied with $\gam=0$. In addition, for  $\alp\in(0,1)$ and  $k\in\mathbb{N}$, assume that
\begin{align*}
&\|\Phi\|_{k,\alp,\R}\le \lambda^{-1},\qquad \Phi(0)=0,\\
&\|(A_{ij}, A)\|_{C^{k,\alp}(\{|{\bf p}|+|z|\le 2K,\, {\bf x}\in \ol{\Om_R}\})}
+\|B\|_{C^{k+1,\alp}(\{|{\bf p}|+|z|\le 2K,\, {\bf x}\in \ol{\Om_R}\})}
\le \lambda^{-1}
\qquad \tx{for $j=1,2$}.
\end{align*}
Then
\begin{equation*}
\|u\|_{k+1,\alp, {{\Om_{R/2}}}}\le C,
\end{equation*}
where $C$ depends only on $(\lambda, K, k, \alp, R)$.
\end{corollary}

\section{Quasilinear Degenerate Elliptic Equations}

Consider a domain $U\subset \R^2$ of the form:
\begin{equation*}
U=\{\rx=(x_1,x_2)\,:\,x_1>0, x_2\in (0, f(x_1))\},
\end{equation*}
where $f\in C^{1}({{\R_+}})$ and $f>0$ on $\R_+$. For a constant $r>0$, denote
\begin{equation*}
\begin{split}
&U_r=U\cap \{x_1<r\},\\
&\Gam_{n,r}=\der U\cap \{(x_1,0)\,:\,0<x_1<r\},\\
&\Gam_{f,r}=\der U \cap\{(x_1, f(x_2))\,:\,0<x_1<r\}.
\end{split}
\end{equation*}
Consider a boundary value problem of the form:
\begin{equation}
\label{bvp-app-eqnxy}
\begin{split}
\sum_{i,j=1}^2 A_{ij}(Du, u, \rx) {{\der_{x_ix_j}}}u+\sum_{i=1}^2 A_i(Du, u, \rx){{\der_{x_i}}} u=0\qquad &\tx{in $U_r$},\\
B(Du, u, \rx)=0\qquad &\tx{on $\Gam_{f,r}$},\\
{{\der_{x_2}}} u=0\qquad &\tx{on $\Gam_{n,r}$},\\
u=0\qquad& \tx{on $\Gam_0=\der U\cap \{x_1=0\}$}.
\end{split}
\end{equation}

\begin{theorem}[Theorem 4.7.4, \cite{CF2}]
\label{elliptic-t7-CF2}
Given constants $r>0$, $M\ge 1$, and $l, \lambda\in(0,1)$, assume that the following conditions are satisfied{\rm :}

\smallskip
\begin{itemize}
\item[(i)] {\emph{Conditions for $\Gam_{f,r}$}}{\rm :}
$f$ is in $C^{1,\beta}([0,r])$ for some $\beta\in(0,1)$ and satisfies
\begin{equation*}
\|f\|_{2,\beta, (0,r)}^{(-1-\beta),\{0\}}\le M,\qquad f\ge l \,\,\,\,\tx{on $\R_+$}.
\end{equation*}

\item[(ii)]  {\emph{Conditions for $(A_{ij}, A_i)$}}{\rm :}
For any $({\bf p}, z, \rx)\in \R^2\times \R\times U_r$ and ${\bm\kappa}=(\kappa_1, \kappa_2)\in \R^2$,
\begin{equation*}
\lambda|\bm\kappa|^2\le \sum_{i,j=1}^2 A_{ij}({\bf p}, z, \rx) \frac{\kappa_i\kappa_j}{x_1^{2-\frac{i+j}{2}}}\le \lambda^{-1}|\bm\kappa|^2.
\end{equation*}
In addition, $(A_{ij}, A_i)$ satisfy the following estimates{\rm :}
\begin{align*}
&\|(A_{11}, A_{12})\|_{0,1, \R^2\times \R\times {{U_r}}}\le M,\\
&|{{\der_{x_2}}} A_{11}({\bf p}, z, \rx)|\le M x_1^{1/2}\qquad \tx{in $\R^2\times \R\times U_r$},\\
&\|(A_{22}, A_1,A_2)\|_{0, \R^2\times \R\times {{U_r}}}+
\|D_{({\bf p}, z)}(A_{22}, A_1, A_2)\|_{0, \R^2\times \R\times {{U_r}}}\le M,\\
&\sup_{({\bf p}, z)\in \R^2\times \R, {\rx \in U_r}}
|(x_1 {{\der}}_{x_1}, x_1^{1/2}{{\der}}_{x_2})(A_{22}, A_1, A_2)({\bf p}, z, \rx)|\le M.
\end{align*}

\item[(iii)]  {\emph{Conditions for $B$}}{\rm :} For any $({\bf p}, z, \rx)\in \R^2\times \R \times \Gam_{f,r}$,
\begin{equation}
\label{condition-b-deriv-p1}
\der_{p_1} B({\bf p}, z, \rx)\le -M^{-1}.
\end{equation}
In addition, $B$ satisfies the following estimates{\rm :}
\begin{align*}
\|B\|_{3, \R^2\times \R\times \Gam_{f,r}}\le M,\qquad
B({\bf 0}, 0, \rx)=0\,\,\,\, \tx{on $\Gam_{f,r}$}.
\end{align*}
\end{itemize}
Let $u\in C(\ol{U_r})\cap C^2(\ol{U_r}\setminus \ol{\Gam_0})$ be a solution of the boundary value
problem \eqref{bvp-app-eqnxy} satisfying
that
\begin{equation*}
|u(\rx)|\le Mx_1^2\qquad\,\tx{in $U_r$}.
\end{equation*}
Then, for any $\alp\in(0,1)$, there exist constants $r_0\in(0,1]$ and $C>0$ depending only on $(M,\lambda, \alp)$ such that,
for $\eps:=\min\{\frac r2, \,r_0,\, l^2\}$,
\begin{equation*}
\|u\|_{2,\alp, U_{\eps}}^{(2),{\rm (par)}}\le C.
\end{equation*}
\end{theorem}

\section{Estimates at a Corner for the Oblique Derivative Boundary Value Problems}
\begin{proposition}[Proposition 4.3.7, \cite{CF2}]
\label{appendix4-proposition1}
Let $R>0$, $\beta\in(0,1)$, $\gam\in[0,1)$, $\lambda>0$, and $K, M\ge 1$.
Let $\Om\subset \R^2$ be a domain with ${\rx}_0\in \der\Om$ and $\der \Om\cap B_R({\rx}_0)=\Gam^1\cup\Gam^2$,
where $\Gam^k$, $k=1,2$, are two Lipschitz curves intersecting only at ${\rx}_0$
and contained within ${\rx}_0+\{{\rx}=(x_1,x_2)\in \R^2:x_2>\tau |x_1|\}$ for some $\tau>0$. Denote
\[
\Om_R:=\Om\cap B_R({\rx}_0).
\]
Assume that $\Gam^2$ is $C^{1,\sigma}$ up to the endpoints for some $\sigma\in(0,1)$ with $\|\ol{\Gam^2}\|_{C^{1,\sigma}}\le M$
in the sense that there exist $c^{(2)}>0$ and $f^{(2)}\in C^{1,\sigma}([0, c^{(2)}])$ such that,
in an appropriate basis in $\R^2$,
$$
\Om_R\subset \{{\rx}\,:\, x_2>f^{(2)}(x_1), 0<x_1<c^{(2)}\}, \quad\,
\ol{\Gam^2}=\{x_2=f^{(2)}(x_1)\,:\, 0<x_1<c^{(2)}\}.
$$
Let $u\in C^1(\ol{\Om_R})\cap C^2(\Om_R\cup \Gam^2)\cap C^3(\Om_R)$ satisfy
\begin{equation}
\label{app-D-u-Lip-est}
\|u\|_{C^{0,1}(\ol{\Om_R})}\le K.
\end{equation}
Assume that $u$ is a solution of
\begin{align}
\label{app-D-Eq}
\sum_{i,j=1}^2 a_{ij}(Du, u, \rx)D_{ij}u+a(Du, u, \rx)=0&\qquad\tx{in $\Om_R$},\\
\label{app-D-bc1}
b^{(1)}(Du, u, \rx)=h(\rx)&\qquad\tx{on $\Gam^1$},\\
\label{app-D-bc2}
b^{(2)}(Du, u, \rx)=0&\qquad \tx{on $\Gam^2$},
\end{align}
where $(a_{ij}, a, b^{(k)})$ are defined in $V=\{({\bf p}, z, \rx)\in \R^2\times \R\times \Om\,:\, |{\bf p}|+|z|\le 2K\}$.
Assume that $(a_{ij}, a)\in C(\ol{V})\cap C^1(\ol{V}\setminus \{\rx=\rx_0\})$, $b^{(1)}\in C^2(\ol{V})$, $b^{(2)}\in C^1(\ol{V})$,
and $h\in C(\ol{\Gam^1})$ with
\begin{align}
\label{appendix4-prop1-cond1}
&\|(a_{ij},a)\|_{C^0(\ol{V})}+\|D_{({\bf p},z)}(a_{ij}, a)\|_{C^0(\ol{V})}\le M,\\
&|D_{\rx}(a_{ij}, a)({\bf p}, z, \rx)|\le M|\rx-\rx_0|^{-\gam}\qquad\tx{for all $({\bf p}, z, \rx)\in V$},\\
&\|b^{(1)}\|_{C^2(\ol{V})}+\|b^{(2)}\|_{C^1(\ol{V})}\le M,\\
\label{appendix4-prop1-cond4}
&|h(\rx)-h(\rx_0)|\le  \frac{1}{\lambda R^{\beta}} 
|\rx-\rx_0|^{\beta}\qquad\tx{for all $\rx\in \Gam^1$}.
\end{align}
In addition to the conditions stated above, assume that the following properties hold{\rm :}

\smallskip
\begin{itemize}
\item[(i)] For any $\rx\in \Om_R$ and ${\bm\kappa}=(\kappa_1, \kappa_2)\in \R^2$,
\[
\lambda|{\bm\kappa}|^2\le \sum_{i,j=1}^2 a_{ij}(Du(\rx), u(\rx), \rx)\kappa_i\kappa_j\le \lambda^{-1}|{\bm\kappa}|^2;
\]

\item[(ii)] For any $\rx\in\Gam^1$,
$
|{{D_{\bf p}b^{(1)}}}(Du(\rx), u(\rx), \rx)|\ge \lambda;
$

\smallskip
\item[(iii)] For any $\rx\in \Gam^2$,
$D_{\bf p}b^{(2)}(Du(\rx), u(\rx), \rx)\cdot{\bm\nu}\ge \lambda$,
where ${\bm\nu}$ is the inner unit normal vector to $\Gam^2${\rm ;}

\smallskip
\item[(iv)] $b^{(1)}$ and $b^{(2)}$ are independent for $u$ on $\Gam^2$ in the sense that,
for any $\rx\in \Gam^2$,
$$
\left|{\rm det} \begin{pmatrix}
D_{\bf p}b^{(1)}(Du(\rx), u(\rx), \rx)\\
D_{\bf p}b^{(2)}(Du(\rx), u(\rx), \rx)
\end{pmatrix}\right|\ge \lambda
\qquad \tx{for any $\rx\in \Gam^2$}.
$$
\end{itemize}
Then there exist $\alp\in(0, \beta]$ and $C$ depending only on $(\lambda, K, M)$,
and $R'\in(0, R]$ depending only on $(\lambda, \gam, K, M, \alp)$ so that,
for any $\rx\in {\ol{{\Om_{R'}}}}$,
\[
|b^{(1)}(Du(\rx), u(\rx), \rx)-b^{(1)}(Du(\rx_0), u(\rx_0), \rx_0)|\le C|\rx-\rx_0|^{\alp}.
\]
\end{proposition}

\begin{proposition}[Proposition 4.3.9, \cite{CF2}]
\label{appendix4-proposition2}
In addition to the assumptions of Proposition {\rm \ref{appendix4-proposition1}}, assume that
\begin{equation}
\label{app4-prop2-cond1}
|b^{(k)}({\bf p}, z, \rx)-b^{(k)}(\til{\bf p}, \til z, \til{\rx})|\le M|({\bf p}, z, \rx)-(\til{\bf p}, \til z, \til{\rx})|
\qquad\,\,\tx{for $k=1,2$},
\end{equation}
for all $({\bf p}, z, \rx)$, $(\til{\bf p}, \til z, \til{\rx})\in V$.
Moreover, denoting $h^{(k)}({\bf p})=b^{(k)}({\bf p}, u(\rx_0), \rx_0)$, $k=1,2$,
and noting that functions $h^{(k)}$ are defined in $B_K(Du(\rx_0))$,
assume that $h^{(k)}\in C^{1,\alp}(\ol{B_K(Du(\rx_0))})$ with $\|h^{(k)}\|_{C^{1,\alp}(\ol{B_{K/2}(Du(\rx_0))})}\le M$ for some $\alp\in(0,1)$, and
\begin{equation}
\label{app4-prop2-cond2}
\left|{\rm det} \begin{pmatrix}
D_{\bf p}h^{(1)}(Du(\rx_0))\\
D_{\bf p}h^{(2)}(Du({{\rx_0}}))
\end{pmatrix}\right|\ge \lambda.
\end{equation}
Let $W\subset \ol{\Om_R}$ satisfy
\begin{equation}
\label{app4-prop2-cond4}
\rx_0\in W,\qquad \emptyset\neq W\cap \der B_r(\rx_0)\subset \ol{W\cap B_r(\rx_0)}\quad\tx{for all $r\in(0,R)$.}
\end{equation}
For each $k=1,2$, let
\begin{equation}
\label{app4-prop2-cond3}
|b^{(k)}(Du(\rx), u(\rx), \rx)-b^{(k)}(Du(\rx_0), u(\rx_0), \rx_0)|\le M|\rx-\rx_0|^{\alp}
\qquad \tx{for all $\rx\in W$}.
\end{equation}
Then there exists a constant $C>0$ depending only on $(K, M, R, \alp)$ such that, for all $\rx\in W$,
\begin{equation*}
|Du(\rx)-Du(\rx_0)|\le C|\rx-\rx_0|^{\alp}.
\end{equation*}
\end{proposition}

\begin{proposition}[Proposition 4.3.11, \cite{CF2}]
\label{appendix4-proposition3}
Let $R, \lambda>0$, $\alp\in(0,1]$, $\gam\in[0,1)$,  and $M\ge 1$.

\smallskip
{\rm (a)} Let $\Om_R$ be as in Proposition {\rm \ref{appendix4-proposition1}}.
Assume that $\Gam^1$ and $\Gam^2$ satisfy that, for each $k=1,2$,

\smallskip
\begin{itemize}
\item[(i)]  $\Gam^k\in C^1$ with $\|\Gam^k\|_{C^{0,1}}\le M$,

\smallskip
\item[(ii)]
$B_{\frac{d(\rx)}{M}}(\rx)\cap \der \Om_R=B_{\frac{d(\rx)}{M}}\cap \Gam^k$ for all $\rx\in \Gam^k\cap B_{\frac{3R}{4}}(\rx_0)$,
for $d(\rx):=|\rx-\rx_0|$.
\end{itemize}
Let $u\in C^1(\ol{\Om_R})\cap C^3(\Om_R)$ be a solution of \eqref{app-D-Eq}--\eqref{app-D-bc2}
with $h\equiv 0$, where $(a_{ij}, a)({\bf p}, z, \rx)$ satisfy all the conditions
stated in Proposition {\rm \ref{appendix4-proposition1}}.
In addition, assume that, for each $k=1,2$,
\begin{equation*}
\begin{split}
&\|b^{(k)}\|_{C^2(\ol{V})}\le M,\\
&|D_{\bf p}b^{(k)}(Du(\rx), u(\rx), \rx)|\ge \lambda\qquad \tx{for all $\rx\in \Om_R$}.
\end{split}
\end{equation*}
Moreover, assume that $u$ satisfies
\begin{equation}
\label{condition-propc14-grad-est}
|Du(\rx)-Du(\rx_0)|\le M|\rx-\rx_0|^{\alp}\qquad\tx{for all $\Om_R$}.
\end{equation}
Then there exist $\beta\in(0,\alp]$ depending only on $(\lambda, K, M, \alp)$
and $C>0$ depending on $(\lambda, K, M, R, \alp)$ such that $u\in C^{1,\beta}(\ol{\Om_{R/2}})$ with
\[
\|u\|_{C^{1,\beta}(\ol{\Om_{R/2}})}\le C.
\]

\smallskip
{\rm (b)} In addition to the previous assumptions,
if $\|\Gam^k\|_{C^{1,\sigma}}\le M$, $k=1,2$, for some $\sigma\in(0,1)$, if $(a_{ij}, a)$ satisfy
\[
\|(a_{ij},a)({\bf 0}, 0,\cdot), D^m_{({\bf p},z)}(a_{ij}, a)({\bf p}, z, \cdot)\|_{1,\delta, \Om_R}^{(-\delta),\{\rx_0\}}\le M
\]
for any $({\bf p}, z)$ satisfying $|{\bf p}|+|z|\le 2K$ and for $m=1,2$,
and if each $b^{(k)}$ satisfies
\[
\|b^{(k)}\|_{C^{2,\delta}(\ol{V})}\le M\qquad\tx{for $k=1,2$},
\]
for some $\delta\in(0,1)$, then there exists a constant $C>0$ depending only on $(\lambda, K, M, R, \alp, \sigma, \delta)$
such that $u$ satisfies
\[
\|u\|_{2,\sigma, \Om_{R/2}}^{(-1-\alp),\{\rx_0\}}\le C.
\]
\end{proposition}

\section{Well-Posedness of a Nonlinear Boundary Value Problem}
For a constant $h>0$ and a function $f_{\rm bd}:[0,h]\rightarrow \R_+$,
denote a bounded domain $\Om\subset \R^2$ as
\begin{equation}
\label{definition-app-C-dom}
\Om:=\{{\bf x}=(x_1,x_2)\in\R^2\,:\, x_1\in(0,h),\,\, x_2\in(0, f_{\rm bd}(x_1))\},
\end{equation}
where  $f_{\rm bd}$ satisfies that, for constants $t_0\ge 0$, $t_1>0$, $t_2>0$, $t_h>0$, $\alp\in(0,1)$, and $M>0$,
\begin{equation}
\label{definition-app-C-dom2}
\begin{split}
&f_{\rm bd}\in C^1([0,h]),\quad f_{\rm bd}(0)=t_0,\quad f_{\rm bd}(h)=t_h,\\
&f_{\rm bd}(x_1)\ge \min\{t_1x_1+t_0, t_2\},\\
&\|f_{\rm bd}\|_{2,\alp,(0,h)}^{(-1-\alp),\{0,h\}}\le M.
\end{split}
\end{equation}

We denote the boundary vertices and segments as follows:
\begin{equation}
\label{definition-app-C-dom3}
\begin{split}
&\lefttop=(0, t_0),\quad \righttop=(h,t_h),\quad \rightbottom=(h,0),\quad \leftbottom=(0,0),\\
&\ol{\Gam_{\rm l}}=\der \Om\cap\{x_1=0\},\quad \ol{\Gam_{\rm r}}=\der \Om\cap\{x_1=h\},\\
&\ol{\Gam_{\rm t}}=\der \Om\cap \{x_2=f_{\rm bd}(x_1)\},\quad \ol{\Gam_{\rm b}}=\der\Om\cap\{x_2=0\};
\end{split}
\end{equation}
and $\Gam_{\rm l}$, $\Gam_{\rm r}$, $\Gam_{\rm t}$, and $\Gam_{\rm b}$ are the relative interiors of the segments defined above.

Let $\phi_0({\bf x})$ be a piecewise smooth function defined in $\R^2$ such that

\smallskip
\begin{itemize}
  \item
  $\phi_0\in C^{\infty}(\{x_1\le\frac h3\})\cap C^{\infty}(\{x_1\ge\frac{2h}{3}\})$ with
      $
      \|\phi_0\|_{C^3(\ol{\Om\setminus\{\frac{h}{3}< x_1< \frac{2h}{3}\}})}\le C_{\phi_0}
      $,

\smallskip
  \item $\phi_0\equiv 0$ in $\{x_1\le\frac h4\}$,

\smallskip
  \item $\phi_0$ is linear in $\{x_1\ge \frac{3h}{4}\}$,

\smallskip
  \item
  $\der_{x_2}\phi_0=0\qquad \tx{on}\,\,\Gam_{\rm b}.$
\end{itemize}
Consider a nonlinear boundary value problem{\rm :}
\begin{equation}
\label{appendix-c-nlbvp}
\begin{split}
\sum_{i,j=1}^2 A_{ij}(Du, {\bf x})D_{ij}u+\sum_{i=1}^2 A_i(Du, {\bf x})D_iu=0\qquad&\tx{in $\Om$},\\
u=\phi_0\qquad&\tx{on $\Gam_{\rm l}\cup \Gam_{\rm r}$},\\
B(Du, u, {\bf x})=0\qquad&\tx{on $\Gam_{\rm t}$},\\
\der_{x_2}u=0\qquad&\tx{on $\Gam_{\rm b}$}.
\end{split}
\end{equation}

Assume that \eqref{appendix-c-nlbvp} satisfies that,
for constants $\lambda\in(0,1)$, $M<\infty$, $\alp\in(0,1)$, $\beta\in[\frac 12, 1)$,
$\sigma\in(0,1)$, and $\eps\in(0, \frac{h}{10})$, the following properties hold{\rm :}

\smallskip
(i) For any ${\bf x}\in \ol{\Om}$, and ${\bf p}, {\bm\kappa}=(\kappa_1, \kappa_2)\in \R^2$,
\[
\lambda\, {\rm dist}({\bf x}, \Gam_{\rm l}\cup\Gam_{\rm r})|\bm\kappa|^2\le \sum_{i,j=1}^2 A_{ij}({\bf p}, {\bf x})\kappa_i\kappa_j\le \lambda^{-1}|\bm\kappa|^2.
\]

(ii)
For any ${\bf x}\in \ol{\Om}\setminus\{\frac{\eps}{2}<x_1<h-\frac{\eps}{2}\}$
and ${\bf p}, {\bm\kappa}=(\kappa_1, \kappa_2)\in \R^2$,
\[
\lambda |\bm\kappa|^2\le \sum_{i,j=1}^2 \frac{A_{ij}({\bf p}, {\bf x})\kappa_i\kappa_j}{(\min\{x_1,h-x_1, \delta\})^{2-\frac{i+j}{2}}}
\le \lambda^{-1}|\bm\kappa|^2.
\]

(iii) $(A_{ij}, A_i)({\bf p}, {\bf x})$ are independent of ${\bf p}$ on $\Om\cap\{\eps \le x_1\le h-\eps \}$ with
\[
\|A_{ij}\|_{L^{\infty}(\Om\cap \{\eps\le x_1\le h-\eps\})}+\|(A_{ij}, A_i)\|_{C^{1,\alp}(\ol{\Om\cap\{\eps\le x\le h-\eps\}})}\le M.
\]

(iv) For any ${\bf p}\in \R^2$,
\[
\|(A_{ij}, A_i)({\bf p},\cdot)\|_{C^{\beta}(\ol{\Om\setminus \{2\eps<x_1<h-2\eps\}})}
+\|(D_{\bf p}A_{ij}, D_{\bf p}A_i)({\bf p,\cdot})\|_{L^{\infty}(\Om\setminus \{2\eps<x_1<h-2\eps\})}\le M.
\]

(v) $(A_{ij}, A_i)\in C^{1,\alp}(\R^2\times (\ol{\Om}\setminus \ol{\Gam_{\rm l}\cup \Gam_{\rm r}}))$ and
\[
\|(A_{ij}, A_i)\|_{C^{1,\alp}(\R^2\times (\ol{\Om}\cap\{s\le x_1\le h-s\}))}\le M\big(\frac hs\big)^M
\qquad\tx{for all $s\in (0, \frac{h}{4})$}.
\]

(vi) For each $({\bf p}, {\bf x})\in \R^2\times \ol{\Om}\setminus\{\frac h4\le x_1\le \frac{3h}{4}\}$, define
\begin{equation*}
 \hat{{\bf p}}={\bf p}-D\phi_0({\bf x}),\qquad
 (a_{ij},a_i)(\hat{{\bf p}}, {\bf x})=(A_{ij}, A_i)({\bf p}, {\bf x}).
\end{equation*}
For each $({\bf p}, (x_1,0))\in \R^2\times (\Gam_{\rm b}\setminus \{\eps\le x_1\le h-\eps\})$,
\[
(a_{11}, a_{22}, a_1)((\hat{p}_1,-\hat{p}_2), (x_1,0))=(a_{11}, a_{22}, a_1)((\hat{p}_1, \hat{p}_2), (x_1,0)),
\]
and, for all $({\bf p}, {\bf x})\in \R^2\times (\Om\setminus \{\eps\le x_1\le h-\eps\})$, $i=1,2$,
\begin{equation*}
\begin{split}
&|a_{ii}({\bf p},(x_1,x_2))-a_{ii}(D\phi_0({0,x_2}), (0,x_2))|\le M|x_1|^{\beta}\qquad\qquad\tx{when $x_1<\eps$},\\
&|a_{ii}({\bf p},(x_1,x_2))-a_{ii}(D\phi_0({h,x_2}), (0,x_2))|\le M|x_1-h|^{\beta}\qquad\tx{when $x_1>h-\eps$}.
\end{split}
\end{equation*}
In $\Om\setminus \{\eps\le x_1< h-\eps\}$, $\phi_0$ satisfies
\[
\sum_{i,j=1}^2 A_{ij}(Du, {\bf x})D_{ij}\phi_0+\sum_{i=1}^2 A_i(Du,{\bf x})D_i\phi_0=0,
\]
so that the equation for $u$ in \eqref{appendix-c-nlbvp} is written as an equation for $\hat u=u-\phi_0$ in the form{\rm :}
\[
\sum_{i,j=1}^2 a_{ij}(D\hat{u}, {\bf x})D_{ij}\hat{u}+\sum_{i=1}^2 a_i(D\hat{u},{\bf x})D_i\hat{u}=0.
\]

(vii) For any ${\bf p}\in \R^2$ and ${\bf x}\in \Gam_{\rm l}\cup \Gam_{\rm r}$,
$(A_{12}, A_{21})({\bf p}, {\bf x})=0.
$

\smallskip
(viii) For any ${\bf p}\in \R^2$ and ${\bf x}\in \Om\setminus\{\frac{\eps}{2}\le x_1\le h-\frac{\eps}{2}\}$,
$A_1({\bf p}, {\bf x})\le -\lambda$.

\smallskip
(ix) For any $({\bf p}, z, {\bf x})\in \R^2\times \R\times \Gam_{\rm t}$,
$D_{\bf p}B({\bf p}, z, {\bf x})\cdot{\bm\nu}^{(1)}(\bf x)\ge \lambda$,
where ${\bm\nu}^{(1)}$ is the inner unit normal vector to $\Gam_{\rm t}$ towards the interior of $\Om$;

\smallskip
(x) For any $({\bf p},z)\in \R^2\times \R$,
\begin{align*}
&\|(B(D\phi_0, \phi_0,\cdot)\|_{C^3(\ol{\Om\setminus\{\frac h3<x_1<\frac{2h}{3}\}})}+\| D^k_{({\bf p}, z)}({\bf p}, z, \cdot))\|_{C^3(\ol{\Om})}
\le M\qquad\mbox{for $k=1,2,3$},\\
&\|D_{\bf p} B({\bf p}, z, \cdot)\|_{C^0(\ol{\Om})}\le \lambda^{-1},\\
&D_zB({\bf p}, z, {\bf x})\le -\lambda\qquad\,\,\tx{for all}\,\,{\bf x}\in \Gam_{\rm t},\\
&D_{p_1}B({\bf p}, z, {\bf x})\le -\lambda\qquad\tx{for all}\,\,\Gam_{\rm t}\setminus\{\eps\le x_1\le h-\eps\}.
\end{align*}

(xi) There exist $v\in C^3(\ol{\Gam_{\rm t}})$ and a nonhomogeneous linear operator:
\[
L({\bf p}, z, {\bf x})={\bf b}^{(1)}({\bf x})\cdot {\bf p}+b_0^{(1)}({\bf x}) z+g_1({\bf x}),
\]
defined for ${\bf x}\in \Gam_{\rm t}$ and $({\bf p}, z)\in \R^2\times \R$, satisfying
\[
\|v\|_{C^3(\ol{\Om})}+\|({\bf b}^{(1)}, b_0^{(1)}, g_1)\|_{C^3(\ol{\Gam_{\rm t}})}\le M
\]
such that, for all $({\bf p}, z, {\bf x})\in \R^2\times \R\times \Gam_{\rm t}$,
\begin{align*}
&|B({\bf p}, z, {\bf x})-L({\bf p}, z, {\bf x})|\le \sigma \big(|{\bf p}-Dv({\bf x})|+|z-v({\bf x})|\big),\\
&|D_{\bf p}B({\bf p}, z, {\bf x})-{\bf b}^{(1)}({\bf x})|+|D_zB({\bf p}, z, {\bf x})-b_0^{(1)}({\bf x})|\le \sigma.
\end{align*}
From \cite[Propositions 4.7.2 and 4.8.7]{CF2}, the following two propositions are obtained:

\begin{proposition}
\label{proposition-app-C-wp}
For fixed constants $\lambda>0$, $M<\infty$, $\alp\in(0,1)$, $\beta\in[\frac 12,1 )$,
and $\eps\in(0, \frac{h}{10})$, there exist constants $\alp_1\in(0,\frac 12)$, $\sigma\in(0,1)$,
and $\delta_0>0$ with $\alp_1$ depending only on $\lambda$,
and $(\sigma, \delta_0)$ depending only on $(\lambda, M, C_{\phi_0}, \alp, \beta, \eps)$
such that the following statement holds{\rm :}
let domain $\Om$ be defined by \eqref{definition-app-C-dom},
and let the nonlinear boundary value problem \eqref{appendix-c-nlbvp} satisfy all the conditions
stated above with $h$, $t_h$, $t_1$, $t_2$, $t_0\ge 0$, $\eps\in(0, \frac{h}{10})$,
and $\delta\in[0, \delta_0)$.
Then the boundary value problem \eqref{appendix-c-nlbvp} has a unique
solution $u\in C(\ol{\Om})\cap C^1(\ol{\Om}\setminus (\ol{\Gam_{\rm l}}\cup \ol{\Gam_{\rm r}}))\cap C^2(\Om)$.
Moreover, $u$ satisfies
\begin{equation}
\label{appendix4-estimate-u-a}
\|u\|_{C^0(\ol{\Om})}\le C,\qquad
|u({\bf x})-\phi_0({\bf x})|\le C\min\{x_1, h-x_1\}\,\,\,\,\,\tx{in $\Om$}
\end{equation}
with a constant $C>0$ depending only on $(\lambda, M, C_{\phi_0}, \eps)$.
Furthermore, $u$ is in $C(\ol{\Om})\cap C^{2,\alp_1}(\ol{\Om}\setminus \ol{\Gam_{\rm l}}\cup \ol{\Gam_{\rm r}})$
and satisfies
\begin{equation}
\label{appendix4-estimate-u-b}
\|u\|_{C^{2,\alp_1} (\ol{\Om\cap \{s<x_1<h-s\}})}\le C_s
\end{equation}
for each $s\in(0, \frac{h}{10})$ with a constant $C_s>0$ depending only on $(\lambda, M, C_{\phi_0}, \alp, \beta, \eps, s)$.
\end{proposition}

\begin{proposition}
\label{proposition-app-C-wp2}
For fixed constants $\lambda>0$, $\delta>0$, $M<\infty$, $\alp\in(0,1)$, $\beta\in[\frac 12, 1)$,
and $\eps\in(0, \frac{h}{10})$, there exist constants $\alp_1\in(0, \frac 12)$, $\sigma\in(0,1)$ with $\alp_1$
depending only on $(\lambda, \delta)$, and $\sigma>0$ depending only
on $(\lambda, \delta, M, C_{\phi_0}, \alp, \beta, \eps)$ such that the following statement holds{\rm :}
let domain $\Om$ be of the structure of \eqref{definition-app-C-dom}--\eqref{definition-app-C-dom3}
with $h>0$, $t_h>0$, $t_1\ge 0$, $t_2\ge 0$, and $t_0=0$, that is,
\[
\lefttop=\leftbottom=(0,0),\qquad \ol{\Gam_{\rm l}}=\{(0,0)\},
\]
and let the nonlinear boundary value problem \eqref{appendix-c-nlbvp} satisfy
conditions {\rm{(iii), (v), and (ix)--(xi)}} above, and the following modified conditions{\rm :}

\smallskip
{\rm{(i*)}} For any ${\bf x}\in \ol{\Om}$ and ${\bf p}, \bm\kappa=(\kappa_1, \kappa_2)\in \R^2$,
\begin{align*}
&\min\{\lambda\, {\rm dist}({\bf x}, \Gam_{\rm l})+\delta, \lambda\, {\rm dist}({\bf x}, \Gam_{\rm r})\}|\bm\kappa|^2
\le \sum_{i,j=1}^2 A_{ij}({\bf p, \bf x})\kappa_i\kappa_j
\le \lambda^{-1}|\bm\kappa|^2,\\
&\|(A_{ij}, A_i)(D\phi_0, \cdot), D^m_{{\bf p}}(A_{ij}, A_i)({\bf p}, \cdot)\|_{1,\alp, \Om\cap \{x_1<2\eps\}}^{(-\alp), \{\lefttop\}}\le M
  \qquad\mbox{for $m=1,2$}.
\end{align*}

\smallskip
{\rm{(ii*)}} Condition {\rm (ii)} holds for any ${\bf x}\in \ol{\Om}\cap\{{\rm dist}({\bf x}, \Gam_{\rm r})< \frac{\eps}{2}\}$
and ${\bf p}, \bm\kappa\in\R^2$.

\smallskip
{\rm{(iv*)}} For any ${\bf p}\in \R^2$,
\[
\|(A_{ij}, A_i)({\bf p},\cdot)\|_{C^{\beta}(\ol{\Om\cap \{x_1\ge h-2\eps\}})}
+\|(D_{\bf p}A_{ij}, D_{\bf p}A_i)({\bf p,\cdot})\|_{L^{\infty}(\Om\cap \{x_1>h-2\eps\})}\le M.
\]

{\rm{(vi*)}}
For each $({\bf p}, (x_1,0))\in \R^2\times (\Gam_{\rm b}\cap \{x_1> h-\eps\})$,
\[
(a_{11}, a_{22}, a_1)((\hat{p}_1,-\hat{p}_2), (x_1,0))=(a_{11}, a_{22}, a_1)((\hat{p}_1, \hat{p}_2), (x_1,0)),
\]
and, for all $({\bf p}, {\bf x})\in \R^2\times (\Om\cap \{
 x_1> h-\eps\})$,
\[
|a_{ii}({\bf p},(x_1,x_2))-a_{ii}(D\phi_0({h,x_2}), (0,x_2))|\le M|x_1-h|^{\beta}, \qquad i=1,2.
\]

{\rm{(vii*)}} Condition {\rm{(vii)}} holds for all ${\bf p}\in\R^2$ and ${\bf x}\in \Gam_{\rm r}$.

\smallskip
{\rm{(viii*)}} Condition {\rm{(viii)}} holds for all ${\bf p}\in \R^2$ and ${\bf x}\in \Om\cap\{x_1>h-\frac{\eps}{2}\}$.

\smallskip
\noindent
Then the boundary value problem \eqref{appendix-c-nlbvp} has a unique
solution $u\in C(\ol{\Om})\cap C^1(\ol{\Om}\setminus (\ol{\Gam_l}\cup \ol{\Gam_r}))\cap C^2(\Om)$.
Moreover, solution $u$
is in $C(\ol{\Om})\cap C^{2,\alp_1}(\ol{\Om}\setminus (\ol{\Gam_l}\cup \ol{\Gam_r}))$
and satisfies \eqref{appendix4-estimate-u-a}--\eqref{appendix4-estimate-u-b}
for $C>0$ in \eqref{appendix4-estimate-u-a} depending only on $(\lambda, \delta, M, C_{\phi_0}, \eps)$,
and $C_s>0$ depending on $(\lambda, \delta, M, C_{\phi_0}, \eps, s)$.
Furthermore, $u$ satisfies
\begin{equation*}
\|u\|_{2,\alp_1, \Om\cap\{x_1<\frac h4\}}^{(-1-\alp_1), \{\lefttop\}}\le \hat C
\end{equation*}
for constant $\hat C>0$ depending only on $(\delta, \lambda, M, \alp, \eps)$.
\end{proposition}

\backmatter


\chapter*{Acknowledgments}

The research of Myoungjean Bae was supported in part by Samsung Science and Technology Foundation
under Project Number SSTF-BA1502-02.
The research of Gui-Qiang G. Chen  was supported in part by
the UK Engineering and Physical Sciences Research Council Award
EP/L015811/1, EP/V008854, and EP/V051121/1,
and the Royal Society--Wolfson Research Merit Award (UK).
The research of Mikhail Feldman was
supported in part by the National Science Foundation under Grants
DMS-1764278 and DMS-2054689,
and the Van Vleck Professorship Research Award by the University of Wisconsin-Madison.
The authors would like to thank the anonymous referees for helpful suggestions
to improve the presentation of this monograph.


\bibliographystyle{amsalpha}



\printindex

\end{document}